\documentclass[leqno,12pt]{article}

\usepackage[dvips]{graphics}
\usepackage{color}
\usepackage{ifthen}
\usepackage{soul}
\usepackage{latexsym}
\usepackage{amsthm,amssymb}
\usepackage{amsbsy,amsfonts,amsmath}

\newcommand{\bysame}{\leavevmode\hbox
to 3em{\hrulefill}\,}

\usepackage{a4}

\renewcommand{\thefootnote}{\fnsymbol{footnote}}

\numberwithin{equation}{section}

\renewcommand{\baselinestretch}{1}
\newtheorem{theorem}{Theorem}[section]
\newtheorem{proposition}[theorem]{Proposition}

\newtheorem{lemma}[theorem]{Lemma}
\newtheorem{example}[theorem]{Example}

\theoremstyle{definition}
\newtheorem{definition}[theorem]{Definition}
\newtheorem{notation}[theorem]{Notation}

\theoremstyle{remark}
\newtheorem{remark}[theorem]{Remark}

\newcommand{\al}{\alpha}
\newcommand{\fkJ}{J}

\newcommand{\bN}{{\mathbb{N}}}
\newcommand{\bZ}{{\mathbb{Z}}}
\newcommand{\bZgeqo}{\bZ_{\geq 0}}

\newcommand{\rmid}{{\mathrm{id}}}

\newcommand{\hatpi}{{\hat \pi}}
\newcommand{\hats}{{\hat s}}
\newcommand{\hath}{{\hat h}}
\newcommand{\hatW}{{\hat W}}
\newcommand{\hatS}{{\hat S}}
\newcommand{\hatvarphi}{{\hat \varphi}}
\newcommand{\hattheta}{{\hat \theta}}

\newcommand{\dotW}{{\dot W}}
\newcommand{\dotvarphi}{{\dot \varphi}}
\newcommand{\dottheta}{{\dot \theta}}

\newcommand{\bC}{{\mathbb{C}}}
\newcommand{\bCt}{\bC^\times}

\newcommand{\bR}{{\mathbb{R}}}

\newcommand{\defOmega}{\Omega}
\newcommand{\defzeta}{\zeta}

\newcommand{\curveSouthWestTHICKsmall}{
\setlength{\unitlength}{0.5mm}
\begin{picture}(10,10)(5.5,0)

\put(0,0){\oval(20.8,20.8)[bl]}
\put(0,0){\oval(20.6,20.6)[bl]}
\put(0,0){\oval(20.4,20.4)[bl]}
\put(0,0){\oval(20.2,20.2)[bl]}
\put(0,0){\oval(20,20)[bl]}
\put(0,0){\oval(19.8,19.8)[bl]}
\put(0,0){\oval(19.6,19.6)[bl]}
\put(0,0){\oval(19.4,19.4)[bl]}
\put(0,0){\oval(19.2,19.2)[bl]}

\end{picture}
}

\newcommand{\curveSouthEastTHICKsmall}{
\setlength{\unitlength}{0.5mm}
\begin{picture}(10,10)(5.5,0)

\put(0,0){\oval(20.8,20.8)[br]}
\put(0,0){\oval(20.6,20.6)[br]}
\put(0,0){\oval(20.4,20.4)[br]}
\put(0,0){\oval(20.2,20.2)[br]}
\put(0,0){\oval(20,20)[br]}
\put(0,0){\oval(19.8,19.8)[br]}
\put(0,0){\oval(19.6,19.6)[br]}
\put(0,0){\oval(19.4,19.4)[br]}
\put(0,0){\oval(19.2,19.2)[br]}

\end{picture}
}

\newcommand{\curveNorthWestTHICKsmall}{
\setlength{\unitlength}{0.5mm}
\begin{picture}(10,10)(5.5,0)

\put(0,0){\oval(20.8,20.8)[tl]}
\put(0,0){\oval(20.6,20.6)[tl]}
\put(0,0){\oval(20.4,20.4)[tl]}
\put(0,0){\oval(20.2,20.2)[tl]}
\put(0,0){\oval(20,20)[tl]}
\put(0,0){\oval(19.8,19.8)[tl]}
\put(0,0){\oval(19.6,19.6)[tl]}
\put(0,0){\oval(19.4,19.4)[tl]}
\put(0,0){\oval(19.2,19.2)[tl]}

\end{picture}
}

\newcommand{\curveNorthEastTHICKsmall}{
\setlength{\unitlength}{0.5mm}
\begin{picture}(10,10)(5.5,0)

\put(0,0){\oval(20.8,20.8)[tr]}
\put(0,0){\oval(20.6,20.6)[tr]}
\put(0,0){\oval(20.4,20.4)[tr]}
\put(0,0){\oval(20.2,20.2)[tr]}
\put(0,0){\oval(20,20)[tr]}
\put(0,0){\oval(19.8,19.8)[tr]}
\put(0,0){\oval(19.6,19.6)[tr]}
\put(0,0){\oval(19.4,19.4)[tr]}
\put(0,0){\oval(19.2,19.2)[tr]}

\end{picture}
}

\newcommand{\RIGHThorizontallineeightmmTHICKsmall}{
\setlength{\unitlength}{0.5mm}
\begin{picture}(10,10)(5.5,0)

\put(0,0.4){\line(1,0){8}}
\put(0,0.2){\line(1,0){8}}
\put(0,0){\line(1,0){8}}
\put(0,-0.2){\line(1,0){8}}
\put(0,-0.4){\line(1,0){8}}

\end{picture}
}

\newcommand{\RIGHThorizontallinetwentyeightmmTHICKsmall}{
\setlength{\unitlength}{0.5mm}
\begin{picture}(10,10)(5.5,0)

\put(0,0.4){\line(1,0){28}}
\put(0,0.2){\line(1,0){28}}
\put(0,0){\line(1,0){28}}
\put(0,-0.2){\line(1,0){28}}
\put(0,-0.4){\line(1,0){28}}

\end{picture}
}

\newcommand{\RIGHTlineTHICKsmall}[1]{
\setlength{\unitlength}{0.5mm}
\begin{picture}(10,10)(5.5,0)

\put(1,0.4){\line(1,0){#1}}
\put(1,0.2){\line(1,0){#1}}
\put(1,0){\line(1,0){#1}}
\put(1,-0.2){\line(1,0){#1}}
\put(1,-0.4){\line(1,0){#1}}

\end{picture}
}

\newcommand{\properRIGHTlineTHICKsmall}[1]{
\setlength{\unitlength}{0.5mm}
\begin{picture}(10,10)(5.5,0)

\put(0,0.4){\line(1,0){#1}}
\put(0,0.2){\line(1,0){#1}}
\put(0,0){\line(1,0){#1}}
\put(0,-0.2){\line(1,0){#1}}
\put(0,-0.4){\line(1,0){#1}}

\end{picture}
}

\newcommand{\properLEFTlineTHICKsmall}[1]{
\setlength{\unitlength}{0.5mm}
\begin{picture}(10,10)(5.5,0)

\put(0,0.4){\line(-1,0){#1}}
\put(0,0.2){\line(-1,0){#1}}
\put(0,0){\line(-1,0){#1}}
\put(0,-0.2){\line(-1,0){#1}}
\put(0,-0.4){\line(-1,0){#1}}

\end{picture}
}

\newcommand{\UPverticallineeightmmTHICKsmall}{
\setlength{\unitlength}{0.5mm}
\begin{picture}(10,10)(5.5,0)

\put(0.4,0){\line(0,1){8}}
\put(0.2,0){\line(0,1){8}}
\put(0,0){\line(0,1){8}}
\put(-0.2,0){\line(0,1){8}}
\put(-0.4,0){\line(0,1){8}}

\end{picture}
}

\newcommand{\UPlineTHICKsmall}[1]{
\setlength{\unitlength}{0.5mm}
\begin{picture}(10,10)(5.5,0)

\put(0.4,1){\line(0,1){8}}
\put(0.2,1){\line(0,1){8}}
\put(0,1){\line(0,1){8}}
\put(-0.2,1){\line(0,1){8}}
\put(-0.4,1){\line(0,1){8}}

\end{picture}
}

\newcommand{\properUPlineTHICKsmall}[1]{
\setlength{\unitlength}{0.5mm}
\begin{picture}(10,10)(5.5,0)

\put(0.4,0){\line(0,1){#1}}
\put(0.2,0){\line(0,1){#1}}
\put(0,0){\line(0,1){#1}}
\put(-0.2,0){\line(0,1){#1}}
\put(-0.4,0){\line(0,1){#1}}

\end{picture}
}

\newcommand{\properDOWNlineTHICKsmall}[1]{
\setlength{\unitlength}{0.5mm}
\begin{picture}(10,10)(5.5,0)

\put(0.4,0){\line(0,-1){#1}}
\put(0.2,0){\line(0,-1){#1}}
\put(0,0){\line(0,-1){#1}}
\put(-0.2,0){\line(0,-1){#1}}
\put(-0.4,0){\line(0,-1){#1}}

\end{picture}
}

\newcommand{\UPRIGHTlineeightmmTHICKsmall}{
\setlength{\unitlength}{0.5mm}
\begin{picture}(10,10)(5.5,0)

\put(0.7,1.1){\line(1,1){8.4}}
\put(0.6,1){\line(1,1){8.4}}
\put(0.5,0.9){\line(1,1){8.4}}
\put(0.8,0.8){\line(1,1){8.4}}
\put(0.9,0.7){\line(1,1){8.4}}
\put(01,0.6){\line(1,1){8.4}}
\put(01.1,0.5){\line(1,1){8.4}}

\end{picture}
}

\newcommand{\UPRIGHTlineTHICKsmall}[1]{
\setlength{\unitlength}{0.5mm}
\begin{picture}(10,10)(5.5,0)

\put(0.7,1.1){\line(1,1){#1}}
\put(0.6,1){\line(1,1){#1}}
\put(0.5,0.9){\line(1,1){#1}}
\put(0.8,0.8){\line(1,1){#1}}
\put(0.9,0.7){\line(1,1){#1}}
\put(01,0.6){\line(1,1){#1}}
\put(01.1,0.5){\line(1,1){#1}}

\end{picture}
}

\newcommand{\UPLEFTlineeightmmTHICKsmall}{
\setlength{\unitlength}{0.5mm}
\begin{picture}(10,10)(5.5,0)

\put(-0.5,1.1){\line(-1,1){8.4}}
\put(-0.6,1){\line(-1,1){8.4}}
\put(-0.7,0.9){\line(-1,1){8.4}}
\put(-0.8,0.8){\line(-1,1){8.4}}
\put(-0.9,0.7){\line(-1,1){8.4}}
\put(-1,0.6){\line(-1,1){8.4}}
\put(-1.1,0.5){\line(-1,1){8.4}}

\end{picture}
}

\newcommand{\UPLEFTlineTHICKsmall}[1]{
\setlength{\unitlength}{0.5mm}
\begin{picture}(10,10)(5.5,0)

\put(-0.5,1.1){\line(-1,1){#1}}
\put(-0.6,1){\line(-1,1){#1}}
\put(-0.7,0.9){\line(-1,1){#1}}
\put(-0.8,0.8){\line(-1,1){#1}}
\put(-0.9,0.7){\line(-1,1){#1}}
\put(-1,0.6){\line(-1,1){#1}}
\put(-1.1,0.5){\line(-1,1){#1}}

\end{picture}
}

\newcommand{\HamiltonWeylBthree}{
\setlength{\unitlength}{0.5mm}
\begin{picture}(230,100)(-30,-20)

\put(9,36){{\tiny{$*$}}}

\multiput(10,50)(40,0){4}{\circle{2}}
\multiput(11,50)(40,0){4}{\line(1,0){8}}
\multiput(20,50)(40,0){4}{\circle{2}}
\multiput(21,50)(40,0){3}{\line(1,0){28}}

\multiput(10,40)(40,0){4}{\circle{2}}
\multiput(10,41)(40,0){4}{\line(0,1){8}}
\multiput(11,40)(40,0){4}{\line(1,0){8}}
\multiput(20,40)(40,0){4}{\circle{2}}
\multiput(20,41)(40,0){4}{\line(0,1){8}}

\multiput(00,30)(40,0){4}{\circle{2}}
\multiput(00.8,30.8)(40,0){4}{\line(1,1){8.4}}
\multiput(30,30)(40,0){4}{\circle{2}}
\multiput(29.2,30.8)(40,0){4}{\line(-1,1){8.4}}
\multiput(31,30)(40,0){3}{\line(1,0){8}}

\multiput(00,20)(40,0){4}{\circle{2}}
\multiput(00,21)(40,0){4}{\line(0,1){8}}
\multiput(30,20)(40,0){4}{\circle{2}}
\multiput(30,21)(40,0){4}{\line(0,1){8}}
\multiput(31,20)(40,0){3}{\line(1,0){8}}

\multiput(10,10)(40,0){4}{\circle{2}}
\multiput(9.2,10.8)(40,0){4}{\line(-1,1){8.4}}
\multiput(11,10)(40,0){4}{\line(1,0){8}}
\multiput(20,10)(40,0){4}{\circle{2}}
\multiput(20.8,10.8)(40,0){4}{\line(1,1){8.4}}

\multiput(10,00)(40,0){4}{\circle{2}}
\multiput(10,01)(40,0){4}{\line(0,1){8}}
\multiput(11,00)(40,0){4}{\line(1,0){8}}
\multiput(20,00)(40,0){4}{\circle{2}}
\multiput(20,01)(40,0){4}{\line(0,1){8}}
\multiput(21,00)(40,0){3}{\line(1,0){28}}

\multiput(10,50)(80,0){2}{\RIGHTlineTHICKsmall{8}}
\multiput(20,50)(40,0){3}{\RIGHTlineTHICKsmall{28}}

\multiput(10,40)(80,0){2}{\RIGHTlineTHICKsmall{8}}
\multiput(50,40)(80,0){2}{\UPlineTHICKsmall{8}}
\multiput(60,40)(80,0){2}{\UPlineTHICKsmall{8}}

\multiput(0,30)(40,0){4}{\UPRIGHTlineTHICKsmall{8.4}}
\multiput(30,30)(40,0){4}{\UPLEFTlineTHICKsmall{8.4}}

\multiput(30,20)(40,0){3}{\UPlineTHICKsmall{8}}
\multiput(40,20)(40,0){3}{\UPlineTHICKsmall{8}}

\multiput(10,10)(40,0){4}{\UPLEFTlineTHICKsmall{8.4}}
\multiput(10,10)(120,0){2}{\RIGHTlineTHICKsmall{8}}
\multiput(20,10)(40,0){4}{\UPRIGHTlineTHICKsmall{8.4}}

\multiput(10,0)(120,0){2}{\RIGHTlineTHICKsmall{8}}
\multiput(20,0)(40,0){3}{\RIGHTlineTHICKsmall{28}}
\multiput(50,0)(40,0){2}{\UPlineTHICKsmall{8}}
\multiput(60,0)(40,0){2}{\UPlineTHICKsmall{8}}

\put(130,51){\curveNorthEastTHICKsmall}
\put(20,61){\properRIGHTlineTHICKsmall{110}}
\put(20,51){\curveNorthWestTHICKsmall}

\put(151,40){\curveSouthEastTHICKsmall}
\put(161,60){\properDOWNlineTHICKsmall{20}}
\put(151,60){\curveNorthEastTHICKsmall}
\put(-1,70){\properRIGHTlineTHICKsmall{152}}
\put(-1,60){\curveNorthWestTHICKsmall}
\put(-11,40){\properUPlineTHICKsmall{20}}
\put(-1,40){\curveSouthWestTHICKsmall}

\put(-1,10){\curveNorthWestTHICKsmall}
\put(-11,10){\properDOWNlineTHICKsmall{20}}
\put(-1,-10){\curveSouthWestTHICKsmall}
\put(-1,-20){\properRIGHTlineTHICKsmall{152}}
\put(151,-10){\curveSouthEastTHICKsmall}
\put(161,-10){\properUPlineTHICKsmall{20}}
\put(151,10){\curveNorthEastTHICKsmall}

\put(20,-1){\curveSouthWestTHICKsmall}
\put(20,-11){\properRIGHTlineTHICKsmall{110}}
\put(130,-1){\curveSouthEastTHICKsmall}

\multiput(74,71)(0,0){1}{{\tiny{3}}}

\multiput(74,56.5)(0,0){1}{{\tiny{2}}}

\multiput(34,51)(40,0){3}{{\tiny{2}}}
\multiput(14,51)(40,0){4}{{\tiny{1}}}

\multiput(6.5,44)(40,0){4}{{\tiny{3}}}
\multiput(21,44)(40,0){4}{{\tiny{3}}}

\multiput(14,36)(40,0){4}{{\tiny{1}}}

\multiput(21,32)(40,0){4}{{\tiny{2}}}
\multiput(6,32)(40,0){4}{{\tiny{2}}}

\multiput(34,31)(40,0){3}{{\tiny{3}}}

\multiput(26.5,24)(40,0){4}{{\tiny{1}}}
\multiput(1,24)(40,0){4}{{\tiny{1}}}

\multiput(34,16)(40,0){3}{{\tiny{3}}}

\multiput(6,15)(40,0){4}{{\tiny{2}}}
\multiput(21,15)(40,0){4}{{\tiny{2}}}

\multiput(14,11)(40,0){4}{{\tiny{1}}}

\multiput(6.5,4)(40,0){4}{{\tiny{3}}}
\multiput(21,4)(40,0){4}{{\tiny{3}}}

\multiput(14,-4)(40,0){4}{{\tiny{1}}}
\multiput(34,-4)(40,0){3}{{\tiny{2}}}

\multiput(74,-10)(0,0){1}{{\tiny{2}}}
\multiput(74,-24)(0,0){1}{{\tiny{3}}}

\end{picture}
}

\newcommand{\VisualExBfour}{
\setlength{\unitlength}{1mm}
\begin{picture}(100,20)(-10,0)

\multiput(0,15)(60,0){2}{\circle{10}}

\multiput(-1,4.8)(60,0){2}{\line(0,1){5.4}}
\multiput(1,4.8)(60,0){2}{\line(0,1){5.4}}

\multiput(-4,7)(60,0){2}{{\tiny{$s_4$}}}
\multiput(1.5,7)(60,0){2}{{\tiny{$s_4$}}}

\multiput(0,0)(15,0){6}{\circle{10}}
\multiput(4.8,1)(15,0){5}{\line(1,0){5.4}}
\multiput(4.8,-1)(15,0){5}{\line(1,0){5.4}}

\multiput(6,2)(15,0){5}{{\tiny{$s_4$}}}
\multiput(6,-2.5)(15,0){5}{{\tiny{$s_4$}}}

\put(-2.5,-1.5){${\hat{W}}_1$}
\put(12.5,-1.5){${\hat{W}}_4$}
\put(27.5,-1.5){${\hat{W}}_7$}
\put(42.5,-1.5){${\hat{W}}_6$}
\put(57.5,-1.5){${\hat{W}}_3$}
\put(72.5,-1.5){${\hat{W}}_2$}

\put(-2.5,13.5){${\hat{W}}_0$}
\put(57.5,13.5){${\hat{W}}_5$}

\end{picture}
}

\newcommand{\TotalCHHecThreeSix}{
\setlength{\unitlength}{0.5mm}
\begin{picture}(230,70)(-30,-10)

\multiput(10,40)(40,0){4}{\circle{2}}
\multiput(11,40)(40,0){4}{\line(1,0){8}}
\multiput(20,40)(40,0){4}{\circle{2}}
\multiput(21,40)(40,0){3}{\line(1,0){28}}

\multiput(10,30)(40,0){4}{\circle{2}}
\multiput(10,31)(40,0){4}{\line(0,1){8}}
\multiput(11,30)(40,0){4}{\line(1,0){8}}
\multiput(20,30)(40,0){4}{\circle{2}}
\multiput(20,31)(40,0){4}{\line(0,1){8}}

\multiput(00,20)(40,0){4}{\circle{2}}
\multiput(00.8,20.8)(40,0){4}{\line(1,1){8.4}}
\multiput(30,20)(40,0){4}{\circle{2}}
\multiput(29.2,20.8)(40,0){4}{\line(-1,1){8.4}}
\multiput(31,20)(40,0){3}{\line(1,0){8}}

\multiput(10,10)(40,0){4}{\circle{2}}
\multiput(9.2,10.8)(40,0){4}{\line(-1,1){8.4}}
\multiput(11,10)(40,0){4}{\line(1,0){8}}
\multiput(20,10)(40,0){4}{\circle{2}}
\multiput(20.8,10.8)(40,0){4}{\line(1,1){8.4}}

\multiput(10,00)(40,0){4}{\circle{2}}
\multiput(10,01)(40,0){4}{\line(0,1){8}}
\multiput(11,00)(40,0){4}{\line(1,0){8}}
\multiput(20,00)(40,0){4}{\circle{2}}
\multiput(20,01)(40,0){4}{\line(0,1){8}}
\multiput(21,00)(40,0){3}{\line(1,0){28}}

\multiput(20,40)(40,0){3}{\RIGHTlineTHICKsmall{28}}

\multiput(10,30)(40,0){4}{\UPlineTHICKsmall{8}}
\multiput(20,30)(40,0){4}{\UPlineTHICKsmall{8}}

\multiput(0,20)(40,0){4}{\UPRIGHTlineTHICKsmall{8.4}}
\multiput(30,20)(40,0){4}{\UPLEFTlineTHICKsmall{8.4}}

\multiput(10,10)(40,0){4}{\UPLEFTlineTHICKsmall{8.4}}
\multiput(20,10)(40,0){4}{\UPRIGHTlineTHICKsmall{8.4}}

\multiput(10,0)(40,0){4}{\RIGHTlineTHICKsmall{8}}
\multiput(10,0)(40,0){4}{\UPlineTHICKsmall{8}}
\multiput(20,0)(40,0){4}{\UPlineTHICKsmall{8}}

\put(130,41){\curveNorthEastTHICKsmall}
\put(20,51){\properRIGHTlineTHICKsmall{110}}
\put(20,41){\curveNorthWestTHICKsmall}

\put(151,30){\oval(20,20)[br]}
\put(161,50){\line(0,-1){20}}
\put(151,50){\oval(20,20)[tr]}
\put(-1,60){\line(1,0){152}}
\put(-1,50){\oval(20,20)[tl]}
\put(-11,30){\line(0,1){20}}
\put(-1,30){\oval(20,20)[bl]}

\put(20,-1){\oval(20,20)[bl]}
\put(20,-11){\line(1,0){110}}
\put(130,-1){\oval(20,20)[br]}

\multiput(74,61)(0,0){1}{{\tiny{1}}}

\multiput(74,46.5)(0,0){1}{{\tiny{2}}}

\multiput(34,41)(40,0){3}{{\tiny{2}}}
\multiput(14,41)(40,0){4}{{\tiny{3}}}

\multiput(6.5,34)(40,0){4}{{\tiny{1}}}
\multiput(21,34)(40,0){4}{{\tiny{1}}}

\multiput(14,26)(40,0){4}{{\tiny{3}}}

\multiput(21,22)(40,0){4}{{\tiny{2}}}
\multiput(6,22)(40,0){4}{{\tiny{2}}}

\multiput(34,21)(40,0){3}{{\tiny{1}}}

\multiput(6,15)(40,0){4}{{\tiny{3}}}
\multiput(21,15)(40,0){4}{{\tiny{3}}}

\multiput(14,11)(40,0){4}{{\tiny{2}}}

\multiput(6.5,4)(40,0){4}{{\tiny{1}}}
\multiput(21,4)(40,0){4}{{\tiny{1}}}

\multiput(14,-4)(40,0){4}{{\tiny{2}}}
\multiput(34,-4)(40,0){3}{{\tiny{3}}}

\multiput(74,-10)(0,0){1}{{\tiny{3}}}

\multiput(5,39)(00,0){1}{{\tiny{$a$}}}
\multiput(18,42)(40,0){3}{{\tiny{$a$}}}
\multiput(48,42)(40,0){3}{{\tiny{$a$}}}
\multiput(142,39)(40,0){1}{{\tiny{$a$}}}

\multiput(5,29)(40,0){4}{{\tiny{$b$}}}
\multiput(22,29)(40,0){4}{{\tiny{$b$}}}

\multiput(-2,22)(40,0){4}{{\tiny{$c$}}}
\multiput(29,22)(40,0){4}{{\tiny{$c$}}}

\multiput(9,12)(40,0){4}{{\tiny{$d$}}}
\multiput(18,12)(40,0){4}{{\tiny{$d$}}}

\multiput(5,-1)(40,0){1}{{\tiny{$e$}}}
\multiput(19,-4)(40,0){3}{{\tiny{$e$}}}
\multiput(49,-4)(40,0){3}{{\tiny{$e$}}}
\multiput(142,-1)(40,0){1}{{\tiny{$e$}}}

\end{picture}
}

\newcommand{\spTotalCHHecThreeSix}{
\setlength{\unitlength}{0.5mm}
\begin{picture}(230,70)(-30,-10)

\multiput(10,40)(40,0){4}{\circle{2}}
\multiput(11,40)(40,0){4}{\line(1,0){8}}
\multiput(20,40)(40,0){4}{\circle{2}}
\multiput(21,40)(40,0){3}{\line(1,0){28}}

\multiput(10,30)(40,0){4}{\circle{2}}
\multiput(10,31)(40,0){4}{\line(0,1){8}}
\multiput(11,30)(40,0){4}{\line(1,0){8}}
\multiput(20,30)(40,0){4}{\circle{2}}
\multiput(20,31)(40,0){4}{\line(0,1){8}}

\multiput(00,20)(40,0){4}{\circle{2}}
\multiput(00.8,20.8)(40,0){4}{\line(1,1){8.4}}
\multiput(30,20)(40,0){4}{\circle{2}}
\multiput(29.2,20.8)(40,0){4}{\line(-1,1){8.4}}
\multiput(31,20)(40,0){3}{\line(1,0){8}}

\multiput(10,10)(40,0){4}{\circle{2}}
\multiput(9.2,10.8)(40,0){4}{\line(-1,1){8.4}}
\multiput(11,10)(40,0){4}{\line(1,0){8}}
\multiput(20,10)(40,0){4}{\circle{2}}
\multiput(20.8,10.8)(40,0){4}{\line(1,1){8.4}}

\multiput(10,00)(40,0){4}{\circle{2}}
\multiput(10,01)(40,0){4}{\line(0,1){8}}
\multiput(11,00)(40,0){4}{\line(1,0){8}}
\multiput(20,00)(40,0){4}{\circle{2}}
\multiput(20,01)(40,0){4}{\line(0,1){8}}
\multiput(21,00)(40,0){3}{\line(1,0){28}}

\multiput(10,40)(80,0){2}{\RIGHTlineTHICKsmall{8}}
\multiput(20,40)(40,0){2}{\RIGHTlineTHICKsmall{28}}
\multiput(130,40)(00,0){1}{\RIGHTlineTHICKsmall{8}}

\multiput(10,30)(40,0){2}{\UPlineTHICKsmall{8}}
\multiput(10,30)(80,0){2}{\RIGHTlineTHICKsmall{8}}
\multiput(60,30)(40,0){3}{\UPlineTHICKsmall{8}}
\multiput(130,30)(40,0){1}{\UPlineTHICKsmall{8}}

\multiput(40,20)(40,0){3}{\UPRIGHTlineTHICKsmall{8.4}}
\multiput(30,20)(40,0){2}{\UPLEFTlineTHICKsmall{8.4}}
\multiput(110,20)(80,0){1}{\RIGHTlineTHICKsmall{8}}
\multiput(150,20)(00,0){1}{\UPLEFTlineTHICKsmall{8.4}}

\multiput(10,10)(40,0){3}{\UPLEFTlineTHICKsmall{8.4}}
\multiput(10,10)(40,0){2}{\RIGHTlineTHICKsmall{8}}
\multiput(20,10)(40,0){3}{\UPRIGHTlineTHICKsmall{8.4}}
\multiput(130,10)(00,0){1}{\RIGHTlineTHICKsmall{8}}

\multiput(10,0)(40,0){2}{\RIGHTlineTHICKsmall{8}}
\multiput(20,0)(40,0){3}{\RIGHTlineTHICKsmall{28}}
\multiput(90,0)(40,0){2}{\UPlineTHICKsmall{8}}
\multiput(100,0)(40,0){2}{\UPlineTHICKsmall{8}}

\put(130,41){\oval(20,20)[tr]}
\put(20,51){\line(1,0){110}}
\put(20,41){\oval(20,20)[tl]}

\put(151,30){\curveSouthEastTHICKsmall}
\put(161,30){\properUPlineTHICKsmall{20}}
\put(151,50){\curveNorthEastTHICKsmall}
\put(-1,60){\properRIGHTlineTHICKsmall{152}}
\put(-1,50){\curveNorthWestTHICKsmall}
\put(-11,30){\properUPlineTHICKsmall{20}}
\put(-1,30){\curveSouthWestTHICKsmall}

\put(20,-1){\curveSouthWestTHICKsmall}
\put(20,-11){\properRIGHTlineTHICKsmall{110}}
\put(130,-1){\curveSouthEastTHICKsmall}

\multiput(74,61)(0,0){1}{{\tiny{1}}}

\multiput(74,46.5)(0,0){1}{{\tiny{2}}}

\multiput(34,41)(40,0){3}{{\tiny{2}}}
\multiput(14,41)(40,0){4}{{\tiny{3}}}

\multiput(6.5,34)(40,0){4}{{\tiny{1}}}
\multiput(21,34)(40,0){4}{{\tiny{1}}}

\multiput(14,26)(40,0){4}{{\tiny{3}}}

\multiput(21,22)(40,0){4}{{\tiny{2}}}
\multiput(6,22)(40,0){4}{{\tiny{2}}}

\multiput(34,21)(40,0){3}{{\tiny{1}}}

\multiput(6,15)(40,0){4}{{\tiny{3}}}
\multiput(21,15)(40,0){4}{{\tiny{3}}}

\multiput(14,11)(40,0){4}{{\tiny{2}}}

\multiput(6.5,4)(40,0){4}{{\tiny{1}}}
\multiput(21,4)(40,0){4}{{\tiny{1}}}

\multiput(14,-4)(40,0){4}{{\tiny{2}}}
\multiput(34,-4)(40,0){3}{{\tiny{3}}}

\multiput(74,-10)(0,0){1}{{\tiny{3}}}

\multiput(5,39)(00,0){1}{{\tiny{$a$}}}
\multiput(18,42)(40,0){3}{{\tiny{$a$}}}
\multiput(48,42)(40,0){3}{{\tiny{$a$}}}
\multiput(142,39)(40,0){1}{{\tiny{$a$}}}

\multiput(5,29)(40,0){4}{{\tiny{$b$}}}
\multiput(22,29)(40,0){4}{{\tiny{$b$}}}

\multiput(-2,22)(40,0){4}{{\tiny{$c$}}}
\multiput(29,22)(40,0){4}{{\tiny{$c$}}}

\multiput(9,12)(40,0){4}{{\tiny{$d$}}}
\multiput(18,12)(40,0){4}{{\tiny{$d$}}}

\multiput(5,-1)(40,0){1}{{\tiny{$e$}}}
\multiput(19,-4)(40,0){3}{{\tiny{$e$}}}
\multiput(49,-4)(40,0){3}{{\tiny{$e$}}}
\multiput(142,-1)(40,0){1}{{\tiny{$e$}}}

\end{picture}
}

\newcommand{\JointForCxeter}{
\setlength{\unitlength}{1mm}
\begin{picture}(120,60)(-10,-20)

\multiput(00,20)(20,0){4}{\circle{2}}
\multiput(00,00)(20,0){4}{\circle{2}}

\multiput(-16,20)(2,0){3}{\circle*{0.6}}

\multiput(72,20)(2,0){3}{\circle*{0.6}}

\multiput(-16,00)(2,0){3}{\circle*{0.6}}

\multiput(72,00)(2,0){3}{\circle*{0.6}}

\put(21,0){\line(1,0){18}}\put(21,20){\line(1,0){18}}

\put(-9,20){\line(1,0){8}}\multiput(-9,20)(0,0.05){4}{\line(1,0){8}}\multiput(-9,20)(0,-0.05){4}{\line(1,0){8}}
\put(61,20){\line(1,0){8}}\multiput(61,20)(0,0.05){4}{\line(1,0){8}}\multiput(61,20)(0,-0.05){4}{\line(1,0){8}}

\put(-9,0){\line(1,0){8}}\multiput(-9,0)(0,0.05){4}{\line(1,0){8}}\multiput(-9,0)(0,-0.05){4}{\line(1,0){8}}
\put(61,0){\line(1,0){8}}\multiput(61,0)(0,0.05){4}{\line(1,0){8}}\multiput(61,0)(0,-0.05){4}{\line(1,0){8}}

\put(1,0){\line(1,0){18}}\multiput(1,0)(0,0.05){4}{\line(1,0){18}}\multiput(1,0)(0,-0.05){4}{\line(1,0){18}}
\put(1,20){\line(1,0){18}}\multiput(1,20)(0,0.05){4}{\line(1,0){18}}\multiput(1,20)(0,-0.05){4}{\line(1,0){18}}

\put(20,1){\line(0,1){18}}\multiput(20,1)(0.05,0){4}{\line(0,1){18}}\multiput(20,1)(-0.05,0){4}{\line(0,1){18}}
\put(40,1){\line(0,1){18}}\multiput(40,1)(0.05,0){4}{\line(0,1){18}}\multiput(40,1)(-0.05,0){4}{\line(0,1){18}}

\put(41,0){\line(1,0){18}}\multiput(41,0)(0,0.05){4}{\line(1,0){18}}\multiput(41,0)(0,-0.05){4}{\line(1,0){18}}
\put(41,20){\line(1,0){18}}\multiput(41,20)(0,0.05){4}{\line(1,0){18}}\multiput(41,20)(0,-0.05){4}{\line(1,0){18}}

\multiput(21,32)(00,0){1}{{\footnotesize{$\hats({\hat{\varphi}}_{u,r}([1])))$}}}
\put(28,30){{\rotatebox{-90}{{\large{$=$}}}}}
\multiput(12,9)(00,0){1}{{\footnotesize{$\hats(n)$}}}\put(17,17){{\footnotesize{$v$}}}
\multiput(42,9)(00,0){1}{{\footnotesize{$\hats(n)$}}}

\put(28,-7){{\rotatebox{-90}{{\large{$=$}}}}}
\multiput(21,-16)(00,0){1}{{\footnotesize{$\hats(\varphi_{l-1}({\overline{b_l}})))$}}}

\multiput(01.5,22)(00,0){1}{{\footnotesize{$\hats({\hat{\varphi}}_{u,r}([2])))$}}}
\multiput(27,22)(00,0){1}{{\footnotesize{$\hats(r)$}}}
\multiput(41.5,22)(00,0){1}{{\footnotesize{$\hats({\hat{\varphi}}_{u,r}([\hath])))$}}}

\multiput(-3,-5)(00,0){1}{{\footnotesize{$\hats(\varphi_{l-1}({\overline{b_l-1}})))$}}}
\multiput(27,-5)(00,0){1}{{\footnotesize{$\hats(r)$}}}
\multiput(40,-5)(00,0){1}{{\footnotesize{$\hats(\varphi_{l-1}({\overline{b_l+1}})))$}}}

\put(80,5){\tiny{$\left(
\begin{array}{l}
\mbox{Here $[t]:=\hatpi_\hath(t)$} \\
\mbox{and ${\overline{t}}:=\hatpi_{(l-1)\hath}(t)$.} \\
\mbox{Recall $r\in\fkJ_{1,n-2}$.}
\end{array}
\right)$}}

\end{picture}
}

\newcommand{\JointForSpecialHighRank}{
\setlength{\unitlength}{1mm}
\begin{picture}(120,70)(-10,-15)

\multiput(00,56)(0,2){3}{\circle*{0.2}}
\multiput(15,56)(0,2){3}{\circle*{0.6}}
\multiput(30,56)(0,2){3}{\circle*{0.2}}
\multiput(45,56)(0,2){3}{\circle*{0.2}}

\multiput(00,-15)(0,2){3}{\circle*{0.2}}
\multiput(15,-15)(0,2){3}{\circle*{0.6}}
\multiput(30,-15)(0,2){3}{\circle*{0.2}}
\multiput(45,-15)(0,2){3}{\circle*{0.2}}

\multiput(-16,45)(2,0){3}{\circle*{0.6}}
\multiput(-16,30)(2,0){3}{\circle*{0.6}}
\multiput(-16,15)(2,0){3}{\circle*{0.6}}
\multiput(-16,00)(2,0){3}{\circle*{0.6}}

\multiput(57,45)(2,0){3}{\circle*{0.6}}
\multiput(57,30)(2,0){3}{\circle*{0.6}}
\multiput(57,15)(2,0){3}{\circle*{0.6}}
\multiput(57,00)(2,0){3}{\circle*{0.6}}

\put(0,46){\line(0,1){8}}
\put(15,46){\line(0,1){8}}\multiput(15,46)(0.05,0){4}{\line(0,1){8}}\multiput(15,46)(-0.05,0){4}{\line(0,1){8}}
\put(30,46){\line(0,1){8}}
\put(45,46){\line(0,1){8}}

\multiput(00,45)(15,0){4}{\circle{2}}
\put(-9,45){\line(1,0){8}}\multiput(-9,45)(0,0.05){4}{\line(1,0){8}}\multiput(-9,45)(0,-0.05){4}{\line(1,0){8}}
\put(1,45){\line(1,0){13}}\multiput(1,45)(0,0.05){4}{\line(1,0){13}}\multiput(1,45)(0,-0.05){4}{\line(1,0){13}}
\put(16,45){\line(1,0){13}}
\put(31,45){\line(1,0){13}}\multiput(31,45)(0,0.05){4}{\line(1,0){13}}\multiput(31,45)(0,-0.05){4}{\line(1,0){13}}
\put(46,45){\line(1,0){8}}\multiput(46,45)(0,0.05){4}{\line(1,0){8}}\multiput(46,45)(0,-0.05){4}{\line(1,0){8}}

\put(0,31){\line(0,1){13}}
\put(15,31){\line(0,1){13}}
\put(30,31){\line(0,1){13}}\multiput(30,31)(0.05,0){4}{\line(0,1){13}}\multiput(30,31)(-0.05,0){4}{\line(0,1){13}}
\put(45,31){\line(0,1){13}}

\multiput(00,30)(15,0){4}{\circle{2}}
\put(-9,30){\line(1,0){8}}\multiput(-9,30)(0,0.05){4}{\line(1,0){8}}\multiput(-9,30)(0,-0.05){4}{\line(1,0){8}}
\put(1,30){\line(1,0){13}}\multiput(1,30)(0,0.05){4}{\line(1,0){13}}\multiput(1,30)(0,-0.05){4}{\line(1,0){13}}
\put(16,30){\line(1,0){13}}
\put(31,30){\line(1,0){13}}\multiput(31,30)(0,0.05){4}{\line(1,0){13}}\multiput(31,30)(0,-0.05){4}{\line(1,0){13}}
\put(46,30){\line(1,0){8}}\multiput(46,30)(0,0.05){4}{\line(1,0){8}}\multiput(46,30)(0,-0.05){4}{\line(1,0){8}}

\put(0,16){\line(0,1){13}}
\put(15,16){\line(0,1){13}}\multiput(15,16)(0.05,0){4}{\line(0,1){13}}\multiput(15,16)(-0.05,0){4}{\line(0,1){13}}
\put(30,16){\line(0,1){13}}
\put(45,16){\line(0,1){13}}

\multiput(00,15)(15,0){4}{\circle{2}}
\put(-9,15){\line(1,0){8}}\multiput(-9,15)(0,0.05){4}{\line(1,0){8}}\multiput(-9,15)(0,-0.05){4}{\line(1,0){8}}
\put(1,15){\line(1,0){13}}\multiput(1,15)(0,0.05){4}{\line(1,0){13}}\multiput(1,15)(0,-0.05){4}{\line(1,0){13}}
\put(16,15){\line(1,0){13}}
\put(31,15){\line(1,0){13}}\multiput(31,15)(0,0.05){4}{\line(1,0){13}}\multiput(31,15)(0,-0.05){4}{\line(1,0){13}}
\put(46,15){\line(1,0){8}}\multiput(46,15)(0,0.05){4}{\line(1,0){8}}\multiput(46,15)(0,-0.05){4}{\line(1,0){8}}

\put(0,1){\line(0,1){13}}
\put(15,1){\line(0,1){13}}
\put(30,1){\line(0,1){13}}\multiput(30,1)(0.05,0){4}{\line(0,1){13}}\multiput(30,1)(-0.05,0){4}{\line(0,1){13}}
\put(45,1){\line(0,1){13}}

\multiput(00,00)(15,0){4}{\circle{2}}
\put(-9,0){\line(1,0){8}}\multiput(-9,0)(0,0.05){4}{\line(1,0){8}}\multiput(-9,0)(0,-0.05){4}{\line(1,0){8}}
\put(1,0){\line(1,0){13}}\multiput(1,0)(0,0.05){4}{\line(1,0){13}}\multiput(1,0)(0,-0.05){4}{\line(1,0){13}}
\put(16,0){\line(1,0){13}}
\put(31,0){\line(1,0){13}}\multiput(31,0)(0,0.05){4}{\line(1,0){13}}\multiput(31,0)(0,-0.05){4}{\line(1,0){13}}
\put(46,0){\line(1,0){8}}\multiput(46,0)(0,0.05){4}{\line(1,0){8}}\multiput(46,0)(0,-0.05){4}{\line(1,0){8}}

\put(0,-9){\line(0,1){8}}
\put(15,-9){\line(0,1){8}}\multiput(15,-9)(0.05,0){4}{\line(0,1){8}}\multiput(15,-9)(-0.05,0){4}{\line(0,1){8}}
\put(30,-9){\line(0,1){8}}
\put(45,-9){\line(0,1){8}}

\multiput(3,2)(0,15){4}{{\tiny{${\hat{k}}_{h_2-1}$}}}
\multiput(18,2)(0,15){4}{{\tiny{${\hat{k}}_{h_2}=k$}}}
\multiput(36,2)(0,15){4}{{\tiny{${\hat{k}}_{1}$}}}

\multiput(-4,36.5)(15,0){4}{{\tiny{${\hat{j}}_{h_1}$}}}
\multiput(-3,21.5)(15,0){4}{{\tiny{${\hat{j}}_1$}}}
\multiput(-3,6.5)(15,0){4}{{\tiny{${\hat{j}}_2$}}}

\put(31,31){{\tiny{$v$}}}

\put(70,5){\tiny{$\left(
\begin{array}{l}
\mbox{Here $v:=(\theta^{\chamchi^\prime}_3\circ\hatpi_{h^{\chamchi^\prime}_3})(0)$,} \\
\mbox{${\hat{k}}_t:
=(\varphi_{\theta^{\chamchi^\prime}_2}\circ\hatpi_{h^{\chamchi^\prime}_2})(t)$,
} \\
\mbox{ 
${\hat{j}}_t:=(\varphi_{\theta^{\chamchi^\prime}_1}\circ$
$\hatpi_{h^{\chamchi^\prime}_1})(t)$,
} \\
\mbox{and $h_1:=h^{\chamchi^\prime}_1$, $h_2:=h^{\chamchi^\prime}_2$.}
\end{array}
\right)$}}

\end{picture}
}

\newcommand{\JointForHecFourEighteen}{
\setlength{\unitlength}{0.5mm}
\begin{picture}(120,70)(50,-10)

\multiput(-6,50)(2,0){3}{\circle*{0.6}}
\multiput(10,50)(10,0){2}{\circle{2}}
\multiput(01,50)(10,0){3}{\line(1,0){8}}
\multiput(32,50)(2,0){3}{\circle*{0.6}}
\multiput(10,41)(10,0){2}{\line(0,1){8}}
\multiput(-6,40)(2,0){3}{\circle*{0.6}}
\multiput(10,40)(10,0){2}{\circle{2}}
\multiput(01,40)(10,0){3}{\line(1,0){8}}
\multiput(32,40)(2,0){3}{\circle*{0.6}}

\multiput(00,50)(10,0){3}{\RIGHTlineTHICKsmall{8}}
\multiput(00,40)(10,0){3}{\RIGHTlineTHICKsmall{8}}

\multiput(-6,10)(2,0){3}{\circle*{0.6}}
\multiput(10,10)(10,0){2}{\circle{2}}
\multiput(01,10)(10,0){3}{\line(1,0){8}}
\multiput(32,10)(2,0){3}{\circle*{0.6}}
\multiput(10,01)(10,0){2}{\line(0,1){8}}
\multiput(-6,00)(2,0){3}{\circle*{0.6}}
\multiput(10,00)(10,0){2}{\circle{2}}
\multiput(01,00)(10,0){3}{\line(1,0){8}}
\multiput(32,00)(2,0){3}{\circle*{0.6}}

\multiput(00,10)(20,0){2}{\RIGHTlineTHICKsmall{8}}
\multiput(10,00)(10,0){2}{\UPlineTHICKsmall{8}}
\multiput(00,00)(20,0){2}{\RIGHTlineTHICKsmall{8}}

\multiput(6,44)(00,0){1}{{\tiny{1}}}
\multiput(21,44)(00,0){1}{{\tiny{1}}}

\multiput(14,51)(00,0){1}{{\tiny{${\bf i}$}}}
\multiput(14,36)(00,0){1}{{\tiny{${\bf i}$}}}

\multiput(9,52)(00,0){1}{{\tiny{$x$}}}
\multiput(19,53)(00,0){1}{{\tiny{$y$}}}

\multiput(9,36)(00,0){1}{{\tiny{$z$}}}
\multiput(19,36)(00,0){1}{{\tiny{$u$}}}

\put(15,30){\vector(0,-1){10}}

\multiput(6,4)(00,0){1}{{\tiny{1}}}
\multiput(21,4)(00,0){1}{{\tiny{1}}}

\multiput(14,11)(00,0){1}{{\tiny{${\bf i}$}}}
\multiput(14,-4)(00,0){1}{{\tiny{${\bf i}$}}}

\multiput(9,12)(00,0){1}{{\tiny{$x$}}}
\multiput(19,13)(00,0){1}{{\tiny{$y$}}}

\multiput(9,-4)(00,0){1}{{\tiny{$z$}}}
\multiput(19,-4)(00,0){1}{{\tiny{$u$}}}

\put(50,15){${\footnotesize{
\left(
\begingroup
\renewcommand{\arraystretch}{1.5}
\begin{array}{l} 
\{({\bf i},x,y,z,u),({\bf i},z,u,x,y),({\bf i},y,x,u,z),({\bf i},u,z,y,x)\}
\\
\cap \{({\bf 3},a,a,a,a), ({\bf 4},a,b,a,b), ({\bf 3},b,c,b,c), ({\bf 4},c,c,c,c), 
\\
({\bf 3},d,d,e,e), ({\bf 4},d,f,e,d^\prime),
({\bf 3},f,f,d^\prime,d^\prime), ({\bf 4},c^\prime,c^\prime,c^\prime,c^\prime), 
\\
({\bf 2},b^\prime,b^\prime,a^\prime,a^\prime),  
({\bf 4},b^\prime,b^\prime,a^\prime,a^\prime)\}\ne\emptyset
\end{array}
\endgroup
\right)}}$} 

\end{picture}
}

\newcommand{\spTotalCHHecThreeSixForHecFourEighteen}{
\setlength{\unitlength}{0.5mm}
\begin{picture}(230,70)(-30,-10)

\multiput(10,40)(40,0){4}{\circle{2}}
\multiput(11,40)(40,0){4}{\line(1,0){8}}
\multiput(20,40)(40,0){4}{\circle{2}}
\multiput(21,40)(40,0){3}{\line(1,0){28}}

\multiput(10,30)(40,0){4}{\circle{2}}
\multiput(10,31)(40,0){4}{\line(0,1){8}}
\multiput(11,30)(40,0){4}{\line(1,0){8}}
\multiput(20,30)(40,0){4}{\circle{2}}
\multiput(20,31)(40,0){4}{\line(0,1){8}}

\multiput(00,20)(40,0){4}{\circle{2}}
\multiput(00.8,20.8)(40,0){4}{\line(1,1){8.4}}
\multiput(30,20)(40,0){4}{\circle{2}}
\multiput(29.2,20.8)(40,0){4}{\line(-1,1){8.4}}
\multiput(31,20)(40,0){3}{\line(1,0){8}}

\multiput(10,10)(40,0){4}{\circle{2}}
\multiput(9.2,10.8)(40,0){4}{\line(-1,1){8.4}}
\multiput(11,10)(40,0){4}{\line(1,0){8}}
\multiput(20,10)(40,0){4}{\circle{2}}
\multiput(20.8,10.8)(40,0){4}{\line(1,1){8.4}}

\multiput(10,00)(40,0){4}{\circle{2}}
\multiput(10,01)(40,0){4}{\line(0,1){8}}
\multiput(11,00)(40,0){4}{\line(1,0){8}}
\multiput(20,00)(40,0){4}{\circle{2}}
\multiput(20,01)(40,0){4}{\line(0,1){8}}
\multiput(21,00)(40,0){3}{\line(1,0){28}}

\multiput(10,40)(80,0){2}{\RIGHTlineTHICKsmall{8}}
\multiput(20,40)(40,0){2}{\RIGHTlineTHICKsmall{28}}
\multiput(130,40)(00,0){1}{\RIGHTlineTHICKsmall{8}}

\multiput(10,30)(40,0){2}{\UPlineTHICKsmall{8}}
\multiput(10,30)(80,0){2}{\RIGHTlineTHICKsmall{8}}
\multiput(60,30)(40,0){3}{\UPlineTHICKsmall{8}}
\multiput(130,30)(40,0){1}{\UPlineTHICKsmall{8}}

\multiput(40,20)(40,0){3}{\UPRIGHTlineTHICKsmall{8.4}}
\multiput(30,20)(40,0){2}{\UPLEFTlineTHICKsmall{8.4}}
\multiput(110,20)(80,0){1}{\RIGHTlineTHICKsmall{8}}
\multiput(150,20)(00,0){1}{\UPLEFTlineTHICKsmall{8.4}}

\multiput(10,10)(40,0){3}{\UPLEFTlineTHICKsmall{8.4}}
\multiput(10,10)(40,0){2}{\RIGHTlineTHICKsmall{8}}
\multiput(20,10)(40,0){3}{\UPRIGHTlineTHICKsmall{8.4}}
\multiput(130,10)(00,0){1}{\RIGHTlineTHICKsmall{8}}

\multiput(10,0)(40,0){2}{\RIGHTlineTHICKsmall{8}}
\multiput(20,0)(40,0){3}{\RIGHTlineTHICKsmall{28}}
\multiput(90,0)(40,0){2}{\UPlineTHICKsmall{8}}
\multiput(100,0)(40,0){2}{\UPlineTHICKsmall{8}}

\put(130,41){\oval(20,20)[tr]}
\put(20,51){\line(1,0){110}}
\put(20,41){\oval(20,20)[tl]}

\put(151,30){\curveSouthEastTHICKsmall}
\put(161,30){\properUPlineTHICKsmall{20}}
\put(151,50){\curveNorthEastTHICKsmall}
\put(-1,60){\properRIGHTlineTHICKsmall{152}}
\put(-1,50){\curveNorthWestTHICKsmall}
\put(-11,30){\properUPlineTHICKsmall{20}}
\put(-1,30){\curveSouthWestTHICKsmall}

\put(20,-1){\curveSouthWestTHICKsmall}
\put(20,-11){\properRIGHTlineTHICKsmall{110}}
\put(130,-1){\curveSouthEastTHICKsmall}

\multiput(74,61)(0,0){1}{{\tiny{{\bf 4}}}}

\multiput(74,46.5)(0,0){1}{{\tiny{3}}}

\multiput(34,41)(40,0){2}{{\tiny{{\bf 3}}}}\multiput(114,41)(00,0){1}{{\tiny{3}}}
\multiput(14,41)(40,0){4}{{\tiny{2}}}

\multiput(6.5,34)(40,0){2}{{\tiny{{\bf 4}}}}\multiput(86.5,34)(00,0){1}{{\tiny{4}}}\multiput(126.5,34)(00,0){1}{{\tiny{{\bf 4}}}}
\multiput(21,34)(00,0){1}{{\tiny{4}}}\multiput(61,34)(40,0){3}{{\tiny{{\bf 4}}}}

\multiput(14,26)(40,0){4}{{\tiny{2}}}

\multiput(21,22)(40,0){2}{{\tiny{{\bf 3}}}}\multiput(101,22)(00,0){1}{{\tiny{3}}}\multiput(141,22)(00,0){1}{{\tiny{{\bf 3}}}}
\multiput(6,22)(00,0){1}{{\tiny{3}}}\multiput(46,22)(40,0){3}{{\tiny{{\bf 3}}}}

\multiput(34,21)(40,0){3}{{\tiny{4}}}

\multiput(6,15)(40,0){4}{{\tiny{2}}}
\multiput(21,15)(40,0){4}{{\tiny{2}}}

\multiput(14,11)(40,0){2}{{\tiny{{\bf 3}}}}\multiput(94,11)(00,0){1}{{\tiny{3}}}\multiput(134,11)(00,0){1}{{\tiny{{\bf 3}}}}

\multiput(6.5,4)(40,0){2}{{\tiny{4}}}\multiput(86.5,4)(40,0){2}{{\tiny{{\bf 4}}}}
\multiput(21,4)(40,0){2}{{\tiny{4}}}\multiput(101,4)(40,0){2}{{\tiny{{\bf 4}}}}

\multiput(14,-4)(40,0){2}{{\tiny{{\bf 3}}}}\multiput(94,-4)(40,0){2}{{\tiny{3}}}
\multiput(34,-4)(40,0){3}{{\tiny{2}}}

\multiput(74,-10)(0,0){1}{{\tiny{2}}}

\multiput(5,39)(00,0){1}{{\tiny{$a$}}}
\multiput(18,42)(40,0){3}{{\tiny{$a$}}}
\multiput(48,42)(40,0){3}{{\tiny{$a$}}}
\multiput(142,39)(40,0){1}{{\tiny{$a$}}}

\multiput(5,29)(40,0){4}{{\tiny{$b$}}}
\multiput(22,29)(40,0){4}{{\tiny{$b$}}}

\multiput(-2,22)(40,0){4}{{\tiny{$c$}}}
\multiput(29,22)(40,0){4}{{\tiny{$c$}}}

\multiput(9,12)(40,0){4}{{\tiny{$d$}}}
\multiput(18,12)(40,0){4}{{\tiny{$d$}}}

\multiput(5,-1)(40,0){1}{{\tiny{$f$}}}
\multiput(19,-5)(40,0){3}{{\tiny{$f$}}}
\multiput(49,-5)(40,0){3}{{\tiny{$f$}}}
\multiput(142,-1)(40,0){1}{{\tiny{$f$}}}

\end{picture}
}

\newcommand{\TotalCHHecThreeSeven}{
\setlength{\unitlength}{0.5mm}
\begin{picture}(230,110)(00,-20)

\multiput(10,70)(40,0){6}{\circle{2}}
\multiput(11,70)(40,0){6}{\line(1,0){8}}
\multiput(20,70)(40,0){6}{\circle{2}}
\multiput(21,70)(40,0){5}{\line(1,0){28}}

\multiput(10,60)(40,0){6}{\circle{2}}
\multiput(10,61)(40,0){6}{\line(0,1){8}}
\multiput(11,60)(40,0){6}{\line(1,0){8}}
\multiput(20,60)(40,0){6}{\circle{2}}
\multiput(20,61)(40,0){6}{\line(0,1){8}}

\multiput(00,50)(40,0){6}{\circle{2}}
\multiput(00.8,50.8)(40,0){6}{\line(1,1){8.4}}
\multiput(30,50)(40,0){6}{\circle{2}}
\multiput(29.2,50.8)(40,0){6}{\line(-1,1){8.4}}
\multiput(31,50)(40,0){5}{\line(1,0){8}}

\multiput(10,40)(40,0){6}{\circle{2}}
\multiput(9.2,40.8)(40,0){6}{\line(-1,1){8.4}}
\multiput(11,40)(40,0){6}{\line(1,0){8}}
\multiput(20,40)(40,0){6}{\circle{2}}
\multiput(20.8,40.8)(40,0){6}{\line(1,1){8.4}}

\multiput(10,30)(40,0){6}{\circle{2}}
\multiput(10,31)(40,0){6}{\line(0,1){8}}
\multiput(11,30)(40,0){6}{\line(1,0){8}}
\multiput(20,30)(40,0){6}{\circle{2}}
\multiput(20,31)(40,0){6}{\line(0,1){8}}

\multiput(00,20)(40,0){6}{\circle{2}}
\multiput(00.8,20.8)(40,0){6}{\line(1,1){8.4}}
\multiput(30,20)(40,0){6}{\circle{2}}
\multiput(29.2,20.8)(40,0){6}{\line(-1,1){8.4}}
\multiput(31,20)(40,0){5}{\line(1,0){8}}

\multiput(10,10)(40,0){6}{\circle{2}}
\multiput(9.2,10.8)(40,0){6}{\line(-1,1){8.4}}
\multiput(11,10)(40,0){6}{\line(1,0){8}}
\multiput(20,10)(40,0){6}{\circle{2}}
\multiput(20.8,10.8)(40,0){6}{\line(1,1){8.4}}

\multiput(10,00)(40,0){6}{\circle{2}}
\multiput(10,01)(40,0){6}{\line(0,1){8}}
\multiput(11,00)(40,0){6}{\line(1,0){8}}
\multiput(20,00)(40,0){6}{\circle{2}}
\multiput(20,01)(40,0){6}{\line(0,1){8}}
\multiput(21,00)(40,0){5}{\line(1,0){28}}

\multiput(20,70)(40,0){5}{\RIGHTlineTHICKsmall{28}}

\multiput(10,60)(40,0){6}{\UPlineTHICKsmall{8}}
\multiput(20,60)(40,0){6}{\UPlineTHICKsmall{8}}

\multiput(0,50)(40,0){6}{\UPRIGHTlineTHICKsmall{8.4}}
\multiput(30,50)(40,0){6}{\UPLEFTlineTHICKsmall{8.4}}

\multiput(10,40)(40,0){6}{\UPLEFTlineTHICKsmall{8.4}}
\multiput(20,40)(40,0){6}{\UPRIGHTlineTHICKsmall{8.4}}

\multiput(10,30)(40,0){6}{\UPlineTHICKsmall{8}}
\multiput(20,30)(40,0){6}{\UPlineTHICKsmall{8}}

\multiput(0,20)(40,0){6}{\UPRIGHTlineTHICKsmall{8.4}}
\multiput(30,20)(40,0){6}{\UPLEFTlineTHICKsmall{8.4}}

\multiput(10,10)(40,0){6}{\UPLEFTlineTHICKsmall{8.4}}
\multiput(20,10)(40,0){6}{\UPRIGHTlineTHICKsmall{8.4}}

\multiput(10,0)(40,0){6}{\RIGHTlineTHICKsmall{8}}
\multiput(10,0)(40,0){6}{\UPlineTHICKsmall{8}}
\multiput(20,0)(40,0){6}{\UPlineTHICKsmall{8}}

\put(210,71){\curveNorthEastTHICKsmall}
\put(20,81){\properRIGHTlineTHICKsmall{190}}
\put(20,71){\curveNorthWestTHICKsmall}

\put(231,60){\oval(20,20)[br]}
\put(241,80){\line(0,-1){20}}
\put(231,80){\oval(20,20)[tr]}
\put(-1,90){\line(1,0){232}}
\put(-1,80){\oval(20,20)[tl]}
\put(-11,60){\line(0,1){20}}
\put(-1,60){\oval(20,20)[bl]}

\put(231,10){\oval(20,20)[tr]}
\put(241,10){\line(0,-1){20}}
\put(231,-10){\oval(20,20)[br]}
\put(-1,-20){\line(1,0){232}}
\put(-1,-10){\oval(20,20)[bl]}
\put(-11,-10){\line(0,1){20}}
\put(-1,10){\oval(20,20)[tl]}

\put(20,-1){\oval(20,20)[bl]}
\put(20,-11){\line(1,0){190}}
\put(210,-1){\oval(20,20)[br]}

\multiput(114,91)(0,0){1}{{\tiny{1}}}

\multiput(114,76.5)(0,0){1}{{\tiny{2}}}

\multiput(34,71)(40,0){5}{{\tiny{2}}}
\multiput(14,71)(40,0){6}{{\tiny{3}}}

\multiput(6.5,63.5)(40,0){6}{{\tiny{1}}}
\multiput(21,63.5)(40,0){6}{{\tiny{1}}}

\multiput(14,61)(40,0){6}{{\tiny{3}}}

\multiput(21,52)(40,0){6}{{\tiny{2}}}
\multiput(6,52)(40,0){6}{{\tiny{2}}}

\multiput(34,46)(40,0){5}{{\tiny{1}}}

\multiput(6,45)(40,0){6}{{\tiny{3}}}
\multiput(21,45)(40,0){6}{{\tiny{3}}}

\multiput(14,36)(40,0){6}{{\tiny{2}}}

\multiput(6.5,34)(40,0){6}{{\tiny{1}}}
\multiput(21,34)(40,0){6}{{\tiny{1}}}

\multiput(14,26)(40,0){6}{{\tiny{2}}}

\multiput(21,22)(40,0){6}{{\tiny{3}}}
\multiput(6,22)(40,0){6}{{\tiny{3}}}

\multiput(34,16)(40,0){5}{{\tiny{1}}}

\multiput(6,15)(40,0){6}{{\tiny{2}}}
\multiput(21,15)(40,0){6}{{\tiny{2}}}

\multiput(14,6)(40,0){6}{{\tiny{3}}}

\multiput(6.5,4)(40,0){6}{{\tiny{1}}}
\multiput(21,4)(40,0){6}{{\tiny{1}}}

\multiput(14,-4)(40,0){6}{{\tiny{3}}}
\multiput(34,-4)(40,0){5}{{\tiny{2}}}

\multiput(114,-10)(0,0){1}{{\tiny{2}}}
\multiput(114,-25)(0,0){1}{{\tiny{1}}}

\multiput(5,69)(40,0){1}{{\tiny{$a$}}}
\multiput(19,72)(40,0){5}{{\tiny{$a$}}}
\multiput(49,72)(40,0){5}{{\tiny{$a$}}}
\multiput(222,69)(40,0){1}{{\tiny{$a$}}}

\multiput(5,59)(40,0){6}{{\tiny{$b$}}}
\multiput(22,59)(40,0){6}{{\tiny{$b$}}}

\multiput(-2,52)(40,0){6}{{\tiny{$c$}}}
\multiput(29,52)(40,0){6}{{\tiny{$c$}}}

\multiput(10,42)(40,0){6}{{\tiny{$d$}}}
\multiput(17,42)(40,0){6}{{\tiny{$d$}}}

\multiput(5,29)(40,0){6}{{\tiny{$d$}}}
\multiput(22,29)(40,0){6}{{\tiny{$d$}}}

\multiput(-2,22)(40,0){6}{{\tiny{$c$}}}
\multiput(29,22)(40,0){6}{{\tiny{$c$}}}

\multiput(9,12)(40,0){6}{{\tiny{$b$}}}
\multiput(18,12)(40,0){6}{{\tiny{$b$}}}

\multiput(5,-1)(40,0){1}{{\tiny{$a$}}}
\multiput(19,-4)(40,0){5}{{\tiny{$a$}}}
\multiput(49,-4)(40,0){5}{{\tiny{$a$}}}
\multiput(222,-1)(40,0){1}{{\tiny{$a$}}}

\end{picture}
}

\newcommand{\sspTotalCHHecThreeSeven}{
\setlength{\unitlength}{0.5mm}
\begin{picture}(230,110)(00,-20)

\multiput(10,70)(40,0){6}{\circle{2}}
\multiput(11,70)(40,0){6}{\line(1,0){8}}
\multiput(20,70)(40,0){6}{\circle{2}}
\multiput(21,70)(40,0){5}{\line(1,0){28}}

\multiput(10,60)(40,0){6}{\circle{2}}
\multiput(10,61)(40,0){6}{\line(0,1){8}}
\multiput(11,60)(40,0){6}{\line(1,0){8}}
\multiput(20,60)(40,0){6}{\circle{2}}
\multiput(20,61)(40,0){6}{\line(0,1){8}}

\multiput(00,50)(40,0){6}{\circle{2}}
\multiput(00.8,50.8)(40,0){6}{\line(1,1){8.4}}
\multiput(30,50)(40,0){6}{\circle{2}}
\multiput(29.2,50.8)(40,0){6}{\line(-1,1){8.4}}
\multiput(31,50)(40,0){5}{\line(1,0){8}}

\multiput(10,40)(40,0){6}{\circle{2}}
\multiput(9.2,40.8)(40,0){6}{\line(-1,1){8.4}}
\multiput(11,40)(40,0){6}{\line(1,0){8}}
\multiput(20,40)(40,0){6}{\circle{2}}
\multiput(20.8,40.8)(40,0){6}{\line(1,1){8.4}}

\multiput(10,30)(40,0){6}{\circle{2}}
\multiput(10,31)(40,0){6}{\line(0,1){8}}
\multiput(11,30)(40,0){6}{\line(1,0){8}}
\multiput(20,30)(40,0){6}{\circle{2}}
\multiput(20,31)(40,0){6}{\line(0,1){8}}

\multiput(00,20)(40,0){6}{\circle{2}}
\multiput(00.8,20.8)(40,0){6}{\line(1,1){8.4}}
\multiput(30,20)(40,0){6}{\circle{2}}
\multiput(29.2,20.8)(40,0){6}{\line(-1,1){8.4}}
\multiput(31,20)(40,0){5}{\line(1,0){8}}

\multiput(10,10)(40,0){6}{\circle{2}}
\multiput(9.2,10.8)(40,0){6}{\line(-1,1){8.4}}
\multiput(11,10)(40,0){6}{\line(1,0){8}}
\multiput(20,10)(40,0){6}{\circle{2}}
\multiput(20.8,10.8)(40,0){6}{\line(1,1){8.4}}

\multiput(10,00)(40,0){6}{\circle{2}}
\multiput(10,01)(40,0){6}{\line(0,1){8}}
\multiput(11,00)(40,0){6}{\line(1,0){8}}
\multiput(20,00)(40,0){6}{\circle{2}}
\multiput(20,01)(40,0){6}{\line(0,1){8}}
\multiput(21,00)(40,0){5}{\line(1,0){28}}

\multiput(20,70)(40,0){5}{\RIGHTlineTHICKsmall{28}}
\multiput(50,70)(80,0){3}{\RIGHTlineTHICKsmall{8}}

\multiput(10,60)(80,0){3}{\UPlineTHICKsmall{8}}
\multiput(20,60)(80,0){3}{\UPlineTHICKsmall{8}}
\multiput(50,60)(80,0){3}{\RIGHTlineTHICKsmall{8}}

\multiput(0,50)(40,0){6}{\UPRIGHTlineTHICKsmall{8.4}}
\multiput(30,50)(40,0){6}{\UPLEFTlineTHICKsmall{8.4}}
\multiput(30,50)(80,0){3}{\RIGHTlineTHICKsmall{8}}

\multiput(10,40)(80,0){3}{\UPLEFTlineTHICKsmall{8.4}}
\multiput(10,40)(40,0){6}{\RIGHTlineTHICKsmall{8}}
\multiput(60,40)(80,0){3}{\UPRIGHTlineTHICKsmall{8.4}}

\multiput(10,30)(40,0){6}{\RIGHTlineTHICKsmall{8}}
\multiput(20,30)(80,0){3}{\UPlineTHICKsmall{8}}
\multiput(50,30)(80,0){3}{\UPlineTHICKsmall{8}}

\multiput(0,20)(80,0){3}{\UPRIGHTlineTHICKsmall{8.4}}
\multiput(30,20)(80,0){3}{\RIGHTlineTHICKsmall{8}}
\multiput(70,20)(80,0){3}{\UPLEFTlineTHICKsmall{8.4}}

\multiput(10,10)(40,0){6}{\UPLEFTlineTHICKsmall{8.4}}
\multiput(20,10)(40,0){6}{\UPRIGHTlineTHICKsmall{8.4}}
\multiput(50,10)(40,0){5}{\RIGHTlineTHICKsmall{8}}

\multiput(10,0)(00,0){1}{\UPlineTHICKsmall{8}}
\multiput(20,0)(00,0){1}{\UPlineTHICKsmall{8}}
\multiput(20,0)(40,0){5}{\RIGHTlineTHICKsmall{28}}
\multiput(50,0)(40,0){5}{\RIGHTlineTHICKsmall{8}}

\put(210,71){\curveNorthEastTHICKsmall}
\put(20,81){\properRIGHTlineTHICKsmall{190}}
\put(20,71){\curveNorthWestTHICKsmall}

\put(231,60){\oval(20,20)[br]}
\put(241,80){\line(0,-1){20}}
\put(231,80){\oval(20,20)[tr]}
\put(-1,90){\line(1,0){232}}
\put(-1,80){\oval(20,20)[tl]}
\put(-11,60){\line(0,1){20}}
\put(-1,60){\oval(20,20)[bl]}

\put(231,10){\oval(20,20)[tr]}
\put(241,10){\line(0,-1){20}}
\put(231,-10){\oval(20,20)[br]}
\put(-1,-20){\line(1,0){232}}
\put(-1,-10){\oval(20,20)[bl]}
\put(-11,-10){\line(0,1){20}}
\put(-1,10){\oval(20,20)[tl]}

\put(20,-1){\curveSouthWestTHICKsmall}
\put(20,-11){\properRIGHTlineTHICKsmall{190}}
\put(210,-1){\curveSouthEastTHICKsmall}

\multiput(114,91)(0,0){1}{{\tiny{1}}}

\multiput(114,76.5)(0,0){1}{{\tiny{2}}}

\multiput(34,71)(40,0){5}{{\tiny{2}}}
\multiput(14,71)(40,0){6}{{\tiny{3}}}

\multiput(6.5,63.5)(40,0){6}{{\tiny{1}}}
\multiput(21,63.5)(40,0){6}{{\tiny{1}}}

\multiput(14,61)(40,0){6}{{\tiny{3}}}

\multiput(21,52)(40,0){6}{{\tiny{2}}}
\multiput(6,52)(40,0){6}{{\tiny{2}}}

\multiput(34,46)(40,0){5}{{\tiny{1}}}

\multiput(6,45)(40,0){6}{{\tiny{3}}}
\multiput(21,45)(40,0){6}{{\tiny{3}}}

\multiput(14,36)(40,0){6}{{\tiny{2}}}

\multiput(6.5,34)(40,0){6}{{\tiny{1}}}
\multiput(21,34)(40,0){6}{{\tiny{1}}}

\multiput(14,26)(40,0){6}{{\tiny{2}}}

\multiput(21,22)(40,0){6}{{\tiny{3}}}
\multiput(6,22)(40,0){6}{{\tiny{3}}}

\multiput(34,16)(40,0){5}{{\tiny{1}}}

\multiput(6,15)(40,0){6}{{\tiny{2}}}
\multiput(21,15)(40,0){6}{{\tiny{2}}}

\multiput(14,6)(40,0){6}{{\tiny{3}}}

\multiput(6.5,4)(40,0){6}{{\tiny{1}}}
\multiput(21,4)(40,0){6}{{\tiny{1}}}

\multiput(14,-4)(40,0){6}{{\tiny{3}}}
\multiput(34,-4)(40,0){5}{{\tiny{2}}}

\multiput(114,-10)(0,0){1}{{\tiny{2}}}
\multiput(114,-25)(0,0){1}{{\tiny{1}}}

\multiput(5,69)(40,0){1}{{\tiny{$a$}}}
\multiput(19,72)(40,0){5}{{\tiny{$a$}}}
\multiput(49,72)(40,0){5}{{\tiny{$a$}}}
\multiput(222,69)(40,0){1}{{\tiny{$a$}}}

\multiput(5,59)(40,0){6}{{\tiny{$b$}}}
\multiput(22,59)(40,0){6}{{\tiny{$b$}}}

\multiput(-2,52)(40,0){6}{{\tiny{$c$}}}
\multiput(29,52)(40,0){6}{{\tiny{$c$}}}

\multiput(10,42)(40,0){6}{{\tiny{$d$}}}
\multiput(17,42)(40,0){6}{{\tiny{$d$}}}

\multiput(5,29)(40,0){6}{{\tiny{$d$}}}
\multiput(22,29)(40,0){6}{{\tiny{$d$}}}

\multiput(-2,22)(40,0){6}{{\tiny{$c$}}}
\multiput(29,22)(40,0){6}{{\tiny{$c$}}}

\multiput(9,12)(40,0){6}{{\tiny{$b$}}}
\multiput(18,12)(40,0){6}{{\tiny{$b$}}}

\multiput(5,-1)(40,0){1}{{\tiny{$a$}}}
\multiput(19,-4)(40,0){5}{{\tiny{$a$}}}
\multiput(49,-4)(40,0){5}{{\tiny{$a$}}}
\multiput(222,-1)(40,0){1}{{\tiny{$a$}}}

\end{picture}
}

\newcommand{\spTotalCHHecThreeSeven}{
\setlength{\unitlength}{0.5mm}
\begin{picture}(230,110)(00,-20)

\multiput(10,70)(40,0){6}{\circle{2}}
\multiput(11,70)(40,0){6}{\line(1,0){8}}
\multiput(20,70)(40,0){6}{\circle{2}}
\multiput(21,70)(40,0){5}{\line(1,0){28}}

\multiput(10,60)(40,0){6}{\circle{2}}
\multiput(10,61)(40,0){6}{\line(0,1){8}}
\multiput(11,60)(40,0){6}{\line(1,0){8}}
\multiput(20,60)(40,0){6}{\circle{2}}
\multiput(20,61)(40,0){6}{\line(0,1){8}}

\multiput(00,50)(40,0){6}{\circle{2}}
\multiput(00.8,50.8)(40,0){6}{\line(1,1){8.4}}
\multiput(30,50)(40,0){6}{\circle{2}}
\multiput(29.2,50.8)(40,0){6}{\line(-1,1){8.4}}
\multiput(31,50)(40,0){5}{\line(1,0){8}}

\multiput(10,40)(40,0){6}{\circle{2}}
\multiput(9.2,40.8)(40,0){6}{\line(-1,1){8.4}}
\multiput(11,40)(40,0){6}{\line(1,0){8}}
\multiput(20,40)(40,0){6}{\circle{2}}
\multiput(20.8,40.8)(40,0){6}{\line(1,1){8.4}}

\multiput(10,30)(40,0){6}{\circle{2}}
\multiput(10,31)(40,0){6}{\line(0,1){8}}
\multiput(11,30)(40,0){6}{\line(1,0){8}}
\multiput(20,30)(40,0){6}{\circle{2}}
\multiput(20,31)(40,0){6}{\line(0,1){8}}

\multiput(00,20)(40,0){6}{\circle{2}}
\multiput(00.8,20.8)(40,0){6}{\line(1,1){8.4}}
\multiput(30,20)(40,0){6}{\circle{2}}
\multiput(29.2,20.8)(40,0){6}{\line(-1,1){8.4}}
\multiput(31,20)(40,0){5}{\line(1,0){8}}

\multiput(10,10)(40,0){6}{\circle{2}}
\multiput(9.2,10.8)(40,0){6}{\line(-1,1){8.4}}
\multiput(11,10)(40,0){6}{\line(1,0){8}}
\multiput(20,10)(40,0){6}{\circle{2}}
\multiput(20.8,10.8)(40,0){6}{\line(1,1){8.4}}

\multiput(10,00)(40,0){6}{\circle{2}}
\multiput(10,01)(40,0){6}{\line(0,1){8}}
\multiput(11,00)(40,0){6}{\line(1,0){8}}
\multiput(20,00)(40,0){6}{\circle{2}}
\multiput(20,01)(40,0){6}{\line(0,1){8}}
\multiput(21,00)(40,0){5}{\line(1,0){28}}

\multiput(20,70)(40,0){5}{\RIGHTlineTHICKsmall{28}}

\multiput(10,60)(40,0){6}{\UPlineTHICKsmall{8}}
\multiput(20,60)(40,0){6}{\UPlineTHICKsmall{8}}
\multiput(90,60)(40,0){4}{\RIGHTlineTHICKsmall{8}}

\multiput(0,50)(40,0){2}{\UPRIGHTlineTHICKsmall{8.4}}
\multiput(30,50)(40,0){2}{\UPLEFTlineTHICKsmall{8.4}}
\multiput(70,50)(40,0){4}{\RIGHTlineTHICKsmall{8}}

\multiput(10,40)(40,0){2}{\RIGHTlineTHICKsmall{8}}
\multiput(20,40)(00,0){1}{\UPRIGHTlineTHICKsmall{8.4}}
\multiput(50,40)(40,0){5}{\UPLEFTlineTHICKsmall{8.4}}
\multiput(100,40)(40,0){4}{\UPRIGHTlineTHICKsmall{8.4}}
\multiput(130,40)(40,0){3}{\RIGHTlineTHICKsmall{8}}

\multiput(10,30)(80,0){2}{\UPlineTHICKsmall{8}}
\multiput(10,30)(40,0){2}{\RIGHTlineTHICKsmall{8}}
\multiput(60,30)(40,0){2}{\UPlineTHICKsmall{8}}
\multiput(130,30)(40,0){3}{\RIGHTlineTHICKsmall{8}}

\multiput(40,20)(40,0){5}{\UPRIGHTlineTHICKsmall{8.4}}
\multiput(30,20)(00,0){1}{\UPLEFTlineTHICKsmall{8.4}}
\multiput(70,20)(40,0){4}{\RIGHTlineTHICKsmall{8}}
\multiput(110,20)(40,0){4}{\UPLEFTlineTHICKsmall{8.4}}

\multiput(10,10)(40,0){2}{\UPLEFTlineTHICKsmall{8.4}}
\multiput(20,10)(40,0){2}{\UPRIGHTlineTHICKsmall{8.4}}
\multiput(50,10)(40,0){5}{\RIGHTlineTHICKsmall{8}}

\multiput(10,0)(00,0){1}{\UPlineTHICKsmall{8}}
\multiput(20,0)(00,0){1}{\UPlineTHICKsmall{8}}
\multiput(20,0)(40,0){5}{\RIGHTlineTHICKsmall{28}}
\multiput(50,00)(00,0){1}{\RIGHTlineTHICKsmall{8}}
\multiput(90,0)(40,0){4}{\UPlineTHICKsmall{8}}
\multiput(100,0)(40,0){4}{\UPlineTHICKsmall{8}}

\put(210,71){\curveNorthEastTHICKsmall}
\put(20,81){\properRIGHTlineTHICKsmall{190}}
\put(20,71){\curveNorthWestTHICKsmall}

\put(231,60){\curveSouthEastTHICKsmall}
\put(241,60){\properUPlineTHICKsmall{20}}
\put(231,80){\curveNorthEastTHICKsmall}
\put(-1,90){\properRIGHTlineTHICKsmall{232}}
\put(-1,80){\curveNorthWestTHICKsmall}
\put(-11,60){\properUPlineTHICKsmall{20}}
\put(-1,60){\curveSouthWestTHICKsmall}

\put(231,10){\curveNorthEastTHICKsmall}
\put(241,-10){\properUPlineTHICKsmall{20}}
\put(231,-10){\curveSouthEastTHICKsmall}
\put(-1,-20){\properRIGHTlineTHICKsmall{232}}
\put(-1,-10){\curveSouthWestTHICKsmall}
\put(-11,-10){\properUPlineTHICKsmall{20}}
\put(-1,10){\curveNorthWestTHICKsmall}

\put(20,-1){\curveSouthWestTHICKsmall}
\put(20,-11){\properRIGHTlineTHICKsmall{190}}
\put(210,-1){\curveSouthEastTHICKsmall}

\multiput(114,91)(0,0){1}{{\tiny{1}}}

\multiput(114,76.5)(0,0){1}{{\tiny{2}}}

\multiput(34,71)(40,0){5}{{\tiny{2}}}
\multiput(14,71)(40,0){6}{{\tiny{3}}}

\multiput(6.5,63.5)(40,0){6}{{\tiny{1}}}
\multiput(21,63.5)(40,0){6}{{\tiny{1}}}

\multiput(14,61)(40,0){6}{{\tiny{3}}}

\multiput(21,52)(40,0){6}{{\tiny{2}}}
\multiput(6,52)(40,0){6}{{\tiny{2}}}

\multiput(34,46)(40,0){5}{{\tiny{1}}}

\multiput(6,45)(40,0){6}{{\tiny{3}}}
\multiput(21,45)(40,0){6}{{\tiny{3}}}

\multiput(14,36)(40,0){6}{{\tiny{2}}}

\multiput(6.5,34)(40,0){6}{{\tiny{1}}}
\multiput(21,34)(40,0){6}{{\tiny{1}}}

\multiput(14,26)(40,0){6}{{\tiny{2}}}

\multiput(21,22)(40,0){6}{{\tiny{3}}}
\multiput(6,22)(40,0){6}{{\tiny{3}}}

\multiput(34,16)(40,0){5}{{\tiny{1}}}

\multiput(6,15)(40,0){6}{{\tiny{2}}}
\multiput(21,15)(40,0){6}{{\tiny{2}}}

\multiput(14,6)(40,0){6}{{\tiny{3}}}

\multiput(6.5,4)(40,0){6}{{\tiny{1}}}
\multiput(21,4)(40,0){6}{{\tiny{1}}}

\multiput(14,-4)(40,0){6}{{\tiny{3}}}
\multiput(34,-4)(40,0){5}{{\tiny{2}}}

\multiput(114,-10)(0,0){1}{{\tiny{2}}}
\multiput(114,-25)(0,0){1}{{\tiny{1}}}

\multiput(5,69)(40,0){1}{{\tiny{$a$}}}
\multiput(19,72)(40,0){5}{{\tiny{$a$}}}
\multiput(49,72)(40,0){5}{{\tiny{$a$}}}
\multiput(222,69)(40,0){1}{{\tiny{$a$}}}

\multiput(5,59)(40,0){6}{{\tiny{$b$}}}
\multiput(22,59)(40,0){6}{{\tiny{$b$}}}

\multiput(-2,52)(40,0){6}{{\tiny{$c$}}}
\multiput(29,52)(40,0){6}{{\tiny{$c$}}}

\multiput(10,42)(40,0){6}{{\tiny{$d$}}}
\multiput(17,42)(40,0){6}{{\tiny{$d$}}}

\multiput(5,29)(40,0){6}{{\tiny{$d$}}}
\multiput(22,29)(40,0){6}{{\tiny{$d$}}}

\multiput(-2,22)(40,0){6}{{\tiny{$c$}}}
\multiput(29,22)(40,0){6}{{\tiny{$c$}}}

\multiput(9,12)(40,0){6}{{\tiny{$b$}}}
\multiput(18,12)(40,0){6}{{\tiny{$b$}}}

\multiput(5,-1)(40,0){1}{{\tiny{$a$}}}
\multiput(19,-4)(40,0){5}{{\tiny{$a$}}}
\multiput(49,-4)(40,0){5}{{\tiny{$a$}}}
\multiput(222,-1)(40,0){1}{{\tiny{$a$}}}

\end{picture}
}

\newcommand{\TotalCHHecThreeSevenForHCHecfourthpartFourTwentytwo}{
\setlength{\unitlength}{0.5mm}
\begin{picture}(230,110)(00,-20)

\multiput(10,70)(40,0){6}{\circle{2}}
\multiput(11,70)(40,0){6}{\line(1,0){8}}
\multiput(20,70)(40,0){6}{\circle{2}}
\multiput(21,70)(40,0){5}{\line(1,0){28}}

\multiput(10,60)(40,0){6}{\circle{2}}
\multiput(10,61)(40,0){6}{\line(0,1){8}}
\multiput(11,60)(40,0){6}{\line(1,0){8}}
\multiput(20,60)(40,0){6}{\circle{2}}
\multiput(20,61)(40,0){6}{\line(0,1){8}}

\multiput(00,50)(40,0){6}{\circle{2}}
\multiput(00.8,50.8)(40,0){6}{\line(1,1){8.4}}
\multiput(30,50)(40,0){6}{\circle{2}}
\multiput(29.2,50.8)(40,0){6}{\line(-1,1){8.4}}
\multiput(31,50)(40,0){5}{\line(1,0){8}}

\multiput(10,40)(40,0){6}{\circle{2}}
\multiput(9.2,40.8)(40,0){6}{\line(-1,1){8.4}}
\multiput(11,40)(40,0){6}{\line(1,0){8}}
\multiput(20,40)(40,0){6}{\circle{2}}
\multiput(20.8,40.8)(40,0){6}{\line(1,1){8.4}}

\multiput(10,30)(40,0){6}{\circle{2}}
\multiput(10,31)(40,0){6}{\line(0,1){8}}
\multiput(11,30)(40,0){6}{\line(1,0){8}}
\multiput(20,30)(40,0){6}{\circle{2}}
\multiput(20,31)(40,0){6}{\line(0,1){8}}

\multiput(00,20)(40,0){6}{\circle{2}}
\multiput(00.8,20.8)(40,0){6}{\line(1,1){8.4}}
\multiput(30,20)(40,0){6}{\circle{2}}
\multiput(29.2,20.8)(40,0){6}{\line(-1,1){8.4}}
\multiput(31,20)(40,0){5}{\line(1,0){8}}

\multiput(10,10)(40,0){6}{\circle{2}}
\multiput(9.2,10.8)(40,0){6}{\line(-1,1){8.4}}
\multiput(11,10)(40,0){6}{\line(1,0){8}}
\multiput(20,10)(40,0){6}{\circle{2}}
\multiput(20.8,10.8)(40,0){6}{\line(1,1){8.4}}

\multiput(10,00)(40,0){6}{\circle{2}}
\multiput(10,01)(40,0){6}{\line(0,1){8}}
\multiput(11,00)(40,0){6}{\line(1,0){8}}
\multiput(20,00)(40,0){6}{\circle{2}}
\multiput(20,01)(40,0){6}{\line(0,1){8}}
\multiput(21,00)(40,0){5}{\line(1,0){28}}

\multiput(20,70)(40,0){5}{\RIGHTlineTHICKsmall{28}}

\multiput(10,60)(40,0){6}{\UPlineTHICKsmall{8}}
\multiput(20,60)(40,0){6}{\UPlineTHICKsmall{8}}
\multiput(90,60)(40,0){4}{\RIGHTlineTHICKsmall{8}}

\multiput(0,50)(40,0){2}{\UPRIGHTlineTHICKsmall{8.4}}
\multiput(30,50)(40,0){2}{\UPLEFTlineTHICKsmall{8.4}}
\multiput(70,50)(40,0){4}{\RIGHTlineTHICKsmall{8}}

\multiput(10,40)(40,0){2}{\RIGHTlineTHICKsmall{8}}
\multiput(20,40)(00,0){1}{\UPRIGHTlineTHICKsmall{8.4}}
\multiput(50,40)(40,0){5}{\UPLEFTlineTHICKsmall{8.4}}
\multiput(100,40)(40,0){4}{\UPRIGHTlineTHICKsmall{8.4}}
\multiput(130,40)(40,0){3}{\RIGHTlineTHICKsmall{8}}

\multiput(10,30)(80,0){2}{\UPlineTHICKsmall{8}}
\multiput(10,30)(40,0){2}{\RIGHTlineTHICKsmall{8}}
\multiput(60,30)(40,0){2}{\UPlineTHICKsmall{8}}
\multiput(130,30)(40,0){3}{\RIGHTlineTHICKsmall{8}}

\multiput(40,20)(40,0){5}{\UPRIGHTlineTHICKsmall{8.4}}
\multiput(30,20)(00,0){1}{\UPLEFTlineTHICKsmall{8.4}}
\multiput(70,20)(40,0){4}{\RIGHTlineTHICKsmall{8}}
\multiput(110,20)(40,0){4}{\UPLEFTlineTHICKsmall{8.4}}

\multiput(10,10)(40,0){2}{\UPLEFTlineTHICKsmall{8.4}}
\multiput(20,10)(40,0){2}{\UPRIGHTlineTHICKsmall{8.4}}
\multiput(50,10)(40,0){5}{\RIGHTlineTHICKsmall{8}}

\multiput(10,0)(00,0){1}{\UPlineTHICKsmall{8}}
\multiput(20,0)(00,0){1}{\UPlineTHICKsmall{8}}
\multiput(20,0)(40,0){5}{\RIGHTlineTHICKsmall{28}}
\multiput(50,00)(00,0){1}{\RIGHTlineTHICKsmall{8}}
\multiput(90,0)(40,0){4}{\UPlineTHICKsmall{8}}
\multiput(100,0)(40,0){4}{\UPlineTHICKsmall{8}}

\put(210,71){\curveNorthEastTHICKsmall}
\put(20,81){\properRIGHTlineTHICKsmall{190}}
\put(20,71){\curveNorthWestTHICKsmall}

\put(231,60){\curveSouthEastTHICKsmall}
\put(241,60){\properUPlineTHICKsmall{20}}
\put(231,80){\curveNorthEastTHICKsmall}
\put(-1,90){\properRIGHTlineTHICKsmall{232}}
\put(-1,80){\curveNorthWestTHICKsmall}
\put(-11,60){\properUPlineTHICKsmall{20}}
\put(-1,60){\curveSouthWestTHICKsmall}

\put(231,10){\curveNorthEastTHICKsmall}
\put(241,-10){\properUPlineTHICKsmall{20}}
\put(231,-10){\curveSouthEastTHICKsmall}
\put(-1,-20){\properRIGHTlineTHICKsmall{232}}
\put(-1,-10){\curveSouthWestTHICKsmall}
\put(-11,-10){\properUPlineTHICKsmall{20}}
\put(-1,10){\curveNorthWestTHICKsmall}

\put(20,-1){\curveSouthWestTHICKsmall}
\put(20,-11){\properRIGHTlineTHICKsmall{190}}
\put(210,-1){\curveSouthEastTHICKsmall}

\multiput(114,91)(0,0){1}{{\tiny{4}}}

\multiput(114,76.5)(0,0){1}{{\bf{{\tiny{3}}}}}

\multiput(34,71)(40,0){5}{{\bf{{\tiny{3}}}}}
\multiput(14,71)(40,0){6}{{\tiny{2}}}

\multiput(6.5,63.5)(40,0){6}{{\tiny{4}}}
\multiput(21,63.5)(40,0){6}{{\tiny{4}}}

\multiput(14,61)(40,0){2}{{\tiny{2}}}\multiput(94,61)(40,0){4}{{\bf{{\tiny{2}}}}}

\multiput(6,52)(40,0){2}{{\bf{{\tiny{3}}}}}\multiput(86,52)(40,0){4}{{\tiny{3}}}
\multiput(21,52)(40,0){2}{{\bf{{\tiny{3}}}}}\multiput(101,52)(40,0){4}{{\tiny{3}}}

\multiput(34,46)(40,0){5}{{\tiny{4}}}

\multiput(6,45)(40,0){1}{{\tiny{2}}}\multiput(46,45)(40,0){5}{{\bf{{\tiny{2}}}}}
\multiput(21,45)(40,0){1}{{\bf{{\tiny{2}}}}}\multiput(61,45)(40,0){1}{{\tiny{2}}}\multiput(101,45)(40,0){4}{{\bf{{\tiny{2}}}}}

\multiput(14,36)(40,0){2}{{\bf{{\tiny{3}}}}}\multiput(94,36)(40,0){1}{{\tiny{3}}}\multiput(134,36)(40,0){3}{{\bf{{\tiny{3}}}}}

\multiput(6.5,34)(40,0){6}{{\tiny{4}}}
\multiput(21,34)(40,0){6}{{\tiny{4}}}

\multiput(14,26)(40,0){2}{{\bf{{\tiny{3}}}}}\multiput(94,26)(40,0){1}{{\tiny{3}}}\multiput(134,26)(40,0){3}{{\bf{{\tiny{3}}}}}

\multiput(6,22)(40,0){1}{{\tiny{2}}}\multiput(46,22)(40,0){5}{{\bf{{\tiny{2}}}}}
\multiput(21,22)(40,0){1}{{\bf{{\tiny{2}}}}}\multiput(61,22)(40,0){1}{{\tiny{2}}}\multiput(101,22)(40,0){4}{{\bf{{\tiny{2}}}}}

\multiput(34,16)(40,0){5}{{\tiny{4}}}

\multiput(6,15)(40,0){2}{{\bf{{\tiny{3}}}}}\multiput(86,15)(40,0){4}{{\tiny{3}}}
\multiput(21,15)(40,0){2}{{\bf{{\tiny{3}}}}}\multiput(101,15)(40,0){4}{{\tiny{3}}}

\multiput(14,6)(40,0){1}{{\tiny{2}}}\multiput(54,6)(40,0){5}{{\bf{{\tiny{2}}}}}

\multiput(6.5,4)(40,0){6}{{\tiny{4}}}
\multiput(21,4)(40,0){6}{{\tiny{4}}}

\multiput(14,-4)(40,0){6}{{\tiny{2}}}
\multiput(34,-4)(40,0){5}{{\bf{{\tiny{3}}}}}

\multiput(114,-10)(0,0){1}{{\bf{{\tiny{3}}}}}
\multiput(114,-25)(0,0){1}{{\tiny{4}}}

\multiput(4,69)(40,0){1}{{\tiny{$a^\prime$}}}
\multiput(19,72)(40,0){5}{{\tiny{$a^\prime$}}}
\multiput(49,72)(40,0){5}{{\tiny{$a^\prime$}}}
\multiput(222,69)(40,0){1}{{\tiny{$a^\prime$}}}

\multiput(4,59)(40,0){6}{{\tiny{$b^\prime$}}}
\multiput(22,59)(40,0){6}{{\tiny{$b^\prime$}}}

\multiput(-2,52)(40,0){6}{{\tiny{$c^\prime$}}}
\multiput(29,52)(40,0){6}{{\tiny{$c^\prime$}}}

\multiput(9.5,42)(40,0){6}{{\tiny{$d^\prime$}}}
\multiput(16.5,42)(40,0){6}{{\tiny{$d^\prime$}}}

\multiput(4,29)(40,0){6}{{\tiny{$d^\prime$}}}
\multiput(22,29)(40,0){6}{{\tiny{$d^\prime$}}}

\multiput(-2,22)(40,0){6}{{\tiny{$c^\prime$}}}
\multiput(29,22)(40,0){6}{{\tiny{$c^\prime$}}}

\multiput(9.5,12)(40,0){6}{{\tiny{$b^\prime$}}}
\multiput(17,12)(40,0){6}{{\tiny{$b^\prime$}}}

\multiput(4,-1)(40,0){1}{{\tiny{$a^\prime$}}}
\multiput(19,-5)(40,0){5}{{\tiny{$a^\prime$}}}
\multiput(49,-5)(40,0){5}{{\tiny{$a^\prime$}}}
\multiput(222,-1)(40,0){1}{{\tiny{$a^\prime$}}}

\end{picture}
}

\newcommand{\TotalCHHecThreeSevenForHCHecfirstpartFourTwentytwo}{
\setlength{\unitlength}{0.5mm}
\begin{picture}(230,110)(00,-20)

\multiput(10,70)(40,0){6}{\circle{2}}
\multiput(11,70)(40,0){6}{\line(1,0){8}}
\multiput(20,70)(40,0){6}{\circle{2}}
\multiput(21,70)(40,0){5}{\line(1,0){28}}

\multiput(10,60)(40,0){6}{\circle{2}}
\multiput(10,61)(40,0){6}{\line(0,1){8}}
\multiput(11,60)(40,0){6}{\line(1,0){8}}
\multiput(20,60)(40,0){6}{\circle{2}}
\multiput(20,61)(40,0){6}{\line(0,1){8}}

\multiput(00,50)(40,0){6}{\circle{2}}
\multiput(00.8,50.8)(40,0){6}{\line(1,1){8.4}}
\multiput(30,50)(40,0){6}{\circle{2}}
\multiput(29.2,50.8)(40,0){6}{\line(-1,1){8.4}}
\multiput(31,50)(40,0){5}{\line(1,0){8}}

\multiput(10,40)(40,0){6}{\circle{2}}
\multiput(9.2,40.8)(40,0){6}{\line(-1,1){8.4}}
\multiput(11,40)(40,0){6}{\line(1,0){8}}
\multiput(20,40)(40,0){6}{\circle{2}}
\multiput(20.8,40.8)(40,0){6}{\line(1,1){8.4}}

\multiput(10,30)(40,0){6}{\circle{2}}
\multiput(10,31)(40,0){6}{\line(0,1){8}}
\multiput(11,30)(40,0){6}{\line(1,0){8}}
\multiput(20,30)(40,0){6}{\circle{2}}
\multiput(20,31)(40,0){6}{\line(0,1){8}}

\multiput(00,20)(40,0){6}{\circle{2}}
\multiput(00.8,20.8)(40,0){6}{\line(1,1){8.4}}
\multiput(30,20)(40,0){6}{\circle{2}}
\multiput(29.2,20.8)(40,0){6}{\line(-1,1){8.4}}
\multiput(31,20)(40,0){5}{\line(1,0){8}}

\multiput(10,10)(40,0){6}{\circle{2}}
\multiput(9.2,10.8)(40,0){6}{\line(-1,1){8.4}}
\multiput(11,10)(40,0){6}{\line(1,0){8}}
\multiput(20,10)(40,0){6}{\circle{2}}
\multiput(20.8,10.8)(40,0){6}{\line(1,1){8.4}}

\multiput(10,00)(40,0){6}{\circle{2}}
\multiput(10,01)(40,0){6}{\line(0,1){8}}
\multiput(11,00)(40,0){6}{\line(1,0){8}}
\multiput(20,00)(40,0){6}{\circle{2}}
\multiput(20,01)(40,0){6}{\line(0,1){8}}
\multiput(21,00)(40,0){5}{\line(1,0){28}}

\multiput(20,70)(40,0){5}{\RIGHTlineTHICKsmall{28}}

\multiput(10,60)(40,0){6}{\UPlineTHICKsmall{8}}
\multiput(20,60)(40,0){6}{\UPlineTHICKsmall{8}}
\multiput(90,60)(40,0){4}{\RIGHTlineTHICKsmall{8}}

\multiput(0,50)(40,0){2}{\UPRIGHTlineTHICKsmall{8.4}}
\multiput(30,50)(40,0){2}{\UPLEFTlineTHICKsmall{8.4}}
\multiput(70,50)(40,0){4}{\RIGHTlineTHICKsmall{8}}

\multiput(10,40)(40,0){2}{\RIGHTlineTHICKsmall{8}}
\multiput(20,40)(00,0){1}{\UPRIGHTlineTHICKsmall{8.4}}
\multiput(50,40)(40,0){5}{\UPLEFTlineTHICKsmall{8.4}}
\multiput(100,40)(40,0){4}{\UPRIGHTlineTHICKsmall{8.4}}
\multiput(130,40)(40,0){3}{\RIGHTlineTHICKsmall{8}}

\multiput(10,30)(80,0){2}{\UPlineTHICKsmall{8}}
\multiput(10,30)(40,0){2}{\RIGHTlineTHICKsmall{8}}
\multiput(60,30)(40,0){2}{\UPlineTHICKsmall{8}}
\multiput(130,30)(40,0){3}{\RIGHTlineTHICKsmall{8}}

\multiput(40,20)(40,0){5}{\UPRIGHTlineTHICKsmall{8.4}}
\multiput(30,20)(00,0){1}{\UPLEFTlineTHICKsmall{8.4}}
\multiput(70,20)(40,0){4}{\RIGHTlineTHICKsmall{8}}
\multiput(110,20)(40,0){4}{\UPLEFTlineTHICKsmall{8.4}}

\multiput(10,10)(40,0){2}{\UPLEFTlineTHICKsmall{8.4}}
\multiput(20,10)(40,0){2}{\UPRIGHTlineTHICKsmall{8.4}}
\multiput(50,10)(40,0){5}{\RIGHTlineTHICKsmall{8}}

\multiput(10,0)(00,0){1}{\UPlineTHICKsmall{8}}
\multiput(20,0)(00,0){1}{\UPlineTHICKsmall{8}}
\multiput(20,0)(40,0){5}{\RIGHTlineTHICKsmall{28}}
\multiput(50,00)(00,0){1}{\RIGHTlineTHICKsmall{8}}
\multiput(90,0)(40,0){4}{\UPlineTHICKsmall{8}}
\multiput(100,0)(40,0){4}{\UPlineTHICKsmall{8}}

\put(210,71){\curveNorthEastTHICKsmall}
\put(20,81){\properRIGHTlineTHICKsmall{190}}
\put(20,71){\curveNorthWestTHICKsmall}

\put(231,60){\curveSouthEastTHICKsmall}
\put(241,60){\properUPlineTHICKsmall{20}}
\put(231,80){\curveNorthEastTHICKsmall}
\put(-1,90){\properRIGHTlineTHICKsmall{232}}
\put(-1,80){\curveNorthWestTHICKsmall}
\put(-11,60){\properUPlineTHICKsmall{20}}
\put(-1,60){\curveSouthWestTHICKsmall}

\put(231,10){\curveNorthEastTHICKsmall}
\put(241,-10){\properUPlineTHICKsmall{20}}
\put(231,-10){\curveSouthEastTHICKsmall}
\put(-1,-20){\properRIGHTlineTHICKsmall{232}}
\put(-1,-10){\curveSouthWestTHICKsmall}
\put(-11,-10){\properUPlineTHICKsmall{20}}
\put(-1,10){\curveNorthWestTHICKsmall}

\put(20,-1){\curveSouthWestTHICKsmall}
\put(20,-11){\properRIGHTlineTHICKsmall{190}}
\put(210,-1){\curveSouthEastTHICKsmall}

\multiput(114,91)(0,0){1}{{\tiny{2}}}

\multiput(114,76.5)(0,0){1}{{\bf{{\tiny{3}}}}}

\multiput(34,71)(40,0){5}{{\bf{{\tiny{3}}}}}
\multiput(14,71)(40,0){6}{{\tiny{4}}}

\multiput(6.5,63.5)(40,0){6}{{\tiny{2}}}
\multiput(21,63.5)(40,0){6}{{\tiny{2}}}

\multiput(14,61)(40,0){2}{{\tiny{4}}}\multiput(94,61)(40,0){4}{{\bf{{\tiny{4}}}}}

\multiput(21,52)(40,0){2}{{\bf{{\tiny{3}}}}}\multiput(101,52)(40,0){4}{{\tiny{3}}}
\multiput(6,52)(40,0){2}{{\bf{{\tiny{3}}}}}\multiput(86,52)(40,0){4}{{\tiny{3}}}

\multiput(34,46)(40,0){5}{{\tiny{2}}}

\multiput(6,45)(40,0){1}{{\tiny{4}}}\multiput(46,45)(40,0){5}{{\bf{{\tiny{4}}}}}
\multiput(21,45)(40,0){1}{{\bf{{\tiny{4}}}}}\multiput(61,45)(40,0){1}{{\tiny{4}}}\multiput(101,45)(40,0){4}{{\bf{{\tiny{4}}}}}

\multiput(14,36)(40,0){2}{{\bf{{\tiny{3}}}}}\multiput(94,36)(40,0){1}{{\tiny{3}}}\multiput(134,36)(40,0){3}{{\bf{{\tiny{3}}}}}

\multiput(6.5,34)(40,0){6}{{\tiny{2}}}
\multiput(21,34)(40,0){6}{{\tiny{2}}}

\multiput(14,26)(40,0){2}{{\bf{{\tiny{3}}}}}\multiput(94,26)(40,0){1}{{\tiny{3}}}\multiput(134,26)(40,0){3}{{\bf{{\tiny{3}}}}}

\multiput(6,22)(40,0){1}{{\tiny{4}}}\multiput(46,22)(40,0){5}{{\bf{{\tiny{4}}}}}
\multiput(21,22)(40,0){1}{{\bf{{\tiny{4}}}}}\multiput(61,22)(40,0){1}{{\tiny{4}}}\multiput(101,22)(40,0){4}{{\bf{{\tiny{4}}}}}

\multiput(34,16)(40,0){5}{{\tiny{2}}}

\multiput(6,15)(40,0){2}{{\bf{{\tiny{3}}}}}\multiput(86,15)(40,0){4}{{\tiny{3}}}
\multiput(21,15)(40,0){2}{{\bf{{\tiny{3}}}}}\multiput(101,15)(40,0){4}{{\tiny{3}}}

\multiput(14,6)(40,0){1}{{\tiny{4}}}\multiput(54,6)(40,0){5}{{\bf{{\tiny{4}}}}}

\multiput(6.5,4)(40,0){6}{{\tiny{2}}}
\multiput(21,4)(40,0){6}{{\tiny{2}}}

\multiput(14,-4)(40,0){1}{{\tiny{4}}}\multiput(54,-4)(40,0){1}{{\bf{{\tiny{4}}}}}\multiput(94,-4)(40,0){4}{{\tiny{4}}}
\multiput(34,-4)(40,0){5}{{\bf{{\tiny{3}}}}}

\multiput(114,-10)(0,0){1}{{\bf{{\tiny{3}}}}}
\multiput(114,-25)(0,0){1}{{\tiny{2}}}

\multiput(5,69)(40,0){1}{{\tiny{$a$}}}
\multiput(19,72)(40,0){5}{{\tiny{$a$}}}
\multiput(49,72)(40,0){5}{{\tiny{$a$}}}
\multiput(222,69)(40,0){1}{{\tiny{$a$}}}

\multiput(5,59)(40,0){6}{{\tiny{$b$}}}
\multiput(22,59)(40,0){6}{{\tiny{$b$}}}

\multiput(-2,52)(40,0){6}{{\tiny{$c$}}}
\multiput(29,52)(40,0){6}{{\tiny{$c$}}}

\multiput(10,42)(40,0){6}{{\tiny{$d$}}}
\multiput(17,42)(40,0){6}{{\tiny{$d$}}}

\multiput(5,29)(40,0){6}{{\tiny{$d$}}}
\multiput(22,29)(40,0){6}{{\tiny{$d$}}}

\multiput(-2,22)(40,0){6}{{\tiny{$c$}}}
\multiput(29,22)(40,0){6}{{\tiny{$c$}}}

\multiput(9,12)(40,0){6}{{\tiny{$b$}}}
\multiput(18,12)(40,0){6}{{\tiny{$b$}}}

\multiput(5,-1)(40,0){1}{{\tiny{$a$}}}
\multiput(19,-4)(40,0){5}{{\tiny{$a$}}}
\multiput(49,-4)(40,0){5}{{\tiny{$a$}}}
\multiput(222,-1)(40,0){1}{{\tiny{$a$}}}

\end{picture}
}

\newcommand{\TotalCHHecThreeFive}{
\setlength{\unitlength}{0.5mm}
\begin{picture}(230,80)(-40,-20)

\multiput(10,50)(40,0){4}{\circle{2}}
\multiput(11,50)(40,0){4}{\line(1,0){8}}
\multiput(20,50)(40,0){4}{\circle{2}}
\multiput(21,50)(40,0){3}{\line(1,0){28}}

\multiput(10,40)(40,0){4}{\circle{2}}
\multiput(10,41)(40,0){4}{\line(0,1){8}}
\multiput(11,40)(40,0){4}{\line(1,0){8}}
\multiput(20,40)(40,0){4}{\circle{2}}
\multiput(20,41)(40,0){4}{\line(0,1){8}}

\multiput(00,30)(40,0){4}{\circle{2}}
\multiput(00.8,30.8)(40,0){4}{\line(1,1){8.4}}
\multiput(30,30)(40,0){4}{\circle{2}}
\multiput(29.2,30.8)(40,0){4}{\line(-1,1){8.4}}
\multiput(31,30)(40,0){3}{\line(1,0){8}}

\multiput(00,20)(40,0){4}{\circle{2}}
\multiput(00,21)(40,0){4}{\line(0,1){8}}
\multiput(30,20)(40,0){4}{\circle{2}}
\multiput(30,21)(40,0){4}{\line(0,1){8}}
\multiput(31,20)(40,0){3}{\line(1,0){8}}

\multiput(10,10)(40,0){4}{\circle{2}}
\multiput(9.2,10.8)(40,0){4}{\line(-1,1){8.4}}
\multiput(11,10)(40,0){4}{\line(1,0){8}}
\multiput(20,10)(40,0){4}{\circle{2}}
\multiput(20.8,10.8)(40,0){4}{\line(1,1){8.4}}

\multiput(10,00)(40,0){4}{\circle{2}}
\multiput(10,01)(40,0){4}{\line(0,1){8}}
\multiput(11,00)(40,0){4}{\line(1,0){8}}
\multiput(20,00)(40,0){4}{\circle{2}}
\multiput(20,01)(40,0){4}{\line(0,1){8}}
\multiput(21,00)(40,0){3}{\line(1,0){28}}

\multiput(20,50)(40,0){3}{\RIGHTlineTHICKsmall{28}}

\multiput(10,40)(40,0){4}{\UPlineTHICKsmall{8}}
\multiput(20,40)(40,0){4}{\UPlineTHICKsmall{8}}

\multiput(0,30)(40,0){4}{\UPRIGHTlineTHICKsmall{8.4}}
\multiput(30,30)(40,0){4}{\UPLEFTlineTHICKsmall{8.4}}

\multiput(00,20)(40,0){4}{\UPlineTHICKsmall{8}}
\multiput(30,20)(40,0){4}{\UPlineTHICKsmall{8}}

\multiput(10,10)(40,0){4}{\UPLEFTlineTHICKsmall{8.4}}
\multiput(20,10)(40,0){4}{\UPRIGHTlineTHICKsmall{8.4}}

\multiput(10,0)(40,0){4}{\RIGHTlineTHICKsmall{8}}
\multiput(10,0)(40,0){4}{\UPlineTHICKsmall{8}}
\multiput(20,0)(40,0){4}{\UPlineTHICKsmall{8}}

\put(130,51){\curveNorthEastTHICKsmall}
\put(20,61){\properRIGHTlineTHICKsmall{110}}
\put(20,51){\curveNorthWestTHICKsmall}

\put(151,40){\oval(20,20)[br]}
\put(161,60){\line(0,-1){20}}
\put(151,60){\oval(20,20)[tr]}
\put(-1,70){\line(1,0){152}}
\put(-1,60){\oval(20,20)[tl]}
\put(-11,40){\line(0,1){20}}
\put(-1,40){\oval(20,20)[bl]}

\put(20,-1){\oval(20,20)[bl]}
\put(20,-11){\line(1,0){110}}
\put(130,-1){\oval(20,20)[br]}

\put(-1,10){\oval(20,20)[tl]}
\put(151,10){\oval(20,20)[tr]}
\multiput(-11,-10)(172,0){2}{\line(0,1){20}}
\put(-1,-10){\oval(20,20)[bl]}
\put(151,-10){\oval(20,20)[br]}
\put(-1,-20){\line(1,0){152}}

\multiput(74,71)(0,0){1}{{\tiny{1}}}

\multiput(74,56.5)(0,0){1}{{\tiny{2}}}

\multiput(34,51)(40,0){3}{{\tiny{2}}}
\multiput(14,51)(40,0){4}{{\tiny{3}}}

\multiput(6.5,43.5)(40,0){4}{{\tiny{1}}}
\multiput(21,43.5)(40,0){4}{{\tiny{1}}}

\multiput(14,36)(40,0){4}{{\tiny{3}}}

\multiput(21,32)(40,0){4}{{\tiny{2}}}
\multiput(6,32)(40,0){4}{{\tiny{2}}}

\multiput(34,31)(40,0){3}{{\tiny{1}}}

\multiput(-3.5,23.5)(00,0){1}{{\tiny{3}}}
\multiput(26.5,23.5)(40,0){3}{{\tiny{3}}}
\multiput(41,23.5)(40,0){3}{{\tiny{3}}}
\multiput(151,23.5)(00,0){1}{{\tiny{3}}}

\multiput(34,16)(40,0){3}{{\tiny{1}}}

\multiput(6,15)(40,0){4}{{\tiny{2}}}
\multiput(21,15)(40,0){4}{{\tiny{2}}}

\multiput(14,11)(40,0){4}{{\tiny{3}}}

\multiput(6.5,3.5)(40,0){4}{{\tiny{1}}}
\multiput(21,3.5)(40,0){4}{{\tiny{1}}}

\multiput(14,-4)(40,0){4}{{\tiny{3}}}
\multiput(34,-4)(40,0){3}{{\tiny{2}}}

\multiput(74,-10)(0,0){1}{{\tiny{2}}}

\multiput(74,-24)(0,0){1}{{\tiny{1}}}

\multiput(5,49)(00,0){1}{{\tiny{$a$}}}
\multiput(18,52)(40,0){3}{{\tiny{$a$}}}
\multiput(48,52)(40,0){3}{{\tiny{$a$}}}
\multiput(142,49)(40,0){1}{{\tiny{$a$}}}

\multiput(5,39)(40,0){4}{{\tiny{$b$}}}
\multiput(22,39)(40,0){4}{{\tiny{$b$}}}

\multiput(-2,32)(40,0){4}{{\tiny{$c$}}}
\multiput(29,32)(40,0){4}{{\tiny{$c$}}}

\multiput(-2,16)(40,0){4}{{\tiny{$c$}}}
\multiput(29,16)(40,0){4}{{\tiny{$c$}}}

\multiput(9,12)(40,0){4}{{\tiny{$b$}}}
\multiput(18,12)(40,0){4}{{\tiny{$b$}}}

\multiput(5,-1)(40,0){1}{{\tiny{$a$}}}
\multiput(19,-4)(40,0){3}{{\tiny{$a$}}}
\multiput(49,-4)(40,0){3}{{\tiny{$a$}}}
\multiput(142,-1)(40,0){1}{{\tiny{$a$}}}

\end{picture}
}

\newcommand{\JointForHecTwentyTwo}{
\setlength{\unitlength}{0.5mm}
\begin{picture}(120,70)(50,-10)

\multiput(-6,50)(2,0){3}{\circle*{0.6}}
\multiput(10,50)(10,0){2}{\circle{2}}
\multiput(01,50)(10,0){3}{\line(1,0){8}}
\multiput(32,50)(2,0){3}{\circle*{0.6}}
\multiput(10,41)(10,0){2}{\line(0,1){8}}
\multiput(-6,40)(2,0){3}{\circle*{0.6}}
\multiput(10,40)(10,0){2}{\circle{2}}
\multiput(01,40)(10,0){3}{\line(1,0){8}}
\multiput(32,40)(2,0){3}{\circle*{0.6}}

\multiput(00,50)(10,0){3}{\RIGHTlineTHICKsmall{8}}
\multiput(00,40)(10,0){3}{\RIGHTlineTHICKsmall{8}}

\multiput(-6,10)(2,0){3}{\circle*{0.6}}
\multiput(10,10)(10,0){2}{\circle{2}}
\multiput(01,10)(10,0){3}{\line(1,0){8}}
\multiput(32,10)(2,0){3}{\circle*{0.6}}
\multiput(10,01)(10,0){2}{\line(0,1){8}}
\multiput(-6,00)(2,0){3}{\circle*{0.6}}
\multiput(10,00)(10,0){2}{\circle{2}}
\multiput(01,00)(10,0){3}{\line(1,0){8}}
\multiput(32,00)(2,0){3}{\circle*{0.6}}

\multiput(00,10)(20,0){2}{\RIGHTlineTHICKsmall{8}}
\multiput(10,00)(10,0){2}{\UPlineTHICKsmall{8}}
\multiput(00,00)(20,0){2}{\RIGHTlineTHICKsmall{8}}

\multiput(6,44)(00,0){1}{{\tiny{1}}}
\multiput(21,44)(00,0){1}{{\tiny{1}}}

\multiput(14,51)(00,0){1}{{\tiny{${\bf i}$}}}
\multiput(14,36)(00,0){1}{{\tiny{${\bf i}$}}}

\multiput(9,52)(00,0){1}{{\tiny{$x$}}}
\multiput(19,53)(00,0){1}{{\tiny{$y$}}}

\multiput(9,36)(00,0){1}{{\tiny{$z$}}}
\multiput(19,36)(00,0){1}{{\tiny{$u$}}}

\put(15,30){\vector(0,-1){10}}

\multiput(6,4)(00,0){1}{{\tiny{1}}}
\multiput(21,4)(00,0){1}{{\tiny{1}}}

\multiput(14,11)(00,0){1}{{\tiny{${\bf i}$}}}
\multiput(14,-4)(00,0){1}{{\tiny{${\bf i}$}}}

\multiput(9,12)(00,0){1}{{\tiny{$x$}}}
\multiput(19,13)(00,0){1}{{\tiny{$y$}}}

\multiput(9,-4)(00,0){1}{{\tiny{$z$}}}
\multiput(19,-4)(00,0){1}{{\tiny{$u$}}}

\put(50,20){${\footnotesize{
\left(
\begingroup
\renewcommand{\arraystretch}{1.5}
\begin{array}{l}
\{({\bf i},x,y,z,u),({\bf i},z,u,x,y),({\bf i},y,x,u,z),({\bf i},u,z,y,x)\}
\\
\cap  
\{({\bf 3},a,a,a,a), ({\bf 4},a,a,a,a), 
({\bf 3},b,c,{\bar b},{\bar a}), ({\bf 4},b,b,{\bar b},{\bar b}), \\
({\bf 4},c,d,{\bar a},{\bar d}), ({\bf 3},d,d,{\bar d},{\bar d}), 
({\bf 3},{\bar c},{\bar c},{\bar c},{\bar c}), ({\bf 3},{\bar c}^\prime,{\bar c}^\prime,{\bar c}^\prime,{\bar c}^\prime),  \\
({\bf 2},{\bar d}^\prime,{\bar a}^\prime,d^\prime,c^\prime), ({\bf 3},{\bar d}^\prime,{\bar d}^\prime,d^\prime,d^\prime),  
({\bf 3},{\bar a}^\prime,{\bar b}^\prime,c^\prime,b^\prime), 
({\bf 2},{\bar b}^\prime,{\bar b}^\prime,b^\prime,b^\prime),   \\
({\bf 2},a^\prime,a^\prime,a^\prime,a^\prime), ({\bf 3},a^\prime,a^\prime,a^\prime,a^\prime)\}\ne\emptyset
\end{array}
\endgroup
\right)}}$} 

\end{picture}
}

\newcommand{\spTotalCHHecThreeFive}{
\setlength{\unitlength}{0.5mm}
\begin{picture}(230,80)(-40,-20)

\multiput(10,50)(40,0){4}{\circle{2}}
\multiput(11,50)(40,0){4}{\line(1,0){8}}
\multiput(20,50)(40,0){4}{\circle{2}}
\multiput(21,50)(40,0){3}{\line(1,0){28}}

\multiput(10,40)(40,0){4}{\circle{2}}
\multiput(10,41)(40,0){4}{\line(0,1){8}}
\multiput(11,40)(40,0){4}{\line(1,0){8}}
\multiput(20,40)(40,0){4}{\circle{2}}
\multiput(20,41)(40,0){4}{\line(0,1){8}}

\multiput(00,30)(40,0){4}{\circle{2}}
\multiput(00.8,30.8)(40,0){4}{\line(1,1){8.4}}
\multiput(30,30)(40,0){4}{\circle{2}}
\multiput(29.2,30.8)(40,0){4}{\line(-1,1){8.4}}
\multiput(31,30)(40,0){3}{\line(1,0){8}}

\multiput(00,20)(40,0){4}{\circle{2}}
\multiput(00,21)(40,0){4}{\line(0,1){8}}
\multiput(30,20)(40,0){4}{\circle{2}}
\multiput(30,21)(40,0){4}{\line(0,1){8}}
\multiput(31,20)(40,0){3}{\line(1,0){8}}

\multiput(10,10)(40,0){4}{\circle{2}}
\multiput(9.2,10.8)(40,0){4}{\line(-1,1){8.4}}
\multiput(11,10)(40,0){4}{\line(1,0){8}}
\multiput(20,10)(40,0){4}{\circle{2}}
\multiput(20.8,10.8)(40,0){4}{\line(1,1){8.4}}

\multiput(10,00)(40,0){4}{\circle{2}}
\multiput(10,01)(40,0){4}{\line(0,1){8}}
\multiput(11,00)(40,0){4}{\line(1,0){8}}
\multiput(20,00)(40,0){4}{\circle{2}}
\multiput(20,01)(40,0){4}{\line(0,1){8}}
\multiput(21,00)(40,0){3}{\line(1,0){28}}

\multiput(10,50)(80,0){2}{\RIGHTlineTHICKsmall{8}}
\multiput(20,50)(40,0){3}{\RIGHTlineTHICKsmall{28}}

\multiput(10,40)(80,0){2}{\RIGHTlineTHICKsmall{8}}
\multiput(50,40)(80,0){2}{\UPlineTHICKsmall{8}}
\multiput(60,40)(80,0){2}{\UPlineTHICKsmall{8}}

\multiput(0,30)(40,0){4}{\UPRIGHTlineTHICKsmall{8.4}}
\multiput(30,30)(40,0){4}{\UPLEFTlineTHICKsmall{8.4}}

\multiput(30,20)(40,0){3}{\UPlineTHICKsmall{8}}
\multiput(40,20)(40,0){3}{\UPlineTHICKsmall{8}}

\multiput(10,10)(40,0){4}{\UPLEFTlineTHICKsmall{8.4}}
\multiput(10,10)(120,0){2}{\RIGHTlineTHICKsmall{8}}
\multiput(20,10)(40,0){4}{\UPRIGHTlineTHICKsmall{8.4}}

\multiput(10,0)(120,0){2}{\RIGHTlineTHICKsmall{8}}
\multiput(20,0)(40,0){3}{\RIGHTlineTHICKsmall{28}}
\multiput(50,0)(40,0){2}{\UPlineTHICKsmall{8}}
\multiput(60,0)(40,0){2}{\UPlineTHICKsmall{8}}

\put(130,51){\curveNorthEastTHICKsmall}
\put(20,61){\properRIGHTlineTHICKsmall{110}}
\put(20,51){\curveNorthWestTHICKsmall}

\put(151,40){\curveSouthEastTHICKsmall}
\put(161,60){\properDOWNlineTHICKsmall{20}}
\put(151,60){\curveNorthEastTHICKsmall}
\put(-1,70){\properRIGHTlineTHICKsmall{152}}
\put(-1,60){\curveNorthWestTHICKsmall}
\put(-11,40){\properUPlineTHICKsmall{20}}
\put(-1,40){\curveSouthWestTHICKsmall}

\put(-1,10){\curveNorthWestTHICKsmall}
\put(-11,10){\properDOWNlineTHICKsmall{20}}
\put(-1,-10){\curveSouthWestTHICKsmall}
\put(-1,-20){\properRIGHTlineTHICKsmall{152}}
\put(151,-10){\curveSouthEastTHICKsmall}
\put(161,-10){\properUPlineTHICKsmall{20}}
\put(151,10){\curveNorthEastTHICKsmall}

\put(20,-1){\curveSouthWestTHICKsmall}
\put(20,-11){\properRIGHTlineTHICKsmall{110}}
\put(130,-1){\curveSouthEastTHICKsmall}

\multiput(74,71)(0,0){1}{{\tiny{1}}}

\multiput(74,56.5)(0,0){1}{{\tiny{2}}}

\multiput(34,51)(40,0){3}{{\tiny{2}}}
\multiput(14,51)(40,0){4}{{\tiny{3}}}

\multiput(6.5,44)(40,0){4}{{\tiny{1}}}
\multiput(21,44)(40,0){4}{{\tiny{1}}}

\multiput(14,36)(40,0){4}{{\tiny{3}}}

\multiput(21,32)(40,0){4}{{\tiny{2}}}
\multiput(6,32)(40,0){4}{{\tiny{2}}}

\multiput(34,31)(40,0){3}{{\tiny{1}}}

\multiput(26.5,24)(40,0){4}{{\tiny{3}}}
\multiput(1,24)(40,0){4}{{\tiny{3}}}

\multiput(34,16)(40,0){3}{{\tiny{1}}}

\multiput(6,15)(40,0){4}{{\tiny{2}}}
\multiput(21,15)(40,0){4}{{\tiny{2}}}

\multiput(14,11)(40,0){4}{{\tiny{3}}}

\multiput(6.5,4)(40,0){4}{{\tiny{1}}}
\multiput(21,4)(40,0){4}{{\tiny{1}}}

\multiput(14,-4)(40,0){4}{{\tiny{3}}}
\multiput(34,-4)(40,0){3}{{\tiny{2}}}

\multiput(74,-10)(0,0){1}{{\tiny{2}}}
\multiput(74,-25)(0,0){1}{{\tiny{1}}}

\multiput(5,49)(00,0){1}{{\tiny{$a$}}}
\multiput(18,52)(40,0){3}{{\tiny{$a$}}}
\multiput(48,52)(40,0){3}{{\tiny{$a$}}}
\multiput(142,49)(40,0){1}{{\tiny{$a$}}}

\multiput(5,39)(40,0){4}{{\tiny{$b$}}}
\multiput(22,39)(40,0){4}{{\tiny{$b$}}}

\multiput(-2,32)(40,0){4}{{\tiny{$c$}}}
\multiput(29,32)(40,0){4}{{\tiny{$c$}}}

\multiput(-2,16)(40,0){4}{{\tiny{$c$}}}
\multiput(29,16)(40,0){4}{{\tiny{$c$}}}

\multiput(9,12)(40,0){4}{{\tiny{$b$}}}
\multiput(18,12)(40,0){4}{{\tiny{$b$}}}

\multiput(5,-1)(40,0){1}{{\tiny{$a$}}}
\multiput(19,-4)(40,0){3}{{\tiny{$a$}}}
\multiput(49,-4)(40,0){3}{{\tiny{$a$}}}
\multiput(142,-1)(40,0){1}{{\tiny{$a$}}}

\put(170,0){$(a^\prime\cong a)$}

\end{picture}
}

\newcommand{\spTotalCHHecThreeFiveForFourEughteen}{
\setlength{\unitlength}{0.5mm}
\begin{picture}(230,80)(-40,-20)

\multiput(10,50)(40,0){4}{\circle{2}}
\multiput(11,50)(40,0){4}{\line(1,0){8}}
\multiput(20,50)(40,0){4}{\circle{2}}
\multiput(21,50)(40,0){3}{\line(1,0){28}}

\multiput(10,40)(40,0){4}{\circle{2}}
\multiput(10,41)(40,0){4}{\line(0,1){8}}
\multiput(11,40)(40,0){4}{\line(1,0){8}}
\multiput(20,40)(40,0){4}{\circle{2}}
\multiput(20,41)(40,0){4}{\line(0,1){8}}

\multiput(00,30)(40,0){4}{\circle{2}}
\multiput(00.8,30.8)(40,0){4}{\line(1,1){8.4}}
\multiput(30,30)(40,0){4}{\circle{2}}
\multiput(29.2,30.8)(40,0){4}{\line(-1,1){8.4}}
\multiput(31,30)(40,0){3}{\line(1,0){8}}

\multiput(00,20)(40,0){4}{\circle{2}}
\multiput(00,21)(40,0){4}{\line(0,1){8}}
\multiput(30,20)(40,0){4}{\circle{2}}
\multiput(30,21)(40,0){4}{\line(0,1){8}}
\multiput(31,20)(40,0){3}{\line(1,0){8}}

\multiput(10,10)(40,0){4}{\circle{2}}
\multiput(9.2,10.8)(40,0){4}{\line(-1,1){8.4}}
\multiput(11,10)(40,0){4}{\line(1,0){8}}
\multiput(20,10)(40,0){4}{\circle{2}}
\multiput(20.8,10.8)(40,0){4}{\line(1,1){8.4}}

\multiput(10,00)(40,0){4}{\circle{2}}
\multiput(10,01)(40,0){4}{\line(0,1){8}}
\multiput(11,00)(40,0){4}{\line(1,0){8}}
\multiput(20,00)(40,0){4}{\circle{2}}
\multiput(20,01)(40,0){4}{\line(0,1){8}}
\multiput(21,00)(40,0){3}{\line(1,0){28}}

\multiput(10,50)(80,0){2}{\RIGHTlineTHICKsmall{8}}
\multiput(20,50)(40,0){3}{\RIGHTlineTHICKsmall{28}}

\multiput(10,40)(80,0){2}{\RIGHTlineTHICKsmall{8}}
\multiput(50,40)(80,0){2}{\UPlineTHICKsmall{8}}
\multiput(60,40)(80,0){2}{\UPlineTHICKsmall{8}}

\multiput(0,30)(40,0){4}{\UPRIGHTlineTHICKsmall{8.4}}
\multiput(30,30)(40,0){4}{\UPLEFTlineTHICKsmall{8.4}}

\multiput(30,20)(40,0){3}{\UPlineTHICKsmall{8}}
\multiput(40,20)(40,0){3}{\UPlineTHICKsmall{8}}

\multiput(10,10)(40,0){4}{\UPLEFTlineTHICKsmall{8.4}}
\multiput(10,10)(120,0){2}{\RIGHTlineTHICKsmall{8}}
\multiput(20,10)(40,0){4}{\UPRIGHTlineTHICKsmall{8.4}}

\multiput(10,0)(120,0){2}{\RIGHTlineTHICKsmall{8}}
\multiput(20,0)(40,0){3}{\RIGHTlineTHICKsmall{28}}
\multiput(50,0)(40,0){2}{\UPlineTHICKsmall{8}}
\multiput(60,0)(40,0){2}{\UPlineTHICKsmall{8}}

\put(130,51){\curveNorthEastTHICKsmall}
\put(20,61){\properRIGHTlineTHICKsmall{110}}
\put(20,51){\curveNorthWestTHICKsmall}

\put(151,40){\curveSouthEastTHICKsmall}
\put(161,60){\properDOWNlineTHICKsmall{20}}
\put(151,60){\curveNorthEastTHICKsmall}
\put(-1,70){\properRIGHTlineTHICKsmall{152}}
\put(-1,60){\curveNorthWestTHICKsmall}
\put(-11,40){\properUPlineTHICKsmall{20}}
\put(-1,40){\curveSouthWestTHICKsmall}

\put(-1,10){\curveNorthWestTHICKsmall}
\put(-11,10){\properDOWNlineTHICKsmall{20}}
\put(-1,-10){\curveSouthWestTHICKsmall}
\put(-1,-20){\properRIGHTlineTHICKsmall{152}}
\put(151,-10){\curveSouthEastTHICKsmall}
\put(161,-10){\properUPlineTHICKsmall{20}}
\put(151,10){\curveNorthEastTHICKsmall}

\put(20,-1){\curveSouthWestTHICKsmall}
\put(20,-11){\properRIGHTlineTHICKsmall{110}}
\put(130,-1){\curveSouthEastTHICKsmall}

\multiput(74,71)(0,0){1}{{\tiny{{\bf 4}}}}

\multiput(74,56.5)(0,0){1}{{\tiny{{\bf 2}}}}

\multiput(34,51)(40,0){3}{{\tiny{{\bf 2}}}}
\multiput(14,51)(40,0){4}{{\tiny{3}}}

\multiput(6.5,44)(80,0){2}{{\tiny{4}}}\multiput(46.5,44)(80,0){2}{{\tiny{{\bf 4}}}}
\multiput(21,44)(80,0){2}{{\tiny{4}}}\multiput(61,44)(80,0){2}{{\tiny{{\bf 4}}}}

\multiput(14,36)(40,0){4}{{\tiny{3}}}

\multiput(21,32)(40,0){4}{{\tiny{{\bf 2}}}}
\multiput(6,32)(40,0){4}{{\tiny{{\bf 2}}}}

\multiput(34,31)(40,0){3}{{\tiny{4}}}

\multiput(26.5,24)(40,0){4}{{\tiny{3}}}
\multiput(1,24)(40,0){4}{{\tiny{3}}}

\multiput(34,16)(40,0){3}{{\tiny{4}}}

\multiput(6,15)(40,0){4}{{\tiny{{\bf 2}}}}
\multiput(21,15)(40,0){4}{{\tiny{{\bf 2}}}}

\multiput(14,11)(40,0){4}{{\tiny{3}}}

\multiput(6.5,4)(00,0){1}{{\tiny{4}}}\multiput(46.5,4)(40,0){2}{{\tiny{{\bf 4}}}}\multiput(126.5,4)(00,0){1}{{\tiny{4}}}
\multiput(21,4)(40,0){1}{{\tiny{4}}}\multiput(61,4)(40,0){2}{{\tiny{{\bf 4}}}}\multiput(141,4)(00,0){1}{{\tiny{4}}}

\multiput(14,-4)(40,0){4}{{\tiny{3}}}
\multiput(34,-4)(40,0){3}{{\tiny{{\bf 2}}}}

\multiput(74,-10)(0,0){1}{{\tiny{{\bf 2}}}}
\multiput(74,-25)(0,0){1}{{\tiny{{\bf 4}}}}

\multiput(4,49)(00,0){1}{{\tiny{$a^\prime$}}}
\multiput(18,52)(40,0){3}{{\tiny{$a^\prime$}}}
\multiput(48,52)(40,0){3}{{\tiny{$a^\prime$}}}
\multiput(142,49)(40,0){1}{{\tiny{$a^\prime$}}}

\multiput(4,39)(40,0){4}{{\tiny{$a^\prime$}}}
\multiput(22,39)(40,0){4}{{\tiny{$a^\prime$}}}

\multiput(-2,32)(40,0){4}{{\tiny{$a^\prime$}}}
\multiput(29,32)(40,0){4}{{\tiny{$a^\prime$}}}

\multiput(-3,14)(40,0){4}{{\tiny{$a^\prime$}}}
\multiput(28,15)(40,0){4}{{\tiny{$a^\prime$}}}

\multiput(9,12)(40,0){4}{{\tiny{$a^\prime$}}}
\multiput(16.5,12)(40,0){4}{{\tiny{$a^\prime$}}}

\multiput(4,-1)(40,0){1}{{\tiny{$a^\prime$}}}
\multiput(19,-5)(40,0){3}{{\tiny{$a^\prime$}}}
\multiput(49,-5)(40,0){3}{{\tiny{$a^\prime$}}}
\multiput(142,-1)(40,0){1}{{\tiny{$a^\prime$}}}

\end{picture}
}

\newcommand{\TotalCHHecThreeEighteen}{
\setlength{\unitlength}{0.5mm}
\begin{picture}(230,130)(-20,-20)

\multiput(20,90)(60,0){4}{\circle{2}}
\multiput(21,90)(60,0){4}{\line(1,0){8}}
\multiput(30,90)(60,0){4}{\circle{2}}
\multiput(31,90)(60,0){3}{\line(1,0){48}}

\multiput(20,80)(60,0){4}{\circle{2}}
\multiput(20,81)(60,0){4}{\line(0,1){8}}
\multiput(21,80)(60,0){4}{\line(1,0){8}}
\multiput(30,80)(60,0){4}{\circle{2}}
\multiput(30,81)(60,0){4}{\line(0,1){8}}

\multiput(10,70)(30,0){8}{\circle{2}}
\multiput(10.8,70.8)(60,0){4}{\line(1,1){8.4}}
\multiput(39.2,70.8)(60,0){4}{\line(-1,1){8.4}}

\multiput(00,60)(30,0){8}{\circle{2}}
\multiput(00.8,60.8)(30,0){8}{\line(1,1){8.4}}
\multiput(20,60)(30,0){8}{\circle{2}}
\multiput(19.2,60.8)(30,0){8}{\line(-1,1){8.4}}
\multiput(21,60)(30,0){7}{\line(1,0){8}}

\multiput(10,50)(30,0){8}{\circle{2}}
\multiput(9.2,50.8)(30,0){8}{\line(-1,1){8.4}}
\multiput(10.8,50.8)(30,0){8}{\line(1,1){8.4}}

\multiput(10,40)(30,0){8}{\circle{2}}
\multiput(10,41)(30,0){8}{\line(0,1){8}}

\multiput(00,30)(30,0){8}{\circle{2}}
\multiput(00.8,30.8)(30,0){8}{\line(1,1){8.4}}
\multiput(20,30)(30,0){8}{\circle{2}}
\multiput(19.2,30.8)(30,0){8}{\line(-1,1){8.4}}
\multiput(21,30)(30,0){7}{\line(1,0){8}}

\multiput(10,20)(30,0){8}{\circle{2}}
\multiput(9.2,20.8)(30,0){8}{\line(-1,1){8.4}}
\multiput(10.8,20.8)(30,0){8}{\line(1,1){8.4}}

\multiput(20,10)(60,0){4}{\circle{2}}
\multiput(19.2,10.8)(60,0){4}{\line(-1,1){8.4}}
\multiput(21,10)(60,0){4}{\line(1,0){8}}
\multiput(30,10)(60,0){4}{\circle{2}}
\multiput(30.8,10.8)(60,0){4}{\line(1,1){8.4}}

\multiput(20,00)(60,0){4}{\circle{2}}
\multiput(20,01)(60,0){4}{\line(0,1){8}}
\multiput(21,00)(60,0){4}{\line(1,0){8}}
\multiput(30,00)(60,0){4}{\circle{2}}
\multiput(30,01)(60,0){4}{\line(0,1){8}}
\multiput(31,00)(60,0){3}{\line(1,0){48}}

\multiput(30,90)(60,0){3}{\RIGHTlineTHICKsmall{48}}

\multiput(20,80)(60,0){4}{\UPlineTHICKsmall{8}}
\multiput(30,80)(60,0){4}{\UPlineTHICKsmall{8}}

\multiput(10,70)(60,0){4}{\UPRIGHTlineTHICKsmall{8.4}}
\multiput(40,70)(60,0){4}{\UPLEFTlineTHICKsmall{8.4}}

\multiput(20,60)(30,0){8}{\UPLEFTlineTHICKsmall{8.4}}
\multiput(20,60)(30,0){7}{\RIGHTlineTHICKsmall{8}}

\multiput(10,50)(30,0){8}{\UPLEFTlineTHICKsmall{8.4}}

\multiput(10,40)(30,0){8}{\UPlineTHICKsmall{8}}

\multiput(00,30)(30,0){8}{\UPRIGHTlineTHICKsmall{8.4}}
\multiput(20,30)(30,0){7}{\RIGHTlineTHICKsmall{8}}

\multiput(10,20)(30,0){8}{\UPRIGHTlineTHICKsmall{8.4}}

\multiput(20,10)(60,0){4}{\UPLEFTlineTHICKsmall{8.4}}
\multiput(30,10)(60,0){4}{\UPRIGHTlineTHICKsmall{8.4}}
\multiput(80,10)(60,0){3}{\RIGHTlineTHICKsmall{8}}

\multiput(20,00)(10,0){2}{\UPlineTHICKsmall{8}}
\multiput(30,00)(60,0){3}{\RIGHTlineTHICKsmall{48}}
\multiput(80,00)(60,0){3}{\RIGHTlineTHICKsmall{8}}

\put(200,91){\curveNorthEastTHICKsmall}
\put(30,101){\properRIGHTlineTHICKsmall{170}}
\put(30,91){\curveNorthWestTHICKsmall}

\put(231,70){\curveSouthEastTHICKsmall}
\put(241,100){\properDOWNlineTHICKsmall{30}}
\put(231,100){\curveNorthEastTHICKsmall}
\put(-1,110){\properRIGHTlineTHICKsmall{232}}
\put(-1,100){\curveNorthWestTHICKsmall}
\put(-11,70){\properUPlineTHICKsmall{30}}
\put(-1,70){\curveSouthWestTHICKsmall}

\put(-1,20){\curveNorthWestTHICKsmall}
\put(-11,20){\properDOWNlineTHICKsmall{30}}
\put(-1,-10){\curveSouthWestTHICKsmall}
\put(-1,-20){\properRIGHTlineTHICKsmall{232}}
\put(231,-10){\curveSouthEastTHICKsmall}
\put(241,-10){\properUPlineTHICKsmall{30}}
\put(231,20){\curveNorthEastTHICKsmall}

\put(30,-1){\curveSouthWestTHICKsmall}
\put(30,-11){\properRIGHTlineTHICKsmall{170}}
\put(200,-1){\curveSouthEastTHICKsmall}

\multiput(114,112)(00,0){1}{{\tiny{2}}}

\multiput(114,96)(00,0){1}{{\tiny{2}}}

\multiput(24,91)(60,0){4}{{\tiny{1}}}
\multiput(54,91)(60,0){3}{{\tiny{2}}}

\multiput(16.5,84)(60,0){4}{{\tiny{3}}}
\multiput(31,84)(60,0){4}{{\tiny{3}}}

\multiput(24,81)(60,0){4}{{\tiny{1}}}

\multiput(11,75)(60,0){4}{{\tiny{2}}}
\multiput(36,75)(60,0){4}{{\tiny{2}}}

\multiput(1,65)(60,0){4}{{\tiny{3}}}
\multiput(16,65)(60,0){4}{{\tiny{1}}}
\multiput(31,65)(60,0){4}{{\tiny{1}}}
\multiput(46,65)(60,0){4}{{\tiny{3}}}

\multiput(24,61)(30,0){7}{{\tiny{2}}}

\multiput(6,55)(60,0){4}{{\tiny{1}}}
\multiput(12,55)(60,0){4}{{\tiny{3}}}
\multiput(36,55)(60,0){4}{{\tiny{3}}}
\multiput(42,55)(60,0){4}{{\tiny{1}}}

\multiput(06,44)(30,0){8}{{\tiny{2}}}

\multiput(6,32)(60,0){4}{{\tiny{1}}}
\multiput(12,32)(60,0){4}{{\tiny{3}}}
\multiput(36,32)(60,0){4}{{\tiny{3}}}
\multiput(42,32)(60,0){4}{{\tiny{1}}}

\multiput(24,25.5)(30,0){7}{{\tiny{2}}}

\multiput(1,22)(60,0){4}{{\tiny{3}}}
\multiput(16,22)(60,0){4}{{\tiny{1}}}
\multiput(31,22)(60,0){4}{{\tiny{1}}}
\multiput(46,22)(60,0){4}{{\tiny{3}}}

\multiput(11,12)(60,0){4}{{\tiny{2}}}
\multiput(36,12)(60,0){4}{{\tiny{2}}}

\multiput(24,6)(60,0){4}{{\tiny{1}}}

\multiput(16.5,4)(60,0){4}{{\tiny{3}}}
\multiput(31,4)(60,0){4}{{\tiny{3}}}

\multiput(24,-4)(60,0){4}{{\tiny{1}}}
\multiput(54,-4)(60,0){3}{{\tiny{2}}}

\multiput(114,-10)(00,0){1}{{\tiny{2}}}
\multiput(114,-25)(00,0){1}{{\tiny{2}}}

\multiput(15,89)(60,0){1}{{\tiny{$b$}}}
\multiput(29,92)(60,0){3}{{\tiny{$b$}}}
\multiput(79,92)(60,0){3}{{\tiny{$b$}}}
\multiput(212,89)(60,0){1}{{\tiny{$b$}}}

\multiput(19,76)(60,0){4}{{\tiny{$a$}}}
\multiput(28,76)(60,0){4}{{\tiny{$a$}}}

\multiput(7,72)(60,0){4}{{\tiny{$a$}}}
\multiput(40,72)(60,0){4}{{\tiny{$a$}}}

\multiput(-2,55)(60,0){4}{{\tiny{$b$}}}
\multiput(19,56)(60,0){4}{{\tiny{$a$}}}
\multiput(28,56)(60,0){4}{{\tiny{$a$}}}
\multiput(49,55)(60,0){4}{{\tiny{$b$}}}

\multiput(12,48)(30,0){8}{{\tiny{$b$}}}

\multiput(12,39)(30,0){8}{{\tiny{$b$}}}

\multiput(-2,32)(60,0){4}{{\tiny{$b$}}}
\multiput(19,32)(60,0){4}{{\tiny{$a$}}}
\multiput(28,32)(60,0){4}{{\tiny{$a$}}}
\multiput(49,32)(60,0){4}{{\tiny{$b$}}}

\multiput(7,16)(60,0){4}{{\tiny{$a$}}}
\multiput(40,16)(60,0){4}{{\tiny{$a$}}}

\multiput(19,12)(60,0){4}{{\tiny{$a$}}}
\multiput(28,12)(60,0){4}{{\tiny{$a$}}}

\multiput(15,-1)(60,0){1}{{\tiny{$b$}}}
\multiput(29,-5)(60,0){3}{{\tiny{$b$}}}
\multiput(79,-5)(60,0){3}{{\tiny{$b$}}}
\multiput(212,-1)(60,0){1}{{\tiny{$b$}}}

\end{picture}
}

\newcommand{\TotalCHHecThreeSixteen}{
\setlength{\unitlength}{0.5mm}
\begin{picture}(180,130)(0,-20)

\put(51,61){\oval(18,18)[tr]}
\put(119,61){\oval(18,18)[tl]}

\multiput(10,70)(140,0){2}{\circle{2}}
\multiput(11,70)(140,0){2}{\line(1,0){8}}
\multiput(20,70)(140,0){2}{\circle{2}}
\multiput(21,70)(98,0){2}{\line(1,0){30}}

\multiput(00,60)(140,0){2}{\circle{2}}
\multiput(00.8,60.8)(140,0){2}{\line(1,1){8.4}}
\multiput(30,60)(140,0){2}{\circle{2}}
\multiput(29.2,60.8)(140,0){2}{\line(-1,1){8.4}}
\multiput(60,60)(40,0){2}{\circle{2}}
\multiput(61,60)(40,0){2}{\line(1,0){8}}
\multiput(70,60)(40,0){2}{\circle{2}}
\multiput(71,60)(00,0){1}{\line(1,0){28}}

\multiput(10,50)(140,0){2}{\circle{2}}
\multiput(11,50)(140,0){2}{\line(1,0){8}}
\multiput(9.2,50.8)(140,0){2}{\line(-1,1){8.4}}
\multiput(20,50)(140,0){2}{\circle{2}}
\multiput(20.8,50.8)(140,0){2}{\line(1,1){8.4}}
\multiput(40,50)(90,0){2}{\circle{2}}
\multiput(39.2,50.8)(00,0){1}{\line(-1,1){8.4}}
\multiput(60,50)(40,0){2}{\circle{2}}
\multiput(60,51)(40,0){2}{\line(0,1){8}}
\multiput(61,50)(40,0){2}{\line(1,0){8}}
\multiput(70,50)(40,0){2}{\circle{2}}
\multiput(70,51)(40,0){2}{\line(0,1){8}}
\multiput(130.8,50.8)(00,0){1}{\line(1,1){8.4}}

\multiput(00,40)(170,0){2}{\circle{2}}
\multiput(00.8,40.8)(00,0){1}{\line(1,1){8.4}}
\multiput(30,40)(110,0){2}{\circle{2}}
\multiput(29.2,40.8)(110,0){2}{\line(-1,1){8.4}}
\multiput(30.2,40.8)(110,0){2}{\line(1,1){8.4}}
\multiput(50,40)(70,0){2}{\circle{2}}
\multiput(49.2,40.8)(70,0){2}{\line(-1,1){8.4}}
\multiput(50.2,40.8)(70,0){2}{\line(1,1){8.4}}
\multiput(80,40)(10,0){2}{\circle{2}}
\multiput(79.2,40.8)(00,0){1}{\line(-1,1){8.4}}
\multiput(81,40)(00,0){1}{\line(1,0){8}}
\multiput(90.8,40.8)(00,0){1}{\line(1,1){8.4}}
\multiput(169.2,40.8)(00,0){1}{\line(-1,1){8.4}}

\multiput(00,30)(170,0){2}{\circle{2}}
\multiput(00,31)(170,0){2}{\line(0,1){8}}
\multiput(30,30)(110,0){2}{\circle{2}}
\multiput(30,31)(110,0){2}{\line(0,1){8}}
\multiput(50,30)(70,0){2}{\circle{2}}
\multiput(50,31)(70,0){2}{\line(0,1){8}}
\multiput(80,30)(10,0){2}{\circle{2}}
\multiput(81,30)(00,0){1}{\line(1,0){8}}
\multiput(80,31)(10,0){2}{\line(0,1){8}}

\multiput(10,20)(140,0){2}{\circle{2}}
\multiput(9.2,20.8)(140,0){2}{\line(-1,1){8.4}}
\multiput(11,20)(140,0){2}{\line(1,0){8}}
\multiput(20,20)(140,0){2}{\circle{2}}
\multiput(20.8,20.8)(140,0){2}{\line(1,1){8.4}}
\multiput(40,20)(90,0){2}{\circle{2}}
\multiput(39.2,20.8)(90,0){2}{\line(-1,1){8.4}}
\multiput(40.8,20.8)(90,0){2}{\line(1,1){8.4}}
\multiput(60,20)(40,0){2}{\circle{2}}
\multiput(59.2,20.8)(40,0){2}{\line(-1,1){8.4}}
\multiput(61,20)(40,0){2}{\line(1,0){8}}
\multiput(70,20)(40,0){2}{\circle{2}}
\multiput(70.8,20.8)(40,0){2}{\line(1,1){8.4}}

\multiput(10,10)(70,0){3}{\circle{2}}
\multiput(10,11)(140,0){2}{\line(0,1){8}}
\multiput(11,10)(70,0){3}{\line(1,0){8}}
\multiput(20,10)(70,0){3}{\circle{2}}
\multiput(20,11)(140,0){2}{\line(0,1){8}}
\multiput(50,10)(40,0){2}{\circle{2}}
\multiput(49.2,10.8)(70,0){2}{\line(-1,1){8.4}}
\multiput(50.8,10.8)(70,0){2}{\line(1,1){8.4}}
\multiput(79.2,10.8)(00,0){1}{\line(-1,1){8.4}}
\multiput(90.8,10.8)(00,0){1}{\line(1,1){8.4}}
\multiput(120,10)(00,0){1}{\circle{2}}

\multiput(29,00)(82,0){2}{\line(1,0){30}}
\multiput(60,00)(40,0){2}{\circle{2}}
\multiput(59.2,00.8)(40,0){2}{\line(-1,1){8.4}}
\multiput(61,00)(40,0){2}{\line(1,0){8}}
\multiput(70,00)(40,0){2}{\circle{2}}
\multiput(70.8,00.8)(40,0){2}{\line(1,1){8.4}}
\multiput(71,00)(00,0){1}{\line(1,0){28}}

\put(29,9){\oval(18,18)[bl]}
\put(141,9){\oval(18,18)[br]}

\put(19,71){\oval(18,18)[tl]}
\multiput(19,80)(00,00){1}{\line(1,0){132}}
\put(151,71){\oval(18,18)[tr]}

\put(-01,69){\oval(18,18)[bl]}
\multiput(-10,69)(190,00){2}{\line(0,1){12}}
\put(-01,81){\oval(18,18)[tl]}
\multiput(-1,90)(00,00){1}{\line(1,0){172}}
\put(171,81){\oval(18,18)[tr]}
\put(171,69){\oval(18,18)[br]}

\put(-11,49){\oval(18,18)[bl]}
\multiput(-20,49)(210,00){2}{\line(0,1){42}}
\put(-11,91){\oval(18,18)[tl]}
\multiput(-11,100)(00,00){1}{\line(1,0){192}}
\put(181,91){\oval(18,18)[tr]}
\put(181,49){\oval(18,18)[br]}

\multiput(-11,40)(182,00){2}{\line(1,0){10}}

\put(-01,21){\oval(18,18)[tl]}
\multiput(-10,21)(190,00){2}{\line(0,-1){32}}
\put(171,21){\oval(18,18)[tr]}
\multiput(-1,-20)(00,00){1}{\line(1,0){172}}
\put(-01,-11){\oval(18,18)[bl]}
\put(171,-11){\oval(18,18)[br]}

\multiput(10,9)(150,00){2}{\line(0,-1){10}}
\put(19,-1){\oval(18,18)[bl]}
\multiput(19,-10)(00,00){1}{\line(1,0){132}}
\put(151,-1){\oval(18,18)[br]}

\multiput(10,70)(140,0){2}{\RIGHTlineTHICKsmall{8}}

\multiput(00,60)(140,0){2}{\UPRIGHTlineTHICKsmall{8.4}}
\multiput(30,60)(140,0){2}{\UPLEFTlineTHICKsmall{8.4}}
\multiput(60,60)(40,0){2}{\RIGHTlineTHICKsmall{8}}

\multiput(10,50)(140,0){2}{\UPLEFTlineTHICKsmall{8.4}}
\multiput(20,50)(140,0){2}{\UPRIGHTlineTHICKsmall{8.4}}
\multiput(60,50)(50,0){2}{\UPlineTHICKsmall{8}}
\multiput(70,50)(30,0){2}{\UPlineTHICKsmall{8}}

\multiput(00,40)(140,0){2}{\UPRIGHTlineTHICKsmall{8.4}}
\multiput(30,40)(140,0){2}{\UPLEFTlineTHICKsmall{8.4}}
\multiput(30,40)(90,0){2}{\UPRIGHTlineTHICKsmall{8.4}}
\multiput(50,40)(90,0){2}{\UPLEFTlineTHICKsmall{8.4}}
\multiput(50,40)(40,0){2}{\UPRIGHTlineTHICKsmall{8.4}}
\multiput(80,40)(40,0){2}{\UPLEFTlineTHICKsmall{8.4}}
\multiput(80,40)(00,0){1}{\RIGHTlineTHICKsmall{8}}

\multiput(00,30)(170,0){2}{\UPlineTHICKsmall{8.4}}
\multiput(80,30)(00,0){1}{\RIGHTlineTHICKsmall{8}}

\multiput(10,20)(140,0){2}{\UPLEFTlineTHICKsmall{8.4}}
\multiput(20,20)(140,0){2}{\UPRIGHTlineTHICKsmall{8.4}}
\multiput(40,20)(60,0){2}{\UPLEFTlineTHICKsmall{8.4}}
\multiput(40,20)(90,0){2}{\UPRIGHTlineTHICKsmall{8.4}}
\multiput(60,20)(70,0){2}{\UPLEFTlineTHICKsmall{8.4}}
\multiput(70,20)(40,0){2}{\UPRIGHTlineTHICKsmall{8.4}}

\multiput(10,10)(140,0){2}{\UPlineTHICKsmall{8}}
\multiput(10,10)(140,0){2}{\RIGHTlineTHICKsmall{8}}
\multiput(20,10)(140,0){2}{\UPlineTHICKsmall{8}}
\multiput(50,10)(40,0){2}{\UPRIGHTlineTHICKsmall{8.4}}
\multiput(80,10)(40,0){2}{\UPLEFTlineTHICKsmall{8.4}}

\multiput(60,00)(40,0){2}{\UPLEFTlineTHICKsmall{8.4}}
\multiput(60,00)(40,0){2}{\RIGHTlineTHICKsmall{8}}
\multiput(70,00)(40,0){2}{\UPRIGHTlineTHICKsmall{8.4}}

\multiput(84,101)(00,0){1}{{\tiny{2}}}

\multiput(84,91)(00,0){1}{{\tiny{3}}}
\multiput(84,81)(00,0){1}{{\tiny{3}}}

\multiput(14,71)(140,0){2}{{\tiny{2}}}
\multiput(39,71)(90,0){2}{{\tiny{3}}}

\multiput(1,65)(140,0){2}{{\tiny{1}}}
\multiput(26,65)(140,0){2}{{\tiny{1}}}

\multiput(64,61)(40,0){2}{{\tiny{2}}}
\multiput(84,61)(00,0){1}{{\tiny{3}}}

\multiput(6,55)(140,0){2}{{\tiny{2}}}
\multiput(21,55)(140,0){2}{{\tiny{2}}}
\multiput(36,55)(00,0){1}{{\tiny{3}}}
\multiput(131,55)(00,0){1}{{\tiny{3}}}

\multiput(56.5,54)(40,0){2}{{\tiny{1}}}
\multiput(71,54)(40,0){2}{{\tiny{1}}}

\multiput(14,51)(140,0){2}{{\tiny{1}}}
\multiput(64,51)(40,0){2}{{\tiny{2}}}

\multiput(01,45)(140,0){2}{{\tiny{3}}}
\multiput(26,45)(140,0){2}{{\tiny{3}}}
\multiput(31,45)(00,0){1}{{\tiny{2}}}
\multiput(46,45)(00,0){1}{{\tiny{1}}}
\multiput(51,45)(40,0){2}{{\tiny{3}}}
\multiput(76,45)(40,0){2}{{\tiny{3}}}
\multiput(121,45)(00,0){1}{{\tiny{1}}}
\multiput(136,45)(00,0){1}{{\tiny{2}}}

\multiput(84,41)(00,0){1}{{\tiny{1}}}

\multiput(-3.5,34)(140,0){2}{{\tiny{1}}}
\multiput(31,34)(140,0){2}{{\tiny{1}}}
\multiput(46.5,34)(40,0){2}{{\tiny{2}}}
\multiput(81,34)(40,0){2}{{\tiny{2}}}

\multiput(84,26)(00,0){1}{{\tiny{1}}}

\multiput(01,22)(140,0){2}{{\tiny{3}}}
\multiput(26,22)(140,0){2}{{\tiny{3}}}
\multiput(31,22)(00,0){1}{{\tiny{2}}}
\multiput(46,22)(00,0){1}{{\tiny{1}}}
\multiput(51,22)(40,0){2}{{\tiny{3}}}
\multiput(76,22)(40,0){2}{{\tiny{3}}}
\multiput(121,22)(00,0){1}{{\tiny{1}}}
\multiput(136,22)(00,0){1}{{\tiny{2}}}

\multiput(14,16)(140,0){2}{{\tiny{1}}}
\multiput(64,16)(40,0){2}{{\tiny{2}}}

\multiput(6.5,14)(140,0){2}{{\tiny{2}}}
\multiput(21,14)(140,0){2}{{\tiny{2}}}

\multiput(41,12)(00,0){1}{{\tiny{3}}}
\multiput(56,12)(40,0){2}{{\tiny{1}}}
\multiput(71,12)(40,0){2}{{\tiny{1}}}
\multiput(126,12)(00,0){1}{{\tiny{3}}}

\multiput(14,06)(140,0){2}{{\tiny{1}}}
\multiput(84,06)(00,0){1}{{\tiny{3}}}

\multiput(51,02)(40,0){2}{{\tiny{2}}}
\multiput(76,02)(40,0){2}{{\tiny{2}}}

\multiput(39,-4)(45,0){3}{{\tiny{3}}}
\multiput(64,-4)(40,0){2}{{\tiny{1}}}

\multiput(84,-14)(00,0){1}{{\tiny{3}}}
\multiput(84,-24)(00,0){1}{{\tiny{2}}}

\multiput(9.5,66)(148,0){2}{{\tiny{$a$}}}
\multiput(18.5,65)(131,0){2}{{\tiny{$b$}}}

\multiput(-3,56.5)(173,0){2}{{\tiny{$c$}}}
\multiput(31.6,60)(104,0){2}{{\tiny{$e$}}}
\multiput(55.5,59)(56.5,0){2}{{\tiny{$f$}}}
\multiput(69,63)(30,0){2}{{\tiny{$g$}}}

\multiput(9.5,45)(148,0){2}{{\tiny{$d$}}}
\multiput(18.5,45)(131,0){2}{{\tiny{$h$}}}
\multiput(41.5,50)(84,0){2}{{\tiny{$e$}}}
\multiput(59,45)(49,0){2}{{\tiny{$f$}}}
\multiput(67.5,46)(32,0){2}{{\tiny{$g$}}}

\multiput(2,36.5)(164,0){2}{{\tiny{$i$}}}
\multiput(25.5,36.5)(116,0){2}{{\tiny{$h$}}}
\multiput(51.5,36.5)(64,0){2}{{\tiny{$b$}}}
\multiput(75.5,36.5)(16,0){2}{{\tiny{$g$}}}

\multiput(2,30.5)(164,0){2}{{\tiny{$j$}}}
\multiput(25.5,30.5)(116,0){2}{{\tiny{$d$}}}
\multiput(51.5,30.5)(64,0){2}{{\tiny{$a$}}}
\multiput(75.5,30.5)(16,0){2}{{\tiny{$f$}}}

\multiput(9.5,23)(149,0){2}{{\tiny{$j$}}}
\multiput(18.5,22)(131,0){2}{{\tiny{$i$}}}
\multiput(35.5,19)(96.5,0){2}{{\tiny{$c$}}}
\multiput(59.5,22)(48.5,0){2}{{\tiny{$a$}}}
\multiput(68.5,22)(31,0){2}{{\tiny{$b$}}}

\multiput(5.5,9)(156.5,0){2}{{\tiny{$j$}}}
\multiput(22,9)(123.5,0){2}{{\tiny{$i$}}}
\multiput(45.5,9)(76.5,0){2}{{\tiny{$c$}}}
\multiput(79.5,12)(8.5,0){2}{{\tiny{$e$}}}

\multiput(58.5,-5)(50,0){2}{{\tiny{$d$}}}
\multiput(68.5,-5)(30,0){2}{{\tiny{$h$}}}

\end{picture}
}

\newcommand{\TotalCHHecThreeSeventeen}{
\setlength{\unitlength}{0.5mm}
\begin{picture}(180,150)(5,-35)

\multiput(00,70)(30,0){4}{\circle{2}}
\multiput(01,70)(30,0){3}{\line(1,0){28}}
\multiput(91,70)(00,0){1}{\line(1,0){8}}
\multiput(100,70)(30,0){4}{\circle{2}}
\multiput(101,70)(30,0){3}{\line(1,0){28}}

\multiput(00,60)(30,0){4}{\circle{2}}
\multiput(00,61)(30,0){4}{\line(0,1){8}}
\multiput(91,60)(00,0){1}{\line(1,0){8}}
\multiput(100,60)(30,0){4}{\circle{2}}
\multiput(100,61)(30,0){4}{\line(0,1){8}}

\multiput(10,50)(30,0){3}{\circle{2}}
\multiput(9.2,50.8)(30,0){3}{\line(-1,1){8.4}} 
\multiput(11,50)(30,0){3}{\line(1,0){8}}
\multiput(20,50)(30,0){3}{\circle{2}}
\multiput(20.8,50.8)(30,0){3}{\line(1,1){8.4}} 
\multiput(110,50)(30,0){3}{\circle{2}}
\multiput(109.2,50.8)(30,0){3}{\line(-1,1){8.4}} 
\multiput(111,50)(30,0){3}{\line(1,0){8}}
\multiput(120,50)(30,0){3}{\circle{2}}
\multiput(120.8,50.8)(30,0){3}{\line(1,1){8.4}} 

\multiput(00,40)(30,0){4}{\circle{2}}
\multiput(00.8,40.8)(30,0){3}{\line(1,1){8.4}}
\multiput(29.2,40.8)(30,0){3}{\line(-1,1){8.4}} 
\multiput(91,40)(00,0){1}{\line(1,0){8}}
\multiput(100,40)(30,0){4}{\circle{2}}
\multiput(100.8,40.8)(30,0){3}{\line(1,1){8.4}}
\multiput(129.2,40.8)(30,0){3}{\line(-1,1){8.4}} 

\multiput(00,30)(30,0){4}{\circle{2}}
\multiput(00,31)(30,0){4}{\line(0,1){8}}
\multiput(91,30)(00,0){1}{\line(1,0){8}}
\multiput(100,30)(30,0){4}{\circle{2}}
\multiput(100,31)(30,0){4}{\line(0,1){8}}

\multiput(10,20)(30,0){3}{\circle{2}}
\multiput(9.2,20.8)(30,0){3}{\line(-1,1){8.4}} 
\multiput(11,20)(30,0){3}{\line(1,0){8}}
\multiput(20,20)(30,0){3}{\circle{2}}
\multiput(20.8,20.8)(30,0){3}{\line(1,1){8.4}} 
\multiput(110,20)(30,0){3}{\circle{2}}
\multiput(109.2,20.8)(30,0){3}{\line(-1,1){8.4}} 
\multiput(111,20)(30,0){3}{\line(1,0){8}}
\multiput(120,20)(30,0){3}{\circle{2}}
\multiput(120.8,20.8)(30,0){3}{\line(1,1){8.4}} 

\multiput(00,10)(30,0){4}{\circle{2}}
\multiput(00.8,10.8)(30,0){3}{\line(1,1){8.4}}
\multiput(29.2,10.8)(30,0){3}{\line(-1,1){8.4}} 
\multiput(91,10)(00,0){1}{\line(1,0){8}}
\multiput(100,10)(30,0){4}{\circle{2}}
\multiput(100.8,10.8)(30,0){3}{\line(1,1){8.4}}
\multiput(129.2,10.8)(30,0){3}{\line(-1,1){8.4}} 

\multiput(00,00)(30,0){4}{\circle{2}}
\multiput(00,01)(30,0){4}{\line(0,1){8}}
\multiput(01,00)(30,0){3}{\line(1,0){28}}
\multiput(91,00)(00,0){1}{\line(1,0){8}}
\multiput(100,00)(30,0){4}{\circle{2}}
\multiput(100,01)(30,0){4}{\line(0,1){8}}
\multiput(101,00)(30,0){3}{\line(1,0){28}}

\multiput(00,70)(100,0){2}{\RIGHTlineTHICKsmall{28}}
\multiput(60,70)(100,0){2}{\RIGHTlineTHICKsmall{28}}

\multiput(00,60)(100,0){2}{\UPlineTHICKsmall{8}}
\multiput(30,60)(100,0){2}{\UPlineTHICKsmall{8}}
\multiput(60,60)(100,0){2}{\UPlineTHICKsmall{8}}
\multiput(90,60)(100,0){2}{\UPlineTHICKsmall{8}}
\multiput(90,60)(00,0){1}{\RIGHTlineTHICKsmall{8}}

\multiput(10,50)(100,0){2}{\RIGHTlineTHICKsmall{8}}
\multiput(40,50)(100,0){2}{\UPLEFTlineTHICKsmall{8.4}}
\multiput(50,50)(100,0){2}{\UPRIGHTlineTHICKsmall{8.4}}
\multiput(70,50)(100,0){2}{\RIGHTlineTHICKsmall{8}}

\multiput(00,40)(100,0){2}{\UPRIGHTlineTHICKsmall{8.4}}
\multiput(30,40)(100,0){2}{\UPLEFTlineTHICKsmall{8.4}}
\multiput(30,40)(100,0){2}{\UPRIGHTlineTHICKsmall{8.4}}
\multiput(60,40)(100,0){2}{\UPLEFTlineTHICKsmall{8.4}}
\multiput(60,40)(100,0){2}{\UPRIGHTlineTHICKsmall{8.4}}
\multiput(90,40)(100,0){2}{\UPLEFTlineTHICKsmall{8.4}}

\multiput(00,30)(100,0){2}{\UPlineTHICKsmall{8}}
\multiput(90,30)(100,0){2}{\UPlineTHICKsmall{8}}

\multiput(10,20)(100,0){2}{\UPLEFTlineTHICKsmall{8.4}}
\multiput(20,20)(100,0){2}{\UPRIGHTlineTHICKsmall{8.4}}
\multiput(40,20)(100,0){2}{\UPLEFTlineTHICKsmall{8.4}}
\multiput(40,20)(100,0){2}{\RIGHTlineTHICKsmall{8}}
\multiput(50,20)(100,0){2}{\UPRIGHTlineTHICKsmall{8.4}}
\multiput(70,20)(100,0){2}{\UPLEFTlineTHICKsmall{8.4}}
\multiput(80,20)(100,0){2}{\UPRIGHTlineTHICKsmall{8.4}}

\multiput(00,10)(100,0){2}{\UPRIGHTlineTHICKsmall{8.4}}
\multiput(30,10)(100,0){2}{\UPLEFTlineTHICKsmall{8.4}}
\multiput(60,10)(100,0){2}{\UPRIGHTlineTHICKsmall{8.4}}
\multiput(90,10)(100,0){2}{\UPLEFTlineTHICKsmall{8.4}}

\multiput(00,00)(100,0){2}{\RIGHTlineTHICKsmall{28}}
\multiput(00,00)(100,0){2}{\UPlineTHICKsmall{8}}
\multiput(30,00)(100,0){2}{\UPlineTHICKsmall{8}}
\multiput(60,00)(100,0){2}{\RIGHTlineTHICKsmall{28}}
\multiput(60,00)(100,0){2}{\UPlineTHICKsmall{8}}
\multiput(90,00)(100,0){2}{\UPlineTHICKsmall{8}}

\put(10,71){\oval(20,20)[tl]}
\put(10,81){\line(1,0){170}}
\put(180,71){\oval(20,20)[tr]}

\put(-1,70){\curveSouthWestTHICKsmall}
\put(191,70){\curveSouthEastTHICKsmall}
\multiput(-11,70)(212,0){2}{\properUPlineTHICKsmall{10}}
\put(-1,80){\curveNorthWestTHICKsmall}
\put(191,80){\curveNorthEastTHICKsmall}
\multiput(-11,70)(212,0){2}{\properUPlineTHICKsmall{10}}
\put(-1,90){\properRIGHTlineTHICKsmall{192}}

\multiput(-11,40)(202,0){2}{\line(1,0){10}}
\put(-11,50){\oval(20,20)[bl]}
\put(201,50){\oval(20,20)[br]}
\multiput(-21,50)(232,0){2}{\line(0,1){40}}
\put(-11,90){\oval(20,20)[tl]}
\put(201,90){\oval(20,20)[tr]}
\put(-11,100){\line(1,0){212}}

\multiput(-11,30)(202,0){2}{\line(1,0){10}}
\put(-11,20){\oval(20,20)[tl]}
\put(201,20){\oval(20,20)[tr]}
\multiput(-21,-20)(232,0){2}{\line(0,1){40}}
\put(-11,-20){\oval(20,20)[bl]}
\put(201,-20){\oval(20,20)[br]}
\put(-11,-30){\line(1,0){212}}

\put(-1,00){\oval(20,20)[tl]}
\put(191,00){\oval(20,20)[tr]}
\multiput(-11,-10)(212,0){2}{\line(0,1){10}}
\put(-1,-10){\oval(20,20)[bl]}
\put(191,-10){\oval(20,20)[br]}
\put(-1,-20){\line(1,0){192}}

\put(10,-01){\oval(20,20)[bl]}
\put(10,-11){\line(1,0){170}}
\put(180,-01){\oval(20,20)[br]}

\multiput(94,101)(00,0){1}{\tiny{2}}

\multiput(94,91)(00,0){1}{\tiny{1}}

\multiput(94,77)(00,0){1}{\tiny{1}}

\multiput(14,71)(100,0){2}{\tiny{2}}
\multiput(44,71)(100,0){2}{\tiny{1}}
\multiput(74,71)(100,0){2}{\tiny{2}}
\multiput(94,71)(00,0){1}{\tiny{1}}

\multiput(-3.5,64)(100,0){2}{{\tiny{3}}}
\multiput(31,64)(100,0){2}{{\tiny{3}}}
\multiput(56.5,64)(100,0){2}{{\tiny{3}}}
\multiput(91,64)(100,0){2}{{\tiny{3}}}

\multiput(94,56)(00,0){1}{\tiny{1}}

\multiput(6,55)(100,0){2}{{\tiny{2}}}
\multiput(21,55)(100,0){2}{{\tiny{2}}}
\multiput(36,55)(100,0){2}{{\tiny{1}}}
\multiput(51,55)(100,0){2}{{\tiny{1}}}
\multiput(66,55)(100,0){2}{{\tiny{2}}}
\multiput(81,55)(100,0){2}{{\tiny{2}}}

\multiput(14,51)(100,0){2}{\tiny{3}}
\multiput(44,51)(100,0){2}{\tiny{3}}
\multiput(74,51)(100,0){2}{\tiny{3}}

\multiput(01,45)(100,0){2}{{\tiny{1}}}
\multiput(26,45)(100,0){2}{{\tiny{1}}}
\multiput(31,45)(100,0){2}{{\tiny{2}}}
\multiput(56,45)(100,0){2}{{\tiny{2}}}
\multiput(61,45)(100,0){2}{{\tiny{1}}}
\multiput(86,45)(100,0){2}{{\tiny{1}}}

\multiput(94,41)(00,0){1}{\tiny{2}}

\multiput(-3.5,34)(100,0){2}{{\tiny{3}}}
\multiput(31,34)(100,0){2}{{\tiny{3}}}
\multiput(56.5,34)(100,0){2}{{\tiny{3}}}
\multiput(91,34)(100,0){2}{{\tiny{3}}}

\multiput(94,26)(00,0){1}{\tiny{2}}

\multiput(01,22)(100,0){2}{{\tiny{1}}}
\multiput(26,22)(100,0){2}{{\tiny{1}}}
\multiput(31,22)(100,0){2}{{\tiny{2}}}
\multiput(56,22)(100,0){2}{{\tiny{2}}}
\multiput(61,22)(100,0){2}{{\tiny{1}}}
\multiput(86,22)(100,0){2}{{\tiny{1}}}

\multiput(14,16)(100,0){2}{\tiny{3}}
\multiput(44,16)(100,0){2}{\tiny{3}}
\multiput(74,16)(100,0){2}{\tiny{3}}

\multiput(06,12)(100,0){2}{{\tiny{2}}}
\multiput(21,12)(100,0){2}{{\tiny{2}}}
\multiput(36,12)(100,0){2}{{\tiny{1}}}
\multiput(51,12)(100,0){2}{{\tiny{1}}}
\multiput(66,12)(100,0){2}{{\tiny{2}}}
\multiput(81,12)(100,0){2}{{\tiny{2}}}

\multiput(94,11)(00,0){1}{\tiny{1}}

\multiput(-3.5,4)(100,0){2}{{\tiny{3}}}
\multiput(31,4)(100,0){2}{{\tiny{3}}}
\multiput(56.5,4)(100,0){2}{{\tiny{3}}}
\multiput(91,4)(100,0){2}{{\tiny{3}}}

\multiput(14,-4)(100,0){2}{\tiny{2}}
\multiput(44,-4)(100,0){2}{\tiny{1}}
\multiput(74,-4)(100,0){2}{\tiny{2}}
\multiput(94,-4)(00,0){1}{\tiny{1}}

\multiput(94,-10)(00,0){1}{\tiny{1}}

\multiput(94,-24)(00,0){1}{\tiny{1}}

\multiput(94,-34)(00,0){1}{\tiny{2}}

\multiput(-4.5,69)(196.5,0){2}{{\tiny{$c$}}}

\multiput(29,72)(130,0){2}{{\tiny{$e$}}}
\multiput(59,72)(70,0){2}{{\tiny{$e$}}}
\multiput(89,72)(10,0){2}{{\tiny{$c$}}}

\multiput(1.5,60.5)(84,0){2}{{\tiny{$a$}}}
\multiput(26,60.5)(35.5,0){2}{{\tiny{$f$}}}
\multiput(101.5,60.5)(84,0){2}{{\tiny{$a$}}}
\multiput(126,60.5)(35.5,0){2}{{\tiny{$f$}}}

\multiput(9.5,45)(68,0){2}{{\tiny{$b$}}}
\multiput(18,45.5)(51,0){2}{{\tiny{$g$}}}
\multiput(39.5,45)(00,0){1}{{\tiny{$h$}}}
\multiput(48,45)(00,0){1}{{\tiny{$h$}}}
\multiput(109.5,45)(68,0){2}{{\tiny{$b$}}}
\multiput(118,45.5)(51,0){2}{{\tiny{$g$}}}
\multiput(139.5,45)(00,0){1}{{\tiny{$h$}}}
\multiput(148,45)(00,0){1}{{\tiny{$h$}}}

\multiput(1.5,36.5)(84,0){2}{{\tiny{$d$}}}
\multiput(26.5,36.5)(35,0){2}{{\tiny{$i$}}}
\multiput(101.5,36.5)(84,0){2}{{\tiny{$d$}}}
\multiput(126.5,36.5)(35,0){2}{{\tiny{$i$}}}

\multiput(1.5,30.5)(84,0){2}{{\tiny{$d$}}}
\multiput(26.5,30.5)(35,0){2}{{\tiny{$i$}}}
\multiput(101.5,30.5)(84.5,0){2}{{\tiny{$d$}}}
\multiput(126.5,30.5)(35,0){2}{{\tiny{$i$}}}

\multiput(9.5,22)(69,0){2}{{\tiny{$b$}}}
\multiput(18.5,23)(51,0){2}{{\tiny{$g$}}}
\multiput(39.5,22)(0,0){2}{{\tiny{$h$}}}
\multiput(48,22)(0,0){2}{{\tiny{$h$}}}
\multiput(109.5,22)(69,0){2}{{\tiny{$b$}}}
\multiput(118.5,23)(51,0){2}{{\tiny{$g$}}}
\multiput(139.5,22)(0,0){2}{{\tiny{$h$}}}
\multiput(148,22)(0,0){2}{{\tiny{$h$}}}

\multiput(1.5,6.5)(84,0){2}{{\tiny{$a$}}}
\multiput(25.5,6.5)(35.5,0){2}{{\tiny{$f$}}}
\multiput(101.5,6.5)(84,0){2}{{\tiny{$a$}}}
\multiput(125.5,6.5)(35.5,0){2}{{\tiny{$f$}}}

\multiput(-4.5,-1)(196.5,0){2}{{\tiny{$c$}}}

\multiput(29,-4)(130,0){2}{{\tiny{$e$}}}
\multiput(59,-4)(70,0){2}{{\tiny{$e$}}}
\multiput(89,-4)(10,0){2}{{\tiny{$c$}}}

\end{picture}
}

\newcommand{\TotalCHHecThreeFifteen}{
\setlength{\unitlength}{0.5mm}
\begin{picture}(210,100)(-10,-15)

\multiput(20,60)(40,0){2}{{\circle{2}}}
\multiput(21,60)(40,0){2}{{\line(1,0){8}}}
\multiput(30,60)(40,0){2}{{\circle{2}}}
\multiput(31,60)(100,0){2}{{\line(1,0){28}}}
\multiput(71,60)(00,0){1}{{\line(1,0){48}}}
\multiput(120,60)(40,0){2}{{\circle{2}}}
\multiput(121,60)(40,0){2}{{\line(1,0){8}}}
\multiput(130,60)(40,0){2}{{\circle{2}}}

\multiput(20,50)(40,0){2}{{\circle{2}}}
\multiput(20,51)(40,0){2}{{\line(0,1){8}}}
\multiput(21,50)(40,0){2}{{\line(1,0){8}}}
\multiput(30,50)(40,0){2}{{\circle{2}}}
\multiput(30,51)(40,0){2}{{\line(0,1){8}}}
\multiput(120,50)(40,0){2}{{\circle{2}}}
\multiput(120,51)(40,0){2}{{\line(0,1){8}}}
\multiput(121,50)(40,0){2}{{\line(1,0){8}}}
\multiput(130,50)(40,0){2}{{\circle{2}}}
\multiput(130,51)(40,0){2}{{\line(0,1){8}}}

\multiput(10,40)(40,0){2}{{\circle{2}}}
\multiput(10.8,40.8)(40,0){2}{{\line(1,1){8.4}}}
\multiput(40,40)(40,0){2}{{\circle{2}}}
\multiput(39.2,40.8)(40,0){2}{{\line(-1,1){8.4}}}
\multiput(41,40)(100,0){2}{{\line(1,0){8}}}
\multiput(110,40)(40,0){2}{{\circle{2}}}
\multiput(110.8,40.8)(40,0){2}{{\line(1,1){8.4}}}
\multiput(140,40)(40,0){2}{{\circle{2}}}
\multiput(139.2,40.8)(40,0){2}{{\line(-1,1){8.4}}}

\multiput(00,30)(100,0){2}{{\circle{2}}}
\multiput(0.8,30.8)(100,0){2}{{\line(1,1){8.4}}}
\multiput(20,30)(40,0){2}{{\circle{2}}}
\multiput(19.2,30.8)(40,0){2}{{\line(-1,1){8.4}}}
\multiput(21,30)(40,0){2}{{\line(1,0){8}}}
\multiput(30,30)(40,0){2}{{\circle{2}}}
\multiput(30.8,30.8)(40,0){2}{{\line(1,1){8.4}}}
\multiput(120,30)(40,0){2}{{\circle{2}}}
\multiput(119.2,30.8)(40,0){2}{{\line(-1,1){8.4}}}
\multiput(121,30)(40,0){2}{{\line(1,0){8}}}
\multiput(130,30)(40,0){2}{{\circle{2}}}
\multiput(130.8,30.8)(40,0){2}{{\line(1,1){8.4}}}
\multiput(90,30)(100,0){2}{{\circle{2}}}
\multiput(89.2,30.8)(100,0){2}{{\line(-1,1){8.4}}}
\multiput(91,30)(00,0){1}{{\line(1,0){8}}}

\multiput(10,20)(40,0){2}{{\circle{2}}}
\multiput(09.2,20.8)(100,0){2}{{\line(-1,1){8.4}}}
\multiput(10.8,20.8)(40,0){2}{{\line(1,1){8.4}}}
\multiput(40,20)(40,0){2}{{\circle{2}}}
\multiput(39.2,20.8)(40,0){2}{{\line(-1,1){8.4}}}
\multiput(41,20)(100,0){2}{{\line(1,0){8}}}
\multiput(80.8,20.8)(100,0){2}{{\line(1,1){8.4}}}
\multiput(110,20)(40,0){2}{{\circle{2}}}
\multiput(110.8,20.8)(40,0){2}{{\line(1,1){8.4}}}
\multiput(140,20)(40,0){2}{{\circle{2}}}
\multiput(139.2,20.8)(40,0){2}{{\line(-1,1){8.4}}}

\multiput(20,10)(40,0){2}{{\circle{2}}}
\multiput(19.2,10.8)(40,0){2}{{\line(-1,1){8.4}}}
\multiput(21,10)(40,0){2}{{\line(1,0){8}}}
\multiput(30,10)(40,0){2}{{\circle{2}}}
\multiput(30.8,10.8)(40,0){2}{{\line(1,1){8.4}}}
\multiput(120,10)(40,0){2}{{\circle{2}}}
\multiput(119.2,10.8)(40,0){2}{{\line(-1,1){8.4}}}
\multiput(121,10)(40,0){2}{{\line(1,0){8}}}
\multiput(130,10)(40,0){2}{{\circle{2}}}
\multiput(130.8,10.8)(40,0){2}{{\line(1,1){8.4}}}

\multiput(20,00)(40,0){2}{{\circle{2}}}
\multiput(20,01)(40,0){2}{{\line(0,1){8}}}
\multiput(21,00)(40,0){2}{{\line(1,0){8}}}
\multiput(30,00)(40,0){2}{{\circle{2}}}
\multiput(30,01)(40,0){2}{{\line(0,1){8}}}
\multiput(31,00)(100,0){2}{{\line(1,0){28}}}
\multiput(71,00)(00,0){1}{{\line(1,0){48}}}
\multiput(120,00)(40,0){2}{{\circle{2}}}
\multiput(120,01)(40,0){2}{{\line(0,1){8}}}
\multiput(121,00)(40,0){2}{{\line(1,0){8}}}
\multiput(130,00)(40,0){2}{{\circle{2}}}
\multiput(130,01)(40,0){2}{{\line(0,1){8}}}

\multiput(20,60)(100,0){2}{\RIGHTlineTHICKsmall{8}}
\multiput(60,60)(100,0){2}{\RIGHTlineTHICKsmall{8}}

\multiput(20,50)(100,0){2}{\UPlineTHICKsmall{8}}
\multiput(30,50)(100,0){2}{\UPlineTHICKsmall{8}}
\multiput(60,50)(100,0){2}{\UPlineTHICKsmall{8}}
\multiput(70,50)(100,0){2}{\UPlineTHICKsmall{8}}

\multiput(10,40)(100,0){2}{\UPRIGHTlineTHICKsmall{8.4}}
\multiput(40,40)(100,0){2}{\UPLEFTlineTHICKsmall{8.4}}
\multiput(50,40)(100,0){2}{\UPRIGHTlineTHICKsmall{8.4}}
\multiput(80,40)(100,0){2}{\UPLEFTlineTHICKsmall{8.4}}

\multiput(0,30)(100,0){2}{\UPRIGHTlineTHICKsmall{8.4}}
\multiput(20,30)(100,0){2}{\RIGHTlineTHICKsmall{8}}

\multiput(30,30)(100,0){2}{\UPRIGHTlineTHICKsmall{8.4}}
\multiput(60,30)(100,0){2}{\UPLEFTlineTHICKsmall{8.4}}
\multiput(60,30)(100,0){2}{\RIGHTlineTHICKsmall{8}}
\multiput(90,30)(100,0){2}{\UPLEFTlineTHICKsmall{8.4}}
\multiput(90,30)(00,0){1}{\RIGHTlineTHICKsmall{8}}

\multiput(10,20)(100,0){2}{\UPRIGHTlineTHICKsmall{8.4}}
\multiput(40,20)(100,0){2}{\RIGHTlineTHICKsmall{8}}
\multiput(80,20)(100,0){2}{\UPLEFTlineTHICKsmall{8.4}}

\multiput(20,10)(100,0){2}{\UPLEFTlineTHICKsmall{8.4}}
\multiput(30,10)(100,0){2}{\UPRIGHTlineTHICKsmall{8.4}}
\multiput(60,10)(100,0){2}{\UPLEFTlineTHICKsmall{8.4}}
\multiput(70,10)(100,0){2}{\UPRIGHTlineTHICKsmall{8.4}}

\multiput(20,00)(100,0){2}{\UPlineTHICKsmall{8}}
\multiput(20,00)(100,0){2}{\RIGHTlineTHICKsmall{8}}
\multiput(30,00)(100,0){2}{\UPlineTHICKsmall{8}}
\multiput(60,00)(100,0){2}{\UPlineTHICKsmall{8}}
\multiput(60,00)(100,0){2}{\RIGHTlineTHICKsmall{8}}
\multiput(70,00)(100,0){2}{\UPlineTHICKsmall{8}}

\put(29,61){\oval(18,18)[tl]}
\put(29,70){\line(1,0){132}}
\put(161,61){\oval(18,18)[tr]}

\put(-1,40){\curveSouthWestTHICKsmall}
\put(191,40){\curveSouthEastTHICKsmall}
\multiput(-11,40)(212,0){2}{\properUPlineTHICKsmall{30}}
\put(-1,70){\curveNorthWestTHICKsmall}
\put(191,70){\curveNorthEastTHICKsmall}
\multiput(-1,80)(00,0){1}{\properRIGHTlineTHICKsmall{192}}

\put(29,-1){\oval(18,18)[bl]}
\put(29,-10){\line(1,0){132}}
\put(161,-1){\oval(18,18)[br]}

\multiput(94,81)(00,0){1}{\tiny{2}}

\multiput(94,71)(00,0){1}{\tiny{2}}

\multiput(24,61)(100,0){2}{\tiny{1}}
\multiput(64,61)(100,0){2}{\tiny{1}}
\multiput(44,61)(100,0){2}{\tiny{2}}
\multiput(94,61)(00,0){1}{\tiny{2}}

\multiput(16.5,54)(100,0){2}{{\tiny{3}}}
\multiput(31,54)(100,0){2}{{\tiny{3}}}
\multiput(56.5,54)(100,0){2}{{\tiny{3}}}
\multiput(71,54)(100,0){2}{{\tiny{3}}}

\multiput(24,51)(100,0){2}{\tiny{1}}
\multiput(64,51)(100,0){2}{\tiny{1}}

\multiput(11,45)(100,0){2}{{\tiny{2}}}
\multiput(36.5,45)(100,0){2}{{\tiny{2}}}
\multiput(51,45)(100,0){2}{{\tiny{2}}}
\multiput(76.5,45)(100,0){2}{{\tiny{2}}}

\multiput(44,36)(100,0){2}{\tiny{3}}

\multiput(01,35)(100,0){2}{{\tiny{3}}}
\multiput(86.5,35)(100,0){2}{{\tiny{3}}}

\multiput(11,32)(100,0){2}{{\tiny{1}}}
\multiput(36,32)(100,0){2}{{\tiny{1}}}
\multiput(51,32)(100,0){2}{{\tiny{1}}}
\multiput(76,32)(100,0){2}{{\tiny{1}}}

\multiput(5.5,25.5)(100,0){2}{{\tiny{1}}}
\multiput(81.5,25.5)(100,0){2}{{\tiny{1}}}

\multiput(24,25.5)(100,0){2}{\tiny{2}}
\multiput(64,25.5)(100,0){2}{\tiny{2}}
\multiput(94,31)(00,0){1}{\tiny{2}}

\multiput(16,22)(100,0){2}{{\tiny{3}}}
\multiput(31,22)(100,0){2}{{\tiny{3}}}
\multiput(56,22)(100,0){2}{{\tiny{3}}}
\multiput(71,22)(100,0){2}{{\tiny{3}}}

\multiput(44,15.5)(100,0){2}{\tiny{1}}

\multiput(11,12)(100,0){2}{{\tiny{2}}}
\multiput(36,12)(100,0){2}{{\tiny{2}}}
\multiput(51,12)(100,0){2}{{\tiny{2}}}
\multiput(76,12)(100,0){2}{{\tiny{2}}}

\multiput(24,6)(100,0){2}{\tiny{3}}
\multiput(64,6)(100,0){2}{\tiny{3}}

\multiput(16.5,4)(100,0){2}{{\tiny{1}}}
\multiput(31,4)(100,0){2}{{\tiny{1}}}
\multiput(56.5,4)(100,0){2}{{\tiny{1}}}
\multiput(71,4)(100,0){2}{{\tiny{1}}}

\multiput(24,-4)(100,0){2}{\tiny{3}}
\multiput(44,-4)(100,0){2}{\tiny{2}}
\multiput(64,-4)(100,0){2}{\tiny{3}}
\multiput(94,-4)(00,0){1}{\tiny{2}}

\multiput(94,-14)(00,0){1}{\tiny{2}}

\multiput(29,62)(100,0){2}{{\tiny{$d$}}}
\multiput(59,62)(100,0){2}{{\tiny{$d$}}}
\multiput(69,62)(00,0){1}{{\tiny{$e$}}}
\multiput(119,62)(00,0){1}{{\tiny{$e$}}}

\multiput(15.5,59)(156.5,0){2}{{\tiny{$e$}}}

\multiput(19.5,46)(100,0){2}{{\tiny{$c$}}}
\multiput(28,45)(100,0){2}{{\tiny{$b$}}}
\multiput(59.5,45)(100,0){2}{{\tiny{$b$}}}
\multiput(68,46)(100,0){2}{{\tiny{$c$}}}

\multiput(39.5,42)(100,0){2}{{\tiny{$a$}}}
\multiput(48,42)(100,0){2}{{\tiny{$a$}}}

\multiput(5.5,40)(100,0){2}{{\tiny{$c$}}}
\multiput(81.5,40)(100,0){2}{{\tiny{$c$}}}

\multiput(19.5,32)(100,0){2}{{\tiny{$b$}}}
\multiput(28,32)(100,0){2}{{\tiny{$a$}}}
\multiput(59.5,32)(100,0){2}{{\tiny{$a$}}}
\multiput(68,32)(100,0){2}{{\tiny{$b$}}}

\multiput(-2,26)(100,0){2}{{\tiny{$e$}}}
\multiput(89,26)(100,0){2}{{\tiny{$e$}}}

\multiput(5.5,18)(100,0){2}{{\tiny{$d$}}}
\multiput(39.5,22)(100,0){2}{{\tiny{$a$}}}
\multiput(48,22)(100,0){2}{{\tiny{$a$}}}
\multiput(81.5,18)(100,0){2}{{\tiny{$d$}}}

\multiput(19.5,12)(100,0){2}{{\tiny{$d$}}}
\multiput(28,12)(100,0){2}{{\tiny{$b$}}}
\multiput(59.5,12)(100,0){2}{{\tiny{$b$}}}
\multiput(68,12)(100,0){2}{{\tiny{$d$}}}

\multiput(15.5,-1)(156.5,0){2}{{\tiny{$e$}}}

\multiput(29,-4)(100,0){2}{{\tiny{$c$}}}
\multiput(59,-4)(100,0){2}{{\tiny{$c$}}}
\multiput(69,-4)(00,0){1}{{\tiny{$e$}}}
\multiput(119,-4)(00,0){1}{{\tiny{$e$}}}

\end{picture}
}

\newcommand{\spTotalCHHecThreeFifteen}{
\setlength{\unitlength}{0.5mm}
\begin{picture}(210,100)(-30,-15)

\multiput(20,60)(40,0){2}{{\circle{2}}}
\multiput(21,60)(40,0){2}{{\line(1,0){8}}}
\multiput(30,60)(40,0){2}{{\circle{2}}}
\multiput(31,60)(100,0){2}{{\line(1,0){28}}}
\multiput(71,60)(00,0){1}{{\line(1,0){48}}}
\multiput(120,60)(40,0){2}{{\circle{2}}}
\multiput(121,60)(40,0){2}{{\line(1,0){8}}}
\multiput(130,60)(40,0){2}{{\circle{2}}}

\multiput(20,50)(40,0){2}{{\circle{2}}}
\multiput(20,51)(40,0){2}{{\line(0,1){8}}}
\multiput(21,50)(40,0){2}{{\line(1,0){8}}}
\multiput(30,50)(40,0){2}{{\circle{2}}}
\multiput(30,51)(40,0){2}{{\line(0,1){8}}}
\multiput(120,50)(40,0){2}{{\circle{2}}}
\multiput(120,51)(40,0){2}{{\line(0,1){8}}}
\multiput(121,50)(40,0){2}{{\line(1,0){8}}}
\multiput(130,50)(40,0){2}{{\circle{2}}}
\multiput(130,51)(40,0){2}{{\line(0,1){8}}}

\multiput(10,40)(40,0){2}{{\circle{2}}}
\multiput(10.8,40.8)(40,0){2}{{\line(1,1){8.4}}}
\multiput(40,40)(40,0){2}{{\circle{2}}}
\multiput(39.2,40.8)(40,0){2}{{\line(-1,1){8.4}}}
\multiput(41,40)(100,0){2}{{\line(1,0){8}}}
\multiput(110,40)(40,0){2}{{\circle{2}}}
\multiput(110.8,40.8)(40,0){2}{{\line(1,1){8.4}}}
\multiput(140,40)(40,0){2}{{\circle{2}}}
\multiput(139.2,40.8)(40,0){2}{{\line(-1,1){8.4}}}

\multiput(00,30)(100,0){2}{{\circle{2}}}
\multiput(0.8,30.8)(100,0){2}{{\line(1,1){8.4}}}
\multiput(20,30)(40,0){2}{{\circle{2}}}
\multiput(19.2,30.8)(40,0){2}{{\line(-1,1){8.4}}}
\multiput(21,30)(40,0){2}{{\line(1,0){8}}}
\multiput(30,30)(40,0){2}{{\circle{2}}}
\multiput(30.8,30.8)(40,0){2}{{\line(1,1){8.4}}}
\multiput(120,30)(40,0){2}{{\circle{2}}}
\multiput(119.2,30.8)(40,0){2}{{\line(-1,1){8.4}}}
\multiput(121,30)(40,0){2}{{\line(1,0){8}}}
\multiput(130,30)(40,0){2}{{\circle{2}}}
\multiput(130.8,30.8)(40,0){2}{{\line(1,1){8.4}}}
\multiput(90,30)(100,0){2}{{\circle{2}}}
\multiput(89.2,30.8)(100,0){2}{{\line(-1,1){8.4}}}
\multiput(91,30)(00,0){1}{{\line(1,0){8}}}

\multiput(10,20)(40,0){2}{{\circle{2}}}
\multiput(09.2,20.8)(100,0){2}{{\line(-1,1){8.4}}}
\multiput(10.8,20.8)(40,0){2}{{\line(1,1){8.4}}}
\multiput(40,20)(40,0){2}{{\circle{2}}}
\multiput(39.2,20.8)(40,0){2}{{\line(-1,1){8.4}}}
\multiput(41,20)(100,0){2}{{\line(1,0){8}}}
\multiput(80.8,20.8)(100,0){2}{{\line(1,1){8.4}}}
\multiput(110,20)(40,0){2}{{\circle{2}}}
\multiput(110.8,20.8)(40,0){2}{{\line(1,1){8.4}}}
\multiput(140,20)(40,0){2}{{\circle{2}}}
\multiput(139.2,20.8)(40,0){2}{{\line(-1,1){8.4}}}

\multiput(20,10)(40,0){2}{{\circle{2}}}
\multiput(19.2,10.8)(40,0){2}{{\line(-1,1){8.4}}}
\multiput(21,10)(40,0){2}{{\line(1,0){8}}}
\multiput(30,10)(40,0){2}{{\circle{2}}}
\multiput(30.8,10.8)(40,0){2}{{\line(1,1){8.4}}}
\multiput(120,10)(40,0){2}{{\circle{2}}}
\multiput(119.2,10.8)(40,0){2}{{\line(-1,1){8.4}}}
\multiput(121,10)(40,0){2}{{\line(1,0){8}}}
\multiput(130,10)(40,0){2}{{\circle{2}}}
\multiput(130.8,10.8)(40,0){2}{{\line(1,1){8.4}}}

\multiput(20,00)(40,0){2}{{\circle{2}}}
\multiput(20,01)(40,0){2}{{\line(0,1){8}}}
\multiput(21,00)(40,0){2}{{\line(1,0){8}}}
\multiput(30,00)(40,0){2}{{\circle{2}}}
\multiput(30,01)(40,0){2}{{\line(0,1){8}}}
\multiput(31,00)(100,0){2}{{\line(1,0){28}}}
\multiput(71,00)(00,0){1}{{\line(1,0){48}}}
\multiput(120,00)(40,0){2}{{\circle{2}}}
\multiput(120,01)(40,0){2}{{\line(0,1){8}}}
\multiput(121,00)(40,0){2}{{\line(1,0){8}}}
\multiput(130,00)(40,0){2}{{\circle{2}}}
\multiput(130,01)(40,0){2}{{\line(0,1){8}}}

\multiput(20,60)(100,0){2}{\RIGHTlineTHICKsmall{8}}
\multiput(30,60)(100,0){2}{\RIGHTlineTHICKsmall{28}}
\multiput(70,60)(00,0){1}{\RIGHTlineTHICKsmall{48}}
\multiput(160,60)(00,0){1}{\RIGHTlineTHICKsmall{8}}

\multiput(20,50)(100,0){2}{\RIGHTlineTHICKsmall{8}}
\multiput(60,50)(00,0){1}{\UPlineTHICKsmall{8}}
\multiput(70,50)(00,0){1}{\UPlineTHICKsmall{8}}
\multiput(160,50)(00,0){1}{\RIGHTlineTHICKsmall{8}}

\multiput(10,40)(100,0){2}{\UPRIGHTlineTHICKsmall{8.4}}
\multiput(40,40)(100,0){2}{\UPLEFTlineTHICKsmall{8.4}}
\multiput(50,40)(100,0){2}{\UPRIGHTlineTHICKsmall{8.4}}
\multiput(80,40)(100,0){2}{\UPLEFTlineTHICKsmall{8.4}}
\multiput(140,40)(00,0){1}{\RIGHTlineTHICKsmall{8}}

\multiput(0,30)(30,0){2}{\UPRIGHTlineTHICKsmall{8.4}}
\multiput(20,30)(100,0){2}{\RIGHTlineTHICKsmall{8}}
\multiput(60,30)(30,0){3}{\UPLEFTlineTHICKsmall{8.4}}
\multiput(60,30)(100,0){2}{\RIGHTlineTHICKsmall{8}}
\multiput(90,30)(00,0){1}{\RIGHTlineTHICKsmall{8}}
\multiput(170,30)(00,0){1}{\UPRIGHTlineTHICKsmall{8.4}}

\multiput(10,20)(140,0){2}{\UPRIGHTlineTHICKsmall{8.4}}
\multiput(40,20)(00,0){1}{\RIGHTlineTHICKsmall{8}}
\multiput(80,20)(30,0){3}{\UPLEFTlineTHICKsmall{8.4}}
\multiput(180,20)(00,0){1}{\UPRIGHTlineTHICKsmall{8.4}}

\multiput(20,10)(100,0){2}{\UPLEFTlineTHICKsmall{8.8}}
\multiput(20,10)(40,0){2}{\RIGHTlineTHICKsmall{8}}
\multiput(30,10)(100,0){2}{\UPRIGHTlineTHICKsmall{8.4}}
\multiput(60,10)(100,0){2}{\UPLEFTlineTHICKsmall{8.4}}
\multiput(70,10)(100,0){2}{\UPRIGHTlineTHICKsmall{8.4}}
\multiput(160,10)(00,0){1}{\RIGHTlineTHICKsmall{8}}

\multiput(20,00)(00,0){1}{\RIGHTlineTHICKsmall{8}}
\multiput(30,00)(100,0){2}{\RIGHTlineTHICKsmall{28}}
\multiput(60,00)(100,0){2}{\RIGHTlineTHICKsmall{8}}
\multiput(70,00)(00,0){1}{\RIGHTlineTHICKsmall{48}}
\multiput(120,00)(00,0){1}{\UPlineTHICKsmall{8}}
\multiput(130,00)(00,0){1}{\UPlineTHICKsmall{8}}

\put(30,61){\curveNorthWestTHICKsmall}
\put(30,71){\properRIGHTlineTHICKsmall{130}}
\put(160,61){\curveNorthEastTHICKsmall}

\put(-1,40){\curveSouthWestTHICKsmall}
\put(191,40){\curveSouthEastTHICKsmall}
\multiput(-11,40)(212,0){2}{\properUPlineTHICKsmall{30}}
\put(-1,70){\curveNorthWestTHICKsmall}
\put(191,70){\curveNorthEastTHICKsmall}
\multiput(-1,80)(00,0){1}{\properRIGHTlineTHICKsmall{192}}

\put(30,-1){\curveSouthWestTHICKsmall}
\put(30,-11){\properRIGHTlineTHICKsmall{130}}
\put(160,-1){\curveSouthEastTHICKsmall}

\multiput(94,81)(00,0){1}{\tiny{2}}

\multiput(94,72)(00,0){1}{\tiny{2}}

\multiput(24,61)(100,0){2}{\tiny{1}}
\multiput(64,61)(100,0){2}{\tiny{1}}
\multiput(44,61)(100,0){2}{\tiny{2}}
\multiput(94,61)(00,0){1}{\tiny{2}}

\multiput(16.5,54)(100,0){2}{{\tiny{3}}}
\multiput(31,54)(100,0){2}{{\tiny{3}}}
\multiput(56.5,54)(100,0){2}{{\tiny{3}}}
\multiput(71,54)(100,0){2}{{\tiny{3}}}

\multiput(24,51)(100,0){2}{\tiny{1}}
\multiput(64,51)(100,0){2}{\tiny{1}}

\multiput(11,45)(100,0){2}{{\tiny{2}}}
\multiput(36.5,45)(100,0){2}{{\tiny{2}}}
\multiput(51,45)(100,0){2}{{\tiny{2}}}
\multiput(76.5,45)(100,0){2}{{\tiny{2}}}

\multiput(44,36)(100,0){2}{\tiny{3}}

\multiput(01,35)(100,0){2}{{\tiny{3}}}
\multiput(86.5,35)(100,0){2}{{\tiny{3}}}

\multiput(11,32)(100,0){2}{{\tiny{1}}}
\multiput(36,32)(100,0){2}{{\tiny{1}}}
\multiput(51,32)(100,0){2}{{\tiny{1}}}
\multiput(76,32)(100,0){2}{{\tiny{1}}}

\multiput(5.5,25.5)(100,0){2}{{\tiny{1}}}
\multiput(81.5,25.5)(100,0){2}{{\tiny{1}}}

\multiput(24,25.5)(100,0){2}{\tiny{2}}
\multiput(64,25.5)(100,0){2}{\tiny{2}}
\multiput(94,31)(00,0){1}{\tiny{2}}

\multiput(16,22)(100,0){2}{{\tiny{3}}}
\multiput(31,22)(100,0){2}{{\tiny{3}}}
\multiput(56,22)(100,0){2}{{\tiny{3}}}
\multiput(71,22)(100,0){2}{{\tiny{3}}}

\multiput(44,15.5)(100,0){2}{\tiny{1}}

\multiput(11,12)(100,0){2}{{\tiny{2}}}
\multiput(36,12)(100,0){2}{{\tiny{2}}}
\multiput(51,12)(100,0){2}{{\tiny{2}}}
\multiput(76,12)(100,0){2}{{\tiny{2}}}

\multiput(24,6)(100,0){2}{\tiny{3}}
\multiput(64,6)(100,0){2}{\tiny{3}}

\multiput(16.5,4)(100,0){2}{{\tiny{1}}}
\multiput(31,4)(100,0){2}{{\tiny{1}}}
\multiput(56.5,4)(100,0){2}{{\tiny{1}}}
\multiput(71,4)(100,0){2}{{\tiny{1}}}

\multiput(24,-4)(100,0){2}{\tiny{3}}
\multiput(44,-4)(100,0){2}{\tiny{2}}
\multiput(64,-4)(100,0){2}{\tiny{3}}
\multiput(94,-4)(00,0){1}{\tiny{2}}

\multiput(94,-15)(00,0){1}{\tiny{2}}

\multiput(29,62)(100,0){2}{{\tiny{$d$}}}
\multiput(59,62)(100,0){2}{{\tiny{$d$}}}
\multiput(69,62)(00,0){1}{{\tiny{$e$}}}
\multiput(119,62)(00,0){1}{{\tiny{$e$}}}

\multiput(15.5,59)(156.5,0){2}{{\tiny{$e$}}}

\multiput(19.5,46)(100,0){2}{{\tiny{$c$}}}
\multiput(28,45)(100,0){2}{{\tiny{$b$}}}
\multiput(59.5,45)(100,0){2}{{\tiny{$b$}}}
\multiput(68,46)(100,0){2}{{\tiny{$c$}}}

\multiput(39.5,42)(100,0){2}{{\tiny{$a$}}}
\multiput(48,42)(100,0){2}{{\tiny{$a$}}}

\multiput(5.5,40)(100,0){2}{{\tiny{$c$}}}
\multiput(81.5,40)(100,0){2}{{\tiny{$c$}}}

\multiput(19.5,32)(100,0){2}{{\tiny{$b$}}}
\multiput(28,32)(100,0){2}{{\tiny{$a$}}}
\multiput(59.5,32)(100,0){2}{{\tiny{$a$}}}
\multiput(68,32)(100,0){2}{{\tiny{$b$}}}

\multiput(-2,26)(100,0){2}{{\tiny{$e$}}}
\multiput(89,26)(100,0){2}{{\tiny{$e$}}}

\multiput(5.5,18)(100,0){2}{{\tiny{$d$}}}
\multiput(39.5,22)(100,0){2}{{\tiny{$a$}}}
\multiput(48,22)(100,0){2}{{\tiny{$a$}}}
\multiput(81.5,18)(100,0){2}{{\tiny{$d$}}}

\multiput(19.5,12)(100,0){2}{{\tiny{$d$}}}
\multiput(28,12)(100,0){2}{{\tiny{$b$}}}
\multiput(59.5,12)(100,0){2}{{\tiny{$b$}}}
\multiput(68,12)(100,0){2}{{\tiny{$d$}}}

\multiput(15.5,-1)(156.5,0){2}{{\tiny{$e$}}}

\multiput(29,-4)(100,0){2}{{\tiny{$c$}}}
\multiput(59,-4)(100,0){2}{{\tiny{$c$}}}
\multiput(69,-4)(00,0){1}{{\tiny{$e$}}}
\multiput(119,-4)(00,0){1}{{\tiny{$e$}}}

\end{picture}
}

\newcommand{\TotalCHHecThreeNine}{
\setlength{\unitlength}{0.5mm}
\begin{picture}(70,100)(-25,-15)

\multiput(50,70)(10,0){2}{\circle{2}}
\multiput(51,70)(00,0){1}{\line(1,0){8}}

\multiput(50,60)(10,0){2}{\circle{2}}
\multiput(51,60)(00,0){1}{\line(1,0){8}}
\multiput(50,61)(00,0){1}{\line(0,1){8}}
\multiput(60,61)(00,0){1}{\line(0,1){8}}

\multiput(10,50)(30,0){4}{\circle{2}}
\multiput(40.8,50.8)(00,0){1}{\line(1,1){8.4}}
\multiput(69.2,50.8)(00,0){1}{\line(-1,1){8.4}}

\multiput(00,40)(30,0){4}{\circle{2}}
\multiput(00.8,40.8)(30,0){4}{\line(1,1){8.4}}
\multiput(20,40)(30,0){4}{\circle{2}}
\multiput(19.2,40.8)(30,0){4}{\line(-1,1){8.4}}
\multiput(21,40)(30,0){3}{\line(1,0){8}}

\multiput(00,30)(30,0){4}{\circle{2}}
\multiput(00,31)(30,0){4}{\line(0,1){8}}
\multiput(20,30)(30,0){4}{\circle{2}}
\multiput(20,31)(30,0){4}{\line(0,1){8}}
\multiput(21,30)(30,0){3}{\line(1,0){8}}

\multiput(10,20)(30,0){4}{\circle{2}}
\multiput(9.2,20.8)(30,0){4}{\line(-1,1){8.4}}
\multiput(10.8,20.8)(30,0){4}{\line(1,1){8.4}}

\multiput(50,10)(10,0){2}{\circle{2}}
\multiput(49.2,10.8)(00,0){1}{\line(-1,1){8.4}}
\multiput(51,10)(00,0){1}{\line(1,0){8}}
\multiput(60.8,10.8)(00,0){1}{\line(1,1){8.4}}

\multiput(50,00)(10,0){2}{\circle{2}}
\multiput(51,00)(00,0){1}{\line(1,0){8}}
\multiput(50,01)(00,0){1}{\line(0,1){8}}
\multiput(60,01)(00,0){1}{\line(0,1){8}}

\multiput(00,00)(110,0){2}{\line(0,1){29}}
\put(10,00){\oval(20,20)[bl]}
\put(100,00){\oval(20,20)[br]}
\multiput(10,-10)(00,0){1}{\line(1,0){90}}

\multiput(10,10)(90,0){2}{\line(0,1){9}}
\put(20,10){\oval(20,20)[bl]}
\put(90,10){\oval(20,20)[br]}
\multiput(20,00)(41,0){2}{\line(1,0){29}}

\multiput(00,41)(110,0){2}{\properUPlineTHICKsmall{29}}
\put(10,70){\curveNorthWestTHICKsmall}
\put(100,70){\curveNorthEastTHICKsmall}
\multiput(10,80)(00,0){1}{\properRIGHTlineTHICKsmall{90}}

\put(20,60){\curveNorthWestTHICKsmall}
\multiput(10,51)(90,0){2}{\properUPlineTHICKsmall{9}}
\put(90,60){\curveNorthEastTHICKsmall}
\multiput(20,70)(41,0){2}{\properRIGHTlineTHICKsmall{29}}

\multiput(50,60)(00,0){1}{\RIGHTlineTHICKsmall{8}}
\multiput(50,60)(10,0){2}{\UPlineTHICKsmall{8}}

\multiput(20,40)(30,0){3}{\UPLEFTlineTHICKsmall{8.4}}
\multiput(30,40)(30,0){3}{\UPRIGHTlineTHICKsmall{8.4}}
\multiput(80,40)(00,0){1}{\RIGHTlineTHICKsmall{8}}

\multiput(00,30)(30,0){3}{\UPlineTHICKsmall{8}}
\multiput(20,30)(30,0){2}{\UPlineTHICKsmall{8}}
\multiput(50,30)(30,0){2}{\RIGHTlineTHICKsmall{8}}
\multiput(110,30)(00,0){1}{\UPlineTHICKsmall{8}}

\multiput(10,20)(30,0){2}{\UPLEFTlineTHICKsmall{8.4}}
\multiput(10,20)(60,0){2}{\UPRIGHTlineTHICKsmall{8.4}}
\multiput(100,20)(00,0){1}{\UPLEFTlineTHICKsmall{8.4}}
\multiput(100,20)(00,0){1}{\UPRIGHTlineTHICKsmall{8.4}}

\multiput(50,10)(00,0){1}{\UPLEFTlineTHICKsmall{8.4}}
\multiput(60,10)(00,0){1}{\UPRIGHTlineTHICKsmall{8.4}}

\multiput(50,00)(10,0){2}{\UPlineTHICKsmall{8}}
\multiput(50,00)(00,0){1}{\RIGHTlineTHICKsmall{8}}

\multiput(54,81)(00,0){1}{\tiny{1}}

\multiput(29,71)(52,0){2}{\tiny{1}}
\multiput(54,71)(00,0){1}{\tiny{2}}

\multiput(46.5,63.5)(00,0){1}{{\tiny{3}}}
\multiput(61,63.5)(00,0){1}{{\tiny{3}}}

\multiput(54,61)(00,0){1}{\tiny{2}}

\multiput(41.5,55)(00,0){1}{{\tiny{1}}}
\multiput(66,55)(00,0){1}{{\tiny{1}}}

\multiput(05,41.5)(60,0){2}{{\tiny{2}}}
\multiput(12,41.5)(60,0){2}{{\tiny{3}}}

\multiput(35,41.5)(60,0){2}{{\tiny{3}}}
\multiput(42,41.5)(60,0){2}{{\tiny{2}}}

\multiput(24,36)(30,0){3}{{\tiny{1}}}

\multiput(1,33.5)(60,0){2}{{\tiny{3}}}
\multiput(16.5,33.5)(60,0){2}{{\tiny{2}}}
\multiput(31,33.5)(60,0){2}{{\tiny{2}}}
\multiput(46.5,33.5)(60,0){2}{{\tiny{3}}}

\multiput(24,31)(30,0){3}{{\tiny{1}}}

\multiput(05,26)(60,0){2}{{\tiny{2}}}
\multiput(42.5,26)(60,0){2}{{\tiny{2}}}

\multiput(11.5,26)(60,0){2}{{\tiny{3}}}
\multiput(35,26)(60,0){2}{{\tiny{3}}}

\multiput(41.5,11.5)(00,0){1}{{\tiny{1}}}
\multiput(66,11.5)(00,0){1}{{\tiny{1}}}

\multiput(54,6)(00,0){1}{\tiny{2}}

\multiput(46.5,3.5)(00,0){1}{{\tiny{3}}}
\multiput(61,3.5)(00,0){1}{{\tiny{3}}}

\multiput(29,-4)(52,0){2}{\tiny{1}}
\multiput(54,-4)(00,0){1}{\tiny{2}}

\multiput(54,-14)(00,0){1}{\tiny{1}}

\multiput(49,72)(10,0){2}{{\tiny{$c$}}}

\multiput(49.5,56)(00,0){1}{{\tiny{$c$}}}
\multiput(58,56)(00,0){1}{{\tiny{$c$}}}

\multiput(11.5,50)(60,0){2}{{\tiny{$a$}}}
\multiput(35.5,50)(60,0){2}{{\tiny{$a$}}}

\multiput(19.5,42)(60,0){2}{{\tiny{$d$}}}
\multiput(27,42)(60,0){2}{{\tiny{$d$}}}
\multiput(49.5,42)(00,0){1}{{\tiny{$b$}}}
\multiput(58,42)(00,0){1}{{\tiny{$b$}}}

\multiput(-4.5,39)(116.5,0){2}{{\tiny{$b$}}}

\multiput(-4.5,29)(116.5,0){2}{{\tiny{$b$}}}

\multiput(19.5,25)(60,0){2}{{\tiny{$d$}}}
\multiput(27,25)(60,0){2}{{\tiny{$d$}}}
\multiput(49.5,25)(00,0){1}{{\tiny{$b$}}}
\multiput(58,25)(00,0){1}{{\tiny{$b$}}}

\multiput(11.5,18)(60,0){2}{{\tiny{$a$}}}
\multiput(35.5,18)(60,0){2}{{\tiny{$a$}}}

\multiput(49.5,12)(00,0){1}{{\tiny{$c$}}}
\multiput(58,12)(00,0){1}{{\tiny{$c$}}}

\multiput(49,-4)(10,0){2}{{\tiny{$c$}}}

\end{picture}
}

\newcommand{\TotalCHHecThreeNineForHecFourThirteen}{
\setlength{\unitlength}{0.5mm}
\begin{picture}(110,100)(00,-15)

\multiput(50,70)(10,0){2}{\circle{2}}
\multiput(51,70)(00,0){1}{\line(1,0){8}}

\multiput(50,60)(10,0){2}{\circle{2}}
\multiput(51,60)(00,0){1}{\line(1,0){8}}
\multiput(50,61)(00,0){1}{\line(0,1){8}}
\multiput(60,61)(00,0){1}{\line(0,1){8}}

\multiput(10,50)(30,0){4}{\circle{2}}
\multiput(40.8,50.8)(00,0){1}{\line(1,1){8.4}}
\multiput(69.2,50.8)(00,0){1}{\line(-1,1){8.4}}

\multiput(00,40)(30,0){4}{\circle{2}}
\multiput(00.8,40.8)(30,0){4}{\line(1,1){8.4}}
\multiput(20,40)(30,0){4}{\circle{2}}
\multiput(19.2,40.8)(30,0){4}{\line(-1,1){8.4}}
\multiput(21,40)(30,0){3}{\line(1,0){8}}

\multiput(00,30)(30,0){4}{\circle{2}}
\multiput(00,31)(30,0){4}{\line(0,1){8}}
\multiput(20,30)(30,0){4}{\circle{2}}
\multiput(20,31)(30,0){4}{\line(0,1){8}}
\multiput(21,30)(30,0){3}{\line(1,0){8}}

\multiput(10,20)(30,0){4}{\circle{2}}
\multiput(9.2,20.8)(30,0){4}{\line(-1,1){8.4}}
\multiput(10.8,20.8)(30,0){4}{\line(1,1){8.4}}

\multiput(50,10)(10,0){2}{\circle{2}}
\multiput(49.2,10.8)(00,0){1}{\line(-1,1){8.4}}
\multiput(51,10)(00,0){1}{\line(1,0){8}}
\multiput(60.8,10.8)(00,0){1}{\line(1,1){8.4}}

\multiput(50,00)(10,0){2}{\circle{2}}
\multiput(51,00)(00,0){1}{\line(1,0){8}}
\multiput(50,01)(00,0){1}{\line(0,1){8}}
\multiput(60,01)(00,0){1}{\line(0,1){8}}

\multiput(00,00)(110,0){2}{\line(0,1){29}}
\put(10,00){\oval(20,20)[bl]}
\put(100,00){\oval(20,20)[br]}
\multiput(10,-10)(00,0){1}{\line(1,0){90}}

\multiput(10,10)(90,0){2}{\line(0,1){9}}
\put(20,10){\oval(20,20)[bl]}
\put(90,10){\oval(20,20)[br]}
\multiput(20,00)(41,0){2}{\line(1,0){29}}

\multiput(00,41)(110,0){2}{\properUPlineTHICKsmall{29}}
\put(10,70){\curveNorthWestTHICKsmall}
\put(100,70){\curveNorthEastTHICKsmall}
\multiput(10,80)(00,0){1}{\properRIGHTlineTHICKsmall{90}}

\put(20,60){\curveNorthWestTHICKsmall}
\multiput(10,51)(90,0){2}{\properUPlineTHICKsmall{9}}
\put(90,60){\curveNorthEastTHICKsmall}
\multiput(20,70)(41,0){2}{\properRIGHTlineTHICKsmall{29}}

\multiput(50,60)(00,0){1}{\RIGHTlineTHICKsmall{8}}
\multiput(50,60)(10,0){2}{\UPlineTHICKsmall{8}}

\multiput(20,40)(30,0){3}{\UPLEFTlineTHICKsmall{8.4}}
\multiput(30,40)(30,0){3}{\UPRIGHTlineTHICKsmall{8.4}}
\multiput(80,40)(00,0){1}{\RIGHTlineTHICKsmall{8}}

\multiput(00,30)(30,0){3}{\UPlineTHICKsmall{8}}
\multiput(20,30)(30,0){2}{\UPlineTHICKsmall{8}}
\multiput(50,30)(30,0){2}{\RIGHTlineTHICKsmall{8}}
\multiput(110,30)(00,0){1}{\UPlineTHICKsmall{8}}

\multiput(10,20)(30,0){2}{\UPLEFTlineTHICKsmall{8.4}}
\multiput(10,20)(60,0){2}{\UPRIGHTlineTHICKsmall{8.4}}
\multiput(100,20)(00,0){1}{\UPLEFTlineTHICKsmall{8.4}}
\multiput(100,20)(00,0){1}{\UPRIGHTlineTHICKsmall{8.4}}

\multiput(50,10)(00,0){1}{\UPLEFTlineTHICKsmall{8.4}}
\multiput(60,10)(00,0){1}{\UPRIGHTlineTHICKsmall{8.4}}

\multiput(50,00)(10,0){2}{\UPlineTHICKsmall{8}}
\multiput(50,00)(00,0){1}{\RIGHTlineTHICKsmall{8}}

\multiput(54,81)(00,0){1}{\tiny{3}}

\multiput(29,71)(52,0){2}{\tiny{3}}
\multiput(54,71)(00,0){1}{\tiny{4}}

\multiput(46.5,63.5)(00,0){1}{{\tiny{2}}}
\multiput(61,63.5)(00,0){1}{{\tiny{2}}}

\multiput(54,61)(00,0){1}{\tiny{4}}

\multiput(41.5,55)(00,0){1}{{\tiny{3}}}
\multiput(66,55)(00,0){1}{{\tiny{3}}}

\multiput(05,41.5)(60,0){2}{{\tiny{4}}}
\multiput(12,41.5)(60,0){2}{{\tiny{2}}}

\multiput(35,41.5)(60,0){2}{{\tiny{2}}}
\multiput(42,41.5)(60,0){2}{{\tiny{2}}}

\multiput(24,36)(30,0){3}{{\tiny{3}}}

\multiput(1,33.5)(60,0){2}{{\tiny{2}}}
\multiput(16.5,33.5)(60,0){2}{{\tiny{4}}}
\multiput(31,33.5)(60,0){2}{{\tiny{4}}}
\multiput(46.5,33.5)(60,0){2}{{\tiny{2}}}

\multiput(24,31)(30,0){3}{{\tiny{3}}}

\multiput(05,26)(60,0){2}{{\tiny{4}}}
\multiput(42.5,26)(60,0){2}{{\tiny{4}}}

\multiput(11.5,26)(60,0){2}{{\tiny{2}}}
\multiput(35,26)(60,0){2}{{\tiny{2}}}

\multiput(41.5,11.5)(00,0){1}{{\tiny{3}}}
\multiput(66,11.5)(00,0){1}{{\tiny{3}}}

\multiput(54,6)(00,0){1}{{\tiny{4}}}

\multiput(46.5,3.5)(00,0){1}{{\tiny{2}}}
\multiput(61,3.5)(00,0){1}{{\tiny{2}}}

\multiput(29,-4)(52,0){2}{\tiny{3}}
\multiput(54,-4)(00,0){1}{\tiny{4}}

\multiput(54,-14)(00,0){1}{\tiny{3}}

\multiput(49,72)(10,0){2}{{\tiny{$a$}}}

\multiput(49.5,56)(00,0){1}{{\tiny{$a$}}}
\multiput(58,56)(00,0){1}{{\tiny{$a$}}}

\multiput(11.5,50)(60,0){2}{{\tiny{$b$}}}
\multiput(35.5,50)(60,0){2}{{\tiny{$b$}}}

\multiput(19.5,42)(60,0){2}{{\tiny{$d$}}}
\multiput(27,42)(60,0){2}{{\tiny{$d$}}}
\multiput(49.5,42)(00,0){1}{{\tiny{$c$}}}
\multiput(58,42)(00,0){1}{{\tiny{$c$}}}

\multiput(-4.5,39)(116.5,0){2}{{\tiny{$c$}}}

\multiput(-4.5,29)(116.5,0){2}{{\tiny{$c$}}}

\multiput(19.5,25)(60,0){2}{{\tiny{$d$}}}
\multiput(27,25)(60,0){2}{{\tiny{$d$}}}
\multiput(49.5,26)(00,0){1}{{\tiny{$c$}}}
\multiput(58,26)(00,0){1}{{\tiny{$c$}}}

\multiput(11.5,18)(60,0){2}{{\tiny{$b$}}}
\multiput(35.5,18)(60,0){2}{{\tiny{$b$}}}

\multiput(49.5,12)(00,0){1}{{\tiny{$a$}}}
\multiput(58,12)(00,0){1}{{\tiny{$a$}}}

\multiput(49,-4)(10,0){2}{{\tiny{$a$}}}

\end{picture}
}

\newcommand{\TotalCHHecThreeNineForHecFourTwetytwo}{
\setlength{\unitlength}{0.5mm}
\begin{picture}(110,100)(00,-15)

\multiput(50,70)(10,0){2}{\circle{2}}
\multiput(51,70)(00,0){1}{\line(1,0){8}}

\multiput(50,60)(10,0){2}{\circle{2}}
\multiput(51,60)(00,0){1}{\line(1,0){8}}
\multiput(50,61)(00,0){1}{\line(0,1){8}}
\multiput(60,61)(00,0){1}{\line(0,1){8}}

\multiput(10,50)(30,0){4}{\circle{2}}
\multiput(40.8,50.8)(00,0){1}{\line(1,1){8.4}}
\multiput(69.2,50.8)(00,0){1}{\line(-1,1){8.4}}

\multiput(00,40)(30,0){4}{\circle{2}}
\multiput(00.8,40.8)(30,0){4}{\line(1,1){8.4}}
\multiput(20,40)(30,0){4}{\circle{2}}
\multiput(19.2,40.8)(30,0){4}{\line(-1,1){8.4}}
\multiput(21,40)(30,0){3}{\line(1,0){8}}

\multiput(00,30)(30,0){4}{\circle{2}}
\multiput(00,31)(30,0){4}{\line(0,1){8}}
\multiput(20,30)(30,0){4}{\circle{2}}
\multiput(20,31)(30,0){4}{\line(0,1){8}}
\multiput(21,30)(30,0){3}{\line(1,0){8}}

\multiput(10,20)(30,0){4}{\circle{2}}
\multiput(9.2,20.8)(30,0){4}{\line(-1,1){8.4}}
\multiput(10.8,20.8)(30,0){4}{\line(1,1){8.4}}

\multiput(50,10)(10,0){2}{\circle{2}}
\multiput(49.2,10.8)(00,0){1}{\line(-1,1){8.4}}
\multiput(51,10)(00,0){1}{\line(1,0){8}}
\multiput(60.8,10.8)(00,0){1}{\line(1,1){8.4}}

\multiput(50,00)(10,0){2}{\circle{2}}
\multiput(51,00)(00,0){1}{\line(1,0){8}}
\multiput(50,01)(00,0){1}{\line(0,1){8}}
\multiput(60,01)(00,0){1}{\line(0,1){8}}

\multiput(00,00)(110,0){2}{\line(0,1){29}}
\put(10,00){\oval(20,20)[bl]}
\put(100,00){\oval(20,20)[br]}
\multiput(10,-10)(00,0){1}{\line(1,0){90}}

\multiput(10,10)(90,0){2}{\line(0,1){9}}
\put(20,10){\oval(20,20)[bl]}
\put(90,10){\oval(20,20)[br]}
\multiput(20,00)(41,0){2}{\line(1,0){29}}

\multiput(00,41)(110,0){2}{\properUPlineTHICKsmall{29}}
\put(10,70){\curveNorthWestTHICKsmall}
\put(100,70){\curveNorthEastTHICKsmall}
\multiput(10,80)(00,0){1}{\properRIGHTlineTHICKsmall{90}}

\put(20,60){\curveNorthWestTHICKsmall}
\multiput(10,51)(90,0){2}{\properUPlineTHICKsmall{9}}
\put(90,60){\curveNorthEastTHICKsmall}
\multiput(20,70)(41,0){2}{\properRIGHTlineTHICKsmall{29}}

\multiput(50,60)(00,0){1}{\RIGHTlineTHICKsmall{8}}
\multiput(50,60)(10,0){2}{\UPlineTHICKsmall{8}}

\multiput(20,40)(30,0){3}{\UPLEFTlineTHICKsmall{8.4}}
\multiput(30,40)(30,0){3}{\UPRIGHTlineTHICKsmall{8.4}}
\multiput(80,40)(00,0){1}{\RIGHTlineTHICKsmall{8}}

\multiput(00,30)(30,0){3}{\UPlineTHICKsmall{8}}
\multiput(20,30)(30,0){2}{\UPlineTHICKsmall{8}}
\multiput(50,30)(30,0){2}{\RIGHTlineTHICKsmall{8}}
\multiput(110,30)(00,0){1}{\UPlineTHICKsmall{8}}

\multiput(10,20)(30,0){2}{\UPLEFTlineTHICKsmall{8.4}}
\multiput(10,20)(60,0){2}{\UPRIGHTlineTHICKsmall{8.4}}
\multiput(100,20)(00,0){1}{\UPLEFTlineTHICKsmall{8.4}}
\multiput(100,20)(00,0){1}{\UPRIGHTlineTHICKsmall{8.4}}

\multiput(50,10)(00,0){1}{\UPLEFTlineTHICKsmall{8.4}}
\multiput(60,10)(00,0){1}{\UPRIGHTlineTHICKsmall{8.4}}

\multiput(50,00)(10,0){2}{\UPlineTHICKsmall{8}}
\multiput(50,00)(00,0){1}{\RIGHTlineTHICKsmall{8}}

\multiput(54,81)(00,0){1}{\tiny{2}}

\multiput(29,71)(52,0){1}{\tiny{2}}\multiput(81,71)(52,0){1}{{\bf{{\tiny{2}}}}}
\multiput(54,71)(00,0){1}{\tiny{4}}

\multiput(46.5,63.5)(00,0){1}{{\bf{{\tiny{3}}}}}
\multiput(61,63.5)(00,0){1}{{\bf{{\tiny{3}}}}}

\multiput(54,61)(00,0){1}{\tiny{4}}

\multiput(41.5,55)(00,0){1}{{\tiny{2}}}
\multiput(66,55)(00,0){1}{{\tiny{2}}}

\multiput(05,41.5)(60,0){2}{{\tiny{4}}}
\multiput(12,41.5)(60,0){2}{{\bf{{\tiny{3}}}}}

\multiput(35,41.5)(60,0){2}{{\bf{{\tiny{3}}}}}
\multiput(42,41.5)(60,0){1}{{\bf{{\tiny{4}}}}}\multiput(102,41.5)(60,0){1}{{\tiny{4}}}

\multiput(24,36)(30,0){3}{{\tiny{2}}}

\multiput(1,33.5)(60,0){2}{{\bf{{\tiny{3}}}}}
\multiput(16.5,33.5)(60,0){1}{{\bf{{\tiny{4}}}}}\multiput(76.5,33.5)(60,0){1}{{\tiny{4}}}
\multiput(31,33.5)(60,0){1}{{\bf{{\tiny{4}}}}}\multiput(91,33.5)(60,0){1}{{\tiny{4}}}
\multiput(46.5,33.5)(60,0){2}{{\bf{{\tiny{3}}}}}

\multiput(24,31)(30,0){3}{{\tiny{2}}}

\multiput(05,26)(60,0){1}{{\bf{{\tiny{4}}}}}\multiput(65,26)(60,0){1}{{\tiny{4}}}
\multiput(42.5,26)(60,0){2}{{\tiny{4}}}

\multiput(11.5,26)(60,0){2}{{\bf{{\tiny{3}}}}}
\multiput(35,26)(60,0){2}{{\bf{{\tiny{3}}}}}

\multiput(41.5,11.5)(00,0){1}{{\tiny{2}}}
\multiput(66,11.5)(00,0){1}{{\bf{{\tiny{2}}}}}

\multiput(54,6)(00,0){1}{\tiny{4}}

\multiput(46.5,3.5)(00,0){1}{{\bf{{\tiny{3}}}}}
\multiput(61,3.5)(00,0){1}{{\bf{{\tiny{3}}}}}

\multiput(29,-4)(52,0){2}{\tiny{2}}
\multiput(54,-4)(00,0){1}{{\bf{{\tiny{4}}}}}

\multiput(54,-14)(00,0){1}{\tiny{2}}

\multiput(49,72)(00,0){1}{{\tiny{${\bar{c}}$}}}
\multiput(59,72)(00,0){1}{{\tiny{${\bar{d}}^\prime$}}}

\multiput(49.5,55)(00,0){1}{{\tiny{${\bar{c}}$}}}
\multiput(57,54)(00,0){1}{{\tiny{${\bar{d}}^\prime$}}}

\multiput(11.5,50)(00,0){1}{{\tiny{${\bar{a}}$}}}
\multiput(35.5,50)(00,0){1}{{\tiny{${\bar{a}}$}}}
\multiput(71.5,50)(00,0){1}{{\tiny{${\bar{a}}^\prime$}}}
\multiput(94,50)(00,0){1}{{\tiny{${\bar{a}}^\prime$}}}

\multiput(19.5,42)(00,0){1}{{\tiny{${\bar{b}}$}}}
\multiput(27,42)(00,0){1}{{\tiny{${\bar{b}}$}}}
\multiput(49.5,42)(00,0){1}{{\tiny{${\bar{d}}$}}}
\multiput(58,42)(00,0){1}{{\tiny{${\bar{c}}^\prime$}}}

\multiput(79.5,42)(00,0){1}{{\tiny{${\bar{b}}^\prime$}}}
\multiput(87,42)(00,0){1}{{\tiny{${\bar{b}}^\prime$}}}

\multiput(-4.5,39)(0,0){1}{{\tiny{${\bar{d}}$}}}
\multiput(112,39)(0,0){1}{{\tiny{${\bar{c}}^\prime$}}}

\multiput(-4.5,29)(0,0){1}{{\tiny{${\bar{d}}$}}}
\multiput(112,29)(0,0){1}{{\tiny{${\bar{c}}^\prime$}}}

\multiput(20,25)(00,0){1}{{\tiny{${\bar{b}}$}}}
\multiput(27,25)(00,0){1}{{\tiny{${\bar{b}}$}}}
\multiput(49.5,25)(00,0){1}{{\tiny{${\bar{d}}$}}}
\multiput(57,25)(00,0){1}{{\tiny{${\bar{c}}^\prime$}}}

\multiput(79.5,25)(00,0){1}{{\tiny{${\bar{b}}^\prime$}}}
\multiput(86,25)(00,0){1}{{\tiny{${\bar{b}}^\prime$}}}

\multiput(11.5,18)(00,0){1}{{\tiny{${\bar{a}}$}}}
\multiput(35.5,18)(00,0){1}{{\tiny{${\bar{a}}$}}}

\multiput(71.5,18)(00,0){1}{{\tiny{${\bar{a}}^\prime$}}}
\multiput(93.5,17)(00,0){1}{{\tiny{${\bar{a}}^\prime$}}}

\multiput(49.5,12)(00,0){1}{{\tiny{${\bar{c}}$}}}
\multiput(58,12)(00,0){1}{{\tiny{${\bar{d}}^\prime$}}}

\multiput(49,-5)(00,0){1}{{\tiny{${\bar{c}}$}}}
\multiput(59,-6)(00,0){1}{{\tiny{${\bar{d}}^\prime$}}}

\end{picture}
}

\newcommand{\TotalCHHecThreeNineForHecFourTwetytwoDash}{
\setlength{\unitlength}{0.5mm}
\begin{picture}(110,100)(00,-15)

\multiput(50,70)(10,0){2}{\circle{2}}
\multiput(51,70)(00,0){1}{\line(1,0){8}}

\multiput(50,60)(10,0){2}{\circle{2}}
\multiput(51,60)(00,0){1}{\line(1,0){8}}
\multiput(50,61)(00,0){1}{\line(0,1){8}}
\multiput(60,61)(00,0){1}{\line(0,1){8}}

\multiput(10,50)(30,0){4}{\circle{2}}
\multiput(40.8,50.8)(00,0){1}{\line(1,1){8.4}}
\multiput(69.2,50.8)(00,0){1}{\line(-1,1){8.4}}

\multiput(00,40)(30,0){4}{\circle{2}}
\multiput(00.8,40.8)(30,0){4}{\line(1,1){8.4}}
\multiput(20,40)(30,0){4}{\circle{2}}
\multiput(19.2,40.8)(30,0){4}{\line(-1,1){8.4}}
\multiput(21,40)(30,0){3}{\line(1,0){8}}

\multiput(00,30)(30,0){4}{\circle{2}}
\multiput(00,31)(30,0){4}{\line(0,1){8}}
\multiput(20,30)(30,0){4}{\circle{2}}
\multiput(20,31)(30,0){4}{\line(0,1){8}}
\multiput(21,30)(30,0){3}{\line(1,0){8}}

\multiput(10,20)(30,0){4}{\circle{2}}
\multiput(9.2,20.8)(30,0){4}{\line(-1,1){8.4}}
\multiput(10.8,20.8)(30,0){4}{\line(1,1){8.4}}

\multiput(50,10)(10,0){2}{\circle{2}}
\multiput(49.2,10.8)(00,0){1}{\line(-1,1){8.4}}
\multiput(51,10)(00,0){1}{\line(1,0){8}}
\multiput(60.8,10.8)(00,0){1}{\line(1,1){8.4}}

\multiput(50,00)(10,0){2}{\circle{2}}
\multiput(51,00)(00,0){1}{\line(1,0){8}}
\multiput(50,01)(00,0){1}{\line(0,1){8}}
\multiput(60,01)(00,0){1}{\line(0,1){8}}

\multiput(00,00)(110,0){2}{\line(0,1){29}}
\put(10,00){\oval(20,20)[bl]}
\put(100,00){\oval(20,20)[br]}
\multiput(10,-10)(00,0){1}{\line(1,0){90}}

\multiput(10,10)(90,0){2}{\line(0,1){9}}
\put(20,10){\oval(20,20)[bl]}
\put(90,10){\oval(20,20)[br]}
\multiput(20,00)(41,0){2}{\line(1,0){29}}

\multiput(00,41)(110,0){2}{\properUPlineTHICKsmall{29}}
\put(10,70){\curveNorthWestTHICKsmall}
\put(100,70){\curveNorthEastTHICKsmall}
\multiput(10,80)(00,0){1}{\properRIGHTlineTHICKsmall{90}}

\put(20,60){\curveNorthWestTHICKsmall}
\multiput(10,51)(90,0){2}{\properUPlineTHICKsmall{9}}
\put(90,60){\curveNorthEastTHICKsmall}
\multiput(20,70)(41,0){2}{\properRIGHTlineTHICKsmall{29}}

\multiput(50,60)(00,0){1}{\RIGHTlineTHICKsmall{8}}
\multiput(50,60)(10,0){2}{\UPlineTHICKsmall{8}}

\multiput(20,40)(30,0){3}{\UPLEFTlineTHICKsmall{8.4}}
\multiput(30,40)(30,0){3}{\UPRIGHTlineTHICKsmall{8.4}}
\multiput(80,40)(00,0){1}{\RIGHTlineTHICKsmall{8}}

\multiput(00,30)(30,0){3}{\UPlineTHICKsmall{8}}
\multiput(20,30)(30,0){2}{\UPlineTHICKsmall{8}}
\multiput(50,30)(30,0){2}{\RIGHTlineTHICKsmall{8}}
\multiput(110,30)(00,0){1}{\UPlineTHICKsmall{8}}

\multiput(10,20)(30,0){2}{\UPLEFTlineTHICKsmall{8.4}}
\multiput(10,20)(60,0){2}{\UPRIGHTlineTHICKsmall{8.4}}
\multiput(100,20)(00,0){1}{\UPLEFTlineTHICKsmall{8.4}}
\multiput(100,20)(00,0){1}{\UPRIGHTlineTHICKsmall{8.4}}

\multiput(50,10)(00,0){1}{\UPLEFTlineTHICKsmall{8.4}}
\multiput(60,10)(00,0){1}{\UPRIGHTlineTHICKsmall{8.4}}

\multiput(50,00)(10,0){2}{\UPlineTHICKsmall{8}}
\multiput(50,00)(00,0){1}{\RIGHTlineTHICKsmall{8}}

\multiput(54,81)(00,0){1}{\tiny{2}}

\multiput(29,71)(52,0){1}{{\bf{{\tiny{2}}}}}\multiput(81,71)(52,0){1}{\tiny{2}}
\multiput(54,71)(00,0){1}{\tiny{4}}

\multiput(46.5,63.5)(00,0){1}{{\bf{{\tiny{3}}}}}
\multiput(61,63.5)(00,0){1}{{\bf{{\tiny{3}}}}}

\multiput(54,61)(00,0){1}{\tiny{4}}

\multiput(41.5,55)(00,0){1}{{\tiny{2}}}
\multiput(66,55)(00,0){1}{{\tiny{2}}}

\multiput(05,41.5)(60,0){1}{{\tiny{4}}}\multiput(65,41.5)(60,0){1}{{\bf{{\tiny{4}}}}}
\multiput(12,41.5)(60,0){2}{{\bf{{\tiny{3}}}}}

\multiput(35,41.5)(60,0){2}{{\bf{{\tiny{3}}}}}
\multiput(42,41.5)(60,0){2}{{\tiny{4}}}

\multiput(24,36)(30,0){3}{{\tiny{2}}}

\multiput(1,33.5)(60,0){2}{{\bf{{\tiny{3}}}}}
\multiput(16.5,33.5)(60,0){2}{{\tiny{4}}}
\multiput(31,33.5)(60,0){2}{{\tiny{4}}}
\multiput(46.5,33.5)(60,0){2}{{\bf{{\tiny{3}}}}}

\multiput(24,31)(30,0){3}{{\tiny{2}}}

\multiput(05,26)(60,0){2}{{\tiny{4}}}
\multiput(42.5,26)(60,0){1}{{\tiny{4}}}\multiput(102.5,26)(60,0){1}{{\bf{{\tiny{4}}}}}

\multiput(11.5,26)(60,0){2}{{\bf{{\tiny{3}}}}}
\multiput(35,26)(60,0){2}{{\bf{{\tiny{3}}}}}

\multiput(41.5,11.5)(00,0){1}{{\bf{{\tiny{2}}}}}
\multiput(66,11.5)(00,0){1}{{\tiny{2}}}

\multiput(54,6)(00,0){1}{\tiny{4}}

\multiput(46.5,3.5)(00,0){1}{{\bf{{\tiny{3}}}}}
\multiput(61,3.5)(00,0){1}{{\bf{{\tiny{3}}}}}

\multiput(29,-4)(52,0){2}{\tiny{2}}
\multiput(54,-4)(00,0){1}{\tiny{4}}

\multiput(54,-14)(00,0){1}{\tiny{2}}

\multiput(49,72)(00,0){1}{{\tiny{${\bar{d}}^\prime$}}}
\multiput(59,72)(00,0){1}{{\tiny{${\bar{c}}$}}}

\multiput(57,55)(00,0){1}{{\tiny{${\bar{c}}$}}}
\multiput(49.5,54)(00,0){1}{{\tiny{${\bar{d}}^\prime$}}}

\multiput(11.5,50)(00,0){1}{{\tiny{${\bar{a}}^\prime$}}}
\multiput(35.5,50)(00,0){1}{{\tiny{${\bar{a}}^\prime$}}}
\multiput(71.5,50)(00,0){1}{{\tiny{${\bar{a}}$}}}
\multiput(94,50)(00,0){1}{{\tiny{${\bar{a}}$}}}

\multiput(19.5,42)(00,0){1}{{\tiny{${\bar{b}}^\prime$}}}
\multiput(27,42)(00,0){1}{{\tiny{${\bar{b}}^\prime$}}}
\multiput(49.5,42)(00,0){1}{{\tiny{${\bar{c}}^\prime$}}}
\multiput(58,42)(00,0){1}{{\tiny{${\bar{d}}$}}}

\multiput(79.5,42)(00,0){1}{{\tiny{${\bar{b}}$}}}
\multiput(87,42)(00,0){1}{{\tiny{${\bar{b}}$}}}

\multiput(-4.5,39)(0,0){1}{{\tiny{${\bar{c}}^\prime$}}}
\multiput(112,39)(0,0){1}{{\tiny{${\bar{d}}$}}}

\multiput(-4.5,29)(0,0){1}{{\tiny{${\bar{c}}^\prime$}}}
\multiput(112,29)(0,0){1}{{\tiny{${\bar{d}}$}}}

\multiput(20,24)(00,0){1}{{\tiny{${\bar{b}}^\prime$}}}
\multiput(27,24)(00,0){1}{{\tiny{${\bar{b}}^\prime$}}}
\multiput(49.5,25)(00,0){1}{{\tiny{${\bar{c}}^\prime$}}}
\multiput(57,24)(00,0){1}{{\tiny{${\bar{d}}$}}}

\multiput(79.5,24)(00,0){1}{{\tiny{${\bar{b}}$}}}
\multiput(87,24)(00,0){1}{{\tiny{${\bar{b}}$}}}

\multiput(11.5,18)(00,0){1}{{\tiny{${\bar{a}}^\prime$}}}
\multiput(33.5,17)(00,0){1}{{\tiny{${\bar{a}}^\prime$}}}

\multiput(71.5,18)(00,0){1}{{\tiny{${\bar{a}}$}}}
\multiput(95.5,18)(00,0){1}{{\tiny{${\bar{a}}$}}}

\multiput(49.5,12)(00,0){1}{{\tiny{${\bar{d}}^\prime$}}}
\multiput(58,12)(00,0){1}{{\tiny{${\bar{c}}$}}}

\multiput(49,-6)(00,0){1}{{\tiny{${\bar{d}}^\prime$}}}
\multiput(59,-5)(00,0){1}{{\tiny{${\bar{c}}$}}}

\end{picture}
}

\newcommand{\TotalCHHecThreeFour}{
\setlength{\unitlength}{0.5mm}
\begin{picture}(110,65)(00,-10)

\multiput(20,31)(00,0){1}{\curveNorthWestTHICKsmall}
\multiput(20,41)(00,0){1}{\properRIGHTlineTHICKsmall{70}}
\multiput(90,31)(00,0){1}{\curveNorthEastTHICKsmall}

\multiput(00,11)(110,0){2}{\line(0,1){30}}
\multiput(09,41)(00,0){1}{\oval(18,18)[tl]}
\multiput(09,50)(00,0){1}{\line(1,0){92}}
\multiput(101,41)(00,0){1}{\oval(18,18)[tr]}

\multiput(10,30)(40,0){3}{\circle{2}}
\multiput(11,30)(40,0){3}{\line(1,0){8}}
\multiput(20,30)(40,0){3}{\circle{2}}
\multiput(21,30)(40,0){2}{\line(1,0){28}}

\multiput(10,20)(40,0){3}{\circle{2}}
\multiput(10,21)(40,0){3}{\line(0,1){8}}
\multiput(11,20)(40,0){3}{\line(1,0){8}}
\multiput(20,20)(40,0){3}{\circle{2}}
\multiput(20,21)(40,0){3}{\line(0,1){8}}

\multiput(00,10)(40,0){3}{\circle{2}}
\multiput(00.8,10.8)(40,0){3}{\line(1,1){8.4}}
\multiput(30,10)(40,0){3}{\circle{2}}
\multiput(29.2,10.8)(40,0){3}{\line(-1,1){8.4}}
\multiput(31,10)(40,0){2}{\line(1,0){8}}

\multiput(00,00)(40,0){3}{\circle{2}}
\multiput(00,01)(40,0){3}{\line(0,1){8}}
\multiput(01,00)(40,0){3}{\line(1,0){28}}
\multiput(30,00)(40,0){3}{\circle{2}}
\multiput(30,01)(40,0){3}{\line(0,1){8}}
\multiput(31,00)(40,0){2}{\line(1,0){8}}

\put(10,-1){\oval(20,20)[bl]}
\multiput(10,-11)(00,0){1}{\line(1,0){90}}
\put(100,-1){\oval(20,20)[br]}

\multiput(20,30)(40,0){2}{\RIGHTlineTHICKsmall{28}}
\multiput(50,30)(40,0){2}{\RIGHTlineTHICKsmall{8}}

\multiput(10,20)(10,0){2}{\UPlineTHICKsmall{8}}
\multiput(50,20)(40,0){2}{\RIGHTlineTHICKsmall{8}}

\multiput(00,10)(40,0){3}{\UPRIGHTlineTHICKsmall{8.4}}
\multiput(30,10)(40,0){2}{\RIGHTlineTHICKsmall{8}}
\multiput(30,10)(40,0){3}{\UPLEFTlineTHICKsmall{8.4}}

\multiput(00,0)(110,0){2}{\UPlineTHICKsmall{8}}

\multiput(00,00)(40,0){3}{\RIGHTlineTHICKsmall{28}}
\multiput(30,00)(40,0){2}{\RIGHTlineTHICKsmall{8}}

\multiput(54,52)(00,0){1}{\tiny{3}}

\multiput(54,42)(00,0){1}{\tiny{2}}

\multiput(14,31)(40,0){3}{\tiny{1}}
\multiput(34,31)(40,0){2}{\tiny{2}}

\multiput(14,16)(40,0){3}{\tiny{1}}

\multiput(11,23.5)(40,0){3}{{\tiny{3}}}
\multiput(16.5,23.5)(40,0){3}{{\tiny{3}}}

\multiput(01.5,15)(40,0){3}{{\tiny{2}}}
\multiput(26,15)(40,0){3}{{\tiny{2}}}

\multiput(34,6)(40,0){2}{\tiny{3}}

\multiput(1,3.5)(40,0){3}{{\tiny{1}}}
\multiput(26.5,3.5)(40,0){3}{{\tiny{1}}}

\multiput(14,2)(40,0){3}{\tiny{2}}
\multiput(34,2)(40,0){2}{\tiny{3}}

\multiput(54,-15.5)(00,0){1}{\tiny{3}}

\multiput(19,32)(70,0){2}{{\tiny{$b^{\prime\prime}$}}}
\multiput(49,33)(10,0){2}{{\tiny{$c^{\prime\prime}$}}}

\multiput(3,29)(99,0){2}{{\tiny{$a^{\prime\prime}$}}}

\multiput(3,19.5)(99,0){2}{{\tiny{$a^{\prime\prime}$}}}
\multiput(22,19.5)(62,0){2}{{\tiny{$b^{\prime\prime}$}}}
\multiput(43,19.5)(19,0){2}{{\tiny{$d^{\prime\prime}$}}}

\multiput(30,11.5)(00,0){1}{{\tiny{$c^{\prime\prime}$}}}
\multiput(82,9)(00,0){1}{{\tiny{$c^{\prime\prime}$}}}
\multiput(42,9)(00,0){1}{{\tiny{$d^{\prime\prime}$}}}
\multiput(69,12)(00,0){1}{{\tiny{$d^{\prime\prime}$}}}

\multiput(-7,09)(119,0){2}{{\tiny{$a^{\prime\prime}$}}}
\multiput(-7,-1)(119,0){2}{{\tiny{$b^{\prime\prime}$}}}

\multiput(29,-6)(50,0){2}{{\tiny{$c^{\prime\prime}$}}}
\multiput(39,-6)(30,0){2}{{\tiny{$d^{\prime\prime}$}}}

\put(120,0){$(a^\prime\cong a^{\prime\prime}\cong a)$}

\end{picture}
}

\newcommand{\TotalCHHecThreeFourForthirdpartHecFourEighteen}{
\setlength{\unitlength}{0.5mm}
\begin{picture}(110,60)(00,-15)

\multiput(20,31)(00,0){1}{\curveNorthWestTHICKsmall}
\multiput(20,41)(00,0){1}{\properRIGHTlineTHICKsmall{70}}
\multiput(90,31)(00,0){1}{\curveNorthEastTHICKsmall}

\multiput(00,11)(110,0){2}{\line(0,1){30}}
\multiput(09,41)(00,0){1}{\oval(18,18)[tl]}
\multiput(09,50)(00,0){1}{\line(1,0){92}}
\multiput(101,41)(00,0){1}{\oval(18,18)[tr]}

\multiput(10,30)(40,0){3}{\circle{2}}
\multiput(11,30)(40,0){3}{\line(1,0){8}}
\multiput(20,30)(40,0){3}{\circle{2}}
\multiput(21,30)(40,0){2}{\line(1,0){28}}

\multiput(10,20)(40,0){3}{\circle{2}}
\multiput(10,21)(40,0){3}{\line(0,1){8}}
\multiput(11,20)(40,0){3}{\line(1,0){8}}
\multiput(20,20)(40,0){3}{\circle{2}}
\multiput(20,21)(40,0){3}{\line(0,1){8}}

\multiput(00,10)(40,0){3}{\circle{2}}
\multiput(00.8,10.8)(40,0){3}{\line(1,1){8.4}}
\multiput(30,10)(40,0){3}{\circle{2}}
\multiput(29.2,10.8)(40,0){3}{\line(-1,1){8.4}}
\multiput(31,10)(40,0){2}{\line(1,0){8}}

\multiput(00,00)(40,0){3}{\circle{2}}
\multiput(00,01)(40,0){3}{\line(0,1){8}}
\multiput(01,00)(40,0){3}{\line(1,0){28}}
\multiput(30,00)(40,0){3}{\circle{2}}
\multiput(30,01)(40,0){3}{\line(0,1){8}}
\multiput(31,00)(40,0){2}{\line(1,0){8}}

\put(10,-1){\oval(20,20)[bl]}
\multiput(10,-11)(00,0){1}{\line(1,0){90}}
\put(100,-1){\oval(20,20)[br]}

\multiput(20,30)(40,0){2}{\RIGHTlineTHICKsmall{28}}
\multiput(50,30)(40,0){2}{\RIGHTlineTHICKsmall{8}}

\multiput(10,20)(10,0){2}{\UPlineTHICKsmall{8}}
\multiput(50,20)(40,0){2}{\RIGHTlineTHICKsmall{8}}

\multiput(00,10)(40,0){3}{\UPRIGHTlineTHICKsmall{8.4}}
\multiput(30,10)(40,0){2}{\RIGHTlineTHICKsmall{8}}
\multiput(30,10)(40,0){3}{\UPLEFTlineTHICKsmall{8.4}}

\multiput(00,0)(110,0){2}{\UPlineTHICKsmall{8}}

\multiput(00,00)(40,0){3}{\RIGHTlineTHICKsmall{28}}
\multiput(30,00)(40,0){2}{\RIGHTlineTHICKsmall{8}}

\multiput(54,52)(00,0){1}{\tiny{4}}

\multiput(54,42)(00,0){1}{\tiny{2}}

\multiput(14,31)(40,0){3}{\tiny{3}}
\multiput(34,31)(40,0){2}{\tiny{2}}

\multiput(14,16)(40,0){3}{\tiny{3}}

\multiput(11,23.5)(40,0){3}{{\tiny{4}}}
\multiput(16.5,23.5)(40,0){3}{{\tiny{4}}}

\multiput(01.5,15)(40,0){3}{{\tiny{2}}}
\multiput(26,15)(40,0){3}{{\tiny{2}}}

\multiput(34,6)(40,0){2}{\tiny{4}}

\multiput(1,3.5)(40,0){3}{{\tiny{3}}}
\multiput(26.5,3.5)(40,0){3}{{\tiny{3}}}

\multiput(14,1)(40,0){3}{\tiny{2}}
\multiput(34,1)(40,0){2}{\tiny{4}}

\multiput(54,-15.5)(00,0){1}{\tiny{4}}

\multiput(19,32)(70,0){2}{{\tiny{$a^\prime$}}}
\multiput(49,32)(10,0){2}{{\tiny{$a^\prime$}}}

\multiput(4,29)(98,0){2}{{\tiny{$a^\prime$}}}

\multiput(4,19.5)(40,0){3}{{\tiny{$a^\prime$}}}
\multiput(22,19.5)(40,0){3}{{\tiny{$a^\prime$}}}

\multiput(30,11.5)(00,0){1}{{\tiny{$a^\prime$}}}
\multiput(82,9)(00,0){1}{{\tiny{$a^\prime$}}}
\multiput(42,9)(00,0){1}{{\tiny{$a^\prime$}}}
\multiput(69,12)(00,0){1}{{\tiny{$a^\prime$}}}

\multiput(-6,09)(118,0){2}{{\tiny{$a^\prime$}}}
\multiput(-6,-1)(118,0){2}{{\tiny{$a^\prime$}}}

\multiput(29,-6)(50,0){2}{{\tiny{$a^\prime$}}}
\multiput(39,-6)(30,0){2}{{\tiny{$a^\prime$}}}

\end{picture}
}

\newcommand{\TotalCHHecThreeEight}{
\setlength{\unitlength}{0.5mm}
\begin{picture}(110,60)(00,-10)

\multiput(19,31)(00,0){1}{\oval(18,18)[tl]}
\multiput(19,40)(00,0){1}{\line(1,0){72}}
\multiput(91,31)(00,0){1}{\oval(18,18)[tr]}

\multiput(00,11)(110,0){2}{\line(0,1){30}}
\multiput(09,41)(00,0){1}{\oval(18,18)[tl]}
\multiput(09,50)(00,0){1}{\line(1,0){92}}
\multiput(101,41)(00,0){1}{\oval(18,18)[tr]}

\multiput(10,30)(40,0){3}{\circle{2}}
\multiput(11,30)(40,0){3}{\line(1,0){8}}
\multiput(20,30)(40,0){3}{\circle{2}}
\multiput(21,30)(40,0){2}{\line(1,0){28}}

\multiput(10,20)(40,0){3}{\circle{2}}
\multiput(10,21)(40,0){3}{\line(0,1){8}}
\multiput(11,20)(40,0){3}{\line(1,0){8}}
\multiput(20,20)(40,0){3}{\circle{2}}
\multiput(20,21)(40,0){3}{\line(0,1){8}}

\multiput(00,10)(40,0){3}{\circle{2}}
\multiput(00.8,10.8)(40,0){3}{\line(1,1){8.4}}
\multiput(30,10)(40,0){3}{\circle{2}}
\multiput(29.2,10.8)(40,0){3}{\line(-1,1){8.4}}
\multiput(31,10)(40,0){2}{\line(1,0){8}}

\multiput(00,00)(40,0){3}{\circle{2}}
\multiput(00,01)(40,0){3}{\line(0,1){8}}
\multiput(01,00)(40,0){3}{\line(1,0){28}}
\multiput(30,00)(40,0){3}{\circle{2}}
\multiput(30,01)(40,0){3}{\line(0,1){8}}
\multiput(31,00)(40,0){2}{\line(1,0){8}}

\put(10,-1){\curveSouthWestTHICKsmall}
\multiput(10,-11)(00,0){1}{\properRIGHTlineTHICKsmall{90}}
\put(100,-1){\curveSouthEastTHICKsmall}

\multiput(10,20)(40,0){3}{\UPlineTHICKsmall{8}}
\multiput(20,20)(40,0){3}{\UPlineTHICKsmall{8}}
\multiput(10,30)(40,0){3}{\RIGHTlineTHICKsmall{8}}

\multiput(0,10)(40,0){3}{\UPRIGHTlineTHICKsmall{8.4}}
\multiput(30,10)(40,0){3}{\UPLEFTlineTHICKsmall{8.4}}

\multiput(0,0)(40,0){3}{\UPlineTHICKsmall{8}}
\multiput(30,0)(40,0){3}{\UPlineTHICKsmall{8}}
\multiput(30,0)(40,0){2}{\RIGHTlineTHICKsmall{8}}

\multiput(54,52)(00,0){1}{\tiny{3}}

\multiput(54,42)(00,0){1}{\tiny{2}}

\multiput(14,31)(40,0){3}{\tiny{1}}
\multiput(34,31)(40,0){2}{\tiny{2}}

\multiput(14,21)(40,0){3}{\tiny{1}}

\multiput(6.5,23.5)(40,0){3}{{\tiny{3}}}
\multiput(21,23.5)(40,0){3}{{\tiny{3}}}

\multiput(01.5,15)(40,0){3}{{\tiny{2}}}
\multiput(26,15)(40,0){3}{{\tiny{2}}}

\multiput(34,6)(40,0){2}{\tiny{3}}

\multiput(1,3.5)(40,0){3}{{\tiny{1}}}
\multiput(26.5,3.5)(40,0){3}{{\tiny{1}}}

\multiput(14,-4)(40,0){3}{\tiny{2}}
\multiput(34,-4)(40,0){2}{\tiny{3}}

\multiput(54,-15.5)(00,0){1}{\tiny{3}}

\multiput(19,32)(70,0){2}{{\tiny{$e$}}}
\multiput(49,33)(10,0){2}{{\tiny{$f$}}}

\multiput(5.5,29)(96.5,0){2}{{\tiny{$d$}}}

\multiput(9.5,15)(00,0){1}{{\tiny{$b$}}}
\multiput(18,16)(00,0){1}{{\tiny{$c$}}}
\multiput(49.5,15)(00,0){1}{{\tiny{$f$}}}
\multiput(57,15)(00,0){1}{{\tiny{$f$}}}
\multiput(89.5,16)(00,0){1}{{\tiny{$c$}}}
\multiput(98,15)(00,0){1}{{\tiny{$b$}}}

\multiput(29.5,12)(00,0){1}{{\tiny{$c$}}}
\multiput(38,12)(00,0){1}{{\tiny{$e$}}}
\multiput(69.5,12)(00,0){1}{{\tiny{$e$}}}
\multiput(78,12)(00,0){1}{{\tiny{$c$}}}

\multiput(-4.5,9)(116.5,0){2}{{\tiny{$a$}}}

\multiput(29,-5)(50,0){2}{{\tiny{$b$}}}
\multiput(39,-5)(30,0){2}{{\tiny{$d$}}}

\multiput(-4.5,-1)(116.5,0){2}{{\tiny{$a$}}}

\put(120,0){$(a^\prime\cong a^{\prime\prime}\cong a)$}

\end{picture}
}

\newcommand{\dashTotalCHHecThreeEightForHecFourThirteen}{
\setlength{\unitlength}{0.5mm}
\begin{picture}(110,60)(00,-15)

\multiput(19,31)(00,0){1}{\oval(18,18)[tl]}
\multiput(19,40)(00,0){1}{\line(1,0){72}}
\multiput(91,31)(00,0){1}{\oval(18,18)[tr]}

\multiput(00,11)(110,0){2}{\line(0,1){30}}
\multiput(09,41)(00,0){1}{\oval(18,18)[tl]}
\multiput(09,50)(00,0){1}{\line(1,0){92}}
\multiput(101,41)(00,0){1}{\oval(18,18)[tr]}

\multiput(10,30)(40,0){3}{\circle{2}}
\multiput(11,30)(40,0){3}{\line(1,0){8}}
\multiput(20,30)(40,0){3}{\circle{2}}
\multiput(21,30)(40,0){2}{\line(1,0){28}}

\multiput(10,20)(40,0){3}{\circle{2}}
\multiput(10,21)(40,0){3}{\line(0,1){8}}
\multiput(11,20)(40,0){3}{\line(1,0){8}}
\multiput(20,20)(40,0){3}{\circle{2}}
\multiput(20,21)(40,0){3}{\line(0,1){8}}

\multiput(00,10)(40,0){3}{\circle{2}}
\multiput(00.8,10.8)(40,0){3}{\line(1,1){8.4}}
\multiput(30,10)(40,0){3}{\circle{2}}
\multiput(29.2,10.8)(40,0){3}{\line(-1,1){8.4}}
\multiput(31,10)(40,0){2}{\line(1,0){8}}

\multiput(00,00)(40,0){3}{\circle{2}}
\multiput(00,01)(40,0){3}{\line(0,1){8}}
\multiput(01,00)(40,0){3}{\line(1,0){28}}
\multiput(30,00)(40,0){3}{\circle{2}}
\multiput(30,01)(40,0){3}{\line(0,1){8}}
\multiput(31,00)(40,0){2}{\line(1,0){8}}

\put(10,-1){\curveSouthWestTHICKsmall}
\multiput(10,-11)(00,0){1}{\properRIGHTlineTHICKsmall{90}}
\put(100,-1){\curveSouthEastTHICKsmall}

\multiput(10,20)(40,0){3}{\UPlineTHICKsmall{8}}
\multiput(20,20)(40,0){3}{\UPlineTHICKsmall{8}}
\multiput(10,30)(40,0){3}{\RIGHTlineTHICKsmall{8}}

\multiput(0,10)(40,0){3}{\UPRIGHTlineTHICKsmall{8.4}}
\multiput(30,10)(40,0){3}{\UPLEFTlineTHICKsmall{8.4}}

\multiput(0,0)(40,0){3}{\UPlineTHICKsmall{8}}
\multiput(30,0)(40,0){3}{\UPlineTHICKsmall{8}}
\multiput(30,0)(40,0){2}{\RIGHTlineTHICKsmall{8}}

\multiput(54,51)(00,0){1}{\tiny{3}}

\multiput(54,41)(00,0){1}{\tiny{2}}

\multiput(14,31)(40,0){3}{\tiny{4}}
\multiput(34,31)(40,0){2}{\tiny{2}}

\multiput(14,21)(40,0){3}{\tiny{4}}

\multiput(6.5,23.5)(40,0){3}{{\tiny{3}}}
\multiput(21,23.5)(40,0){3}{{\tiny{3}}}

\multiput(01.5,15)(40,0){3}{{\tiny{2}}}
\multiput(26,15)(40,0){3}{{\tiny{2}}}

\multiput(34,6)(40,0){2}{\tiny{3}}

\multiput(1,3.5)(40,0){3}{{\tiny{4}}}
\multiput(26.5,3.5)(40,0){3}{{\tiny{4}}}

\multiput(14,-4)(40,0){3}{\tiny{2}}
\multiput(34,-4)(40,0){2}{\tiny{3}}

\multiput(54,-15.5)(00,0){1}{\tiny{2}}

\multiput(19,32)(70,0){2}{{\tiny{$e$}}}
\multiput(49,32)(10,0){2}{{\tiny{$e$}}}

\multiput(5.5,29)(96.5,0){2}{{\tiny{$e$}}}

\multiput(9.5,16)(40,0){3}{{\tiny{$e$}}}
\multiput(18,16)(40,0){3}{{\tiny{$e$}}}

\multiput(29.5,12)(00,0){1}{{\tiny{$e$}}}
\multiput(38,12)(00,0){1}{{\tiny{$e$}}}
\multiput(69.5,12)(00,0){1}{{\tiny{$e$}}}
\multiput(78,12)(00,0){1}{{\tiny{$e$}}}

\multiput(-4.5,9)(116.5,0){2}{{\tiny{$e$}}}

\multiput(29,-4)(50,0){2}{{\tiny{$e$}}}
\multiput(39,-4)(30,0){2}{{\tiny{$e$}}}

\multiput(-4.5,-1)(116.5,0){2}{{\tiny{$e$}}}

\end{picture}
}

\newcommand{\TotalCHHecThreeEightForHecFourThirteen}{
\setlength{\unitlength}{0.5mm}
\begin{picture}(110,60)(00,-15)

\multiput(20,31)(00,0){1}{\curveNorthWestTHICKsmall}
\multiput(20,41)(00,0){1}{\properRIGHTlineTHICKsmall{70}}
\multiput(90,31)(00,0){1}{\curveNorthEastTHICKsmall}

\multiput(00,11)(110,0){2}{\line(0,1){30}}
\multiput(09,41)(00,0){1}{\oval(18,18)[tl]}
\multiput(09,50)(00,0){1}{\line(1,0){92}}
\multiput(101,41)(00,0){1}{\oval(18,18)[tr]}

\multiput(10,30)(40,0){3}{\circle{2}}
\multiput(11,30)(40,0){3}{\line(1,0){8}}
\multiput(20,30)(40,0){3}{\circle{2}}
\multiput(21,30)(40,0){2}{\line(1,0){28}}

\multiput(10,20)(40,0){3}{\circle{2}}
\multiput(10,21)(40,0){3}{\line(0,1){8}}
\multiput(11,20)(40,0){3}{\line(1,0){8}}
\multiput(20,20)(40,0){3}{\circle{2}}
\multiput(20,21)(40,0){3}{\line(0,1){8}}

\multiput(00,10)(40,0){3}{\circle{2}}
\multiput(00.8,10.8)(40,0){3}{\line(1,1){8.4}}
\multiput(30,10)(40,0){3}{\circle{2}}
\multiput(29.2,10.8)(40,0){3}{\line(-1,1){8.4}}
\multiput(31,10)(40,0){2}{\line(1,0){8}}

\multiput(00,00)(40,0){3}{\circle{2}}
\multiput(00,01)(40,0){3}{\line(0,1){8}}
\multiput(01,00)(40,0){3}{\line(1,0){28}}
\multiput(30,00)(40,0){3}{\circle{2}}
\multiput(30,01)(40,0){3}{\line(0,1){8}}
\multiput(31,00)(40,0){2}{\line(1,0){8}}

\put(10,-1){\oval(20,20)[bl]}
\multiput(10,-11)(00,0){1}{\line(1,0){90}}
\put(100,-1){\oval(20,20)[br]}

\multiput(20,30)(40,0){2}{\RIGHTlineTHICKsmall{28}}
\multiput(50,30)(40,0){2}{\RIGHTlineTHICKsmall{8}}

\multiput(10,20)(10,0){2}{\UPlineTHICKsmall{8}}
\multiput(50,20)(40,0){2}{\RIGHTlineTHICKsmall{8}}

\multiput(00,10)(40,0){3}{\UPRIGHTlineTHICKsmall{8.4}}
\multiput(30,10)(40,0){2}{\RIGHTlineTHICKsmall{8}}
\multiput(30,10)(40,0){3}{\UPLEFTlineTHICKsmall{8.4}}

\multiput(00,0)(110,0){2}{\UPlineTHICKsmall{8}}

\multiput(00,00)(40,0){3}{\RIGHTlineTHICKsmall{28}}
\multiput(30,00)(40,0){2}{\RIGHTlineTHICKsmall{8}}

\multiput(54,51)(00,0){1}{\tiny{3}}

\multiput(54,42.5)(00,0){1}{\tiny{2}}

\multiput(14,31)(40,0){3}{\tiny{4}}
\multiput(34,31)(40,0){2}{\tiny{2}}

\multiput(14,21)(40,0){3}{\tiny{4}}

\multiput(6.5,23.5)(40,0){3}{{\tiny{3}}}
\multiput(21,23.5)(40,0){3}{{\tiny{3}}}

\multiput(01.5,15)(40,0){3}{{\tiny{2}}}
\multiput(26,15)(40,0){3}{{\tiny{2}}}

\multiput(34,6)(40,0){2}{\tiny{3}}

\multiput(1,3.5)(40,0){3}{{\tiny{4}}}
\multiput(26.5,3.5)(40,0){3}{{\tiny{4}}}

\multiput(14,-4)(40,0){3}{\tiny{2}}
\multiput(34,-4)(40,0){2}{\tiny{3}}

\multiput(54,-15.5)(00,0){1}{\tiny{2}}

\multiput(19,32)(70,0){2}{{\tiny{$e$}}}
\multiput(49,32)(10,0){2}{{\tiny{$e$}}}

\multiput(5.5,29)(96.5,0){2}{{\tiny{$e$}}}

\multiput(9.5,16)(40,0){3}{{\tiny{$e$}}}
\multiput(18,16)(40,0){3}{{\tiny{$e$}}}

\multiput(29.5,12)(00,0){1}{{\tiny{$e$}}}
\multiput(38,12)(00,0){1}{{\tiny{$e$}}}
\multiput(69.5,12)(00,0){1}{{\tiny{$e$}}}
\multiput(78,12)(00,0){1}{{\tiny{$e$}}}

\multiput(-4.5,9)(116.5,0){2}{{\tiny{$e$}}}

\multiput(29,-4)(50,0){2}{{\tiny{$e$}}}
\multiput(39,-4)(30,0){2}{{\tiny{$e$}}}

\multiput(-4.5,-1)(116.5,0){2}{{\tiny{$e$}}}

\end{picture}
}

\newcommand{\TotalCHHecThreeEightForHecFourEighteen}{
\setlength{\unitlength}{0.5mm}
\begin{picture}(110,60)(00,-15)

\multiput(19,31)(00,0){1}{\oval(18,18)[tl]}
\multiput(19,40)(00,0){1}{\line(1,0){72}}
\multiput(91,31)(00,0){1}{\oval(18,18)[tr]}

\multiput(00,11)(110,0){2}{\properUPlineTHICKsmall{30}}
\multiput(10,41)(00,0){1}{\curveNorthWestTHICKsmall}
\multiput(10,51)(00,0){1}{\properRIGHTlineTHICKsmall{90}}
\multiput(100,41)(00,0){1}{\curveNorthEastTHICKsmall}

\multiput(10,30)(40,0){3}{\circle{2}}
\multiput(11,30)(40,0){3}{\line(1,0){8}}
\multiput(20,30)(40,0){3}{\circle{2}}
\multiput(21,30)(40,0){2}{\line(1,0){28}}

\multiput(10,20)(40,0){3}{\circle{2}}
\multiput(10,21)(40,0){3}{\line(0,1){8}}
\multiput(11,20)(40,0){3}{\line(1,0){8}}
\multiput(20,20)(40,0){3}{\circle{2}}
\multiput(20,21)(40,0){3}{\line(0,1){8}}

\multiput(00,10)(40,0){3}{\circle{2}}
\multiput(00.8,10.8)(40,0){3}{\line(1,1){8.4}}
\multiput(30,10)(40,0){3}{\circle{2}}
\multiput(29.2,10.8)(40,0){3}{\line(-1,1){8.4}}
\multiput(31,10)(40,0){2}{\line(1,0){8}}

\multiput(00,00)(40,0){3}{\circle{2}}
\multiput(00,01)(40,0){3}{\line(0,1){8}}
\multiput(01,00)(40,0){3}{\line(1,0){28}}
\multiput(30,00)(40,0){3}{\circle{2}}
\multiput(30,01)(40,0){3}{\line(0,1){8}}
\multiput(31,00)(40,0){2}{\line(1,0){8}}

\put(10,-1){\oval(20,20)[bl]}
\multiput(10,-11)(00,0){1}{\line(1,0){90}}
\put(100,-1){\oval(20,20)[br]}

\multiput(10,30)(80,0){2}{\RIGHTlineTHICKsmall{8}}
\multiput(20,30)(40,0){2}{\RIGHTlineTHICKsmall{28}}

\multiput(10,20)(40,0){2}{\UPlineTHICKsmall{8}}
\multiput(10,20)(40,0){3}{\RIGHTlineTHICKsmall{8}}
\multiput(60,20)(40,0){2}{\UPlineTHICKsmall{8}}

\multiput(30,10)(00,0){1}{\UPLEFTlineTHICKsmall{8.4}}
\multiput(30,10)(40,0){2}{\RIGHTlineTHICKsmall{8}}
\multiput(80,10)(00,0){1}{\UPRIGHTlineTHICKsmall{8.4}}

\multiput(00,0)(70,0){2}{\UPlineTHICKsmall{8}}
\multiput(00,0)(80,0){2}{\RIGHTlineTHICKsmall{28}}
\multiput(40,0)(70,0){2}{\UPlineTHICKsmall{8}}
\multiput(30,0)(40,0){2}{\RIGHTlineTHICKsmall{8}}

\multiput(54,53)(00,0){1}{\tiny{{\bf 4}}}

\multiput(54,41)(00,0){1}{\tiny{2}}

\multiput(14,31)(00,0){1}{\tiny{{\bf 3}}}\multiput(54,31)(00,0){1}{\tiny{3}}\multiput(94,31)(00,0){1}{\tiny{{\bf 3}}}
\multiput(34,31)(40,0){2}{\tiny{2}}

\multiput(14,21)(40,0){3}{\tiny{{\bf 3}}}

\multiput(6.5,23.5)(40,0){2}{{\tiny{{\bf 4}}}}\multiput(86.5,23.5)(00,0){1}{{\tiny{4}}}
\multiput(21,23.5)(00,0){1}{{\tiny{4}}}\multiput(61,23.5)(40,0){2}{{\tiny{{\bf 4}}}}

\multiput(01.5,15)(40,0){3}{{\tiny{2}}}
\multiput(26,15)(40,0){3}{{\tiny{2}}}

\multiput(34,6)(40,0){2}{\tiny{{\bf 4}}}

\multiput(1,3.5)(40,0){3}{{\tiny{3}}}
\multiput(26.5,3.5)(40,0){3}{{\tiny{3}}}

\multiput(14,-4)(00,0){1}{\tiny{{\bf 2}}}\multiput(54,-4)(00,0){1}{\tiny{2}}\multiput(94,-4)(00,0){1}{\tiny{{\bf 2}}}
\multiput(34,-4)(40,0){2}{\tiny{{\bf 4}}}

\multiput(54,-15.5)(00,0){1}{\tiny{2}}

\multiput(19,32)(70,0){2}{{\tiny{$e$}}}
\multiput(49,32)(10,0){2}{{\tiny{$e$}}}

\multiput(5.5,29)(96.5,0){2}{{\tiny{$e$}}}

\multiput(9,14.5)(40,0){3}{{\tiny{$d^\prime$}}}
\multiput(16,14.5)(40,0){3}{{\tiny{$d^\prime$}}}

\multiput(29.5,12)(40,0){2}{{\tiny{$c^\prime$}}}
\multiput(37.5,12)(40,0){2}{{\tiny{$c^\prime$}}}

\multiput(-5.5,9)(117.5,0){2}{{\tiny{$c^\prime$}}}

\multiput(29,-5.5)(50,0){2}{{\tiny{$b^\prime$}}}
\multiput(39,-5.5)(30,0){2}{{\tiny{$b^\prime$}}}

\multiput(-5.5,-1)(117.5,0){2}{{\tiny{$b^\prime$}}}

\end{picture}
}

\newcommand{\TotalCHHecThreeOne}{
\setlength{\unitlength}{0.5mm}
\begin{picture}(110,40)(00,-15)

\multiput(19,31)(00,0){1}{\oval(18,18)[tl]}
\multiput(19,40)(00,0){1}{\line(1,0){72}}
\multiput(91,31)(00,0){1}{\oval(18,18)[tr]}

\multiput(00,11)(110,0){2}{\line(0,1){30}}
\multiput(09,41)(00,0){1}{\oval(18,18)[tl]}
\multiput(09,50)(00,0){1}{\line(1,0){92}}
\multiput(101,41)(00,0){1}{\oval(18,18)[tr]}

\multiput(10,30)(40,0){3}{\circle{2}}
\multiput(11,30)(40,0){3}{\line(1,0){8}}
\multiput(20,30)(40,0){3}{\circle{2}}
\multiput(21,30)(40,0){2}{\line(1,0){28}}

\multiput(10,20)(40,0){3}{\circle{2}}
\multiput(10,21)(40,0){3}{\line(0,1){8}}
\multiput(11,20)(40,0){3}{\line(1,0){8}}
\multiput(20,20)(40,0){3}{\circle{2}}
\multiput(20,21)(40,0){3}{\line(0,1){8}}

\multiput(00,10)(40,0){3}{\circle{2}}
\multiput(00.8,10.8)(40,0){3}{\line(1,1){8.4}}
\multiput(30,10)(40,0){3}{\circle{2}}
\multiput(29.2,10.8)(40,0){3}{\line(-1,1){8.4}}
\multiput(31,10)(40,0){2}{\line(1,0){8}}

\multiput(00,00)(40,0){3}{\circle{2}}
\multiput(00,01)(40,0){3}{\line(0,1){8}}
\multiput(01,00)(40,0){3}{\line(1,0){28}}
\multiput(30,00)(40,0){3}{\circle{2}}
\multiput(30,01)(40,0){3}{\line(0,1){8}}
\multiput(31,00)(40,0){2}{\line(1,0){8}}

\put(10,-1){\curveSouthWestTHICKsmall}
\multiput(10,-11)(00,0){1}{\properRIGHTlineTHICKsmall{90}}
\put(100,-1){\curveSouthEastTHICKsmall}

\multiput(10,20)(40,0){3}{\UPlineTHICKsmall{8}}
\multiput(20,20)(40,0){3}{\UPlineTHICKsmall{8}}
\multiput(10,30)(40,0){3}{\RIGHTlineTHICKsmall{8}}

\multiput(0,10)(40,0){3}{\UPRIGHTlineTHICKsmall{8.4}}
\multiput(30,10)(40,0){3}{\UPLEFTlineTHICKsmall{8.4}}

\multiput(0,0)(40,0){3}{\UPlineTHICKsmall{8}}
\multiput(30,0)(40,0){3}{\UPlineTHICKsmall{8}}
\multiput(30,0)(40,0){2}{\RIGHTlineTHICKsmall{8}}

\multiput(54,51)(00,0){1}{\tiny{3}}

\multiput(54,41)(00,0){1}{\tiny{2}}

\multiput(14,31)(40,0){3}{\tiny{1}}
\multiput(34,31)(40,0){2}{\tiny{2}}

\multiput(14,21)(40,0){3}{\tiny{1}}

\multiput(6.5,23.5)(40,0){3}{{\tiny{3}}}
\multiput(21,23.5)(40,0){3}{{\tiny{3}}}

\multiput(01.5,15)(40,0){3}{{\tiny{2}}}
\multiput(26,15)(40,0){3}{{\tiny{2}}}

\multiput(34,6)(40,0){2}{\tiny{3}}

\multiput(1,3.5)(40,0){3}{{\tiny{1}}}
\multiput(26.5,3.5)(40,0){3}{{\tiny{1}}}

\multiput(14,-4)(40,0){3}{\tiny{2}}
\multiput(34,-4)(40,0){2}{\tiny{3}}

\multiput(54,-15.5)(00,0){1}{\tiny{3}}

\end{picture}
}

\newcommand{\TotalCHHecThreeOneDASH}{
\setlength{\unitlength}{0.5mm}
\begin{picture}(110,60)(00,-15)

\multiput(20,31)(00,0){1}{\curveNorthWestTHICKsmall}
\multiput(20,41)(00,0){1}{\properRIGHTlineTHICKsmall{70}}
\multiput(90,31)(00,0){1}{\curveNorthEastTHICKsmall}

\multiput(00,11)(110,0){2}{\line(0,1){30}}
\multiput(09,41)(00,0){1}{\oval(18,18)[tl]}
\multiput(09,50)(00,0){1}{\line(1,0){92}}
\multiput(101,41)(00,0){1}{\oval(18,18)[tr]}

\multiput(10,30)(40,0){3}{\circle{2}}
\multiput(11,30)(40,0){3}{\line(1,0){8}}
\multiput(20,30)(40,0){3}{\circle{2}}
\multiput(21,30)(40,0){2}{\line(1,0){28}}

\multiput(10,20)(40,0){3}{\circle{2}}
\multiput(10,21)(40,0){3}{\line(0,1){8}}
\multiput(11,20)(40,0){3}{\line(1,0){8}}
\multiput(20,20)(40,0){3}{\circle{2}}
\multiput(20,21)(40,0){3}{\line(0,1){8}}

\multiput(00,10)(40,0){3}{\circle{2}}
\multiput(00.8,10.8)(40,0){3}{\line(1,1){8.4}}
\multiput(30,10)(40,0){3}{\circle{2}}
\multiput(29.2,10.8)(40,0){3}{\line(-1,1){8.4}}
\multiput(31,10)(40,0){2}{\line(1,0){8}}

\multiput(00,00)(40,0){3}{\circle{2}}
\multiput(00,01)(40,0){3}{\line(0,1){8}}
\multiput(01,00)(40,0){3}{\line(1,0){28}}
\multiput(30,00)(40,0){3}{\circle{2}}
\multiput(30,01)(40,0){3}{\line(0,1){8}}
\multiput(31,00)(40,0){2}{\line(1,0){8}}

\put(10,-1){\oval(20,20)[bl]}
\multiput(10,-11)(00,0){1}{\line(1,0){90}}
\put(100,-1){\oval(20,20)[br]}

\multiput(20,30)(40,0){2}{\RIGHTlineTHICKsmall{28}}
\multiput(50,30)(40,0){2}{\RIGHTlineTHICKsmall{8}}

\multiput(10,20)(10,0){2}{\UPlineTHICKsmall{8}}
\multiput(50,20)(40,0){2}{\RIGHTlineTHICKsmall{8}}

\multiput(00,10)(40,0){3}{\UPRIGHTlineTHICKsmall{8.4}}
\multiput(30,10)(40,0){2}{\RIGHTlineTHICKsmall{8}}
\multiput(30,10)(40,0){3}{\UPLEFTlineTHICKsmall{8.4}}

\multiput(00,0)(110,0){2}{\UPlineTHICKsmall{8}}

\multiput(00,00)(40,0){3}{\RIGHTlineTHICKsmall{28}}
\multiput(30,00)(40,0){2}{\RIGHTlineTHICKsmall{8}}

\multiput(54,51)(00,0){1}{\tiny{3}}

\multiput(54,42)(00,0){1}{\tiny{2}}

\multiput(14,31)(40,0){3}{\tiny{1}}
\multiput(34,31)(40,0){2}{\tiny{2}}

\multiput(14,21)(40,0){3}{\tiny{1}}

\multiput(6.5,23.5)(40,0){3}{{\tiny{3}}}
\multiput(21,23.5)(40,0){3}{{\tiny{3}}}

\multiput(01.5,15)(40,0){3}{{\tiny{2}}}
\multiput(26,15)(40,0){3}{{\tiny{2}}}

\multiput(34,6)(40,0){2}{\tiny{3}}

\multiput(1,3.5)(40,0){3}{{\tiny{1}}}
\multiput(26.5,3.5)(40,0){3}{{\tiny{1}}}

\multiput(14,-4)(40,0){3}{\tiny{2}}
\multiput(34,-4)(40,0){2}{\tiny{3}}

\multiput(54,-15.5)(00,0){1}{\tiny{3}}

\end{picture}
}

\newcommand{\TotalCHAtwoAone}{
\setlength{\unitlength}{0.5mm}
\begin{picture}(50,40)(00,00)

\multiput(31,40)(00,0){1}{\line(1,0){9}}
\multiput(40,30)(00,0){1}{\oval(20,20)[tr]}
\multiput(50,21)(00,0){1}{\line(0,1){9}}

\multiput(20,40)(10,0){2}{\circle{2}}
\multiput(21,40)(00,0){1}{\line(1,0){8}}

\multiput(20,30)(10,0){2}{\circle{2}}
\multiput(20,31)(10,0){2}{\line(0,1){8}}
\multiput(21,30)(00,0){1}{\line(1,0){8}}

\multiput(00,20)(40,0){2}{\circle{2}}
\multiput(01,20)(40,0){2}{\line(1,0){8}}
\multiput(10,20)(40,0){2}{\circle{2}}
\multiput(10.8,20.8)(00,0){1}{\line(1,1){8.4}}
\multiput(39.2,20.8)(00,0){1}{\line(-1,1){8.4}}

\multiput(20,10)(10,0){2}{\circle{2}}
\multiput(19.2,10.8)(00,0){1}{\line(-1,1){8.4}}
\multiput(21,10)(00,0){1}{\line(1,0){8}}
\multiput(30.8,10.8)(00,0){1}{\line(1,1){8.4}}

\multiput(20,00)(10,0){2}{\circle{2}}
\multiput(20,01)(10,0){2}{\line(0,1){8}}
\multiput(21,00)(00,0){1}{\line(1,0){8}}

\multiput(20,40)(00,0){1}{\RIGHTlineTHICKsmall{8}}

\multiput(20,30)(00,0){1}{\RIGHTlineTHICKsmall{8}}
\multiput(30,30)(00,0){1}{\UPlineTHICKsmall{8}}

\multiput(10,20)(00,0){1}{\UPRIGHTlineTHICKsmall{8.4}}
\multiput(40,20)(00,0){1}{\RIGHTlineTHICKsmall{8}}

\multiput(20,10)(00,0){1}{\UPLEFTlineTHICKsmall{8.4}}
\multiput(20,10)(00,0){1}{\RIGHTlineTHICKsmall{8}}
\multiput(30,10)(00,0){1}{\UPRIGHTlineTHICKsmall{8.4}}

\multiput(20,00)(00,0){1}{\RIGHTlineTHICKsmall{8}}

\multiput(10,40)(00,0){1}{\properRIGHTlineTHICKsmall{9}}

\put(10,30){\curveNorthWestTHICKsmall}
\multiput(20,30)(00,0){1}{\RIGHTlineTHICKsmall{8}}

\multiput(00,21)(00,0){1}{\properUPlineTHICKsmall{9}}

\multiput(00,10)(50,0){2}{\properUPlineTHICKsmall{9}}

\put(10,10){\curveSouthWestTHICKsmall}
\multiput(10,00)(21,0){2}{\properRIGHTlineTHICKsmall{9}}
\put(40,10){\curveSouthEastTHICKsmall}

\multiput(9,41.5)(30,0){2}{{\tiny{1}}}
\multiput(24,41.5)(00,0){1}{{\tiny{2}}}

\multiput(16.5,33.5)(00,0){1}{{\tiny{3}}}
\multiput(31,33.5)(00,0){1}{{\tiny{3}}}

\multiput(24,31)(00,0){1}{{\tiny{2}}}

\multiput(16,22)(00,0){1}{{\tiny{1}}}
\multiput(31,22)(00,0){1}{{\tiny{1}}}

\multiput(4,21)(40,0){2}{{\tiny{3}}}

\multiput(16,15)(00,0){1}{{\tiny{2}}}
\multiput(31,15)(00,0){1}{{\tiny{2}}}

\multiput(24,5.5)(00,0){1}{{\tiny{1}}}

\multiput(16.5,3.5)(00,0){1}{{\tiny{3}}}
\multiput(31,3.5)(00,0){1}{{\tiny{3}}}

\multiput(9,-4.5)(30,0){2}{{\tiny{2}}}
\multiput(24,-4.5)(00,0){1}{{\tiny{1}}}

\end{picture}
}

\newcommand{\whitecircleempty}{
\setlength{\unitlength}{1mm}
\begin{picture}(10,10)(2.74,0)
\put(0,0){\circle{2}}
\end{picture}
}

\newcommand{\whitecircleEasta}{
\setlength{\unitlength}{1mm}
\begin{picture}(10,10)(2.74,0)
\put(0,0){\circle{2}}
\put(1.5,-0.7){$a$}
\end{picture}
}

\newcommand{\whitecircleWesta}{
\setlength{\unitlength}{1mm}
\begin{picture}(10,10)(2.74,0)
\put(0,0){\circle{2}}
\put(-3.8,-0.7){$a$}
\end{picture}
}

\newcommand{\whitecircleNortha}{
\setlength{\unitlength}{1mm}
\begin{picture}(10,10)(2.74,0)
\put(0,0){\circle{2}}
\put(-0.9,1.5){$a$}
\end{picture}
}

\newcommand{\whitecircleSoutha}{
\setlength{\unitlength}{1mm}
\begin{picture}(10,10)(2.74,0)
\put(0,0){\circle{2}}
\put(-0.9,-3.4){$a$}
\end{picture}
}

\newcommand{\whitecircleNorthWesta}{
\setlength{\unitlength}{1mm}
\begin{picture}(10,10)(2.74,0)
\put(0,0){\circle{2}}
\put(-3.8,1.5){$a$}
\end{picture}
}

\newcommand{\whitecircleNorthEasta}{
\setlength{\unitlength}{1mm}
\begin{picture}(10,10)(2.74,0)
\put(0,0){\circle{2}}
\put(1.5,1.5){$a$}
\end{picture}
}

\newcommand{\whitecircleSouthWesta}{
\setlength{\unitlength}{1mm}
\begin{picture}(10,10)(2.74,0)
\put(0,0){\circle{2}}
\put(-3.8,-3.4){$a$}
\end{picture}
}

\newcommand{\whitecircleSouthEasta}{
\setlength{\unitlength}{1mm}
\begin{picture}(10,10)(2.74,0)
\put(0,0){\circle{2}}
\put(1.5,-3.4){$a$}
\end{picture}
}

\newcommand{\whitecircleEastb}{
\setlength{\unitlength}{1mm}
\begin{picture}(10,10)(2.74,0)
\put(0,0){\circle{2}}
\put(1.5,-0.7){$b$}
\end{picture}
}

\newcommand{\whitecircleWestb}{
\setlength{\unitlength}{1mm}
\begin{picture}(10,10)(2.74,0)
\put(0,0){\circle{2}}
\put(-3.8,-0.7){$b$}
\end{picture}
}

\newcommand{\whitecircleNorthb}{
\setlength{\unitlength}{1mm}
\begin{picture}(10,10)(2.74,0)
\put(0,0){\circle{2}}
\put(-0.9,1.5){$b$}
\end{picture}
}

\newcommand{\whitecircleSouthb}{
\setlength{\unitlength}{1mm}
\begin{picture}(10,10)(2.74,0)
\put(0,0){\circle{2}}
\put(-0.9,-4.8){$b$}
\end{picture}
}

\newcommand{\whitecircleNorthWestb}{
\setlength{\unitlength}{1mm}
\begin{picture}(10,10)(2.74,0)
\put(0,0){\circle{2}}
\put(-3.8,1.5){$b$}
\end{picture}
}

\newcommand{\whitecircleNorthEastb}{
\setlength{\unitlength}{1mm}
\begin{picture}(10,10)(2.74,0)
\put(0,0){\circle{2}}
\put(1.5,1.5){$b$}
\end{picture}
}

\newcommand{\whitecircleSouthWestb}{
\setlength{\unitlength}{1mm}
\begin{picture}(10,10)(2.74,0)
\put(0,0){\circle{2}}
\put(-3.8,-3.4){$b$}
\end{picture}
}

\newcommand{\whitecircleSouthEastb}{
\setlength{\unitlength}{1mm}
\begin{picture}(10,10)(2.74,0)
\put(0,0){\circle{2}}
\put(1.5,-3.4){$b$}
\end{picture}
}

\newcommand{\whitecircleEastc}{
\setlength{\unitlength}{1mm}
\begin{picture}(10,10)(2.74,0)
\put(0,0){\circle{2}}
\put(1.5,-0.7){$c$}
\end{picture}
}

\newcommand{\whitecircleWestc}{
\setlength{\unitlength}{1mm}
\begin{picture}(10,10)(2.74,0)
\put(0,0){\circle{2}}
\put(-3.8,-0.7){$c$}
\end{picture}
}

\newcommand{\whitecircleNorthc}{
\setlength{\unitlength}{1mm}
\begin{picture}(10,10)(2.74,0)
\put(0,0){\circle{2}}
\put(-0.9,1.5){$c$}
\end{picture}
}

\newcommand{\whitecircleSouthc}{
\setlength{\unitlength}{1mm}
\begin{picture}(10,10)(2.74,0)
\put(0,0){\circle{2}}
\put(-0.9,-3.4){$c$}
\end{picture}
}

\newcommand{\whitecircleNorthWestc}{
\setlength{\unitlength}{1mm}
\begin{picture}(10,10)(2.74,0)
\put(0,0){\circle{2}}
\put(-3.8,1.5){$c$}
\end{picture}
}

\newcommand{\whitecircleNorthEastc}{
\setlength{\unitlength}{1mm}
\begin{picture}(10,10)(2.74,0)
\put(0,0){\circle{2}}
\put(1.5,1.5){$c$}
\end{picture}
}

\newcommand{\whitecircleSouthWestc}{
\setlength{\unitlength}{1mm}
\begin{picture}(10,10)(2.74,0)
\put(0,0){\circle{2}}
\put(-3.8,-3.4){$c$}
\end{picture}
}

\newcommand{\whitecircleSouthEastc}{
\setlength{\unitlength}{1mm}
\begin{picture}(10,10)(2.74,0)
\put(0,0){\circle{2}}
\put(1.5,-3.4){$c$}
\end{picture}
}

\newcommand{\whitecircleEastd}{
\setlength{\unitlength}{1mm}
\begin{picture}(10,10)(2.74,0)
\put(0,0){\circle{2}}
\put(1.5,-0.7){$d$}
\end{picture}
}

\newcommand{\whitecircleWestd}{
\setlength{\unitlength}{1mm}
\begin{picture}(10,10)(2.74,0)
\put(0,0){\circle{2}}
\put(-3.8,-0.7){$d$}
\end{picture}
}

\newcommand{\whitecircleNorthd}{
\setlength{\unitlength}{1mm}
\begin{picture}(10,10)(2.74,0)
\put(0,0){\circle{2}}
\put(-0.9,1.5){$d$}
\end{picture}
}

\newcommand{\whitecircleSouthd}{
\setlength{\unitlength}{1mm}
\begin{picture}(10,10)(2.74,0)
\put(0,0){\circle{2}}
\put(-0.9,-4.8){$d$}
\end{picture}
}

\newcommand{\whitecircleNorthWestd}{
\setlength{\unitlength}{1mm}
\begin{picture}(10,10)(2.74,0)
\put(0,0){\circle{2}}
\put(-3.8,1.5){$d$}
\end{picture}
}

\newcommand{\whitecircleNorthEastd}{
\setlength{\unitlength}{1mm}
\begin{picture}(10,10)(2.74,0)
\put(0,0){\circle{2}}
\put(1.5,1.5){$d$}
\end{picture}
}

\newcommand{\whitecircleSouthWestd}{
\setlength{\unitlength}{1mm}
\begin{picture}(10,10)(2.74,0)
\put(0,0){\circle{2}}
\put(-3.8,-3.4){$d$}
\end{picture}
}

\newcommand{\whitecircleSouthEastd}{
\setlength{\unitlength}{1mm}
\begin{picture}(10,10)(2.74,0)
\put(0,0){\circle{2}}
\put(1.5,-3.4){$d$}
\end{picture}
}

\newcommand{\whitecircleEaste}{
\setlength{\unitlength}{1mm}
\begin{picture}(10,10)(2.74,0)
\put(0,0){\circle{2}}
\put(1.5,-0.7){$e$}
\end{picture}
}

\newcommand{\whitecircleWeste}{
\setlength{\unitlength}{1mm}
\begin{picture}(10,10)(2.74,0)
\put(0,0){\circle{2}}
\put(-3.8,-0.7){$e$}
\end{picture}
}

\newcommand{\whitecircleNorthe}{
\setlength{\unitlength}{1mm}
\begin{picture}(10,10)(2.74,0)
\put(0,0){\circle{2}}
\put(-0.9,1.5){$e$}
\end{picture}
}

\newcommand{\whitecircleSouthe}{
\setlength{\unitlength}{1mm}
\begin{picture}(10,10)(2.74,0)
\put(0,0){\circle{2}}
\put(-0.9,-3.4){$e$}
\end{picture}
}

\newcommand{\whitecircleNorthWeste}{
\setlength{\unitlength}{1mm}
\begin{picture}(10,10)(2.74,0)
\put(0,0){\circle{2}}
\put(-3.8,1.5){$e$}
\end{picture}
}

\newcommand{\whitecircleNorthEaste}{
\setlength{\unitlength}{1mm}
\begin{picture}(10,10)(2.74,0)
\put(0,0){\circle{2}}
\put(1.5,1.5){$e$}
\end{picture}
}

\newcommand{\whitecircleSouthWeste}{
\setlength{\unitlength}{1mm}
\begin{picture}(10,10)(2.74,0)
\put(0,0){\circle{2}}
\put(-3.8,-3.4){$e$}
\end{picture}
}

\newcommand{\whitecircleSouthEaste}{
\setlength{\unitlength}{1mm}
\begin{picture}(10,10)(2.74,0)
\put(0,0){\circle{2}}
\put(1.5,-3.4){$e$}
\end{picture}
}

\newcommand{\whitecircleEastf}{
\setlength{\unitlength}{1mm}
\begin{picture}(10,10)(2.74,0)
\put(0,0){\circle{2}}
\put(1.5,-0.7){$f$}
\end{picture}
}

\newcommand{\whitecircleWestf}{
\setlength{\unitlength}{1mm}
\begin{picture}(10,10)(2.74,0)
\put(0,0){\circle{2}}
\put(-3.8,-0.7){$f$}
\end{picture}
}

\newcommand{\whitecircleNorthf}{
\setlength{\unitlength}{1mm}
\begin{picture}(10,10)(2.74,0)
\put(0,0){\circle{2}}
\put(-0.9,2.5){$f$}
\end{picture}
}

\newcommand{\whitecircleSouthf}{
\setlength{\unitlength}{1mm}
\begin{picture}(10,10)(2.74,0)
\put(0,0){\circle{2}}
\put(-0.9,-4.8){$f$}
\end{picture}
}

\newcommand{\whitecircleNorthWestf}{
\setlength{\unitlength}{1mm}
\begin{picture}(10,10)(2.74,0)
\put(0,0){\circle{2}}
\put(-3.8,1.5){$f$}
\end{picture}
}

\newcommand{\whitecircleNorthEastf}{
\setlength{\unitlength}{1mm}
\begin{picture}(10,10)(2.74,0)
\put(0,0){\circle{2}}
\put(1.5,1.5){$f$}
\end{picture}
}

\newcommand{\whitecircleSouthWestf}{
\setlength{\unitlength}{1mm}
\begin{picture}(10,10)(2.74,0)
\put(0,0){\circle{2}}
\put(-3.8,-3.4){$f$}
\end{picture}
}

\newcommand{\whitecircleSouthEastf}{
\setlength{\unitlength}{1mm}
\begin{picture}(10,10)(2.74,0)
\put(0,0){\circle{2}}
\put(1.5,-3.4){$f$}
\end{picture}
}

\newcommand{\whitecircleEastg}{
\setlength{\unitlength}{1mm}
\begin{picture}(10,10)(2.74,0)
\put(0,0){\circle{2}}
\put(1.5,-0.7){$g$}
\end{picture}
}

\newcommand{\whitecircleWestg}{
\setlength{\unitlength}{1mm}
\begin{picture}(10,10)(2.74,0)
\put(0,0){\circle{2}}
\put(-3.8,-0.7){$g$}
\end{picture}
}

\newcommand{\whitecircleNorthg}{
\setlength{\unitlength}{1mm}
\begin{picture}(10,10)(2.74,0)
\put(0,0){\circle{2}}
\put(-0.5,2.2){$g$}
\end{picture}
}

\newcommand{\whitecircleSouthg}{
\setlength{\unitlength}{1mm}
\begin{picture}(10,10)(2.74,0)
\put(0,0){\circle{2}}
\put(-1.3,-3.7){$g$}
\end{picture}
}

\newcommand{\whitecircleNorthWestg}{
\setlength{\unitlength}{1mm}
\begin{picture}(10,10)(2.74,0)
\put(0,0){\circle{2}}
\put(-3.8,1.5){$g$}
\end{picture}
}

\newcommand{\whitecircleNorthEastg}{
\setlength{\unitlength}{1mm}
\begin{picture}(10,10)(2.74,0)
\put(0,0){\circle{2}}
\put(1.5,1.5){$g$}
\end{picture}
}

\newcommand{\whitecircleSouthWestg}{
\setlength{\unitlength}{1mm}
\begin{picture}(10,10)(2.74,0)
\put(0,0){\circle{2}}
\put(-3.8,-3.4){$g$}
\end{picture}
}

\newcommand{\whitecircleSouthEastg}{
\setlength{\unitlength}{1mm}
\begin{picture}(10,10)(2.74,0)
\put(0,0){\circle{2}}
\put(1.5,-3.4){$g$}
\end{picture}
}

\newcommand{\whitecircleEasth}{
\setlength{\unitlength}{1mm}
\begin{picture}(10,10)(2.74,0)
\put(0,0){\circle{2}}
\put(1.5,-0.7){$h$}
\end{picture}
}

\newcommand{\whitecircleWesth}{
\setlength{\unitlength}{1mm}
\begin{picture}(10,10)(2.74,0)
\put(0,0){\circle{2}}
\put(-3.8,-0.7){$h$}
\end{picture}
}

\newcommand{\whitecircleNorthh}{
\setlength{\unitlength}{1mm}
\begin{picture}(10,10)(2.74,0)
\put(0,0){\circle{2}}
\put(-0.9,1.5){$h$}
\end{picture}
}

\newcommand{\whitecircleSouthh}{
\setlength{\unitlength}{1mm}
\begin{picture}(10,10)(2.74,0)
\put(0,0){\circle{2}}
\put(-0.9,-4.8){$h$}
\end{picture}
}

\newcommand{\whitecircleNorthWesth}{
\setlength{\unitlength}{1mm}
\begin{picture}(10,10)(2.74,0)
\put(0,0){\circle{2}}
\put(-3.8,1.5){$h$}
\end{picture}
}

\newcommand{\whitecircleNorthEasth}{
\setlength{\unitlength}{1mm}
\begin{picture}(10,10)(2.74,0)
\put(0,0){\circle{2}}
\put(1.5,1.5){$h$}
\end{picture}
}

\newcommand{\whitecircleSouthWesth}{
\setlength{\unitlength}{1mm}
\begin{picture}(10,10)(2.74,0)
\put(0,0){\circle{2}}
\put(-3.8,-3.4){$h$}
\end{picture}
}

\newcommand{\whitecircleSouthEasth}{
\setlength{\unitlength}{1mm}
\begin{picture}(10,10)(2.74,0)
\put(0,0){\circle{2}}
\put(1.5,-3.4){$h$}
\end{picture}
}

\newcommand{\whitecircleEasti}{
\setlength{\unitlength}{1mm}
\begin{picture}(10,10)(2.74,0)
\put(0,0){\circle{2}}
\put(1.5,-0.7){$i$}
\end{picture}
}

\newcommand{\whitecircleWesti}{
\setlength{\unitlength}{1mm}
\begin{picture}(10,10)(2.74,0)
\put(0,0){\circle{2}}
\put(-3.8,-0.7){$i$}
\end{picture}
}

\newcommand{\whitecircleNorthi}{
\setlength{\unitlength}{1mm}
\begin{picture}(10,10)(2.74,0)
\put(0,0){\circle{2}}
\put(-0.9,1.5){$i$}
\end{picture}
}

\newcommand{\whitecircleSouthi}{
\setlength{\unitlength}{1mm}
\begin{picture}(10,10)(2.74,0)
\put(0,0){\circle{2}}
\put(-0.9,-4.8){$i$}
\end{picture}
}

\newcommand{\whitecircleNorthWesti}{
\setlength{\unitlength}{1mm}
\begin{picture}(10,10)(2.74,0)
\put(0,0){\circle{2}}
\put(-3.8,1.5){$i$}
\end{picture}
}

\newcommand{\whitecircleNorthEasti}{
\setlength{\unitlength}{1mm}
\begin{picture}(10,10)(2.74,0)
\put(0,0){\circle{2}}
\put(1.5,1.5){$i$}
\end{picture}
}

\newcommand{\whitecircleSouthWesti}{
\setlength{\unitlength}{1mm}
\begin{picture}(10,10)(2.74,0)
\put(0,0){\circle{2}}
\put(-3.8,-3.4){$i$}
\end{picture}
}

\newcommand{\whitecircleSouthEasti}{
\setlength{\unitlength}{1mm}
\begin{picture}(10,10)(2.74,0)
\put(0,0){\circle{2}}
\put(1.5,-3.4){$i$}
\end{picture}
}

\newcommand{\whitecircleEastj}{
\setlength{\unitlength}{1mm}
\begin{picture}(10,10)(2.74,0)
\put(0,0){\circle{2}}
\put(1.5,-0.7){$j$}
\end{picture}
}

\newcommand{\whitecircleWestj}{
\setlength{\unitlength}{1mm}
\begin{picture}(10,10)(2.74,0)
\put(0,0){\circle{2}}
\put(-3.8,-0.7){$j$}
\end{picture}
}

\newcommand{\whitecircleNorthj}{
\setlength{\unitlength}{1mm}
\begin{picture}(10,10)(2.74,0)
\put(0,0){\circle{2}}
\put(-0.9,2.2){$j$}
\end{picture}
}

\newcommand{\whitecircleSouthj}{
\setlength{\unitlength}{1mm}
\begin{picture}(10,10)(2.74,0)
\put(0,0){\circle{2}}
\put(-0.9,-4.8){$j$}
\end{picture}
}

\newcommand{\whitecircleNorthWestj}{
\setlength{\unitlength}{1mm}
\begin{picture}(10,10)(2.74,0)
\put(0,0){\circle{2}}
\put(-3.8,1.5){$j$}
\end{picture}
}

\newcommand{\whitecircleNorthEastj}{
\setlength{\unitlength}{1mm}
\begin{picture}(10,10)(2.74,0)
\put(0,0){\circle{2}}
\put(1.5,1.5){$j$}
\end{picture}
}

\newcommand{\whitecircleSouthWestj}{
\setlength{\unitlength}{1mm}
\begin{picture}(10,10)(2.74,0)
\put(0,0){\circle{2}}
\put(-3.8,-3.4){$j$}
\end{picture}
}

\newcommand{\whitecircleSouthEastj}{
\setlength{\unitlength}{1mm}
\begin{picture}(10,10)(2.74,0)
\put(0,0){\circle{2}}
\put(1.5,-3.4){$j$}
\end{picture}
}

\newcommand{\holizontallinetenmm}{
\setlength{\unitlength}{1mm}
\begin{picture}(10,10)(3.74,0)
\put(1,0){\line(1,0){10}}
\end{picture}
}

\newcommand{\holizontallinetenmmTHICK}{
\setlength{\unitlength}{1mm}
\begin{picture}(10,10)(3.74,0)
\put(1,0){\line(1,0){10}}\multiput(1,0)(0,0.05){4}{\line(1,0){10}}\multiput(1,0)(0,-0.05){4}{\line(1,0){10}}
\end{picture}
}

\newcommand{\holizontallineAboveone}{
\setlength{\unitlength}{1mm}
\begin{picture}(10,10)(2.74,0)
\put(1,0){\line(1,0){8}}
\put(4.5,0.5){{\tiny{1}}}
\end{picture}
}

\newcommand{\holizontallineAboveonebar}{
\setlength{\unitlength}{1mm}
\begin{picture}(10,10)(2.74,0)
\put(1,0){\line(1,0){8}}
\put(4.5,0.5){{\tiny{${\bar{1}}$}}}
\end{picture}
}

\newcommand{\holizontallineAboveoneTHICK}{
\setlength{\unitlength}{1mm}
\begin{picture}(10,10)(2.74,0)
\put(1,0){\line(1,0){8}}\multiput(1,0)(0,0.05){4}{\line(1,0){8}}\multiput(1,0)(0,-0.05){4}{\line(1,0){8}}
\put(4.5,0.5){{\tiny{1}}}
\end{picture}
}

\newcommand{\holizontallineAboveoneTHICKbar}{
\setlength{\unitlength}{1mm}
\begin{picture}(10,10)(2.74,0)
\put(1,0){\line(1,0){8}}\multiput(1,0)(0,0.05){4}{\line(1,0){8}}\multiput(1,0)(0,-0.05){4}{\line(1,0){8}}
\put(4.5,0.5){{\tiny{${\bar{1}}$}}}
\end{picture}
}

\newcommand{\holizontallineAboveoneEast}{
\setlength{\unitlength}{1mm}
\begin{picture}(10,10)(2.74,0)
\put(1,0){\vector(1,0){8}}
\put(4.5,0.5){{\tiny{1}}}
\end{picture}
}

\newcommand{\holizontallineAboveoneWest}{
\setlength{\unitlength}{1mm}
\begin{picture}(10,10)(2.74,0)
\put(9,0){\vector(-1,0){8}}
\put(4.5,0.5){{\tiny{1}}}
\end{picture}
}

\newcommand{\holizontallineUnderone}{
\setlength{\unitlength}{1mm}
\begin{picture}(10,10)(2.74,0)
\put(1,0){\line(1,0){8}}
\put(4.5,-2){{\tiny{1}}}
\end{picture}
}

\newcommand{\verticallineLeftone}{
\setlength{\unitlength}{1mm}
\begin{picture}(10,10)(2.74,0)
\put(0,1){\line(0,1){8}}
\put(-1.7,4.5){{\tiny{1}}}
\end{picture}
}

\newcommand{\verticallineLeftoneUp}{
\setlength{\unitlength}{1mm}
\begin{picture}(10,10)(2.74,0)
\put(0,1){\vector(0,1){8}}
\put(-1.7,4.5){{\tiny{1}}}
\end{picture}
}

\newcommand{\verticallineLeftoneDown}{
\setlength{\unitlength}{1mm}
\begin{picture}(10,10)(2.74,0)
\put(0,9){\line(0,-1){8}}
\put(-1.7,4.5){{\tiny{1}}}
\end{picture}
}

\newcommand{\verticallineRightone}{
\setlength{\unitlength}{1mm}
\begin{picture}(10,10)(2.74,0)
\put(0,1){\line(0,1){8}}
\put(0.5,4.5){{\tiny{1}}}
\end{picture}
}

\newcommand{\verticallineRightonebar}{
\setlength{\unitlength}{1mm}
\begin{picture}(10,10)(2.74,0)
\put(0,1){\line(0,1){8}}
\put(0.5,4.5){{\tiny{${\bar{1}}$}}}
\end{picture}
}

\newcommand{\verticallineRightoneTHICK}{
\setlength{\unitlength}{1mm}
\begin{picture}(10,10)(2.74,0)
\put(0,1){\line(0,1){8}}\multiput(0,1)(0.05,0){4}{\line(0,1){8}}\multiput(0,1)(-0.05,0){4}{\line(0,1){8}}
\put(0.5,4.5){{\tiny{1}}}
\end{picture}
}

\newcommand{\verticallineRightoneTHICKbar}{
\setlength{\unitlength}{1mm}
\begin{picture}(10,10)(2.74,0)
\put(0,1){\line(0,1){8}}\multiput(0,1)(0.05,0){4}{\line(0,1){8}}\multiput(0,1)(-0.05,0){4}{\line(0,1){8}}
\put(0.5,4.5){{\tiny{${\bar{1}}$}}}
\end{picture}
}

\newcommand{\verticallineRightoneUp}{
\setlength{\unitlength}{1mm}
\begin{picture}(10,10)(2.74,0)
\put(0,1){\vector(0,1){8}}
\put(0.5,4.5){{\tiny{1}}}
\end{picture}
}

\newcommand{\verticallineRightoneDown}{
\setlength{\unitlength}{1mm}
\begin{picture}(10,10)(2.74,0)
\put(0,9){\vector(0,-1){8}}
\put(0.5,4.5){{\tiny{1}}}
\end{picture}
}

\newcommand{\EastUplineaboveone}{
\setlength{\unitlength}{1mm}
\begin{picture}(10,10)(2.74,0)
\put(0.8,0.8){\line(1,1){8.4}}
\put(4,5.5){{\tiny{1}}}
\end{picture}
}

\newcommand{\EastUplineaboveoneTHICK}{
\setlength{\unitlength}{1mm}
\begin{picture}(10,10)(2.74,0)
\put(0.8,0.8){\line(1,1){8.4}}\multiput(0.8,0.8)(0.03,-0.03){4}{\line(1,1){8.4}}\multiput(0.8,0.8)(-0.03,0.03){4}{\line(1,1){8.4}}
\put(4,5.5){{\tiny{1}}}
\end{picture}
}

\newcommand{\EastUplineaboveoneUp}{
\setlength{\unitlength}{1mm}
\begin{picture}(10,10)(2.74,0)
\put(0.8,0.8){\vector(1,1){8.4}}
\put(4,5.5){{\tiny{1}}}
\end{picture}
}

\newcommand{\EastUplineaboveoneDown}{
\setlength{\unitlength}{1mm}
\begin{picture}(10,10)(2.74,0)
\put(9.2,9.2){\vector(-1,-1){8.4}}
\put(4,5.5){{\tiny{1}}}
\end{picture}
}

\newcommand{\EastUplineunderone}{
\setlength{\unitlength}{1mm}
\begin{picture}(10,10)(2.74,0)
\put(0.8,0.8){\line(1,1){8.4}}
\put(5,3.5){{\tiny{1}}}
\end{picture}
}

\newcommand{\EastDownlineaboveone}{
\setlength{\unitlength}{1mm}
\begin{picture}(10,10)(2.74,0)
\put(0.8,-0.8){\line(1,-1){8.4}}
\put(5.5,-5){{\tiny{1}}}
\end{picture}
}

\newcommand{\EastDownlineaboveoneTHICK}{
\setlength{\unitlength}{1mm}
\begin{picture}(10,10)(2.74,0)
\put(0.8,-0.8){\line(1,-1){8.4}}\multiput(0.8,-0.8)(0.03,0.03){4}{\line(1,-1){8.4}}\multiput(0.8,-0.8)(-0.03,-0.03){4}{\line(1,-1){8.4}}
\put(5.5,-5){{\tiny{1}}}
\end{picture}
}

\newcommand{\EastDownlineaboveoneUp}{
\setlength{\unitlength}{1mm}
\begin{picture}(10,10)(2.74,0)
\put(9.2,-9.2){\vector(-1,1){8.4}}
\put(5.5,-5){{\tiny{1}}}
\end{picture}
}

\newcommand{\EastDownlineaboveoneDown}{
\setlength{\unitlength}{1mm}
\begin{picture}(10,10)(2.74,0)
\put(0.8,-0.8){\vector(1,-1){8.4}}
\put(5.5,-5){{\tiny{1}}}
\end{picture}
}

\newcommand{\EastDownlineubderone}{
\setlength{\unitlength}{1mm}
\begin{picture}(10,10)(2.74,0)
\put(0.8,-0.8){\line(1,-1){8.4}}
\put(4,-6.5){{\tiny{1}}}
\end{picture}
}

\newcommand{\holizontallineAbovetwo}{
\setlength{\unitlength}{1mm}
\begin{picture}(10,10)(2.74,0)
\put(1,0){\line(1,0){8}}
\put(4.5,0.5){{\tiny{2}}}
\end{picture}
}

\newcommand{\holizontallineAbovetwoTHICK}{
\setlength{\unitlength}{1mm}
\begin{picture}(10,10)(2.74,0)
\put(1,0){\line(1,0){8}}\multiput(1,0)(0,0.05){4}{\line(1,0){8}}\multiput(1,0)(0,-0.05){4}{\line(1,0){8}}
\put(4.5,0.5){{\tiny{2}}}
\end{picture}
}

\newcommand{\holizontallineAbovetwoTHICKbar}{
\setlength{\unitlength}{1mm}
\begin{picture}(10,10)(2.74,0)
\put(1,0){\line(1,0){8}}\multiput(1,0)(0,0.05){4}{\line(1,0){8}}\multiput(1,0)(0,-0.05){4}{\line(1,0){8}}
\put(4.5,0.5){{\tiny{${\bar{2}}$}}}
\end{picture}
}

\newcommand{\holizontallineAbovetwoEast}{
\setlength{\unitlength}{1mm}
\begin{picture}(10,10)(2.74,0)
\put(1,0){\vector(1,0){8}}
\put(4.5,0.5){{\tiny{2}}}
\end{picture}
}

\newcommand{\holizontallineAbovetwoWest}{
\setlength{\unitlength}{1mm}
\begin{picture}(10,10)(2.74,0)
\put(9,0){\vector(-1,0){8}}
\put(4.5,0.5){{\tiny{2}}}
\end{picture}
}

\newcommand{\holizontallineUndertwo}{
\setlength{\unitlength}{1mm}
\begin{picture}(10,10)(2.74,0)
\put(1,0){\line(1,0){8}}
\put(4.5,-2){{\tiny{2}}}
\end{picture}
}

\newcommand{\verticallineLefttwo}{
\setlength{\unitlength}{1mm}
\begin{picture}(10,10)(2.74,0)
\put(0,1){\line(0,1){8}}
\put(-1.7,4.5){{\tiny{2}}}
\end{picture}
}

\newcommand{\verticallineLefttwoUp}{
\setlength{\unitlength}{1mm}
\begin{picture}(10,10)(2.74,0)
\put(0,1){\vector(0,1){8}}
\put(-1.7,4.5){{\tiny{2}}}
\end{picture}
}

\newcommand{\verticallineLefttwoDown}{
\setlength{\unitlength}{1mm}
\begin{picture}(10,10)(2.74,0)
\put(0,9){\line(0,-1){8}}
\put(-1.7,4.5){{\tiny{2}}}
\end{picture}
}

\newcommand{\verticallineRighttwo}{
\setlength{\unitlength}{1mm}
\begin{picture}(10,10)(2.74,0)
\put(0,1){\line(0,1){8}}
\put(0.5,4.5){{\tiny{2}}}
\end{picture}
}

\newcommand{\verticallineRighttwoTHICK}{
\setlength{\unitlength}{1mm}
\begin{picture}(10,10)(2.74,0)
\put(0,1){\line(0,1){8}}\multiput(0,1)(0.05,0){4}{\line(0,1){8}}\multiput(0,1)(-0.05,0){4}{\line(0,1){8}}
\put(0.5,4.5){{\tiny{2}}}
\end{picture}
}

\newcommand{\verticallineRighttwoUp}{
\setlength{\unitlength}{1mm}
\begin{picture}(10,10)(2.74,0)
\put(0,1){\vector(0,1){8}}
\put(0.5,4.5){{\tiny{2}}}
\end{picture}
}

\newcommand{\verticallineRighttwoDown}{
\setlength{\unitlength}{1mm}
\begin{picture}(10,10)(2.74,0)
\put(0,9){\vector(0,-1){8}}
\put(0.5,4.5){{\tiny{2}}}
\end{picture}
}

\newcommand{\EastUplineabovetwo}{
\setlength{\unitlength}{1mm}
\begin{picture}(10,10)(2.74,0)
\put(0.8,0.8){\line(1,1){8.4}}
\put(4,5.5){{\tiny{2}}}
\end{picture}
}

\newcommand{\EastUplineabovetwoTHICK}{
\setlength{\unitlength}{1mm}
\begin{picture}(10,10)(2.74,0)
\put(0.8,0.8){\line(1,1){8.4}}\multiput(0.8,0.8)(0.03,-0.03){4}{\line(1,1){8.4}}\multiput(0.8,0.8)(-0.03,0.03){4}{\line(1,1){8.4}}
\put(4,5.5){{\tiny{2}}}
\end{picture}
}

\newcommand{\EastUplineabovetwoUp}{
\setlength{\unitlength}{1mm}
\begin{picture}(10,10)(2.74,0)
\put(0.8,0.8){\vector(1,1){8.4}}
\put(4,5.5){{\tiny{2}}}
\end{picture}
}

\newcommand{\EastUplineabovetwoDown}{
\setlength{\unitlength}{1mm}
\begin{picture}(10,10)(2.74,0)
\put(9.2,9.2){\vector(-1,-1){8.4}}
\put(4,5.5){{\tiny{2}}}
\end{picture}
}

\newcommand{\EastUplineundertwo}{
\setlength{\unitlength}{1mm}
\begin{picture}(10,10)(2.74,0)
\put(0.8,0.8){\line(1,1){8.4}}
\put(5,3.5){{\tiny{2}}}
\end{picture}
}

\newcommand{\EastDownlineabovetwo}{
\setlength{\unitlength}{1mm}
\begin{picture}(10,10)(2.74,0)
\put(0.8,-0.8){\line(1,-1){8.4}}
\put(5.5,-5){{\tiny{2}}}
\end{picture}
}

\newcommand{\EastDownlineabovetwoTHICK}{
\setlength{\unitlength}{1mm}
\begin{picture}(10,10)(2.74,0)
\put(0.8,-0.8){\line(1,-1){8.4}}\multiput(0.8,-0.8)(0.03,0.03){4}{\line(1,-1){8.4}}\multiput(0.8,-0.8)(-0.03,-0.03){4}{\line(1,-1){8.4}}
\put(5.5,-5){{\tiny{2}}}
\end{picture}
}

\newcommand{\EastDownlineabovetwoUp}{
\setlength{\unitlength}{1mm}
\begin{picture}(10,10)(2.74,0)
\put(9.2,-9.2){\vector(-1,1){8.4}}
\put(5.5,-5){{\tiny{2}}}
\end{picture}
}

\newcommand{\EastDownlineabovetwoDown}{
\setlength{\unitlength}{1mm}
\begin{picture}(10,10)(2.74,0)
\put(0.8,-0.8){\vector(1,-1){8.4}}
\put(5.5,-5){{\tiny{2}}}
\end{picture}
}

\newcommand{\EastDownlineubdertwo}{
\setlength{\unitlength}{1mm}
\begin{picture}(10,10)(2.74,0)
\put(0.8,-0.8){\line(1,-1){8.4}}
\put(4,-6.5){{\tiny{2}}}
\end{picture}
}

\newcommand{\holizontallineAbovethree}{
\setlength{\unitlength}{1mm}
\begin{picture}(10,10)(2.74,0)
\put(1,0){\line(1,0){8}}
\put(4.5,0.5){{\tiny{3}}}
\end{picture}
}

\newcommand{\holizontallineAbovethreeTHICK}{
\setlength{\unitlength}{1mm}
\begin{picture}(10,10)(2.74,0)
\put(1,0){\line(1,0){8}}\multiput(1,0)(0,0.05){4}{\line(1,0){8}}\multiput(1,0)(0,-0.05){4}{\line(1,0){8}}
\put(4.5,0.5){{\tiny{3}}}
\end{picture}
}

\newcommand{\holizontallineAbovethreeTHICKbar}{
\setlength{\unitlength}{1mm}
\begin{picture}(10,10)(2.74,0)
\put(1,0){\line(1,0){8}}\multiput(1,0)(0,0.05){4}{\line(1,0){8}}\multiput(1,0)(0,-0.05){4}{\line(1,0){8}}
\put(4.5,0.5){{\tiny{${\bar{3}}$}}}
\end{picture}
}

\newcommand{\holizontallineAbovethreeEast}{
\setlength{\unitlength}{1mm}
\begin{picture}(10,10)(2.74,0)
\put(1,0){\vector(1,0){8}}
\put(4.5,0.5){{\tiny{3}}}
\end{picture}
}

\newcommand{\holizontallineAbovethreeWest}{
\setlength{\unitlength}{1mm}
\begin{picture}(10,10)(2.74,0)
\put(9,0){\vector(-1,0){8}}
\put(4.5,0.5){{\tiny{3}}}
\end{picture}
}

\newcommand{\holizontallineUnderthree}{
\setlength{\unitlength}{1mm}
\begin{picture}(10,10)(2.74,0)
\put(1,0){\line(1,0){8}}
\put(4.5,-2){{\tiny{3}}}
\end{picture}
}

\newcommand{\verticallineLeftthree}{
\setlength{\unitlength}{1mm}
\begin{picture}(10,10)(2.74,0)
\put(0,1){\line(0,1){8}}
\put(-1.7,4.5){{\tiny{3}}}
\end{picture}
}

\newcommand{\verticallineLeftthreeUp}{
\setlength{\unitlength}{1mm}
\begin{picture}(10,10)(2.74,0)
\put(0,1){\vector(0,1){8}}
\put(-1.7,4.5){{\tiny{3}}}
\end{picture}
}

\newcommand{\verticallineLeftthreeDown}{
\setlength{\unitlength}{1mm}
\begin{picture}(10,10)(2.74,0)
\put(0,9){\line(0,-1){8}}
\put(-1.7,4.5){{\tiny{3}}}
\end{picture}
}

\newcommand{\verticallineRightthree}{
\setlength{\unitlength}{1mm}
\begin{picture}(10,10)(2.74,0)
\put(0,1){\line(0,1){8}}
\put(0.5,4.5){{\tiny{3}}}
\end{picture}
}

\newcommand{\verticallineRightthreeTHICK}{
\setlength{\unitlength}{1mm}
\begin{picture}(10,10)(2.74,0)
\put(0,9){\line(0,-1){8}}\multiput(0,9)(0.05,0){4}{\line(0,-1){8}}\multiput(0,9)(-0.05,0){4}{\line(0,-1){8}}
\put(0.75,4.5){{\tiny{3}}}
\end{picture}
}

\newcommand{\verticallineRightthreeUp}{
\setlength{\unitlength}{1mm}
\begin{picture}(10,10)(2.74,0)
\put(0,1){\vector(0,1){8}}
\put(0.5,4.5){{\tiny{3}}}
\end{picture}
}

\newcommand{\verticallineRightthreeDown}{
\setlength{\unitlength}{1mm}
\begin{picture}(10,10)(2.74,0)
\put(0,9){\vector(0,-1){8}}
\put(0.5,4.5){{\tiny{3}}}
\end{picture}
}

\newcommand{\EastUplineabovethree}{
\setlength{\unitlength}{1mm}
\begin{picture}(10,10)(2.74,0)
\put(0.8,0.8){\line(1,1){8.4}}
\put(4,5.5){{\tiny{3}}}
\end{picture}
}

\newcommand{\EastUplineabovethreeTHICK}{
\setlength{\unitlength}{1mm}
\begin{picture}(10,10)(2.74,0)
\put(0.8,0.8){\line(1,1){8.4}}\multiput(0.8,0.8)(0.03,-0.03){4}{\line(1,1){8.4}}\multiput(0.8,0.8)(-0.03,0.03){4}{\line(1,1){8.4}}
\put(4,5.5){{\tiny{3}}}
\end{picture}
}

\newcommand{\EastUplineabovethreeUp}{
\setlength{\unitlength}{1mm}
\begin{picture}(10,10)(2.74,0)
\put(0.8,0.8){\vector(1,1){8.4}}
\put(4,5.5){{\tiny{3}}}
\end{picture}
}

\newcommand{\EastUplineabovethreeDown}{
\setlength{\unitlength}{1mm}
\begin{picture}(10,10)(2.74,0)
\put(9.2,9.2){\vector(-1,-1){8.4}}
\put(4,5.5){{\tiny{3}}}
\end{picture}
}

\newcommand{\EastUplineunderthree}{
\setlength{\unitlength}{1mm}
\begin{picture}(10,10)(2.74,0)
\put(0.8,0.8){\line(1,1){8.4}}
\put(5,3.5){{\tiny{3}}}
\end{picture}
}

\newcommand{\EastDownlineabovethree}{
\setlength{\unitlength}{1mm}
\begin{picture}(10,10)(2.74,0)
\put(0.8,-0.8){\line(1,-1){8.4}}
\put(5.5,-5){{\tiny{3}}}
\end{picture}
}

\newcommand{\EastDownlineabovethreeTHICK}{
\setlength{\unitlength}{1mm}
\begin{picture}(10,10)(2.74,0)
\put(0.8,-0.8){\line(1,-1){8.4}}\multiput(0.8,-0.8)(0.03,0.03){4}{\line(1,-1){8.4}}\multiput(0.8,-0.8)(-0.03,-0.03){4}{\line(1,-1){8.4}}
\put(5.5,-5){{\tiny{3}}}
\end{picture}
}

\newcommand{\EastDownlineabovethreeUp}{
\setlength{\unitlength}{1mm}
\begin{picture}(10,10)(2.74,0)
\put(9.2,-9.2){\vector(-1,1){8.4}}
\put(5.5,-5){{\tiny{3}}}
\end{picture}
}

\newcommand{\EastDownlineabovethreeDown}{
\setlength{\unitlength}{1mm}
\begin{picture}(10,10)(2.74,0)
\put(0.8,-0.8){\vector(1,-1){8.4}}
\put(5.5,-5){{\tiny{3}}}
\end{picture}
}

\newcommand{\EastDownlineubderthree}{
\setlength{\unitlength}{1mm}
\begin{picture}(10,10)(2.74,0)
\put(0.8,-0.8){\line(1,-1){8.4}}
\put(4,-6.5){{\tiny{3}}}
\end{picture}
}

\newcommand{\curveNorthWest}{
\setlength{\unitlength}{1mm}
\begin{picture}(10,10)(2.74,0)
\put(0,0){\oval(20,20)[tl]}
\end{picture}
}

\newcommand{\curveNorthEast}{
\setlength{\unitlength}{1mm}
\begin{picture}(10,10)(2.74,0)
\put(0,0){\oval(20,20)[tr]}
\end{picture}
}

\newcommand{\curveNorthEastSE}{
\setlength{\unitlength}{1mm}
\begin{picture}(10,10)(2.74,0)
\put(0,1){\oval(20,18)[tr]}
\end{picture}
}

\newcommand{\curveSouthWest}{
\setlength{\unitlength}{1mm}
\begin{picture}(10,10)(2.74,0)
\put(0,0){\oval(20,20)[bl]}
\end{picture}
}

\newcommand{\curveSouthWestNW}{
\setlength{\unitlength}{1mm}
\begin{picture}(10,10)(2.74,0)
\put(0,-1){\oval(20,18)[bl]}
\end{picture}
}

\newcommand{\curveSouthEast}{
\setlength{\unitlength}{1mm}
\begin{picture}(10,10)(2.74,0)
\put(0,0){\oval(20,20)[br]}
\end{picture}
}

\newcommand{\curveSouthWestTHICK}{
\setlength{\unitlength}{1mm}
\begin{picture}(10,10)(2.74,0)

\put(0,0){\oval(20.20,20.40)[bl]}
\put(0,0){\oval(20.15,20.30)[bl]}
\put(0,0){\oval(20.10,20.20)[bl]}
\put(0,0){\oval(20.05,20.10)[bl]}
\put(0,0){\oval(20,20)[bl]}
\put(0,0){\oval(19.90,19.90)[bl]}
\put(0,0){\oval(19.80,19.80)[bl]}
\put(0,0){\oval(19.70,19.70)[bl]}
\put(0,0){\oval(19.70,19.60)[bl]}

\end{picture}
}

\newcommand{\curveSouthEastTHICK}{
\setlength{\unitlength}{1mm}
\begin{picture}(10,10)(2.74,0)

\put(0,0){\oval(20.20,20.40)[br]}
\put(0,0){\oval(20.15,20.30)[br]}
\put(0,0){\oval(20.10,20.20)[br]}
\put(0,0){\oval(20.05,20.10)[br]}
\put(0,0){\oval(20,20)[br]}
\put(0,0){\oval(19.90,19.90)[br]}
\put(0,0){\oval(19.80,19.80)[br]}
\put(0,0){\oval(19.70,19.70)[br]}
\put(0,0){\oval(19.70,19.60)[br]}

\end{picture}
}

\newcommand{\curveNorthWestTHICK}{
\setlength{\unitlength}{1mm}
\begin{picture}(10,10)(2.74,0)
\put(0,0){\oval(20.20,20.40)[tl]}
\put(0,0){\oval(20.15,20.30)[tl]}
\put(0,0){\oval(20.10,20.20)[tl]}
\put(0,0){\oval(20.05,20.10)[tl]}
\put(0,0){\oval(20,20)[tl]}
\put(0,0){\oval(19.90,19.90)[tl]}
\put(0,0){\oval(19.80,19.80)[tl]}
\put(0,0){\oval(19.70,19.70)[tl]}
\put(0,0){\oval(19.70,19.60)[tl]}
\end{picture}
}

\newcommand{\curveNorthEastTHICK}{
\setlength{\unitlength}{1mm}
\begin{picture}(10,10)(2.74,0)
\put(0,0){\oval(20.40,20.40)[tr]}
\put(0,0){\oval(20.30,20.30)[tr]}
\put(0,0){\oval(20.20,20.20)[tr]}
\put(0,0){\oval(20.10,20.10)[tr]}
\put(0,0){\oval(20,20)[tr]}
\put(0,0){\oval(19.90,19.90)[tr]}
\put(0,0){\oval(19.80,19.80)[tr]}
\put(0,0){\oval(19.70,19.70)[tr]}
\put(0,0){\oval(19.70,19.60)[tr]}
\end{picture}
}

\newcommand{\acn}{n}
\newcommand{\acl}{l}
\newcommand{\acchi}{\mathcal{D}}
\newcommand{\acal}{\al}

\newcommand{\EXAMPLEHecNinea}{
\setlength{\unitlength}{1mm}
\begin{picture}(40,20)(2,10)
\put(-12,10){$\acchi^{(2,9)}_2=$}
\put(10,0){\circle{2}}\put(25,0){\circle{2}}\put(17.5,15){\circle{2}}
\put(11,0){\line(1,0){13}}\put(24.6,1){\line(-1,2){6.55}}\put(10.4,1){\line(1,2){6.55}}
\put(16.5,-4.5){$s$}\put(10,7.5){$q$}\put(22.5,7.5){$r$}
\put(2,-1){$-1$}\put(28,-1){$-1$}\put(15,18){{\small{$-1$}}}
\put(20,14){$\acal_2$}\put(9,-5){$\acal_1$}\put(24,-5){$\acal_3$}
\end{picture}
}

\newcommand{\EXAMPLEHecNineb}{
\setlength{\unitlength}{1mm}
\begin{picture}(40,20)(5,-5)
\put(-10,-1){$\acchi^{(2,9)}_1=$}
\put(10,0){\circle{2}}\put(25,0){\circle{2}}\put(40,0){\circle{2}}
\put(11,0){\line(1,0){13}}\put(26,0){\line(1,0){13}}
\put(9,3){{\small{$q$}}}\put(22,3){{\small{$-1$}}}\put(39,3){{\small{$r$}}}
\put(15.5,1.5){{\small{$q^{-1}$}}}\put(30.5,1.5){{\small{$r^{-1}$}}}
\put(8,-6){$\acal_1$}\put(23,-6){$\acal_2$}\put(38,-6){$\acal_3$}
\end{picture}
}

\newcommand{\EXAMPLEHecNinec}{
\setlength{\unitlength}{1mm}
\begin{picture}(40,20)(5,-5)
\put(-10,-1){$\acchi^{(2,9)}_3=$}
\put(10,0){\circle{2}}\put(25,0){\circle{2}}\put(40,0){\circle{2}}
\put(11,0){\line(1,0){13}}\put(26,0){\line(1,0){13}}
\put(9,3){{\small{$q$}}}\put(22,3){{\small{$-1$}}}\put(39,3){{\small{$s$}}}
\put(15.5,1.5){{\small{$q^{-1}$}}}\put(30.5,1.5){{\small{$s^{-1}$}}}
\put(8,-6){$\acal_1$}\put(23,-6){$\acal_2$}\put(38,-6){$\acal_3$}
\end{picture}
}

\newcommand{\EXAMPLEHecNined}{
\setlength{\unitlength}{1mm}
\begin{picture}(40,20)(5,-5)
\put(-10,-1){$\acchi^{(2,9)}_4=$}
\put(10,0){\circle{2}}\put(25,0){\circle{2}}\put(40,0){\circle{2}}
\put(11,0){\line(1,0){13}}\put(26,0){\line(1,0){13}}
\put(9,3){{\small{$s$}}}\put(22,3){{\small{$-1$}}}\put(39,3){{\small{$r$}}}
\put(15.5,1.5){{\small{$s^{-1}$}}}\put(30.5,1.5){{\small{$r^{-1}$}}}
\put(8,-6){$\acal_1$}\put(23,-6){$\acal_2$}\put(38,-6){$\acal_3$}
\end{picture}
}

\newcommand{\EXAMPLEHecThreeNine}{
\setlength{\unitlength}{1mm}
\begin{picture}(120,45)(0,0)
\put(0,35){\EXAMPLEHecNineb} 
\put(70,35){\EXAMPLEHecNinea} 
\put(0,15){\EXAMPLEHecNinec}
\put(0,0){\EXAMPLEHecNined}
\put(45,5){$(q,r,s\in\bK\setminus\{1\}$, $qrs=1$, $q\ne r\ne s\ne q)$}
\end{picture}
}

\newcommand{\EXAMPLEHecFourEightThree}{
\setlength{\unitlength}{1mm}
\begin{picture}(40,35)(0,-5)
\put(-12,13){$\acchi^{(3,18)}_3=$}
\put(10,15){\circle{2}}\put(25,15){\circle{2}}\put(40,0){\circle{2}}\put(40,30){\circle{2}}
\put(11,15){\line(1,0){13}}
\put(25.7,14.3){\line(1,-1){13.6}}\put(40,1){\line(0,1){28}}\put(25.7,15.7){\line(1,1){13.6}}

\put(7,18){{\small{$\defzeta^{-1}$}}}\put(22,18){{\small{$-1$}}}
\put(43,29){{\small{$-1$}}}\put(43,-1){{\small{$\defzeta$}}}

\put(16.5,16.5){{\small{$\defzeta$}}}
\put(28,24.5){{\small{$\defzeta^{-1}$}}}
\put(28,1.5){{\small{$\defzeta^{-1}$}}}
\put(42,14){{{\small{$\defzeta^{-1}$}}}}

\put(38,33){$\acal_3$}
\put(8,10){$\acal_1$}
\put(23,10){$\acal_2$}
\put(38,-5){$\acal_4$}
\end{picture}
}

\newcommand{\EXAMPLEHecFourEightThreeDASH}{
\setlength{\unitlength}{1mm}
\begin{picture}(40,40)(0,-5)
\put(-29,13){$\acchi^{(3,18)}_3\circ
\bigl[{{1,2,3,4}\atop{4,2,3,1}}\bigr]
=$}
\put(10,15){\circle{2}}\put(25,15){\circle{2}}\put(40,0){\circle{2}}\put(40,30){\circle{2}}
\put(11,15){\line(1,0){13}}
\put(25.7,14.3){\line(1,-1){13.6}}\put(40,1){\line(0,1){28}}\put(25.7,15.7){\line(1,1){13.6}}

\put(7,18){{\small{$\defzeta^{-1}$}}}\put(22,18){{\small{$-1$}}}
\put(43,29){{\small{$-1$}}}\put(43,-1){{\small{$\defzeta$}}}

\put(16.5,16.5){{\small{$\defzeta$}}}
\put(28,24.5){{\small{$\defzeta^{-1}$}}}
\put(28,1.5){{\small{$\defzeta^{-1}$}}}
\put(42,14){{{\small{$\defzeta^{-1}$}}}}

\put(38,33){$\acal_3$}
\put(8,10){$\acal_4$}
\put(23,10){$\acal_2$}
\put(38,-5){$\acal_1$}
\end{picture}
}

\newcommand{\ExampleCayleyA}{
\setlength{\unitlength}{1mm}
\begin{picture}(60,60)(25,7)

\put(10,19){\curveSouthWestTHICK}
\multiput(10,9.05)(0,0.05){3}{\line(1,0){90}}
\multiput(10,8.95)(0,-0.05){3}{\line(1,0){90}}
\put(10,9){\line(1,0){90}}
\put(100,19){\curveSouthEastTHICK}
\put(54.5,10){{\tiny{3}}}

\put(20,51){\curveNorthWest}
\put(20,61){\line(1,0){70}}
\put(90,51){\curveNorthEast}
\put(54.5,58){{\tiny{2}}}

\put(10,61){\curveNorthWest}
\put(10,71){\line(1,0){90}}
\put(100,61){\curveNorthEast}
\multiput(0,31)(110,0){2}{\line(0,1){30}}
\put(54.5,68){{\tiny{3}}}

\put(0,20){\whitecircleempty}\put(40,20){\whitecircleempty}\put(80,20){\whitecircleempty}
\put(0,20){\verticallineRightoneTHICK}
\put(40,20){\verticallineRightoneTHICK}\put(80,20){\verticallineRightoneTHICK}
\put(1,20){\line(1,0){10}}\put(41,20){\line(1,0){10}}\put(81,20){\line(1,0){10}}
\multiput(10,20)(40,0){3}{\holizontallineAbovetwo}
\multiput(19,20)(40,0){3}{\line(1,0){10}}
\multiput(30,20)(40,0){3}{\whitecircleempty}
\multiput(30,20)(40,0){2}{\verticallineRightoneTHICK}\put(110,20){\verticallineRightoneTHICK}
\multiput(30,20)(40,0){2}{\holizontallineAbovethreeTHICK}

\multiput(0,30)(40,0){3}{\whitecircleempty}
\multiput(0,30)(40,0){3}{\EastUplineabovetwoTHICK}
\multiput(30,30)(40,0){3}{\whitecircleempty}
\multiput(30,30)(40,0){2}{\holizontallineAbovethree}

\multiput(10,40)(40,0){3}{\whitecircleempty}
\multiput(10,40)(40,0){3}{\verticallineRightthreeTHICK}
\multiput(10,40)(40,0){3}{\holizontallineAboveone}
\multiput(20,40)(40,0){3}{\whitecircleempty}
\multiput(20,40)(40,0){3}{\verticallineRightthreeTHICK}
\multiput(20,40)(40,0){3}{\EastDownlineabovetwoTHICK}

\multiput(10,50)(40,0){3}{\whitecircleempty}
\multiput(10,50)(40,0){3}{\holizontallineAboveoneTHICK}
\multiput(20,50)(40,0){3}{\whitecircleempty}
\multiput(21,50)(40,0){2}{\line(1,0){10}}
\multiput(30,50)(40,0){2}{\holizontallineAbovetwo}
\multiput(39,50)(40,0){2}{\line(1,0){10}}

\end{picture}
}

\newcommand{\AtwoAone}{
\setlength{\unitlength}{1mm}
\begin{picture}(80,20)(5,-5)
\put(8,10){$a:=[\acchi^{(1,1)}_1\times\acchi^{(1,0)}_1]$}
\put(10,0){\circle{2}}\put(25,0){\circle{2}}\put(40,0){\circle{2}}
\put(11,0){\line(1,0){13}}
\put(9,3){{\small{$q$}}}\put(24,3){{\small{$q$}}}\put(39,3){{\small{$r$}}}
\put(15.5,1.5){{\small{$q^{-1}$}}}
\put(8,-6){$\al_1$}\put(23,-6){$\al_2$}\put(38,-6){$\al_3$}
\put(50,-3){$(\chamq, r\in\bKt, \chamq\ne 1)$}
\end{picture}
}

\newcommand{\aTotalCHAtwoAone}{
\setlength{\unitlength}{0.5mm}
\begin{picture}(50,45)(00,00)

\multiput(31,40)(00,0){1}{\line(1,0){9}}
\multiput(40,30)(00,0){1}{\oval(20,20)[tr]}
\multiput(50,21)(00,0){1}{\line(0,1){9}}

\multiput(20,40)(10,0){2}{\circle{2}}
\multiput(21,40)(00,0){1}{\line(1,0){8}}

\multiput(20,30)(10,0){2}{\circle{2}}
\multiput(20,31)(10,0){2}{\line(0,1){8}}
\multiput(21,30)(00,0){1}{\line(1,0){8}}

\multiput(00,20)(40,0){2}{\circle{2}}
\multiput(01,20)(40,0){2}{\line(1,0){8}}
\multiput(10,20)(40,0){2}{\circle{2}}
\multiput(10.8,20.8)(00,0){1}{\line(1,1){8.4}}
\multiput(39.2,20.8)(00,0){1}{\line(-1,1){8.4}}

\multiput(20,10)(10,0){2}{\circle{2}}
\multiput(19.2,10.8)(00,0){1}{\line(-1,1){8.4}}
\multiput(21,10)(00,0){1}{\line(1,0){8}}
\multiput(30.8,10.8)(00,0){1}{\line(1,1){8.4}}

\multiput(20,00)(10,0){2}{\circle{2}}
\multiput(20,01)(10,0){2}{\line(0,1){8}}
\multiput(21,00)(00,0){1}{\line(1,0){8}}

\multiput(20,40)(00,0){1}{\RIGHTlineTHICKsmall{8}}

\multiput(20,30)(00,0){1}{\RIGHTlineTHICKsmall{8}}
\multiput(30,30)(00,0){1}{\UPlineTHICKsmall{8}}

\multiput(10,20)(00,0){1}{\UPRIGHTlineTHICKsmall{8.4}}
\multiput(40,20)(00,0){1}{\RIGHTlineTHICKsmall{8}}

\multiput(20,10)(00,0){1}{\UPLEFTlineTHICKsmall{8.4}}
\multiput(20,10)(00,0){1}{\RIGHTlineTHICKsmall{8}}
\multiput(30,10)(00,0){1}{\UPRIGHTlineTHICKsmall{8.4}}

\multiput(20,00)(00,0){1}{\RIGHTlineTHICKsmall{8}}

\multiput(10,40)(00,0){1}{\properRIGHTlineTHICKsmall{9}}

\put(10,30){\curveNorthWestTHICKsmall}
\multiput(20,30)(00,0){1}{\RIGHTlineTHICKsmall{8}}

\multiput(00,21)(00,0){1}{\properUPlineTHICKsmall{9}}

\multiput(00,10)(50,0){2}{\properUPlineTHICKsmall{9}}

\put(10,10){\curveSouthWestTHICKsmall}
\multiput(10,00)(21,0){2}{\properRIGHTlineTHICKsmall{9}}
\put(40,10){\curveSouthEastTHICKsmall}

\multiput(9,41.5)(30,0){2}{{\tiny{2}}}
\multiput(24,41.5)(00,0){1}{{\tiny{1}}}

\multiput(16.5,33.5)(00,0){1}{{\tiny{3}}}
\multiput(31,33.5)(00,0){1}{{\tiny{3}}}

\multiput(24,31)(00,0){1}{{\tiny{1}}}

\multiput(16,22)(00,0){1}{{\tiny{2}}}
\multiput(31,22)(00,0){1}{{\tiny{2}}}

\multiput(4,21)(40,0){2}{{\tiny{3}}}

\multiput(16,15)(00,0){1}{{\tiny{1}}}
\multiput(31,15)(00,0){1}{{\tiny{1}}}

\multiput(24,5.5)(00,0){1}{{\tiny{2}}}

\multiput(16.5,3.5)(00,0){1}{{\tiny{3}}}
\multiput(31,3.5)(00,0){1}{{\tiny{3}}}

\multiput(9,-4.5)(30,0){2}{{\tiny{1}}}
\multiput(24,-4.5)(00,0){1}{{\tiny{2}}}

\multiput(19,42)(10,0){2}{{\tiny{$a$}}}
\multiput(19,26)(9,0){2}{{\tiny{$a$}}}
\multiput(8,16)(31.5,0){2}{{\tiny{$a$}}}
\multiput(-5,19)(57,0){2}{{\tiny{$a$}}}

\multiput(19,12)(9,0){2}{{\tiny{$a$}}}
\multiput(19,-4)(10,0){2}{{\tiny{$a$}}}
\end{picture}
}

\newcommand{\HecOnea}{
\setlength{\unitlength}{1mm}
\begin{picture}(80,20)(5,-5)
\put(8,10){$a^\prime:=[\acchi^{(2,1)}_1]$}
\put(10,0){\circle{2}}\put(25,0){\circle{2}}\put(40,0){\circle{2}}
\put(11,0){\line(1,0){13}}\put(26,0){\line(1,0){13}}
\put(9,3){{\small{$q$}}}\put(24,3){{\small{$q$}}}\put(39,3){{\small{$q$}}}
\put(15.5,1.5){{\small{$q^{-1}$}}}\put(30.5,1.5){{\small{$q^{-1}$}}}
\put(8,-6){$\al_1$}\put(23,-6){$\al_2$}\put(38,-6){$\al_3$}
\put(50,-3){$(q\in\bKt, q\ne 1)$}
\end{picture}
}

\newcommand{\HecFoura}{
\setlength{\unitlength}{1mm}
\begin{picture}(40,20)(5,-5)
\put(8,10){$a^{\prime\prime}:=[\acchi^{(2,4)}_1]$}
\put(10,0){\circle{2}}\put(25,0){\circle{2}}\put(40,0){\circle{2}}
\put(11,0){\line(1,0){13}}\put(26,0){\line(1,0){13}}
\put(8,3){{\small{$-1$}}}\put(24,3){{\small{$q$}}}\put(39,3){{\small{$q$}}}
\put(15.5,1.5){{\small{$q^{-1}$}}}\put(30.5,1.5){{\small{$q^{-1}$}}}
\put(8,-6){$\al_1$}\put(23,-6){$\al_2$}\put(38,-6){$\al_3$}
\end{picture}
}

\newcommand{\HecFourb}{
\setlength{\unitlength}{1mm}
\begin{picture}(40,20)(5,-5)
\put(8,10){$b^{\prime\prime}:=[\acchi^{(2,4)}_2]$}
\put(10,0){\circle{2}}\put(25,0){\circle{2}}\put(40,0){\circle{2}}
\put(11,0){\line(1,0){13}}\put(26,0){\line(1,0){13}}
\put(8,3){{\small{$-1$}}}\put(23,3){{\small{$-1$}}}\put(39,3){{\small{$q$}}}
\put(16.5,1.5){{\small{$q$}}}\put(30.5,1.5){{\small{$q^{-1}$}}}
\put(8,-6){$\al_1$}\put(23,-6){$\al_2$}\put(38,-6){$\al_3$}
\end{picture}
}

\newcommand{\HecFourc}{
\setlength{\unitlength}{1mm}
\begin{picture}(40,20)(5,-5)
\put(8,10){$c^{\prime\prime}:=b^{\prime\prime}\circ\left[{{123}\atop{321}}\right]$}
\put(10,0){\circle{2}}\put(25,0){\circle{2}}\put(40,0){\circle{2}}
\put(11,0){\line(1,0){13}}\put(26,0){\line(1,0){13}}
\put(9,3){{\small{$q$}}}\put(23,3){{\small{$-1$}}}\put(38,3){{\small{$-1$}}}
\put(15.5,1.5){{\small{$q^{-1}$}}}\put(31.5,1.5){{\small{$q$}}}
\put(8,-6){$\al_1$}\put(23,-6){$\al_2$}\put(38,-6){$\al_3$}
\end{picture}
}

\newcommand{\HecFourd}{
\setlength{\unitlength}{1mm}
\begin{picture}(40,20)(5,-5)
\put(8,10){$d^{\prime\prime}:=a^{\prime\prime}\circ\left[{{123}\atop{321}}\right]$}
\put(10,0){\circle{2}}\put(25,0){\circle{2}}\put(40,0){\circle{2}}
\put(11,0){\line(1,0){13}}\put(26,0){\line(1,0){13}}
\put(9,3){{\small{$q$}}}\put(24,3){{\small{$q$}}}\put(38,3){{\small{$-1$}}}
\put(15.5,1.5){{\small{$q^{-1}$}}}\put(30.5,1.5){{\small{$q^{-1}$}}}
\put(8,-6){$\al_1$}\put(23,-6){$\al_2$}\put(38,-6){$\al_3$}
\end{picture}
}

\newcommand{\HecFourab}
{\setlength{\unitlength}{1mm}
\begin{picture}(90,55)(20,-5)

\put(00,25){\HecFoura}
\put(45,25){\HecFourb}
\put(90,25){\HecFourc}
\put(00,00){\HecFourd}
\put(60,5){$(q\in\bKt,\,q\ne 1)$}

\end{picture}
}

\newcommand{\HecEighta}{
\setlength{\unitlength}{1mm}
\begin{picture}(40,20)(5,-5)
\put(8,10){$a:=[\acchi^{(2,8)}_1]$}
\put(10,0){\circle{2}}\put(25,0){\circle{2}}\put(40,0){\circle{2}}
\put(11,0){\line(1,0){13}}\put(26,0){\line(1,0){13}}
\put(9,3){{\small{$q$}}}\put(23,3){{\small{$-1$}}}\put(38,3){{\small{$q^{-1}$}}}
\put(15.5,1.5){{\small{$q^{-1}$}}}\put(31.5,1.5){{\small{$q$}}}
\put(8,-6){$\al_1$}\put(23,-6){$\al_2$}\put(38,-6){$\al_3$}
\end{picture}
}

\newcommand{\HecEightb}{
\setlength{\unitlength}{1mm}
\begin{picture}(40,20)(5,-5)
\put(8,10){$b:=[\acchi^{(2,8)}_2]$}
\put(10,0){\circle{2}}\put(25,0){\circle{2}}\put(40,0){\circle{2}}
\put(11,0){\line(1,0){13}}\put(26,0){\line(1,0){13}}
\put(8,3){{\small{$-1$}}}\put(23,3){{\small{$-1$}}}\put(38,3){{\small{$-1$}}}
\put(16.5,1.5){{\small{$q$}}}\put(30.5,1.5){{\small{$q^{-1}$}}}
\put(8,-6){$\al_1$}\put(23,-6){$\al_2$}\put(38,-6){$\al_3$}
\end{picture}
}

\newcommand{\HecEightc}{
\setlength{\unitlength}{1mm}
\begin{picture}(40,20)(5,-5)
\put(8,10){$c:=[\acchi^{(2,8)}_3]$}
\put(10,0){\circle{2}}\put(25,0){\circle{2}}\put(40,0){\circle{2}}
\put(11,0){\line(1,0){13}}\put(26,0){\line(1,0){13}}
\put(8,3){{\small{$-1$}}}\put(24,3){{\small{$q$}}}\put(38,3){{\small{$-1$}}}
\put(15.5,1.5){{\small{$q^{-1}$}}}\put(30.5,1.5){{\small{$q^{-1}$}}}
\put(8,-6){$\al_1$}\put(23,-6){$\al_2$}\put(38,-6){$\al_3$}
\end{picture}
}

\newcommand{\HecEightd}{
\setlength{\unitlength}{1mm}
\begin{picture}(40,20)(5,-5)
\put(8,10){$d:=[\acchi^{(2,8)}_4]$}
\put(10,0){\circle{2}}\put(25,0){\circle{2}}\put(40,0){\circle{2}}
\put(11,0){\line(1,0){13}}\put(26,0){\line(1,0){13}}
\put(8,3){{\small{$-1$}}}\put(23,3){{\small{$q^{-1}$}}}\put(38,3){{\small{$-1$}}}
\put(16.5,1.5){{\small{$q$}}}\put(31.5,1.5){{\small{$q$}}}
\put(8,-6){$\al_1$}\put(23,-6){$\al_2$}\put(38,-6){$\al_3$}
\end{picture}
}

\newcommand{\HecEighte}{
\setlength{\unitlength}{1mm}
\begin{picture}(40,20)(5,-5)
\put(8,10){$e:=[\acchi^{(2,8)}_2]\circ\left[{{123}\atop{321}}\right]$}
\put(10,0){\circle{2}}\put(25,0){\circle{2}}\put(40,0){\circle{2}}
\put(11,0){\line(1,0){13}}\put(26,0){\line(1,0){13}}
\put(8,3){{\small{$-1$}}}\put(23,3){{\small{$-1$}}}\put(38,3){{\small{$-1$}}}
\put(16.5,1.5){{\small{$q^{-1}$}}}\put(30.5,1.5){{\small{$q$}}}
\put(8,-6){$\al_1$}\put(23,-6){$\al_2$}\put(38,-6){$\al_3$}
\end{picture}
}

\newcommand{\HecEightf}{
\setlength{\unitlength}{1mm}
\begin{picture}(40,20)(5,-5)
\put(8,10){$f:=[\acchi^{(2,8)}_1]\circ\left[{{123}\atop{321}}\right]$}
\put(10,0){\circle{2}}\put(25,0){\circle{2}}\put(40,0){\circle{2}}
\put(11,0){\line(1,0){13}}\put(26,0){\line(1,0){13}}
\put(8,3){{\small{$q^{-1}$}}}\put(23,3){{\small{$-1$}}}\put(39,3){{\small{$q$}}}
\put(16.5,1.5){{\small{$q$}}}\put(30.5,1.5){{\small{$q^{-1}$}}}
\put(8,-6){$\al_1$}\put(23,-6){$\al_2$}\put(38,-6){$\al_3$}
\end{picture}
}

\newcommand{\HecEightabcd}{
\setlength{\unitlength}{1mm}
\begin{picture}(150,50)(-5,-10)
\put(0,25){$\HecEighta$}\put(45,25){$\HecEightb$}\put(90,25){$\HecEightc$}
\put(0,0){$\HecEightd$}\put(45,0){$\HecEighte$}\put(90,0){$\HecEightf$}
\put(5,-10){$(q\in\bK\setminus\{0,1\})$}
\end{picture}
}

\newcommand{\CayleyA}{
\setlength{\unitlength}{1mm}
\begin{picture}(60,80)(25,-3)

\put(10,19){\curveSouthWest}
\put(10,9){\line(1,0){90}}
\put(100,19){\curveSouthEast}
\put(54.5,10){{\tiny{3}}}

\put(20,51){\curveNorthWest}
\put(20,61){\line(1,0){70}}
\put(90,51){\curveNorthEast}
\put(54.5,58){{\tiny{2}}}

\put(10,61){\curveNorthWest}
\put(10,71){\line(1,0){90}}
\put(100,61){\curveNorthEast}
\multiput(0,31)(110,0){2}{\line(0,1){30}}
\put(54.5,68){{\tiny{3}}}

\put(0,20){\whitecircleempty}\put(40,20){\whitecircleempty}\put(80,20){\whitecircleempty}
\put(0,20){\verticallineRightoneTHICK}
\put(40,20){\verticallineRightone}\put(80,20){\verticallineRightone}
\put(1,20){\line(1,0){10}}\put(41,20){\line(1,0){10}}\put(81,20){\line(1,0){10}}
\multiput(1,20.05)(40,0){3}{\line(1,0){10}}\multiput(1,20.10)(40,0){3}{\line(1,0){10}}\multiput(1,20.15)(40,0){3}{\line(1,0){10}}
\multiput(1,19.95)(40,0){3}{\line(1,0){10}}\multiput(1,19.90)(40,0){3}{\line(1,0){10}}\multiput(1,19.85)(40,0){3}{\line(1,0){10}}
\multiput(10,20)(40,0){3}{\holizontallineAbovetwoTHICK}
\multiput(19,20)(40,0){3}{\line(1,0){10}}
\multiput(19,20.05)(40,0){3}{\line(1,0){10}}\multiput(19,20.10)(40,0){3}{\line(1,0){10}}\multiput(19,20.15)(40,0){3}{\line(1,0){10}}
\multiput(19,19.95)(40,0){3}{\line(1,0){10}}\multiput(19,19.90)(40,0){3}{\line(1,0){10}}\multiput(19,19.85)(40,0){3}{\line(1,0){10}}
\multiput(30,20)(40,0){3}{\whitecircleempty}
\multiput(30,20)(40,0){2}{\verticallineRightone}\put(110,20){\verticallineRightoneTHICK}
\multiput(30,20)(40,0){2}{\holizontallineAbovethreeTHICK}

\multiput(0,30)(40,0){3}{\whitecircleempty}
\multiput(0,30)(40,0){3}{\EastUplineabovetwoTHICK}
\multiput(30,30)(40,0){3}{\whitecircleempty}
\multiput(30,30)(40,0){2}{\holizontallineAbovethreeTHICK}

\multiput(10,40)(40,0){3}{\whitecircleempty}
\multiput(10,40)(40,0){3}{\verticallineRightthreeTHICK}
\multiput(10,40)(40,0){3}{\holizontallineAboveone}
\multiput(20,40)(40,0){3}{\whitecircleempty}
\multiput(20,40)(40,0){3}{\verticallineRightthreeTHICK}
\multiput(20,40)(40,0){3}{\EastDownlineabovetwoTHICK}

\multiput(10,50)(40,0){3}{\whitecircleempty}
\multiput(10,50)(40,0){3}{\holizontallineAboveoneTHICK}
\multiput(20,50)(40,0){3}{\whitecircleempty}
\multiput(21,50)(40,0){2}{\line(1,0){10}}
\multiput(30,50)(40,0){2}{\holizontallineAbovetwo}
\multiput(39,50)(40,0){2}{\line(1,0){10}}

\put(-6.5, 19){$v_{19}$}
\put(-6.5, 29){$v_{13}$}
\put(4.5, 39){$v_7$}\put(44.5, 39){$v_9$}\put(83, 39){$v_{11}$}
\put(4.5, 49){$v_1$}

\put(111.5, 19){$v_{24}$}
\put(111.5, 29){$v_{18}$}
\put(101.5, 39){$v_{12}$}\put(61.5, 39){$v_{10}$}\put(21.5, 39){$v_8$}
\put(101.5, 49){$v_6$}

\put(28, 16){$v_{20}$}\put(38, 16){$v_{21}$}\put(68, 16){$v_{22}$}\put(78, 16){$v_{23}$}
\put(28, 32.5){$v_{14}$}\put(36.5, 32.5){$v_{15}$}\put(68, 32){$v_{16}$}\put(76.5, 32.5){$v_{17}$}
\put(18, 52.5){$v_2$}\put(48, 52.5){$v_3$}\put(58, 52){$v_4$}\put(88, 52.5){$v_5$}

\end{picture}
}

\newcommand{\HecTwod}{
\setlength{\unitlength}{1mm}
\begin{picture}(40,20)(5,10)
\put(8,10){$a^\prime:=[\acchi^{(2,2)}_1]$}
\put(10,0){\circle{2}}\put(25,0){\circle{2}}\put(40,0){\circle{2}}
\put(11,0){\line(1,0){13}}\put(26,0){\line(1,0){13}}
\put(9,3){{\small{$q^2$}}}\put(24,3){{\small{$q^2$}}}\put(39,3){{\small{$q$}}}
\put(15.5,2.0){{\small{$q^{-2}$}}}\put(30.5,2.0){{\small{$q^{-2}$}}}
\put(8,-4){$\al_1$}\put(23,-4){$\al_2$}\put(38,-4){$\al_3$}

\put(8,-10){\footnotesize{$(q\in\bKt$, $q^2\ne 1)$}}
\end{picture}
}

\newcommand{\HecThreee}{
\setlength{\unitlength}{1mm}
\begin{picture}(40,20)(5,10)
\put(8,10){$a^\prime:=[\acchi^{(2,3)}_1]$}
\put(10,0){\circle{2}}\put(25,0){\circle{2}}\put(40,0){\circle{2}}
\put(11,0){\line(1,0){13}}\put(26,0){\line(1,0){13}}
\put(9,3){{\small{$q$}}}\put(24,3){{\small{$q$}}}\put(39,3){{\small{$q^2$}}}
\put(15.5,2.0){{\small{$q^{-1}$}}}\put(30.5,2.0){{\small{$q^{-2}$}}}
\put(8,-4){$\al_1$}\put(23,-4){$\al_2$}\put(38,-4){$\al_3$}

\put(8,-10){\footnotesize{$(q\in\bKt$, $q^2\ne 1)$}}
\end{picture}
}

\newcommand{\HecTwelvef}{
\setlength{\unitlength}{1mm}
\begin{picture}(40,20)(5,10)
\put(8,10){$a^\prime:=[\acchi^{(2,12)}_1]$} 
\put(10,0){\circle{2}}\put(25,0){\circle{2}}\put(40,0){\circle{2}}
\put(11,0){\line(1,0){13}}\put(26,0){\line(1,0){13}}
\put(5.5,3){{\small{$-\defzeta^{-1}$}}}\put(21,3){{\small{$-\defzeta^{-1}$}}}\put(39,3){{\small{$\defzeta$}}}
\put(15,2.0){{\small{$-\defzeta$}}}\put(30,2.0){{\small{$-\defzeta$}}}
\put(8,-4){$\al_1$}\put(23,-4){$\al_2$}\put(38,-4){$\al_3$}

\put(8,-10){\footnotesize{$(\zeta\in\bKt$, $\zeta^2+\zeta+1=0)$}}
\end{picture}
}

\newcommand{\HecTwoThreeTwelvedef}
{\setlength{\unitlength}{1mm}
\begin{picture}(130,20)(0,-20)

\put(0,0){\HecTwod}
\put(45,0){\HecThreee}
\put(90,0){\HecTwelvef}

\end{picture}
}

\newcommand{\HecFivea}{
\setlength{\unitlength}{1mm}
\begin{picture}(40,20)(5,10)
\put(8,10){$a:=[\acchi^{(2,5)}_1]$}
\put(10,0){\circle{2}}\put(25,0){\circle{2}}\put(40,0){\circle{2}}
\put(11,0){\line(1,0){13}}\put(26,0){\line(1,0){13}}
\put(8,3){{\small{$-1$}}}\put(24,3){{\small{$q^2$}}}\put(39,3){{\small{$q$}}}
\put(15.5,2.0){{\small{$q^{-2}$}}}\put(30.5,2.0){{\small{$q^{-2}$}}}
\put(8,-4){$\al_1$}\put(23,-4){$\al_2$}\put(38,-4){$\al_3$}
\end{picture}
}

\newcommand{\HecFiveb}{
\setlength{\unitlength}{1mm}
\begin{picture}(40,20)(5,10)
\put(8,10){$b:=[\acchi^{(2,5)}_2]$}
\put(10,0){\circle{2}}\put(25,0){\circle{2}}\put(40,0){\circle{2}}
\put(11,0){\line(1,0){13}}\put(26,0){\line(1,0){13}}
\put(8,3){{\small{$-1$}}}\put(23,3){{\small{$-1$}}}\put(39,3){{\small{$q$}}}
\put(16,2.0){{\small{$q^2$}}}\put(30.5,2.0){{\small{$q^{-2}$}}}
\put(8,-4){$\al_1$}\put(23,-4){$\al_2$}\put(38,-4){$\al_3$}
\end{picture}
}

\newcommand{\HecFivec}{
\setlength{\unitlength}{1mm}
\begin{picture}(40,20)(5,10)
\put(8,10){$c:=[\acchi^{(2,5)}_3]$}
\put(10,0){\circle{2}}\put(25,0){\circle{2}}\put(40,0){\circle{2}}
\put(11,0){\line(1,0){13}}\put(26,0){\line(1,0){13}}
\put(9,3){{\small{$q^2$}}}\put(23,3){{\small{$-1$}}}\put(38,3){{\small{$-q^{-1}$}}}
\put(15.5,2.0){{\small{$q^{-2}$}}}\put(31,2.0){{\small{$q^2$}}}
\put(8,-4){$\al_1$}\put(23,-4){$\al_2$}\put(38,-4){$\al_3$}
\end{picture}
}

\newcommand{\HecFiveabc}
{\setlength{\unitlength}{1mm}
\begin{picture}(130,15)(0,-20)

\put(0,0){\HecFivea}
\put(45,0){\HecFiveb}
\put(90,0){\HecFivec}

\put(8,-20){\footnotesize{$(q\in\bKt$, $q^2\ne 1)$}}
\end{picture}
}

\newcommand{\HecSixa}{
\setlength{\unitlength}{1mm}
\begin{picture}(40,20)(5,10)
\put(8,10){$a:=[\acchi^{(2,6)}_1]$}
\put(10,0){\circle{2}}\put(25,0){\circle{2}}\put(40,0){\circle{2}}
\put(11,0){\line(1,0){13}}\put(26,0){\line(1,0){13}}
\put(8,3){{\small{$-1$}}}\put(24,3){{\small{$q$}}}\put(39,3){{\small{$q^2$}}}
\put(15.5,2.0){{\small{$q^{-1}$}}}\put(30.5,2.0){{\small{$q^{-2}$}}}
\put(8,-4){$\al_1$}\put(23,-4){$\al_2$}\put(38,-4){$\al_3$}
\end{picture}
}

\newcommand{\HecSixb}{
\setlength{\unitlength}{1mm}
\begin{picture}(40,20)(5,10)
\put(8,10){$b:=[\acchi^{(2,6)}_3]$}
\put(10,0){\circle{2}}\put(25,0){\circle{2}}\put(40,0){\circle{2}}
\put(11,0){\line(1,0){13}}\put(26,0){\line(1,0){13}}
\put(8,3){{\small{$-1$}}}\put(23,3){{\small{$-1$}}}\put(39,3){{\small{$q^2$}}}
\put(17,2.0){{\small{$q$}}}\put(30.5,2.0){{\small{$q^{-2}$}}}
\put(8,-4){$\al_1$}\put(23,-4){$\al_2$}\put(38,-4){$\al_3$}
\end{picture}
}

\newcommand{\HecSixc}{
\setlength{\unitlength}{1mm}
\begin{picture}(40,20)(5,10)
\put(8,15){$c:=[\acchi^{(2,6)}_2]$}
\put(10,0){\circle{2}}\put(25,0){\circle{2}}\put(40,0){\circle{2}}
\put(11,0){\line(1,0){13}}\put(26,0){\line(1,0){13}}
\put(6.5,3){{\small{$q$}}}\put(23,3){{\small{$-1$}}}\put(41,3){{\small{$-1$}}}
\put(15.5,2.0){{\small{$q^{-1}$}}}\put(31,2.0){{\small{$q^2$}}}
\put(8,-4){$\al_1$}\put(23,-4){$\al_2$}\put(38,-4){$\al_3$}
\put(25,1){\oval(30,14)[t]}\put(24,10){{\small{$q^{-1}$}}}
\end{picture}
}

\newcommand{\HecSixd}{
\setlength{\unitlength}{1mm}
\begin{picture}(40,20)(5,10)
\put(8,10){$d:=[\acchi^{(2,6)}_1]\circ\left[{{123}\atop{132}}\right]$}
\put(10,0){\circle{2}}\put(25,0){\circle{2}}\put(40,0){\circle{2}}
\put(11,0){\line(1,0){13}}\put(26,0){\line(1,0){13}}
\put(8,3){{\small{$-1$}}}\put(23,3){{\small{$-1$}}}\put(39,3){{\small{$q^2$}}}
\put(17,2.0){{\small{$q$}}}\put(30.5,2.0){{\small{$q^{-2}$}}}
\put(8,-4){$\al_1$}\put(23,-4){$\al_3$}\put(38,-4){$\al_2$}
\end{picture}
}

\newcommand{\HecSixe}{
\setlength{\unitlength}{1mm}
\begin{picture}(40,20)(5,10)
\put(8,10){$e:=[\acchi^{(2,6)}_3]\circ\left[{{123}\atop{132}}\right]$}
\put(10,0){\circle{2}}\put(25,0){\circle{2}}\put(40,0){\circle{2}}
\put(11,0){\line(1,0){13}}\put(26,0){\line(1,0){13}}
\put(8,3){{\small{$-1$}}}\put(24,3){{\small{$q$}}}\put(39,3){{\small{$q^2$}}}
\put(15.5,2.0){{\small{$q^{-1}$}}}\put(30.5,2.0){{\small{$q^{-2}$}}}
\put(8,-4){$\al_1$}\put(23,-4){$\al_3$}\put(38,-4){$\al_2$}
\end{picture}
}

\newcommand{\HecSixabcde}{
\setlength{\unitlength}{1mm}
\begin{picture}(130,50)(0,-20)
\put(0,25){$\HecSixa$}\put(45,25){$\HecSixb$}\put(90,25){$\HecSixc$}
\put(0,0){$\HecSixd$}
\put(45,0){$\HecSixe$}
\put(90,-10){$(q\in\bKt, q^2\ne 1)$}
\end{picture}
}

\newcommand{\CayleySix}{
\setlength{\unitlength}{1mm}
\begin{picture}(110,60)(-10,0)

\multiput(10,0)(2,0){6}{\line(1,0){0.5}}\put(13,1){$C_i$}
\put(20,0){\holizontallineAbovethree}
\put(30,0){\whitecircleSouthe}
\put(30,0){\holizontallineAbovetwoTHICK}
\put(30,0){\verticallineRightoneTHICK}
\put(40,0){\whitecircleSouthe}
\put(40,0){\holizontallineAbovethree}
\put(40,0){\verticallineRightoneTHICK}
\multiput(50,0)(2,0){6}{\line(1,0){0.5}}\put(51,1){$C_{i+1}$}

\put(30,10){\whitecircleNorthd}
\put(30,10){\holizontallineAbovetwo}
\put(40,10){\whitecircleNorthd}
\put(40,10){\EastUplineabovethreeTHICK}

\multiput(0,20)(2,0){6}{\line(1,0){0.5}}\put(2,21){$B_i$}
\put(10,20){\holizontallineAboveone}
\put(20,20){\whitecircleNorthc}
\put(20,20){\EastDownlineabovethreeTHICK}
\put(20,20){\EastUplineabovetwoTHICK}
\put(50,20){\whitecircleNorthc}
\put(50,20){\holizontallineAboveone}
\multiput(60,20)(2,0){6}{\line(1,0){0.5}}\put(61,21){$B_{i+1}$}

\put(76,21){$(i\in\bZ/4\bZ)$}

\put(30,30){\whitecircleSouthb}
\put(30,30){\holizontallineAbovethree}
\put(30,30){\verticallineRightoneTHICK}
\put(40,30){\whitecircleSouthb}
\put(40,30){\EastDownlineabovetwoTHICK}
\put(40,30){\verticallineRightoneTHICK}

\multiput(10,40)(2,0){6}{\line(1,0){0.5}}\put(13,41){$A_i$}
\put(20,40){\holizontallineAbovetwoTHICK}
\put(30,40){\whitecircleNortha}
\put(30,40){\holizontallineAbovethree}
\put(40,40){\whitecircleNortha}
\put(40,40){\holizontallineAbovetwoTHICK}
\multiput(50,40)(2,0){6}{\line(1,0){0.5}}\put(51,41){$A_{i+1}$}

\end{picture}
}

\newcommand{\WholeCayleySix}{
\setlength{\unitlength}{1mm}
\begin{picture}(120,160)(-5,0)

\put(20,10){\curveSouthWest}

\put(20,80){\curveNorthWest}

\put(20,0){\line(1,0){9}}
\put(50,0){\holizontallineAbovethree}
\multiput(30,0)(40,0){2}{\whitecircleSouthe}
\multiput(30,0)(40,0){2}{\holizontallineAbovetwoTHICK}
\multiput(30,0)(40,0){2}{\verticallineRightoneTHICK}
\multiput(40,0)(40,0){2}{\whitecircleSouthe}
\multiput(41,0)(18,0){2}{\holizontallinetenmm}
\multiput(40,0)(40,0){2}{\verticallineRightoneTHICK}

\multiput(10,10)(0,39){2}{\line(0,1){31}}
\multiput(30,10)(40,0){2}{\whitecircleNorthd}
\multiput(30,10)(40,0){2}{\holizontallineAbovetwo}
\multiput(40,10)(40,0){2}{\whitecircleNorthd}
\multiput(40,10)(40,0){2}{\EastUplineabovethreeTHICK}

\put(81,0){\line(1,0){9}}
\multiput(100,10)(0,39){2}{\line(0,1){31}}

\put(90,10){\curveSouthEast}

\put(90,80){\curveNorthEast}

\multiput(20,21)(0,28){2}{\line(0,1){20}}
\put(20,20){\whitecircleWestc}
\multiput(20,20)(40,0){2}{\EastDownlineabovethreeTHICK}
\multiput(20,20)(40,0){2}{\EastUplineabovetwoTHICK}
\multiput(50,20)(40,0){2}{\whitecircleNorthc}
\put(50,20){\holizontallineAboveone}
\put(60,20){\whitecircleNorthc}
\multiput(90,21)(0,28){2}{\line(0,1){20}}

\multiput(30,30)(40,0){2}{\whitecircleSouthb}
\multiput(30,30)(40,0){2}{\holizontallineAbovethree}
\multiput(30,30)(40,0){2}{\verticallineRightoneTHICK}
\multiput(40,30)(40,0){2}{\whitecircleSouthb}
\multiput(40,30)(40,0){2}{\EastDownlineabovetwoTHICK}
\multiput(40,30)(40,0){2}{\verticallineRightoneTHICK}

\put(10,40){\verticallineRightthree}
\put(20,40){\verticallineRightone}
\put(30,40){\verticallineRighttwoTHICK}
\put(30,40){\whitecircleWesta}
\multiput(30,40)(40,0){2}{\holizontallineAbovethree}
\put(40,40){\whitecircleNortha}
\multiput(41,40)(18,0){2}{\holizontallinetenmmTHICK}
\put(50,40){\holizontallineAbovetwoTHICK}
\put(70,40){\whitecircleNortha}
\put(80,40){\whitecircleEasta}
\put(80,40){\verticallineRighttwoTHICK}
\put(90,40){\verticallineRightone}
\put(100,40){\verticallineRightthree}

\put(30,50){\whitecircleWesta}
\put(30,50){\verticallineRightoneTHICK}
\put(30,50){\holizontallineAbovethree}
\put(40,50){\whitecircleSoutha}
\put(40,50){\verticallineRightoneTHICK}
\put(50,50){\holizontallineAbovetwoTHICK}
\put(70,50){\whitecircleSoutha}
\put(70,50){\verticallineRightoneTHICK}
\put(70,50){\holizontallineAbovethree}
\put(80,50){\whitecircleEasta}
\put(80,50){\verticallineRightoneTHICK}
\multiput(41,50)(18,0){2}{\holizontallinetenmmTHICK}

\multiput(30,60)(40,0){2}{\whitecircleNorthb}
\multiput(30,60)(40,0){2}{\holizontallineAbovethree}
\multiput(40,60)(40,0){2}{\whitecircleNorthb}
\multiput(40,60)(40,0){2}{\EastUplineabovetwoTHICK}

\put(20,70){\whitecircleWestc}
\multiput(20,70)(40,0){2}{\EastDownlineabovetwoTHICK}
\multiput(20,70)(40,0){2}{\EastUplineabovethreeTHICK}
\put(50,70){\whitecircleNorthc}
\put(50,70){\holizontallineAboveone}
\put(60,70){\whitecircleNorthc}
\put(90,70){\whitecircleEastc}

\multiput(30,80)(40,0){2}{\whitecircleSouthd}
\multiput(30,80)(40,0){2}{\verticallineRightoneTHICK}
\multiput(30,80)(40,0){2}{\holizontallineAbovetwo}
\multiput(40,80)(40,0){2}{\whitecircleSouthd}
\multiput(40,80)(40,0){2}{\verticallineRightoneTHICK}
\multiput(40,80)(40,0){2}{\EastDownlineabovethreeTHICK}

\put(20,90){\line(1,0){9}}
\put(50,90){\holizontallineAbovethree}
\multiput(30,90)(40,0){2}{\whitecircleNorthe}
\multiput(30,90)(40,0){2}{\holizontallineAbovetwoTHICK}
\multiput(40,90)(40,0){2}{\whitecircleNorthe}
\put(81,90){\line(1,0){10}}
\multiput(41,90)(18,0){2}{\holizontallinetenmm}

\end{picture}
}

\newcommand{\HecSevena}{
\setlength{\unitlength}{1mm}
\begin{picture}(40,20)(5,10)
\put(8,10){$a:=[\acchi^{(2,7)}_1]$}
\put(10,0){\circle{2}}\put(25,0){\circle{2}}\put(40,0){\circle{2}}
\put(11,0){\line(1,0){13}}\put(26,0){\line(1,0){13}}
\put(8,3){{\small{$-1$}}}\put(24,3){{\small{$q$}}}\put(39,3){{\small{$q^3$}}}
\put(15.5,2.0){{\small{$q^{-1}$}}}\put(30.5,2.0){{\small{$q^{-3}$}}}
\put(8,-4){$\al_1$}\put(23,-4){$\al_2$}\put(38,-4){$\al_3$}
\end{picture}
}

\newcommand{\HecSevenb}{
\setlength{\unitlength}{1mm}
\begin{picture}(40,20)(5,10)
\put(8,10){$b:=[\acchi^{(2,7)}_3]$}
\put(10,0){\circle{2}}\put(25,0){\circle{2}}\put(40,0){\circle{2}}
\put(11,0){\line(1,0){13}}\put(26,0){\line(1,0){13}}
\put(8,3){{\small{$-1$}}}\put(23,3){{\small{$-1$}}}\put(39,3){{\small{$q^3$}}}
\put(17,2.0){{\small{$q$}}}\put(30.5,2.0){{\small{$q^{-3}$}}}
\put(8,-4){$\al_1$}\put(23,-4){$\al_2$}\put(38,-4){$\al_3$}
\end{picture}
}

\newcommand{\HecSevenc}{
\setlength{\unitlength}{1mm}
\begin{picture}(40,20)(5,10)
\put(8,15){$c:=[\acchi^{(2,7)}_3]\circ\left[{{123}\atop{132}}\right]$}
\put(10,0){\circle{2}}\put(25,0){\circle{2}}\put(40,0){\circle{2}}
\put(11,0){\line(1,0){13}}\put(26,0){\line(1,0){13}}
\put(7,3){{\small{$q$}}}\put(23,3){{\small{$-1$}}}\put(41,3){{\small{$-1$}}}
\put(16,2.0){{\small{$q^{-1}$}}}\put(31,2.0){{\small{$q^3$}}}
\put(8,-4){$\al_1$}\put(23,-4){$\al_2$}\put(38,-4){$\al_3$}
\put(25,1){\oval(30,14)[t]}\put(24,10){{\small{$q^{-2}$}}}
\end{picture}
}

\newcommand{\HecSevend}{
\setlength{\unitlength}{1mm}
\begin{picture}(40,20)(5,10)
\put(8,10){$d:=[\acchi^{(2,7)}_3]$}
\put(10,0){\circle{2}}\put(25,0){\circle{2}}\put(40,0){\circle{2}}
\put(11,0){\line(1,0){13}}\put(26,0){\line(1,0){13}}
\put(9,3){{\small{$q^3$}}}\put(23,3){{\small{$-1$}}}\put(38,3){{\small{$-q^{-1}$}}}
\put(15.5,2.0){{\small{$q^{-3}$}}}\put(31,2.0){{\small{$q^2$}}}
\put(8,-4){$\al_2$}\put(23,-4){$\al_3$}\put(38,-4){$\al_1$}
\end{picture}
}

\newcommand{\HecSevenabcd}
{\setlength{\unitlength}{1mm}
\begin{picture}(130,40)(0,-15)
\put(0,20){\HecSevena}
\put(45,20){\HecSevenb}
\put(90,20){\HecSevenc}

\put(0,0){\HecSevend}

\put(60,-10){$(q\in\bKt, q^2\ne 1\ne q^3)$}

\end{picture}
}

\newcommand{\CayleySeven}{
\setlength{\unitlength}{1mm}
\begin{picture}(110,80)(-10,0)

\multiput(10,0)(2,0){6}{\line(1,0){0.5}}\put(13,1){$D_i$}
\put(20,0){\holizontallineAbovetwo}
\put(30,0){\whitecircleSoutha}
\put(30,0){\holizontallineAbovethreeTHICK}
\put(30,0){\verticallineRightoneTHICK}
\put(40,0){\whitecircleSoutha}
\put(40,0){\holizontallineAbovetwo}
\put(40,0){\verticallineRightoneTHICK}
\multiput(50,0)(2,0){6}{\line(1,0){0.5}}\put(51,1){$D_{i+1}$}

\put(30,10){\whitecircleNorthb}
\put(30,10){\holizontallineAbovethree}
\put(40,10){\whitecircleNorthb}
\put(40,10){\EastUplineabovetwoTHICK}

\multiput(0,20)(2,0){6}{\line(1,0){0.5}}\put(2,21){$C_i$}
\put(10,20){\holizontallineAboveone}
\put(20,20){\whitecircleSouthc}
\put(20,20){\EastDownlineabovetwoTHICK}
\put(20,20){\EastUplineabovethreeTHICK}
\put(50,20){\whitecircleSouthc}
\put(50,20){\holizontallineAboveone}
\multiput(60,20)(2,0){6}{\line(1,0){0.5}}\put(61,21){$C_{i+1}$}

\put(30,30){\whitecircleSouthd}
\put(30,30){\holizontallineAbovetwo}
\put(30,30){\verticallineRightoneTHICK}
\put(40,30){\whitecircleSouthd}
\put(40,30){\EastDownlineabovethreeTHICK}
\put(40,30){\verticallineRightoneTHICK}

\put(30,40){\whitecircleNorthd}
\put(30,40){\holizontallineAbovetwo}
\put(40,40){\whitecircleNorthd}
\put(40,40){\EastUplineabovethreeTHICK}

\multiput(0,50)(2,0){6}{\line(1,0){0.5}}\put(2,51){$B_i$}
\put(10,50){\holizontallineAboveone}
\put(20,50){\whitecircleSouthc}
\put(20,50){\EastDownlineabovethreeTHICK}
\put(20,50){\EastUplineabovetwoTHICK}
\put(50,50){\whitecircleSouthc}
\put(50,50){\holizontallineAboveone}
\multiput(60,50)(2,0){6}{\line(1,0){0.5}}\put(61,51){$B_{i+1}$}

\put(30,60){\whitecircleSouthb}
\put(30,60){\holizontallineAbovethree}
\put(30,60){\verticallineRightoneTHICK}
\put(40,60){\whitecircleSouthb}
\put(40,60){\EastDownlineabovetwoTHICK}
\put(40,60){\verticallineRightoneTHICK}

\multiput(10,70)(2,0){6}{\line(1,0){0.5}}\put(13,71){$A_i$}
\put(20,70){\holizontallineAbovetwoTHICK}
\put(30,70){\whitecircleNortha}
\put(30,70){\holizontallineAbovethree}
\put(40,70){\whitecircleNortha}
\put(40,70){\holizontallineAbovetwoTHICK}
\multiput(50,70)(2,0){6}{\line(1,0){0.5}}\put(51,71){$A_{i+1}$}

\end{picture}
}

\newcommand{\HecNinea}{
\setlength{\unitlength}{1mm}
\begin{picture}(40,20)(5,10)
\put(8,17){$a$}
\put(10,0){\circle{2}}\put(25,0){\circle{2}}\put(40,0){\circle{2}}
\put(11,0){\line(1,0){13}}\put(26,0){\line(1,0){13}}
\put(4,3){{\small{$-1$}}}\put(23,3){{\small{$-1$}}}\put(41,3){{\small{$-1$}}}
\put(17,2.0){{\small{$q$}}}\put(32,2.0){{\small{$r$}}}
\put(8,-4){$\al_1$}\put(23,-4){$\al_2$}\put(38,-4){$\al_3$}
\put(25,1){\oval(30,14)[t]}\put(22,10){{\small{$(qr)^{-1}$}}}
\end{picture}
}

\newcommand{\HecNineb}{
\setlength{\unitlength}{1mm}
\begin{picture}(40,20)(5,10)
\put(8,10){$b$}
\put(10,0){\circle{2}}\put(25,0){\circle{2}}\put(40,0){\circle{2}}
\put(11,0){\line(1,0){13}}\put(26,0){\line(1,0){13}}
\put(9,3){{\small{$q$}}}\put(23,3){{\small{$-1$}}}\put(39,3){{\small{$r$}}}
\put(15.5,1.5){{\small{$q^{-1}$}}}\put(30.5,1.5){{\small{$r^{-1}$}}}
\put(8,-6){$\al_1$}\put(23,-6){$\al_2$}\put(38,-6){$\al_3$}
\end{picture}
}

\newcommand{\HecNinec}{
\setlength{\unitlength}{1mm}
\begin{picture}(40,20)(5,10)
\put(8,10){$c$}
\put(10,0){\circle{2}}\put(25,0){\circle{2}}\put(40,0){\circle{2}}
\put(11,0){\line(1,0){13}}\put(26,0){\line(1,0){13}}
\put(9,3){{\small{$q$}}}\put(23,3){{\small{$-1$}}}\put(37,3){{\small{$(qr)^{-1}$}}}
\put(15.5,1.5){{\small{$q^{-1}$}}}\put(31,1.5){{\small{$qr$}}}

\put(8,-6){$\al_2$}\put(23,-6){$\al_1$}\put(38,-6){$\al_3$}
\end{picture}
}

\newcommand{\HecNined}{
\setlength{\unitlength}{1mm}
\begin{picture}(40,20)(5,10)
\put(8,10){$d$}
\put(10,0){\circle{2}}\put(25,0){\circle{2}}\put(40,0){\circle{2}}
\put(11,0){\line(1,0){13}}\put(26,0){\line(1,0){13}}
\put(7,3){{\small{$(qr)^{-1}$}}}\put(23,3){{\small{$-1$}}}\put(39,3){{\small{$r$}}}
\put(16,1.5){{\small{$qr$}}}\put(30.5,1.5){{\small{$r^{-1}$}}}

\put(8,-6){$\al_1$}\put(23,-6){$\al_3$}\put(38,-6){$\al_2$}
\end{picture}
}

\newcommand{\HecNineabcd}{
\setlength{\unitlength}{1mm}
\begin{picture}(130,60)(0,-20)
\put(0,25){$\HecNinea$}
\put(0,0){$\HecNineb$}
\put(45,0){$\HecNinec$}
\put(90,0){$\HecNined$}
\put(45,20){${\tiny{
\left(
\begingroup
\renewcommand{\arraystretch}{1.5}
\begin{array}{l} 
q,r\in\bKt, q\ne 1,  r\ne 1, qr\ne 1,\\
\mbox{If $q\ne r\ne (qr)^{-1}\ne q$, then
$a:=[\acchi^{(2,9)}_2]$, $b:=[\acchi^{(2,9)}_1]$,} \\
\mbox{$c:=[\acchi^{(2,9)}_3]\circ\left[{{123}\atop{213}}\right]$ and
$d:=[\acchi^{(2,9)}_4]\circ\left[{{123}\atop{312}}\right]$.} \\
\mbox{If $q=r\ne (qr)^{-1}$, then
$a:=[\acchi^{(2,10)}_2]$, $b:=[\acchi^{(2,10)}_1]$,} \\
\mbox{$c:=[\acchi^{(2,10)}_3]\circ\left[{{123}\atop{213}}\right]$ and
$d:=c\circ\left[{{123}\atop{312}}\right]$.} \\
\mbox{If $q=r= (qr)^{-1}$, then
$a:=[\acchi^{(2,11)}_2]$, $b:=[\acchi^{(2,11)}_1]$,} \\
\mbox{$c:=b\circ\left[{{123}\atop{213}}\right]$ and
$d:=c\circ\left[{{123}\atop{312}}\right]$.} 
\end{array}
\endgroup
\right)}}$} 
\end{picture}
}

\newcommand{\CayleyNine}{
\setlength{\unitlength}{1mm}
\begin{picture}(60,85)(25,0)
\put(50,10){\whitecircleSouthc}
\put(50,10){\holizontallineAbovetwoTHICK}
\put(60,10){\whitecircleSouthc}
\put(50,10){\verticallineRightthreeTHICK}
\put(60,10){\verticallineRightthreeTHICK}

\put(50,20){\whitecircleNorthc}
\put(50,20){\holizontallineAbovetwo}
\put(60,20){\whitecircleNorthc}
\put(60,20){\EastUplineaboveoneTHICK}

\put(10,30){\whitecircleNortha}
\put(10,30){\EastUplineabovethreeTHICK}
\put(40,30){\whitecircleSoutha}
\put(40,30){\EastUplineabovetwo}
\put(40,30){\EastDownlineaboveoneTHICK}
\put(70,30){\whitecircleSoutha}
\put(70,30){\EastUplineabovethreeTHICK}
\put(100,30){\whitecircleNortha}
\put(100,30){\EastUplineabovetwoTHICK}

\put(0,40){\whitecircleWestb}
\put(0,40){\verticallineRightthreeTHICK}
\put(0,40){\EastDownlineabovetwoTHICK}
\put(20,40){\whitecircleSouthd}
\put(20,40){\verticallineRighttwoTHICK}
\put(20,40){\holizontallineAboveone}
\put(30,40){\whitecircleSouthd}
\put(30,40){\verticallineRighttwoTHICK}
\put(30,40){\EastDownlineabovethreeTHICK}
\put(50,40){\whitecircleSouthb}
\put(50,40){\verticallineRightthreeTHICK}
\put(50,40){\holizontallineAboveoneTHICK}
\put(60,40){\whitecircleSouthb}
\put(60,40){\verticallineRightthreeTHICK}
\put(60,40){\EastDownlineaboveone}
\put(80,40){\whitecircleSouthd}
\put(80,40){\verticallineRighttwo}
\put(80,40){\holizontallineAboveoneTHICK}
\put(90,40){\whitecircleSouthd}
\put(90,40){\verticallineRighttwo}
\put(90,40){\EastDownlineabovethreeTHICK}
\put(110,40){\whitecircleEastb}
\put(110,40){\verticallineRightthreeTHICK}

\put(0,50){\whitecircleWestb}
\put(0,50){\EastUplineabovetwo}
\put(20,50){\whitecircleNorthd}
\put(20,50){\holizontallineAboveone}
\put(30,50){\whitecircleNorthd}
\put(30,50){\EastUplineabovethreeTHICK}
\put(50,50){\whitecircleNorthb}
\put(50,50){\holizontallineAboveone}
\put(60,50){\whitecircleNorthb}
\put(60,50){\EastUplineabovetwoTHICK}
\put(80,50){\whitecircleNorthd}
\put(80,50){\holizontallineAboveoneTHICK}
\put(90,50){\whitecircleNorthd}
\put(90,50){\EastUplineabovethreeTHICK}
\put(110,50){\whitecircleEastb}

\put(10,60){\whitecircleSoutha}
\put(10,60){\EastDownlineabovethreeTHICK}
\put(40,60){\whitecircleNortha}
\put(40,60){\EastUplineaboveone}
\put(40,60){\EastDownlineabovetwoTHICK}
\put(70,60){\whitecircleSoutha}
\put(70,60){\EastDownlineabovethreeTHICK}
\put(100,60){\whitecircleSoutha}
\put(100,60){\EastDownlineabovetwo}

\put(50,70){\whitecircleSouthc}
\put(50,70){\verticallineRightthreeTHICK}
\put(50,70){\holizontallineAbovetwoTHICK}
\put(60,70){\whitecircleSouthc}
\put(60,70){\verticallineRightthreeTHICK}
\put(60,70){\EastDownlineaboveone}

\put(50,80){\whitecircleNorthc}
\put(50,80){\holizontallineAbovetwo}
\put(60,80){\whitecircleNorthc}

\put(10,0){\line(1,0){90}}

\put(10,10){\curveSouthWest}
\put(0,10){\line(0,1){29}}
\put(100,10){\curveSouthEast}
\put(110,10){\line(0,1){29}}
\put(20,10){\line(1,0){29}}
\put(61,10){\line(1,0){29}}

\put(20,20){\curveSouthWest}
\put(10,20){\line(0,1){9}}
\put(90,20){\curveSouthEast}
\put(100,20){\line(0,1){9}}

\put(0,80){\line(0,-1){29}}\multiput(0,80)(0.05,0){4}{\line(0,-1){29}}\multiput(0,80)(-0.05,0){4}{\line(0,-1){29}}
\put(110,51){\line(0,1){29}}\multiput(110,51)(0.05,0){4}{\line(0,1){29}}\multiput(110,51)(-0.05,0){4}{\line(0,1){29}}

\put(10,61){\line(0,1){9}}\multiput(10,61)(0.05,0){4}{\line(0,1){9}}\multiput(10,61)(-0.05,0){4}{\line(0,1){9}}
\put(100,70){\line(0,-1){9}}\multiput(100,70)(0.05,0){4}{\line(0,-1){9}}\multiput(100,70)(-0.05,0){4}{\line(0,-1){9}}

\put(20,70){\curveNorthWestTHICK}
\put(90,70){\curveNorthEastTHICK}

\put(10,80){\curveNorthWestTHICK}
\put(20,80){\line(1,0){29}}\multiput(20,80)(0,0.05){4}{\line(1,0){29}}\multiput(20,80)(0,-0.05){4}{\line(1,0){29}}
\put(61,80){\line(1,0){29}}\multiput(61,80)(0,0,05){4}{\line(1,0){29}}\multiput(61,80)(0,-0,05){4}{\line(1,0){29}}
\put(100,80){\curveNorthEastTHICK}

\put(10,90){\line(1,0){90}}\multiput(10,90)(0,0,05){4}{\line(1,0){90}}\multiput(10,90)(0,-0,05){4}{\line(1,0){90}}

\put(1,70){{\tiny{$1$}}}
\put(11,70){{\tiny{$1$}}}
\put(101,70){{\tiny{$1$}}}

\put(1,15){{\tiny{$1$}}}
\put(11,15){{\tiny{$1$}}}
\put(101,15){{\tiny{$1$}}}

\end{picture}
}

\newcommand{\HecThirteena}{
\setlength{\unitlength}{1mm}
\begin{picture}(40,20)(5,10)
\put(8,10){$a:=[\acchi^{(3,13)}_1]$}
\put(10,0){\circle{2}}\put(25,0){\circle{2}}\put(40,0){\circle{2}}
\put(11,0){\line(1,0){13}}\put(26,0){\line(1,0){13}}
\put(9,3){{\small{$\defzeta$}}}\put(24,3){{\small{$\defzeta$}}}\put(38,3){{\small{$-1$}}}
\put(15.5,2.0){{\small{$\defzeta^{-1}$}}}\put(30.5,2.0){{\small{$\defzeta^{-2}$}}}
\put(8,-4){$\al_1$}\put(23,-4){$\al_2$}\put(38,-4){$\al_3$}
\end{picture}
}

\newcommand{\HecThirteenb}{
\setlength{\unitlength}{1mm}
\begin{picture}(40,20)(5,10)
\put(8,10){$b:=[\acchi^{(3,13)}_2]$}
\put(10,0){\circle{2}}\put(25,0){\circle{2}}\put(40,0){\circle{2}}
\put(11,0){\line(1,0){13}}\put(26,0){\line(1,0){13}}
\put(9,3){{\small{$\defzeta$}}}\put(22.5,3){{\small{$-\defzeta^{-1}$}}}\put(38,3){{\small{$-1$}}}
\put(15.5,2.0){{\small{$\defzeta^{-1}$}}}\put(32,2.0){{\small{$\defzeta^2$}}}
\put(8,-4){$\al_1$}\put(23,-4){$\al_2$}\put(38,-4){$\al_3$}
\end{picture}
}

\newcommand{\HecThirteenab}
{\setlength{\unitlength}{1mm}
\begin{picture}(70,20)(15,-15)

\put(0,0){\HecThirteena}
\put(45,0){\HecThirteenb}

\put(90,-10){$(\zeta\in\bKt,\,\zeta^4+\zeta^2+1=0)$}
\end{picture}
}

\newcommand{\HecFourteena}{
\setlength{\unitlength}{1mm}
\begin{picture}(40,20)(5,10)
\put(8,10){$a:=[\acchi^{(2,14)}_1]$} 
\put(10,0){\circle{2}}\put(25,0){\circle{2}}\put(40,0){\circle{2}}
\put(11,0){\line(1,0){13}}\put(26,0){\line(1,0){13}}
\put(8,3){{\small{$-1$}}}\put(21.5,3){{\small{$-\defzeta^{-1}$}}}\put(39,3){{\small{$\defzeta$}}}
\put(15,2.0){{\small{$-\defzeta$}}}\put(30,2.0){{\small{$-\defzeta$}}}
\put(8,-4){$\al_1$}\put(23,-4){$\al_2$}\put(38,-4){$\al_3$}
\end{picture}
}

\newcommand{\HecFourteenb}{
\setlength{\unitlength}{1mm}
\begin{picture}(40,20)(5,10)
\put(8,10){$b:=[\acchi^{(2,14)}_2]$} 
\put(10,0){\circle{2}}\put(25,0){\circle{2}}\put(40,0){\circle{2}}
\put(11,0){\line(1,0){13}}\put(26,0){\line(1,0){13}}
\put(8,3){{\small{$-1$}}}\put(23,3){{\small{$-1$}}}\put(39,3){{\small{$\defzeta$}}}
\put(14,2.0){{\small{$-\defzeta^{-1}$}}}\put(30,2.0){{\small{$-\defzeta$}}}
\put(8,-4){$\al_1$}\put(23,-4){$\al_2$}\put(38,-4){$\al_3$}
\end{picture}
}

\newcommand{\HecFourteenc}{
\setlength{\unitlength}{1mm}
\begin{picture}(40,20)(5,10)
\put(8,10){$c:=[\acchi^{(2,14)}_3]$} 
\put(10,0){\circle{2}}\put(25,0){\circle{2}}\put(40,0){\circle{2}}
\put(11,0){\line(1,0){13}}\put(26,0){\line(1,0){13}}
\put(5.5,3){{\small{$-\defzeta^{-1}$}}}\put(23,3){{\small{$-1$}}}\put(39,3){{\small{$\defzeta$}}}
\put(15,2.0){{\small{$-\defzeta$}}}\put(29,2.0){{\small{$-\defzeta^{-1}$}}}
\put(8,-4){$\al_1$}\put(23,-4){$\al_2$}\put(38,-4){$\al_3$}
\end{picture}
}

\newcommand{\HecFourteenabc}
{\setlength{\unitlength}{1mm}
\begin{picture}(130,25)(0,-20)

\put(0,0){\HecFourteena}
\put(45,0){\HecFourteenb}
\put(90,0){\HecFourteenc}

\put(8,-20){\footnotesize{$(\zeta\in\bKt$, $\zeta^2+\zeta+1=0)$}}
\end{picture}
}

\newcommand{\CayleyFourteen}{
\setlength{\unitlength}{1mm}
\begin{picture}(110,50)(-10,0)

\multiput(10,0)(2,0){6}{\line(1,0){0.5}}\put(13,1){$D_i$}
\put(20,0){\holizontallineAbovetwo}
\put(30,0){\whitecircleSoutha}
\put(30,0){\holizontallineAbovethreeTHICK}
\put(30,0){\verticallineRightoneTHICK}
\put(40,0){\whitecircleSoutha}
\put(40,0){\holizontallineAbovetwo}
\put(40,0){\verticallineRightoneTHICK}
\multiput(50,0)(2,0){6}{\line(1,0){0.5}}\put(51,1){$D_{i+1}$}

\put(30,10){\whitecircleNorthb}
\put(30,10){\holizontallineAbovethree}
\put(40,10){\whitecircleNorthb}
\put(40,10){\EastUplineabovetwoTHICK}

\multiput(0,20)(2,0){6}{\line(1,0){0.5}}\put(2,21){$C_i$}
\put(10,20){\holizontallineAboveone}
\put(20,20){\whitecircleSouthc}
\put(20,20){\EastDownlineabovetwoTHICK}
\put(20,20){\verticallineRightthreeTHICK}
\put(50,20){\whitecircleSouthc}
\put(50,20){\holizontallineAboveone}
\put(50,20){\verticallineRightthreeTHICK}
\multiput(60,20)(2,0){6}{\line(1,0){0.5}}\put(61,21){$C_{i+1}$}

\multiput(0,30)(2,0){6}{\line(1,0){0.5}}\put(2,31){$B_i$}
\put(10,30){\holizontallineAboveone}
\put(20,30){\whitecircleNorthc}
\put(20,30){\EastUplineabovetwoTHICK}
\put(50,30){\whitecircleNorthc}
\put(50,30){\holizontallineAboveone}
\multiput(60,30)(2,0){6}{\line(1,0){0.5}}\put(61,31){$B_{i+1}$}

\put(30,40){\whitecircleSouthb}
\put(30,40){\holizontallineAbovethree}
\put(30,40){\verticallineRightoneTHICK}
\put(40,40){\whitecircleSouthb}
\put(40,40){\EastDownlineabovetwoTHICK}
\put(40,40){\verticallineRightoneTHICK}

\multiput(10,50)(2,0){6}{\line(1,0){0.5}}\put(13,51){$A_i$}
\put(20,50){\holizontallineAbovetwoTHICK}
\put(30,50){\whitecircleNortha}
\put(30,50){\holizontallineAbovethree}
\put(40,50){\whitecircleNortha}
\put(40,50){\holizontallineAbovetwoTHICK}
\multiput(50,50)(2,0){6}{\line(1,0){0.5}}\put(51,51){$A_{i+1}$}

\end{picture}
}

\newcommand{\HecFifteena}{
\setlength{\unitlength}{1mm}
\begin{picture}(40,20)(5,10)
\put(8,17){$a:=[\acchi^{(2,15)}_2]$}
\put(10,0){\circle{2}}\put(25,0){\circle{2}}\put(40,0){\circle{2}}
\put(11,0){\line(1,0){13}}\put(26,0){\line(1,0){13}}
\put(6.5,3){{\small{$\defzeta$}}}\put(23,3){{\small{$-1$}}}\put(41,3){{\small{$\defzeta$}}}
\put(15.5,2.0){{\small{$\defzeta^{-1}$}}}\put(31,2.0){{\small{$\defzeta^{-1}$}}}
\put(8,-4){$\al_1$}\put(23,-4){$\al_2$}\put(38,-4){$\al_3$}
\put(25,1){\oval(30,14)[t]}\put(24,10){{\small{$\defzeta^{-1}$}}}
\end{picture}
}

\newcommand{\HecFifteenb}{
\setlength{\unitlength}{1mm}
\begin{picture}(40,20)(5,10)
\put(8,10){$b:=[\acchi^{(2,15)}_3]$}
\put(10,0){\circle{2}}\put(25,0){\circle{2}}\put(40,0){\circle{2}}
\put(11,0){\line(1,0){13}}\put(26,0){\line(1,0){13}}
\put(8,3){{\small{$-1$}}}\put(23,3){{\small{$-1$}}}\put(38,3){{\small{$-1$}}}
\put(17,2.0){{\small{$\defzeta$}}}\put(32,2.0){{\small{$\defzeta$}}}
\put(8,-4){$\al_1$}\put(23,-4){$\al_2$}\put(38,-4){$\al_3$}
\end{picture}
}

\newcommand{\HecFifteenc}{
\setlength{\unitlength}{1mm}
\begin{picture}(40,20)(5,10)
\put(8,10){$c:=[\acchi^{(2,15)}_1]$}
\put(10,0){\circle{2}}\put(25,0){\circle{2}}\put(40,0){\circle{2}}
\put(11,0){\line(1,0){13}}\put(26,0){\line(1,0){13}}
\put(8,3){{\small{$-1$}}}\put(24,3){{\small{$\defzeta$}}}\put(38,3){{\small{$-1$}}}
\put(15.5,2.0){{\small{$\defzeta^{-1}$}}}\put(32,2.0){{\small{$\defzeta$}}}
\put(8,-4){$\al_1$}\put(23,-4){$\al_2$}\put(38,-4){$\al_3$}
\end{picture}
}

\newcommand{\HecFifteend}{
\setlength{\unitlength}{1mm}
\begin{picture}(40,20)(5,10)
\put(8,10){$d:=c\circ\left[{{123}\atop{321}}\right]$}
\put(10,0){\circle{2}}\put(25,0){\circle{2}}\put(40,0){\circle{2}}
\put(11,0){\line(1,0){13}}\put(26,0){\line(1,0){13}}
\put(8,3){{\small{$-1$}}}\put(24,3){{\small{$\defzeta$}}}\put(38,3){{\small{$-1$}}}
\put(17,2.0){{\small{$\defzeta$}}}\put(30.5,2.0){{\small{$\defzeta^{-1}$}}}
\put(8,-4){$\al_1$}\put(23,-4){$\al_2$}\put(38,-4){$\al_3$}
\end{picture}
}

\newcommand{\HecFifteene}{
\setlength{\unitlength}{1mm}
\begin{picture}(40,20)(5,10)
\put(8,10){$e:=[\acchi^{(2,15)}_4]$}
\put(10,0){\circle{2}}\put(25,0){\circle{2}}\put(40,0){\circle{2}}
\put(11,0){\line(1,0){13}}\put(26,0){\line(1,0){13}}
\put(8,3){{\small{$-1$}}}\put(22,3){{\small{$-\defzeta^{-1}$}}}\put(38,3){{\small{$-1$}}}
\put(15.5,2.0){{\small{$\defzeta^{-1}$}}}\put(30.5,2.0){{\small{$\defzeta^{-1}$}}}
\put(8,-4){$\al_1$}\put(23,-4){$\al_2$}\put(38,-4){$\al_3$}
\end{picture}
}

\newcommand{\HecFifteenabcd}
{\setlength{\unitlength}{1mm}
\begin{picture}(130,60)(10,-20)
\put(0,25){\HecFifteena}
\put(45,25){\HecFifteenb}
\put(90,25){\HecFifteenc}

\put(0,0){\HecFifteend}
\put(45,0){\HecFifteene}
\put(90,-10){$(\zeta\in\bKt,\,\zeta^2+\zeta+1=0)$}
\end{picture}
}

\newcommand{\CayleyFifteen}{
\setlength{\unitlength}{1mm}
\begin{picture}(80,80)(-20,0)

\put(0,0){\holizontallineAbovetwo}
\put(10,0){\whitecircleSouthe}
\put(10,0){\holizontallineAbovethreeTHICK}
\put(10,0){\verticallineRightoneTHICK}
\put(20,0){\whitecircleSouthc}
\put(20,0){\verticallineRightoneTHICK}
\put(20,0){\holizontallineAbovetwo}

\put(10,10){\whitecircleNortha}
\put(10,10){\holizontallineAbovethree}
\put(20,10){\whitecircleNorthb}
\put(20,10){\EastUplineabovetwoTHICK}

\put(0,20){\whitecircleSouthd}
\put(0,20){\EastDownlineabovetwoTHICK}
\put(0,20){\EastUplineabovethreeTHICK}
\put(30,20){\whitecircleSoutha}
\put(30,20){\holizontallineAboveoneTHICK}

\put(-20,30){\holizontallineAbovetwoTHICK}
\put(-10,30){\whitecircleEastd}
\put(-10,30){\EastUplineabovethreeTHICK}
\put(-10,30){\EastDownlineaboveone}
\put(10,30){\whitecircleSouthb}
\put(10,30){\holizontallineAbovetwoTHICK}
\put(20,30){\whitecircleSoutha}
\put(20,30){\EastDownlineabovethree}
\put(20,30){\EastUplineaboveoneTHICK}

\put(0,40){\whitecircleNorthc}
\put(0,40){\EastDownlineaboveone}
\put(0,40){\EastUplineabovetwoTHICK}
\put(30,40){\whitecircleNortha}
\put(30,40){\holizontallineAbovethree}

\put(10,50){\whitecircleSouthc}
\put(10,50){\holizontallineAboveone}
\put(10,50){\verticallineRightthreeTHICK}
\put(20,50){\whitecircleSouthb}
\put(20,50){\verticallineRightthreeTHICK}
\put(20,50){\EastDownlineabovetwoTHICK}

\put(0,60){\holizontallineAbovethree}
\put(10,60){\whitecircleNorthe}
\put(10,60){\holizontallineAboveoneTHICK}
\put(20,60){\whitecircleNorthd}
\put(20,60){\holizontallineAbovethree}

\multiput(-10,60)(2,0){6}{\line(1,0){0.5}}\put(-6.5,60.5){{\tiny{A(i)}}}
\multiput(-30,30)(2,0){6}{\line(1,0){0.5}}\put(-26.5,30.5){{\tiny{B(i)}}}
\multiput(-10,0)(2,0){6}{\line(1,0){0.5}}\put(-6.5,0.5){{\tiny{C(i)}}}

\multiput(29.5,60)(2,0){6}{\line(1,0){0.5}}\put(33.5,60.5){{\tiny{D(i)}}}
\multiput(39.5,40)(2,0){6}{\line(1,0){0.5}}\put(43.5,40.5){{\tiny{E(i)}}}
\multiput(39.5,20)(2,0){6}{\line(1,0){0.5}}\put(43.5,20.5){{\tiny{F(i)}}}
\multiput(29.5,00)(2,0){6}{\line(1,0){0.5}}\put(33.5,00.5){{\tiny{G(i)}}}

\end{picture}
}

\newcommand{\HecSixteenhati}{
\setlength{\unitlength}{1mm}
\begin{picture}(40,20)(5,10)
\put(8,10){$a:=[\acchi^{(2,16)}_1]$}
\put(10,0){\circle{2}}\put(25,0){\circle{2}}\put(40,0){\circle{2}}
\put(11,0){\line(1,0){13}}\put(26,0){\line(1,0){13}}
\put(8,3){{\small{$-1$}}}\put(24,3){{\small{$\defzeta$}}}\put(38,3){{\small{$-\defzeta$}}}
\put(15.5,2.0){{\small{$\defzeta^{-1}$}}}\put(27.5,2.0){{\small{$-\defzeta^{-1}$}}}
\put(8,-4){$\al_1$}\put(23,-4){$\al_2$}\put(38,-4){$\al_3$}
\end{picture}
}

\newcommand{\HecSixteenhatii}{
\setlength{\unitlength}{1mm}
\begin{picture}(40,20)(2,10)
\put(8,17){$b:=[\acchi^{(2,16)}_2]\circ\left[{{123}\atop{213}}\right]$}
\put(10,0){\circle{2}}\put(25,0){\circle{2}}\put(40,0){\circle{2}}
\put(11,0){\line(1,0){13}}\put(26,0){\line(1,0){13}}
\put(6.5,3){{\small{$\defzeta$}}}\put(23,3){{\small{$-1$}}}\put(41,3){{\small{$-1$}}}
\put(14.5,2.0){{\small{$-1$}}}\put(30,2.0){{\small{$-\defzeta$}}}
\put(8,-4){$\al_1$}\put(23,-4){$\al_3$}\put(38,-4){$\al_2$}
\put(25,1){\oval(30,14)[t]}\put(24,10){{\small{$\defzeta^{-1}$}}}
\end{picture}
}

\newcommand{\HecSixteenhatiii}{
\setlength{\unitlength}{1mm}
\begin{picture}(40,20)(5,10)
\put(8,10){$c:=[\acchi^{(2,16)}_3]$}
\put(10,0){\circle{2}}\put(25,0){\circle{2}}\put(40,0){\circle{2}}
\put(11,0){\line(1,0){13}}\put(26,0){\line(1,0){13}}
\put(8,3){{\small{$-1$}}}\put(23,3){{\small{$-1$}}}\put(38,3){{\small{$-\defzeta$}}}
\put(17,2.0){{\small{$\defzeta$}}}\put(28.5,2.0){{\small{$-\defzeta^{-1}$}}}
\put(8,-4){$\al_1$}\put(23,-4){$\al_2$}\put(38,-4){$\al_3$}
\end{picture}
}

\newcommand{\HecSixteenhativ}{
\setlength{\unitlength}{1mm}
\begin{picture}(40,20)(5,10)
\put(8,10){$f:=[\acchi^{(2,16)}_4]\circ\left[{{123}\atop{312}}\right]$}
\put(10,0){\circle{2}}\put(25,0){\circle{2}}\put(40,0){\circle{2}}
\put(11,0){\line(1,0){13}}\put(26,0){\line(1,0){13}}
\put(9,3){{\small{$\defzeta$}}}\put(23,3){{\small{$-1$}}}\put(38,3){{\small{$-\defzeta$}}}
\put(15,2.0){{\small{$-1$}}}\put(28.5,2.0){{\small{$-\defzeta^{-1}$}}}
\put(8,-4){$\al_1$}\put(23,-4){$\al_2$}\put(38,-4){$\al_3$}
\end{picture}
}

\newcommand{\HecSixteenhatv}{
\setlength{\unitlength}{1mm}
\begin{picture}(40,20)(5,10)
\put(8,10){$g:=[\acchi^{(2,16)}_5]\circ\left[{{123}\atop{312}}\right]$}
\put(10,0){\circle{2}}\put(25,0){\circle{2}}\put(40,0){\circle{2}}
\put(11,0){\line(1,0){13}}\put(26,0){\line(1,0){13}}
\put(9,3){{\small{$\defzeta$}}}\put(23,3){{\small{$-\defzeta$}}}\put(38,3){{\small{$-\defzeta$}}}
\put(13.5,2.0){{\small{$-\defzeta^{-1}$}}}\put(28.5,2.0){{\small{$-\defzeta^{-1}$}}}
\put(8,-4){$\al_2$}\put(23,-4){$\al_3$}\put(38,-4){$\al_1$}
\end{picture}
}

\newcommand{\HeckSixteenabcdefghij}
{\setlength{\unitlength}{1mm}
\begin{picture}(130,70)(10,-30)
\put(0,25){\HecSixteenhati}
\put(45,25){\HecSixteenhatii}
\put(90,25){\HecSixteenhatiii}

\put(0,0){\HecSixteenhativ}
\put(45,0){\HecSixteenhatv}
\put(90,-10){$(\zeta\in\bKt,\,\zeta^2+\zeta+1=0)$}

\put(5,-25){$d:=b\circ\left[{{123}\atop{213}}\right]$}
\put(30,-25){$e:=c\circ\left[{{123}\atop{213}}\right]$}
\put(55,-25){$h:=a\circ\left[{{123}\atop{213}}\right]$}
\put(80,-25){$i:=f\circ\left[{{123}\atop{213}}\right]$}
\put(105,-25){$j:=g\circ\left[{{123}\atop{213}}\right]$}
\end{picture}
}
\newcommand{\CayleySixteen}{
\setlength{\unitlength}{1mm}
\begin{picture}(110,80)(-10,0)

\put(70,0){\whitecircleSouthd}
\put(70,0){\holizontallineAboveoneTHICK}
\put(80,0){\whitecircleSouthh}
\put(80,0){\EastUplineabovetwoTHICK}

\put(20,10){\whitecircleSouthj}
\put(20,10){\verticallineRighttwoTHICK}
\put(20,10){\holizontallineAboveoneTHICK}
\put(30,10){\whitecircleEasti}
\put(30,10){\verticallineRighttwoTHICK}
\put(60,10){\whitecircleSouthc}
\put(60,10){\EastUplineaboveoneTHICK}
\put(60,10){\EastDownlineabovetwoTHICK}
\put(90,10){\whitecircleSouthe}
\put(90,10){\holizontallineAbovethree}

\put(20,20){\whitecircleNorthj}
\put(20,20){\holizontallineAboveone}
\put(30,20){\whitecircleNorthi}
\put(30,20){\EastUplineabovethreeTHICK}
\put(50,20){\whitecircleNorthc}
\put(50,20){\EastUplineaboveoneTHICK}
\put(50,20){\EastDownlineabovethree}
\put(70,20){\whitecircleSoutha}
\put(70,20){\holizontallineAbovetwo}
\put(80,20){\whitecircleSouthb}
\put(80,20){\EastDownlineaboveoneTHICK}
\put(80,20){\EastUplineabovethreeTHICK}

\put(0,30){\holizontallineAbovetwo}
\put(10,30){\whitecircleSouthj}
\put(10,30){\EastDownlineabovethreeTHICK}
\put(10,30){\verticallineRightoneTHICK}
\put(40,30){\whitecircleSouthd}
\put(40,30){\EastDownlineabovetwoTHICK}
\put(40,30){\verticallineRightone}
\put(60,30){\whitecircleSoutha}
\put(60,30){\EastDownlineabovethreeTHICK}
\put(60,30){\verticallineRighttwo}
\put(90,30){\whitecircleSouthf}
\put(90,30){\holizontallineAboveoneTHICK}
\put(90,30){\verticallineRighttwo}

\put(0,40){\holizontallineAbovetwo}
\put(10,40){\whitecircleNorthi}
\put(10,40){\EastUplineabovethreeTHICK}
\put(40,40){\whitecircleNorthh}
\put(40,40){\EastUplineabovetwoTHICK}
\put(60,40){\whitecircleNorthb}
\put(60,40){\EastUplineabovethreeTHICK}
\put(90,40){\whitecircleNorthg}
\put(90,40){\holizontallineAboveoneTHICK}

\put(20,50){\whitecircleNorthd}
\put(20,50){\holizontallineAboveone}
\put(30,50){\whitecircleNorthh}
\put(30,50){\EastDownlineabovethreeTHICK}
\put(30,50){\EastUplineabovetwoTHICK}
\put(50,50){\whitecircleNorthe}
\put(50,50){\EastDownlineaboveoneTHICK}
\put(70,50){\whitecircleSouthf}
\put(70,50){\holizontallineAbovetwo}
\put(70,50){\verticallineRightoneTHICK}
\put(80,50){\whitecircleSouthg}
\put(80,50){\EastDownlineabovethreeTHICK}
\put(80,50){\verticallineRightoneTHICK}

\put(0,60){\holizontallineAbovethree}
\put(10,60){\whitecircleNorthc}
\put(10,60){\EastDownlineabovetwoTHICK}
\put(10,60){\EastUplineaboveoneTHICK}
\put(40,60){\whitecircleNorthe}
\put(40,60){\EastDownlineabovethree}
\put(70,60){\whitecircleWestf}
\put(70,60){\holizontallineAbovetwoTHICK}
\put(80,60){\whitecircleNorthg}

\put(20,70){\whitecircleNortha}
\put(20,70){\holizontallineAbovetwoTHICK}
\put(30,70){\whitecircleNorthb}
\put(30,70){\EastDownlineaboveoneTHICK}

\put(40,10){\curveSouthWestNW}
\put(40,00){\line(1,0){29}}\put(44.5,0.5){{\tiny{3}}}

\put(60,60){\curveNorthEastSE}
\put(31,70){\line(1,0){29}}\put(44.5,70.5){{\tiny{3}}}

\put(1,10){\line(1,0){18}}\put(4.5,10.5){{\tiny{3}}}
\put(1,70){\line(1,0){18}}\put(4.5,70.5){{\tiny{3}}}

\put(81,00){\line(1,0){18}}\put(94.5,0.5){{\tiny{3}}}
\put(81,60){\line(1,0){18}}\put(94.5,60.5){{\tiny{3}}}

\multiput(-10,10)(2,0){6}{\line(1,0){0.5}}\put(-5.5,10.5){{\tiny{E}}}
\multiput(-10,30)(2,0){6}{\line(1,0){0.5}}\put(-5.5,30.5){{\tiny{D}}}
\multiput(-10,40)(2,0){6}{\line(1,0){0.5}}\put(-5.5,40.5){{\tiny{C}}}
\multiput(-10,60)(2,0){6}{\line(1,0){0.5}}\put(-5.5,60.5){{\tiny{B}}}
\multiput(-10,70)(2,0){6}{\line(1,0){0.5}}\put(-5.5,70.5){{\tiny{A}}}

\multiput(99.5,0)(2,0){6}{\line(1,0){0.5}}\put(104.5,0.5){{\tiny{J}}}
\multiput(99.5,10)(2,0){6}{\line(1,0){0.5}}\put(104.5,10.5){{\tiny{I}}}
\multiput(99.5,30)(2,0){6}{\line(1,0){0.5}}\put(104.5,30.5){{\tiny{H}}}
\multiput(99.5,40)(2,0){6}{\line(1,0){0.5}}\put(104.5,40.5){{\tiny{G}}}
\multiput(99.5,60)(2,0){6}{\line(1,0){0.5}}\put(104.5,60.5){{\tiny{F}}}

\end{picture}
}

\newcommand{\HecSeventeena}{
\setlength{\unitlength}{1mm}
\begin{picture}(40,20)(5,10)
\put(8,10){$a:=[\acchi^{(2,17)}_1]$}
\put(10,0){\circle{2}}\put(25,0){\circle{2}}\put(40,0){\circle{2}}
\put(11,0){\line(1,0){13}}\put(26,0){\line(1,0){13}}
\put(8,3){{\small{$-1$}}}\put(23,3){{\small{$-1$}}}\put(38,3){{\small{$-1$}}}
\put(15,2.0){{\small{$-1$}}}\put(32,2.0){{\small{$\defzeta$}}}
\put(8,-4){$\al_1$}\put(23,-4){$\al_2$}\put(38,-4){$\al_3$}
\end{picture}
}

\newcommand{\HecSeventeenb}{
\setlength{\unitlength}{1mm}
\begin{picture}(40,20)(5,10)
\put(8,15){$b:=[\acchi^{(2,17)}_7]$}
\put(10,0){\circle{2}}\put(25,0){\circle{2}}\put(40,0){\circle{2}}
\put(11,0){\line(1,0){13}}\put(26,0){\line(1,0){13}}
\put(4,3){{\small{$-1$}}}\put(23,3){{\small{$-1$}}}\put(41,3){{\small{$\defzeta$}}}
\put(15,2.0){{\small{$-1$}}}\put(31,2.0){{\small{$\defzeta^{-1}$}}}
\put(8,-4){$\al_1$}\put(23,-4){$\al_2$}\put(38,-4){$\al_3$}
\put(25,1){\oval(30,14)[t]}\put(23,10){{\small{$-\defzeta$}}}
\end{picture}
}

\newcommand{\HecSeventeenc}{
\setlength{\unitlength}{1mm}
\begin{picture}(40,20)(5,10)
\put(8,10){$c:=[\acchi^{(2,17)}_3]$}
\put(10,0){\circle{2}}\put(25,0){\circle{2}}\put(40,0){\circle{2}}
\put(11,0){\line(1,0){13}}\put(26,0){\line(1,0){13}}
\put(8,3){{\small{$-1$}}}\put(24,3){{\small{$\defzeta$}}}\put(38,3){{\small{$-1$}}}
\put(15,2.0){{\small{$-1$}}}\put(30.5,2.0){{\small{$\defzeta^{-1}$}}}
\put(8,-4){$\al_1$}\put(23,-4){$\al_2$}\put(38,-4){$\al_3$}
\end{picture}
}

\newcommand{\HecSeventeend}{
\setlength{\unitlength}{1mm}
\begin{picture}(40,20)(5,10)
\put(8,10){$d:=[\acchi^{(2,17)}_9]\circ\left[{{123}\atop{213}}\right]$}
\put(10,0){\circle{2}}\put(25,0){\circle{2}}\put(40,0){\circle{2}}
\put(11,0){\line(1,0){13}}\put(26,0){\line(1,0){13}}
\put(8,3){{\small{$-1$}}}\put(23,3){{\small{$-1$}}}\put(39,3){{\small{$\defzeta^{-1}$}}}
\put(15,2.0){{\small{$-1$}}}\put(29,2.0){{\small{$-\defzeta^{-1}$}}}
\put(8,-4){$\al_2$}\put(23,-4){$\al_1$}\put(38,-4){$\al_3$}
\end{picture}
}

\newcommand{\HecSeventeene}{
\setlength{\unitlength}{1mm}
\begin{picture}(40,20)(5,10)
\put(8,15){$e:=[\acchi^{(2,17)}_2]$}
\put(10,0){\circle{2}}\put(25,0){\circle{2}}\put(40,0){\circle{2}}
\put(11,0){\line(1,0){13}}\put(26,0){\line(1,0){13}}
\put(4,3){{\small{$-\defzeta$}}}\put(24,3){{\small{$\defzeta$}}}\put(41,3){{\small{$-1$}}}
\put(14,2.0){{\small{$-\defzeta^{-1}$}}}\put(31,2.0){{\small{$\defzeta^{-1}$}}}
\put(8,-4){$\al_1$}\put(23,-4){$\al_2$}\put(38,-4){$\al_3$}
\put(25,1){\oval(30,14)[t]}\put(22,10){{\small{$-\defzeta^{-1}$}}}
\end{picture}
}

\newcommand{\HecSeventeenf}{
\setlength{\unitlength}{1mm}
\begin{picture}(40,20)(5,10)
\put(8,10){$f:=[\acchi^{(2,17)}_4]\circ\left[{{123}\atop{312}}\right]$}
\put(10,0){\circle{2}}\put(25,0){\circle{2}}\put(40,0){\circle{2}}
\put(11,0){\line(1,0){13}}\put(26,0){\line(1,0){13}}
\put(8,3){{\small{$-1$}}}\put(23,3){{\small{$-1$}}}\put(38,3){{\small{$-1$}}}
\put(17,2.0){{\small{$\defzeta$}}}\put(30,2.0){{\small{$-\defzeta$}}}
\put(8,-4){$\al_2$}\put(23,-4){$\al_3$}\put(38,-4){$\al_1$}
\end{picture}
}

\newcommand{\HecSeventeeng}{
\setlength{\unitlength}{1mm}
\begin{picture}(40,20)(5,10)
\put(8,10){$g:=[\acchi^{(2,17)}_8]\circ\left[{{123}\atop{312}}\right]$}
\put(10,0){\circle{2}}\put(25,0){\circle{2}}\put(40,0){\circle{2}}
\put(11,0){\line(1,0){13}}\put(26,0){\line(1,0){13}}
\put(8,3){{\small{$-1$}}}\put(24,3){{\small{$\defzeta$}}}\put(38,3){{\small{$-1$}}}
\put(15.5,2.0){{\small{$\defzeta^{-1}$}}}\put(30,2.0){{\small{$-\defzeta$}}}
\put(8,-4){$\al_2$}\put(23,-4){$\al_3$}\put(38,-4){$\al_1$}
\end{picture}
}

\newcommand{\HecSeventeenh}{
\setlength{\unitlength}{1mm}
\begin{picture}(40,20)(5,10)
\put(8,10){$h:=[\acchi^{(2,17)}_6]\circ\left[{{123}\atop{312}}\right]$}
\put(10,0){\circle{2}}\put(25,0){\circle{2}}\put(40,0){\circle{2}}
\put(11,0){\line(1,0){13}}\put(26,0){\line(1,0){13}}
\put(8,3){{\small{$-1$}}}\put(22,3){{\small{$-\defzeta$}}}\put(38,3){{\small{$-1$}}}
\put(17,2.0){{\small{$\defzeta$}}}\put(28.5,2.0){{\small{$-\defzeta^{-1}$}}}
\put(8,-4){$\al_2$}\put(23,-4){$\al_3$}\put(38,-4){$\al_1$}
\end{picture}
}

\newcommand{\HecSeventeeni}{
\setlength{\unitlength}{1mm}
\begin{picture}(40,20)(5,10)
\put(8,10){$i:=[\acchi^{(2,17)}_5]\circ\left[{{123}\atop{312}}\right]$}
\put(10,0){\circle{2}}\put(25,0){\circle{2}}\put(40,0){\circle{2}}
\put(11,0){\line(1,0){13}}\put(26,0){\line(1,0){13}}
\put(8,3){{\small{$-1$}}}\put(22.5,3){{\small{$\defzeta^{-1}$}}}\put(38,3){{\small{$-1$}}}
\put(15.5,2.0){{\small{$\defzeta^{-1}$}}}\put(28.5,2.0){{\small{$-\defzeta^{-1}$}}}
\put(8,-4){$\al_2$}\put(23,-4){$\al_3$}\put(38,-4){$\al_1$}
\end{picture}
}

\newcommand{\HecSeventeenabcdefghi}
{\setlength{\unitlength}{1mm}
\begin{picture}(130,90)(15,-20)
\put(0,60){\HecSeventeena}
\put(45,60){\HecSeventeenb}
\put(90,60){\HecSeventeenc}

\put(0,30){\HecSeventeend}
\put(45,30){\HecSeventeene}
\put(90,30){\HecSeventeenf}

\put(0,5){\HecSeventeeng}
\put(45,5){\HecSeventeenh}
\put(90,5){\HecSeventeeni}

\put(5,-20){$(\zeta\in\bKt,\,\zeta^2+\zeta+1=0)$}

\end{picture}
}

\newcommand{\CayleySeventeen}{
\setlength{\unitlength}{1mm}
\begin{picture}(110,65)(-15,-5)

\put(10,0){\holizontallineAboveone}
\put(20,0){\whitecircleSouthe}
\put(20,0){\verticallineRightthreeTHICK}
\put(30,0){\holizontallineAbovetwoTHICK}\put(21,0){\holizontallinetenmmTHICK}\put(39,0){\holizontallinetenmmTHICK}
\put(50,0){\whitecircleSouthc}
\put(50,0){\verticallineRightthreeTHICK}
\put(50,0){\holizontallineAboveone}

\put(20,10){\whitecircleNorthf}
\put(20,10){\EastUplineabovetwoTHICK}
\put(50,10){\whitecircleNortha}
\put(50,10){\holizontallineAboveone}

\put(0,20){\holizontallineAbovethreeTHICK}
\put(10,20){\whitecircleNorthh}
\put(10,20){\EastUplineabovetwoTHICK}
\put(10,20){\EastDownlineaboveone}
\put(30,20){\whitecircleNorthg}
\put(30,20){\holizontallineAbovethree}
\put(40,20){\whitecircleNorthb}
\put(40,20){\EastUplineaboveoneTHICK}
\put(40,20){\EastDownlineabovetwoTHICK}

\put(20,30){\whitecircleSouthi}
\put(20,30){\verticallineRightthree}
\put(20,30){\EastDownlineaboveoneTHICK}
\put(50,30){\whitecircleSouthd}
\put(50,30){\verticallineRightthreeTHICK}
\put(50,30){\holizontallineAbovetwo}

\put(20,40){\whitecircleNorthi}
\put(20,40){\EastUplineaboveoneTHICK}
\put(50,40){\whitecircleNorthd}
\put(50,40){\holizontallineAbovethree}

\put(0,50){\holizontallineAbovetwo}
\put(10,50){\whitecircleNorthh}
\put(10,50){\EastUplineaboveoneTHICK}
\put(10,50){\EastDownlineabovetwoTHICK}
\put(30,50){\whitecircleNorthg}
\put(30,50){\holizontallineAbovethreeTHICK}
\put(40,50){\whitecircleNorthb}
\put(40,50){\EastUplineabovetwo}
\put(40,50){\EastDownlineaboveoneTHICK}

\put(20,60){\whitecircleSouthf}
\put(20,60){\EastDownlineabovetwo}
\put(20,60){\verticallineRightthreeTHICK}
\put(50,60){\whitecircleSoutha}
\put(50,60){\verticallineRightthreeTHICK}
\put(50,60){\holizontallineAboveoneTHICK}

\put(10,70){\holizontallineAboveone}
\put(20,70){\whitecircleNorthe}
\put(30,70){\holizontallineAbovetwoTHICK}\put(21,70){\holizontallinetenmmTHICK}\put(39,70){\holizontallinetenmmTHICK}
\put(50,70){\whitecircleNorthc}
\put(50,70){\holizontallineAboveone}

\multiput(0,70)(2,0){6}{\line(1,0){0.5}}\put(-16,71){$A_{i-(-1)^i}=A_i$}
\multiput(-10,50)(2,0){6}{\line(1,0){0.5}}\put(-26,51){$B_{i-(-1)^i}=B_i$}
\multiput(-10,20)(2,0){6}{\line(1,0){0.5}}\put(-26,21){$C_{i-(-1)^i}=C_i$}
\multiput(0,0)(2,0){6}{\line(1,0){0.5}}\put(-16,1){$D_{i-(-1)^i}=D_i$}

\multiput(60,70)(2,0){6}{\line(1,0){0.5}}\put(61,71){$E_i=E_{i+(-1)^i}$}
\multiput(60,60)(2,0){6}{\line(1,0){0.5}}\put(61,61){$F_i=F_{i+(-1)^i}$}
\multiput(60,40)(2,0){6}{\line(1,0){0.5}}\put(61,41){$G_i=G_{i+(-1)^i}$}\put(90,34){$(i\in\bZ/4\bZ)$}
\multiput(60,30)(2,0){6}{\line(1,0){0.5}}\put(61,31){$H_i=H_{i+(-1)^i}$}
\multiput(60,10)(2,0){6}{\line(1,0){0.5}}\put(61,11){$I_i=I_{i+(-1)^i}$}
\multiput(60,0)(2,0){6}{\line(1,0){0.5}}\put(61,1){$J_i=J_{i+(-1)^i}$}

\end{picture}
}

\newcommand{\HecEighteena}{
\setlength{\unitlength}{1mm}
\begin{picture}(40,20)(5,10)
\put(8,10){$a:=[\acchi^{(3,18)}_1]$}
\put(10,0){\circle{2}}\put(25,0){\circle{2}}\put(40,0){\circle{2}}
\put(11,0){\line(1,0){13}}\put(26,0){\line(1,0){13}}
\put(9,3){{\small{$\zeta$}}}\put(24,3){{\small{$\zeta$}}}\put(39,3){{\small{$\zeta^{-3}$}}}
\put(15.5,2.0){{\small{$\zeta^{-1}$}}}\put(30.5,2.0){{\small{$\zeta^{-1}$}}}
\put(8,-4){$\al_1$}\put(23,-4){$\al_2$}\put(38,-4){$\al_3$}
\end{picture}
}

\newcommand{\HecEighteenb}{
\setlength{\unitlength}{1mm}
\begin{picture}(40,20)(5,10)
\put(8,10){$b:=[\acchi^{(3,18)}_2]$}
\put(10,0){\circle{2}}\put(25,0){\circle{2}}\put(40,0){\circle{2}}
\put(11,0){\line(1,0){13}}\put(26,0){\line(1,0){13}}
\put(9,3){{\small{$\zeta$}}}\put(24,3){{\small{$\zeta^{-4}$}}}\put(39,3){{\small{$\zeta^{-3}$}}}
\put(15.5,2.0){{\small{$\zeta^{-1}$}}}\put(32,2.0){{\small{$\zeta^4$}}}
\put(8,-4){$\al_1$}\put(23,-4){$\al_2$}\put(38,-4){$\al_3$}
\end{picture}
}

\newcommand{\HecEighteenab}
{\setlength{\unitlength}{1mm}
\begin{picture}(70,20)(15,-15)
\put(0,0){\HecEighteena}
\put(45,0){\HecEighteenb}

\put(90,-10){$(\zeta\in\bKt,\,\zeta^6+\zeta^3+1=0)$}
\end{picture}
}

\newcommand{\CayleyEighteen}{
\setlength{\unitlength}{1mm}
\begin{picture}(110,100)(-10,-10)
\put(0,0){\holizontallineAbovetwo}
\put(10,0){\whitecircleSouthb}
\put(10,0){\verticallineRightthreeTHICK}
\put(10,0){\holizontallineAboveoneTHICK}
\put(20,0){\whitecircleSouthb}
\put(20,0){\verticallineRightthreeTHICK}
\put(20,0){\holizontallineAbovetwo}

\put(10,10){\whitecircleNortha}
\put(10,10){\holizontallineAboveone}
\put(20,10){\whitecircleNortha}
\put(20,10){\EastUplineabovetwoTHICK}

\put(0,20){\whitecircleNortha}
\put(0,20){\EastUplineaboveoneTHICK}
\put(0,20){\EastDownlineabovetwoTHICK}
\put(30,20){\whitecircleNortha}
\put(30,20){\EastUplineabovethreeTHICK}

\put(-10,30){\EastDownlineabovethree}
\put(10,30){\whitecircleNortha}
\put(10,30){\holizontallineAbovetwoTHICK}
\put(20,30){\whitecircleNortha}
\put(20,30){\EastUplineabovethreeTHICK}
\put(20,30){\EastDownlineaboveone}
\put(40,30){\whitecircleNorthb}
\put(40,30){\holizontallineAbovetwoTHICK}
\put(50,30){\whitecircleNorthb}
\put(50,30){\EastUplineaboveoneTHICK}
\put(50,30){\EastDownlineabovethree}

\put(0,40){\EastDownlineabovethree}
\put(30,40){\whitecircleSouthb}
\put(30,40){\verticallineRighttwoTHICK}
\put(30,40){\EastDownlineaboveone}
\put(60,40){\whitecircleSouthb}
\put(60,40){\verticallineRighttwoTHICK}
\put(60,40){\EastDownlineabovethree}

\put(0,50){\EastUplineabovethree}
\put(30,50){\whitecircleNorthb}
\put(30,50){\EastUplineaboveone}
\put(60,50){\whitecircleNorthb}
\put(60,50){\EastUplineabovethree}

\put(-10,60){\EastUplineabovethree}
\put(10,60){\whitecircleNortha}
\put(10,60){\holizontallineAbovetwoTHICK}
\put(20,60){\whitecircleNortha}
\put(20,60){\EastUplineaboveone}
\put(20,60){\EastDownlineabovethreeTHICK}
\put(40,60){\whitecircleNorthb}
\put(40,60){\holizontallineAbovetwoTHICK}
\put(50,60){\whitecircleNorthb}
\put(50,60){\EastUplineabovethree}
\put(50,60){\EastDownlineaboveoneTHICK}

\put(0,70){\whitecircleNortha}
\put(0,70){\EastUplineabovetwoTHICK}
\put(0,70){\EastDownlineaboveoneTHICK}
\put(30,70){\whitecircleNortha}
\put(30,70){\EastDownlineabovethreeTHICK}

\put(10,80){\whitecircleSoutha}
\put(10,80){\verticallineRightthreeTHICK}
\put(10,80){\holizontallineAboveone}
\put(20,80){\whitecircleSoutha}
\put(20,80){\verticallineRightthreeTHICK}
\put(20,80){\EastDownlineabovetwoTHICK}

\put(0,90){\holizontallineAbovetwoTHICK}
\put(10,90){\whitecircleNorthb}
\put(10,90){\holizontallineAboveone}
\put(20,90){\whitecircleNorthb}
\put(20,90){\holizontallineAbovetwoTHICK}

\multiput(-10,90)(2,0){6}{\line(1,0){0.5}}\put(-7,91){$A_i$}
\multiput(-20,60)(2,0){6}{\line(1,0){0.5}}\put(-17,61){$B_i$}
\multiput(-10,50)(2,0){6}{\line(1,0){0.5}}\put(-7,51){$C_i$}
\multiput(-10,40)(2,0){6}{\line(1,0){0.5}}\put(-7,41){$D_i$}
\multiput(-20,30)(2,0){6}{\line(1,0){0.5}}\put(-17,31){$E_i$}
\multiput(-10,0)(2,0){6}{\line(1,0){0.5}}\put(-7,1){$F_i$}

\multiput(30,90)(2,0){6}{\line(1,0){0.5}}\put(32,91){$A_{i+1}$}
\multiput(60,70)(2,0){6}{\line(1,0){0.5}}\put(62,71){$B_{i+1}$}
\multiput(70,60)(2,0){6}{\line(1,0){0.5}}\put(72,61){$C_{i+1}$}
\multiput(70,30)(2,0){6}{\line(1,0){0.5}}\put(72,31){$D_{i+1}$}
\multiput(60,20)(2,0){6}{\line(1,0){0.5}}\put(62,21){$E_{i+1}$}
\multiput(30,0)(2,0){6}{\line(1,0){0.5}}\put(32,1){$F_{i+1}$}

\end{picture}
}

\newcommand{\TopHecFourQUCX}{
\setlength{\unitlength}{1mm}
\begin{picture}(220,70)(00,-10)
\put(10,50){$\smallHecFourOneTop$}
\put(40,50){$\smallHecFourTwoTop$}
\put(70,50){$\smallHecFourThreeTop$}
\put(100,50){$\smallHecFourFourTop$}
\put(10,25){$\smallHecFourFiveTop$}
\put(40,25){$\smallHecFourSixTop$}
\put(70,25){$\smallHecFourSevenTop$}
\put(100,25){$\smallHecFourTenTop$}
\put(10,00){$\smallHecFourElevenTop$}
\put(40,00){$\smallHecFourFifteenTop$}
\put(70,00){$\smallHecFourSixteenTop$}
\put(100,00){$\smallHecFourNineteenTop$}
\end{picture}}

\newcommand{\smallHecFourOneTop}{
\setlength{\unitlength}{0.5mm}
\begin{picture}(30,20)(00,00)

\put(-5,16){\tiny{$\acchi^{(3,1)}_1$}}

\multiput(00,00)(15,00){4}{\circle{2}}
\multiput(01,00)(15,00){3}{\line(1,0){13}}

\put(-3,-15){\tiny{$(q\in\bKt,\,q\ne 1)$}}

\put(-1,3){\tiny{$q$}}
\put(14,3){\tiny{$q$}}
\put(29,3){\tiny{$q$}}
\put(43,3){\tiny{$q$}}

\put(3.5,2){\tiny{$q^{-1}$}}
\put(18.5,2){\tiny{$q^{-1}$}}
\put(33.5,2){\tiny{$q^{-1}$}}

\put(-2,-5){\tiny{$\al_1$}}
\put(13,-5){\tiny{$\al_2$}}
\put(28,-5){\tiny{$\al_3$}}
\put(43,-5){\tiny{$\al_4$}}
\end{picture}}

\newcommand{\smallHecFourTwoTop}{
\setlength{\unitlength}{0.5mm}
\begin{picture}(30,20)(00,00)

\put(-5,16){\tiny{$\acchi^{(3,2)}_1$}}

\multiput(00,00)(15,00){4}{\circle{2}}
\multiput(01,00)(15,00){3}{\line(1,0){13}}

\put(-3,-15){\tiny{$(q\in\bKt,\,q^2\ne 1)$}}

\put(-1,3){\tiny{$q^2$}}
\put(14,3){\tiny{$q^2$}}
\put(29,3){\tiny{$q^2$}}
\put(43,3){\tiny{$q$}}

\put(3.5,2){\tiny{$q^{-2}$}}
\put(18.5,2){\tiny{$q^{-2}$}}
\put(33.5,2){\tiny{$q^{-2}$}}

\put(-2,-5){\tiny{$\al_1$}}
\put(13,-5){\tiny{$\al_2$}}
\put(28,-5){\tiny{$\al_3$}}
\put(43,-5){\tiny{$\al_4$}}
\end{picture}}

\newcommand{\smallHecFourThreeTop}{
\setlength{\unitlength}{0.5mm}
\begin{picture}(30,20)(00,00)

\put(-5,16){\tiny{$\acchi^{(3,3)}_1$}}

\multiput(00,00)(15,00){4}{\circle{2}}
\multiput(01,00)(15,00){3}{\line(1,0){13}}

\put(-3,-15){\tiny{$(q\in\bKt,\,q^2\ne 1)$}}

\put(-1,3){\tiny{$q$}}
\put(14,3){\tiny{$q$}}
\put(29,3){\tiny{$q$}}
\put(43,3){\tiny{$q^2$}}

\put(3.5,2){\tiny{$q^{-1}$}}
\put(18.5,2){\tiny{$q^{-1}$}}
\put(33.5,2){\tiny{$q^{-2}$}}

\put(-2,-5){\tiny{$\al_1$}}
\put(13,-5){\tiny{$\al_2$}}
\put(28,-5){\tiny{$\al_3$}}
\put(43,-5){\tiny{$\al_4$}}
\end{picture}}

\newcommand{\smallHecFourFourTop}{
\setlength{\unitlength}{0.5mm}
\begin{picture}(30,20)(00,00)

\put(-5,16){\tiny{$\acchi^{(3,4)}_1$}}

\multiput(00,00)(15,00){4}{\circle{2}}
\multiput(01,00)(15,00){3}{\line(1,0){13}}

\put(-3,-15){\tiny{$(q\in\bKt,\,q^2\ne 1)$}}

\put(-1,3){\tiny{$q^2$}}
\put(14,3){\tiny{$q^2$}}
\put(29,3){\tiny{$q$}}
\put(43,3){\tiny{$q$}}

\put(3.5,2){\tiny{$q^{-2}$}}
\put(18.5,2){\tiny{$q^{-2}$}}
\put(33.5,2){\tiny{$q^{-1}$}}

\put(-2,-5){\tiny{$\al_1$}}
\put(13,-5){\tiny{$\al_2$}}
\put(28,-5){\tiny{$\al_3$}}
\put(43,-5){\tiny{$\al_4$}}
\end{picture}}

\newcommand{\smallHecFourFiveTop}{
\setlength{\unitlength}{0.5mm}
\begin{picture}(30,20)(00,00)

\put(-5,16){\tiny{$\acchi^{(3,5)}_1$}}

\multiput(15.8,3.8)(00,00){1}{\oval(6.4,6.4)[br]}
\multiput(32,3.8)(00,00){1}{\oval(26,13)[t]}
\multiput(45,01)(00,00){1}{\line(0,1){2.8}}

\multiput(00,00)(15,00){4}{\circle{2}}
\multiput(01,00)(15,00){2}{\line(1,0){13}}

\put(-3,-15){\tiny{$(q\in\bKt,\,q\ne 1)$}}

\put(-1,3){\tiny{$q$}}
\put(14,3){\tiny{$q$}}
\put(29.5,3){\tiny{$q$}}
\put(47,2){\tiny{$q$}}

\put(3.5,2){\tiny{$q^{-1}$}}
\put(20,2){\tiny{$q^{-1}$}}
\put(31,12){\tiny{$q^{-1}$}}

\put(-2,-5){\tiny{$\al_1$}}
\put(13,-5){\tiny{$\al_2$}}
\put(28,-5){\tiny{$\al_3$}}
\put(43,-5){\tiny{$\al_4$}}
\end{picture}}

\newcommand{\smallHecFourSixTop}{
\setlength{\unitlength}{0.5mm}
\begin{picture}(30,20)(00,00)

\put(-5,16){\tiny{$\acchi^{(3,6)}_1$}}

\multiput(00,00)(15,00){4}{\circle{2}}
\multiput(01,00)(15,00){3}{\line(1,0){13}}

\put(-3,-15){\tiny{$(q\in\bKt,\,q^2\ne 1)$}}

\put(-4,3){\tiny{$-1$}}
\put(14,3){\tiny{$q$}}
\put(29,3){\tiny{$q$}}
\put(43,3){\tiny{$q$}}

\put(3.5,2){\tiny{$q^{-1}$}}
\put(18.5,2){\tiny{$q^{-1}$}}
\put(33.5,2){\tiny{$q^{-1}$}}

\put(-2,-5){\tiny{$\al_1$}}
\put(13,-5){\tiny{$\al_2$}}
\put(28,-5){\tiny{$\al_3$}}
\put(43,-5){\tiny{$\al_4$}}
\end{picture}}

\newcommand{\smallHecFourSevenTop}{
\setlength{\unitlength}{0.5mm}
\begin{picture}(30,20)(00,00)

\put(-5,16){\tiny{$\acchi^{(3,7)}_1$}}

\multiput(00,00)(15,00){4}{\circle{2}}
\multiput(01,00)(15,00){3}{\line(1,0){13}}

\put(-3,-15){\tiny{$(q\in\bKt,\,q^4\ne 1)$}}

\put(-4,3){\tiny{$-1$}}
\put(13.5,3){\tiny{$q^2$}}
\put(28.5,3){\tiny{$q^2$}}
\put(43,3){\tiny{$q$}}

\put(3.5,2){\tiny{$q^{-2}$}}
\put(18.5,2){\tiny{$q^{-2}$}}
\put(33.5,2){\tiny{$q^{-2}$}}

\put(-2,-5){\tiny{$\al_1$}}
\put(13,-5){\tiny{$\al_2$}}
\put(28,-5){\tiny{$\al_3$}}
\put(43,-5){\tiny{$\al_4$}}
\end{picture}}

\newcommand{\smallHecFourTenTop}{
\setlength{\unitlength}{0.5mm}
\begin{picture}(30,20)(00,00)

\put(-5,16){\tiny{$\acchi^{(3,10)}_1$}}

\multiput(00,00)(15,00){4}{\circle{2}}
\multiput(01,00)(15,00){3}{\line(1,0){13}}

\put(-3,-15){\tiny{$(q\in\bKt,\,q^4\ne 1)$}}

\put(-4,3){\tiny{$q^{-1}$}}
\put(11,3){\tiny{$-1$}}
\put(28.5,3){\tiny{$q$}}
\put(43,3){\tiny{$q$}}

\put(6,2){\tiny{$q$}}
\put(18.5,2){\tiny{$q^{-1}$}}
\put(33.5,2){\tiny{$q^{-1}$}}

\put(-2,-5){\tiny{$\al_1$}}
\put(13,-5){\tiny{$\al_2$}}
\put(28,-5){\tiny{$\al_3$}}
\put(43,-5){\tiny{$\al_4$}}
\end{picture}}

\newcommand{\smallHecFourElevenTop}{
\setlength{\unitlength}{0.5mm}
\begin{picture}(30,20)(00,00)

\put(-5,16){\tiny{$\acchi^{(3,11)}_1$}}

\multiput(00,00)(15,00){4}{\circle{2}}
\multiput(01,00)(15,00){3}{\line(1,0){13}}

\put(-3,-15){\tiny{$(q\in\bKt,\,q^4\ne 1)$}}

\put(-4,3){\tiny{$q^{-2}$}}
\put(11,3){\tiny{$-1$}}
\put(28,3){\tiny{$q^2$}}
\put(43,3){\tiny{$q$}}

\put(5,2){\tiny{$q^2$}}
\put(18.5,2){\tiny{$q^{-2}$}}
\put(33.5,2){\tiny{$q^{-2}$}}

\put(-2,-5){\tiny{$\al_1$}}
\put(13,-5){\tiny{$\al_2$}}
\put(28,-5){\tiny{$\al_3$}}
\put(43,-5){\tiny{$\al_4$}}
\end{picture}}

\newcommand{\smallHecFourFifteenTop}{
\setlength{\unitlength}{0.5mm}
\begin{picture}(30,20)(00,00)

\put(-5,16){\tiny{$\acchi^{(3,15)}_1$}}

\multiput(00,00)(15,00){4}{\circle{2}}
\multiput(01,00)(15,00){3}{\line(1,0){13}}

\put(-3,-15){\tiny{$(\zeta\in\bKt,\,\zeta^2+\zeta+1=0)$}}

\put(-5,4.5){\tiny{$-{\frac 1 \zeta}$}}
\put(10,4.5){\tiny{$-{\frac 1 \zeta}$}}
\put(25,4.5){\tiny{$-{\frac 1 \zeta}$}}
\put(43,3){\tiny{$\zeta$}}

\put(3,2){\tiny{$-\zeta$}}
\put(18,2){\tiny{$-\zeta$}}
\put(33,2){\tiny{$-\zeta$}}

\put(-2,-5){\tiny{$\al_1$}}
\put(13,-5){\tiny{$\al_2$}}
\put(28,-5){\tiny{$\al_3$}}
\put(43,-5){\tiny{$\al_4$}}
\end{picture}}

\newcommand{\smallHecFourSixteenTop}{
\setlength{\unitlength}{0.5mm}
\begin{picture}(30,20)(00,00)

\put(-5,16){\tiny{$\acchi^{(3,16)}_1$}}

\multiput(00,00)(15,00){4}{\circle{2}}
\multiput(01,00)(15,00){3}{\line(1,0){13}}

\put(-3,-15){\tiny{$(\zeta\in\bKt,\,\zeta^2+\zeta+1=0)$}}

\put(-3,4.5){\tiny{$-1$}}
\put(10,4.5){\tiny{$-{\frac 1 \zeta}$}}
\put(25,4.5){\tiny{$-{\frac 1 \zeta}$}}
\put(43,3){\tiny{$\zeta$}}

\put(3,2){\tiny{$-\zeta$}}
\put(18,2){\tiny{$-\zeta$}}
\put(33,2){\tiny{$-\zeta$}}

\put(-2,-5){\tiny{$\al_1$}}
\put(13,-5){\tiny{$\al_2$}}
\put(28,-5){\tiny{$\al_3$}}
\put(43,-5){\tiny{$\al_4$}}
\end{picture}}

\newcommand{\smallHecFourNineteenTop}{
\setlength{\unitlength}{0.5mm}
\begin{picture}(30,20)(00,00)

\put(-5,16){\tiny{$\acchi^{(3,19)}_1$}}

\multiput(00,00)(15,00){4}{\circle{2}}
\multiput(01,00)(15,00){3}{\line(1,0){13}}

\put(-3,-15){\tiny{$(\zeta\in\bKt,\,\zeta^2+\zeta+1=0)$}}

\put(-6,3){\tiny{$-\zeta$}}
\put(11,3){\tiny{$-1$}}
\put(25,4.5){\tiny{$-{\frac 1 \zeta}$}}
\put(43,3){\tiny{$\zeta$}}

\put(2.5,3.5){\tiny{$-{\frac 1 \zeta}$}}
\put(18,2){\tiny{$-\zeta$}}
\put(33,2){\tiny{$-\zeta$}}

\put(-2,-5){\tiny{$\al_1$}}
\put(13,-5){\tiny{$\al_2$}}
\put(28,-5){\tiny{$\al_3$}}
\put(43,-5){\tiny{$\al_4$}}
\end{picture}}

\newcommand{\smallHecFourNineTop}{
\setlength{\unitlength}{0.5mm}
\begin{picture}(30,20)(00,00)

\put(-5,16){\tiny{$\acchi^{(3,9)}_1$}}

\multiput(00,00)(15,00){4}{\circle{2}}
\multiput(01,00)(15,00){3}{\line(1,0){13}}

\put(-3,-15){\tiny{$(q\in\bKt,\,q^2\ne 1\ne q^3)$}}

\put(-1,3){\tiny{$q^2$}}
\put(14,3){\tiny{$q^2$}}
\put(29,3){\tiny{$q$}}
\put(43,3){\tiny{$-1$}}

\put(3.5,2){\tiny{$q^{-2}$}}
\put(21,2){\tiny{$q^2$}}
\put(33.5,2){\tiny{$q^{-1}$}}

\put(-2,-5){\tiny{$\al_1$}}
\put(13,-5){\tiny{$\al_2$}}
\put(28,-5){\tiny{$\al_3$}}
\put(43,-5){\tiny{$\al_4$}}
\end{picture}}

\newcommand{\smallHecFourTwelveTop}{
\setlength{\unitlength}{0.5mm}
\begin{picture}(30,20)(00,00)

\put(-5,16){\tiny{$\acchi^{(3,12)}_1$}}

\multiput(00,00)(15,00){4}{\circle{2}}
\multiput(01,00)(15,00){3}{\line(1,0){13}}

\put(-1,-15){\tiny{$(q\in\bKt, q^2\ne 1)$}}

\put(-3,3){\tiny{$q^{-1}$}}
\put(11,3){\tiny{$-1$}}
\put(29,3){\tiny{$q$}}
\put(44,3){\tiny{$q^2$}}

\put(6,2){\tiny{$q$}}
\put(18.5,2){\tiny{$q^{-1}$}}
\put(33.5,2){\tiny{$q^{-2}$}}

\put(-2,-5){\tiny{$\al_1$}}
\put(13,-5){\tiny{$\al_2$}}
\put(28,-5){\tiny{$\al_3$}}
\put(43,-5){\tiny{$\al_4$}}
\end{picture}}

\newcommand{\smallHecFourFourteenTop}{
\setlength{\unitlength}{0.5mm}
\begin{picture}(30,20)(00,00)

\put(-1,-15){\tiny{$(q\in\bKt, q^2\ne 1)$}}

\put(-5,16){\tiny{$\acchi^{(3,14)}_1$}}

\multiput(00,00)(15,00){4}{\circle{2}}
\multiput(01,00)(15,00){3}{\line(1,0){13}}

\put(-1,3){\tiny{$q$}}
\put(14,3){\tiny{$q$}}
\put(28,3){\tiny{$-1$}}
\put(40,-5){\tiny{$-q^{-1}$}}

\put(3.5,2){\tiny{$q^{-1}$}}
\put(18.5,2){\tiny{$q^{-1}$}}
\put(35,2){\tiny{$-q$}}

\put(-2,-5){\tiny{$\al_1$}}
\put(13,-5){\tiny{$\al_2$}}
\put(28,-5){\tiny{$\al_3$}}
\put(43,3){\tiny{$\al_4$}}
\end{picture}}

\newcommand{\TopHecFiveEightQUCX}{
\setlength{\unitlength}{1mm}
\begin{picture}(220,125)(00,-10)
\put(10,100){$\smallHecFiveElevenTop$}
\put(50,100){$\smallHecFiveTwelveTop$}
\put(90,100){$\smallHecFiveThirteenTop$}

\put(10,75){$\smallHecFiveFourteenTop$}
\put(55,75){$\smallHecFiveFifteenTop$}
\put(100,75){$\smallHecSixSixteenTop$}

\put(10,50){$\smallHecSixSeventeenTop$}
\put(55,50){$\smallHecSixEighteenTop$}
\put(100,50){$\smallHecSixNineteenTop$}

\put(10,25){$\smallHecSevenTwentyTop$}
\put(70,25){$\smallHecSevenTwetyoneTop$}

\put(10,00){$\smallHecSevenTwentyTwoTop$}
\end{picture}}

\newcommand{\smallHecFiveElevenTop}{
\setlength{\unitlength}{0.5mm}
\begin{picture}(30,20)(00,00)

\put(-5,16){\tiny{$\acchi^{(4,11)}_1$}}

\multiput(00,00)(15,00){5}{\circle{2}}
\multiput(01,00)(15,00){4}{\line(1,0){13}}

\put(-3,-15){\tiny{$(\zeta\in\bKt,\,\zeta^2+\zeta+1=0)$}}

\put(-1,3){\tiny{$\zeta$}}
\put(14,3){\tiny{$\zeta$}}
\put(27,3){\tiny{$-1$}}
\put(44,3){\tiny{$\zeta$}}
\put(59,3){\tiny{$\zeta$}}

\put(3.5,2){\tiny{$\zeta^{-1}$}}
\put(18,2){\tiny{$\zeta^{-1}$}}
\put(33.5,2){\tiny{$\zeta^{-1}$}}
\put(48.5,2){\tiny{$\zeta^{-1}$}}

\put(-2,-5){\tiny{$\al_1$}}
\put(13,-5){\tiny{$\al_2$}}
\put(28,-5){\tiny{$\al_3$}}
\put(43,-5){\tiny{$\al_4$}}
\put(58,-5){\tiny{$\al_5$}}
\end{picture}}

\newcommand{\smallHecFiveTwelveTop}{
\setlength{\unitlength}{0.5mm}
\begin{picture}(30,20)(00,00)

\put(-5,16){\tiny{$\acchi^{(4,12)}_1$}}

\multiput(00,00)(15,00){5}{\circle{2}}
\multiput(01,00)(15,00){4}{\line(1,0){13}}

\put(-3,-15){\tiny{$(\zeta\in\bKt,\,\zeta^2+\zeta+1=0)$}}

\put(-1,3){\tiny{$\zeta$}}
\put(14,3){\tiny{$\zeta$}}
\put(29,3){\tiny{$\zeta$}}
\put(43,5){\tiny{${\frac 1 \zeta}$}}
\put(57,3){\tiny{$-1$}}

\put(3.5,2){\tiny{$\zeta^{-1}$}}
\put(18.5,2){\tiny{$\zeta^{-1}$}}
\put(33.5,2){\tiny{$\zeta^{-1}$}}
\put(51,2){\tiny{$\zeta$}}

\put(-2,-5){\tiny{$\al_1$}}
\put(13,-5){\tiny{$\al_2$}}
\put(28,-5){\tiny{$\al_3$}}
\put(43,-5){\tiny{$\al_4$}}
\put(58,-5){\tiny{$\al_5$}}
\end{picture}}

\newcommand{\smallHecFiveThirteenTop}{
\setlength{\unitlength}{0.5mm}
\begin{picture}(30,20)(00,00)

\put(-5,16){\tiny{$\acchi^{(4,13)}_1$}}

\multiput(00,00)(15,00){5}{\circle{2}}
\multiput(01,00)(15,00){4}{\line(1,0){13}}

\put(-3,-15){\tiny{$(\zeta\in\bKt,\,\zeta^2+\zeta+1=0)$}}

\put(-3,3){\tiny{$-1$}}
\put(14,3){\tiny{$\zeta$}}
\put(29,3){\tiny{$\zeta$}}
\put(43,5){\tiny{${\frac 1 \zeta}$}}
\put(57,3){\tiny{$-1$}}

\put(3.5,2){\tiny{$\zeta^{-1}$}}
\put(18.5,2){\tiny{$\zeta^{-1}$}}
\put(33.5,2){\tiny{$\zeta^{-1}$}}
\put(51,2){\tiny{$\zeta$}}

\put(-2,-5){\tiny{$\al_1$}}
\put(13,-5){\tiny{$\al_2$}}
\put(28,-5){\tiny{$\al_3$}}
\put(43,-5){\tiny{$\al_4$}}
\put(58,-5){\tiny{$\al_5$}}
\end{picture}}

\newcommand{\smallHecFiveFourteenTop}{
\setlength{\unitlength}{0.5mm}
\begin{picture}(30,20)(00,00)

\put(-5,16){\tiny{$\acchi^{(4,14)}_1$}}

\multiput(00,00)(15,00){5}{\circle{2}}
\multiput(01,00)(15,00){4}{\line(1,0){13}}

\put(-3,-15){\tiny{$(\zeta\in\bKt,\,\zeta^2+1=0)$}}

\put(-1,3){\tiny{$\zeta$}}
\put(14,3){\tiny{$\zeta$}}
\put(27,3){\tiny{$-1$}}
\put(42,3){\tiny{$-1$}}
\put(57,3){\tiny{$-1$}}

\put(4,2){\tiny{$-\zeta$}}
\put(19,2){\tiny{$-\zeta$}}
\put(34.5,2){\tiny{$-1$}}
\put(49,2){\tiny{$-\zeta$}}

\put(-2,-5){\tiny{$\al_1$}}
\put(13,-5){\tiny{$\al_2$}}
\put(28,-5){\tiny{$\al_3$}}
\put(43,-5){\tiny{$\al_4$}}
\put(58,-5){\tiny{$\al_5$}}
\end{picture}}

\newcommand{\smallHecFiveFifteenTop}{
\setlength{\unitlength}{0.5mm}
\begin{picture}(30,20)(00,00)

\put(-5,16){\tiny{$\acchi^{(4,15)}_1$}}

\multiput(00,00)(15,00){5}{\circle{2}}
\multiput(01,00)(15,00){4}{\line(1,0){13}}

\put(-8,-15){\tiny{$(\zeta\in\bKt,\,\zeta^4+\zeta^3+\zeta^2+\zeta+1=0)$}}

\put(-2,3){\tiny{$\zeta^2$}}
\put(13,3){\tiny{$\zeta^2$}}
\put(27,3){\tiny{$-1$}}
\put(44,3){\tiny{$\zeta$}}
\put(58,3){\tiny{$\zeta^2$}}

\put(3.5,2){\tiny{$\zeta^{-2}$}}
\put(18,2){\tiny{$\zeta^{-2}$}}
\put(33.5,2){\tiny{$\zeta^{-1}$}}
\put(48.5,2){\tiny{$\zeta^{-2}$}}

\put(-2,-5){\tiny{$\al_1$}}
\put(13,-5){\tiny{$\al_2$}}
\put(28,-5){\tiny{$\al_3$}}
\put(43,-5){\tiny{$\al_4$}}
\put(58,-5){\tiny{$\al_5$}}
\end{picture}}

\newcommand{\smallHecSixSixteenTop}{
\setlength{\unitlength}{0.5mm}
\begin{picture}(30,20)(00,00)

\put(-5,16){\tiny{$\acchi^{(4,16)}_1$}}

\multiput(30.8,3.8)(00,00){1}{\oval(6.4,6.4)[br]}
\multiput(54.5,3.8)(00,00){1}{\oval(41,13)[t]}
\multiput(75,01)(00,00){1}{\line(0,1){2.8}}

\multiput(00,00)(15,00){6}{\circle{2}}
\multiput(01,00)(15,00){4}{\line(1,0){13}}

\put(-3,-15){\tiny{$(q\in\bKt,\,q\ne 1)$}}

\multiput(-1,3)(15,00){5}{\tiny{$q$}}

\put(77,2){\tiny{$q$}}

\multiput(4.5,2)(15,00){4}{\tiny{$q^{-1}$}}

\put(53,12){\tiny{$q^{-1}$}}

\put(-2,-5){\tiny{$\al_1$}}
\put(13,-5){\tiny{$\al_2$}}
\put(28,-5){\tiny{$\al_3$}}
\put(43,-5){\tiny{$\al_4$}}
\put(58,-5){\tiny{$\al_5$}}
\put(73,-5){\tiny{$\al_6$}}

\end{picture}}

\newcommand{\smallHecSixSeventeenTop}{
\setlength{\unitlength}{0.5mm}
\begin{picture}(30,20)(00,00)

\put(-5,16){\tiny{$\acchi^{(4,17)}_1$}}

\multiput(00,00)(15,00){6}{\circle{2}}
\multiput(01,00)(15,00){5}{\line(1,0){13}}

\put(-3,-15){\tiny{$(\zeta\in\bKt,\,\zeta^2+\zeta+1=0)$}}

\put(-1,3){\tiny{$\zeta$}}
\put(14,3){\tiny{$\zeta$}}
\put(29,3){\tiny{$\zeta$}}
\put(42,3){\tiny{$-1$}}
\put(59,3){\tiny{$\zeta$}}
\put(74,3){\tiny{$\zeta$}}

\put(3.5,2){\tiny{$\zeta^{-1}$}}
\put(18.5,2){\tiny{$\zeta^{-1}$}}
\put(33,2){\tiny{$\zeta^{-1}$}}
\put(48.5,2){\tiny{$\zeta^{-1}$}}
\put(63.5,2){\tiny{$\zeta^{-1}$}}

\put(-2,-5){\tiny{$\al_1$}}
\put(13,-5){\tiny{$\al_2$}}
\put(28,-5){\tiny{$\al_3$}}
\put(43,-5){\tiny{$\al_4$}}
\put(58,-5){\tiny{$\al_5$}}
\put(73,-5){\tiny{$\al_6$}}
\end{picture}}

\newcommand{\smallHecSixEighteenTop}{
\setlength{\unitlength}{0.5mm}
\begin{picture}(30,20)(00,00)

\put(-5,16){\tiny{$\acchi^{(4,18)}_1$}}

\multiput(00,00)(15,00){6}{\circle{2}}
\multiput(01,00)(15,00){5}{\line(1,0){13}}

\put(-3,-15){\tiny{$(\zeta\in\bKt,\,\zeta^2+\zeta+1=0)$}}

\put(-1,3){\tiny{$\zeta$}}
\put(14,3){\tiny{$\zeta$}}
\put(29,3){\tiny{$\zeta$}}
\put(44,3){\tiny{$\zeta$}}
\put(58,5){\tiny{${\frac 1 \zeta}$}}
\put(73,3){\tiny{$-1$}}

\put(3.5,2){\tiny{$\zeta^{-1}$}}
\put(18.5,2){\tiny{$\zeta^{-1}$}}
\put(33,2){\tiny{$\zeta^{-1}$}}
\put(48.5,2){\tiny{$\zeta^{-1}$}}
\put(63.5,2){\tiny{$\zeta^{-1}$}}

\put(-2,-5){\tiny{$\al_1$}}
\put(13,-5){\tiny{$\al_2$}}
\put(28,-5){\tiny{$\al_3$}}
\put(43,-5){\tiny{$\al_4$}}
\put(58,-5){\tiny{$\al_5$}}
\put(73,-5){\tiny{$\al_6$}}
\end{picture}}

\newcommand{\smallHecSixNineteenTop}{
\setlength{\unitlength}{0.5mm}
\begin{picture}(30,20)(00,00)

\put(-5,16){\tiny{$\acchi^{(4,19)}_1$}}

\multiput(00,00)(15,00){6}{\circle{2}}
\multiput(01,00)(15,00){5}{\line(1,0){13}}

\put(-3,-15){\tiny{$(\zeta\in\bKt,\,\zeta^2+\zeta+1=0)$}}

\put(-1,3){\tiny{$\zeta$}}
\put(14,3){\tiny{$\zeta$}}
\put(29,3){\tiny{$\zeta$}}
\put(42,3){\tiny{$-1$}}
\put(57,3){\tiny{$-1$}}
\put(74,3){\tiny{$\zeta$}}

\put(4,2){\tiny{$-\zeta$}}
\put(19,2){\tiny{$-\zeta$}}
\put(34,2){\tiny{$-\zeta$}}
\put(49,2){\tiny{$-1$}}
\put(64,2){\tiny{$-\zeta$}}

\put(-2,-5){\tiny{$\al_1$}}
\put(13,-5){\tiny{$\al_2$}}
\put(28,-5){\tiny{$\al_3$}}
\put(43,-5){\tiny{$\al_4$}}
\put(58,-5){\tiny{$\al_5$}}
\put(73,-5){\tiny{$\al_6$}}

\end{picture}}

\newcommand{\smallHecSevenTwentyTop}{
\setlength{\unitlength}{0.5mm}
\begin{picture}(30,20)(00,00)

\put(-5,16){\tiny{$\acchi^{(4,20)}_1$}}

\multiput(45.8,3.8)(00,00){1}{\oval(6.4,6.4)[br]}
\multiput(69.5,3.8)(00,00){1}{\oval(41,13)[t]}
\multiput(90,01)(00,00){1}{\line(0,1){2.8}}

\multiput(00,00)(15,00){7}{\circle{2}}
\multiput(01,00)(15,00){5}{\line(1,0){13}}

\put(-3,-15){\tiny{$(q\in\bKt,\,q\ne 1)$}}

\multiput(-1,3)(15,00){6}{\tiny{$q$}}

\put(92,2){\tiny{$q$}}

\multiput(4.5,2)(15,00){5}{\tiny{$q^{-1}$}}

\put(68,12){\tiny{$q^{-1}$}}

\put(-2,-5){\tiny{$\al_1$}}
\put(13,-5){\tiny{$\al_2$}}
\put(28,-5){\tiny{$\al_3$}}
\put(43,-5){\tiny{$\al_4$}}
\put(58,-5){\tiny{$\al_5$}}
\put(73,-5){\tiny{$\al_6$}}
\put(88,-5){\tiny{$\al_7$}}

\end{picture}}

\newcommand{\smallHecSevenTwetyoneTop}{
\setlength{\unitlength}{0.5mm}
\begin{picture}(30,20)(00,00)

\put(-5,16){\tiny{$\acchi^{(4,21)}_1$}}

\multiput(00,00)(15,00){7}{\circle{2}}
\multiput(01,00)(15,00){6}{\line(1,0){13}}

\put(-3,-15){\tiny{$(\zeta\in\bKt,\,\zeta^2+\zeta+1=0)$}}

\put(-1,3){\tiny{$\zeta$}}
\put(14,3){\tiny{$\zeta$}}
\put(29,3){\tiny{$\zeta$}}
\put(44,3){\tiny{$\zeta$}}
\put(57,3){\tiny{$-1$}}
\put(74,3){\tiny{$\zeta$}}
\put(89,3){\tiny{$\zeta$}}

\put(3.5,2){\tiny{$\zeta^{-1}$}}
\put(18.5,2){\tiny{$\zeta^{-1}$}}
\put(33,2){\tiny{$\zeta^{-1}$}}
\put(48.5,2){\tiny{$\zeta^{-1}$}}
\put(63.5,2){\tiny{$\zeta^{-1}$}}
\put(78.5,2){\tiny{$\zeta^{-1}$}}

\put(-2,-5){\tiny{$\al_1$}}
\put(13,-5){\tiny{$\al_2$}}
\put(28,-5){\tiny{$\al_3$}}
\put(43,-5){\tiny{$\al_4$}}
\put(58,-5){\tiny{$\al_5$}}
\put(73,-5){\tiny{$\al_6$}}
\put(88,-5){\tiny{$\al_7$}}
\end{picture}}

\newcommand{\smallHecSevenTwentyTwoTop}{
\setlength{\unitlength}{0.5mm}
\begin{picture}(30,20)(00,00)

\put(-5,16){\tiny{$\acchi^{(4,22)}_1$}}

\multiput(60.8,3.8)(00,00){1}{\oval(6.4,6.4)[br]}
\multiput(84.5,3.8)(00,00){1}{\oval(41,13)[t]}
\multiput(105,01)(00,00){1}{\line(0,1){2.8}}

\multiput(00,00)(15,00){8}{\circle{2}}
\multiput(01,00)(15,00){6}{\line(1,0){13}}

\put(-3,-15){\tiny{$(q\in\bKt,\,q\ne 1)$}}

\multiput(-1,3)(15,00){7}{\tiny{$q$}}

\put(107,2){\tiny{$q$}}

\multiput(4.5,2)(15,00){6}{\tiny{$q^{-1}$}}

\put(83,12){\tiny{$q^{-1}$}}

\put(-2,-5){\tiny{$\al_1$}}
\put(13,-5){\tiny{$\al_2$}}
\put(28,-5){\tiny{$\al_3$}}
\put(43,-5){\tiny{$\al_4$}}
\put(58,-5){\tiny{$\al_5$}}
\put(73,-5){\tiny{$\al_6$}}
\put(88,-5){\tiny{$\al_7$}}
\put(103,-5){\tiny{$\al_8$}}

\end{picture}}

\newcommand{\allHecFourEight}{
\setlength{\unitlength}{1mm}
\begin{picture}(220,30)(00,-10)
\put(10,10){$\smallHecFourEighta$}
\put(40,10){$\smallHecFourEightb$}
\put(70,10){$\smallHecFourEightc$}
\put(100,10){$\smallHecFourEightd$}
\put(05,00){$(q\in\bKt,\,q\ne -1)$}
\put(10,-10){\tiny{$e:=c\circ\left[{{1234}\atop{1243}}\right]$}}
\put(30,-10){\tiny{$f:=b\circ\left[{{1234}\atop{1243}}\right]$}}
\put(50,-10){\tiny{$g:=a\circ\left[{{1234}\atop{1243}}\right]$}}
\end{picture}}

\newcommand{\smallHecFourEighta}{
\setlength{\unitlength}{0.5mm}
\begin{picture}(30,20)(00,00)

\put(-5,16){\tiny{$a:=[\acchi^{(3,8)}_1]$}}

\multiput(00,00)(15,00){4}{\circle{2}}
\multiput(01,00)(15,00){3}{\line(1,0){13}}

\put(-4,3){\tiny{$-1$}}
\put(14,3){\tiny{$q$}}
\put(29,3){\tiny{$q$}}
\put(44,3){\tiny{$q^2$}}

\put(3.5,2){\tiny{$q^{-1}$}}
\put(18.5,2){\tiny{$q^{-1}$}}
\put(33.5,2){\tiny{$q^{-2}$}}

\put(-2,-5){\tiny{$\al_1$}}
\put(13,-5){\tiny{$\al_2$}}
\put(28,-5){\tiny{$\al_3$}}
\put(43,-5){\tiny{$\al_4$}}
\end{picture}}

\newcommand{\smallHecFourEightb}{
\setlength{\unitlength}{0.5mm}
\begin{picture}(30,20)(00,00)

\put(-5,16){\tiny{$b:=[\acchi^{(3,8)}_2]$}}

\multiput(00,00)(15,00){4}{\circle{2}}
\multiput(01,00)(15,00){3}{\line(1,0){13}}

\put(-4.5,3){\tiny{$-1$}}
\put(13,3){\tiny{$-1$}}
\put(29,3){\tiny{$q$}}
\put(44,3){\tiny{$q^2$}}

\put(4.5,2){\tiny{$q^{-1}$}}
\put(19.5,2){\tiny{$q^{-1}$}}
\put(33.5,2){\tiny{$q^{-2}$}}

\put(-2,-5){\tiny{$\al_1$}}
\put(13,-5){\tiny{$\al_2$}}
\put(28,-5){\tiny{$\al_3$}}
\put(43,-5){\tiny{$\al_4$}}
\end{picture}}

\newcommand{\smallHecFourEightc}{
\setlength{\unitlength}{0.5mm}
\begin{picture}(30,20)(00,00)

\put(-5,16){\tiny{$c:=[\acchi^{(3,8)}_3]$}}

\multiput(00,00)(15,00){4}{\circle{2}}
\multiput(01,00)(15,00){3}{\line(1,0){13}}

\put(-1,3){\tiny{$q$}}
\put(13,3){\tiny{$-1$}}
\put(26.5,3){\tiny{$-1$}}
\put(44,3){\tiny{$q^2$}}

\put(3.5,2){\tiny{$q^{-1}$}}
\put(21,2){\tiny{$q$}}
\put(33.5,2){\tiny{$q^{-2}$}}

\put(-2,-5){\tiny{$\al_1$}}
\put(13,-5){\tiny{$\al_2$}}
\put(28,-5){\tiny{$\al_3$}}
\put(43,-5){\tiny{$\al_4$}}
\end{picture}}

\newcommand{\smallHecFourEightd}{
\setlength{\unitlength}{0.5mm}
\begin{picture}(30,20)(00,00)

\put(-5,20){\tiny{$d:=[\acchi^{(3,8)}_4]$}}

\multiput(15.8,3.8)(00,00){1}{\oval(6.4,6.4)[br]}
\multiput(32,3.8)(00,00){1}{\oval(26,13)[t]}
\multiput(45,01)(00,00){1}{\line(0,1){2.8}}

\multiput(00,00)(15,00){4}{\circle{2}}
\multiput(01,00)(15,00){3}{\line(1,0){13}}

\put(-1,3){\tiny{$q$}}
\put(14,3){\tiny{$q$}}
\put(29,3){\tiny{$-1$}}
\put(46,2){\tiny{$-1$}}

\put(3.5,2){\tiny{$q^{-1}$}}
\put(20,2){\tiny{$q^{-1}$}}
\put(37,2){\tiny{$q^2$}}
\put(29,12){\tiny{$q^{-1}$}}

\put(-2,-5){\tiny{$\al_1$}}
\put(13,-5){\tiny{$\al_2$}}
\put(28,-5){\tiny{$\al_3$}}
\put(43,-5){\tiny{$\al_4$}}
\end{picture}}

\newcommand{\RankFourHecEightTau}{
\setlength{\unitlength}{0.5mm}
\begin{picture}(140,20)(0,0)

\put(00,00){\footnotesize{$a$}}
\put(23,00){\footnotesize{$b$}}
\put(46,00){\footnotesize{$c$}}
\put(69,00){\footnotesize{$d$}}
\put(92,00){\footnotesize{$e$}}
\put(115,00){\footnotesize{$f$}}
\put(138,00){\footnotesize{$g$}}

\multiput(6.5,00)(23,0){6}{\vector(1,0){13}}
\multiput(18.5,00)(23,0){6}{\vector(-1,0){13}}

\put(10,03){\footnotesize{${\bar{\tau}}_1$}}
\put(33,03){\footnotesize{${\bar{\tau}}_2$}}
\put(56,03){\footnotesize{${\bar{\tau}}_3$}}
\put(79,03){\footnotesize{${\bar{\tau}}_4$}}
\put(102,03){\footnotesize{${\bar{\tau}}_2$}}
\put(125,03){\footnotesize{${\bar{\tau}}_1$}}

\end{picture}}

\newcommand{\allHecFourNine}{
\setlength{\unitlength}{1mm}
\begin{picture}(220,40)(00,-5)
\put(10,20){$\smallHecFourNinea$}
\put(40,20){$\smallHecFourNineb$}
\put(70,20){$\smallHecFourNinec$}
\put(100,20){$\smallHecFourNined$}
\put(10,00){$\smallHecFourNinee$}
\put(40,00){$\smallHecFourNinef$}
\put(75,00){$(q\in\bKt,\,q^2\ne 1\ne q^3)$}
\end{picture}}

\newcommand{\smallHecFourNinea}{
\setlength{\unitlength}{0.5mm}
\begin{picture}(30,20)(00,00)

\put(-5,16){\tiny{$a:=[\acchi^{(3,9)}_1]\circ\left[{{1234}\atop{3421}}\right]$}}

\multiput(00,00)(15,00){4}{\circle{2}}
\multiput(01,00)(15,00){3}{\line(1,0){13}}

\put(-1,3){\tiny{$q^2$}}
\put(14,3){\tiny{$q^2$}}
\put(29,3){\tiny{$q$}}
\put(43,3){\tiny{$-1$}}

\put(3.5,2){\tiny{$q^{-2}$}}
\put(21,2){\tiny{$q^2$}}
\put(33.5,2){\tiny{$q^{-1}$}}

\put(-2,-5){\tiny{$\al_4$}}
\put(13,-5){\tiny{$\al_3$}}
\put(28,-5){\tiny{$\al_1$}}
\put(43,-5){\tiny{$\al_2$}}
\end{picture}}

\newcommand{\smallHecFourNineb}{
\setlength{\unitlength}{0.5mm}
\begin{picture}(30,20)(00,00)

\put(-5,16){\tiny{$b:=[\acchi^{(3,9)}_2]\circ\left[{{1234}\atop{3421}}\right]$}}

\multiput(00,00)(15,00){4}{\circle{2}}
\multiput(01,00)(15,00){3}{\line(1,0){13}}

\put(-1,3){\tiny{$q^2$}}
\put(14,3){\tiny{$q^2$}}
\put(28,3){\tiny{$-1$}}
\put(43,3){\tiny{$-1$}}

\put(3.5,2){\tiny{$q^{-2}$}}
\put(21,2){\tiny{$q^2$}}
\put(36,2){\tiny{$q$}}

\put(-2,-5){\tiny{$\al_4$}}
\put(13,-5){\tiny{$\al_3$}}
\put(28,-5){\tiny{$\al_1$}}
\put(43,-5){\tiny{$\al_2$}}
\end{picture}}

\newcommand{\smallHecFourNinec}{
\setlength{\unitlength}{0.5mm}
\begin{picture}(30,20)(00,00)

\put(-5,20){\tiny{$c:=[\acchi^{(3,9)}_3]\circ\left[{{1234}\atop{3421}}\right]$}}

\multiput(15.8,3.8)(00,00){1}{\oval(6.4,6.4)[br]}
\multiput(32,3.8)(00,00){1}{\oval(26,13)[t]}
\multiput(45,01)(00,00){1}{\line(0,1){2.8}}

\multiput(00,00)(15,00){4}{\circle{2}}
\multiput(01,00)(15,00){3}{\line(1,0){13}}

\put(-3,3){\tiny{$q^2$}}
\put(11.5,3){\tiny{$-1$}}
\put(27,3){\tiny{$-1$}}
\put(46,2){\tiny{$q$}}

\put(2,2){\tiny{$q^{-2}$}}
\put(20,2){\tiny{$q^2$}}
\put(35,2){\tiny{$q^{-1}$}}
\put(29,12){\tiny{$q^{-1}$}}

\put(-2,-5){\tiny{$\al_4$}}
\put(13,-5){\tiny{$\al_3$}}
\put(28,-5){\tiny{$\al_1$}}
\put(43,-5){\tiny{$\al_2$}}
\end{picture}}

\newcommand{\smallHecFourNined}{
\setlength{\unitlength}{0.5mm}
\begin{picture}(30,20)(00,00)

\put(-5,20){\tiny{$d:=[\acchi^{(3,9)}_4]\circ\left[{{1234}\atop{1423}}\right]$}}

\multiput(15.8,3.8)(00,00){1}{\oval(6.4,6.4)[br]}
\multiput(32,3.8)(00,00){1}{\oval(26,13)[t]}
\multiput(45,01)(00,00){1}{\line(0,1){2.8}}

\multiput(00,00)(15,00){4}{\circle{2}}
\multiput(01,00)(15,00){3}{\line(1,0){13}}

\put(-3,3){\tiny{$q^2$}}
\put(11.5,3){\tiny{$-1$}}
\put(27,3){\tiny{$-1$}}
\put(46,2){\tiny{$-1$}}

\put(2,2){\tiny{$q^{-2}$}}
\put(20,2){\tiny{$q^2$}}
\put(35,2){\tiny{$q^{-3}$}}
\put(29,12){\tiny{$q$}}

\put(-2,-5){\tiny{$\al_1$}}
\put(13,-5){\tiny{$\al_3$}}
\put(28,-5){\tiny{$\al_4$}}
\put(43,-5){\tiny{$\al_2$}}
\end{picture}}

\newcommand{\smallHecFourNinee}{
\setlength{\unitlength}{0.5mm}
\begin{picture}(30,20)(00,00)

\put(-5,16){\tiny{$e:=[\acchi^{(3,9)}_5]\circ\left[{{1234}\atop{1423}}\right]$}}

\multiput(00,00)(15,00){4}{\circle{2}}
\multiput(01,00)(15,00){3}{\line(1,0){13}}

\put(-4,3){\tiny{$q^{-3}$}}
\put(12,3){\tiny{$-1$}}
\put(29,3){\tiny{$q^2$}}
\put(44,3){\tiny{$q^2$}}

\put(6,2){\tiny{$q^3$}}
\put(18.5,2){\tiny{$q^{-2}$}}
\put(34,2){\tiny{$q^{-2}$}}

\put(-2,-5){\tiny{$\al_2$}}
\put(13,-5){\tiny{$\al_4$}}
\put(28,-5){\tiny{$\al_3$}}
\put(43,-5){\tiny{$\al_1$}}

\end{picture}}

\newcommand{\smallHecFourNinef}{
\setlength{\unitlength}{0.5mm}
\begin{picture}(30,20)(00,00)

\put(-5,16){\tiny{$f:=[\acchi^{(3,9)}_6]\circ\left[{{1234}\atop{1342}}\right]$}}

\multiput(00,00)(15,00){4}{\circle{2}}
\multiput(01,00)(15,00){3}{\line(1,0){13}}

\put(-4,3){\tiny{$q^{-3}$}}
\put(12,3){\tiny{$-1$}}
\put(29,3){\tiny{$q$}}
\put(44,3){\tiny{$q^2$}}

\put(6,2){\tiny{$q^3$}}
\put(18.5,2){\tiny{$q^{-1}$}}
\put(33.5,2){\tiny{$q^{-2}$}}

\put(-2,-5){\tiny{$\al_1$}}
\put(13,-5){\tiny{$\al_3$}}
\put(28,-5){\tiny{$\al_2$}}
\put(43,-5){\tiny{$\al_4$}}

\end{picture}}

\newcommand{\RankFourHecNineTau}{
\setlength{\unitlength}{0.5mm}
\begin{picture}(140,30)(0,0)

\put(00,23){\footnotesize{$a$}}
\put(23,23){\footnotesize{$b$}}
\put(46,23){\footnotesize{$c$}}
\put(69,23){\footnotesize{$d$}}
\put(92,23){\footnotesize{$e$}}

\multiput(6.5,23)(23,0){4}{\vector(1,0){13}}
\multiput(18.5,23)(23,0){4}{\vector(-1,0){13}}

\multiput(70,20)(00,0){1}{\vector(0,-1){13}}
\multiput(70,8)(00,0){1}{\vector(0,1){13}}

\put(69,00){\footnotesize{$f$}}

\put(10,26){\footnotesize{${\bar{\tau}}_2$}}
\put(33,26){\footnotesize{${\bar{\tau}}_1$}}
\put(56,26){\footnotesize{${\bar{\tau}}_3$}}
\put(79,26){\footnotesize{${\bar{\tau}}_4$}}

\put(63,12){\footnotesize{${\bar{\tau}}_2$}}

\end{picture}}

\newcommand{\allHecFourTwelve}{
\setlength{\unitlength}{1mm}
\begin{picture}(220,40)(00,-10)
\put(10,20){$\smallHecFourTwelvea$}
\put(40,20){$\smallHecFourTwelveb$}
\put(70,20){$\smallHecFourTwelvec$}
\put(100,20){$\smallHecFourTwelved$}
\put(10,00){$\smallHecFourTwelvee$}
\put(40,00){$\smallHecFourTwelvef$}
\put(75,00){$(q\in\bKt,\,q\ne -1)$}
\put(10,-10){\tiny{$g:=c\circ\left[{{1234}\atop{1243}}\right]$}}
\put(30,-10){\tiny{$h:=b\circ\left[{{1234}\atop{1243}}\right]$}}
\put(50,-10){\tiny{$i:=a\circ\left[{{1234}\atop{1243}}\right]$}}
\end{picture}}

\newcommand{\smallHecFourTwelvea}{
\setlength{\unitlength}{0.5mm}
\begin{picture}(30,20)(00,00)

\put(-5,16){\tiny{$a:=[\acchi^{(3,12)}_1]$}}

\multiput(00,00)(15,00){4}{\circle{2}}
\multiput(01,00)(15,00){3}{\line(1,0){13}}

\put(-3,3){\tiny{$q^{-1}$}}
\put(11.5,3){\tiny{$-1$}}
\put(29,3){\tiny{$q$}}
\put(44,3){\tiny{$q^2$}}

\put(6,2){\tiny{$q$}}
\put(18.5,2){\tiny{$q^{-1}$}}
\put(33.5,2){\tiny{$q^{-2}$}}

\put(-2,-5){\tiny{$\al_1$}}
\put(13,-5){\tiny{$\al_2$}}
\put(28,-5){\tiny{$\al_3$}}
\put(43,-5){\tiny{$\al_4$}}
\end{picture}}

\newcommand{\smallHecFourTwelveb}{
\setlength{\unitlength}{0.5mm}
\begin{picture}(30,20)(00,00)

\put(-5,16){\tiny{$b:=[\acchi^{(3,12)}_2]$}}

\multiput(00,00)(15,00){4}{\circle{2}}
\multiput(01,00)(15,00){3}{\line(1,0){13}}

\put(-4.5,3){\tiny{$-1$}}
\put(13,3){\tiny{$-1$}}
\put(26.5,3){\tiny{$-1$}}
\put(44,3){\tiny{$q^2$}}

\put(3.5,2){\tiny{$q^{-1}$}}
\put(21,2){\tiny{$q$}}
\put(33.5,2){\tiny{$q^{-2}$}}

\put(-2,-5){\tiny{$\al_1$}}
\put(13,-5){\tiny{$\al_2$}}
\put(28,-5){\tiny{$\al_3$}}
\put(43,-5){\tiny{$\al_4$}}
\end{picture}}

\newcommand{\smallHecFourTwelvec}{
\setlength{\unitlength}{0.5mm}
\begin{picture}(30,20)(00,00)

\put(-5,16){\tiny{$c:=[\acchi^{(3,12)}_3]$}}

\multiput(00,00)(15,00){4}{\circle{2}}
\multiput(01,00)(15,00){3}{\line(1,0){13}}

\put(-4.5,3){\tiny{$-1$}}
\put(13,3){\tiny{$q^{-1}$}}
\put(26.5,3){\tiny{$-1$}}
\put(44,3){\tiny{$q^2$}}

\put(6,2){\tiny{$q$}}
\put(21,2){\tiny{$q$}}
\put(33.5,2){\tiny{$q^{-2}$}}

\put(-2,-5){\tiny{$\al_1$}}
\put(13,-5){\tiny{$\al_2$}}
\put(28,-5){\tiny{$\al_3$}}
\put(43,-5){\tiny{$\al_4$}}
\end{picture}}

\newcommand{\smallHecFourTwelved}{
\setlength{\unitlength}{0.5mm}
\begin{picture}(30,20)(00,00)

\put(-5,20){\tiny{$d:=[\acchi^{(3,12)}_4]$}}

\multiput(15.8,3.8)(00,00){1}{\oval(6.4,6.4)[br]}
\multiput(32,3.8)(00,00){1}{\oval(26,13)[t]}
\multiput(45,01)(00,00){1}{\line(0,1){2.8}}

\multiput(00,00)(15,00){4}{\circle{2}}
\multiput(01,00)(15,00){3}{\line(1,0){13}}

\put(-3,3){\tiny{$-1$}}
\put(14,3){\tiny{$q$}}
\put(29,3){\tiny{$-1$}}
\put(46,2){\tiny{$-1$}}

\put(5,2){\tiny{$q^{-1}$}}
\put(20,2){\tiny{$q^{-1}$}}
\put(37,2){\tiny{$q^2$}}
\put(29,12){\tiny{$q^{-1}$}}

\put(-2,-5){\tiny{$\al_1$}}
\put(13,-5){\tiny{$\al_2$}}
\put(28,-5){\tiny{$\al_3$}}
\put(43,-5){\tiny{$\al_4$}}
\end{picture}}

\newcommand{\smallHecFourTwelvee}{
\setlength{\unitlength}{0.5mm}
\begin{picture}(30,20)(00,00)

\put(-5,20){\tiny{$e:=[\acchi^{(3,12)}_5]$}}

\multiput(15.8,3.8)(00,00){1}{\oval(6.4,6.4)[br]}
\multiput(32,3.8)(00,00){1}{\oval(26,13)[t]}
\multiput(45,01)(00,00){1}{\line(0,1){2.8}}

\multiput(00,00)(15,00){4}{\circle{2}}
\multiput(01,00)(15,00){3}{\line(1,0){13}}

\put(-3,3){\tiny{$-1$}}
\put(11,3){\tiny{$-1$}}
\put(29,3){\tiny{$-1$}}
\put(46,2){\tiny{$-1$}}

\put(6,2){\tiny{$q$}}
\put(20,2){\tiny{$q^{-1}$}}
\put(37,2){\tiny{$q^2$}}
\put(29,12){\tiny{$q^{-1}$}}

\put(-2,-5){\tiny{$\al_1$}}
\put(13,-5){\tiny{$\al_2$}}
\put(28,-5){\tiny{$\al_3$}}
\put(43,-5){\tiny{$\al_4$}}
\end{picture}}

\newcommand{\smallHecFourTwelvef}{
\setlength{\unitlength}{0.5mm}
\begin{picture}(30,20)(00,00)

\put(-5,20){\tiny{$f:=[\acchi^{(3,12)}_6]$}}

\multiput(15.8,3.8)(00,00){1}{\oval(6.4,6.4)[br]}
\multiput(32,3.8)(00,00){1}{\oval(26,13)[t]}
\multiput(45,01)(00,00){1}{\line(0,1){2.8}}

\multiput(00,00)(15,00){4}{\circle{2}}
\multiput(01,00)(15,00){2}{\line(1,0){13}}

\put(-1,3){\tiny{$q$}}
\put(12,3){\tiny{$-1$}}
\put(29,3){\tiny{$q^{-1}$}}
\put(46,2){\tiny{$q^{-1}$}}

\put(3,2){\tiny{$q^{-1}$}}
\put(22,2){\tiny{$q$}}
\put(29,12){\tiny{$q$}}

\put(-2,-5){\tiny{$\al_1$}}
\put(13,-5){\tiny{$\al_2$}}
\put(28,-5){\tiny{$\al_3$}}
\put(43,-5){\tiny{$\al_4$}}
\end{picture}}

\newcommand{\RankFourHecTwelveTau}{
\setlength{\unitlength}{0.5mm}
\begin{picture}(140,60)(0,0)

\put(00,46){\footnotesize{$a$}}
\put(23,46){\footnotesize{$b$}}
\put(46,46){\footnotesize{$c$}}
\multiput(6.5,47)(23,0){2}{\vector(1,0){13}}
\multiput(18.5,47)(23,0){2}{\vector(-1,0){13}}

\multiput(24,43)(23,0){2}{\vector(0,-1){13}}
\multiput(24,31)(23,0){2}{\vector(0,1){13}}

\put(23,23){\footnotesize{$d$}}
\put(46,23){\footnotesize{$e$}}
\put(69,23){\footnotesize{$f$}}

\multiput(29.5,24)(23,0){2}{\vector(1,0){13}}
\multiput(41.5,24)(23,0){2}{\vector(-1,0){13}}

\multiput(24,20)(23,0){2}{\vector(0,-1){13}}
\multiput(24,8)(23,0){2}{\vector(0,1){13}}

\put(00,00){\footnotesize{$i$}}
\put(23,00){\footnotesize{$h$}}
\put(46,00){\footnotesize{$g$}}
\multiput(6.5,1)(23,0){2}{\vector(1,0){13}}
\multiput(18.5,1)(23,0){2}{\vector(-1,0){13}}

\put(10,49){\footnotesize{${\bar{\tau}}_2$}}
\put(33,49){\footnotesize{${\bar{\tau}}_1$}}

\put(17,35){\footnotesize{${\bar{\tau}}_3$}}
\put(49,35){\footnotesize{${\bar{\tau}}_3$}}

\put(33,26){\footnotesize{${\bar{\tau}}_1$}}
\put(56,26){\footnotesize{${\bar{\tau}}_2$}}

\put(17,12){\footnotesize{${\bar{\tau}}_4$}}
\put(48,12){\footnotesize{${\bar{\tau}}_4$}}

\put(10,03){\footnotesize{${\bar{\tau}}_2$}}
\put(33,03){\footnotesize{${\bar{\tau}}_1$}}

\end{picture}}

\newcommand{\smallHecFourThirteenOne}{
\setlength{\unitlength}{0.5mm}
\begin{picture}(30,20)(00,00)

\put(-5,16){\tiny{$a:=[\acchi^{(3,13)}_1]$}}

\multiput(00,00)(15,00){4}{\circle{2}}
\multiput(01,00)(15,00){3}{\line(1,0){13}}

\put(-1,3){\tiny{$q$}}
\put(14,3){\tiny{$q$}}
\put(28,3){\tiny{$-1$}}
\put(44,3){\tiny{$q^{-2}$}}

\put(6,2){\tiny{$q$}}
\put(18.5,2){\tiny{$q^{-1}$}}
\put(36,2){\tiny{$q^2$}}

\put(-2,-5){\tiny{$\al_1$}}
\put(13,-5){\tiny{$\al_2$}}
\put(28,-5){\tiny{$\al_3$}}
\put(43,-5){\tiny{$\al_4$}}
\end{picture}}

\newcommand{\smallHecFourThirteenOnedash}{
\setlength{\unitlength}{0.5mm}
\begin{picture}(30,20)(00,00)

\put(-5,16){\tiny{$c:=[\acchi^{(3,13)}_1]\circ\left[{{1234}\atop{1243}}\right]$}}

\multiput(00,00)(15,00){4}{\circle{2}}
\multiput(01,00)(15,00){3}{\line(1,0){13}}

\put(-1,3){\tiny{$q$}}
\put(14,3){\tiny{$q$}}
\put(28,3){\tiny{$-1$}}
\put(44,3){\tiny{$q^{-2}$}}

\put(6,2){\tiny{$q$}}
\put(18.5,2){\tiny{$q^{-1}$}}
\put(36,2){\tiny{$q^2$}}

\put(-2,-5){\tiny{$\al_1$}}
\put(13,-5){\tiny{$\al_2$}}
\put(28,-5){\tiny{$\al_4$}}
\put(43,-5){\tiny{$\al_3$}}
\end{picture}}

\newcommand{\smallHecFourThirteenTwo}{
\setlength{\unitlength}{0.5mm}
\begin{picture}(30,20)(00,00)

\put(-5,16){\tiny{$b:=[\acchi^{(3,13)}_2]$}}

\multiput(15.8,3.8)(00,00){1}{\oval(6.4,6.4)[br]}
\multiput(32,3.8)(00,00){1}{\oval(26,13)[t]}
\multiput(45,01)(00,00){1}{\line(0,1){2.8}}

\multiput(00,00)(15,00){4}{\circle{2}}
\multiput(01,00)(15,00){3}{\line(1,0){13}}

\put(-1,3){\tiny{$q$}}
\put(11.5,3){\tiny{$-1$}}
\put(28,3){\tiny{$-1$}}
\put(46,2){\tiny{$-1$}}

\put(3,2){\tiny{$q^{-1}$}}
\put(19.5,2){\tiny{$q^{-1}$}}
\put(35,2){\tiny{$q^{-2}$}}
\put(31,12){\tiny{$q^{-1}$}}

\put(-2,-5){\tiny{$\al_1$}}
\put(13,-5){\tiny{$\al_2$}}
\put(28,-5){\tiny{$\al_3$}}
\put(43,-5){\tiny{$\al_4$}}
\end{picture}}

\newcommand{\smallHecFourThirteenThree}{
\setlength{\unitlength}{0.5mm}
\begin{picture}(30,20)(00,00)

\put(-5,16){\tiny{$d:=[\acchi^{(3,13)}_3]$}}

\multiput(15.8,3.8)(00,00){1}{\oval(6.4,6.4)[br]}
\multiput(32,3.8)(00,00){1}{\oval(26,13)[t]}
\multiput(45,01)(00,00){1}{\line(0,1){2.8}}

\multiput(00,00)(15,00){4}{\circle{2}}
\multiput(01,00)(15,00){2}{\line(1,0){13}}

\put(-2,3){\tiny{$-1$}}
\put(11.5,3){\tiny{$-1$}}
\put(30,3){\tiny{$q$}}
\put(46,2){\tiny{$q$}}

\put(6,2){\tiny{$q$}}
\put(20.5,2){\tiny{$q^{-1}$}}
\put(31,12){\tiny{$q^{-1}$}}

\put(-2,-5){\tiny{$\al_1$}}
\put(13,-5){\tiny{$\al_2$}}
\put(28,-5){\tiny{$\al_3$}}
\put(43,-5){\tiny{$\al_4$}}
\end{picture}}

\newcommand{\smallHecFourThirteenFour}{
\setlength{\unitlength}{0.5mm}
\begin{picture}(30,20)(00,00)

\put(-5,16){\tiny{$e:=[\acchi^{(3,13)}_4]$}}

\multiput(15.8,3.8)(00,00){1}{\oval(6.4,6.4)[br]}
\multiput(32,3.8)(00,00){1}{\oval(26,13)[t]}
\multiput(45,01)(00,00){1}{\line(0,1){2.8}}

\multiput(00,00)(15,00){4}{\circle{2}}
\multiput(01,00)(15,00){2}{\line(1,0){13}}

\put(-4,3){\tiny{$-1$}}
\put(13,3){\tiny{$q$}}
\put(30,3){\tiny{$q$}}
\put(46,2){\tiny{$q$}}

\put(3,2){\tiny{$q^{-1}$}}
\put(20.5,2){\tiny{$q^{-1}$}}
\put(31,12){\tiny{$q^{-1}$}}

\put(-2,-5){\tiny{$\al_1$}}
\put(13,-5){\tiny{$\al_2$}}
\put(28,-5){\tiny{$\al_3$}}
\put(43,-5){\tiny{$\al_4$}}
\end{picture}}

\newcommand{\allHecFourThirteen}{
\setlength{\unitlength}{1mm}
\begin{picture}(220,40)(00,-10)
\put(10,20){$\smallHecFourThirteenOne$}
\put(10,00){$\smallHecFourThirteenOnedash$}
\put(40,20){$\smallHecFourThirteenTwo$}
\put(70,20){$\smallHecFourThirteenThree$}
\put(100,20){$\smallHecFourThirteenFour$}

\put(90,00){$(q\in\bKt,\,q^2\ne 1)$}

\end{picture}}

\newcommand{\RankForHecFourThirteenTau}{
\setlength{\unitlength}{1mm}
\begin{picture}(50,30)(40,-10)
\put(50,7){$a$}\put(65,7){$b$}\put(80,7){$d$}\put(95,7){$e$}
\put(65,-8){$c$}
\multiput(55.5,8)(15,00){3}{\vector(1,0){8}}
\multiput(61.5,8)(15,00){3}{\vector(-1,0){8}}
\multiput(66,3.5)(00,00){1}{\vector(0,-1){8}}
\multiput(66,-2.5)(00,00){1}{\vector(0,1){8}}
\put(56.5,10){$\chambartau_3$}
\put(71.5,10){$\chambartau_2$}
\put(86.5,10){$\chambartau_1$}
\put(61.5,-1){$\chambartau_4$}
\end{picture}}

\newcommand{\allHecFourFourteen}{
\setlength{\unitlength}{1mm}
\begin{picture}(220,60)(00,-20)
\put(10,20){$\smallHecFourFourteena$}
\put(40,20){$\smallHecFourFourteenb$}
\put(70,20){$\smallHecFourFourteenf$}
\put(100,20){$\smallHecFourFourteend$}
\put(10,00){$\smallHecFourFourteenj$}
\put(55,00){$(q\in\bKt,\,q^2\ne 1)$}
\put(10,-10){\tiny{$c:=b\circ\left[{{1234}\atop{3214}}\right]$}}
\put(30,-10){\tiny{$e:=d\circ\left[{{1234}\atop{3214}}\right]$}}
\put(50,-10){\tiny{$h:=a\circ\left[{{1234}\atop{3214}}\right]$}}
\put(70,-10){\tiny{$g:=f\circ\left[{{1234}\atop{3214}}\right]$}}
\put(90,-10){\tiny{$i:=j\circ\left[{{1234}\atop{3214}}\right]$}}
\end{picture}}

\newcommand{\smallHecFourFourteena}{
\setlength{\unitlength}{0.5mm}
\begin{picture}(30,20)(00,00)

\put(-5,16){\tiny{$a:=[\acchi^{(3,14)}_1]$}}

\multiput(00,00)(15,00){4}{\circle{2}}
\multiput(01,00)(15,00){3}{\line(1,0){13}}

\put(-1,3){\tiny{$q$}}
\put(14,3){\tiny{$q$}}
\put(28,3){\tiny{$-1$}}
\put(40,-5){\tiny{$-q^{-1}$}}

\put(3.5,2){\tiny{$q^{-1}$}}
\put(18.5,2){\tiny{$q^{-1}$}}
\put(35,2){\tiny{$-q$}}

\put(-2,-5){\tiny{$\al_1$}}
\put(13,-5){\tiny{$\al_2$}}
\put(28,-5){\tiny{$\al_3$}}
\put(43,3){\tiny{$\al_4$}}
\end{picture}}

\newcommand{\smallHecFourFourteenb}{
\setlength{\unitlength}{0.5mm}
\begin{picture}(30,20)(00,00)

\put(-5,16){\tiny{$b:=[\acchi^{(3,14)}_2]$}}

\multiput(15.8,3.8)(00,00){1}{\oval(6.4,6.4)[br]}
\multiput(32,3.8)(00,00){1}{\oval(26,13)[t]}
\multiput(45,01)(00,00){1}{\line(0,1){2.8}}

\multiput(00,00)(15,00){4}{\circle{2}}
\multiput(01,00)(15,00){3}{\line(1,0){13}}

\put(-1,3){\tiny{$q$}}
\put(11.5,3){\tiny{$-1$}}
\put(26,4){\tiny{$-1$}}
\put(46,2){\tiny{$-1$}}

\put(2.5,2){\tiny{$q^{-1}$}}
\put(20,1){\tiny{$-1$}}
\put(31.5,1.5){\tiny{$-q^{-1}$}}
\put(31,12){\tiny{$q$}}

\put(-2,-5){\tiny{$\al_1$}}
\put(13,-5){\tiny{$\al_2$}}
\put(28,-5){\tiny{$\al_3$}}
\put(43,-5){\tiny{$\al_4$}}
\end{picture}}

\newcommand{\smallHecFourFourteenf}{
\setlength{\unitlength}{0.5mm}
\begin{picture}(30,20)(00,00)

\put(-5,20){\tiny{$f:=[\acchi^{(3,14)}_3]\circ\left[{{1234}\atop{3412}}\right]$}}

\multiput(15.8,3.8)(00,00){1}{\oval(6.4,6.4)[br]}
\multiput(32,3.8)(00,00){1}{\oval(26,13)[t]}
\multiput(45,01)(00,00){1}{\line(0,1){2.8}}

\multiput(00,00)(15,00){4}{\circle{2}}
\multiput(01,00)(15,00){3}{\line(1,0){13}}

\put(-3,3){\tiny{$q^{-1}$}}
\put(11.5,3){\tiny{$-1$}}
\put(27,3){\tiny{$-1$}}
\put(46,2){\tiny{$-1$}}

\put(3.5,2){\tiny{$-q$}}
\put(20,1){\tiny{$-1$}}
\put(37,2){\tiny{$q$}}
\put(27,12){\tiny{$-q^{-1}$}}

\put(-2,-5){\tiny{$\al_3$}}
\put(13,-5){\tiny{$\al_4$}}
\put(28,-5){\tiny{$\al_2$}}
\put(43,-5){\tiny{$\al_1$}}
\end{picture}}

\newcommand{\smallHecFourFourteend}{
\setlength{\unitlength}{0.5mm}
\begin{picture}(30,20)(00,00)

\put(-5,16){\tiny{$d:=[\acchi^{(3,14)}_4]\circ\left[{{1234}\atop{1243}}\right]$}}

\multiput(00,00)(15,00){4}{\circle{2}}
\multiput(01,00)(15,00){3}{\line(1,0){13}}

\put(-1,3){\tiny{$q$}}
\put(13,3){\tiny{$-1$}}
\put(28,3){\tiny{$-1$}}
\put(40,-5){\tiny{$-q^{-1}$}}

\put(3.5,2){\tiny{$q^{-1}$}}
\put(20,2){\tiny{$-1$}}
\put(35,2){\tiny{$-q$}}

\put(-2,-5){\tiny{$\al_1$}}
\put(13,-5){\tiny{$\al_2$}}
\put(28,-5){\tiny{$\al_4$}}
\put(43,3){\tiny{$\al_3$}}
\end{picture}}

\newcommand{\smallHecFourFourteenj}{
\setlength{\unitlength}{0.5mm}
\begin{picture}(30,20)(00,00)

\put(-5,16){\tiny{$j:=[\acchi^{(3,14)}_5]\circ\left[{{1234}\atop{3412}}\right]$}}

\multiput(00,00)(15,00){4}{\circle{2}}
\multiput(01,00)(15,00){3}{\line(1,0){13}}

\put(-1,3){\tiny{$q$}}
\put(12,-5){\tiny{$-1$}}
\put(25,-6){\tiny{$-q^{-1}$}}
\put(41,-6){\tiny{$-q^{-1}$}}

\put(3.5,2){\tiny{$q^{-1}$}}
\put(19.5,2){\tiny{$-q$}}
\put(35,2){\tiny{$-q$}}

\put(-2,-5){\tiny{$\al_2$}}
\put(13,3){\tiny{$\al_1$}}
\put(28,3){\tiny{$\al_4$}}
\put(43,3){\tiny{$\al_3$}}
\end{picture}}

\newcommand{\RankFourHecFourteenTau}{
\setlength{\unitlength}{0.5mm}
\begin{picture}(140,30)(0,0)

\put(23,23){\footnotesize{$a$}}
\put(46,23){\footnotesize{$b$}}
\put(69,23){\footnotesize{$d$}}
\put(92,23){\footnotesize{$f$}}
\put(115,23){\footnotesize{$j$}}

\multiput(29.5,24)(23,0){4}{\vector(1,0){13}}
\multiput(41.5,24)(23,0){4}{\vector(-1,0){13}}

\multiput(47,20)(46,0){2}{\vector(0,-1){13}}
\multiput(47,8)(46,0){2}{\vector(0,1){13}}

\put(23,00){\footnotesize{$h$}}
\put(46,00){\footnotesize{$c$}}
\put(69,00){\footnotesize{$e$}}
\put(92,00){\footnotesize{$g$}}
\put(115,00){\footnotesize{$i$}}

\multiput(29.5,1)(23,0){4}{\vector(1,0){13}}
\multiput(41.5,1)(23,0){4}{\vector(-1,0){13}}

\put(33,26){\footnotesize{${\bar{\tau}}_3$}}
\put(56,26){\footnotesize{${\bar{\tau}}_4$}}
\put(79,26){\footnotesize{${\bar{\tau}}_2$}}
\put(102,26){\footnotesize{${\bar{\tau}}_1$}}

\put(40,12){\footnotesize{${\bar{\tau}}_4$}}
\multiput(94,12)(00,0){1}{\footnotesize{${\bar{\tau}}_4$}}

\put(33,03){\footnotesize{${\bar{\tau}}_1$}}
\put(56,03){\footnotesize{${\bar{\tau}}_4$}}
\put(79,03){\footnotesize{${\bar{\tau}}_2$}}
\put(102,03){\footnotesize{${\bar{\tau}}_3$}}

\end{picture}}

\newcommand{\allHecFourSeventeen}{
\setlength{\unitlength}{1mm}
\begin{picture}(220,60)(00,-20)
\put(10,20){$\smallHecFourSeventeenZero$}
\put(40,20){$\smallHecFourSeventeenOne$}
\put(70,20){$\smallHecFourSeventeenTwo$}
\put(100,20){$\smallHecFourSeventeenThree$}
\put(10,00){$\smallHecFourSeventeenFour$}
\put(40,00){$\smallHecFourSeventeenEleven$}
\put(75,00){$(\zeta\in\bKt,\,\zeta^2+\zeta+1=0)$}
\put(10,-10){\tiny{$h:=a\circ\left[{{1234}\atop{3214}}\right]$}}
\put(30,-10){\tiny{$c:=a\circ\left[{{1234}\atop{1432}}\right]$}}
\put(50,-10){\tiny{$i:=f\circ\left[{{1234}\atop{1432}}\right]$}}
\put(70,-10){\tiny{$d:=e\circ\left[{{1234}\atop{3214}}\right]$}}
\put(90,-10){\tiny{$k:=f\circ\left[{{1234}\atop{3214}}\right]$}}
\put(110,-10){\tiny{$j:=g\circ\left[{{1234}\atop{1432}}\right]$}}
\end{picture}}

\newcommand{\smallHecFourSeventeenZero}{
\setlength{\unitlength}{0.5mm}
\begin{picture}(30,20)(00,00)

\put(-5,16){\tiny{$g:=[\acchi^{(3,17)}_1]$}}

\multiput(00,00)(15,00){4}{\circle{2}}
\multiput(01,00)(15,00){3}{\line(1,0){13}}

\put(-6,3){\tiny{$-\zeta$}}
\put(11,3){\tiny{$-\zeta$}}
\put(27,3){\tiny{$-\zeta$}}
\put(44,3){\tiny{$\zeta$}}

\put(2,3.5){\tiny{$-{\frac 1 \zeta}$}}
\put(18.5,3.5){\tiny{$-{\frac 1 \zeta}$}}
\put(35,3.5){\tiny{$-{\frac 1 \zeta}$}}

\put(-2,-5){\tiny{$\al_1$}}
\put(13,-5){\tiny{$\al_2$}}
\put(28,-5){\tiny{$\al_3$}}
\put(43,-5){\tiny{$\al_4$}}
\end{picture}}

\newcommand{\smallHecFourSeventeenOne}{
\setlength{\unitlength}{0.5mm}
\begin{picture}(30,20)(00,00)

\put(-5,16){\tiny{$f:=[\acchi^{(3,17)}_2]$}}

\multiput(00,00)(15,00){4}{\circle{2}}
\multiput(01,00)(15,00){3}{\line(1,0){13}}

\put(-6,3){\tiny{$-\zeta$}}
\put(11,3){\tiny{$-\zeta$}}
\put(27,3){\tiny{$-1$}}
\put(44,3){\tiny{$\zeta$}}

\put(2,3.5){\tiny{$-{\frac 1 \zeta}$}}
\put(18.5,3.5){\tiny{$-{\frac 1 \zeta}$}}
\put(35,2){\tiny{$-1$}}

\put(-2,-5){\tiny{$\al_1$}}
\put(13,-5){\tiny{$\al_2$}}
\put(28,-5){\tiny{$\al_3$}}
\put(43,-5){\tiny{$\al_4$}}
\end{picture}}

\newcommand{\smallHecFourSeventeenTwo}{
\setlength{\unitlength}{0.5mm}
\begin{picture}(30,20)(00,00)

\put(-5,16){\tiny{$b:=[\acchi^{(3,17)}_3]$}}

\multiput(15.8,3.8)(00,00){1}{\oval(6.4,6.4)[br]}
\multiput(32,3.8)(00,00){1}{\oval(26,13)[t]}
\multiput(45,01)(00,00){1}{\line(0,1){2.8}}

\multiput(00,00)(15,00){4}{\circle{2}}
\multiput(01,00)(15,00){3}{\line(1,0){13}}

\put(-6,3){\tiny{$-\zeta$}}
\put(11.5,3){\tiny{$-1$}}
\put(28,3){\tiny{$-1$}}
\put(46,2){\tiny{$\zeta$}}

\put(2,3.5){\tiny{$-{\frac 1 \zeta}$}}
\put(19.5,2){\tiny{$-\zeta$}}
\put(35,2){\tiny{$-1$}}
\put(31,12){\tiny{$\zeta^{-1}$}}

\put(-2,-5){\tiny{$\al_1$}}
\put(13,-5){\tiny{$\al_2$}}
\put(28,-5){\tiny{$\al_3$}}
\put(43,-5){\tiny{$\al_4$}}
\end{picture}}

\newcommand{\smallHecFourSeventeenThree}{
\setlength{\unitlength}{0.5mm}
\begin{picture}(30,20)(00,00)

\put(-5,20){\tiny{$e:=[\acchi^{(3,17)}_4]\circ\left[{{1234}\atop{3214}}\right]$}}

\multiput(15.8,3.8)(00,00){1}{\oval(6.4,6.4)[br]}
\multiput(32,3.8)(00,00){1}{\oval(26,13)[t]}
\multiput(45,01)(00,00){1}{\line(0,1){2.8}}

\multiput(00,00)(15,00){4}{\circle{2}}
\multiput(01,00)(15,00){3}{\line(1,0){13}}

\put(-6,3){\tiny{$-\zeta$}}
\put(11.5,3){\tiny{$-1$}}
\put(28,3){\tiny{$-1$}}
\put(46,2){\tiny{$-1$}}

\put(2,3.5){\tiny{$-{\frac 1 \zeta}$}}
\put(19.5,2){\tiny{$-\zeta$}}
\put(35,2){\tiny{$-\zeta$}}
\put(32,12){\tiny{$\zeta$}}

\put(-2,-5){\tiny{$\al_3$}}
\put(13,-5){\tiny{$\al_2$}}
\put(28,-5){\tiny{$\al_1$}}
\put(43,-5){\tiny{$\al_4$}}
\end{picture}}

\newcommand{\smallHecFourSeventeenFour}{
\setlength{\unitlength}{0.5mm}
\begin{picture}(30,20)(00,00)

\put(-5,16){\tiny{$a:=[\acchi^{(3,17)}_5]\circ\left[{{1234}\atop{4312}}\right]$}}

\multiput(00,00)(15,00){4}{\circle{2}}
\multiput(01,00)(15,00){3}{\line(1,0){13}}

\put(-6,3){\tiny{$-\zeta$}}
\put(12,3){\tiny{$-1$}}
\put(28.5,3){\tiny{$\zeta$}}
\put(43,3){\tiny{$-\zeta$}}

\put(2,3.5){\tiny{$-{\frac 1 \zeta}$}}
\put(19,2){\tiny{$\zeta^{-1}$}}
\put(33.5,3.5){\tiny{$-{\frac 1 \zeta}$}}

\put(-2,-5){\tiny{$\al_1$}}
\put(13,-5){\tiny{$\al_2$}}
\put(28,-5){\tiny{$\al_4$}}
\put(43,-5){\tiny{$\al_3$}}
\end{picture}}

\newcommand{\smallHecFourSeventeenEleven}{
\setlength{\unitlength}{0.5mm}
\begin{picture}(30,20)(00,00)

\put(-5,16){\tiny{$l:=[\acchi^{(3,17)}_6]\circ\left[{{1234}\atop{3214}}\right]$}}

\multiput(00,00)(15,00){4}{\circle{2}}
\multiput(01,00)(15,00){3}{\line(1,0){13}}

\put(-6,3){\tiny{$-\zeta$}}
\put(11,3){\tiny{$-\zeta$}}
\put(27,3){\tiny{$-1$}}
\put(44,3){\tiny{$-\zeta$}}

\put(2,3.5){\tiny{$-{\frac 1 \zeta}$}}
\put(18.5,3.5){\tiny{$-{\frac 1 \zeta}$}}
\put(35,3.5){\tiny{$-{\frac 1 \zeta}$}}

\put(-2,-5){\tiny{$\al_3$}}
\put(13,-5){\tiny{$\al_4$}}
\put(28,-5){\tiny{$\al_1$}}
\put(43,-5){\tiny{$\al_2$}}
\end{picture}}

\newcommand{\RankFourHecSeventeenTau}{
\setlength{\unitlength}{0.5mm}
\begin{picture}(140,30)(0,0)

\put(00,23){\footnotesize{$j$}}
\put(23,23){\footnotesize{$i$}}
\put(46,23){\footnotesize{$c$}}
\put(69,23){\footnotesize{$a$}}
\put(92,23){\footnotesize{$b$}}
\put(115,23){\footnotesize{$f$}}
\put(138,23){\footnotesize{$g$}}

\multiput(6.5,23)(23,0){6}{\vector(1,0){13}}
\multiput(18.5,23)(23,0){6}{\vector(-1,0){13}}

\multiput(47,20)(46,0){2}{\vector(0,-1){13}}
\multiput(47,8)(46,0){2}{\vector(0,1){13}}

\put(23,00){\footnotesize{$k$}}
\put(46,00){\footnotesize{$d$}}
\put(69,00){\footnotesize{$h$}}
\put(92,00){\footnotesize{$e$}}
\put(115,00){\footnotesize{$l$}}

\multiput(29.5,01)(23,0){4}{\vector(1,0){13}}
\multiput(41.5,01)(23,0){4}{\vector(-1,0){13}}

\put(10,26){\footnotesize{${\bar{\tau}}_2$}}
\put(33,26){\footnotesize{${\bar{\tau}}_3$}}
\put(56,26){\footnotesize{${\bar{\tau}}_2$}}
\put(79,26){\footnotesize{${\bar{\tau}}_4$}}
\put(102,26){\footnotesize{${\bar{\tau}}_3$}}
\put(125,26){\footnotesize{${\bar{\tau}}_4$}}

\put(40,12){\footnotesize{${\bar{\tau}}_4$}}
\put(95,12){\footnotesize{${\bar{\tau}}_2$}}

\put(33,03){\footnotesize{${\bar{\tau}}_1$}}
\put(56,03){\footnotesize{${\bar{\tau}}_2$}}
\put(79,03){\footnotesize{${\bar{\tau}}_4$}}
\put(102,03){\footnotesize{${\bar{\tau}}_1$}}

\end{picture}}

\newcommand{\allHecFourEighteen}{
\setlength{\unitlength}{1mm}
\begin{picture}(220,40)(00,-10)
\put(10,20){$\smallHecFourEighteenOne$}
\put(40,20){$\smallHecFourEighteenTwo$}
\put(70,20){$\smallHecFourEighteenThree$}
\put(100,20){$\smallHecFourEighteenFour$}
\put(10,00){$\smallHecFourEighteenSix$}
\put(40,00){$\smallHecFourEighteenFive$}
\put(75,00){$(\zeta\in\bKt,\,\zeta^2+\zeta+1=0)$}
\end{picture}}

\newcommand{\smallHecFourEighteenOne}{
\setlength{\unitlength}{0.5mm}
\begin{picture}(30,20)(00,00)

\put(-5,16){\tiny{$a:=[\acchi^{(3,18)}_1]$}}

\multiput(00,00)(15,00){4}{\circle{2}}
\multiput(01,00)(15,00){3}{\line(1,0){13}}

\put(-3,3){\tiny{$\zeta^{-1}$}}
\put(11.5,3){\tiny{$\zeta^{-1}$}}
\put(28,3){\tiny{$\zeta$}}
\put(42,3){\tiny{$-1$}}

\put(6,2){\tiny{$\zeta$}}
\put(21.5,2){\tiny{$\zeta$}}
\put(33,2){\tiny{$\zeta^{-1}$}}

\put(-2,-5){\tiny{$\al_1$}}
\put(13,-5){\tiny{$\al_2$}}
\put(28,-5){\tiny{$\al_3$}}
\put(43,-5){\tiny{$\al_4$}}
\end{picture}}

\newcommand{\smallHecFourEighteenTwo}{
\setlength{\unitlength}{0.5mm}
\begin{picture}(30,20)(00,00)

\put(-5,16){\tiny{$b:=[\acchi^{(3,18)}_2]$}}

\multiput(00,00)(15,00){4}{\circle{2}}
\multiput(01,00)(15,00){3}{\line(1,0){13}}

\put(-3,3){\tiny{$\zeta^{-1}$}}
\put(11.5,3){\tiny{$\zeta^{-1}$}}
\put(27,3){\tiny{$-1$}}
\put(42,3){\tiny{$-1$}}

\put(6,2){\tiny{$\zeta$}}
\put(21.5,2){\tiny{$\zeta$}}
\put(36,2){\tiny{$\zeta$}}

\put(-2,-5){\tiny{$\al_1$}}
\put(13,-5){\tiny{$\al_2$}}
\put(28,-5){\tiny{$\al_3$}}
\put(43,-5){\tiny{$\al_4$}}
\end{picture}}

\newcommand{\smallHecFourEighteenThree}{
\setlength{\unitlength}{0.5mm}
\begin{picture}(30,20)(00,00)

\put(-5,16){\tiny{$c:=[\acchi^{(3,18)}_3]$}}

\multiput(15.8,3.8)(00,00){1}{\oval(6.4,6.4)[br]}
\multiput(32,3.8)(00,00){1}{\oval(26,13)[t]}
\multiput(45,01)(00,00){1}{\line(0,1){2.8}}

\multiput(00,00)(15,00){4}{\circle{2}}
\multiput(01,00)(15,00){3}{\line(1,0){13}}

\put(-3,3){\tiny{$\zeta^{-1}$}}
\put(11.5,3){\tiny{$-1$}}
\put(28,3){\tiny{$-1$}}
\put(46,2){\tiny{$\zeta$}}

\put(6,2){\tiny{$\zeta$}}
\put(19.5,2){\tiny{$\zeta^{-1}$}}
\put(35,2){\tiny{$\zeta^{-1}$}}
\put(31,12){\tiny{$\zeta^{-1}$}}

\put(-2,-5){\tiny{$\al_1$}}
\put(13,-5){\tiny{$\al_2$}}
\put(28,-5){\tiny{$\al_3$}}
\put(43,-5){\tiny{$\al_4$}}
\end{picture}}

\newcommand{\smallHecFourEighteenFour}{
\setlength{\unitlength}{0.5mm}
\begin{picture}(30,20)(00,00)

\put(-5,20){\tiny{$d:=[\acchi^{(3,18)}_4]\circ\left[{{1234}\atop{3214}}\right]$}}

\multiput(15.8,3.8)(00,00){1}{\oval(6.4,6.4)[br]}
\multiput(32,3.8)(00,00){1}{\oval(26,13)[t]}
\multiput(45,01)(00,00){1}{\line(0,1){2.8}}

\multiput(00,00)(15,00){4}{\circle{2}}
\multiput(01,00)(15,00){2}{\line(1,0){13}}

\put(-5,3){\tiny{$-1$}}
\put(11.5,3){\tiny{$-1$}}
\put(28,3){\tiny{$\zeta^{-1}$}}
\put(46,2){\tiny{$-1$}}

\put(3,2){\tiny{$\zeta^{-1}$}}
\put(22.5,2){\tiny{$\zeta$}}
\put(31,12){\tiny{$\zeta$}}

\put(-2,-5){\tiny{$\al_1$}}
\put(13,-5){\tiny{$\al_2$}}
\put(28,-5){\tiny{$\al_3$}}
\put(43,-5){\tiny{$\al_4$}}
\end{picture}}

\newcommand{\smallHecFourEighteenFive}{
\setlength{\unitlength}{0.5mm}
\begin{picture}(30,20)(00,00)

\put(-5,20){\tiny{$f:=[\acchi^{(3,18)}_5]\circ\left[{{1234}\atop{3214}}\right]$}}

\multiput(15.8,3.8)(00,00){1}{\oval(6.4,6.4)[br]}
\multiput(32,3.8)(00,00){1}{\oval(26,13)[t]}
\multiput(45,01)(00,00){1}{\line(0,1){2.8}}

\multiput(00,00)(15,00){4}{\circle{2}}
\multiput(01,00)(15,00){2}{\line(1,0){13}}

\put(-5,3){\tiny{$-1$}}
\put(14,3){\tiny{$\zeta$}}
\put(28,3){\tiny{$\zeta^{-1}$}}
\put(46,2){\tiny{$-1$}}

\put(5,2){\tiny{$\zeta^{-1}$}}
\put(22.5,2){\tiny{$\zeta$}}
\put(31,12){\tiny{$\zeta^{-1}$}}

\put(-2,-5){\tiny{$\al_1$}}
\put(13,-5){\tiny{$\al_2$}}
\put(28,-5){\tiny{$\al_3$}}
\put(43,-5){\tiny{$\al_4$}}
\end{picture}}

\newcommand{\smallHecFourEighteenSix}{
\setlength{\unitlength}{0.5mm}
\begin{picture}(30,20)(00,00)

\put(-5,22){\tiny{$e:=[\acchi^{(3,18)}_6]\circ\left[{{1234}\atop{3214}}\right]$}}

\multiput(15.8,3.8)(00,00){1}{\oval(6.4,6.4)[br]}
\multiput(32,3.8)(00,00){1}{\oval(26,13)[t]}
\multiput(45,01)(00,00){1}{\line(0,1){2.8}}

\multiput(00,00)(15,00){4}{\circle{2}}
\multiput(01,00)(15,00){2}{\line(1,0){13}}

\put(-5,3){\tiny{$-1$}}
\put(10,3){\tiny{$\zeta^{-1}$}}
\put(28,3){\tiny{$\zeta^{-1}$}}
\put(46,2){\tiny{$-1$}}

\put(5,2){\tiny{$\zeta$}}
\put(22.5,2){\tiny{$\zeta$}}
\put(31,12){\tiny{$\zeta$}}

\put(-2,-5){\tiny{$\al_1$}}
\put(13,-5){\tiny{$\al_2$}}
\put(28,-5){\tiny{$\al_3$}}
\put(43,-5){\tiny{$\al_4$}}
\end{picture}}

\newcommand{\RankFourHecEighteena}{
\setlength{\unitlength}{1mm}
\begin{picture}(40,20)(5,10)
\put(8,10){$\chi_a$}
\put(10,0){\circle{2}}\put(25,0){\circle{2}}\put(40,0){\circle{2}}\put(55,0){\circle{2}}
\put(11,0){\line(1,0){13}}\put(26,0){\line(1,0){13}}\put(41,0){\line(1,0){13}}
\put(9,3){{\small{$\defzeta^{-1}$}}}\put(24,3){{\small{$\defzeta^{-1}$}}}\put(39,3){{\small{$\defzeta$}}}\put(53,3){{\small{$-1$}}}
\put(17,2.0){{\small{$\defzeta$}}}\put(32,2.0){{\small{$\defzeta$}}}\put(46,2.0){{\small{$\defzeta^{-1}$}}}
\put(8,-4){$\al_1$}\put(23,-4){$\al_2$}\put(38,-4){$\al_3$}\put(53,-4){$\al_4$}
\end{picture}
}

\newcommand{\RankFourHecEighteenb}{
\setlength{\unitlength}{1mm}
\begin{picture}(40,20)(5,10)
\put(8,10){$\chi_b$}
\put(10,0){\circle{2}}\put(25,0){\circle{2}}\put(40,0){\circle{2}}\put(55,0){\circle{2}}
\put(11,0){\line(1,0){13}}\put(26,0){\line(1,0){13}}\put(41,0){\line(1,0){13}}
\put(9,3){{\small{$\defzeta^{-1}$}}}\put(24,3){{\small{$\defzeta^{-1}$}}}\put(38,3){{\small{$-1$}}}\put(53,3){{\small{$-1$}}}
\put(17,2.0){{\small{$\defzeta$}}}\put(32,2.0){{\small{$\defzeta$}}}\put(47,2.0){{\small{$\defzeta$}}}
\put(8,-4){$\al_1$}\put(23,-4){$\al_2$}\put(38,-4){$\al_3$}\put(53,-4){$\al_4$}
\end{picture}
}

\newcommand{\RankFourHecEighteenc}{
\setlength{\unitlength}{1mm}
\begin{picture}(40,20)(5,10)
\put(8,10){$\chi_c$}
\put(10,0){\circle{2}}\put(25,0){\circle{2}}\put(40,0){\circle{2}}\put(55,0){\circle{2}}
\put(11,0){\line(1,0){13}}\put(26,0){\line(1,0){13}}\put(41,0){\line(1,0){13}}
\put(9,3){{\small{$\defzeta^{-1}$}}}\put(19,3){{\small{$-1$}}}\put(38,3){{\small{$-1$}}}\put(56,3){{\small{$\defzeta$}}}
\put(15,2.0){{\small{$\defzeta$}}}\put(31,2.0){{\small{$\defzeta^{-1}$}}}\put(46,2.0){{\small{$\defzeta^{-1}$}}}
\put(8,-4){$\al_1$}\put(23,-4){$\al_2$}\put(38,-4){$\al_3$}\put(53,-4){$\al_4$}
\put(40,1){\oval(30,14)[t]}\put(39,10){{\small{$\defzeta^{-1}$}}}
\end{picture}
}

\newcommand{\RankFourHecEighteend}{
\setlength{\unitlength}{1mm}
\begin{picture}(40,20)(5,10)
\put(8,10){$\chi_d$}
\put(10,0){\circle{2}}\put(25,0){\circle{2}}\put(40,0){\circle{2}}\put(55,0){\circle{2}}
\put(11,0){\line(1,0){13}}\put(26,0){\line(1,0){13}}
\put(7,3){{\small{$-1$}}}\put(19.5,3){{\small{$-1$}}}\put(39,3){{\small{$\defzeta^{-1}$}}}\put(56,3){{\small{$-1$}}}
\put(13,2.0){{\small{$\defzeta^{-1}$}}}\put(31.5,2.0){{\small{$\defzeta$}}}
\put(8,-4){$\al_1$}\put(23,-4){$\al_2$}\put(38,-4){$\al_3$}\put(53,-4){$\al_4$}
\put(40,1){\oval(30,14)[t]}\put(39,10){{\small{$\defzeta$}}}
\end{picture}
}

\newcommand{\RankFourHecEighteene}{
\setlength{\unitlength}{1mm}
\begin{picture}(40,20)(5,10)
\put(8,10){$\chi_e$}
\put(10,0){\circle{2}}\put(25,0){\circle{2}}\put(40,0){\circle{2}}\put(55,0){\circle{2}}
\put(11,0){\line(1,0){13}}\put(26,0){\line(1,0){13}}
\put(7,3){{\small{$-1$}}}\put(19.5,3){{\small{$\defzeta^{-1}$}}}\put(39,3){{\small{$\defzeta^{-1}$}}}\put(56,3){{\small{$-1$}}}
\put(15,2.0){{\small{$\defzeta$}}}\put(31.5,2.0){{\small{$\defzeta$}}}
\put(8,-4){$\al_1$}\put(23,-4){$\al_2$}\put(38,-4){$\al_3$}\put(53,-4){$\al_4$}
\put(40,1){\oval(30,14)[t]}\put(39,10){{\small{$\defzeta$}}}
\end{picture}
}

\newcommand{\RankFourHecEighteenf}{
\setlength{\unitlength}{1mm}
\begin{picture}(40,20)(5,10)
\put(8,10){$\chi_f$}
\put(10,0){\circle{2}}\put(25,0){\circle{2}}\put(40,0){\circle{2}}\put(55,0){\circle{2}}
\put(11,0){\line(1,0){13}}\put(26,0){\line(1,0){13}}
\put(7,3){{\small{$-1$}}}\put(21.5,3){{\small{$\defzeta$}}}\put(39,3){{\small{$\defzeta^{-1}$}}}\put(56,3){{\small{$-1$}}}
\put(14,2.0){{\small{$\defzeta^{-1}$}}}\put(31.5,2.0){{\small{$\defzeta$}}}
\put(8,-4){$\al_1$}\put(23,-4){$\al_2$}\put(38,-4){$\al_3$}\put(53,-4){$\al_4$}
\put(40,1){\oval(30,14)[t]}\put(39,10){{\small{$\defzeta^{-1}$}}}
\end{picture}
}

\newcommand{\RankFourHecEighteenabcde}
{\setlength{\unitlength}{1mm}
\begin{picture}(120,75)(0,-15)
\put(0,50){\RankFourHecEighteena}
\put(60,50){\RankFourHecEighteenb}

\put(0,25){\RankFourHecEighteenc}
\put(60,25){\RankFourHecEighteend}

\put(0,0){\RankFourHecEighteene}
\put(60,0){\RankFourHecEighteenf}
\end{picture}
}

\newcommand{\RankFourHecEighteenTau}
{\setlength{\unitlength}{1mm}
\begin{picture}(115,20)(0,0)

\multiput(6,16)(15,00){4}{\vector(1,0){7}}
\multiput(11,16)(15,00){4}{\vector(-1,0){7}}

\put(0,15){$a$}\put(7,17.5){${\bar{\tau}}_4$}
\put(15,15){$b$}\put(22,17.5){${\bar{\tau}}_3$}
\put(30,15){$c$}\put(37,17.5){${\bar{\tau}}_2$}
\put(45,15){$d$}\put(52,17.5){${\bar{\tau}}_1$}
\put(60,15){$e$}

\multiput(46,6.5)(15,00){2}{\vector(0,1){7}}
\multiput(46,11.5)(15,00){2}{\vector(0,-1){7}}
\put(42,8){${\bar{\tau}}_4$}
\put(62,8){${\bar{\tau}}_4$}

\multiput(51,1)(15,00){4}{\vector(1,0){7}}
\multiput(56,1)(15,00){4}{\vector(-1,0){7}}

\put(45,0){$f$}\put(52,2.5){${\bar{\tau}}_1$}
\put(60,0){$d^\prime$}\put(67,2.5){${\bar{\tau}}_3$}
\put(75,0){$c^\prime$}\put(82,2.5){${\bar{\tau}}_2$}
\put(90,0){$b^\prime$}\put(97,2.5){${\bar{\tau}}_1$}
\put(105,0){$a^\prime$}

\end{picture}
}

\newcommand{\CayleyAdash}{
\setlength{\unitlength}{1mm}
\begin{picture}(60,80)(25,-3)

\put(10,19){\curveSouthWest}
\put(10,9){\line(1,0){90}}
\put(100,19){\curveSouthEast}
\put(54.5,10){{\tiny{${\bar{3}}$}}}

\put(20,51){\curveNorthWest}
\put(20,61){\line(1,0){70}}
\put(90,51){\curveNorthEast}
\put(54.5,58){{\tiny{${\bar{2}}$}}}

\put(10,61){\curveNorthWestTHICK}
\put(10,71){\line(1,0){90}}
\put(10,70.95){\line(1,0){90}}\put(10,70.90){\line(1,0){90}}\put(10,70.85){\line(1,0){90}}
\put(10,71.05){\line(1,0){90}}\put(10,71.10){\line(1,0){90}}\put(10,71.15){\line(1,0){90}}
\put(100,61){\curveNorthEastTHICK}
\put(0,31){\line(0,1){30}}
\put(-0.05,31){\line(0,1){30}}\put(-0.10,31){\line(0,1){30}}\put(-0.15,31){\line(0,1){30}}
\put(0.05,31){\line(0,1){30}}\put(0.10,31){\line(0,1){30}}\put(0.15,31){\line(0,1){30}}

\put(110,31){\line(0,1){30}}
\put(109.95,31){\line(0,1){30}}\put(109.90,31){\line(0,1){30}}\put(109.85,31){\line(0,1){30}}
\put(110.05,31){\line(0,1){30}}\put(110.10,31){\line(0,1){30}}\put(110.15,31){\line(0,1){30}}

\put(54.5,68){{\tiny{${\bar{3}}$}}}

\put(0,20){\whitecircleempty}\put(40,20){\whitecircleempty}\put(80,20){\whitecircleempty}
\put(0,20){\verticallineRightoneTHICKbar}

\put(40,20){\verticallineRightoneTHICKbar}\put(80,20){\verticallineRightonebar}

\put(1,20){\line(1,0){10}}\put(41,20){\line(1,0){10}}\put(81,20){\line(1,0){10}}
\multiput(1,20.05)(80,0){2}{\line(1,0){10}}\multiput(1,20.10)(80,0){2}{\line(1,0){10}}\multiput(1,20.15)(80,0){2}{\line(1,0){10}}
\multiput(1,19.95)(80,0){2}{\line(1,0){10}}\multiput(1,19.90)(80,0){2}{\line(1,0){10}}\multiput(1,19.85)(80,0){2}{\line(1,0){10}}
\multiput(10,20)(80,0){2}{\holizontallineAbovetwoTHICK}\put(50,20){\holizontallineAbovetwo}
\multiput(19,20)(40,0){3}{\line(1,0){10}}
\multiput(19,20.05)(80,0){2}{\line(1,0){10}}\multiput(19,20.10)(80,0){2}{\line(1,0){10}}\multiput(19,20.15)(80,0){2}{\line(1,0){10}}
\multiput(19,19.95)(80,0){2}{\line(1,0){10}}\multiput(19,19.90)(80,0){2}{\line(1,0){10}}\multiput(19,19.85)(80,0){2}{\line(1,0){10}}
\multiput(30,20)(40,0){3}{\whitecircleempty}
\put(30,20){\verticallineRightonebar}\put(70,20){\verticallineRightoneTHICKbar}
\put(110,20){\verticallineRightoneTHICKbar}
\multiput(30,20)(40,0){2}{\holizontallineAbovethreeTHICKbar}

\multiput(0,30)(40,0){3}{\whitecircleempty}
\multiput(0,30)(40,0){2}{\EastUplineabovetwo}\put(80,30){\EastUplineabovetwoTHICK}
\multiput(30,30)(40,0){3}{\whitecircleempty}
\multiput(30,30)(40,0){2}{\holizontallineAbovethreeTHICKbar}

\multiput(10,40)(40,0){3}{\whitecircleempty}
\multiput(10,40)(40,0){2}{\verticallineRightthreeTHICK}\put(90,40){\verticallineRightthree}
\multiput(10,40)(40,0){3}{\holizontallineAboveoneTHICKbar}
\multiput(20,40)(40,0){3}{\whitecircleempty}
\put(20,40){\verticallineRightthree}\multiput(60,40)(40,0){2}{\verticallineRightthreeTHICK}

\put(20,40){\EastDownlineabovetwoTHICK}
\multiput(60,40)(40,0){2}{\EastDownlineabovetwo}

\multiput(10,50)(40,0){3}{\whitecircleempty}
\multiput(10,50)(80,0){2}{\holizontallineAboveoneTHICKbar}\put(50,50){\holizontallineAboveonebar}
\multiput(20,50)(40,0){3}{\whitecircleempty}

\multiput(21,50)(40,0){2}{\line(1,0){10}}
\multiput(21,49.95)(40,0){2}{\line(1,0){10}}\multiput(21,49.90)(40,0){2}{\line(1,0){10}}\multiput(21,49.85)(40,0){2}{\line(1,0){10}}
\multiput(21,50.05)(40,0){2}{\line(1,0){10}}\multiput(21,50.10)(40,0){2}{\line(1,0){10}}\multiput(21,50.15)(40,0){2}{\line(1,0){10}}
\multiput(30,50)(40,0){2}{\holizontallineAbovetwoTHICKbar}
\multiput(39,50)(40,0){2}{\line(1,0){10}}
\multiput(39,49.95)(40,0){2}{\line(1,0){10}}\multiput(39,49.90)(40,0){2}{\line(1,0){10}}\multiput(39,49.85)(40,0){2}{\line(1,0){10}}
\multiput(39,50.05)(40,0){2}{\line(1,0){10}}\multiput(39,50.10)(40,0){2}{\line(1,0){10}}\multiput(39,50.15)(40,0){2}{\line(1,0){10}}

\put(-4, 19){$e$}
\put(-4, 29){$d$}
\put(6, 39){$c$}\put(46, 39){$c$}\put(86, 39){$c$}
\put(6, 49){$b$}

\put(111.5, 19){$e$}
\put(111.5, 29){$d$}
\put(101.5, 39){$c$}\put(61.5, 39){$c$}\put(21.5, 39){$c$}
\put(101.5, 49){$b$}

\put(29, 16){$e$}\put(39, 16){$e$}\put(69, 16){$e$}\put(79, 16){$e$}
\put(29, 32.5){$d$}\put(39, 32.5){$d$}\put(69, 32){$d$}\put(79, 32){$d$}
\put(19, 52){$b$}\put(49, 52){$b$}\put(59, 52){$b$}\put(89, 52){$b$}

\end{picture}
}

\newcommand{\CayleyAdashsmall}{
\setlength{\unitlength}{0.5mm}
\begin{picture}(60,80)(25,-3)

\put(10,19){\curveSouthWest}
\put(10,9){\line(1,0){90}}
\put(100,19){\curveSouthEast}
\put(54.5,10){{\tiny{${\bar{3}}$}}}

\put(20,51){\curveNorthWest}
\put(20,61){\line(1,0){70}}
\put(90,51){\curveNorthEast}
\put(54.5,58){{\tiny{${\bar{2}}$}}}

\put(10,61){\curveNorthWestTHICK}
\put(10,71){\line(1,0){90}}
\put(10,70.95){\line(1,0){90}}\put(10,70.90){\line(1,0){90}}\put(10,70.85){\line(1,0){90}}
\put(10,71.05){\line(1,0){90}}\put(10,71.10){\line(1,0){90}}\put(10,71.15){\line(1,0){90}}
\put(100,61){\curveNorthEastTHICK}
\put(0,31){\line(0,1){30}}
\put(-0.05,31){\line(0,1){30}}\put(-0.10,31){\line(0,1){30}}\put(-0.15,31){\line(0,1){30}}
\put(0.05,31){\line(0,1){30}}\put(0.10,31){\line(0,1){30}}\put(0.15,31){\line(0,1){30}}

\put(110,31){\line(0,1){30}}
\put(109.95,31){\line(0,1){30}}\put(109.90,31){\line(0,1){30}}\put(109.85,31){\line(0,1){30}}
\put(110.05,31){\line(0,1){30}}\put(110.10,31){\line(0,1){30}}\put(110.15,31){\line(0,1){30}}

\put(54.5,68){{\tiny{${\bar{3}}$}}}

\put(0,20){\whitecircleempty}\put(40,20){\whitecircleempty}\put(80,20){\whitecircleempty}
\put(0,20){\verticallineRightoneTHICKbar}

\put(40,20){\verticallineRightoneTHICKbar}\put(80,20){\verticallineRightonebar}

\put(1,20){\line(1,0){10}}\put(41,20){\line(1,0){10}}\put(81,20){\line(1,0){10}}
\multiput(1,20.05)(80,0){2}{\line(1,0){10}}\multiput(1,20.10)(80,0){2}{\line(1,0){10}}\multiput(1,20.15)(80,0){2}{\line(1,0){10}}
\multiput(1,19.95)(80,0){2}{\line(1,0){10}}\multiput(1,19.90)(80,0){2}{\line(1,0){10}}\multiput(1,19.85)(80,0){2}{\line(1,0){10}}
\multiput(10,20)(80,0){2}{\holizontallineAbovetwoTHICK}\put(50,20){\holizontallineAbovetwo}
\multiput(19,20)(40,0){3}{\line(1,0){10}}
\multiput(19,20.05)(80,0){2}{\line(1,0){10}}\multiput(19,20.10)(80,0){2}{\line(1,0){10}}\multiput(19,20.15)(80,0){2}{\line(1,0){10}}
\multiput(19,19.95)(80,0){2}{\line(1,0){10}}\multiput(19,19.90)(80,0){2}{\line(1,0){10}}\multiput(19,19.85)(80,0){2}{\line(1,0){10}}
\multiput(30,20)(40,0){3}{\whitecircleempty}
\put(30,20){\verticallineRightonebar}\put(70,20){\verticallineRightoneTHICKbar}
\put(110,20){\verticallineRightoneTHICKbar}
\multiput(30,20)(40,0){2}{\holizontallineAbovethreeTHICKbar}

\multiput(0,30)(40,0){3}{\whitecircleempty}
\multiput(0,30)(40,0){2}{\EastUplineabovetwo}\put(80,30){\EastUplineabovetwoTHICK}
\multiput(30,30)(40,0){3}{\whitecircleempty}
\multiput(30,30)(40,0){2}{\holizontallineAbovethreeTHICKbar}

\multiput(10,40)(40,0){3}{\whitecircleempty}
\multiput(10,40)(40,0){2}{\verticallineRightthreeTHICK}\put(90,40){\verticallineRightthree}
\multiput(10,40)(40,0){3}{\holizontallineAboveoneTHICKbar}
\multiput(20,40)(40,0){3}{\whitecircleempty}
\put(20,40){\verticallineRightthree}\multiput(60,40)(40,0){2}{\verticallineRightthreeTHICK}

\put(20,40){\EastDownlineabovetwoTHICK}
\multiput(60,40)(40,0){2}{\EastDownlineabovetwo}

\multiput(10,50)(40,0){3}{\whitecircleempty}
\multiput(10,50)(80,0){2}{\holizontallineAboveoneTHICKbar}\put(50,50){\holizontallineAboveonebar}
\multiput(20,50)(40,0){3}{\whitecircleempty}

\multiput(21,50)(40,0){2}{\line(1,0){10}}
\multiput(21,49.95)(40,0){2}{\line(1,0){10}}\multiput(21,49.90)(40,0){2}{\line(1,0){10}}\multiput(21,49.85)(40,0){2}{\line(1,0){10}}
\multiput(21,50.05)(40,0){2}{\line(1,0){10}}\multiput(21,50.10)(40,0){2}{\line(1,0){10}}\multiput(21,50.15)(40,0){2}{\line(1,0){10}}
\multiput(30,50)(40,0){2}{\holizontallineAbovetwoTHICKbar}
\multiput(39,50)(40,0){2}{\line(1,0){10}}
\multiput(39,49.95)(40,0){2}{\line(1,0){10}}\multiput(39,49.90)(40,0){2}{\line(1,0){10}}\multiput(39,49.85)(40,0){2}{\line(1,0){10}}
\multiput(39,50.05)(40,0){2}{\line(1,0){10}}\multiput(39,50.10)(40,0){2}{\line(1,0){10}}\multiput(39,50.15)(40,0){2}{\line(1,0){10}}

\put(-4, 19){$e$}
\put(-4, 29){$d$}
\put(6, 39){$c$}\put(46, 39){$c$}\put(86, 39){$c$}
\put(6, 49){$b$}

\put(111.5, 19){$e$}
\put(111.5, 29){$d$}
\put(101.5, 39){$c$}\put(61.5, 39){$c$}\put(21.5, 39){$c$}
\put(101.5, 49){$b$}

\put(29, 16){$e$}\put(39, 16){$e$}\put(69, 16){$e$}\put(79, 16){$e$}
\put(29, 32.5){$d$}\put(39, 32.5){$d$}\put(69, 32){$d$}\put(79, 32){$d$}
\put(19, 52){$b$}\put(49, 52){$b$}\put(59, 52){$b$}\put(89, 52){$b$}

\end{picture}
}

\newcommand{\allHecFourTewnty}{
\setlength{\unitlength}{1mm}
\begin{picture}(220,60)(-5,-10)
\put(10,40){$\smallHecFourTwentyZero$}
\put(40,40){$\smallHecFourTwentyOne$}
\put(70,40){$\smallHecFourTwentyTwo$}
\put(100,40){$\smallHecFourTwentyThree$}
\put(10,20){$\smallHecFourTwentyFour$}
\put(40,20){$\smallHecFourTwentyFive$}
\put(70,20){$\smallHecFourTwentySix$}
\put(100,20){$\smallHecFourTwentySeven$}
\put(10,00){$\smallHecFourTwentyEight$}
\put(40,00){$\smallHecFourTwentyNine$}
\put(75,00){$(\zeta\in\bKt,\,\zeta^2+\zeta+1=0)$}
\end{picture}}

\newcommand{\smallHecFourTwentyZero}{
\setlength{\unitlength}{0.5mm}
\begin{picture}(30,20)(00,00)

\put(-5,16){\tiny{$a:=[\acchi^{(3,20)}_7]$}}

\multiput(15.8,3.8)(00,00){1}{\oval(6.4,6.4)[br]}
\multiput(32,3.8)(00,00){1}{\oval(26,13)[t]}
\multiput(45,01)(00,00){1}{\line(0,1){2.8}}

\multiput(00,00)(15,00){4}{\circle{2}}
\multiput(01,00)(15,00){3}{\line(1,0){13}}

\put(-3,3){\tiny{$-1$}}
\put(14,3){\tiny{$\zeta$}}
\put(29,3){\tiny{$-1$}}
\put(46,3){\tiny{$\zeta$}}

\put(5,2){\tiny{$\zeta^{-1}$}}
\put(20,2){\tiny{$\zeta^{-1}$}}
\put(36,2){\tiny{$\zeta^{-1}$}}
\put(31,12){\tiny{$\zeta^{-1}$}}

\put(-2,-5){\tiny{$\al_1$}}
\put(13,-5){\tiny{$\al_2$}}
\put(28,-5){\tiny{$\al_3$}}
\put(43,-5){\tiny{$\al_4$}}
\end{picture}}

\newcommand{\smallHecFourTwentyOne}{
\setlength{\unitlength}{0.5mm}
\begin{picture}(30,20)(00,00)

\put(-5,16){\tiny{$b:=[\acchi^{(3,20)}_3]$}}

\multiput(00,00)(15,00){4}{\circle{2}}
\multiput(01,00)(15,00){3}{\line(1,0){13}}

\put(-3,3){\tiny{$-1$}}
\put(14,3){\tiny{$-1$}}
\put(27,3){\tiny{$-1$}}
\put(42,3){\tiny{$-1$}}

\put(5,2){\tiny{$\zeta^{-1}$}}
\put(21,2){\tiny{$\zeta$}}
\put(36,2){\tiny{$\zeta$}}

\put(-2,-5){\tiny{$\al_1$}}
\put(13,-5){\tiny{$\al_2$}}
\put(28,-5){\tiny{$\al_3$}}
\put(43,-5){\tiny{$\al_4$}}
\end{picture}}

\newcommand{\smallHecFourTwentyTwo}{
\setlength{\unitlength}{0.5mm}
\begin{picture}(30,20)(00,00)

\put(-5,16){\tiny{$c:=[\acchi^{(3,20)}_1]$}}

\multiput(00,00)(15,00){4}{\circle{2}}
\multiput(01,00)(15,00){3}{\line(1,0){13}}

\put(-3,3){\tiny{$\zeta^{-1}$}}
\put(12,3){\tiny{$-1$}}
\put(29,3){\tiny{$\zeta$}}
\put(42,3){\tiny{$-1$}}

\put(6,2){\tiny{$\zeta$}}
\put(19.5,2){\tiny{$\zeta^{-1}$}}
\put(36,2){\tiny{$\zeta$}}

\put(-2,-5){\tiny{$\al_1$}}
\put(13,-5){\tiny{$\al_2$}}
\put(28,-5){\tiny{$\al_3$}}
\put(43,-5){\tiny{$\al_4$}}
\end{picture}}

\newcommand{\smallHecFourTwentyThree}{
\setlength{\unitlength}{0.5mm}
\begin{picture}(30,20)(00,00)

\put(-5,16){\tiny{$d:=[\acchi^{(3,20)}_5]$}}

\multiput(00,00)(15,00){4}{\circle{2}}
\multiput(01,00)(15,00){3}{\line(1,0){13}}

\put(-3,3){\tiny{$\zeta^{-1}$}}
\put(12,3){\tiny{$-1$}}
\put(29,3){\tiny{$\zeta$}}
\put(42,3){\tiny{$-1$}}

\put(6,2){\tiny{$\zeta$}}
\put(19.5,2){\tiny{$\zeta^{-1}$}}
\put(36,2){\tiny{$\zeta$}}

\put(-2,-5){\tiny{$\al_1$}}
\put(13,-5){\tiny{$\al_2$}}
\put(28,-5){\tiny{$\al_3$}}
\put(43,-5){\tiny{$\al_4$}}
\end{picture}}

\newcommand{\smallHecFourTwentyFour}{
\setlength{\unitlength}{0.5mm}
\begin{picture}(30,20)(00,00)

\put(-5,16){\tiny{$e:=[\acchi^{(3,20)}_2]$}}

\multiput(00,00)(15,00){4}{\circle{2}}
\multiput(01,00)(15,00){3}{\line(1,0){13}}

\put(-3,3){\tiny{$\zeta^{-1}$}}
\put(12,3){\tiny{$-1$}}
\put(29,3){\tiny{$\al_3$}}
\put(45,3){\tiny{$-1$}}

\put(6,2){\tiny{$\zeta$}}
\put(19.5,2){\tiny{$\zeta^{-1}$}}
\put(36,2){\tiny{$\zeta^{-1}$}}

\put(-2,-5){\tiny{$\al_1$}}
\put(13,-5){\tiny{$\al_2$}}
\put(26,-6){\tiny{$-\zeta^{-1}$}}
\put(43,-5){\tiny{$\al_4$}}
\end{picture}}

\newcommand{\smallHecFourTwentyFive}{
\setlength{\unitlength}{0.5mm}
\begin{picture}(30,20)(00,00)

\put(-5,16){\tiny{${\bar{a}}:=[\acchi^{(3,20)}_{10}]$}}

\multiput(15.8,3.8)(00,00){1}{\oval(6.4,6.4)[br]}
\multiput(32,3.8)(00,00){1}{\oval(26,13)[t]}
\multiput(45,01)(00,00){1}{\line(0,1){2.8}}

\multiput(00,00)(15,00){4}{\circle{2}}
\multiput(01,00)(15,00){2}{\line(1,0){13}}

\put(-1,3){\tiny{$\zeta$}}
\put(14,3){\tiny{$\zeta$}}
\put(29,3){\tiny{$\zeta^{-1}$}}
\put(46,3){\tiny{$-1$}}

\put(5,2){\tiny{$\zeta^{-1}$}}
\put(21,2){\tiny{$\zeta$}}
\put(31,12){\tiny{$\zeta^{-1}$}}

\put(-2,-5){\tiny{$\al_1$}}
\put(13,-5){\tiny{$\al_2$}}
\put(28,-5){\tiny{$\al_3$}}
\put(43,-5){\tiny{$\al_4$}}
\end{picture}}

\newcommand{\smallHecFourTwentySix}{
\setlength{\unitlength}{0.5mm}
\begin{picture}(30,20)(00,00)

\put(-5,16){\tiny{${\bar{b}}:=[\acchi^{(3,20)}_9]$}}

\multiput(15.8,3.8)(00,00){1}{\oval(6.4,6.4)[br]}
\multiput(32,3.8)(00,00){1}{\oval(26,13)[t]}
\multiput(45,01)(00,00){1}{\line(0,1){2.8}}

\multiput(00,00)(15,00){4}{\circle{2}}
\multiput(01,00)(15,00){2}{\line(1,0){13}}

\put(-1,3){\tiny{$\zeta$}}
\put(12,3){\tiny{$-1$}}
\put(29,3){\tiny{$\zeta^{-1}$}}
\put(46,3){\tiny{$-1$}}

\put(3.5,2){\tiny{$\zeta^{-1}$}}
\put(21,2){\tiny{$\zeta$}}
\put(31,12){\tiny{$\zeta$}}

\put(-2,-5){\tiny{$\al_1$}}
\put(13,-5){\tiny{$\al_2$}}
\put(28,-5){\tiny{$\al_3$}}
\put(43,-5){\tiny{$\al_4$}}
\end{picture}}

\newcommand{\smallHecFourTwentySeven}{
\setlength{\unitlength}{0.5mm}
\begin{picture}(30,20)(00,00)

\put(-5,16){\tiny{${\bar{c}}:=[\acchi^{(3,20)}_8]$}}

\multiput(15.8,3.8)(00,00){1}{\oval(6.4,6.4)[br]}
\multiput(32,3.8)(00,00){1}{\oval(26,13)[t]}
\multiput(45,01)(00,00){1}{\line(0,1){2.8}}

\multiput(00,00)(15,00){4}{\circle{2}}
\multiput(01,00)(15,00){3}{\line(1,0){13}}

\put(-3,3){\tiny{$-1$}}
\put(12,3){\tiny{$-1$}}
\put(29,3){\tiny{$-1$}}
\put(46,3){\tiny{$\zeta$}}

\put(6,2){\tiny{$\zeta$}}
\put(20,2){\tiny{$\zeta^{-1}$}}
\put(36,2){\tiny{$\zeta^{-1}$}}
\put(31,12){\tiny{$\zeta^{-1}$}}

\put(-2,-5){\tiny{$\al_1$}}
\put(13,-5){\tiny{$\al_2$}}
\put(28,-5){\tiny{$\al_3$}}
\put(43,-5){\tiny{$\al_4$}}
\end{picture}}

\newcommand{\smallHecFourTwentyEight}{
\setlength{\unitlength}{0.5mm}
\begin{picture}(30,20)(00,00)

\put(-5,16){\tiny{${\bar{d}}:=[\acchi^{(3,20)}_4]$}}

\multiput(00,00)(15,00){4}{\circle{2}}
\multiput(01,00)(15,00){3}{\line(1,0){13}}

\put(-3,3){\tiny{$-1$}}
\put(14,3){\tiny{$-1$}}
\put(27,3){\tiny{$-1$}}
\put(42,3){\tiny{$-1$}}

\put(5,2){\tiny{$\zeta^{-1}$}}
\put(21,2){\tiny{$\zeta$}}
\put(36,2){\tiny{$\zeta$}}

\put(-2,-5){\tiny{$\al_1$}}
\put(13,-5){\tiny{$\al_2$}}
\put(28,-5){\tiny{$\al_3$}}
\put(43,-5){\tiny{$\al_4$}}
\end{picture}}

\newcommand{\smallHecFourTwentyNine}{
\setlength{\unitlength}{0.5mm}
\begin{picture}(30,20)(00,00)

\put(-5,16){\tiny{${\bar{e}}:=[\acchi^{(3,20)}_6]$}}

\multiput(00,00)(15,00){4}{\circle{2}}
\multiput(01,00)(15,00){3}{\line(1,0){13}}

\put(-4,3){\tiny{$-1$}}
\put(13,3){\tiny{$\zeta^{-1}$}}
\put(29,3){\tiny{$\zeta$}}
\put(43,3){\tiny{$-1$}}

\put(3,2){\tiny{$\zeta^{-1}$}}
\put(22,2){\tiny{$\zeta$}}
\put(34,2){\tiny{$\zeta^{-1}$}}

\put(-2,-5){\tiny{$\al_1$}}
\put(13,-5){\tiny{$\al_2$}}
\put(28,-5){\tiny{$\al_3$}}
\put(43,-5){\tiny{$\al_4$}}
\end{picture}}

\newcommand{\RankFourHecTwentyTau}{
\setlength{\unitlength}{0.5mm}
\begin{picture}(140,70)(0,0)

\put(46,46){\footnotesize{$a$}}
\put(69,46){\footnotesize{${\bar{c}}$}}
\put(92,46){\footnotesize{${\bar{b}}$}}
\put(115,46){\footnotesize{${\bar{a}}$}}
\multiput(52.5,47)(23,0){3}{\vector(1,0){13}}
\multiput(64.5,47)(23,0){3}{\vector(-1,0){13}}

\multiput(47,43)(23,0){2}{\vector(0,-1){13}}
\multiput(47,31)(23,0){2}{\vector(0,1){13}}

\put(23,23){\footnotesize{$c$}}
\put(46,23){\footnotesize{$b$}}
\put(69,23){\footnotesize{${\bar{d}}$}}

\multiput(29.5,24)(23,0){2}{\vector(1,0){13}}
\multiput(41.5,24)(23,0){2}{\vector(-1,0){13}}

\multiput(24,20)(23,0){3}{\vector(0,-1){13}}
\multiput(24,8)(23,0){3}{\vector(0,1){13}}

\put(23,00){\footnotesize{$e$}}
\put(46,00){\footnotesize{$d$}}
\put(69,00){\footnotesize{${\bar{e}}$}}
\multiput(29.5,1)(23,0){2}{\vector(1,0){13}}
\multiput(41.5,1)(23,0){2}{\vector(-1,0){13}}

\put(56,49){\footnotesize{${\bar{\tau}}_1$}}
\put(79,49){\footnotesize{${\bar{\tau}}_2$}}
\put(102,49){\footnotesize{${\bar{\tau}}_4$}}

\put(40,35){\footnotesize{${\bar{\tau}}_3$}}
\put(72,35){\footnotesize{${\bar{\tau}}_3$}}

\put(33,26){\footnotesize{${\bar{\tau}}_2$}}
\put(56,26){\footnotesize{${\bar{\tau}}_1$}}

\put(17,12){\footnotesize{${\bar{\tau}}_4$}}
\multiput(48,12)(23,0){2}{\footnotesize{${\bar{\tau}}_4$}}

\put(33,03){\footnotesize{${\bar{\tau}}_2$}}
\put(56,03){\footnotesize{${\bar{\tau}}_1$}}

\end{picture}}

\newcommand{\allHecFourTwentyone}{
\setlength{\unitlength}{1mm}
\begin{picture}(220,40)(00,-10)
\put(10,20){$\smallHecFourTwentyoneOne$}
\put(40,20){$\smallHecFourTwentyoneTwo$}
\put(70,20){$\smallHecFourTwentyoneThree$}
\put(100,20){$\smallHecFourTwentyoneFour$}
\put(10,00){$\smallHecFourTwentyoneFive$}
\put(40,00){$\smallHecFourTwentyoneSix$}
\put(70,00){$\smallHecFourTwentyoneSeven$}
\put(10,-10){$(\zeta\in\bKt,\,\zeta^2+\zeta+1=0)$}
\end{picture}}

\newcommand{\smallHecFourTwentyoneOne}{
\setlength{\unitlength}{0.5mm}
\begin{picture}(30,20)(00,00)

\put(-5,16){\tiny{${\bar{a}}:=[\acchi^{(3,21)}_1]$}}

\multiput(00,00)(15,00){4}{\circle{2}}
\multiput(01,00)(15,00){3}{\line(1,0){13}}

\put(-4,3){\tiny{$-1$}}
\put(13.5,3){\tiny{$\zeta$}}
\put(28.5,3){\tiny{$\zeta$}}
\put(42,3){\tiny{$-1$}}

\put(3,2){\tiny{$\zeta^{-1}$}}
\put(18,2){\tiny{$\zeta^{-1}$}}
\put(36,2){\tiny{$\zeta$}}

\put(-2,-5){\tiny{$\al_1$}}
\put(13,-5){\tiny{$\al_2$}}
\put(28,-5){\tiny{$\al_3$}}
\put(43,-5){\tiny{$\al_4$}}
\end{picture}}

\newcommand{\smallHecFourTwentyoneTwo}{
\setlength{\unitlength}{0.5mm}
\begin{picture}(30,20)(00,00)

\put(-5,16){\tiny{${\bar{b}}:=[\acchi^{(3,21)}_2]$}}

\multiput(00,00)(15,00){4}{\circle{2}}
\multiput(01,00)(15,00){3}{\line(1,0){13}}

\put(-4,3){\tiny{$-1$}}
\put(13.5,3){\tiny{$\zeta$}}
\put(25,-6){\tiny{$-\zeta^{-1}$}}
\put(43,3){\tiny{$-1$}}

\put(3,2){\tiny{$\zeta^{-1}$}}
\put(18,2){\tiny{$\zeta^{-1}$}}
\put(33.5,2){\tiny{$\zeta^{-1}$}}

\put(-2,-5){\tiny{$\al_1$}}
\put(13,-5){\tiny{$\al_2$}}
\put(27.5,3){\tiny{$\al_3$}}
\put(43,-5){\tiny{$\al_4$}}
\end{picture}}

\newcommand{\smallHecFourTwentyoneThree}{
\setlength{\unitlength}{0.5mm}
\begin{picture}(30,20)(00,00)

\put(-5,16){\tiny{$c:=[\acchi^{(3,21)}_3]$}}

\multiput(00,00)(15,00){4}{\circle{2}}
\multiput(01,00)(15,00){3}{\line(1,0){13}}

\put(-4,3){\tiny{$-1$}}
\put(11,3){\tiny{$-1$}}
\put(28.5,3){\tiny{$\zeta$}}
\put(42,3){\tiny{$-1$}}

\put(6,2){\tiny{$\zeta$}}
\put(18,2){\tiny{$\zeta^{-1}$}}
\put(36,2){\tiny{$\zeta$}}

\put(-2,-5){\tiny{$\al_1$}}
\put(13,-5){\tiny{$\al_2$}}
\put(28,-5){\tiny{$\al_3$}}
\put(43,-5){\tiny{$\al_4$}}
\end{picture}}

\newcommand{\smallHecFourTwentyoneFour}{
\setlength{\unitlength}{0.5mm}
\begin{picture}(30,20)(00,00)

\put(-5,16){\tiny{$e:=[\acchi^{(3,21)}_4]$}}

\multiput(00,00)(15,00){4}{\circle{2}}
\multiput(01,00)(15,00){3}{\line(1,0){13}}

\put(-4,3){\tiny{$-1$}}
\put(11,3){\tiny{$-1$}}
\put(25,-6){\tiny{$-\zeta^{-1}$}}
\put(43.5,3){\tiny{$-1$}}

\put(6,2){\tiny{$\zeta$}}
\put(18,2){\tiny{$\zeta^{-1}$}}
\put(34,2){\tiny{$\zeta^{-1}$}}

\put(-2,-5){\tiny{$\al_1$}}
\put(13,-5){\tiny{$\al_2$}}
\put(27.5,3){\tiny{$\al_3$}}
\put(43,-5){\tiny{$\al_4$}}
\end{picture}}

\newcommand{\smallHecFourTwentyoneFive}{
\setlength{\unitlength}{0.5mm}
\begin{picture}(30,20)(00,00)

\put(-5,16){\tiny{$b:=[\acchi^{(3,21)}_6]$}}

\multiput(00,00)(15,00){4}{\circle{2}}
\multiput(01,00)(15,00){3}{\line(1,0){13}}

\put(-1,3){\tiny{$\zeta$}}
\put(13,3){\tiny{$-1$}}
\put(27,3){\tiny{$-1$}}
\put(42,3){\tiny{$-1$}}

\put(3,2){\tiny{$\zeta^{-1}$}}
\put(21,2){\tiny{$\zeta$}}
\put(36,2){\tiny{$\zeta$}}

\put(-2,-5){\tiny{$\al_1$}}
\put(13,-5){\tiny{$\al_2$}}
\put(28,-5){\tiny{$\al_3$}}
\put(43,-5){\tiny{$\al_4$}}
\end{picture}}

\newcommand{\smallHecFourTwentyoneSix}{
\setlength{\unitlength}{0.5mm}
\begin{picture}(30,20)(00,00)

\put(-5,16){\tiny{$d:=[\acchi^{(3,21)}_5]$}}

\multiput(00,00)(15,00){4}{\circle{2}}
\multiput(01,00)(15,00){3}{\line(1,0){13}}

\put(-1,3){\tiny{$\zeta$}}
\put(13,3){\tiny{$-1$}}
\put(29,3){\tiny{$\zeta$}}
\put(43,3){\tiny{$-1$}}

\put(3,2){\tiny{$\zeta^{-1}$}}
\put(22,2){\tiny{$\zeta$}}
\put(33,2){\tiny{$\zeta^{-1}$}}

\put(-2,-5){\tiny{$\al_1$}}
\put(13,-5){\tiny{$\al_2$}}
\put(28,-5){\tiny{$\al_3$}}
\put(43,-5){\tiny{$\al_4$}}
\end{picture}}

\newcommand{\smallHecFourTwentyoneSeven}{
\setlength{\unitlength}{0.5mm}
\begin{picture}(30,20)(00,00)

\put(-5,16){\tiny{$a:=[\acchi^{(3,18)}_7]$}}

\multiput(15.8,3.8)(00,00){1}{\oval(6.4,6.4)[br]}
\multiput(32,3.8)(00,00){1}{\oval(26,13)[t]}
\multiput(45,01)(00,00){1}{\line(0,1){2.8}}

\multiput(00,00)(15,00){4}{\circle{2}}
\multiput(01,00)(15,00){3}{\line(1,0){13}}

\put(-1,3){\tiny{$\zeta$}}
\put(14,3){\tiny{$\zeta$}}
\put(28,3){\tiny{$-1$}}
\put(46,2){\tiny{$\zeta$}}

\put(4,2){\tiny{$\zeta^{-1}$}}
\put(19.5,2){\tiny{$\zeta^{-1}$}}
\put(35,2){\tiny{$\zeta^{-1}$}}
\put(31,12){\tiny{$\zeta^{-1}$}}

\put(-2,-5){\tiny{$\al_1$}}
\put(13,-5){\tiny{$\al_2$}}
\put(28,-5){\tiny{$\al_3$}}
\put(43,-5){\tiny{$\al_4$}}
\end{picture}}

\newcommand{\RankFourHecTwentyoneTau}{
\setlength{\unitlength}{0.5mm}
\begin{picture}(140,30)(0,0)

\put(23,23){\footnotesize{${\bar{a}}$}}
\put(46,23){\footnotesize{$c$}}
\put(69,23){\footnotesize{$b$}}
\put(92,23){\footnotesize{$a$}}

\multiput(29.5,24)(23,0){3}{\vector(1,0){13}}
\multiput(41.5,24)(23,0){3}{\vector(-1,0){13}}

\multiput(24,20)(23,0){3}{\vector(0,-1){13}}
\multiput(24,8)(23,0){3}{\vector(0,1){13}}

\put(23,00){\footnotesize{${\bar{b}}$}}
\put(46,00){\footnotesize{$e$}}
\put(69,00){\footnotesize{$d$}}
\multiput(29.5,1)(23,0){2}{\vector(1,0){13}}
\multiput(41.5,1)(23,0){2}{\vector(-1,0){13}}

\put(33,26){\footnotesize{${\bar{\tau}}_1$}}
\put(56,26){\footnotesize{${\bar{\tau}}_2$}}
\put(79,26){\footnotesize{${\bar{\tau}}_3$}}

\put(17,12){\footnotesize{${\bar{\tau}}_4$}}
\multiput(48,12)(23,0){2}{\footnotesize{${\bar{\tau}}_4$}}

\put(33,03){\footnotesize{${\bar{\tau}}_1$}}
\put(56,03){\footnotesize{${\bar{\tau}}_2$}}

\end{picture}}

\newcommand{\allHecFourTewntytwo}{
\setlength{\unitlength}{1mm}
\begin{picture}(220,60)(-5,10)
\put(10,40){$\smallHecFourTwentytwoZero$}
\put(40,40){$\smallHecFourTwentytwoOne$}
\put(70,40){$\smallHecFourTwentytwoTwo$}
\put(100,40){$\smallHecFourTwentytwoThree$}
\put(10,20){$\smallHecFourTwentytwoFour$}
\put(40,20){$\smallHecFourTwentytwoFive$}
\put(70,20){$\smallHecFourTwentytwoSix$}
\put(100,20){$\smallHecFourTwentytwoSeven$}
\put(10,10){$t^\prime:=t\circ\left[{{1234}\atop{1432}}\right]$} 
\put(80,10){$(\zeta\in\bKt,\,\zeta^2+1=0)$}
\end{picture}}

\newcommand{\smallHecFourTwentytwoZero}{
\setlength{\unitlength}{0.5mm}
\begin{picture}(30,20)(00,00)

\put(-5,16){\tiny{$a:=[\acchi^{(3,22)}_1]$}}

\multiput(00,00)(15,00){4}{\circle{2}}
\multiput(01,00)(15,00){3}{\line(1,0){13}}

\put(-3,3){\tiny{$-\zeta$}}
\put(11.5,3){\tiny{$-1$}}
\put(28,3){\tiny{$\zeta$}}
\put(42,3){\tiny{$-\zeta$}}

\put(7,2){\tiny{$\zeta$}}
\put(20,2){\tiny{$-\zeta$}}
\put(36,2){\tiny{$\zeta$}}

\put(-2,-5){\tiny{$\al_1$}}
\put(13,-5){\tiny{$\al_2$}}
\put(28,-5){\tiny{$\al_3$}}
\put(43,-5){\tiny{$\al_4$}}
\end{picture}}

\newcommand{\smallHecFourTwentytwoOne}{
\setlength{\unitlength}{0.5mm}
\begin{picture}(30,20)(00,00)

\put(-5,16){\tiny{$b:=[\acchi^{(3,22)}_2]$}}

\multiput(00,00)(15,00){4}{\circle{2}}
\multiput(01,00)(15,00){3}{\line(1,0){13}}

\put(-3,3){\tiny{$-1$}}
\put(11.5,3){\tiny{$-1$}}
\put(26.5,3){\tiny{$-1$}}
\put(42,3){\tiny{$-\zeta$}}

\put(5,2){\tiny{$-\zeta$}}
\put(22,2){\tiny{$\zeta$}}
\put(36,2){\tiny{$\zeta$}}

\put(-2,-5){\tiny{$\al_1$}}
\put(13,-5){\tiny{$\al_2$}}
\put(28,-5){\tiny{$\al_3$}}
\put(43,-5){\tiny{$\al_4$}}
\end{picture}}

\newcommand{\smallHecFourTwentytwoTwo}{
\setlength{\unitlength}{0.5mm}
\begin{picture}(30,20)(00,00)

\put(-5,16){\tiny{${\bar{b}}:=[\acchi^{(3,22)}_3]$}}

\multiput(00,00)(15,00){4}{\circle{2}}
\multiput(01,00)(15,00){3}{\line(1,0){13}}

\put(-3,3){\tiny{$-1$}}
\put(12,3){\tiny{$-\zeta$}}
\put(28,3){\tiny{$-1$}}
\put(42,3){\tiny{$-\zeta$}}

\put(6,2){\tiny{$\zeta$}}
\put(21,2){\tiny{$\zeta$}}
\put(36,2){\tiny{$\zeta$}}

\put(-2,-5){\tiny{$\al_1$}}
\put(13,-5){\tiny{$\al_2$}}
\put(28,-5){\tiny{$\al_3$}}
\put(43,-5){\tiny{$\al_4$}}
\end{picture}}

\newcommand{\smallHecFourTwentytwoThree}{
\setlength{\unitlength}{0.5mm}
\begin{picture}(30,20)(00,00)

\put(-5,16){\tiny{$c:=[\acchi^{(3,22)}_4]$}}

\multiput(15.8,3.8)(00,00){1}{\oval(6.4,6.4)[br]}
\multiput(32,3.8)(00,00){1}{\oval(26,13)[t]}
\multiput(45,01)(00,00){1}{\line(0,1){2.8}}

\multiput(00,00)(15,00){4}{\circle{2}}
\multiput(01,00)(15,00){3}{\line(1,0){13}}

\put(-3,3){\tiny{$-1$}}
\put(13,3){\tiny{$\zeta$}}
\put(28,3){\tiny{$-1$}}
\put(46,2){\tiny{$-1$}}

\put(6,2){\tiny{$-\zeta$}}
\put(19.5,2){\tiny{$-\zeta$}}
\put(35,2){\tiny{$-\zeta$}}
\put(29,12){\tiny{$-1$}}

\put(-2,-5){\tiny{$\al_1$}}
\put(13,-5){\tiny{$\al_2$}}
\put(28,-5){\tiny{$\al_3$}}
\put(43,-5){\tiny{$\al_4$}}
\end{picture}}

\newcommand{\smallHecFourTwentytwoFour}{
\setlength{\unitlength}{0.5mm}
\begin{picture}(30,20)(00,00)

\put(-5,16){\tiny{${\bar{a}}:=[\acchi^{(3,22)}_5]$}}

\multiput(15.8,3.8)(00,00){1}{\oval(6.4,6.4)[br]}
\multiput(32,3.8)(00,00){1}{\oval(26,13)[t]}
\multiput(45,01)(00,00){1}{\line(0,1){2.8}}

\multiput(00,00)(15,00){4}{\circle{2}}
\multiput(01,00)(15,00){3}{\line(1,0){13}}

\put(-3,3){\tiny{$-1$}}
\put(11.5,3){\tiny{$-1$}}
\put(28,3){\tiny{$-1$}}
\put(46,2){\tiny{$-1$}}

\put(7,2){\tiny{$\zeta$}}
\put(19.5,2){\tiny{$-\zeta$}}
\put(35,2){\tiny{$-\zeta$}}
\put(29,12){\tiny{$-1$}}

\put(-2,-5){\tiny{$\al_1$}}
\put(13,-5){\tiny{$\al_2$}}
\put(28,-5){\tiny{$\al_3$}}
\put(43,-5){\tiny{$\al_4$}}
\end{picture}}

\newcommand{\smallHecFourTwentytwoFive}{
\setlength{\unitlength}{0.5mm}
\begin{picture}(30,20)(00,00)

\put(-5,16){\tiny{$d:=[\acchi^{(3,22)}_6]\circ\left[{{1234}\atop{1243}}\right]$}}

\multiput(00,00)(15,00){4}{\circle{2}}
\multiput(01,00)(15,00){3}{\line(1,0){13}}

\put(-3,3){\tiny{$-1$}}
\put(13,3){\tiny{$\zeta$}}
\put(28,3){\tiny{$-1$}}
\put(42,3){\tiny{$-\zeta$}}

\put(5,2){\tiny{$-\zeta$}}
\put(20,2){\tiny{$-1$}}
\put(36,2){\tiny{$\zeta$}}

\put(-2,-5){\tiny{$\al_1$}}
\put(13,-5){\tiny{$\al_2$}}
\put(28,-5){\tiny{$\al_4$}}
\put(43,-5){\tiny{$\al_3$}}
\end{picture}}

\newcommand{\smallHecFourTwentytwoSix}{
\setlength{\unitlength}{0.5mm}
\begin{picture}(30,20)(00,00)

\put(-5,16){\tiny{${\bar{d}}:=[\acchi^{(3,22)}_7]$}}

\multiput(00,00)(15,00){4}{\circle{2}}
\multiput(01,00)(15,00){3}{\line(1,0){13}}

\put(-3,3){\tiny{$-1$}}
\put(11.5,3){\tiny{$-1$}}
\put(26.5,3){\tiny{$-1$}}
\put(42,3){\tiny{$-\zeta$}}

\put(7,2){\tiny{$\zeta$}}
\put(19.5,2){\tiny{$-1$}}
\put(36,2){\tiny{$\zeta$}}

\put(-2,-5){\tiny{$\al_1$}}
\put(13,-5){\tiny{$\al_2$}}
\put(28,-5){\tiny{$\al_3$}}
\put(43,-5){\tiny{$\al_4$}}
\end{picture}}

\newcommand{\smallHecFourTwentytwoSeven}{
\setlength{\unitlength}{0.5mm}
\begin{picture}(30,20)(00,00)

\put(-5,19){\tiny{${\bar{c}}:=[\acchi^{(3,22)}_8]\circ\left[{{1234}\atop{3214}}\right]$}}

\multiput(15.8,3.8)(00,00){1}{\oval(6.4,6.4)[br]}
\multiput(32,3.8)(00,00){1}{\oval(26,13)[t]}
\multiput(45,01)(00,00){1}{\line(0,1){2.8}}

\multiput(00,00)(15,00){4}{\circle{2}}
\multiput(01,00)(15,00){3}{\line(1,0){13}}

\put(-3,3){\tiny{$-\zeta$}}
\put(11.5,3){\tiny{$-1$}}
\put(29,3){\tiny{$\zeta$}}
\put(46,2){\tiny{$-1$}}

\put(7,2){\tiny{$\zeta$}}
\put(20,2){\tiny{$-\zeta$}}
\put(34,2){\tiny{$-\zeta$}}
\put(29,12){\tiny{$-1$}}

\put(-2,-5){\tiny{$\al_3$}}
\put(13,-5){\tiny{$\al_2$}}
\put(28,-5){\tiny{$\al_1$}}
\put(43,-5){\tiny{$\al_4$}}
\end{picture}}

\newcommand{\RankFourHecTwentytwoTau}{
\setlength{\unitlength}{0.5mm}
\begin{picture}(140,70)(0,0)
\put(0,69){\footnotesize{$a$}}
\put(23,69){\footnotesize{$b$}}
\put(46,69){\footnotesize{${\bar{b}}$}}
\multiput(6.5,70)(23,0){2}{\vector(1,0){13}}
\multiput(18.5,70)(23,0){2}{\vector(-1,0){13}}

\multiput(24,66)(23,0){2}{\vector(0,-1){13}}
\multiput(24,54)(23,0){2}{\vector(0,1){13}}

\put(23,46){\footnotesize{$c$}}
\put(46,46){\footnotesize{${\bar{a}}$}}
\put(69,46){\footnotesize{${\bar{c}}$}}
\put(92,46){\footnotesize{${\bar{d}}^\prime$}}
\put(115,46){\footnotesize{$d^\prime$}}
\multiput(29.5,47)(23,0){4}{\vector(1,0){13}}
\multiput(41.5,47)(23,0){4}{\vector(-1,0){13}}

\multiput(24,43)(23,0){2}{\vector(0,-1){13}}
\multiput(24,31)(23,0){2}{\vector(0,1){13}}
\multiput(93,43)(23,0){2}{\vector(0,-1){13}}
\multiput(93,31)(23,0){2}{\vector(0,1){13}}

\put(23,23){\footnotesize{$d$}}
\put(46,23){\footnotesize{${\bar{d}}$}}
\put(69,23){\footnotesize{${\bar{c}}^\prime$}}
\put(92,23){\footnotesize{${\bar{a}}^\prime$}}
\put(115,23){\footnotesize{$c^\prime$}}
\multiput(29.5,24)(23,0){4}{\vector(1,0){13}}
\multiput(41.5,24)(23,0){4}{\vector(-1,0){13}}

\multiput(93,20)(23,0){2}{\vector(0,-1){13}}
\multiput(93,8)(23,0){2}{\vector(0,1){13}}

\put(92,00){\footnotesize{${\bar{b}}^\prime$}}
\put(115,00){\footnotesize{$b^\prime$}}
\put(138,00){\footnotesize{$a^\prime$}}
\multiput(98.5,1)(23,0){2}{\vector(1,0){13}}
\multiput(110.5,1)(23,0){2}{\vector(-1,0){13}}

\put(10,72){\footnotesize{${\bar{\tau}}_2$}}
\put(33,72){\footnotesize{${\bar{\tau}}_1$}}

\put(17,58){\footnotesize{${\bar{\tau}}_3$}}
\put(48,58){\footnotesize{${\bar{\tau}}_3$}}

\put(33,49){\footnotesize{${\bar{\tau}}_1$}}
\put(56,49){\footnotesize{${\bar{\tau}}_2$}}
\put(79,49){\footnotesize{${\bar{\tau}}_4$}}
\put(102,49){\footnotesize{${\bar{\tau}}_1$}}

\put(17,35){\footnotesize{${\bar{\tau}}_4$}}
\put(48,35){\footnotesize{${\bar{\tau}}_4$}}

\put(86,35){\footnotesize{${\bar{\tau}}_2$}}
\put(117,35){\footnotesize{${\bar{\tau}}_2$}}

\put(33,26){\footnotesize{${\bar{\tau}}_1$}}
\put(56,26){\footnotesize{${\bar{\tau}}_2$}}
\put(79,26){\footnotesize{${\bar{\tau}}_4$}}
\put(102,26){\footnotesize{${\bar{\tau}}_1$}}

\put(86,12){\footnotesize{${\bar{\tau}}_3$}}
\put(117,12){\footnotesize{${\bar{\tau}}_3$}}

\put(102,03){\footnotesize{${\bar{\tau}}_1$}}
\put(125,03){\footnotesize{${\bar{\tau}}_4$}}

\end{picture}}

\newcommand{\allHecHigherquasiCartan}{
\setlength{\unitlength}{1mm}
\begin{picture}(220,110)(-25,-10)

\put(0,105){$(1)$}\put(20,105){$\smallHecHigherOneZero$}
\put(0,90){$(2)$}\put(20,90){$\smallHecHigherTwoZero$}
\put(0,75){$(3)$}\put(20,75){$\smallHecHigherThreeZero$}
\put(0,60){$(4)$}\put(20,60){$\smallHecHigherFourZero$}
\put(0,45){$(5)$}\put(20,45){$\smallHecHigherFiveZero$}
\put(0,30){$(6)$}\put(20,30){$\smallHecHigherSixZero$}
\put(0,15){$(7)$}\put(20,15){$\smallHecHigherSevenZero$}
\put(0,00){$(8)$}\put(20,00){$\smallHecHigherEightZero$}
\end{picture}}

\newcommand{\smallHecHigherOneZero}{
\setlength{\unitlength}{0.5mm}
\begin{picture}(230,20)(00,00)

\multiput(00,00)(15,00){2}{\circle{2}}
\multiput(01,00)(15,00){2}{\line(1,0){13}}
\multiput(32.5,00)(2,00){6}{\circle*{0.6}}
\multiput(46,00)(15,00){2}{\line(1,0){13}}
\multiput(60,00)(00,00){2}{\circle{2}}
\multiput(77.5,00)(2,00){6}{\circle*{0.6}}
\multiput(91,00)(15,00){3}{\line(1,0){13}}
\multiput(105,00)(15,00){3}{\circle{2}}

\put(-1,3){\tiny{$\chamq$}}
\put(3.5,2){\tiny{$\chamq^{-1}$}}
\put(14,3){\tiny{$\chamq$}}
\put(18.5,2){\tiny{$\chamq^{-1}$}}
\put(48.5,2){\tiny{$\chamq^{-1}$}}
\put(59,3){\tiny{$\chamq$}}
\put(63.5,2){\tiny{$\chamq^{-1}$}}
\put(93.5,2){\tiny{$\chamq^{-1}$}}
\put(104,3){\tiny{$\chamq$}}
\put(108.5,2){\tiny{$\chamq^{-1}$}}
\put(119,3){\tiny{$\chamq$}}
\put(123.5,2){\tiny{$\chamq^{-1}$}}
\put(134,3){\tiny{$\chamq$}}

\put(-2,-5){\tiny{$\al_1$}}
\put(13,-5){\tiny{$\al_2$}}
\put(58,-5){\tiny{$\al_m$}}
\put(103,-5){\tiny{$\al_{\chamhatN-2}$}}
\put(118,-5){\tiny{$\al_{\chamhatN-1}$}}
\put(133,-5){\tiny{$\al_{\chamhatN}$}}

\put(-2,-13){\tiny{$(\chamhatN\geq 5,\,\chamq\in\bKt,\,\chamq\ne 1\,(1\leq m\leq\chamhatN))$}}

\end{picture}}

\newcommand{\smallHecHigherTwoZero}{
\setlength{\unitlength}{0.5mm}
\begin{picture}(230,20)(00,00)

\multiput(00,00)(15,00){2}{\circle{2}}
\multiput(01,00)(15,00){2}{\line(1,0){13}}
\multiput(32.5,00)(2,00){6}{\circle*{0.6}}
\multiput(46,00)(15,00){2}{\line(1,0){13}}
\multiput(60,00)(00,00){2}{\circle{2}}
\multiput(77.5,00)(2,00){6}{\circle*{0.6}}
\multiput(91,00)(15,00){3}{\line(1,0){13}}
\multiput(105,00)(15,00){3}{\circle{2}}

\put(-1,3){\tiny{$\chamq$}}
\put(3.5,2){\tiny{$\chamq^{-1}$}}
\put(14,3){\tiny{$\chamq$}}
\put(18.5,2){\tiny{$\chamq^{-1}$}}
\put(48.5,2){\tiny{$\chamq^{-1}$}}
\put(58,3){\tiny{$-1$}}
\put(66,2){\tiny{$\chamq$}}
\put(96,2){\tiny{$\chamq$}}
\put(102,3){\tiny{$\chamq^{-1}$}}
\put(111,2){\tiny{$\chamq$}}
\put(117,3){\tiny{$\chamq^{-1}$}}
\put(126,2){\tiny{$\chamq$}}
\put(132,3){\tiny{$\chamq^{-1}$}}

\put(-2,-5){\tiny{$\al_1$}}
\put(13,-5){\tiny{$\al_2$}}
\put(58,-5){\tiny{$\al_m$}}
\put(103,-5){\tiny{$\al_{\chamhatN-2}$}}
\put(118,-5){\tiny{$\al_{\chamhatN-1}$}}
\put(133,-5){\tiny{$\al_{\chamhatN}$}}

\put(-2,-13){\tiny{$(\chamhatN\geq 5,\,\chamq\in\bKt,\,\chamq^2\ne 1,\,1\leq m\leq{\frac {\chamhatN+1} 2})$}}

\end{picture}}

\newcommand{\smallHecHigherThreeZero}{
\setlength{\unitlength}{0.5mm}
\begin{picture}(230,20)(00,00)

\multiput(00,00)(15,00){2}{\circle{2}}
\multiput(01,00)(15,00){2}{\line(1,0){13}}
\multiput(32.5,00)(2,00){6}{\circle*{0.6}}
\multiput(46,00)(15,00){2}{\line(1,0){13}}
\multiput(60,00)(00,00){2}{\circle{2}}
\multiput(77.5,00)(2,00){6}{\circle*{0.6}}
\multiput(91,00)(15,00){3}{\line(1,0){13}}
\multiput(105,00)(15,00){3}{\circle{2}}

\put(-1,3){\tiny{$\chamq^2$}}
\put(3.5,2){\tiny{$\chamq^{-2}$}}
\put(14,3){\tiny{$\chamq^2$}}
\put(18.5,2){\tiny{$\chamq^{-2}$}}
\put(48.5,2){\tiny{$\chamq^{-2}$}}
\put(59,3){\tiny{$\chamq^2$}}
\put(63.5,2){\tiny{$\chamq^{-2}$}}
\put(93.5,2){\tiny{$\chamq^{-2}$}}
\put(104,3){\tiny{$\chamq^2$}}
\put(108.5,2){\tiny{$\chamq^{-2}$}}
\put(119,3){\tiny{$\chamq^2$}}
\put(123.5,2){\tiny{$\chamq^{-2}$}}
\put(134,3){\tiny{$\chamq$}}

\put(-2,-5){\tiny{$\al_1$}}
\put(13,-5){\tiny{$\al_2$}}
\put(58,-5){\tiny{$\al_m$}}
\put(103,-5){\tiny{$\al_{\chamhatN-2}$}}
\put(118,-5){\tiny{$\al_{\chamhatN-1}$}}
\put(133,-5){\tiny{$\al_{\chamhatN}$}}

\put(-2,-13){\tiny{$(\chamhatN\geq 5,\,\chamq\in\bKt,\,\chamq^2\ne 1,
\,(1\leq m\leq\chamhatN-1))$}}

\end{picture}}

\newcommand{\smallHecHigherFourZero}{
\setlength{\unitlength}{0.5mm}
\begin{picture}(230,20)(00,00)

\multiput(00,00)(15,00){2}{\circle{2}}
\multiput(01,00)(15,00){2}{\line(1,0){13}}
\multiput(32.5,00)(2,00){6}{\circle*{0.6}}
\multiput(46,00)(15,00){2}{\line(1,0){13}}
\multiput(60,00)(00,00){2}{\circle{2}}
\multiput(77.5,00)(2,00){6}{\circle*{0.6}}
\multiput(91,00)(15,00){3}{\line(1,0){13}}
\multiput(105,00)(15,00){3}{\circle{2}}

\put(-3,3){\tiny{$\chamq^{-2}$}}
\put(6,2){\tiny{$\chamq^2$}}
\put(12,3){\tiny{$\chamq^{-2}$}}
\put(21,2){\tiny{$\chamq^2$}}
\put(51,2){\tiny{$\chamq^2$}}
\put(57,3){\tiny{$-1$}}
\put(64.5,2){\tiny{$\chamq^{-2}$}}
\put(93.5,2){\tiny{$\chamq^{-2}$}}
\put(104,3){\tiny{$\chamq^2$}}
\put(108.5,2){\tiny{$\chamq^{-2}$}}
\put(119,3){\tiny{$\chamq^2$}}
\put(123.5,2){\tiny{$\chamq^{-2}$}}
\put(134,3){\tiny{$\chamq$}}

\put(-2,-5){\tiny{$\al_1$}}
\put(13,-5){\tiny{$\al_2$}}
\put(58,-5){\tiny{$\al_m$}}
\put(103,-5){\tiny{$\al_{\chamhatN-2}$}}
\put(118,-5){\tiny{$\al_{\chamhatN-1}$}}
\put(133,-5){\tiny{$\al_{\chamhatN}$}}

\put(-2,-13){\tiny{$(\chamhatN\geq 5,\,\chamq\in\bKt,\,\chamq^4\ne 1,\,\,1\leq m\leq{\frac {\chamhatN+1} 2})$}}
\end{picture}}

\newcommand{\smallHecHigherFiveZero}{
\setlength{\unitlength}{0.5mm}
\begin{picture}(230,20)(00,00)

\multiput(00,00)(15,00){2}{\circle{2}}
\multiput(01,00)(15,00){2}{\line(1,0){13}}
\multiput(32.5,00)(2,00){6}{\circle*{0.6}}
\multiput(46,00)(15,00){2}{\line(1,0){13}}
\multiput(60,00)(00,00){2}{\circle{2}}
\multiput(77.5,00)(2,00){6}{\circle*{0.6}}
\multiput(91,00)(15,00){3}{\line(1,0){13}}
\multiput(105,00)(15,00){3}{\circle{2}}

\put(-5,4.5){\tiny{$-{\frac 1 \zeta}$}}
\put(3,2){\tiny{$-\zeta$}}
\put(10,4.5){\tiny{$-{\frac 1 \zeta}$}}
\put(18,2){\tiny{$-\zeta$}}
\put(48,2){\tiny{$-\zeta$}}
\put(55,4.5){\tiny{$-{\frac 1 \zeta}$}}
\put(63,2){\tiny{$-\zeta$}}
\put(93,2){\tiny{$-\zeta$}}
\put(100,4.5){\tiny{$-{\frac 1 \zeta}$}}
\put(108,2){\tiny{$-\zeta$}}
\put(115,4.5){\tiny{$-{\frac 1 \zeta}$}}
\put(123,2){\tiny{$-\zeta$}}
\put(133,3){\tiny{$\zeta$}}

\put(-2,-5){\tiny{$\al_1$}}
\put(13,-5){\tiny{$\al_2$}}
\put(58,-5){\tiny{$\al_m$}}
\put(103,-5){\tiny{$\al_{\chamhatN-2}$}}
\put(118,-5){\tiny{$\al_{\chamhatN-1}$}}
\put(133,-5){\tiny{$\al_{\chamhatN}$}}

\put(-2,-13){\tiny{$(\chamhatN\geq 5,\,\zeta\in\bKt,\,\zeta^2+\zeta+1=0\,(1\leq m\leq\chamhatN-1))$}}
\end{picture}}

\newcommand{\smallHecHigherSixZero}{
\setlength{\unitlength}{0.5mm}
\begin{picture}(230,20)(00,00)

\multiput(00,00)(15,00){2}{\circle{2}}
\multiput(01,00)(15,00){2}{\line(1,0){13}}
\multiput(32.5,00)(2,00){6}{\circle*{0.6}}
\multiput(46,00)(15,00){2}{\line(1,0){13}}
\multiput(60,00)(00,00){2}{\circle{2}}
\multiput(77.5,00)(2,00){6}{\circle*{0.6}}
\multiput(91,00)(15,00){3}{\line(1,0){13}}
\multiput(105,00)(15,00){3}{\circle{2}}

\put(-4,3){\tiny{$-\zeta$}}
\put(2.5,3.5){\tiny{$-{\frac 1 \zeta}$}}
\put(11,3){\tiny{$-\zeta$}}
\put(17.5,3.5){\tiny{$-{\frac 1 \zeta}$}}
\put(47.5,3.5){\tiny{$-{\frac 1 \zeta}$}}
\put(56,3){\tiny{$-1$}}
\put(63,2){\tiny{$-\zeta$}}
\put(93,2){\tiny{$-\zeta$}}
\put(100,4.5){\tiny{$-{\frac 1 \zeta}$}}
\put(108,2){\tiny{$-\zeta$}}
\put(115,4.5){\tiny{$-{\frac 1 \zeta}$}}
\put(123,2){\tiny{$-\zeta$}}
\put(133,3){\tiny{$\zeta$}}

\put(-2,-5){\tiny{$\al_1$}}
\put(13,-5){\tiny{$\al_2$}}
\put(58,-5){\tiny{$\al_m$}}
\put(103,-5){\tiny{$\al_{\chamhatN-2}$}}
\put(118,-5){\tiny{$\al_{\chamhatN-1}$}}
\put(133,-5){\tiny{$\al_{\chamhatN}$}}

\put(-2,-13){\tiny{$(\chamhatN\geq 5,\,\zeta\in\bKt,\,\zeta^2+\zeta+1=0,\,1\leq m\leq\chamhatN-1)$}}
\end{picture}}

\newcommand{\smallHecHigherSevenZero}{
\setlength{\unitlength}{0.5mm}
\begin{picture}(230,20)(00,00)

\multiput(00,00)(15,00){2}{\circle{2}}
\multiput(01,00)(15,00){2}{\line(1,0){13}}
\multiput(32.5,00)(2,00){6}{\circle*{0.6}}
\multiput(46,00)(15,00){2}{\line(1,0){13}}
\multiput(60,00)(00,00){2}{\circle{2}}
\multiput(77.5,00)(2,00){6}{\circle*{0.6}}
\multiput(91,00)(15,00){3}{\line(1,0){13}}
\multiput(105,00)(15,00){3}{\circle{2}}

\put(-1,3){\tiny{$\chamq$}}
\put(3.5,2){\tiny{$\chamq^{-1}$}}
\put(14,3){\tiny{$\chamq$}}
\put(18.5,2){\tiny{$\chamq^{-1}$}}
\put(48.5,2){\tiny{$\chamq^{-1}$}}
\put(59,3){\tiny{$\chamq$}}
\put(63.5,2){\tiny{$\chamq^{-1}$}}
\put(93.5,2){\tiny{$\chamq^{-1}$}}
\put(104,3){\tiny{$\chamq$}}
\put(108.5,2){\tiny{$\chamq^{-1}$}}
\put(119,3){\tiny{$\chamq$}}
\put(123.5,2){\tiny{$\chamq^{-2}$}}
\put(134,3){\tiny{$\chamq^2$}}

\put(-2,-5){\tiny{$\al_1$}}
\put(13,-5){\tiny{$\al_2$}}
\put(58,-5){\tiny{$\al_m$}}
\put(103,-5){\tiny{$\al_{\chamhatN-2}$}}
\put(118,-5){\tiny{$\al_{\chamhatN-1}$}}
\put(133,-5){\tiny{$\al_{\chamhatN}$}}

\put(-2,-13){\tiny{$(\chamhatN\geq 5,\,\chamq\in\bKt,\,\chamq^2\ne 1\,(1\leq m\leq\chamhatN-2))$}}
\end{picture}}

\newcommand{\smallHecHigherEightZero}{
\setlength{\unitlength}{0.5mm}
\begin{picture}(230,20)(00,00)

\multiput(105.8,3.8)(00,00){1}{\oval(6.4,6.4)[br]}
\multiput(122,3.8)(00,00){1}{\oval(26,13)[t]}
\multiput(135,01)(00,00){1}{\line(0,1){2.8}}

\multiput(00,00)(15,00){2}{\circle{2}}
\multiput(01,00)(15,00){2}{\line(1,0){13}}
\multiput(32.5,00)(2,00){6}{\circle*{0.6}}
\multiput(46,00)(15,00){2}{\line(1,0){13}}
\multiput(60,00)(00,00){2}{\circle{2}}
\multiput(77.5,00)(2,00){6}{\circle*{0.6}}
\multiput(91,00)(15,00){2}{\line(1,0){13}}
\multiput(105,00)(15,00){3}{\circle{2}}

\put(-1,3){\tiny{$\chamq$}}
\put(3.5,2){\tiny{$\chamq^{-1}$}}
\put(14,3){\tiny{$\chamq$}}
\put(18.5,2){\tiny{$\chamq^{-1}$}}
\put(48.5,2){\tiny{$\chamq^{-1}$}}
\put(59,3){\tiny{$\chamq$}}
\put(63.5,2){\tiny{$\chamq^{-1}$}}
\put(93.5,2){\tiny{$\chamq^{-1}$}}
\put(104,3){\tiny{$\chamq$}}
\put(110,2){\tiny{$\chamq^{-1}$}}
\put(119,3){\tiny{$\chamq$}}

\put(136,3){\tiny{$\chamq$}}

\put(119,12){\tiny{$\chamq^{-1}$}}

\put(-2,-5){\tiny{$\al_1$}}
\put(13,-5){\tiny{$\al_2$}}
\put(58,-5){\tiny{$\al_m$}}
\put(103,-5){\tiny{$\al_{\chamhatN-2}$}}
\put(118,-5){\tiny{$\al_{\chamhatN-1}$}}
\put(133,-5){\tiny{$\al_{\chamhatN}$}}

\put(-2,-13){\tiny{$(\chamhatN\geq 5,\,\chamq\in\bKt,\,\chamq\ne 1\,(1\leq m\leq\chamhatN-3))$}}

\end{picture}}

\newcommand{\allHecHigherNineTen}{
\setlength{\unitlength}{1mm}
\begin{picture}(220,70)(-15,-10)
\put(10,60){$\smallHecHigherNineTenZero$}
\put(10,45){$\smallHecHigherNineTenOne$}
\put(10,30){$\smallHecHigherNineTenTwo$}
\put(10,15){$\smallHecHigherNineTenThree$}
\put(10,00){$\smallHecHigherNineTenFour$}

\multiput(105,29)(00,00){1}{\oval(5,28)[r]}
\multiput(105,15)(00,28){2}{\vector(-1,0){3}}
\multiput(109,35)(00,00){1}{{\tiny{$\chamhattau_{\chamhatN-1}$
$(\chamq_{\chamhatN-1}=-1)$}}}

\multiput(115,22)(00,00){1}{\oval(5,14)[r]}
\multiput(115,15)(00,14){2}{\vector(-1,0){3}}
\multiput(119,22)(00,00){1}{{\tiny{$\chamhattau_{\chamhatN-2}$
$(\chamq_{\chamhatN-2}=-1)$}}}

\multiput(88,7)(00,00){1}{\oval(5,14)[r]}
\multiput(88,00)(00,14){2}{\vector(-1,0){3}}
\multiput(92,7)(00,00){1}{{\tiny{$\chamhattau_\chamhatN$
$(\chamq_{\chamhatN-1}=-1)$}}}
\end{picture}}

\newcommand{\smallHecHigherNineTenZero}{
\setlength{\unitlength}{0.5mm}
\begin{picture}(230,20)(00,00)

\put(-5,13){\tiny{$a:=\champichi(\chi)$, $\chamchi:=\chamchi^{\chamhatpzero}$,
$r_i=q^{-1}$ $(i\in\fkJ_{0,m})$, $r_j=q$ $(j\in\fkJ_{m+1,\chamhatN-1})$}}

\multiput(00,00)(15,00){2}{\circle{2}}
\multiput(01,00)(15,00){2}{\line(1,0){13}}
\multiput(32.5,00)(2,00){6}{\circle*{0.6}}
\multiput(46,00)(15,00){2}{\line(1,0){13}}
\multiput(60,00)(00,00){2}{\circle{2}}
\multiput(77.5,00)(2,00){6}{\circle*{0.6}}
\multiput(91,00)(15,00){3}{\line(1,0){13}}
\multiput(105,00)(15,00){3}{\circle{2}}

\put(-1,3){\tiny{$\chamq$}}
\put(3.5,2){\tiny{$\chamq^{-1}$}}
\put(14,3){\tiny{$\chamq$}}
\put(18.5,2){\tiny{$\chamq^{-1}$}}
\put(48.5,2){\tiny{$\chamq^{-1}$}}
\put(58,3){\tiny{$-1$}}
\put(66,2){\tiny{$\chamq$}}
\put(96,2){\tiny{$\chamq$}}
\put(102,3){\tiny{$\chamq^{-1}$}}
\put(111,2){\tiny{$\chamq$}}
\put(117,3){\tiny{$\chamq^{-1}$}}
\put(126,2){\tiny{$\chamq^2$}}
\put(132,3){\tiny{$\chamq^{-2}$}}

\put(-2,-5){\tiny{$\al_1$}}
\put(13,-5){\tiny{$\al_2$}}
\put(58,-5){\tiny{$\al_m$}}
\put(103,-5){\tiny{$\al_{\chamhatN-2}$}}
\put(118,-5){\tiny{$\al_{\chamhatN-1}$}}
\put(133,-5){\tiny{$\al_{\chamhatN}$}}
\end{picture}}

\newcommand{\smallHecHigherNineTenOne}{
\setlength{\unitlength}{0.5mm}
\begin{picture}(230,20)(00,00)

\put(-5,15){\tiny{$\chamhatp(\chamhatN)=1$,
$r_{\chamhatN-1}=\chamq$}}

\multiput(00,00)(15,00){2}{\circle{2}}
\multiput(01,00)(15,00){2}{\line(1,0){13}}
\multiput(32.5,00)(2,00){6}{\circle*{0.6}}
\multiput(46,00)(15,00){2}{\line(1,0){13}}
\multiput(60,00)(00,00){2}{\circle{2}}
\multiput(77.5,00)(2,00){6}{\circle*{0.6}}
\multiput(91,00)(15,00){3}{\line(1,0){13}}
\multiput(105,00)(15,00){3}{\circle{2}}

\put(-3,3){\tiny{$\chamq_1$}}
\put(5,2){\tiny{$r_1$}}
\put(13,3){\tiny{$\chamq_2$}}
\put(20,2){\tiny{$r_2$}}
\put(47,2){\tiny{$r_{i-1}$}}
\put(57,3){\tiny{$\chamq_i$}}
\put(62,2){\tiny{$r_{i+1}$}}
\put(91,4){\tiny{$r_{\chamhatN-2}$}}
\put(102,9){\tiny{$\chamq_{\chamhatN-2}$}}
\put(106,4){\tiny{$r_{\chamhatN-1}$}}
\put(117,9){\tiny{$\chamq_{\chamhatN-1}$}}
\put(125,2){\tiny{$\chamq^2$}}
\put(132,3){\tiny{$\chamq^{-2}$}}

\put(-2,-5){\tiny{$\al_1$}}
\put(13,-5){\tiny{$\al_2$}}
\put(58,-5){\tiny{$\al_i$}}
\put(103,-5){\tiny{$\al_{\chamhatN-2}$}}
\put(118,-5){\tiny{$\al_{\chamhatN-1}$}}
\put(133,-5){\tiny{$\al_{\chamhatN}$}}

\end{picture}}

\newcommand{\smallHecHigherNineTenTwo}{
\setlength{\unitlength}{0.5mm}
\begin{picture}(230,20)(00,00)

\put(-5,16){\tiny{$\chamhatp(\chamhatN-1)=\chamhatp(\chamhatN)=\chamhatp(\chamhatN+1)=0$, 
$r_{\chamhatN-2}=r_{\chamhatN-1}=\chamq^{-1}$}}

\multiput(105.8,3.8)(00,00){1}{\oval(6.4,6.4)[br]}
\multiput(122,3.8)(00,00){1}{\oval(26,13)[t]}
\multiput(135,01)(00,00){1}{\line(0,1){2.8}}

\multiput(00,00)(15,00){2}{\circle{2}}
\multiput(01,00)(15,00){2}{\line(1,0){13}}
\multiput(32.5,00)(2,00){6}{\circle*{0.6}}
\multiput(46,00)(15,00){2}{\line(1,0){13}}
\multiput(60,00)(00,00){2}{\circle{2}}
\multiput(77.5,00)(2,00){6}{\circle*{0.6}}
\multiput(91,00)(15,00){2}{\line(1,0){13}}
\multiput(105,00)(15,00){3}{\circle{2}}

\put(-3,3){\tiny{$\chamq_1$}}
\put(5,2){\tiny{$r_1$}}
\put(13,3){\tiny{$\chamq_2$}}
\put(20,2){\tiny{$r_2$}}
\put(47,2){\tiny{$r_{i-1}$}}
\put(57,3){\tiny{$\chamq_i$}}
\put(62,2){\tiny{$r_{i+1}$}}
\put(91,4){\tiny{$r_{\chamhatN-2}$}}
\put(96,9){\tiny{$\chamq_{\chamhatN-2}$}}
\put(109.5,2){\tiny{$\chamq^{-1}$}}
\put(119,3){\tiny{$\chamq$}}
\put(136,2){\tiny{$\chamq$}}

\put(119,12){\tiny{$\chamq^{-1}$}}

\put(-2,-5){\tiny{$\al_1$}}
\put(13,-5){\tiny{$\al_2$}}
\put(58,-5){\tiny{$\al_i$}}
\put(103,-5){\tiny{$\al_{\chamhatN-2}$}}
\put(118,-5){\tiny{$\al_{\chamhatN-1}$}}
\put(133,-5){\tiny{$\al_{\chamhatN}$}}

\end{picture}}

\newcommand{\smallHecHigherNineTenThree}{
\setlength{\unitlength}{0.5mm}
\begin{picture}(230,20)(00,00)

\put(-5,16){\tiny{$\chamhatp(\chamhatN-1)=1$, 
$\chamhatp(\chamhatN)=\chamhatp(\chamhatN+1)=0$, 
$r_{\chamhatN-2}=\chamq$, $r_{\chamhatN-1}=\chamq^{-1}$}}

\multiput(105.8,3.8)(00,00){1}{\oval(6.4,6.4)[br]}
\multiput(122,3.8)(00,00){1}{\oval(26,13)[t]}
\multiput(135,01)(00,00){1}{\line(0,1){2.8}}

\multiput(00,00)(15,00){2}{\circle{2}}
\multiput(01,00)(15,00){2}{\line(1,0){13}}
\multiput(32.5,00)(2,00){6}{\circle*{0.6}}
\multiput(46,00)(15,00){2}{\line(1,0){13}}
\multiput(60,00)(00,00){2}{\circle{2}}
\multiput(77.5,00)(2,00){6}{\circle*{0.6}}
\multiput(91,00)(15,00){3}{\line(1,0){13}}
\multiput(105,00)(15,00){3}{\circle{2}}

\put(-3,3){\tiny{$\chamq_1$}}
\put(5,2){\tiny{$r_1$}}
\put(13,3){\tiny{$\chamq_2$}}
\put(20,2){\tiny{$r_2$}}
\put(47,2){\tiny{$r_{i-1}$}}
\put(57,3){\tiny{$\chamq_i$}}
\put(62,2){\tiny{$r_{i+1}$}}
\put(91,4){\tiny{$r_{\chamhatN-2}$}}
\put(96,9){\tiny{$\chamq_{\chamhatN-2}$}}
\put(111,2){\tiny{$\chamq$}}
\put(117,3){\tiny{$-1$}}
\put(124,2){\tiny{$\chamq^{-2}$}}
\put(136,2){\tiny{$-1$}}

\put(119,12){\tiny{$\chamq$}}

\put(-2,-5){\tiny{$\al_1$}}
\put(13,-5){\tiny{$\al_2$}}
\put(58,-5){\tiny{$\al_i$}}
\put(103,-5){\tiny{$\al_{\chamhatN-2}$}}
\put(118,-5){\tiny{$\al_{\chamhatN-1}$}}
\put(133,-5){\tiny{$\al_{\chamhatN}$}}

\end{picture}}

\newcommand{\smallHecHigherNineTenFour}{
\setlength{\unitlength}{0.5mm}
\begin{picture}(230,20)(00,00)

\put(-5,15){\tiny{$\chamhatp(\chamhatN)=1$,
$r_{\chamhatN-1}=\chamq$}}

\multiput(00,00)(15,00){2}{\circle{2}}
\multiput(01,00)(15,00){2}{\line(1,0){13}}
\multiput(32.5,00)(2,00){6}{\circle*{0.6}}
\multiput(46,00)(15,00){2}{\line(1,0){13}}
\multiput(60,00)(00,00){2}{\circle{2}}
\multiput(77.5,00)(2,00){6}{\circle*{0.6}}
\multiput(91,00)(15,00){3}{\line(1,0){13}}
\multiput(105,00)(15,00){3}{\circle{2}}

\put(-3,3){\tiny{$\chamq_1$}}
\put(5,2){\tiny{$r_1$}}
\put(13,3){\tiny{$\chamq_2$}}
\put(20,2){\tiny{$r_2$}}
\put(47,2){\tiny{$r_{i-1}$}}
\put(57,3){\tiny{$\chamq_i$}}
\put(62,2){\tiny{$r_{i+1}$}}
\put(91,4){\tiny{$r_{\chamhatN-2}$}}
\put(102,9){\tiny{$\chamq_{\chamhatN-2}$}}
\put(106,4){\tiny{$r_{\chamhatN-1}$}}
\put(117,9){\tiny{$\chamq_{\chamhatN-1}$}}
\put(125,2){\tiny{$\chamq^2$}}
\put(132,3){\tiny{$\chamq^{-2}$}}

\put(-2,-5){\tiny{$\al_1$}}
\put(13,-5){\tiny{$\al_2$}}
\put(58,-5){\tiny{$\al_i$}}
\put(103,-5){\tiny{$\al_{\chamhatN-2}$}}
\put(118,-5){\tiny{$\al_{\chamhatN}$}}
\put(133,-5){\tiny{$\al_{\chamhatN-1}$}}

\end{picture}}

\newcommand{\cartanAone}{
\setlength{\unitlength}{0.5mm}
\begin{picture}(230,70)(-30,-10)
\put(0,30){\circle{2}}
\put(-10,20){\UPRIGHTlineTHICKsmall{8.4}}\put(10,20){\UPLEFTlineTHICKsmall{8.4}}
\put(-10,20){\circle{2}}\put(10,20){\circle{2}}
\put(-10,10){\UPlineTHICKsmall{8}}\put(10,10){\UPlineTHICKsmall{8}}
\put(-10,10){\circle{2}}\put(10,10){\circle{2}}
\put(0,0){\UPLEFTlineTHICKsmall{8.4}}\put(0,0){\UPRIGHTlineTHICKsmall{8.4}}
\put(0,0){\circle{2}}

\put(-2,33){{\tiny{$a$}}}
\put(-16,20){{\tiny{$a$}}}\put(12,20){{\tiny{$a$}}}
\put(-16,10){{\tiny{$a$}}}\put(12,10){{\tiny{$a$}}}
\put(-2,-5){{\tiny{$a$}}}

\put(-8,26){{\tiny{$1$}}}\put(6,26){{\tiny{$2$}}}
\put(-8,14){{\tiny{$2$}}}\put(6,14){{\tiny{$1$}}}
\put(-8,1){{\tiny{$1$}}}\put(6,1){{\tiny{$2$}}}

\put(45,28){$a$} 
\put(60,30){\circle{4}}\put(62,30){\line(1,0){26}}\put(90,30){\circle{4}}
\put(58,20){$\al_1$}\put(88,20){$\al_2$}
\put(58,35){${\dot{q}}$}\put(88,35){${\dot{q}}^{-1}$}
\put(70,33){${\dot{q}}^{-1}$}

\put(50,00){$({\dot{q}}\in\bKt\setminus\{1\})$}

\end{picture}
}

\newcommand{\supersltwoone}{
\setlength{\unitlength}{0.5mm}
\begin{picture}(230,70)(-30,-10)
\put(0,30){\circle{2}}
\put(-10,20){\UPRIGHTlineTHICKsmall{8.4}}\put(10,20){\UPLEFTlineTHICKsmall{8.4}}
\put(-10,20){\circle{2}}\put(10,20){\circle{2}}
\put(-10,10){\UPlineTHICKsmall{8}}\put(10,10){\UPlineTHICKsmall{8}}
\put(-10,10){\circle{2}}\put(10,10){\circle{2}}
\put(0,0){\UPLEFTlineTHICKsmall{8.4}}\put(0,0){\UPRIGHTlineTHICKsmall{8.4}}
\put(0,0){\circle{2}}

\put(-2,33){{\tiny{$a$}}}
\put(-16,20){{\tiny{$a$}}}\put(12,20){{\tiny{$b$}}}
\put(-16,10){{\tiny{$b$}}}\put(12,10){{\tiny{$c$}}}
\put(-2,-5){{\tiny{$c$}}}

\put(-8,26){{\tiny{$1$}}}\put(6,26){{\tiny{$2$}}}
\put(-8,14){{\tiny{$2$}}}\put(6,14){{\tiny{$1$}}}
\put(-8,1){{\tiny{$1$}}}\put(6,1){{\tiny{$2$}}}

\put(45,28){$a$} 
\put(60,30){\circle{4}}\put(62,30){\line(1,0){26}}\put(90,30){\circle{4}}
\put(58,20){$\al_1$}\put(88,20){$\al_2$}
\put(58,35){${\dot{q}}$}\put(86,35){$-1$}
\put(70,33){${\dot{q}}^{-1}$}

\put(125,28){$b$} 
\put(140,30){\circle{4}}\put(142,30){\line(1,0){26}}\put(170,30){\circle{4}}
\put(138,20){$\al_1$}\put(168,20){$\al_2$}
\put(136,35){$-1$}\put(166,35){$-1$}
\put(153,33){${\dot{q}}$}

\put(45,-2){$c$} 
\put(60,00){\circle{4}}\put(62,00){\line(1,0){26}}\put(90,00){\circle{4}}
\put(58,-10){$\al_1$}\put(88,-10){$\al_2$}
\put(56,5){$-1$}\put(88,5){${\dot{q}}$}
\put(70,3){${\dot{q}}^{-1}$}

\put(110,00){$({\dot{q}}\in\bKt\setminus\{-1,1\})$}

\end{picture}
}

\newcommand{\smallHecHigherEighteenZero}{
\setlength{\unitlength}{0.5mm}
\begin{picture}(230,20)(00,00)

\put(-5,18){\tiny{
$q_{\chamhatN-3}=q_{\chamhatN-5}^{-1}\in\{-1,\zeta^{-1}\}$, 
$r_{\chamhatN-4}\chamq_{\chamhatN-3}^2\zeta=1$}}

\multiput(105.8,3.8)(00,00){1}{\oval(6.4,6.4)[br]}
\multiput(129.8,3.8)(00,00){1}{\oval(41.6,13)[t]}
\multiput(150.6,00.8)(00,00){1}{\line(0,1){3.0}}

\multiput(120.8,3.8)(00,00){1}{\oval(6.4,6.4)[br]}
\multiput(136.7,3.8)(00,00){1}{\oval(25.4,11)[t]}
\multiput(149.4,00.8)(00,00){1}{\line(0,1){3.0}}

\multiput(00,00)(15,00){2}{\circle{2}}
\multiput(01,00)(15,00){2}{\line(1,0){13}}
\multiput(32.5,00)(2,00){6}{\circle*{0.6}}
\multiput(46,00)(15,00){2}{\line(1,0){13}}
\multiput(60,00)(00,00){2}{\circle{2}}
\multiput(77.5,00)(2,00){6}{\circle*{0.6}}
\multiput(91,00)(15,00){3}{\line(1,0){13}}
\multiput(105,00)(15,00){4}{\circle{2}}

\put(-3,3){\tiny{$\chamq_1$}}
\put(5,2){\tiny{$r_1$}}
\put(13,3){\tiny{$\chamq_2$}}
\put(20,2){\tiny{$r_2$}}
\put(47,2){\tiny{$r_{i-1}$}}
\put(57,3){\tiny{$\chamq_i$}}
\put(62,2){\tiny{$r_{i+1}$}}
\put(91,4){\tiny{$r_{\chamhatN-4}$}}
\put(96,9){\tiny{$\chamq_{\chamhatN-3}$}}
\put(111,2){\tiny{$\zeta$}}
\put(117,3){\tiny{$-1$}}
\put(125,2){\tiny{$\zeta^{-1}$}}
\put(135,3){\tiny{$\zeta$}}
\put(152,0){\tiny{$-1$}}

\put(142,12.5){\tiny{$\zeta$}}
\put(142,5){\tiny{$\zeta$}}

\put(-2,-5){\tiny{$\al_1$}}
\put(13,-5){\tiny{$\al_2$}}
\put(58,-5){\tiny{$\al_i$}}
\put(103,-5){\tiny{$\al_{\chamhatN-3}$}}
\put(118,-5){\tiny{$\al_{\chamhatN-2}$}}
\put(133,-5){\tiny{$\al_{\chamhatN}$}}
\put(148,-5){\tiny{$\al_{\chamhatN-1}$}}

\end{picture}}

\newcommand{\FigOneA}{{\rm{Fig.1}}}
\newcommand{\FigTwoA}{{\rm{Fig.2}}}
\newcommand{\FigThreeA}{{\rm{Fig.3}}}
\newcommand{\FigFourA}{{\rm{Fig.4}}}
\newcommand{\FigFiveA}{{\rm{Fig.5}}}
\newcommand{\FigSixA}{{\rm{Fig.6}}}
\newcommand{\FigSevenA}{{\rm{Fig.7}}}
\newcommand{\FigEightA}{{\rm{Fig.8}}} \newcommand{\NFigEightA}{8}
\newcommand{\FigNineA}{{\rm{Fig.9}}}
\newcommand{\FigTenA}{{\rm{Fig.10}}}

\newcommand{\FigElevenA}{{\rm{Fig.11}}}
\newcommand{\FigTwelveA}{{\rm{Fig.12}}}
\newcommand{\FigThirteenA}{{\rm{Fig.13}}}
\newcommand{\FigFourteenA}{{\rm{Fig.14}}}
\newcommand{\FigFifteenA}{{\rm{Fig.15}}}
\newcommand{\FigSixteenA}{{\rm{Fig.16}}}\newcommand{\NFigSixteenA}{16}
\newcommand{\FigSeventeenA}{{\rm{Fig.17}}}\newcommand{\NFigSeventeenA}{17}
\newcommand{\FigEighteenA}{{\rm{Fig.18}}}
\newcommand{\FigNineteenA}{{\rm{Fig.19}}}
\newcommand{\FigTwentyA}{{\rm{Fig.20}}}

\newcommand{\FigTwentyOneA}{{\rm{Fig.21}}} \newcommand{\NFigTwentyOneA}{21}
\newcommand{\FigTwentyTwoA}{{\rm{Fig.22}}}
\newcommand{\FigTwentyThreeA}{{\rm{Fig.23}}}
\newcommand{\FigTwentyFourA}{{\rm{Fig.24}}}
\newcommand{\FigTwentyFiveA}{{\rm{Fig.25}}} \newcommand{\NFigTwentyFiveA}{25}
\newcommand{\FigTwentySixA}{{\rm{Fig.26}}}
\newcommand{\FigTwentySevenA}{{\rm{Fig.27}}} \newcommand{\NFigTwentySevenA}{{\rm{27}}}
\newcommand{\FigTwentyEightA}{{\rm{Fig.28}}}
\newcommand{\FigTwentyNineA}{{\rm{Fig.29}}}
\newcommand{\FigThirtyA}{{\rm{Fig.30}}}

\newcommand{\FigThirtyOneA}{{\rm{Fig.31}}}
\newcommand{\FigThirtyTwoA}{{\rm{Fig.32}}} \newcommand{\NFigThirtyTwoA}{32}
\newcommand{\FigThirtyThreeA}{{\rm{Fig.33}}}
\newcommand{\FigThirtyFourA}{{\rm{Fig.34}}}
\newcommand{\FigThirtyFiveA}{{\rm{Fig.35}}} \newcommand{\NFigThirtyFiveA}{35}
\newcommand{\FigThirtySixA}{{\rm{Fig.36}}}
\newcommand{\FigThirtySevenA}{{\rm{Fig.37}}} \newcommand{\NFigThirtySevenA}{37}
\newcommand{\FigThirtyEightA}{{\rm{Fig.38}}}
\newcommand{\FigThirtyNineA}{{\rm{Fig.39}}}
\newcommand{\FigFourtyA}{{\rm{Fig.40}}}

\newcommand{\FigFourtyOneA}{{\rm{Fig.41}}}
\newcommand{\FigFourtyTwoA}{{\rm{Fig.42}}}
\newcommand{\FigFourtyThreeA}{{\rm{Fig.43}}} \newcommand{\NFigFourtyThreeA}{43}
\newcommand{\FigFourtyFourA}{{\rm{Fig.44}}} \newcommand{\NFigFourtyFourA}{44}
\newcommand{\FigFourtyFiveA}{{\rm{Fig.45}}}
\newcommand{\FigFourtySixA}{{\rm{Fig.46}}} \newcommand{\NFigFourtySixA}{46}
\newcommand{\FigFourtySevenA}{{\rm{Fig.47}}} \newcommand{\NFigFourtySevenA}{47}
\newcommand{\FigFourtyEightA}{{\rm{Fig.48}}}
\newcommand{\FigFourtyNineA}{{\rm{Fig.49}}} \newcommand{\NFigFourtyNineA}{49}
\newcommand{\FigFiftyA}{{\rm{Fig.50}}} \newcommand{\NFigFiftyA}{50}

\newcommand{\FigFiftyOneA}{{\rm{Fig.51}}} \newcommand{\NFigFiftyOneA}{51}
\newcommand{\FigFiftyTwoA}{{\rm{Fig.52}}} \newcommand{\NFigFiftyTwoA}{52}
\newcommand{\FigFiftyThreeA}{{\rm{Fig.53}}} \newcommand{\NFigFiftyThreeA}{53}
\newcommand{\FigFiftyFourA}{{\rm{Fig.54}}} \newcommand{\NFigFiftyFourA}{54}
\newcommand{\FigFiftyFiveA}{{\rm{Fig.55}}} \newcommand{\NFigFiftyFiveA}{55}
\newcommand{\FigFiftySixA}{{\rm{Fig.56}}} \newcommand{\NFigFiftySixA}{56}
\newcommand{\FigFiftySevenA}{{\rm{Fig.57}}} \newcommand{\NFigFiftySevenA}{57}
\newcommand{\FigFiftyEightA}{{\rm{Fig.58}}}
\newcommand{\FigFiftyNineA}{{\rm{Fig.59}}} \newcommand{\NFigFiftyNineA}{59}
\newcommand{\FigSixtyA}{{\rm{Fig.60}}} \newcommand{\NFigSixtyA}{60}

\newcommand{\FigSixtyOneA}{{\rm{Fig.61}}} \newcommand{\NFigSixtyOneA}{61}
\newcommand{\FigSixtyTwoA}{{\rm{Fig.62}}} \newcommand{\NFigSixtyTwoA}{62}
\newcommand{\FigSixtyThreeA}{{\rm{Fig.63}}} \newcommand{\NFigSixtyThreeA}{63}
\newcommand{\FigSixtyFourA}{{\rm{Fig.64}}} \newcommand{\NFigSixtyFourA}{64}
\newcommand{\FigSixtyFiveA}{{\rm{Fig.65}}} 
\newcommand{\FigSixtySixA}{{\rm{Fig.66}}} \newcommand{\NFigSixtyFiveA}{66}
\newcommand{\FigSixtySevenA}{{\rm{Fig.67}}}
\newcommand{\FigSixtyEightA}{{\rm{Fig.68}}}
\newcommand{\FigSixtyNineA}{{\rm{Fig.69}}}
\newcommand{\FigSeventyA}{{\rm{Fig.70}}}

\newcommand{\FigSeventyOneA}{{\rm{Fig.71}}}
\newcommand{\FigSeventyTwoA}{{\rm{Fig.72}}}
\newcommand{\FigSeventyThreeA}{{\rm{Fig.73}}}
\newcommand{\FigSeventyFourA}{{\rm{Fig.74}}} \newcommand{\NFigSeventyFourA}{74}
\newcommand{\FigSeventyFiveA}{{\rm{Fig.75}}}
\newcommand{\FigSeventySixA}{{\rm{Fig.76}}}
\newcommand{\FigSeventySevenA}{{\rm{Fig.77}}}
\newcommand{\FigSeventyEightA}{{\rm{Fig.78}}}
\newcommand{\FigSeventyNineA}{{\rm{Fig.79}}}
\newcommand{\FigEightyA}{{\rm{Fig.80}}}

\newcommand{\FigEightyOneA}{{\rm{Fig.81}}}
\newcommand{\FigEightyTwoA}{{\rm{Fig.82}}}
\newcommand{\FigEightyThreeA}{{\rm{Fig.83}}}
\newcommand{\FigEightyFourA}{{\rm{Fig.84}}}
\newcommand{\FigEightyFiveA}{{\rm{Fig.85}}}
\newcommand{\FigEightySixA}{{\rm{Fig.86}}}
\newcommand{\FigEightySevenA}{{\rm{Fig.87}}}
\newcommand{\FigEightyEightA}{{\rm{Fig.88}}}
\newcommand{\FigEightyNineA}{{\rm{Fig.89}}}
\newcommand{\FigNinetyA}{{\rm{Fig.90}}}

\newcommand{\FigNinetyOneA}{{\rm{Fig.91}}}
\newcommand{\FigNinetyTwoA}{{\rm{Fig.92}}}
\newcommand{\FigNinetyThreeA}{{\rm{Fig.93}}}
\newcommand{\FigNinetyFourA}{{\rm{Fig.94}}}
\newcommand{\FigNinetyFiveA}{{\rm{Fig.95}}}
\newcommand{\FigNinetySixA}{{\rm{Fig.96}}}
\newcommand{\FigNinetySevenA}{{\rm{Fig.97}}}
\newcommand{\FigNinetyEightA}{{\rm{Fig.98}}}
\newcommand{\FigNinetyNineA}{{\rm{Fig.99}}}
\newcommand{\FigNeHundredA}{{\rm{Fig.100}}}

\newcommand{\tableone}{{\rm{1}}}
\newcommand{\tabletwo}{{\rm{2}}}
\newcommand{\tablethree}{{\rm{3}}}
\newcommand{\tablefour}{{\rm{4}}}

\begin{document}
\begingroup
\renewcommand{\arraystretch}{1.1}
{\Huge{
\begin{tabular}{c}
Hamilton Circuits  of  Cayley Graphs \\ of Weyl Groupoids \\
of Generalized Quantum Groups
\end{tabular}}}
\endgroup
\newline
\begin{center}
{\large{Hiroyuki Yamane}}
\end{center}

\vspace{1cm}

Abstract. We 
study Hamilton circuits of the Cayley graphs of the Weyl groupoids of 
the generalized quantum groups, or the quantum double of the Nichols algebras of diagonal-type,
with finite root systems. 
We prove the existence of a Hamilton circuit for any of them
and explicitly  draw one of them for rank 3 and 4 cases.

\section{Introduction}
This paper studies a Hamilton circuit of the Cayley graphs of the Weyl groupoids
associated
with the generalized root systems, or the quantum doubles of the finite-type
Nichols algebras of diagonal type, over an field of characteristic zero. The
importance of the Weyl groupoids for study of Hopf algebras and Lie 
superalgebras is consistently written in a recent book \cite{HS20}. In particular, they
play important role in study of classification of Nichols algebras. Moreover
since the quantum double of the Nichols algebras of diagonal type is a
generalization of the quantum groups, we call them
{\it{the generalized quantum groups}}.

Let ${\mathbb{K}}$ be a field, and ${\mathbb{K}}^\times:={\mathbb{K}}\setminus\{0\}$.
Let $\bZ\Pi$ be a finite-rank free $\bZ$-module with a basis $\Pi$.
Let $\chi:\bZ\Pi\times\bZ\Pi\to{\mathbb{K}}^\times$ be the map satisfying
\begin{equation*}
\chi(\lambda+\mu,\nu)=\chi(\lambda,\nu)\chi(\mu,\nu),\,\,
\chi(\lambda,\mu+\nu)=\chi(\lambda,\mu)\chi(\lambda,\nu)\quad(\lambda,\mu,\nu\in\bZ\Pi).
\end{equation*}
The Nichols algebras $U^+(\chi)$ of diagonal type has a Kharchenko's PBW-basis described by
the generaliezed root system $R(\chi)=R^+(\chi)\cup(-R^+(\chi))$. 
I.~Heckenberger \cite{Hec09} classified $\chi$ with $|R(\chi)|<\infty$ for ${\mathbb{K}}$ of characteristic zero.
I.~Heckenberger and the author~\cite{HY08} introduced an axiomatic definition of the generalized root systems,
which include  $R(\chi)$'s.
M.~Cuntz and I.~Heckenbeger improved the axiomatic definition by \cite{CH09} 
and completely classified the finite Weyl groupoids by \cite{CH15}.
The result of \cite{CH15} tells that there are so many finite Weyl groupoids other than 
those associated with the Nichols algebras of diagonal type.
In this paper, we concentrate on studying a Hamilton circuit of the Cayley graph
of the Weyl groupoid of the Nichols algebras $U^+(\chi)$ of diagonal type with $|R(\chi)|<\infty$
over a field ${\mathbb{K}}$ of characteristic zero.
 
On the other hand, J.H.~Conway, N.J.A.~Sloane and Allan~R.~Wilks \cite{CSW89} proved an existence of a Hamilton circuit for the Cayley graph
of every
finite Coxeter group in a very concise way. Therefore it is natural to ask whether the Cayley graph of a finite Weyl groupoid has a Hamilton circuit. 
In general, to determine whether a Hamiton path or a Hamilton 
circuit 
exists for a given graph is one of the   
major theme in graph theory, and called the 
Hamiltonian path problem.
For basic terminologies, see \cite{diestel}. 
 One important open problem is whether  every  Cayley graph    
  contains a Hamilton path or a Hamilton circuit, see \cite{lanel-pallage-ratnayake-thevasha-welihinda}.

We also mention that the author with his collaborators achieved
results \cite{AYY15}, \cite{AY15}, \cite{AY18}, \cite{BY18}, 
\cite{BY20}, \cite{HY10}, \cite{JMY18}, \cite{Y09}, \cite{Y16}, \cite{Y21} concerning representation theory 
the generalized quantum groups
$U(\chi)$ defined as the quantum double of $U^+(\chi)$.

Since there  are so many finite Weyl groupoids 
\cite{CH15} (see also \cite[($|A|$ of) Table~1 and Table~2]{CH12}) other than
those of the generalized quantum groups, it would be interesting if in future
we are able
to obtain Hamilton circuits for them using computers (if any).

\newcommand{\NVMpi}{\pi}

\begin{notation}
{\rm{(1)}}
For $n\in\bN$, we 
let $\hatpi:\bZ\to\bZ/n\bZ$ mean the canonical group epimorphism,
i.e.,
$\hatpi_n(t):=t+n\bZ$ $(t\in\bZ)$.
\newline
{\rm{(2)}} For $x$, $y\in\bR\cup\{\infty,-\infty\}$, let
$\fkJ_{x,y}:=\{k\in\bZ|x\leq k\leq y\}$.
\end{notation}

\section{Hamilton circuits of Cayley graphs of finite Coxeter Groups}

In this section we recall the main ideas of the proof given in \cite{CSW89} of the
existence of a Hamilton circuit for the Cayley graph of the finite Coxeter
group.
An adapted argument will be used later for finite Coxeter groupoids.

Recall that a Coxeter system consists of the following data:
\newline\newline
(a) A non-empty finite set $S$.
\newline
(b) A map $\Psi=\Psi^{(W,S)}:S\times S\to\bN$ such that
\begin{equation*}
\Psi(s,s)=1,\,\,\Psi(s,t)=\Psi(t,s)\geq 2,\quad s,t\in S,\,s\ne t.
\end{equation*}
(c) A group $W$ defined by the generators $s\in S$ and relations
\begin{equation*}
(st)^{\Psi(s,t)}=e,\quad s,t\in S,
\end{equation*}
\newline  
Let $(W, S)$ be a finite Coxeter system, i.e., 
the following $({\rm{Cx}}1)$-$({\rm{Cx}}4)$ are fulfilled.
\newline\newline
$({\rm{Cx}}1)$ $S$ is a non-empty finite set. \newline
$({\rm{Cx}}2)$ There exists a map $\Psi^{(W, S)}:S\times S\to\bN$ 
such that $\Psi^{(W, S)}(s,s)=1$ $(s\in S)$,
$\Psi^{(W, S)}(s^\prime,s^{\prime\prime})=\Psi^{(W, S)}(s^{\prime\prime},s^\prime)
\geq 2$ $(s^\prime,s^{\prime\prime}\in S,\,
s^\prime\ne s^{\prime\prime})$. \newline
$({\rm{Cx}}3)$ 
$W$ is the group defined by the generators $s\in S$
and 
\begin{equation*}
\mbox{the relations $(s^\prime s^{\prime\prime})^{\Psi^{(W, S)}(s^\prime,s^{\prime\prime})}
=e$ for all $s^\prime$, $s^{\prime\prime}\in S$.} 
\end{equation*}
where $e$ means the unit element of $W$.
\newline
$({\rm{Cx}}4)$ $|W|<\infty$.
\newline\newline
Let $X$ be an $|S|$-dimensional $\bR$-linear space with a basis $\{x_s|s\in S\}$.
Define the symmetric bilinear map
$(\,,\,): X\times X\to\bR$ by 
\begin{equation*}
(x_s,x_{s^\prime}):=
-2\cos{\frac {\pi} {\Psi^{(W, S)}(s^\prime,s^{\prime\prime})}}\quad
(s^\prime,s^{\prime\prime}\in S).
\end{equation*}
Then we have a group monomorphism 
$f:W\to{\mathrm{GL}}(X)$ defined by
$f(s)(y):=y-(y,x_s)x_s$
$(s\in S,y\in X)$. Hence
\begin{equation*}
\forall s, \forall s^\prime\in S,
\Psi^{(W, S)}(s,s^\prime)
=\min\{m\in\bN|(ss^\prime)^m=e\}.
\end{equation*}

In this section, we explain the idea of \cite{CSW89} to obtain a Hamilton circuit ${\mathcal{H}}_{W, S}$
of the Cayley graph ${\mathcal{C}}_{W, S}$ of $(W, S)$.
Here ${\mathcal{C}}_{W, S}$ is the graph composed of 
the vertices bijectively corresponding to
the elements of $W$ and the edges labeled
by $s\in S$
connecting two vertices corresponding to 
$x$, $y\in W$ with $y=sx$.
Let $k:=|W|$.
If there exists a bijection $\theta:\bZ/k\bZ\to W$
such that $\theta(\hatpi_k(t+1))\theta(\hatpi_k(t))^{-1}\in S$ for all $t\in\bZ$,
then {\it{the Hamilton circuit ${\mathcal{H}}_{W, S}$ of ${\mathcal{C}}_{W, S}$ defined by $\theta$}}
is the subgraph of ${\mathcal{C}}_{W, S}$
composed of all vertices of ${\mathcal{C}}_{W, S}$ 
and the $k$-edges each of which connects two vertices
labeled by  $\theta(\hatpi_k(t))$ and $\theta(\hatpi_k(t+1))$ for some $t\in\bZ$.

We clearly have the following lemma.

\begin{lemma} \label{lemma:ranktwoCoxeterGroup}
Assume $|S|=2$. Let $\{s_1,s_2\}:=S$.
Then a Hamilton circuit ${\mathcal{H}}_{W, S}$ uniquely exists.
More precisely we have the following.
Let $m:=\Psi^{(W, S)}(s_1,s_2)$.
Then $|W|=2m$. 
Define the bijection $\theta:\bZ/2m\bZ\to W$ by
$\theta(\hatpi_k(0)):=e$, and
$\theta(\hatpi_k(2t+1)):=s_1\theta(\hatpi_k(2t))$,
$\theta(\hatpi_k(2t+2)):=s_2\theta(\hatpi_k(2t+1))$
$(t\in\bZ)$.
Then ${\mathcal{H}}_{W, S}$ defined by $\theta$ coincides with 
${\mathcal{C}}_{W, S}$.
\end{lemma}

We shall mostly use the argument of the proof  of the following theorem.

\begin{theorem}\label{theorem:Main:CSW} {\rm{(\cite[Theorem]{CSW89})}}
A Hamilton circuit ${\mathcal{H}}_{W, S}$ exists
for any finite Coxeter system $(W,S)$.
\end{theorem}
\noindent
{\it{Proof.}}
Let $n:=|S|$.
If $n=1$, the claim is obvious.
If $n=2$, the claim follows from Lemma~\ref{lemma:ranktwoCoxeterGroup}.

Assume $n\geq 3$. 
Then we can easily see from the classification of the irreducible
finite Coxeter systems 
that there exists a bijection 
$\hats:\fkJ_{1,n}\to S$ fulfilling that
\begin{equation} \label{eqn:sisn}
\forall i\in\fkJ_{1,n-2},\,\, \hats(i)\hats(n)=\hats(n)\hats(i).
\end{equation}
Let $\hatS:=\{\hats(i)|i\in\fkJ_{1,n-1}\}$.
Let $\hatW$ be the subgroup of $W$ generated by
$\hatS$. Note that $(\hatW,\hatS)$ is a Coxeter system.
Let $h:=|W|$, $\hath:=|\hatW|$, and
$p:={\frac {h} {\hath}}\,(\in\fkJ_{2,\infty})$. 
Then there exists $w_u\in W$ $(u\in\fkJ_{1,p})$
such that $w_1=e$, 
$W=\cup_{x=1}^p\hatW w_x$
and $\hatW w_x\cap \hatW w_{x^\prime}=\emptyset$
$(x\ne x^\prime)$ , i.e.,
$W=\amalg_{x=1}^p\hatW w_x$ is the right coset decomposition
of $W$ by $\hatW$. 

By induction, we have
the Hamilton circuit ${\mathcal{H}}_{\hatW, \hatS}$ 
of ${\mathcal{C}}_{\hatW, \hatS}$ defined by 
a bijection $\hattheta:\bZ/\hath\bZ\to\hatW$.
Since
the elements of a non-empty proper subset of $\hatS$
can not generate $\hatW$,
we see that 
\begin{equation} \label{eqn:existsmk}
\forall k\in\fkJ_{1,n-1},\,\exists m_k\in\bZ,\,
\hattheta(\hatpi_\hath(m_k+1))=\hats(k)\hattheta(\hatpi_\hath(m_k)).
\end{equation}
For $u\in\hatW$ and $k\in\fkJ_{1,n-1}$, fixing $m_k$ of \eqref{eqn:existsmk},  
define the bijection
$\hattheta_{u,k}:\bZ/\hath\bZ\to\hatW$ by
\begin{equation*}
\hattheta_{u,k}(\hatpi_\hath(t)):=
\hattheta(\hatpi_\hath(m_k+t))\hattheta(\hatpi_\hath(m_k))^{-1}u
\quad (t\in\bZ),
\end{equation*} and define the surjection
${\hat{\varphi}}_{u,k}:\bZ/\hath\bZ\to\fkJ_{1,n-1}$
by 
\begin{equation*}
(\hats\circ{\hat{\varphi}}_{u,k})(\hatpi_\hath(t)):
=\hattheta(\hatpi_\hath(m_k+t))\hattheta(\hatpi_\hath(m_k+t-1))^{-1}
\quad (t\in\bZ).
\end{equation*}
Then we have
the Hamilton circuit ${\mathcal{H}}^{u,k}_{\hatW, \hatS}$ 
of ${\mathcal{C}}_{\hatW, \hatS}$ defined by $\hattheta_{u,k}$.
Note that 
\begin{equation}\label{eqn:hvphiukoi}
\mbox{$\hattheta_{u,k}(\hatpi_\hath(0))=u$
and $\hattheta_{u,k}(\hatpi_\hath(1))=\hats(k)u=(\hats\circ{\hat{\varphi}}_{u,k})(\hatpi_\hath(1))u$.}
\end{equation}

{\it{We shall construct a Hamilton circuit in an inductive way.}}
Let $l\in\fkJ_{1,p}$.
We shall show by induction on $l$ that there exist
an injection $g_l:\fkJ_{1,l}\to\fkJ_{1,p}$ and 
a bijection
$\dottheta_l:\bZ/l\hath\bZ\to\dotW_l$, 
where $\dotW_l:=\amalg_{z=1}^l\hatW w_{g_l(z)}$,
such that 
$\dottheta_l(\hatpi_{l\hath}(t+1))\dottheta_l(\hatpi_{l\hath}(t))^{-1}\in S$
for all $t\in\bZ$. Recall $p\hath=|W|$.
If $l=1$, we may let $g_l(1):=1$ and $\dottheta_1:=\hattheta$.
Let $l\in\fkJ_{2,p}$, and assume that the hypothesis of the induction for $l-1$
in place of $l$ holds. 	
Define the surjection $\varphi_{l-1}:\bZ/(l-1)\hath\bZ\to\fkJ_{1,n-1}$
by 
\begin{equation*}
\hats(\varphi_{l-1}(\hatpi_{(l-1)\hath}(t)))\dottheta_{l-1}(\hatpi_{(l-1)\hath}(t-1))
=\dottheta_{l-1}(\hatpi_{(l-1)\hath}(t))\quad (t\in\bZ).
\end{equation*}
Since $W$ is a finite group generated by $S$,
there exist $b_l\in\fkJ_{1,(l-1)\hath}$ and $c\in\fkJ_{1,p}\setminus g_{l-1}(\fkJ_{1,l-1})$
such that
$\hats(n)\dottheta_{l-1}(\hatpi_{(l-1)\hath}(b_l-1))\in\hatW w_c$.
If $\varphi_{l-1}(\hatpi_{(l-1)\hath}(b_l))=n-1$, then $\varphi_{l-1}(\hatpi_{(l-1)\hath}(b_l-1))\in\fkJ_{1,n-2}$,
which implies
$\hats(n)\dottheta_{l-1}(\hatpi_{(l-1)\hath}(b_l-2))\in\hatW w_c$
by \eqref{eqn:sisn}.
Hence we may assume $\varphi_{l-1}(\hatpi_{(l-1)\hath}(b_l))\in\fkJ_{1,n-2}$. 
Let $r:=\varphi_{l-1}(\hatpi_{(l-1)\hath}(b_l))\,(\in\fkJ_{1,n-2})$. 
Let $v:=\hats(n)\dottheta_{l-1}(\hatpi_{(l-1)\hath}(b_l-1))w_c^{-1}
\in\hatW$. 
Define the map $g_l:\fkJ_{1,l}\to\fkJ_{1,p}$ by by $(g_l)_{|\fkJ_{1,l-1}}:=g_{l-1}$
and  $g_l(l):=c$.
Define the map $\dottheta_l:\bZ/l\bZ\to \dotW_l(:=\amalg_{z=1}^l\hatW w_{g_l(z)})$ by
\begin{equation*}
\dottheta_l(\hatpi_{l\hath}(t)):=
\left\{\begin{array}{ll}
\dottheta_{l-1}(\hatpi_{(l-1)\hath}(t)) & 
\quad\mbox{if $t\in\fkJ_{1,b_l-1}$,} \\
vw_c\,(=\hats(n)\dottheta_{l-1}(\hatpi_{(l-1)\hath}(b_l-1))) & \quad\mbox{if $t=b_l$,} \\
\hattheta_{\hats(r)v,r}(\hatpi_\hath(t+1-b_l))w_c &
\quad\mbox{if $t\in\fkJ_{b_l+1, b_l+\hath-1}$,} \\
\dottheta_{l-1}(\hatpi_{(l-1)\hath}(t-\hath)) & 
\quad\mbox{if $t\in\fkJ_{b_l+\hath,l\hath}$,}
\end{array}\right.
\end{equation*} for $t\in\fkJ_{1,l\hath}$.
Let $u:=\hats(r)v\,(\in\hatW w_c)$. 
For $t\in\fkJ_{b_l+1, b_l+\hath-1}$, 
we have
\begin{equation*}
\begingroup
\renewcommand{\arraystretch}{1.2}
\begin{array}{lcl}
\dottheta_l(\hatpi_{l\hath}(t)) &  
= & \hattheta_{u,r}(\hatpi_\hath(t+1-b_l))w_c \\
& = & \underbrace{\hats({\hat{\varphi}}_{u,r}(\hatpi_\hath(t+1-b_l)))\cdots
\hats({\hat{\varphi}}_{u,r}(\hatpi_\hath(2)))
\hats({\hat{\varphi}}_{u,r}(\hatpi_\hath(1)))}_{t+1-b_l}\cdot  
uw_c
\\
& = & 
\underbrace{\hats({\hat{\varphi}}_{u,r}(\hatpi_\hath(t+1-b_l)))\cdots
\hats({\hat{\varphi}}_{u,r}(\hatpi_\hath(1)))}_{t-b_l+1}\cdot  
\hats(r)vw_c \\
& = & 
\underbrace{\hats({\hat{\varphi}}_{u,r}(\hatpi_\hath(t+1-b_l)))\cdots
\hats({\hat{\varphi}}_{u,r}(\hatpi_\hath(3)))
\hats({\hat{\varphi}}_{u,r}(\hatpi_\hath(2)))}_{t-b_l}\cdot vw_c.
\end{array}
\endgroup
\end{equation*}
We also have
\begin{equation*}
\begin{array}{l}
\hats(n)\dottheta_l(\hatpi_{l\hath}(b_l+\hath-1)) \\
\quad = \hats(n){\hat{\theta}}_{\hats(r)v,r}(\hatpi_\hath(\hath))w_c \\
\quad = \hats(n){\hat{\theta}}_{\hats(r)v,r}(\hatpi_\hath(0))w_c \\
\quad = \hats(n)\hats(r)vw_c \quad\mbox{(by \eqref{eqn:hvphiukoi})} \\
\quad = \hats(n)\hats(r)\hats(n)\dottheta_{l-1}(\hatpi_{(l-1)\hath}(b_l-1)) \\
\quad = \hats(r)\dottheta_{l-1}(\hatpi_{(l-1)\hath}(b_l-1)) \\
\quad = \dottheta_{l-1}(\hatpi_{(l-1)\hath}(b_l)) \\
\quad = \dottheta_l(\hatpi_{l\hath}(b_l+\hath)).
\end{array}
\end{equation*}
Then $g_l$ and $\dottheta_l$ are the desired ones,
where see also $\FigOneA$. 
\hfill $\Box$

\begin{equation*}
\begin{array}{c}       
\mbox{$\JointForCxeter$}
\\ 
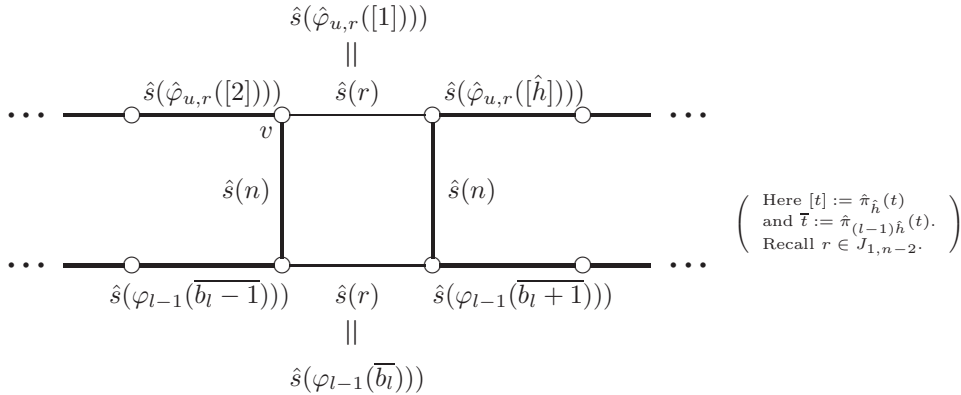
\mbox{$\FigOneA$: Joint to get a Hamilton circuit in Proof of Theorem~\ref{theorem:Main:CSW}}
\end{array}
\end{equation*}

\begin{example}\label{example:HCWeylBfour}
 {\rm{Here let $(W,S)$ be a Coxeter system of type $B_n$ $(n\in\fkJ_{2,\infty})$.
It means that $W$ is the group 
defined by generators $s_i$ $(i\in\fkJ_{1,n})$ and relations 
$s_i^2=e$, $s_is_j=s_js_i$ $(|i-j|\geq 2)$, 
$s_is_{i+1}s_i=s_{i+1}s_is_{i+1}$ $(i\in\fkJ_{2,n-1})$,
and $s_1s_2s_1s_2=s_2s_1s_2s_1$,
where we let $S=\{s_i|i\in\fkJ_{1,n}\}$.
It is well-known that $|W|=2^n\cdot n!$.
If $n=2$, a Hamilton circuit ${\mathcal{H}}_{W, S}$ is
$s_1 s_2 s_1 s_2 s_1 s_2 s_1 s_2$ by Lemma~\ref{lemma:ranktwoCoxeterGroup}. 
If $n=3$, a Hamilton circuit ${\mathcal{H}}_{W, S}$ is:
\begin{equation*} 
\begin{array}{c}
\begin{array}{l}
s_1 s_2 s_1 s_2 s_1 s_2( s_3 s_2 s_1 s_2 s_1 s_2( s_3 s_2 s_1 s_2( s_3 s_2 s_1 s_2 (s_3 s_2 (s_3 
s_2 s_1 s_2 s_1 s_2 s_1 s_2 s_3) s_2 s_1 s_2 s_1 s_2 \\ 
\cdot s_3) s_2 s_1 s_2 s_3) s_2 s_1 s_2 s_3) s_2 s_3) s_2,
\end{array} 
\\ \quad \\
\mbox{$\FigTwoA$: Hamilton circuit of the Weyl group of type $B_3$. See also 
$\FigThreeA$
below.}
\end{array}
\end{equation*} where we use the argument of Proof of Theorem~\ref{theorem:Main:CSW}
for $m_{k_l}=0$, $\hats(k_l)=s_1$ $(l\in\fkJ_{2,6})$ and
$b_2=7$, $b_3=13$, $b_4=17$, $b_5=21$, $b_6=23$.
If $n=4$, see $\FigFourA$ 
below.
\begin{equation*}
\begin{array}{l}
\hspace{1cm}    $\HamiltonWeylBthree$
\\ \quad \\
\mbox{$\FigThreeA$: The ($2$-convenient (see Definition~\ref{definition:specialgraph})) Hamilton circuit 
$\FigTwoA$} \\ 
\mbox{of the Weyl group of type $B_3$,
where $*$ is the initial and end point of $\FigTwoA$} 
\end{array}
\end{equation*}

\begin{equation*}
\begin{array}{l}\hspace{0.5cm}
\mbox{\footnotesize{$
\begingroup
\renewcommand{\arraystretch}{1.5}
\begin{array}{l}
s^{{\bar{0}}}_2  s^{{\bar{0}}}_3  s^{{\bar{0}}}_2  s^{{\bar{0}}}_3  s^{{\bar{0}}}_2  s^{{\bar{0}}}_1  s^{{\bar{0}}}_2  s^{{\bar{0}}}_3  s^{{\bar{0}}}_2  
s^{{\bar{0}}}_1  s^{{\bar{0}}}_2  s^{{\bar{0}}}_3  s^{{\bar{0}}}_2  s^{{\bar{0}}}_1  s^{{\bar{0}}}_2  s^{{\bar{0}}}_1  s^{{\bar{0}}}_2  s^{{\bar{0}}}_3  s^{{\bar{0}}}_2  
s^{{\bar{0}}}_1  s^{{\bar{0}}}_2  s^{{\bar{0}}}_1  s^{{\bar{0}}}_2  s^{{\bar{0}}}_1  s^{{\bar{0}}}_2  s^{{\bar{0}}}_3  s^{{\bar{0}}}_2  s^{{\bar{0}}}_3  s^{{\bar{0}}}_2  
s^{{\bar{0}}}_1  s^{{\bar{0}}}_2  s^{{\bar{0}}}_3  s^{{\bar{0}}}_2  s^{{\bar{0}}}_1  s^{{\bar{0}}}_2  s^{{\bar{0}}}_3  s^{{\bar{0}}}_2  s^{{\bar{0}}}_1  s^{{\bar{0}}}_2  
s^{{\bar{0}}}_1  \\ \cdot  s^{{\bar{0}}}_2  s^{{\bar{0}}}_3  s^{{\bar{0}}}_2  s^{{\bar{0}}}_1  s^{{\bar{0}}}_2  s^{{\bar{0}}}_1  s^{{\bar{0}}}_2  
s^{{\bar{0}}}_4  s^{{\bar{1}}}_2  s^{{\bar{1}}}_3  s^{{\bar{1}}}_2  s^{{\bar{1}}}_3  s^{{\bar{1}}}_2  s^{{\bar{1}}}_1  s^{{\bar{1}}}_2  s^{{\bar{1}}}_3  s^{{\bar{1}}}_2  
s^{{\bar{1}}}_1  s^{{\bar{1}}}_2  s^{{\bar{1}}}_3  s^{{\bar{1}}}_2  s^{{\bar{1}}}_1  s^{{\bar{1}}}_2  s^{{\bar{1}}}_1  s^{{\bar{1}}}_2  s^{{\bar{1}}}_3  s^{{\bar{1}}}_2  
s^{{\bar{1}}}_1  s^{{\bar{1}}}_2  s^{{\bar{1}}}_1  s^{{\bar{1}}}_2  s^{{\bar{1}}}_1  s^{{\bar{1}}}_2  s^{{\bar{1}}}_3  s^{{\bar{1}}}_2  s^{{\bar{1}}}_3  s^{{\bar{1}}}_2  
s^{{\bar{1}}}_1  s^{{\bar{1}}}_2  s^{{\bar{1}}}_3  \\ \cdot  s^{{\bar{1}}}_2  s^{{\bar{1}}}_1  s^{{\bar{1}}}_2  s^{{\bar{1}}}_3  s^{{\bar{1}}}_2  
s^{{\bar{1}}}_1  s^{{\bar{1}}}_2  s^{{\bar{1}}}_1  s^{{\bar{1}}}_4  s^{{\bar{4}}}_1  s^{{\bar{4}}}_2  s^{{\bar{4}}}_3  s^{{\bar{4}}}_2  s^{{\bar{4}}}_3  s^{{\bar{4}}}_2  
s^{{\bar{4}}}_1  s^{{\bar{4}}}_2  s^{{\bar{4}}}_3  s^{{\bar{4}}}_2  s^{{\bar{4}}}_1  s^{{\bar{4}}}_2  s^{{\bar{4}}}_3  s^{{\bar{4}}}_2  s^{{\bar{4}}}_1  s^{{\bar{4}}}_2  
s^{{\bar{4}}}_1  s^{{\bar{4}}}_2  s^{{\bar{4}}}_3  s^{{\bar{4}}}_2  s^{{\bar{4}}}_1  s^{{\bar{4}}}_2  s^{{\bar{4}}}_1  s^{{\bar{4}}}_2  s^{{\bar{4}}}_1  s^{{\bar{4}}}_2  
s^{{\bar{4}}}_3  s^{{\bar{4}}}_2  s^{{\bar{4}}}_3  s^{{\bar{4}}}_2  s^{{\bar{4}}}_1  \\ \cdot  s^{{\bar{4}}}_2  s^{{\bar{4}}}_3  s^{{\bar{4}}}_2  
s^{{\bar{4}}}_1  s^{{\bar{4}}}_2  s^{{\bar{4}}}_3  s^{{\bar{4}}}_2  s^{{\bar{4}}}_1  s^{{\bar{4}}}_2  s^{{\bar{4}}}_1  s^{{\bar{4}}}_4  s^{{\bar{7}}}_1  s^{{\bar{7}}}_2  
s^{{\bar{7}}}_3  s^{{\bar{7}}}_2  s^{{\bar{7}}}_3  s^{{\bar{7}}}_2  s^{{\bar{7}}}_1  s^{{\bar{7}}}_2  s^{{\bar{7}}}_3  s^{{\bar{7}}}_2  s^{{\bar{7}}}_1  s^{{\bar{7}}}_2  
s^{{\bar{7}}}_3  s^{{\bar{7}}}_2  s^{{\bar{7}}}_1  s^{{\bar{7}}}_2  s^{{\bar{7}}}_1  s^{{\bar{7}}}_2  s^{{\bar{7}}}_3  s^{{\bar{7}}}_2  s^{{\bar{7}}}_1  s^{{\bar{7}}}_2  
s^{{\bar{7}}}_1  s^{{\bar{7}}}_2  s^{{\bar{7}}}_1  s^{{\bar{7}}}_2  s^{{\bar{7}}}_3  s^{{\bar{7}}}_2  s^{{\bar{7}}}_3  \\ \cdot  s^{{\bar{7}}}_2  
s^{{\bar{7}}}_1  s^{{\bar{7}}}_2  s^{{\bar{7}}}_3  s^{{\bar{7}}}_2  s^{{\bar{7}}}_1  s^{{\bar{7}}}_2  s^{{\bar{7}}}_3  s^{{\bar{7}}}_2  s^{{\bar{7}}}_1  s^{{\bar{7}}}_2  
s^{{\bar{7}}}_1  s^{{\bar{7}}}_4  s^{{\bar{6}}}_1  s^{{\bar{6}}}_2  s^{{\bar{6}}}_3  s^{{\bar{6}}}_2  s^{{\bar{6}}}_3  s^{{\bar{6}}}_2  s^{{\bar{6}}}_1  s^{{\bar{6}}}_2  
s^{{\bar{6}}}_3  s^{{\bar{6}}}_2  s^{{\bar{6}}}_1  s^{{\bar{6}}}_2  s^{{\bar{6}}}_3  s^{{\bar{6}}}_2  s^{{\bar{6}}}_1  s^{{\bar{6}}}_2  s^{{\bar{6}}}_1  s^{{\bar{6}}}_2  
s^{{\bar{6}}}_3  s^{{\bar{6}}}_2  s^{{\bar{6}}}_1  s^{{\bar{6}}}_2  s^{{\bar{6}}}_1  s^{{\bar{6}}}_2  s^{{\bar{6}}}_1  s^{{\bar{6}}}_2  s^{{\bar{6}}}_3  
\\ \cdot  s^{{\bar{6}}}_2  s^{{\bar{6}}}_3  s^{{\bar{6}}}_2  s^{{\bar{6}}}_1  s^{{\bar{6}}}_2  s^{{\bar{6}}}_3  s^{{\bar{6}}}_2  s^{{\bar{6}}}_1  
s^{{\bar{6}}}_2  s^{{\bar{6}}}_3  s^{{\bar{6}}}_2  s^{{\bar{6}}}_1  s^{{\bar{6}}}_2  s^{{\bar{6}}}_1  s^{{\bar{6}}}_4  s^{{\bar{3}}}_1  s^{{\bar{3}}}_2  s^{{\bar{3}}}_3  
s^{{\bar{3}}}_2  s^{{\bar{3}}}_3  s^{{\bar{3}}}_2  s^{{\bar{3}}}_1  s^{{\bar{3}}}_2  s^{{\bar{3}}}_3  s^{{\bar{3}}}_2  s^{{\bar{3}}}_1  s^{{\bar{3}}}_2  s^{{\bar{3}}}_3  
s^{{\bar{3}}}_2  s^{{\bar{3}}}_1  s^{{\bar{3}}}_2  s^{{\bar{3}}}_1  s^{{\bar{3}}}_2  s^{{\bar{3}}}_3  s^{{\bar{3}}}_2  s^{{\bar{3}}}_1  s^{{\bar{3}}}_2  s^{{\bar{3}}}_1  
s^{{\bar{3}}}_2  s^{{\bar{3}}}_1  \\ \cdot  s^{{\bar{3}}}_2  s^{{\bar{3}}}_3  s^{{\bar{3}}}_4  s^{{\bar{5}}}_1  s^{{\bar{5}}}_2  s^{{\bar{5}}}_3  
s^{{\bar{5}}}_2  s^{{\bar{5}}}_3  s^{{\bar{5}}}_2  s^{{\bar{5}}}_1  s^{{\bar{5}}}_2  s^{{\bar{5}}}_3  s^{{\bar{5}}}_2  s^{{\bar{5}}}_1  s^{{\bar{5}}}_2  s^{{\bar{5}}}_3  
s^{{\bar{5}}}_2  s^{{\bar{5}}}_1  s^{{\bar{5}}}_2  s^{{\bar{5}}}_1  s^{{\bar{5}}}_2  s^{{\bar{5}}}_3  s^{{\bar{5}}}_2  s^{{\bar{5}}}_1  s^{{\bar{5}}}_2  s^{{\bar{5}}}_1  
s^{{\bar{5}}}_2  s^{{\bar{5}}}_1  s^{{\bar{5}}}_2  s^{{\bar{5}}}_3  s^{{\bar{5}}}_2  s^{{\bar{5}}}_3  s^{{\bar{5}}}_2  s^{{\bar{5}}}_1  s^{{\bar{5}}}_2  s^{{\bar{5}}}_3  
s^{{\bar{5}}}_2  s^{{\bar{5}}}_1  s^{{\bar{5}}}_2  s^{{\bar{5}}}_3  \\ \cdot  s^{{\bar{5}}}_2  s^{{\bar{5}}}_1  s^{{\bar{5}}}_2  s^{{\bar{5}}}_1  
s^{{\bar{5}}}_2  s^{{\bar{5}}}_3  s^{{\bar{5}}}_2  s^{{\bar{5}}}_1  s^{{\bar{5}}}_2  s^{{\bar{5}}}_1  s^{{\bar{5}}}_4  s^{{\bar{3}}}_3  s^{{\bar{3}}}_2  s^{{\bar{3}}}_1  
s^{{\bar{3}}}_2  s^{{\bar{3}}}_3  s^{{\bar{3}}}_2  s^{{\bar{3}}}_1  s^{{\bar{3}}}_4  s^{{\bar{2}}}_1  s^{{\bar{2}}}_2  s^{{\bar{2}}}_3  s^{{\bar{2}}}_2  s^{{\bar{2}}}_3  
s^{{\bar{2}}}_2  s^{{\bar{2}}}_1  s^{{\bar{2}}}_2  s^{{\bar{2}}}_3  s^{{\bar{2}}}_2  s^{{\bar{2}}}_1  s^{{\bar{2}}}_2  s^{{\bar{2}}}_3  s^{{\bar{2}}}_2  s^{{\bar{2}}}_1  
s^{{\bar{2}}}_2  s^{{\bar{2}}}_1  s^{{\bar{2}}}_2  s^{{\bar{2}}}_3  s^{{\bar{2}}}_2  s^{{\bar{2}}}_1  \\ \cdot  s^{{\bar{2}}}_2  s^{{\bar{2}}}_1  
s^{{\bar{2}}}_2  s^{{\bar{2}}}_1  s^{{\bar{2}}}_2  s^{{\bar{2}}}_3  s^{{\bar{2}}}_2  s^{{\bar{2}}}_3  s^{{\bar{2}}}_2  s^{{\bar{2}}}_1  s^{{\bar{2}}}_2  s^{{\bar{2}}}_3  
s^{{\bar{2}}}_2  s^{{\bar{2}}}_1  s^{{\bar{2}}}_2  s^{{\bar{2}}}_3  s^{{\bar{2}}}_2  s^{{\bar{2}}}_1  s^{{\bar{2}}}_2  s^{{\bar{2}}}_1  s^{{\bar{2}}}_2  s^{{\bar{2}}}_3  
s^{{\bar{2}}}_2  s^{{\bar{2}}}_1  s^{{\bar{2}}}_2  s^{{\bar{2}}}_1  s^{{\bar{2}}}_4  s^{{\bar{3}}}_3  s^{{\bar{3}}}_2  s^{{\bar{3}}}_1  s^{{\bar{3}}}_2  s^{{\bar{3}}}_1  
s^{{\bar{3}}}_2  s^{{\bar{3}}}_3  s^{{\bar{3}}}_2  s^{{\bar{3}}}_1  s^{{\bar{3}}}_2  s^{{\bar{3}}}_1  s^{{\bar{3}}}_4  s^{{\bar{6}}}_3  \\ \cdot  
s^{{\bar{6}}}_2  s^{{\bar{6}}}_1  s^{{\bar{6}}}_2  s^{{\bar{6}}}_1  s^{{\bar{6}}}_4  s^{{\bar{7}}}_3  s^{{\bar{7}}}_2  s^{{\bar{7}}}_1  s^{{\bar{7}}}_2  s^{{\bar{7}}}_1  
s^{{\bar{7}}}_4  s^{{\bar{4}}}_3  s^{{\bar{4}}}_2  s^{{\bar{4}}}_1  s^{{\bar{4}}}_2  s^{{\bar{4}}}_1  s^{{\bar{4}}}_4  s^{{\bar{1}}}_3  s^{{\bar{1}}}_2  s^{{\bar{1}}}_1  
s^{{\bar{1}}}_2  s^{{\bar{1}}}_1  s^{{\bar{1}}}_2  s^{{\bar{1}}}_4 
\end{array}
\endgroup
$}}
\\ \quad \\
\mbox{$\FigFourA$: This is a Hamilton circuit of Weyl group of type $B_4$
with $s_i^2=e$, } \\
\mbox{$(s_1s_2)^4=(s_2s_3)^3=(s_3s_4)^3=e$, $(s_is_j)^2=2$ $(|i-j|\geq 2)$.
(Length $=384$)} \\
\mbox{This is obtained by the argument of Proof of Theorem~\ref{theorem:Main:CSW} for 
$\FigTwoA$} 
\\
\mbox{(that is $\FigThreeA$) 
and ${\hat{W}}:=\langle s_1,s_2,s_3\rangle$.
Here $s^{{\overline{k(r)}}}_{i_r}\cdots.s^{{\overline{k(2)}}}_{i_2}s^{{\overline{k(1)}}}_{i_1}$
means} \\
\mbox{$s_{i_r}\cdots s_{i_2}s_{i_1}\in{\hat{W}}_{k(r)}$ $(i_1=4,i_2=2,i_3=1,\ldots)$,
where ${\hat{W}}_0:={\hat{W}}$, 
${\hat{W}}_1:={\hat{W}}_0s_4$,} \\
\mbox{${\hat{W}}_2:={\hat{W}}_1s_3$,
${\hat{W}}_3:={\hat{W}}_2s_2$, ${\hat{W}}_4:={\hat{W}}_3s_1$,
${\hat{W}}_5:={\hat{W}}_4s_2$,
${\hat{W}}_6:={\hat{W}}_5s_3$ and} \\
\mbox{${\hat{W}}_7:={\hat{W}}_5s_4$.
Notice $W=\amalg_{z=0}^7\hatW_z$. See also $\FigFiveA$ 
below.}
\end{array}
\end{equation*}

\begin{equation*}
\begin{array}{l}
\hspace{1.5cm}    $\VisualExBfour$
\\ \quad \\
\mbox{$\FigFiveA$: Visual explanation to obtain the Hamilton circuit of 
$\FigFourA$
of} \\
\mbox{the Weyl group of type $B_4$}
\end{array}
\end{equation*}
}}
\end{example}

%
%


\begin{equation*}
\begin{array}{l}
\hspace{4cm} $\TotalCHHecThreeOneDASH$
\\
\mbox{$\FigSixA$: The Caylay graph and a ($2$-convenient (see Definition~\ref{definition:specialgraph})) Hamilton} \\
\mbox{circuit of the Weyl group of type-$A_3$}
\end{array}
\end{equation*}

\begin{equation*}
\begin{array}{l}
\hspace{6cm} $\TotalCHAtwoAone$
\\ \quad \\
\mbox{$\FigSevenA$: The Caylay graph and a ($2$-convenient (see Definition~\ref{definition:specialgraph})) Hamilton} \\
\mbox{circuit of the Weyl group of type-$A_2\times A_1$}
\end{array}
\end{equation*}

\begin{equation*}
\begin{array}{l}
\mbox{{\footnotesize{$\begin{array}{l}
s_3 s_2 s_3 s_2 s_3 s_4 s_3 s_4 s_3 
s_2 s_3 s_4 s_3 s_2 s_3 s_2 s_3 s_4 s_3 
s_4 s_3 s_2 s_1 s_4 s_3 s_2 s_3 s_2 s_3 
s_4 s_3 s_4 s_3 s_2 s_3 s_4 s_3 s_2 s_3 
s_2 \\ \cdot s_3 s_4 s_3 s_4 s_1 s_4 s_3 
s_2 s_3 s_2 s_3 s_4 s_3 s_4 s_3 s_2 s_3 
s_4 s_1 s_4 s_3 s_2 s_3 s_2 s_3 s_4 s_3 
s_4 s_3 s_2 s_3 s_4 s_1 s_4 s_3 s_2 s_3 
s_2 s_3 s_4 \\ \cdot s_3 s_4 s_3 s_2 s_3 
s_4 s_3 s_2 s_3 s_2 s_3 s_4 s_3 s_4 s_3 
s_2 s_1 s_2 s_3 s_2 s_3 s_4 s_3 s_4 s_3 
s_2 s_1 s_2 s_3 s_2 s_3 s_4 s_3 s_4 s_3 
s_2 s_1 s_2 s_1 s_4
\end{array}
$}}}
\\ \quad
\\ \mbox{$\FigEightA$: Hamilton circuit of Weyl group of type $A_4$
with $s_i^2=e$,}
\\ \mbox{$s_is_{i+1}s_i=s_{i+1}s_is_{i+1}$, $s_is_{i+2}=s_{i+2}s_i$ and $s_1s_4=s_4s_1$
(Length $=120$)
}
\end{array}
\end{equation*}

\begin{equation*}
\begin{array}{l}
\mbox{{\tiny{$\begin{array}{l}
s_2 s_3 s_2 s_3 s_2 s_4 s_2 s_4 s_2 
s_3 s_2 s_4 s_2 s_3 s_2 s_3 s_2 s_4 s_2 
s_4 s_2 s_3 s_2 s_1 s_2 s_3 s_2 s_3 s_2 
s_4 s_2 s_4 s_2 s_3 s_2 s_4 s_2 s_3 s_2 
s_3 s_2 s_4 s_2 s_4 s_2 s_1 s_2 s_4 s_2 
s_3 \\ \cdot  s_2 s_3 s_2 s_4 s_2 s_4 s_2 
s_3 s_2 s_4 s_2 s_3 s_2 s_3 s_2 s_1 s_2 
s_3 s_2 s_3 s_2 s_4 s_2 s_4 s_2 s_3 s_2 
s_4 s_2 s_3 s_2 s_3 s_2 s_4 s_2 s_1 s_2 
s_3 s_2 s_3 s_2 s_4 s_2 s_4 s_2 s_3 s_2 
s_4 s_2 s_3 \\ \cdot  s_2 s_3 s_2 s_1 s_2 
s_3 s_2 s_3 s_2 s_4 s_2 s_4 s_2 s_3 s_2 
s_4 s_2 s_3 s_2 s_3 s_2 s_4 s_2 s_1 s_2 
s_3 s_2 s_3 s_2 s_4 s_2 s_4 s_2 s_1 s_2 
 s_4 s_2 s_3 s_2 s_3 s_2 s_4 s_2 s_4 
s_2 s_3 s_2 s_4 s_2 s_3 \\ \cdot  s_2 s_3 
s_2 s_4 s_2 s_4 s_2 s_1 s_2 s_4 s_2 s_3 
s_2 s_3 s_2 s_4 s_2 s_4 s_2 s_3 s_2 s_1 
s_2 s_3 s_2 s_1 s_2 s_4 s_2 s_3 s_2 s_1 
s_2 s_3 s_2 s_1 s_2 s_4 s_2 s_1 s_2 s_1 
\end{array}
$}}}
\\ \quad
\\ \mbox{$\FigNineA$: Hamilton circuit of Weyl group of type $D_4$
with $s_i^2=e$, $s_is_js_i=s_js_is_j$}
\\ \mbox{$((i,j)\in\{(1,2),(2,3),(2,4)\})$ and
$s_is_j=s_js_i$ $((i,j)\in\{(1,3),(1,4),(3,4)\})$}
\\ 
\mbox{(Length $=192$)}
\end{array}
\end{equation*}

\begin{equation*}
\begin{array}{l} \hspace{0.7cm}
\mbox{{\tiny{$\begin{array}{l}
s_3 s_4 s_3 s_4 s_3 s_2 s_3 s_4 s_3 
s_2 s_3 s_4 s_3 s_2 s_3 s_2 s_3 s_4 s_3 
s_2 s_3 s_2 s_3 s_2 s_3 s_4 s_3 s_4 s_3 
s_2 s_3 s_4 s_3 s_2 s_3 s_4 s_3 s_2 s_3 
s_2 s_3 s_4 s_3 s_2 s_3 s_2 s_1 s_2 s_3 
s_4 \\ \cdot s_3 s_4 s_3 s_2 s_3 s_4 s_3 
s_2 s_3 s_4 s_3 s_2 s_3 s_2 s_3 s_4 s_3 
s_2 s_3 s_2 s_3 s_2 s_3 s_4 s_3 s_4 s_3 
s_2 s_3 s_4 s_3 s_2 s_3 s_4 s_3 s_2 s_3 
s_2 s_3 s_4 s_3 s_2 s_1 s_2 s_3 s_4 s_3 
s_4 s_3 s_2 \\ \cdot s_3 s_4 s_3 s_2 s_3 
s_4 s_3 s_2 s_3 s_2 s_3 s_4 s_3 s_2 s_3 
s_2 s_3 s_2 s_3 s_4 s_3 s_4 s_3 s_2 s_3 
s_4 s_3 s_2 s_3 s_4 s_3 s_2 s_3 s_2 s_3 
 s_4 s_1 s_2 s_3 s_4 s_3 s_4 s_3 s_2 
s_3 s_4 s_3 s_2 s_3 s_4 \\ \cdot s_3 s_2 
s_3 s_2 s_3 s_4 s_3 s_2 s_3 s_2 s_3 s_2 
s_3 s_4 s_3 s_4 s_3 s_2 s_3 s_4 s_3 s_2 
s_3 s_4 s_3 s_2 s_3 s_2 s_3 s_4 s_1 s_2 
s_3 s_4 s_3 s_4 s_3 s_2 s_3 s_4 s_3 s_2 
s_3 s_4 s_3 s_2 s_3 s_2 s_3 s_4 \\ \cdot 
s_3 s_2 s_3 s_2 s_3 s_2 s_3 s_4 s_3 s_4 
s_3 s_2 s_3 s_4 s_3 s_2 s_3 s_4 s_3 s_2 
s_1 s_2 s_3 s_4 s_3 s_4 s_3 s_2 s_3 s_4 
s_3 s_2 s_3 s_4 s_3 s_2 s_3 s_2 s_3 s_4 
s_3 s_2 s_3 s_2 s_3 s_2 s_3 s_4 s_3 s_4 
\\ \cdot s_3 s_2 s_3 s_4 s_3 s_2 s_3 s_4 
s_3 s_2 s_3 s_2 s_3 s_4 s_1 s_2 s_3 s_4 
s_3 s_4 s_3 s_2 s_3 s_4 s_3 s_2 s_3 s_4 
s_3 s_2 s_3 s_2 s_3 s_4 s_3 s_2 s_3 s_2 
s_3 s_2 s_3 s_4 s_3 s_4 s_3 s_2 s_3 s_4 
s_3 s_2 \\ \cdot s_3 s_4 s_3 s_2 s_3 s_2 
s_3 s_4 s_1 s_2 s_3 s_4 s_3 s_4 s_3 s_2 
s_3 s_4 s_3 s_2 s_3 s_4 s_3 s_2 s_3 s_2 
s_3 s_4 s_3 s_2 s_3 s_2 s_3 s_2 s_3 s_4 
s_3 s_4 s_3 s_2 s_3 s_4 s_3 s_2 s_3 s_4 
s_1 s_2 s_3 s_4 \\ \cdot s_3 s_4 s_3 s_2 
s_3 s_4 s_3 s_2 s_3 s_4 s_3 s_2 s_3 s_2 
s_3 s_4 s_3 s_2 s_3 s_2 s_3 s_2 s_3 s_4 
s_3 s_4 s_3 s_2 s_3 s_4 s_3 s_2 s_3 s_4 
s_3 s_2 s_3 s_2 s_3 s_4 s_1 s_2 s_3 s_4 
s_3 s_4 s_3 s_2 s_3 s_4 \\ \cdot s_3 s_2 
s_3 s_4 s_3 s_2 s_3 s_2 s_3 s_4 s_3 s_2 
s_3 s_2 s_3 s_2 s_3 s_4 s_3 s_4 s_3 s_2 
s_3 s_4 s_3 s_2 s_3 s_4 s_3 s_2 s_1 s_2 
s_3 s_4 s_3 s_4 s_3 s_2 s_3 s_4 s_3 s_2 
s_3 s_4 s_3 s_2 s_3 s_2 s_3 s_4 \\ \cdot 
s_3 s_2 s_3 s_2 s_3 s_2 s_3 s_4 s_3 s_4 
s_3 s_2 s_3 s_4 s_3 s_2 s_3 s_4 s_3 s_2 
s_1 s_2 s_3 s_4 s_3 s_4 s_3 s_2 s_3 s_4 
s_3 s_2 s_3 s_4 s_3 s_2 s_3 s_2 s_3 s_4 
s_3 s_2 s_3 s_2 s_3 s_2 s_3 s_4 s_3 s_4 
\\ \cdot s_3 s_2 s_3 s_4 s_3 s_2 s_3 s_4 
s_3 s_2 s_3 s_2 s_3 s_4 s_1 s_2 s_3 s_4 
s_3 s_4 s_3 s_2 s_3 s_4 s_3 s_2 s_3 s_4 
 s_3 s_2 s_3 s_2 s_3 s_4 s_3 s_2 s_3 
s_2 s_3 s_2 s_3 s_4 s_3 s_4 s_3 s_2 s_3 
s_4 s_3 s_2 \\ \cdot s_3 s_4 s_3 s_2 s_3 
s_2 s_3 s_4 s_1 s_2 s_3 s_4 s_3 s_4 s_3 
s_2 s_3 s_4 s_3 s_2 s_3 s_4 s_3 s_2 s_3 
s_2 s_3 s_4 s_3 s_2 s_3 s_2 s_3 s_2 s_3 
s_4 s_3 s_4 s_3 s_2 s_3 s_4 s_3 s_2 s_3 
s_4 s_1 s_2 s_3 s_4 \\ \cdot s_3 s_4 s_3 
s_2 s_3 s_4 s_3 s_2 s_3 s_4 s_3 s_2 s_3 
s_2 s_3 s_4 s_3 s_2 s_3 s_2 s_3 s_2 s_3 
s_4 s_3 s_4 s_3 s_2 s_3 s_4 s_3 s_2 s_3 
s_4 s_3 s_2 s_3 s_2 s_3 s_4 s_1 s_2 s_3 
s_4 s_3 s_4 s_3 s_2 s_3 s_4 \\ \cdot s_3 
s_2 s_3 s_4 s_3 s_2 s_3 s_2 s_3 s_4 s_3 
s_2 s_3 s_2 s_3 s_2 s_3 s_4 s_3 s_4 s_3 
s_2 s_3 s_4 s_3 s_2 s_3 s_4 s_1 s_2 s_3 
s_4 s_3 s_4 s_3 s_2 s_3 s_4 s_3 s_2 s_3 
s_4 s_3 s_2 s_3 s_2 s_3 s_4 s_3 s_2 
\\ \cdot s_3 s_2 s_3 s_2 s_3 s_4 s_3 s_4 
s_3 s_2 s_3 s_4 s_3 s_2 s_3 s_4 s_3 s_2 
s_1 s_2 s_3 s_4 s_3 s_4 s_3 s_2 s_3 s_4 
s_3 s_2 s_3 s_4 s_3 s_2 s_3 s_2 s_3 s_4 
s_3 s_2 s_3 s_2 s_3 s_2 s_3 s_4 s_3 s_4 
s_3 s_2 \\ \cdot s_3 s_4 s_3 s_2 s_3 s_4 
s_1 s_2 s_3 s_4 s_3 s_4 s_3 s_2 s_3 s_4 
s_3 s_2 s_3 s_4 s_3 s_2 s_3 s_2 s_3 s_4 
s_3 s_2 s_3 s_2 s_3 s_2 s_3 s_4 s_3 s_4 
s_3 s_2 s_3 s_4 s_3 s_2 s_3 s_4 s_3 s_2
  s_3 s_2 s_3 s_4 \\ \cdot s_1 s_2 s_3 s_4 
s_3 s_4 s_1 s_2 s_3 s_4 s_3 s_4 s_3 s_2 
s_3 s_4 s_3 s_2 s_3 s_4 s_3 s_2 s_3 s_2 
s_3 s_4 s_3 s_2 s_3 s_2 s_3 s_2 s_3 s_4 
s_3 s_4 s_3 s_2 s_3 s_4 s_3 s_2 s_3 s_4 
s_1 s_2 s_3 s_4 s_3 s_4 \\ \cdot s_3 s_2 
s_3 s_4 s_3 s_2 s_3 s_4 s_3 s_2 s_3 s_2 
s_3 s_4 s_3 s_2 s_3 s_2 s_3 s_2 s_3 s_4 
s_3 s_4 s_3 s_2 s_3 s_4 s_3 s_2 s_3 s_4 
s_3 s_2 s_3 s_2 s_3 s_4 s_3 s_2 s_3 s_2 
s_1 s_2 s_3 s_2 s_3 s_4 s_1 s_2 \\ \cdot 
s_3 s_4 s_3 s_4 s_3 s_2 s_3 s_4 s_3 s_2 
s_3 s_4 s_3 s_2 s_3 s_2 s_3 s_4 s_3 s_2 
s_3 s_2 s_3 s_2 s_3 s_4 s_3 s_4 s_3 s_2 
s_3 s_4 s_3 s_2 s_3 s_4 s_3 s_2 s_3 s_2 
s_3 s_4 s_3 s_2 s_3 s_2 s_1 s_2 s_3 s_2 
\\ \cdot s_1 s_2 s_3 s_4 s_3 s_2 s_3 s_4 
s_3 s_2 s_3 s_2 s_3 s_4 s_3 s_2 s_3 s_2 
s_3 s_2 s_3 s_4 s_3 s_4 s_1 s_2 s_3 s_4 
s_3 s_4 s_3 s_2 s_3 s_4 s_3 s_2 s_3 s_4 
s_3 s_2 s_3 s_2 s_3 s_4 s_3 s_2 s_3 s_2 
s_3 s_2 \\ \cdot s_3 s_4 s_3 s_4 s_3 s_2 
s_3 s_4 s_3 s_2 s_3 s_4 s_3 s_2 s_3 s_2 
s_3 s_4 s_3 s_2 s_3 s_2 s_1 s_2 s_3 s_4 
s_3 s_2 s_3 s_4 s_3 s_2 s_3 s_2 s_3 s_4 
s_3 s_2 s_3 s_2 s_1 s_2 s_3 s_2 s_1 s_2 
s_3 s_2 s_3 s_4 \\ \cdot s_3 s_2 s_3 s_2 
s_1 s_2 s_3 s_4 s_3 s_2 s_3 s_2 s_1 s_2 
s_3 s_2 s_3 s_4 s_3 s_2 s_3 s_2 s_1 s_2 
s_3 s_2 s_1 s_2 s_3 s_2 s_3 s_4 s_3 s_2 
s_3 s_2 s_1 s_2 s_3 s_2 s_1 s_2 s_3 s_2 
s_1 s_2 s_3 s_4 s_3 s_2 \\ \cdot s_3 s_2 
s_1 s_2 s_3 s_4 s_3 s_2 s_3 s_2 s_1 s_2 
s_3 s_2 s_1 s_2 s_3 s_2 s_3 s_4 s_3 s_2 
s_3 s_2 s_1 s_2 s_3 s_2 s_1 s_2 s_3 s_2 
s_1 s_2 s_3 s_4 s_3 s_2 s_3 s_2 s_1 s_2 
s_3 s_2 s_1 s_2 s_3 s_2 s_1 s_2 \\ \cdot
s_1 s_2
\end{array}
$}}}
\\ \quad
\\ \mbox{$\FigTenA$: Hamilton circuit of Weyl group of type $F_4$
with $s_i^2=e$,}
\\ \mbox{$(s_1s_2)^3=(s_2s_3)^4=(s_3s_4)^3=e$ and
$(s_is_j)^2=2$ $(|i-j|\geq 2)$
(Length $=1152$)}
\end{array}
\end{equation*}

\section{Weyl groupoids associated with Nichols algebras of diagonal type}
\subsection{Basic of Nichols algebras of diagonal-type}
\newcommand{\bK}{{\mathbb{K}}}
\newcommand{\bKt}{\bK^\times}
\newcommand{\chamfkI}{I}
\newcommand{\chamPi}{{\Pi}}
\newcommand{\chambZPi}{\bZ\chamPi}
\newcommand{\chamZgeqoPi}{\bZgeqo\chamPi}
\newcommand{\chamchi}{\chi}
\newcommand{\chamAutbZPi}{{\mathrm{Aut}_\bZ(\chambZPi)}}

\newcommand{\chamcalX}{{\mathcal{X}}}
\newcommand{\chamcalXPi}{\chamcalX_\chamPi}

\newcommand{\chamE}{E}
\newcommand{\chamUp}{U^+}
\newcommand{\chamUpchi}{\chamUp(\chamchi)}
\newcommand{\champartialK}{\partial^K}
\newcommand{\champartialL}{\partial^L}

\newcommand{\chamR}{R}
\newcommand{\chamRp}{\chamR^+}
\newcommand{\chamRpchi}{\chamRp(\chi)}
\newcommand{\champhipchi}{\varphi^\chi_+}

\newcommand{\chamq}{q}

\newcommand{\chamN}{N}
\newcommand{\chamNchi}{\chamN^\chi}

\newcommand{\chams}{s}
\newcommand{\chamschi}{\chams^\chi}

\newcommand{\chamtau}{\tau}
\newcommand{\chambartau}{{\bar{\chamtau}}}
\newcommand{\chamtauichi}{\chamtau_i\chi}
\newcommand{\chamtauiNchi}{\chamN^{\chamtauichi}}
\newcommand{\chamstauichi}{\chams^{\chamtauichi}}

\newcommand{\chamT}{T}
\newcommand{\chamTchi}{\chamT^\chi}

\newcommand{\chamUptauichi}{\chamUp(\chamtauichi)}

\newcommand{\chamRptauichi}{\chamRp(\chamtau_i\chi)}
\newcommand{\champhiptauichi}{\varphi^{\chamtau_i\chi}_+}

Let $\bK$ be a field of characteristic zero, and $\bKt:=\bK\setminus\{0\}$. 
For $x,y\in\bK$ and $k\in\bZgeqo$,
let 
$(k)_x:=\sum_{r=1}^kq^{r-1}$, $(k)_x!:=\prod_{r=1}^k(r)_x$, 
$(k;x,y):=1-x^{k-1}y$ and $(k;x,y)!:=\prod_{r=1}^k(r;x,y)$.
Let  $\chamfkI$ be an non-empty finite set.
Let $\chambZPi(=\oplus_{i\in\chamfkI}\bZ\al_i)$ be a free $\bZ$-module with a basis 
$\chamPi:=\{\al_i|i\in\chamfkI\}$,
where $|\chamPi|=n$.
Let $\chamZgeqoPi:=\oplus_{i\in\chamfkI}\bZgeqo\al_i$
$(\subset\chambZPi)$.
Let $\chamAutbZPi$ be the group formed by $\bZ$-module automorphism of 
$\chambZPi$.

Let $\chamcalXPi$ be the set formed by all maps
$\chamchi^\prime:\chambZPi\times\chambZPi\to\bKt$ with
\begin{equation}\label{eqn;condchi}
\chamchi^\prime(\lambda,\mu+\nu)
=\chamchi^\prime(\lambda,\mu)\chamchi^\prime(\lambda,\nu),
\,
\chamchi^\prime(\lambda+\mu,\nu)=
\chamchi^\prime(\lambda,\nu)\chamchi^\prime(\mu,\nu)
\quad (\lambda, \mu, \nu\in\chambZPi).
\end{equation}
We call an element of $\chamcalXPi$
{\it{a bi-character associated with $\chamfkI$}}.

In this subsection, let $\chamchi\in\chamcalXPi$, 
and let $\chamq_{ij}:=\chamchi(\al_i,\al_j)$ for $i$, $j\in\chamfkI$.

Let $\chamUpchi$ be the Nichols algebra of diagonal-type associated with $\chamchi$, that is to say, $\chamUpchi$ is the unital associative 
$\bK$-algebra with the generators $\chamE_i$ $(i\in\chamfkI)$
satisfying the conditions ({\rm{Nic}1)-({\rm{Nic}2) below.
\newline\newline
({\rm{Nic}1) There exist subspaces $\chamUpchi_\lambda$
$(\lambda\in\chamZgeqoPi)$ such that 
\begin{equation*}
\begin{array}{l}
\mbox{$\chamUpchi=\oplus_{\lambda\in\chamZgeqoPi}\chamUpchi_\lambda$
(as a $\bK$-linear space),} \\
\mbox{$\chamUpchi_\lambda\chamUpchi_\mu\subset
\chamUpchi_{\lambda+\mu}$
$(\lambda, \mu\in\chamZgeqoPi)$,} \\
\mbox{$\bK1=\chamUpchi_0$, $\dim\chamUpchi_0=1$, 
$\bK\chamE_i=\chamUpchi_{\al_i}$, $\dim\chamUpchi_{\al_i}=1$
$(i\in\chamfkI)$.} 
\end{array}
\end{equation*}
\newline
({\rm{Nic}2) There exist $\bK$-linear maps
$\champartialL_i:\chamUpchi\to\chamUpchi$
$(i\in\chamfkI)$ such that
\begin{equation*}
\begin{array}{cl}
\bullet & \champartialL_i(1)=0, \\
\bullet & \forall i,j\in\chamfkI, 
\forall\lambda\in\chamZgeqoPi,
\forall X\in\chamUpchi_\lambda, \\
&\quad\quad
\champartialL_i(\chamE_jX)=
\delta_{i,j}X
+\chamq_{ji}\chamE_j\champartialL_i(X), \\
\bullet
&
\cap_{i\in\chamfkI}\ker\champartialL_i=\chamUpchi_0.
\end{array}
\end{equation*}
Namely $\chamUpchi$ is
{\it{the Nichols algebra of diagonal-type
associated with $\chamchi$}}, see
\cite[Proposition~5.4]{Hec10}, \cite[Section~3]{HY10}
(see also \cite[CHAPTER~38]{Lus10}).

\begin{theorem}\label{theorem:KhaTh}
{\rm{(}}Kharchenko PBW Theorem {\rm{\cite{Kha99}}}{\rm{)}}
{\rm{(}}We emphasis that $\chamchi$ be anyone satisfying 
\eqref{eqn;condchi}.{\rm{)}}
There exists a unique pair $(\chamRpchi,\champhipchi)$
of a subset $\chamRpchi$ of $\chamZgeqoPi\setminus\{0\}$
and a map $\champhipchi:\chamRpchi\to\bN$ satisfying the following 
condition $(*)$, where notice  $\chamPi\subset
\chamRpchi\subset\chamZgeqoPi\setminus\{0\}$.
\newline\newline
$(*)$ 
There exist $\chamE_{\beta,t}\in\chamUpchi_\beta$ 
$(\beta\in\chamRpchi,t\in\fkJ_{1,\champhipchi(\beta)})$
and a total order $\preceq$ on $X:=\{(\beta,t)|\beta\in\chamRpchi,t\in\fkJ_{1,\champhipchi(\beta)}\}(\subset\chamRpchi\times\bN)$ such that 
for every $\lambda\in\chamZgeqoPi\setminus\{0\}$, 
a $\bK$-basis of $\chamUpchi_\lambda$ is formed by
the elements 
$\chamE_{\beta_1,t_1}^{k_1}\chamE_{\beta_2,t_2}^{k_2}
\cdots\chamE_{\beta_2,t_m}^{k_m}$
with
$m\in\bN$,
$(\beta_1,t_1)\prec(\beta_2,t_2)\prec\cdots\prec(\beta_m,t_m)$,
$(k_t)_{\chamchi(\beta_t,\beta_t)}!\ne 0$
$(t\in\fkJ_{1,m})$ and $\sum_{t=1}^mk_t\beta_t=\lambda$.
\end{theorem}

\begin{definition}\label{definition:DefUpi}
(1) Let $i,j\in\chamfkI$ with $i\ne j$. 
For $k\in\bZgeqo$,
let $\chamE^{\chamchi,+}_{0;i,j}:=\chamE^{\chamchi,-}_{0;i,j}:=\chamE_i$ (if $k=0$)
and $\chamE^{\chamchi,+}_{k;i,j}:=\chamE_i\chamE^{\chamchi,+}_{k-1;i,j}
-\chamq_{ii}^{k-1}\chamq_{ij}\chamE^{\chamchi,+}_{k-1;i,j}\chamE_i$,
$\chamE^{\chamchi,-}_{k;i,j}:=\chamE^{\chamchi,-}_{k-1;i,j}\chamE_i
-\chamq_{ii}^{k-1}\chamq_{ji}\chamE_i\chamE^{\chamchi,-}_{k-1;i,j}$
 (if $k\geq 1$).
Let $\chamNchi_{ij}
:=|\{k\in\bN|\chamE^{\chamchi,+}_{k;i,j}\ne 0\}|\,(\in\bZgeqo\cup\{\infty\})$.
\newline
(2) For $i\in\chamfkI$ and $\epsilon\in\{+,-\}$, let $\chamUpchi[\epsilon;i]$
be the non-unital  $\bK$-subalgebra of 
$\oplus_{\lambda\in\chamZgeqoPi\setminus\{0\}}\chamUpchi_\lambda$
generated by $\chamE^{\chamchi,\epsilon}_{k;i,j}$
with $j\in\chamfkI\setminus\{i\}$ and $k\in\bZgeqo$. 
\newline
(3) Let 
$\chamUpchi[\epsilon;i]_\lambda:=\chamUpchi_\lambda\cap\chamUpchi[\epsilon;i]$
for $\lambda\in\chamZgeqoPi$,
where notice \newline
$\chamUpchi[\epsilon;i]_{k\al_i+\al_j}=\bK\chamE^{\chamchi,\epsilon}_{k;i,j}$
$(\epsilon\in\{+,-\}, \{i,j\}\subset\chamfkI, i\ne j,k\in\bZgeqo)$.
\end{definition}

\begin{lemma}\label{lemma:easystring}
Let $i$, $j\in\chamfkI$ with $i\ne j$. Let $k\in\bZgeqo$. Then
\begin{equation*}
(k)_{\chamq_{ii}}!(k;\chamq_{ii},\chamq_{ij}\chamq_{ji})!=0
\Leftrightarrow \chamE^{\chamchi,+}_{k;i,j}=0 
\Leftrightarrow \chamE^{\chamchi,-}_{k;i,j}=0 
\Leftrightarrow k\al_i+\al_j\in\chamRpchi.
\end{equation*}
\end{lemma}

\begin{lemma}\label{lemma:kerpartiali}
{\rm{(\cite[Proposition~5.10]{Hec10}, \cite[Lemma 3.6~(ii)]{HY10})}}
Let $i\in\chamfkI$. Let $Y:=\{y\in\bZgeqo|(y)_{\chamq_{ii}}!\ne 0\}$,
and let $Z$ be the $|Y|$-dimensional
$\bK$-linear space with a basis $\{z_y|y\in Y\}$.
Then we have{\rm{:}}
\newline\newline
{\rm{(1)}} We have $Y=\{y\in\bZgeqo|\chamE_i^y\ne 0\}$.
\newline
{\rm{(2)}} We have 
$\ker\champartialL_i=\chamUpchi[-;i]$.
\newline
{\rm{(3)}} For $\epsilon\in\{+,-\}$, we have the $\bK$-linear isomorphism{\rm{:}}
\begin{equation*}
m:Z\otimes\chamUpchi[\epsilon;i]\to\chamUpchi\,\,
\mbox{defined by}\,\,
m(z_y\otimes x):=\chamE_i^yx
\,(y\in Y,\,x\in\chamUpchi[\epsilon;i]).
\end{equation*}

\end{lemma}

\begin{definition}\label{definition:tauichi}
Let $i\in\chamfkI$.
Assume $\chamNchi_{ij}<\infty$ for all $j\in\chamfkI\setminus\{i\}$.
Define 
$\chamschi_i\in\chamAutbZPi$
by $\chamschi_i(\al_i):=-\al_i$ and 
$\chamschi_i(\al_j):=\al_j+\chamNchi_{ij}\al_i$
$(j\in\chamfkI\setminus\{i\})$.
Define $\chamtauichi\in\chamcalXPi$ by
$\chamtauichi(\lambda,\mu):=\chi(\chamschi_i(\lambda),
\chamschi_i(\mu))$ $(\lambda,\mu\in\chambZPi)$.
\end{definition}

\begin{lemma}\label{lemma:Lusiso}
Let $i\in\chamfkI$.
Assume $\chamNchi_{ij}<\infty$ for all $j\in\chamfkI\setminus\{i\}$.
Then we have the following.
\newline\newline
{\rm{(1)}} We have $\chamtauiNchi_{ij}=\chamNchi_{ij}$ for all $j\in\chamfkI\setminus\{i\}$. In particular,
$\chamstauichi_i=\chamschi_i=(\chamschi_i)^{-1}$
and $\chamtau_i\chamtauichi=\chamchi$.
\newline\newline
{\rm{(2)}} There exists a unique {\rm{(}}non-unital{\rm{)}} 
$\bK$-algebra isomorphism 
$\chamTchi_i:\chamUpchi[-;i]\to\chamUptauichi[+;i]$
such that 
\begin{equation*}
\chamTchi_i(\chamE^{\chamchi,-}_{k;i,j})=
{\acute{\chamq}}_{ji}^{-k}\chamq_{ii}^{\frac {k(k-1-2\chamNchi_{ij})} 2}
(\prod_{r=0}^k(\chamNchi_{ij}-r+1)_{\chamq_{ii}}
(\chamNchi_{ij}-r+1;\chamq_{ii},
{\acute{\chamq}}_{ij}{\acute{\chamq}}_{ji}))\chamE^{\chamtauichi,+}_{\chamNchi_{ij}-k;i,j}
\end{equation*}
for $j\in\chamfkI\setminus\{i\}$ and $k\in\fkJ_{0,\chamNchi_{ij}}$,
where ${\acute{\chamq}}_{xy}:=\chamtauichi(\al_x,\al_y)$
$(x,y\in\chamfkI)$.
\end{lemma}

By Theorem~\ref{theorem:KhaTh} and Lemmas~\ref{lemma:kerpartiali}~(3)
and \ref{lemma:Lusiso}~(2), we have:
\begin{proposition}\label{proposition:preservemulti}
Let $i\in\chamfkI$.
Assume $\chamNchi_{ij}<\infty$ for all $j\in\chamfkI\setminus\{i\}$.
Then we have
\begin{equation*}
\chamschi_i(\chamRpchi\setminus\{\al_i\})=\chamRptauichi\setminus\{\al_i\}
\,\,\mbox{and}\,\,
(\champhiptauichi\circ\chamschi_i)_{|\chamRpchi\setminus\{\al_i\}}
={\champhipchi}_{|\chamRpchi\setminus\{\al_i\}}.
\end{equation*}

\end{proposition}

\subsection{Cayley graph associated with $\chamchi$ satisfying $|\chamRpchi|<\infty$}
\newcommand{\chamcalG}{{\mathcal{G}}}
\newcommand{\chamcalGchi}{\chamcalG(\chi)}
\newcommand{\chambarcalG}{{\bar{\chamcalG}}}
\newcommand{\chambarcalGchi}{\chambarcalG(\chi)}

\newcommand{\chamonechis}{{{1^\chamchi}\chams}}

\newcommand{\chamcalV}{{\mathcal{V}}}
\newcommand{\chamcalVchi}{\chamcalV(\chi)}
\newcommand{\chamcalE}{{\mathcal{E}}}
\newcommand{\chamcalEchi}{\chamcalE(\chi)}

\newcommand{\champi}{\pi}
\newcommand{\champichi}{\champi_\chamchi}

\newcommand{\chamGamma}{\Gamma}
\newcommand{\chamGammachi}{\chamGamma(\chi)}
In this subsection, let $\chamchi\in\chamcalXPi$, and
assume $|\chamRpchi|<\infty$. 
For a finite sequence
$i_1,i_2,\ldots,i_k$ in $\chamfkI$, 
let 
$\chamchi_{(i_k,\ldots,i_2,i_1)}:=
\chamtau_{i_k}\cdots\chamtau_{i_2}\chamtau_{i_1}\chamchi\,(\in\chamcalXPi)$,
and let
$\chams^\chamchi_{(i_k,\cdots,i_2,i_1)}
:=\chams^{\chamchi_k}_{i_k}
\cdots\chams^{\chamchi_2}_{i_2}\chams^{\chamchi_1}_{i_1}\,(\in\chamAutbZPi)$,
where $\chamchi_1:=\chamchi$
and 
$\chamchi_t:=\chamchi_{(i_{t-1},\ldots i_2,i_1)}$
$(t\in\fkJ_{2,k})$.
Let $\chamchi_{(\,)}:=\chamchi$ and 
$\chams^\chamchi_{(\,)}:=\rmid_{\chambZPi}$.
If $k=0$,
let $\chamchi_{(i_k,\cdots,i_2,i_1)}$ and
$\chams^\chamchi_{(i_k,\cdots,i_2,i_1)}$
mean $\chamchi_{(\,)}$ and $\chams^\chamchi_{(\,)}$
respectively.
\begin{definition}\label{definition:DefCAYLEYchi}
(Cayley graph associated with $\chamchi$) 
Let
\begin{equation*}
\begin{array}{lcl}
\chamcalVchi & := &
\{\chams^\chamchi_{(i_k,\cdots,i_2,i_1)}\in\chamAutbZPi\,|\,
k\in\bZgeqo,\,i_t\in\chamfkI\,(t\in\fkJ_{1,k})\}, \\
\chamcalGchi & := & 
\{\chamchi_{(i_k,\ldots,i_2,i_1)}\,\in\chamcalXPi\,|\,
k\in\bZgeqo,\,i_t\in\chamfkI\,(t\in\fkJ_{1,k})\}.
\end{array}
\end{equation*} 
Let ${\mathfrak{p}}_2(\chi)$
be the subset of the power set ${\mathfrak{P}}(\chamcalVchi)$ 
of $\chamcalVchi$ formed 
by all the sets having two different elements of $\chamcalVchi$.
Let 
\begin{equation*}
\chamcalEchi:=\{\{\chams^\chi_{(i_k,\ldots,i_2,i_1)},
\chams^\chi_{(i_{k+1},i_k,\ldots,i_2,i_1)}\}\in{\mathfrak{p}}_2(\chi)
\,|\,k\in\bZgeqo,\,i_t\in\chamfkI\,(t\in\fkJ_{1,k+1})\}.
\end{equation*}
For $\chamchi^\prime$, $\chamchi^{\prime\prime}\in\chamcalGchi$,
we write $\chamchi^\prime\equiv\chamchi^{\prime\prime}$
if 
\begin{equation*}
\mbox{$\chamchi^\prime(\al_i,\al_i)=\chamchi^{\prime\prime}(\al_i,\al_i)$
and $\chamchi^\prime(\al_i,\al_j)\chamchi^\prime(\al_j,\al_i)=\chamchi^{\prime\prime}(\al_i,\al_j)\chamchi^{\prime\prime}(\al_j,\al_i)$}
\end{equation*}
for all $i$, $j\in\chamfkI$. Clearly
$\equiv$ is an equivalent relation on $\chamcalGchi$.
Then for $\chamchi^\prime$, $\chamchi^{\prime\prime}
\in\chamcalGchi$, we have
\begin{equation}\label{eqn:prpequiv}
\chamchi^\prime\equiv\chamchi^{\prime\prime}
\,\,\Rightarrow\,\,
\chamRp(\chamchi^\prime)=\chamRp(\chamchi^{\prime\prime}),
\chams^{\chamchi^\prime}_i=\chams^{\chamchi^{\prime\prime}}_i,
\chamtau_i\chamchi^\prime\equiv\chamtau_i\chamchi^{\prime\prime}
\end{equation}
(see \cite[(4.28), Lemma~4.22]{AYY15} for example).
Let $\chambarcalGchi$ be the quotient set
of $\chamcalGchi$ by $\equiv$.
Let $\champichi:\chamcalGchi\to\chambarcalGchi$
be the canonical projection map, i.e.,
\begin{equation*}
\mbox{$\champichi(\chamchi^\prime):=\{\chamchi^{\prime\prime}
\in\chamcalGchi|\chamchi^{\prime\prime}
\equiv\chamchi^\prime\}\quad (\chamchi^\prime\in
\chamcalGchi)$.}
\end{equation*}
Let $\chamGammachi$ be the graph 
$(\chamcalVchi,\chamcalEchi)$, 
see also $(\star)$ below.
When we draw the graph $\chamGammachi$ in a visible way, we 
do it in the following way. The subgraph of $\chamGammachi$ formed by 
the the two vertices 
$\chams^\chamchi_{(i_k,\ldots,i_2,i_1)}$, $\chams^\chamchi_{(i_{k+1},i_k,\ldots,i_2,i_1)}\in\chamcalVchi$
and the edge $\{\chams^\chamchi_{(i_k,\ldots,i_2,i_1)}, \chams^\chamchi_{(i_{k+1},i_k,\ldots,i_2,i_1)}\}$ $\in\chamcalEchi$
is written below, where
$a:=\champichi(\chamchi_{(i_k,\ldots,i_2,i_1)})$
and 
$b:=\champichi(\chamchi_{(i_{k+1},i_k,\ldots,i_2,i_1)})$.
\begin{center}
\setlength{\unitlength}{1mm}
\begin{picture}(30,5)(0,0)
\put(-01,02){$a$}
\put(07,01.5){$i_{k+1}$}
\put(19,02){$b$}
\put(00,00){\circle{2}}
\put(01,00){\line(1,0){18}}
\put(20,00){\circle{2}}
\end{picture}
\end{center}
\end{definition}
\noindent
$(\star)$ In general,  we do not care 
(resp. do care)
that each vertex of the graph
$\chamGammachi$ exactly corresponds to an element of $\chamcalVchi$
(resp. is assigned by an element of 
$\chambarcalGchi$).
\newline\par
For $i\in\chamfkI$, define the map
$\chamtau_i:\chamcalGchi\to\chamcalGchi$
by $\chamtau_i(\chamchi^\prime):=\chamtau_i\chamchi^\prime$
$(\chamchi^\prime\in\chamcalGchi)$. 

\begin{definition}\label{definition:genrtsysOfChi}
Let $\chamchi^\prime\in\chamcalGchi$,
$a:=\champichi(\chamchi^\prime)\in\chambarcalGchi$,
$\chamRp(a):=\chamRp(\chamchi^\prime)$ and
$\chamR(a):=\chamRp(a)\cup(-\chamRp(a))$.
Let
$\chams^a_i:=\chams^{\chamchi^\prime}_i$
$(i\in\chamfkI)$, and define $|I|\times |I|$-matrix
$C^a=(c^a_{ij})_{i,j\in\chamfkI}$ over $\bZ$ by
$c^a_{ii}:=-2$ and $c^a_{ij}:=-\chamN^{\chamchi^\prime}_{ij}$
$(i\ne j)$.
For $i\in\chamfkI$, 
define the map
$\chambartau_i:\chambarcalGchi\to\chambarcalGchi$
by $\chambartau_i(\champichi(\chamchi^\prime))
=\champichi(\chamtau_i\chamchi^\prime)$
$(\chamchi^\prime\in\chamcalGchi)$.
For $a\in\chambarcalGchi$ and 
$k\in\bZgeqo$, $i_t\in\chamfkI$
$(t\in\fkJ_{1,k})$, 
define $a_{(i_k,\ldots,i_2,i_1)}\in\chambarcalGchi$
{\rm{(}}resp. $\chams^a_{(i_k,\ldots,i_2,i_1)}\in\chamAutbZPi${\rm{)}} in the same way 
as those for $\chamchi_{(i_k,\ldots,i_2,i_1)}$
{\rm{(}}resp. $\chams^\chamchi_{(i_k,\ldots,i_2,i_1)}${\rm{)}}
with $a$ and $\chambartau_i$ in place of
$\chamchi$ and $\chamtau_i$.
\newline\newline
{\rm{(1)}} For all $a\in\chambarcalGchi$ and all 
$i$, $j\in\chamfkI$, we have 
$\chambartau_i^2=\rmid_\chambarcalGchi$
and $c^{\chambartau_i(a)}_{ij}=c^a_{ij}$.
Namely ${\mathcal{C}}(\chi)={\mathcal{C}}(\chamfkI,\chambarcalGchi,
(\chambartau_i)_{i\in\chamfkI}, 
(C^a)_{a\in\chambarcalGchi})$ is a Cartan scheme,
see 
{\rm{\cite[Definition~2.1]{CH09}}},
{\rm{\cite[Definition~2.1]{HY10}}}
and {\rm{\cite[\S2.1]{AYY15}}} for this term.
\newline
{\rm{(2)}} We have the following $({\mathrm{R}}1)$-$({\mathrm{R}}4)$.
Namely ${\mathcal{R}}(\chi)={\mathcal{R}}({\mathcal{C}}(\chi),
(\chamR(a))_{a\in\chambarcalGchi})$ is
a generalized root system of type ${\mathcal{C}}(\chi)$,
see {\rm{\cite[Definition~2.2]{CH09}}},
{\rm{\cite[Definition~2.3]{HY10}}}
and
{\rm{\cite[Definition~1.2, Lemma~1.5, Remark~1.6]{AYY15}}}
for this term.
\newline\newline
$({\mathrm{R}}1)$ \quad $\chamR(a)=\chamRp(a)\cup(-\chamRp(a))$
$(a\in\chambarcalGchi)$. 
\newline
$({\mathrm{R}}2)$ \quad 
$\chamR(a)\cap\bZ\al_i=\{\al_i,-\al_i\}$
$(a\in\chambarcalGchi,i\in\chamfkI)$.
\newline
$({\mathrm{R}}3)$ \quad 
$\chams^a_i(\chamR(a))=\chamR({\bar{\chamtau_i}}(a))$
$(a\in\chambarcalGchi,i\in\chamfkI)$.
\newline
$({\mathrm{R}}4)$ \quad
For $a\in\chambarcalGchi$, $k\in\bN$, $i_t\in\chamfkI$ $(t\in\fkJ_{1,k})$,
if $\chams^a_{(i_k,\ldots,i_2,i_1)}=\rmid_{\chambZPi}$,
then $a_{(i_1,i_2\ldots,i_k)}=a$.
It also follows
that for  $a\in\chambarcalGchi$
and $i,j\in\chamfkI$ with $i\ne j$,
letting $m^a_{ij}:=|\chamRp(a)\cap(\bZgeqo\al_i
\oplus\bZgeqo\al_j)|$, we have
\begin{equation*}
\mbox{$a_{(\underbrace{\ldots,i,j,i}_{m^a_{ij}})}=a_{(\underbrace{\ldots,j,i,j}_{m^a_{ij}})}$
and
$\chams^a_{(\underbrace{\ldots,i,j,i}_{m^a_{ij}})}=
\chams^a_{(\underbrace{\ldots,j,i,j}_{m^a_{ij}})}$.}
\end{equation*}
\end{definition}

\newcommand{\chamvarpi}{\varpi}
\newcommand{\chamvarpitheta}{\chamvarpi_\theta}
\newcommand{\chamvarphitheta}{{\varphi_\theta}}
\newcommand{\chamvartheta}{\vartheta}
\newcommand{\chamvarthetatheta}{\chamvartheta_\theta}

\begin{definition}\label{definition:specialgraph}
Let $k:=|\chamcalVchi|\,(\in\bN)$.
Let $a:=\champichi(\chamchi)\in\chambarcalGchi$.
\newline\newline
(1) We say that a map $\theta:\bZ/k\bZ\to\chamcalVchi$ 
is {\it{a Hamilton circuit map}} of $\chamGammachi$
if $\theta$ is a bijection for which
$\{(\theta\circ\hatpi_k)(t),(\theta\circ\hatpi_k)(t+1)\}\in\chamcalEchi$ for all
$t\in\bZ$. \newline
(2) Let $\theta$ be a Hamilton circuit map of $\chamGammachi$. 
Let $a:=\champichi(\chamchi)\in\chambarcalGchi$,
and let $x\in\fkJ_{0,k-1}$ be such that $(\theta\circ\hatpi_k)(x)
=\rmid_{\chambZPi}\,(\in\chamcalVchi)$.
Let $j_t\in\chamfkI$ $(t\in\fkJ_{1,k})$ be such that 
$\chams^\chamchi_{(j_t,\ldots,j_2,j_1)}=(\theta\circ\hatpi_k)(x+t)$.
Define the map $\chamvarphitheta:\bZ/k\bZ\to\chamfkI$
(resp. $\chamvarthetatheta:\bZ/k\bZ\to\chambarcalGchi$)
by 
\begin{equation*}
\mbox{$(\chamvarphitheta\circ\hatpi_k)(x+t):=j_t$
(resp. $(\chamvarthetatheta\circ\hatpi_k)(x+t):=a_{(j_t,\ldots,j_2,j_1)}$)
for $t\in\fkJ_{1,k}$.}
\end{equation*}
\newline
(3) Let $\theta$ be a Hamilton circuit map of $\chamGammachi$. 
We call $\theta$ {\it{special}}
if for every $a\in\chambarcalGchi$ and 
every $i\in\chamfkI$, there exists $t\in\bZ$
such that $(\chamvarthetatheta\circ\hatpi_k)(t)=a$
and $(\chamvarphitheta\circ\hatpi_k)(t)=i$.
\newline
(4) Let $\theta$ be a Hamilton circuit map of $\chamGammachi$.
For $i\in\chamfkI$, we say that $\theta$ is
{\it{$i$-convenient}}
if there exists $r\in\fkJ_{0,1}$
such that 
$(\chamvarphitheta\circ\hatpi_k)(2t+r)=i$
for all $t\in\bZ$.
\end{definition}

By \eqref{eqn:prpequiv} and the definition of $\chambarcalGchi$,
we see:
\begin{lemma}\label{lemma:lemspecialgraph}
{\rm{(1)}} If there exists a special Hamilton map of $\chamGammachi$,
then for every $v\in\chamcalVchi$ and every $i\in\chamfkI$,
there exists a Hamilton map $\theta$ of $\chamGammachi$
such that $(\theta\circ\hatpi_k)(t)=v$
and $(\chamvarphitheta\circ\hatpi_k)(t)=i$
for some $t\in\bZ$. \newline
{\rm{(2)}} Let $i\in\chamfkI$. If there exists an $i$-convenient Hamilton map of $\chamGammachi$,
then for every $v\in\chamcalVchi$,
there exists a Hamilton map $\theta$ of $\chamGammachi$
such that $(\theta\circ\hatpi_k)(t)=v$
and $(\chamvarphitheta\circ\hatpi_k)(t)=i$
for some $t\in\bZ$.
\end{lemma}
\newcommand{\chamcalW}{{\mathcal{W}}}
\newcommand{\chamcalWchi}{\chamcalW(\chamchi)}

Let $\chamcalWchi$ be the semigroup 
defined by generators 
\begin{equation*}
0,\,\,e^a\,(a\in\chambarcalGchi),\,\,
z^a_i\,(a\in\chambarcalGchi,i\in\chamfkI)
\end{equation*} and relations
\begin{equation*}
\begin{array}{c}
0=00=0e^a=e^a0=0z^a_i=z^a_i0,\,\,e^ae^a=e^a,\,\,
e^ae^b=0\,(a\ne b),
\\
z^a_{(\underbrace{\ldots,i,j,i}_{m^a_{ij}})}=
z^a_{(\underbrace{\ldots,j,i,j}_{m^a_{ij}})}
\quad(i\ne j),
\end{array}
\end{equation*}
where $z^a_{(i_k,\ldots,i_2,i_1)}$
is defined in the same way as that for $\chams^a_{(i_k,\ldots,i_2,i_1)}$
with $z^a_i$ in place of $\chams^a_i$.
Let $z^a_{(\,)}:=e^a$.
Notice that by the definition \cite[Section~3]{HY08}
of a semigroup defined by generators and relations, 
we see:
\begin{lemma}\label{lemma:reconfirmW}
$z^a_{(i_k,\ldots,i_2,i_1)}=z^a_{(j_r,\ldots,j_2,j_1)}$
if and only if there exists $m\in\bN$, $x_t\in\bZgeqo$,
and $(i^{(t)}_{x_t},\ldots,i^{(t)}_2,i^{(t)}_1)\in\chamfkI^{x_t}$ $(t\in\fkJ_{0,m})$,
{\rm{(}}we let $\chamfkI^0$
mean $\{(\,)\}${\rm{)}}, such that the following 
$(r1)$-$(r2)$ are fulfilled.
\newline\newline
$(r1)$ $(i^{(0)}_{x_0},\ldots,i^{(0)}_2,i^{(0)}_1)
=(i_k,\ldots,i_2,i_1)$
and $(i^{(m)}_{x_m},\ldots,i^{(m)}_2,i^{(m)}_1)
=(j_r,\ldots,j_2,j_1)$.
\newline
$(r2)$ For every $t\in\fkJ_{0,m-1}$, 
there exist $t^\prime,t^{\prime\prime}\in\fkJ_{0,m}$
with $(t^\prime,t^{\prime\prime})\in\{(t,t+1),(t+1,t)\}$ such that
one of  $(r2\mbox{-}1)$-$(r2\mbox{-}2)$ below is fulfilled.
\newline
$(r2\mbox{-}1)$ We have $x_{t^{\prime\prime}}=x_{t^\prime}-2$
and there exist
$y\in\fkJ_{1,x_{t^\prime}-1}$ such that
\begin{equation*}
(i^{(t^{\prime\prime})}_{x_{t^{\prime\prime}}},
\ldots,i^{(t^{\prime\prime})}_{y+2},i^{(t^\prime)}_y,i^{(t^\prime)}_y,
i^{(t^{\prime\prime})}_{y-1},\ldots,i^{(t^{\prime\prime})}_2,i^{(t^{\prime\prime})}_1)
=(i^{(t^\prime)}_{x_{t^\prime}},\ldots,i^{(t^\prime)}_2,i^{(t^\prime)}_1).
\end{equation*}
\newline
$(r2\mbox{-}1)$ 
Let $a^\prime_0:=a$ and $a^\prime_u:=a_{(i^{(t^\prime)}_u,\ldots,i^{(t^\prime)}_2,i^{(t^\prime)}_1)}$
$(u\in\fkJ_{1,x_{t^\prime}})$.
For $u\in\fkJ_{1,x_{t^\prime}-1}$, if $i^{(t^\prime)}_u\ne i^{(t^\prime)}_{u+1}$,
let $m^\prime_u:=m^{a^\prime_{u-1}}_{i^{(t^\prime)}_{u+1}, i^{(t^\prime)}_u}$.
We have $x_{t^{\prime\prime}}=x_{t^\prime}$
and there exists
$y\in\fkJ_{1,x_{t^\prime}-1}$ with 
$i^{(t^\prime)}_y\ne i^{(t^\prime)}_{y+1}$ such that
\begin{equation*}
(i^{(t^\prime)}_{x_{t^\prime}},
\ldots,i^{(t^\prime)}_2,i^{(t^\prime)}_1)=(i^{(t^\prime)}_{x_{t^\prime}},
\ldots,
i^{(t^\prime)}_{y+m^\prime_{y-1}},
\underbrace{\ldots,i^{(t^\prime)}_y,i^{(t^\prime)}_{y+1},i^{(t^\prime)}_y}_{m^\prime_y},
i^{(t^\prime)}_{y-1},
\ldots,i^{(t^\prime)}_2,i^{(t^\prime)}_1)
\end{equation*} and 
\begin{equation*}
(i^{(t^{\prime\prime})}_{x_{t^{\prime\prime}}},
\ldots,i^{(t^{\prime\prime})}_2,i^{(t^{\prime\prime})}_1)=(i^{(t^\prime)}_{x_{t^\prime}},
\ldots,
i^{(t^\prime)}_{y+m^\prime_{y-1}},
\underbrace{\ldots,i^{(t^\prime)}_{y+1},i^{(t^\prime)}_y,i^{(t^\prime)}_{y+1}}_{m^\prime_y},
i^{(t^\prime)}_{y-1},
\ldots,i^{(t^\prime)}_2,i^{(t^\prime)}_1),
\end{equation*} where
in fact, $x_{t^\prime}\geq m^\prime_y$
and $m^\prime_{y^\prime}=m^\prime_y$
$(y^\prime\in\fkJ_{y,y+m^\prime_y-1})$.
\end{lemma}

\begin{theorem}\label{theorem:presentationW}
{\rm{(\cite[Theorem~1]{HY08})}}
We have
\begin{equation*}
\chams^a_{(i_k,\ldots,i_2,i_1)}=\chams^a_{(j_r,\ldots,j_2,j_1)}
\,\,
\Leftrightarrow\,\,
z^a_{(i_k,\ldots,i_2,i_1)}=z^a_{(j_r,\ldots,j_2,j_1)}
\end{equation*}
In particular,
$\chams^a_{(i_k,\ldots,i_2,i_1)}=\rmid_{\chambZPi}$
$\Leftrightarrow$ 
$\chams^a_{(i_k,\ldots,i_2,i_1)}=e^a$.
\end{theorem} We have
\begin{equation}\label{eqn:VchitwobN}
|\chamcalVchi|\in 2\bN.
\end{equation}
{\it{Proof of \eqref{eqn:VchitwobN}.}}
Let $\chamcalVchi^\prime:= 
\{\chams^\chamchi_{(i_{2k},\cdots,i_2,i_1)}\in\chamAutbZPi\,|\,
k\in\bZgeqo,\,i_t\in\chamfkI\,(t\in\fkJ_{1,2k})\}$
and $\chamcalVchi^{\prime\prime}:= 
\{\chams^\chamchi_{(i_{2k-1},\cdots,i_2,i_1)}\in\chamAutbZPi\,|\,
k\in\bN,\,i_t\in\chamfkI\,(t\in\fkJ_{1,2k-1})\}$.
By Lemma~\ref{lemma:reconfirmW} and
Theorem~\ref{theorem:presentationW}, we
have $\chamcalVchi=\chamcalVchi^\prime\cup\chamcalVchi^{\prime\prime}$
and $\chamcalVchi^\prime\cap\chamcalVchi^{\prime\prime}=\emptyset$.
Fix $j\in\chamfkI$ and define the map $f:\chamcalVchi^\prime\to\chamcalVchi^{\prime\prime}$
by $f(\chams^\chamchi_{(i_{2k},\cdots,i_2,i_1)}):=\chams^\chamchi_{(j,i_{2k},\cdots,i_2,i_1)}$.
By Lemma~\ref{lemma:reconfirmW} and
Theorem~\ref{theorem:presentationW}, we see that $f$ is bijective.
\hfill $\Box$
\begin{definition} \label{definition:haveSameCH}
Let $\chamchi^\prime\in\chamcalXPi$, and
assume $|\chamRp(\chi^\prime)|<\infty$.
We write $\chamchi\cong\chamchi^\prime$ if
there exists a bijection $f:\chamfkI\to \chamfkI$ such that 
\begin{equation*}
\forall k\in\bN,\,
\forall i_t\in\chamfkI\,(t\in\fkJ_{1,k}),\,
\chams^\chamchi_{(i_k,\ldots,i_2,i_1)}=
\chams^{\chamchi^\prime}_{(f(i_k),\ldots,f(i_2),f(i_1))}\,(\in\chamAutbZPi).
\end{equation*}
\end{definition}

Clearly we have:
\begin{lemma} \label{lemma:haveSameCHlemma}
Assume $\chamchi\cong\chamchi^\prime$.
Then there exists a unique bijection $f:\chamcalVchi
\to\chamcalV(\chamchi^\prime)$
such that $f(\rmid_{\chambZPi})=\rmid_{\chambZPi}$
and 
$\{\{f(x),f(y)\}|\{x,y\}\in\chamcalEchi\}
=\chamcalE(\chamchi^\prime)\,(\subset{\mathfrak{p}}_2(\chi^\prime))$.
In particular, if there exists a Hamilton circuit map of $\chamGammachi$,
then it is also the case for
$\chamGamma(\chamchi^\prime)$.
\end{lemma}

\begin{definition} \label{definition:quasiCartan}
Let $\chamq_{ij}:=\chamchi(\al_i,\al_j)$ $(i,j\in\chamfkI)$.
\newline\newline
(1) We say that $\chamchi$ is of (finite) Cartan-type
if $\chamq_{ii}\ne 1$ $(i\in\chamfkI)$
and $\chamq_{ii}^{\chamNchi_{ij}}\chamq_{ij}\chamq_{ii}=1$
$(i,j\in\chamfkI,\,i\ne j)$.
In addition, if there exists a positive definite 
symmetric bi-linear map 
$(\,,\,):(\bR\otimes_\bZ\chambZPi)\times(\bR\otimes_\bZ\chambZPi)
\to\bR$ such that $[{\frac {2(\al_i,\al_j)} {(\al_i,\al_i)}}]_{i,j\in\chamfkI}$
is a {\rm{(}}finite{\rm{)}} Cartan matrix of type $X$
and $\chamNchi_{ij}=-{\frac {2(\al_i,\al_j)} {(\al_i,\al_i)}}$
$(i,j\in\chamfkI,\,i\ne j)$,
then we say that $\chamchi$ is 
{\it{a {\rm{(}}finite{\rm{)}} Cartan-type bi-character of type $X$}}.
(Such $X$ really exists for any (finite) Cartan-type $\chamchi$, which can be directly seen by \cite{Hec09}.)
\newline
(2) We say that $\chamchi$ is of quasi-Cartan-type
if there exists $\chamchi^\prime\in\chamcalXPi$
with $|\chamRp(\chamchi^\prime)|<\infty$
such that $\chamchi^\prime\cong\chamchi$
and $\chamchi^\prime$ is of Cartan-type.
In addition, if  $\chamchi^\prime$
is  a (finite) Cartan-type
of type $X$, we say that 
{\it{a {\rm{(}}finite{\rm{)}} quasi-Cartan-type bi-character  of type $X$}}.
\end{definition}

Using Theorem~\ref{theorem:Main:CSW},
by \eqref{eqn:prpequiv} and the definitions of $\chambarcalGchi$
and $\cong$,
we see:
\begin{lemma} \label{lemma:HCquasiCartan}
Assume that $|\chamfkI|\geq 2$. 
Let $k:=|\chamcalVchi|\,(\in\bN)$
and $a:=\champichi(\chamchi)\in\chambarcalGchi$.
\newline\newline
{\rm{(1)}} Let $\chamchi$ be of Cartan-type of type $X$.
Then $|\chambarcalGchi|=1$, and the graph $\chamGammachi$ is isomorphic to the Cayley graph ${\mathcal{C}}_{W, S}$ 
of the Coxeter system $(W, S)$ of type $X$. 
In particular, a Hamilton circuit of $\chamGammachi$ exists and it
is special.
\newline
{\rm{(2)}} Let $\chamchi$ be of quasi-Cartan-type of type $X$.
Then  the graph $\chamGammachi$ is isomorphic to the Cayley graph ${\mathcal{C}}_{W, S}$ 
of the Coxeter system $(W, S)$ of type $X$. 
In particular, 
for every $v\in\chamcalVchi$ and every $i\in\chamfkI$,
there exists a Hamilton circuit map 
$\theta:\bZ/k\bZ\to\chamcalVchi$
such that $(\chamvarthetatheta\circ\hatpi_k)(t)=v$
and $(\chamvarphitheta\circ\hatpi_k)(t)=i$
for some $t\in\bZ$. 
\end{lemma}

The following lemma seems interesting although we will not use it.

\begin{lemma} \label{lemma:HCquasiCartanb}
Assume that $|\chamfkI|\geq 3$. 
Let $\chamchi$ be of quasi-Cartan-type.
Let $\chamq_{i^\prime,j^\prime}:=
\chamchi(\al_{i^\prime},\al_{j^\prime})$
$({i^\prime},{j^\prime}\in\chamfkI)$.
Let $k:=|\chamcalVchi|$.
{\rm{(}}Notice $k\in 2\bN$ since  
$\chamchi$ is of quasi-Cartan-type.{\rm{)}}
Let $i_1$, $i_2$, $i_3\in\chamfkI$
be such that $i_1\ne i_2\ne i_3\ne i_1$. 
Assume that $\chamq_{i_1,i_2}\chamq_{i_2,i_1}\ne 1$,
$\chamq_{i_2,i_3}\chamq_{i_3,i_2}\ne 1$
and $\chamq_{i_1,i_3}\chamq_{i_3,i_1}=1$.
Then there exists a Hamilton circuit map
$\theta:\bZ/k\bZ\to\chamcalVchi$ 
of $\chamGammachi$
such that 
\begin{equation*}
\forall t\in\bZ,\,\,
(\chamvarphitheta\circ\hatpi_k)(2t-1)
\notin\{i_1,i_3\}.
\end{equation*}
In particular, if $|\chamfkI|=3$, $\theta$ is an $i_2$-convenient Hamilton map.
\end{lemma}
\noindent
{\it{Proof.}}
We use the Dynkin diagrams of the finite Coxeter systems.
If $|\chamfkI|=3$, the claim follows from 
$\FigTwoA$ and $\FigSixA$.
If $|\chamfkI|\geq 4$, the claim is inductively proved 
using the
same argument as that of Proof of 
Theorem~\ref{theorem:Main:CSW}.
\hfill $\Box$

\begin{definition} \label{definition:verygoodchi}
Assume that $|\chamfkI|\geq 2$. 
We call $\chamchi$ {\it{special}} 
if there exists a special Hamilton circuit map of 
$\chamGammachi$.
\end{definition}

\subsection{Main tools}\label{subsection:maintools}
In this subsection, let $\chamchi\in\chamcalXPi$, and
assume $|\chamRpchi|<\infty$. 

By an argument similar to that of Proof of
Theorem~\ref{theorem:Main:CSW}, we have:
\begin{proposition}\label{proposition:maintool}
Assume $|I|\geq  3$. 
Fix $i\in\chamfkI$.
For $a\in\chambarcalGchi$,
let 
\begin{equation*}
\begin{array}{lcl}
W^a & := & \{\chams^a_{(u_k,\ldots,u_2,u_1)}|
k\in\bZgeqo,\,u_t\in\chamfkI\,(t\in\fkJ_{1,k})\}, 
\,(\subset\chamAutbZPi), \\
r & := & |W^a|\,\,
{\rm{(}}\mbox{Notice that $r=|W^{\grave{a}}|$ for all ${\grave{a}}\in\chambarcalGchi$}{\rm{)}}, \\
\chambarcalGchi^{\wedge,a} & := & \{a_{(j_k,\ldots,j_2,j_1)}|
k\in\bZgeqo,\,j_t\in\chamfkI\setminus\{i\}\,(t\in\fkJ_{1,k})\}
\,(\subset\chambarcalGchi), \\
{\hat{W}}^a & := & \{\chams^a_{(j_k,\ldots,j_2,j_1)}|
k\in\bZgeqo,\,j_t\in\chamfkI\setminus\{i\}\,(t\in\fkJ_{1,k})\}
\,(\subset\chamAutbZPi), \\
{\hat{r}}^a & := & |{\hat{W}}^a|.
\end{array}
\end{equation*}
{\rm{(}}${\hat{r}}^a\ne {\hat{r}}^{a^\prime}$ may happen for some 
$a$, $a^\prime\in\chambarcalGchi$ with $a\ne a^\prime$.{\rm{)}}
Let $B$ be a subset of $\chambarcalGchi$
such that $\cup_{b\in B}\chambarcalGchi^{\wedge,b}=\chambarcalGchi$
and $\chambarcalGchi^{\wedge,b}\cap\chambarcalGchi^{\wedge,b^\prime}=\emptyset$
$(b,b^\prime\in B,\,b\ne b^\prime)$.
Assume that there 
exist $x^b\in\bN$ $(b\in B)$
and surjections 
${\hat{\varphi}}^b_y:\bZ/{\hat{r}}^b\bZ\to\chamfkI\setminus\{i\}$ 
$(b\in B,\,y\in\fkJ_{1,x^b})$ satisfying the following
$(B1)$-$(B3)$.
In the following,
for $b\in B$, $y\in\fkJ_{1,x^b}$ and $t\in\bZgeqo$, let
\begin{equation*}
\mbox{$j_{b,y;t}:={\hat{\varphi}}^b_y(\hatpi_{{\hat{r}}^b}(t))$,
$b_{y;t}:=b_{(j_{b,y;t},\ldots,j_{b,y;2},j_{b,y;1})}$
and $\chams_{b,y;t}:=\chams^b_{(j_{b,y;t},\ldots,j_{b,y;2},j_{b,y;1})}$,}
\end{equation*} where $b_{y;0}:=b\,(=b_{(\,)})$ and  
$\chams_{b,y;0}:=\rmid_{\chambZPi}(=\chams^b_{(\,)})$.
\newline\newline
$(B1)$ For $b\in B$ and $y\in\fkJ_{1,x^b}$,
the map ${\hat{\theta}}^b_y:\bZ/{\hat{r}}^b\bZ\to{\hat{W}}^b$
defined by 
${\hat{\theta}}^b_y(\hatpi_{{\hat{r}}^b}(t)):=\chams_{b,y;t}$
$(t\in\bZgeqo)$ 
is bijective.
\newline
$(B2)$ 
For $b\in B$, $y\in\fkJ_{1,x^b}$ and $t\in\fkJ_{1,{\hat{r}}^b}$, 
there exist ${\dot t}\in\fkJ_{t,t+1}$,
$b^\prime\in B$, $y^\prime\in\fkJ_{1,x^{b^\prime}}$
and $t^\prime\in\fkJ_{1,{\hat{r}}^{b^\prime}}$ such that 
\begin{equation*}
\mbox{$\chambartau_i(b_{y;{\dot t}-1})=b^\prime_{y^\prime;
t^\prime-1}$
$j_{b,y;{\dot t}}=j_{b^\prime,y^\prime;t^\prime}$
and $m^{b_{y;{\dot t}-1}}_{i,j_{b,y;{\dot t}-1}}=2$,}
\end{equation*}
where notice
$\chams^{b^\prime_{y^\prime;t^\prime}}_{j_{b,y^\prime;t^\prime-1}}\chams^{b_{y;{\dot t}-1}}_i
=\chams^{b_{y;{\dot t}}}_i\chams^{b_{y;{\dot t}-1}}_{j_{b,y;{\dot t}}}$.
Then for $a\in\chambarcalGchi$,
there exists a surjection $\varphi^a:\bZ/r\bZ\to\chamfkI$
such that the map $\theta^a:\fkJ_{1,r}\to W^a$ 
defined by 
$\theta^a(h):=\chams^a_{(\varphi^a(\hatpi_r(h),\ldots,\varphi^a(\hatpi_r(2)),\varphi^a(\hatpi_r(1)))}$
$(h\in\fkJ_{1,r})$ is bijective,
where notice that 
$\theta^a(\hatpi_r(0))
=\theta^a(\hatpi_r(r))
=\rmid_{\chambZPi}(=\chams^a_{(\,)})$.
In other words, there exists a Hamilton circuit map of $\chamGammachi$.
\end{proposition}

\begin{proposition}\label{proposition:EdgeSpecial}
{\rm{(}}Recall \eqref{eqn:VchitwobN}.{\rm{)}}
Let $\chamchi\in\chamcalXPi$, and
assume $|\chamRpchi|<\infty$. 
Let $N\in\fkJ_{3,\infty}$ and assume $\chamfkI=\fkJ_{1,N}$.
Let $i\in\fkJ_{1,N-2}$.
Let $\chamfkI_1:=\fkJ_{1,i-1}$ and
$\chamfkI_2:=\fkJ_{i+1,N}$.
Let $k:=i+1$.
\begin{equation*}
\begin{array}{l}
\mbox{Assume that $\forall \chamchi^\prime\in\chamcalGchi$,
$\forall (i_1,i_2)\in(\chamfkI_1\times\chamfkI_2)
\cup(\{i\}\times(\chamfkI_2\setminus\{k\}))$,} \\
\mbox{$\chamchi^\prime(\al_{i_1},\al_{i_2})\chamchi^\prime(\al_{i_2},\al_{i_1})=1$.}
\end{array}
\end{equation*}
For $\chamchi^\prime\in\chamcalGchi$
and $r\in\fkJ_{1,2}$, 
if $\chamfkI_r\ne\emptyset$,
let $\chamPi_r:=\{\al_g|g\in\chamfkI_r\}$, 
$Y_r:=\oplus_{g\in\chamfkI_r}\bZ\al_g$, 
and $\chamchi^\prime_r:
=\chamchi^\prime_{|Y_r\times Y_r}$,  
regard $\chamchi^\prime_r\in\chamcalX_{\chamPi_r}$ in a natural way,
and let $h^{\chamchi^\prime}_r:=|\chamcalV(\chamchi^\prime_r)|$,
where notice $\chamfkI_2\ne\emptyset$.
For $\chamchi^\prime\in\chamcalGchi$
and $r\in\fkJ_{1,2}$, assume that if
$h^{\chamchi^\prime}_r\geq 4$,
there exists a Hamilton map 
$
\theta^{\chamchi^\prime}_r:\bZ/h^{\chamchi^\prime}_r\bZ\to\chamcalV(\chamchi^\prime_r)
$ of 
$\chamGamma(\chamchi^\prime_r)$,
where notice $h^{\chamchi^\prime}_2\geq 4$.
\begin{equation}\label{eqn:assEdgeSpecial}
\begin{array}{l}
\mbox{Assume that  $\theta^{\chamchi^\prime}_2$ is special for  
every $\chamchi^\prime\in\chamcalGchi$.} \\
\mbox{{\rm{(}}See also Remark~{\rm{\ref{remark:restrictionremark}}} below.{\rm{)}}}
\end{array}
\end{equation}
Then there exists a Hamilton map of $\chamGammachi$.
\end{proposition}
\noindent
{\it{Proof.}} 
Let $\chamchi^\prime\in\chamcalGchi$.
By Lemma~\ref{lemma:lemspecialgraph}~(1),
we have a special Hamilton map 
$\theta^{\chamchi^\prime}_2:
\bZ/h^{\chamchi^\prime}_2\bZ\to\chamcalV(\chamchi^\prime_2)$
of $\chamGamma(\chamchi^\prime_2)$
such that 
$(\theta^{\chamchi^\prime}_2
\circ\hatpi_{h^{\chamchi^\prime}_2})(0)=\rmid_{Y_2}$
and $(\varphi_{\theta^{\chamchi^\prime}_2}\circ
\hatpi_{h^{\chamchi^\prime}_2})(h^{\chamchi^\prime}_2)=k$
(See also Remark~{\rm{\ref{remark:restrictionremark}}} below.)
Let $\chamfkI_3:=\chamfkI\setminus\{i\}\,(=\chamfkI_1\cup\chamfkI_2)$,
$\chamPi_3:=\chamPi\setminus\{\al_i\}\,(=\chamPi_1\cup\chamPi_2)$,
$Y_3:=Y_1\oplus Y_2$, 
$\chamchi^\prime_3:=\chamchi^\prime_{|Y_3\times Y_3}$
and $h^{\chamchi^\prime}_3:=|\chamcalV(\chamchi^\prime_3)|$.
Then $h^{\chamchi^\prime}_3=h^{\chamchi^\prime}_1\cdot h^{\chamchi^\prime}_2$
if $\chamfkI_1\ne\emptyset$,
and $h^{\chamchi^\prime}_3=h^{\chamchi^\prime}_2$
if $\chamfkI_1=\emptyset$.
Define the Hamilton map 
$\theta^{\chamchi^\prime}_3
:\bZ/h^{\chamchi^\prime}_3\bZ\to\chamcalV(\chamchi^\prime_3)$
of $\chamGamma(\chamchi^\prime_3)$ by 
$(\theta^{\chamchi^\prime}_3\circ\hatpi_{h^{\chamchi^\prime}_3})(0):=\rmid_{Y_3}$ and
\begin{equation*}
\begin{array}{l}
(\varphi_{\theta^{\chamchi^\prime}_3}\circ
\hatpi_{h^{\chamchi^\prime}_3})(t) \\
\quad =
\left\{\begin{array}{l}
(\varphi_{\theta^{\chamchi^\prime}_2}\circ
\hatpi_{h^{\chamchi^\prime}_2})(t^\prime) \\
\quad \mbox{for $t=t^\prime+u\cdot h^{\chamchi^\prime}_2$ with $t^\prime\in\fkJ_{1,h^{\chamchi^\prime}_2-1}$
and $u\in 2\bZ$,} \\
(\varphi_{\theta^{\chamchi^\prime}_2}\circ
\hatpi_{h^{\chamchi^\prime}_2})(-t^\prime) \\
\quad \mbox{for $t=t^\prime+u\cdot h^{\chamchi^\prime}_2$ with $t^\prime\in\fkJ_{1,h^{\chamchi^\prime}_2-1}$
and $u\in 2\bZ+1$ if $|\chamfkI_1|\geq 1$,} \\
k \quad \mbox{for $t=h^{\chamchi^\prime}_2$ 
if $|\chamfkI_1|=0$,}
\\
1 \quad \mbox{for $t\in\{h^{\chamchi^\prime}_2,2h^{\chamchi^\prime}_2\}$ 
if $|\chamfkI_1|=1$,}
\\
(\varphi_{\theta^{\chamchi^\prime}_1}\circ
\hatpi_{h^{\chamchi^\prime}_1})(r) \quad 
\mbox{for $t=r\cdot h^{\chamchi^\prime}_2$
with $r\in\fkJ_{1,h^{\chamchi^\prime}_1}$
if $|\chamfkI_1|\geq 2$.} 
\end{array}\right.
\end{array}
\end{equation*}
$(t\in\fkJ_{1,h^{\chamchi^\prime}_3})$, where 
let $(\varphi_{\theta^{\chamchi^\prime}_1}\circ
\hatpi_{h^{\chamchi^\prime}_1})(1):=
(\varphi_{\theta^{\chamchi^\prime}_1}\circ
\hatpi_{h^{\chamchi^\prime}_1})(2):=1$ if $|\chamfkI_1|=1$, 
and see also $\FigElevenA$.
Notice that for  $\chamchi^{\prime\prime}\in\chamcalGchi$
and $t\in\fkJ_{1,h^{\chamchi^{\prime\prime}}_3}$,
there exists $r\in\fkJ_{0,1}$ such that
$(\varphi_{\theta^{\chamchi^{\prime\prime}}_3}\circ
\hatpi_{h^{\chamchi^{\prime\prime}}_3})(t+r)
\in\chamfkI_2\setminus\{k\}$.
Then the claim follows from Proposition~\ref{proposition:maintool}. 
\newline
\hfill $\Box$

\begin{remark} \label{remark:restrictionremark}
We can define the injection 
$x:\chamcalV(\chamchi^\prime_r)\to\chamcalV(\chamchi^\prime)$ 
by $x(\chams^{\chamchi^\prime_r}_{(i_p,\ldots,i_2,i_1)}):=\chams^{\chamchi^\prime_r}_{(i_p,\ldots,i_2,i_1)}$
$(p\in\bZgeqo$, $i_t\in\chamfkI_r$ $(t\in\fkJ_{1,p}))$.
Let $a^\prime:=\champichi(\chamchi^\prime)\,
(=\champi_{\chamchi^\prime}(\chamchi^\prime))\in\chambarcalGchi\,
(=\chambarcalG(\chamchi^\prime))$
and $a^\prime_r:=\champi_{\chamchi^\prime_r}(\chamchi^\prime_r)
\in\chambarcalG(\chamchi^\prime_r)$.
Let $Y:=\{a^\prime_{(i_p,\ldots,i_2,i_1)}|p\in\bZgeqo,\,i_t\in\chamfkI_r\,(t\in\fkJ_{1,p})\}$.
We can define the surjection
$y:Y\to\chamcalV(\chamchi^\prime)$
by $y(a^\prime_{(i_p,\ldots,i_2,i_1)}):=(a^\prime_r)_{(i_p,\ldots,i_2,i_1)}$
$(p\in\bZgeqo$, $i_t\in\chamfkI_r$ $(t\in\fkJ_{1,p}))$.
However $y$ is not necessarily injective.
See $\FigSeventyEightA$. 
\end{remark}

\begin{equation*}
\begin{array}{l}
$\JointForSpecialHighRank$
\\ \quad
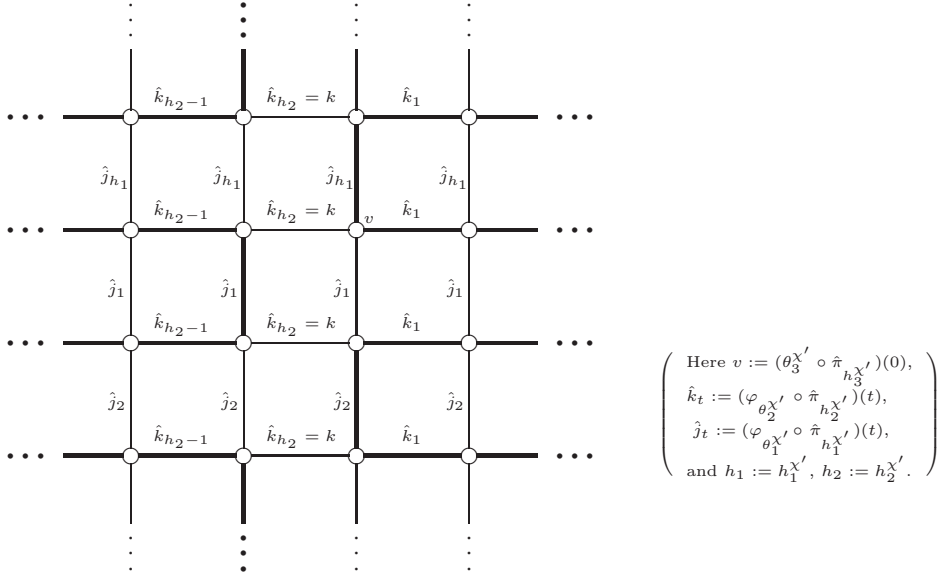
\\ \mbox{$\FigElevenA$: Hamilton circuit of $\chamGamma(\chamchi^\prime_3)$ in Proof of Proposition~\ref{proposition:EdgeSpecial}}
\end{array}
\end{equation*}

Using an argument similar to that of Proof of Proposition~\ref{proposition:EdgeAlmostSpecial},
we also have:
\begin{proposition}\label{proposition:EdgeAlmostSpecial}
Replace the assumption \eqref{eqn:assEdgeSpecial} of 
Proposition~{\rm{\ref{proposition:EdgeSpecial}}}
by that of \eqref{eqn:assEdgeAlmostSpecial} below.
\begin{equation}\label{eqn:assEdgeAlmostSpecial}
\begin{array}{l}
\mbox{For every $\chamchi^\prime\in\chamcalGchi$, there exists a Hamilton map} \\
\mbox{ 
$
(\theta^{\chamchi^\prime}_2)^\prime:\bZ/h^{\chamchi^\prime}_2\bZ\to\chamcalV(\chamchi^\prime_2)
$ of 
$\chamGamma(\chamchi^\prime_2)$ such that} \\
\mbox{
$\chambarcalG(\chamchi^\prime_2)\times(\chamfkI_2\setminus\{k\})$
$\subset\{((\chamvartheta_{\theta^{\chamchi^\prime}_2}\circ\hatpi_{h^{\chamchi^\prime}_2})(t),
(\varphi_{\theta^{\chamchi^\prime}_2}\circ
\hatpi_{h^{\chamchi^\prime}_2})(t))
|t\in\bZ\}$} \\
\mbox{\quad\quad\quad\quad\quad\quad\quad\quad\quad\quad
 $\cup\{((\chamvartheta_{(\theta^{\chamchi^\prime}_2)^\prime}\circ\hatpi_{h^{\chamchi^\prime}_2})
(t),(\varphi_{(\theta^{\chamchi^\prime}_2)^\prime}\circ
\hatpi_{h^{\chamchi^\prime}_2})(t))
|t\in\bZ
\}$.}
\end{array}
\end{equation} 
Then there exists a Hamilton map of $\chamGammachi$.
\end{proposition}

\newcommand{\chambarcalGi}{\chambarcalG_{\langle i\rangle}}
\newcommand{\chambarcalGisharpOne}[1]{\chambarcalGi({#1})}
\newcommand{\chambarcalGichi}{\chambarcalGisharpOne{\chamchi}}

\begin{definition} \label{definition:subGammachi}
Assume $|I|\geq  3$, and let $i\in\chamfkI$ and
$(\chambZPi)^\prime:=\oplus_{j\in\chamfkI\setminus\{i\}}\bZ\al_j$.
Let $\chamchi_{\langle i\rangle}:=
\chamchi_{|(\chambZPi)^\prime\times(\chambZPi)^\prime}\in\chamcalX_{\chamPi\setminus\{\al_i\}}$.
Notice that $\chamGamma(\chamchi_{\langle i\rangle})
=(\chamcalV(\chamchi_{\langle i\rangle}),
\chamcalE(\chamchi_{\langle i\rangle}))$ means
the one defined in the same way as that for 
$\chamGammachi=(\chamcalVchi,\chamcalEchi)$
with $\chamchi_{\langle i\rangle}$ 
and $\chamPi\setminus\{\al_i\}$
in place of 
$\chamchi$ 
and $\chamPi$ respectively. 
(Notice that $\chamRp(\chamchi_{\langle i\rangle})$
can be identified with
$\chamRpchi\cap(\chambZPi)^\prime$,
and that $\chamGamma(\chamchi_{\langle i\rangle})$
is the full subgraph of $\chamGammachi$ with
the vertices labeled by $\al_j$ with $j\in\chamfkI\setminus\{i\}$.)
Let $a:=\champichi(\chamchi)\in\chambarcalGchi$.
Let 
\begin{equation*}
\chambarcalGichi:=
\{\chambartau_{j_1}\chambartau_{j_2}\cdots\chambartau_{j_r}a
|r\in\bN,i_t\in\chamfkI\setminus\{i\}
\,(t\in\fkJ_{1,r})\}\,(\subset
\chambarcalGchi), 
\end{equation*} where notice 
$a\in\chambarcalGichi$.
For a non-empty subset $G$ of $\chamcalGchi$,
We say that  $G$ (resp. $\{\champichi(\chamchi^\prime)|\chamchi^\prime\in G\}$) is 
{\it{an $i$-complete $\chambartau$-representative subset of $\chamcalGchi$
{\rm{(}}resp. $\chambarcalGchi${\rm{)}}}} 
if $\chambarcalGchi
=\cup_{\chamchi^\prime\in G}\chambarcalGisharpOne{\chamchi^\prime}$
and
$\chambarcalGisharpOne{\chamchi^\prime}
\cap\chambarcalGisharpOne{\chamchi^{\prime\prime}}
=\emptyset$ $(\chamchi^\prime, \chamchi^{\prime\prime}\in G,\,
\chamchi^\prime\ne\chamchi^{\prime\prime})$.
\end{definition}

\begin{definition} \label{definition:handy}
Assume $|\chamfkI|\geq 3$. 
\newline\newline
{\rm{(1)}}
We call $\chamchi$ {\it{$(i,j)$-special}}
if there exists $i,j\in\chamfkI$ with $i\ne j$ such that 
for every $\chamchi^\prime\in\chamcalGchi$,
the equations $\chamchi^\prime(\al_i,\al_k)\chamchi^\prime(\al_k,\al_i)=1$
$(k\in\chamfkI\setminus\{i,j\})$ hold and
$\chamGamma(\chamchi^\prime_{\langle i\rangle})$
is special.
\newline
{\rm{(2)}}
We call $\chamchi$ {\it{$(i,j)$-convenient}}
if there exists $i,j\in\chamfkI$ with $i\ne j$ such that 
for every $\chamchi^\prime\in\chamcalGchi$, 
the equation $\chamchi^\prime(\al_i,\al_j)\chamchi^\prime(\al_j,\al_i)=1$
holds and $\chamGamma(\chamchi^\prime_{\langle i\rangle})$
is $j$-convenient.
\end{definition}

By Lemma~\ref{lemma:lemspecialgraph} and 
Proposition~\ref{proposition:maintool},
we have:
\begin{proposition}\label{proposition:maintoolcor}
Assume $|I|\geq  3$. 
If $\chamchi$ is $(i,j)$-special or $(i,j)$-convenient for some 
$i,j\in\chamfkI$ with $i\ne j$,
then there exists a Hamilton circuit map of
$\chamGammachi$.
\end{proposition}

\subsection{Further notation, terminology and examples}

\begin{definition} \label{definition:genDynkindiagram}
Let $\chamchi\in\chamcalXPi$.
Let $\chamq_{ij}:=\chamchi(\al_i,\al_j)$
$(\in\bKt)$
$(i,j\in\chamfkI)$.
{\it{The generalized Dynkin diagram}} of $\chamchi$ 
is
the graph composed of the $|\chamfkI|$-verticies 
one-to-one corresponding to $\al_i$'s $(i\in\chamfkI)$
and the edges connected vertices corresponding to
$\al_i$ and $\al_i$ $(i\ne j)$ with
$\chamq_{ij}\chamq_{ji}\ne 1$; moreover
$\al_i$ and $\chamq_{ii}$ are attached to the vertex corresponding to $\al_i$,
and $\chamq_{ij}\chamq_{ji}(\ne 1)$ is attached 
to the edge connecting 
$\al_i$ and $\al_j$.
For example, the generalized Dynkin diagram 
$\acchi^{(2,9)}_1$
of 
$\FigTwelveA$
means that $\chamfkI=\fkJ_{1,3}$
and 
$\chamq_{11}=q$, $\chamq_{22}=-1$, $\chamq_{33}=r$, 
$\chamq_{12}\chamq_{21}=q^{-1}$,
$\chamq_{23}\chamq_{32}=r^{-1}$ and   
$\chamq_{13}\chamq_{31}=1$.
\end{definition}

Let $\acn_{i,t}$ denote the cardinal number
of the generalized Dynkin diagrams drawn in 
\cite[Row~t of Table~i]{Hec09}.
For example, 
$\acn_{1,t_1}=1$ $(t_1\in\fkJ_{0,1})$, $\acn_{1,2}=2$, $\acn_{1,3}=1$, $\ldots$, 
$\acn_{2,t_2}=1$ $(t_2\in\fkJ_{1,3})$, $\acn_{2,4}=2$, $\acn_{2,3}=3$, $\ldots$, 
$\acn_{3,t_3}=1$ $(t_3\in\fkJ_{1,5})$, $\acn_{3,6}=3$, $\acn_{3,7}=4$, $\ldots$,
and  
$\acn_{4,11}=4$, $\acn_{4,12}=15$, $\acn_{4,13}=2$, $\ldots$.

Let $\acchi^{(i,t)}_k$ $(k\in\fkJ_{1,\acn_{i,t}})$ denote 
the the generalized Dynkin diagrams drawn in 
\cite[Row~t of Table~i]{Hec09}.
Let $\acchi^{(i,t)}_\bullet$ mean $\acchi^{(i,t)}_{k^\prime}$ for any $k^\prime\in\fkJ_{1,\acn_{i,t}}$.
Let $\acl^{(i,t)}$ denote the cardinal number of the vertices of $\acchi^{(i,t)}_k$.
To each vertex of $\acchi^{(i,t)}_k$, we attach the symbol 
$\acal_p$ $(p\in\fkJ_{1,\acl^{(i,t)}})$.
Define the order of $\acchi^{(i,t)}_k$'s 
and the order of vertices of $\acchi^{(i,t)}_k$
by the following $(\sharp1)$-$(\sharp3)$.
\newline\newline
$(\sharp1)$ Let $u$ and $v$ be two different vertices
of $\acchi^{(i,t)}_k$'s
in the middle box of \cite[Row~t of Table~i]{Hec09}.
They may belong to different $\acchi^{(i,t)}_k$ and $\acchi^{(i,t)}_{k^\prime}$.
We say that {\it{$u$ locates at a north-east place to $v$}}
if $u$ locates at a higher place than $v$,
or if they locate at the same hight but $u$ locates on the left side of $v$.
We also say that {\it{$u$ locates at an east-north place to $v$}}
if $u$ locates on the left side of $v$,
or if they locate at the same column but $u$ locates at a higher place than $v$.
\newline\newline
$(\sharp2)$ If the leftmost vertex of $\acchi^{(i,t)}_{k_1}$
locates at a north-east place to the one of $\acchi^{(i,t)}_{k_2}$, 
then $k_1<k_2$.
\newline\newline
$(\sharp3)$ 
For two different vertices of $\acchi^{(i,t)}_k$ attached by 
$\acal_{p_1}$ and $\acal_{p_2}$, if 
the vertex attached by $\acal_{p_1}$
locates at an east-north place to the one attached by $\acal_{p_2}$,
then $p_1<p_2$.
\newline\newline
For example, see $\FigTwelveA$ and $\FigThirteenA$.
Moreover for a bijection $f:\fkJ_{1,\acl^{(i,t)}}\to\fkJ_{1,\acl^{(i,t)}}$,
let $\acchi^{(i,t)}_k\circ f$ mean the generalized Dynkin diagram
obtained from $\acchi^{(i,t)}_k$ by replacing $\acal_j$ by $\acal_{f(j)}$
$(j\in\fkJ_{1,\acl^{(i,t)}})$, i.e., 
$\acchi^{(i,t)}_k\circ f$ is obtained from $\acchi^{(i,t)}_k$
by changing the order of its vertices by $j\to f(j)$ for
$j\in\fkJ_{1,\acl^{(i,t)}}$.
We also denote a bijection $f:\fkJ_{1,\acl^{(i,t)}}\to\fkJ_{1,\acl^{(i,t)}}$
by 
$\left[
\begin{array}{cccc}
1 & 2 & \cdots & \acl^{(i,t)} \\
f(1) & f(2) & \cdots & f(\acl^{(i,t)})
\end{array}\right]$.
See $\FigFourteenA$
for example.

\begin{equation*}
\begin{array}{l}
\hspace{1cm} $\EXAMPLEHecThreeNine$
\\ \quad
\\ \mbox{$\FigTwelveA$: Orders of the generalized Dynkin diagrams and}
\\ \mbox{their vertices for \cite[Row~9 of Table~2]{Hec09}}
\end{array}
\end{equation*}

\begin{equation*}
\begin{array}{l} \hspace{4cm}
$\EXAMPLEHecFourEightThree$
\\ \quad
\\ \mbox{$\FigThirteenA$: Order of vertices of \cite[3rd graph of Row~18 of Table~3]{Hec09}}
\end{array}
\end{equation*}

\begin{equation*}
\begin{array}{l}
\hspace{3cm} $\EXAMPLEHecFourEightThreeDASH$
\\ \quad
\\ \mbox{$\FigFourteenA$: Order changing of vertices of $\acchi^{(3,18)}_3$}
\end{array}
\end{equation*}

\begin{definition} \label{definition:irrchi}
Let $\chamchi\in\chamcalXPi$.
We say that $\chamchi$
is {\it{irreducibe}}
if the generalized Dynkin diagram of $\chamchi$ is connected.
\end{definition}

\section{Hamilton Circuits of Rank-3 Cases}\label{section:HCRkThree}

In this section, we directly draw a Hamilton circuit map 
of $\chamGammachi$ of each irreducible $\chamchi\in\chamcalXPi$
with $|\chamRpchi|<\infty$ and $|\chamfkI|=3$.

Proposition~\ref{proposition:areaorder} below
tells that for $\chamchi\in\chamcalXPi$ with $|\chamRpchi|<\infty$,
and $|\chamfkI|=3$, if we draw a graph inside a bounded 2-dimensional region fulfilling
such an appropriate property as $(d1)$-$(d5)$ of its statement,
it must be isomorphic to
$\chamGammachi$ as a graph.

\newcommand{\chamsigmaX}{\sigma^X}

\begin{proposition} \label{proposition:areaorder}
Let $\chamchi\in\chamcalXPi$
be such that $|\chamfkI|=3$ and $|\chamRpchi|<\infty$.
Let $X$ be a non-empty finite subset of $\bR^2$
fulfilling the following conditions 
$(d1)$-$(d6)$.
\newline\newline
$(d1)$ For $i\in\chamfkI$, we have 
$\chamsigmaX_i\in{\mathfrak{S}}_X$
with $(\chamsigmaX_i)^2=\rmid_X$
and $\chamsigmaX_i(x)\ne x$ $(x\in X)$,
where ${\mathfrak{S}}_X$ is the group formed by
all bijections from $X$ to $X$.
\newline
$(d2)$ We have
a surjection $\eta:X\to\chambarcalGchi$
with $\eta\circ\chamsigmaX_i=\chambartau_i\circ\eta$
for all $i\in\chamfkI$.
\newline
$(d3)$ For $x\in X$ and $i\in\chamfkI$, 
we have
a piecewise smooth curve
$\gamma_{x,i}:[0,1]\to\bR^2$ such that
$\gamma_{x,i}(0)=x$, $\gamma_{x,i}(1)=\chamsigmaX_i(x)$,
$\gamma_{\chamsigmaX_i(x),i}(v)=\gamma_{x,i}(1-v)$ $(v\in[0,1])$
and $\gamma_{x,i}$ is injective.
Moreover,
for $x\in X$ and $i,j\in\chamfkI$ with $i\ne j$, we have
$\gamma_{x,i}([0,1])\cap\gamma_{x,j}([0,1])=\{x\}$.
\newline
$(d4)$ Let $x\in X$, and let $i$, $j\in\chamfkI$ with $i\ne j$.
Let $m:=m^{\eta(x)}_{i,j}$. 
Then $(\chamsigmaX_j\chamsigmaX_i)^m(x)=x$.
Let $i_{2t-1}:=i$
and $i_{2t}:=j$ for $t\in\fkJ_{1,m}$.
Let $x_0:=x$ and $x_u:=\chamsigmaX_u(x_{u-1})$ 
for $u\in\fkJ_{1,2m}$ {\rm{(}}whence $x_{2m}=x${\rm{)}}.
Define the piecewise smooth closed curve
${\hat{\gamma}}_{x,i,j}:[0,2m]\to\bR^2$
by ${\hat{\gamma}}_{x,i,j}(t)=\gamma_{x_u,i_u}(t+u)$
$(u\in\fkJ_{0,m-1}$, $t\in[u,u+1])$.
Then ${\hat{\gamma}}_{x,i,j}$ is injective.
\newline
$(d5)$ 
For $k\in\fkJ_{1,2}$, let $x_k\in X$, $i_k$, $j_k\in\chamfkI$ with $i_k\ne j_k$
and let ${\hat{m}}_k:=m^{\eta(x_k)}_{i_k,j_k}$ and $H_{x_1,i_1,j_1}:={\hat{\gamma}}_{x,i,j}([0,2{\hat{m}}_k])$.
Then exactly one of the following $(d5)_1$-$(d5)_4$ holds.
\newline\newline
$(d5)_1$ \quad $H_{x_1,i_1,j_1}\cap H_{x_2,i_2,j_2}=\emptyset$. \newline
$(d5)_2$ \quad $H_{x_1,i_1,j_1}=H_{x_2,i_2,j_2}$, $\{i_1,j_1\}=\{i_2,j_2\}$ and 
${\hat{\gamma}}_{x_1,i_1,j_1}(u)=x_2$ for some $u\in\fkJ_{0,{\hat{m}}_1-1}$. \newline
$(d5)_3$ \quad $\exists u_1\in\fkJ_{0,{\hat{m}}_1-1}$, $\exists u_2\in\fkJ_{0,{\hat{m}}_2-1}$,
$H_{x_1,i_1,j_1}\cap H_{x_2,i_2,j_2}=\{{\hat{\gamma}}_{x_1,i_1,j_1}(u_1)\}=\{{\hat{\gamma}}_{x_2,i_2,j_2}(u_2)\}$. \newline
$(d5)_4$ \quad $\exists u_1\in\fkJ_{0,{\hat{m}}_1-1}$, $\exists u_2\in\fkJ_{0,{\hat{m}}_2-1}$,
$H_{x_1,i_1,j_1}\cap H_{x_2,i_2,j_2}={\hat{\gamma}}_{x_1,i_1,j_1}([u_1,u_1+1])={\hat{\gamma}}_{x_2,i_2,j_2}([u_2,u_2+1])$. 
\newline\newline
$(d6)$ 
For $x\in X$ and $i$, $j\in\chamfkI$ with $i\ne j$, 
let ${\widehat{H}}_{x,i,j}$
be the bounded simply connected closed subset of $\bR^2$
whose boundary is ${\hat{\gamma}}_{x,i,j}([0,2m])$.
Then $\cup_{x\in X}\cup_{i,j\in\chamfkI,i\ne j}{\widehat{H}}_{x,i,j}$ is simply connected.
\newline
$(d7)$ Let ${\hat{x}}\in X$ be such that $\champichi(\chamchi)=\eta({\hat{x}})$.
For $l\in\bZgeqo$ and $i_{t^\prime}\in\chamfkI$  
$(t^\prime\in\fkJ_{1,l})$,
let ${\hat{x}}_{(i_l,\ldots,i_2,i_1)}:=\chamsigmaX_{i_l}\cdots\chamsigmaX_{i_2}\chamsigmaX_{i_1}({\hat{x}})$,
where ${\hat{x}}_{(i_l,\ldots,i_2,i_1)}:={\hat{x}}_{(\,)}:=x$ if $l=0$.
Then $X=\{{\hat{x}}_{(i_t,\ldots,i_2,i_1)}|l\in\bZgeqo,
i_{t^\prime}\in\chamfkI\,(t^\prime\in\fkJ_{1,l})\}$.
\newline\newline
Then there exists a bijection $\varphi:\chamcalVchi\to X$ such that
$\varphi(\chams^\chamchi_{(i_l,\ldots,i_2,i_1)})={\hat{x}}_{(i_l,\ldots,i_2,i_1)}$
for all $l\in\bZgeqo$ and $i_{t^\prime}\in\chamfkI$  
$(t^\prime\in\fkJ_{1,l})$,
where notice $\varphi(\rmid_{\chambZPi})={\hat{x}}$.
\end{proposition}
\noindent
{\it{Proof.}}
Let $\bR X$ be an $|X|$-dimensional 
$\bR$-linear space whose basis is $X$.
For $b\in\chambarcalGchi$ and $i\in\chamfkI$, 
define $f^b_i\in{\mathrm{End}}_\bR(X)$ 
by $f^b_i(x):=\delta_{b,\eta(x)}\chamsigmaX_i(x)$
$(x\in X)$.
By $(d1)$ and $(d4)$,
we have the semigroup homomorphism 
$\zeta:\chamcalWchi\to{\mathrm{End}}_\bR(X)$
with $\zeta(z^b_i):=f^b_i$
$(b\in\chambarcalGchi,i\in\chamfkI)$.
By Theorem~\ref{theorem:presentationW},
there exists a surjecition $\varphi:\chamcalVchi\to X$
of the statement.

Let $l\in\bN$ and $i_t\in\chamfkI$ $(t\in\fkJ_{1,l})$
be such that ${\hat{x}}_{(i_l,\ldots,i_2,i_1)}={\hat{x}}$.
Let ${\hat{x}}_0:={\hat{x}}$
and ${\hat{x}}_t:={\hat{x}}_{(i_t,\ldots,i_2,i_1)}$
$(t\in\fkJ_{1,l})$. Then ${\hat{x}}_l={\hat{x}}$.
We have $l\geq 2$ by $(d1)$.
Let $h_1$, $h_2\in\fkJ_{0,l-1}$ be such that
$h_1<h_2$,
${\hat{x}}_{h_1}={\hat{x}}_{h_2}$
and ${\hat{x}}_{h^\prime}\ne{\hat{x}}_{h_1}$
for $h^\prime\in\fkJ_{h_1+1,h_2-1}$.
We have $h_1\leq h_2-2$ by $(d1)$.
If $h_1=h_2-2$, we have $i_{h_1+1}=i_{h_2}$
by $(d1)$ and $(d3)$,
whence $\chams^\chamchi_{(i_l,\ldots,i_2,i_1)}=
\chams^\chamchi_{(i_l,\ldots,i_{h_1+3},i_{h_1},\ldots,i_2,i_1)}$,
where if $h_2=l$, let $\chams^\chamchi_{(i_l,\ldots,i_{h_1+3},i_{h_1},\ldots,i_2,i_1)}$
mean $\chams^\chamchi_{(i_{l-2},\ldots,i_2,i_1)}$.
Assume $h_2\geq h_1+3$.
Let $Z$ be the bounded simply connected closed subset of $\bR^2$
whose boundary is composed of
$\gamma_{x_{t-1},i_t}([0,1])$ 
$(t\in\fkJ_{h_1,h_2-1})$.
Let $||Z||$ mean the area of $Z$.
Then $0<||Z||<\infty$. By $(d4)$, $(d5)$ and $(d6)$, using an induction on $||Z||$,
we see that $\chams^\chamchi_{(i_l,\ldots,i_2,i_1)}=
\chams^\chamchi_{(i_l,\ldots,i_{h_2+1},i_{h_1},\ldots,i_2,i_1)}$.
Thus we prove $\chams^\chamchi_{(i_l,\ldots,i_2,i_1)}=\rmid_{\chambZPi}$.
Hence by $(d7)$, the inverse map $\varphi^{-1}$ of $\varphi$ can be defined.
In particular, $\varphi$ is bijective.
\hfill $\Box$
\newline\par
Using Proposition~\ref{proposition:areaorder}, we directly have:

\begin{theorem} \label{theorem:existRkThree}
Assume that $|\chamfkI|=3$.
Then for every $\chamchi$ of {\rm{\cite[Table~2]{Hec09}}},
there exists a Hamilton circuit map of 
$\chamGammachi$.
In fact, 
{\rm{(}}up to elements of $\chambarcalGchi${\rm{)}}
it is drawn by 
$\FigTwentyA$ {\rm{(}}$\acchi^{(2,1)}_\bullet$, $\acchi^{(2,4)}_\bullet$, $\acchi^{(2,8)}_\bullet${\rm{)}},
$\FigTwentyFiveA$ {\rm{(}}$\acchi^{(2,2)}_\bullet$, $\acchi^{(2,3)}_\bullet$, $\acchi^{(2,5)}_\bullet$,
$\acchi^{(2,12)}_\bullet$, $\acchi^{(2,14)}_\bullet${\rm{)}},
$\FigTwentySevenA$ {\rm{(}}$\acchi^{(2,6)}_\bullet${\rm{)}},
$\FigTwentyNineA$ {\rm{(}}$\acchi^{(2,7)}_\bullet${\rm{)}},
$\FigThirtyTwoA$ {\rm{(}}$\acchi^{(2,9)}_\bullet$, $\acchi^{(2,10)}_\bullet$, $\acchi^{(2,11)}_\bullet${\rm{)}},
$\FigThirtyFiveA$ {\rm{(}}$\acchi^{(2,13)}_\bullet$, $\acchi^{(2,18)}_\bullet${\rm{)}},
$\FigThirtySevenA$ {\rm{(}}$\acchi^{(2,15)}_\bullet${\rm{)}},
$\FigThirtyNineA$ {\rm{(}}$\acchi^{(2,16)}_\bullet${\rm{)}} 
and $\FigFourtyOneA$ {\rm{(}}$\acchi^{(2,17)}_\bullet${\rm{)}}.
\end{theorem}

\noindent
{\bf{Rank-$3$-Case-$0$:}} We shall also need 
the Hamilton circuit map of 
$\FigSixteenA$ 
for $\chamchi$ of $A_2\times A_1$-case.  
\begin{equation*}
\begin{array}{l}
\hspace{1cm} $\AtwoAone$
\\ \quad
\\ \mbox{$\FigFifteenA$: the generalized Dynkin diagrams for $A_2\times A_1$}
\end{array}
\end{equation*}

\begin{equation*}
\begin{array}{l}
\hspace{4cm} $\aTotalCHAtwoAone$
\\ \quad
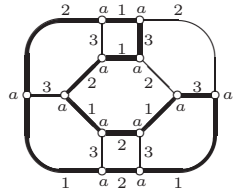
\\ \mbox{$\FigSixteenA$: Special and $1$-convenient Hamilton circuit for $A_2\times A_1$}
\end{array}
\end{equation*}

\noindent
{\bf{Rank-$3$-Case-$1$:}} Let $\chamchi$ be of $\acchi^{(2,8)}_1$.
We draw a Hamilton circuit map of $\chamGammachi$
by $\FigTwentyA$.
Notice that
$\chamchi$ is of \cite[Table~2, Row~8]{Hec09}
and that $\chamchi_1\cong\chamchi_2$
for any
two $\chamchi_1$ and $\chamchi_2$
of \cite[Table~2, Row~1,4,8]{Hec09}.
\begin{equation*}
\begin{array}{l}
$\HecOnea$
\\ \quad
\\ \mbox{$\FigSeventeenA$: The generalized Dynkin diagrams for \cite[Table~2, Row~1]{Hec09}}
\end{array}
\end{equation*}

\begin{equation*}
\begin{array}{l} \hspace{0.5cm}
$\HecFourab$
\\ \quad
\\ \mbox{$\FigEighteenA$: Generalized Dynkin diagrams for \cite[Table~2, Row~4]{Hec09}}
\end{array}
\end{equation*}

\begin{equation*}
\begin{array}{l}
$\HecEightabcd$
\\ \quad
\\ \mbox{$\FigNineteenA$: The generalized Dynkin diagrams for \cite[Table~2, Row~8]{Hec09}}
\end{array}
\end{equation*}


\begin{equation*}
\begin{array}{l} \hspace{1.5cm}
$\TotalCHHecThreeEight$
\\ \quad
\\ \mbox{$\FigTwentyA$:  Hamilton circuit for \cite[Table~2, Row~8 ($\cong$ Rows~1,4)]{Hec09}}
\end{array}
\end{equation*}

\begin{equation*}
\begin{array}{l} \hspace{1.5cm}
$\TotalCHHecThreeFour$
\\ \quad
\\ \mbox{$\FigTwentyOneA$: $2$-convenient and special Hamilton circuit} 
\\ \mbox{for \cite[Table~2, Row~4 ($\cong$ Rows~1,8)]{Hec09}}
\end{array}
\end{equation*}

\noindent
{\bf{Rank-$3$-Case-$2$:}} Let $\chamchi$ be of $\acchi^{(2,5)}_1$ or $\acchi^{(2,14)}_1$.
We draw Hamilton circuit maps of $\chamGammachi$
by $\FigTwentyFiveA$. 
Notice that
$\chamchi$ is of \cite[Table~2, Rows~5,14]{Hec09}
and that $\chamchi_1\cong\chamchi_2$
for any
two $\chamchi_1$ and $\chamchi_2$ of
\cite[Table~2, Row~2,3,5,12,14]{Hec09}.

\begin{equation*}
\begin{array}{l}
$\HecTwoThreeTwelvedef$
\\ \quad
\\ \mbox{$\FigTwentyTwoA$: Generalized Dynkin diagrams of \cite[Table~2, Rows~2,3,12]{Hec09}}
\end{array}
\end{equation*}

\begin{equation*}
\begin{array}{l}
$\HecFiveabc$
\\ \quad
\\ \mbox{$\FigTwentyThreeA$: Generalized Dynkin diagrams for \cite[Table~2,\,Row~5]{Hec09}}
\end{array}
\end{equation*}

\begin{equation*}
\begin{array}{l}
$\HecFourteenabc$
\\ \quad
\\ \mbox{$\FigTwentyFourA$: The generalized Dynkin diagrams for \cite[Table~2,\,Row~14]{Hec09}}
\end{array}
\end{equation*}



\begin{equation*}
\begin{array}{l}
$\spTotalCHHecThreeFive$
\\ \quad
\\ \mbox{$\FigTwentyFiveA$: Special Hamilton circuit of \cite[Table~2, Rows~5,14 (and Rows~2,3,12)]{Hec09}}
\end{array}
\end{equation*}

\noindent
{\bf{Rank-$3$-Case-$3$:}} Let $\chamchi$ be of $\acchi^{(2,6)}_1$.
We draw a Hamilton circuit map of $\chamGammachi$
by $\FigTwentySevenA$. 
Notice that
$\chamchi$ is of \cite[Table~2, Row~6]{Hec09}.

\begin{equation*}
\begin{array}{l}
$\HecSixabcde$
\\ \quad
\\ \mbox{$\FigTwentySixA$: Generalized Dynkin diagrams for \cite[Table~2, Row~6]{Hec09}}
\end{array}
\end{equation*}


\begin{equation*}
\begin{array}{l}
$\spTotalCHHecThreeSix$
\\ \quad
\\ \mbox{$\FigTwentySevenA$: Special Hamilton circuit of \cite[Table~2, Row~6]{Hec09}}
\end{array}
\end{equation*}

\noindent
{\bf{Rank-$3$-Case-$4$:}} Let $\chamchi$ be of $\acchi^{(2,7)}_1$.
We draw Hamilton circuit maps of $\chamGammachi$
by $\FigTwentyNineA$ 
and $\FigThirtyA$. 
Notice that
$\chamchi$ is of \cite[Table~2, Rows~7]{Hec09}.

\begin{equation*}
\begin{array}{l}
$\HecSevenabcd$
\\ \quad
\\ \mbox{$\FigTwentyEightA$: Generalized Dynkin diagrams of \cite[Table~2, Row~7]{Hec09}}
\end{array}
\end{equation*}



\begin{equation*}
\begin{array}{l}
$\sspTotalCHHecThreeSeven$
\\ \quad
\\ \mbox{$\FigTwentyNineA$: $2$-convenient Hamilton circuit of \cite[Table~2, Row~7]{Hec09}}
\end{array}
\end{equation*}

\begin{equation*}
\begin{array}{l}
$\spTotalCHHecThreeSeven$
\\ \quad
\\ \mbox{$\FigThirtyA$: Special Hamilton circuit of \cite[Table~2, Row~7]{Hec09}}
\end{array}
\end{equation*}

\noindent
{\bf{Rank-$3$-Case-$5$:}} Let $\chamchi$ be of $\acchi^{(2,9)}_1$,
$\acchi^{(2,10)}_1$ or $\acchi^{(2,11)}_1$.
We draw a Hamilton circuit map of $\chamGammachi$
by $\FigThirtyTwoA$. 
Notice that
$\chamchi$ is of \cite[Table~2, Rows~9,10,11]{Hec09}.

\begin{equation*}
\begin{array}{l}
$\HecNineabcd$
\\ \quad
\\ \mbox{$\FigThirtyOneA$: Generalized Dynkin diagrams for \cite[Table~2, Rows~9,10,11]{Hec09}}
\end{array}
\end{equation*}


\begin{equation*}
\begin{array}{l}
$\TotalCHHecThreeNine$
\\ \quad
\\ \mbox{$\FigThirtyTwoA$: Special and $3$-convenient Hamilton  circuit}
\\ \mbox{for \cite[Table~2, Rows~9, 10, 11]{Hec09}}
\end{array}
\end{equation*}

\noindent
{\bf{Rank-$3$-Case-$6$:}} Let $\chamchi$ be of $\acchi^{(2,13)}_1$ or $\acchi^{(2,18)}_1$.
We draw a Hamilton circuit map of $\chamGammachi$
by  $\FigThirtyFiveA$. 
Notice that
$\chamchi$ is of \cite[Table~2, Rows~13,18]{Hec09}.

\begin{equation*}
\begin{array}{l}
$\HecThirteenab$
\\ \quad
\\ \mbox{$\FigThirtyThreeA$: Generalized Dynkin diagrams for  \cite[Table~2, Rows~13]{Hec09}}
\end{array}
\end{equation*}

\begin{equation*}
\begin{array}{l}
$\HecEighteenab$
\\ \quad
\\ \mbox{$\FigThirtyFourA$: Generalized Dynkin diagrams for \cite[Table~2, Rows~18]{Hec09}}
\end{array}
\end{equation*}

\begin{equation*}
\begin{array}{l}
    $\TotalCHHecThreeEighteen$
\\ \quad
\\ \mbox{$\FigThirtyFiveA$: Special and $2$-convenient Hamilton circuit of \cite[Table~3, Rows~13,18]{Hec09}}
\end{array}
\end{equation*}

\noindent
{\bf{Rank-$3$-Case-$7$:}} Let $\chamchi$ be of $\acchi^{(2,15)}_1$,
where let $\defzeta\in\bKt$ be such that $\defzeta^2+\defzeta+1=0$.
We draw a Hamilton circuit map of $\chamGammachi$
by $\FigThirtySevenA$. 
Notice that
$\chamchi$ is of \cite[Table~2, Rows~15]{Hec09}.

\begin{equation*}
\begin{array}{l}
$\HecFifteenabcd$
\\ \quad
\\ \mbox{$\FigThirtySixA$: The generalized Dynkin diagrams for  \cite[Table~2, Row15]{Hec09}}
\end{array}
\end{equation*}




\begin{equation*}
\begin{array}{l}
$\spTotalCHHecThreeFifteen$
\\ \quad
\\ \mbox{$\FigThirtySevenA$: Special and $2$-convenient Hamilton circuit of \cite[Table~2, Row~15]{Hec09}}
\end{array}
\end{equation*}

\noindent
{\bf{Rank-$3$-Case-$8$:}} Let $\chamchi$ be of $\acchi^{(2,16)}_1$,
where let $\defzeta\in\bKt$ be such that $\defzeta^2+\defzeta+1=0$.
We draw a Hamilton circuit map of $\chamGammachi$
by $\FigThirtyEightA$. 
Notice that
$\chamchi$ is of \cite[Table~2, Rows~16]{Hec09}.

\begin{equation*}
\begin{array}{l}
$\HeckSixteenabcdefghij$
\\ \quad
\\ \mbox{$\FigThirtyEightA$: The generalized Dynkin diagrams for  \cite[Table~2, Row16]{Hec09}}
\end{array}
\end{equation*}



\begin{equation*}
\begin{array}{l}
    $\TotalCHHecThreeSixteen$
\\ \quad
\\ \mbox{$\FigThirtyNineA$: Hamilton circuit for \cite[Table~2, Row~16]{Hec09}}
\end{array}
\end{equation*}

\noindent
{\bf{Rank-$3$-Case-$9$:}}  
Let $\chamchi$ be of $\acchi^{(2,17)}_1$,
where let $\defzeta\in\bKt$ be such that $\defzeta^2+\defzeta+1=0$.
We draw a Hamilton circuit map of $\chamGammachi$
by $\FigFourtyOneA$. 
Notice that
$\chamchi$ is of \cite[Table~2, Rows~17]{Hec09}.

\begin{equation*}
\begin{array}{l}
$\HecSeventeenabcdefghi$
\\ \quad
\\ \mbox{$\FigFourtyA$: Generalized Dynkin diagrams in \cite[Table~2, Row~17]{Hec09}}
\end{array}
\end{equation*}


\begin{equation*}
\begin{array}{l}
    $\TotalCHHecThreeSeventeen$
\\ \quad
\\ \mbox{$\FigFourtyOneA$: Hamilton Circuit for \cite[Table~2, Row~17]{Hec09}}
\end{array}
\end{equation*}

\section{Hamilton Circuits of Rank-4 Cases}
In this section, analyzing carefully in each case,
we shall see:
\begin{theorem} \label{theorem:existRkFour}
Assume that $|\chamfkI|=4$.
Then for every $\chamchi$ of {\rm{\cite[Table~3]{Hec09}}},
a Hamilton circuit map of 
$\chamGammachi$ exists.
In fact, {\rm{(}}up to elements of $\chambarcalGchi${\rm{)}} it is drawn by 
$\FigEightA$ {\rm{(}}$\acchi^{(3,1)}_\bullet$, $\acchi^{(3,6)}_\bullet$, $\acchi^{(3,10)}_\bullet${\rm{)}},
$\FigFourA$ {\rm{(}}$\acchi^{(3,2)}_\bullet$, $\acchi^{(3,3)}_\bullet$, $\acchi^{(3,7)}_\bullet$, $\acchi^{(3,11)}_\bullet$,
$\acchi^{(3,15)}_\bullet$, $\acchi^{(3,16)}_\bullet$, $\acchi^{(3,19)}_\bullet${\rm{)}},
$\FigTenA$ {\rm{(}}$\acchi^{(3,4)}_\bullet${\rm{)}}, 
$\FigNineA$ {\rm{(}}$\acchi^{(3,5)}_\bullet${\rm{)}}, 
$\FigFourtyFiveA$ {\rm{(}}$\acchi^{(3,8)}_\bullet${\rm{)}}, 
$\FigFourtyEightA$ {\rm{(}}$\acchi^{(3,12)}_\bullet${\rm{)}}, 
$\FigFiftyTwoA$ {\rm{(}}$\acchi^{(3,13)}_\bullet${\rm{)}}, 
$\FigFiftyFiveA$ {\rm{(}}$\acchi^{(3,20)}_\bullet${\rm{)}}, 
$\FigFiftyEightA$ {\rm{(}}$\acchi^{(3,21)}_\bullet${\rm{)}}, 
$\FigSixtyOneA$ {\rm{(}}$\acchi^{(3,9)}_\bullet${\rm{)}}, 
$\FigSixtyFourA$ {\rm{(}}$\acchi^{(3,14)}_\bullet${\rm{)}}, 
$\FigSixtySevenA$ {\rm{(}}$\acchi^{(3,17)}_\bullet${\rm{)}}, 
$\FigSeventyFourA$ {\rm{(}}$\acchi^{(3,18)}_\bullet${\rm{)}} and
$\FigEightyOneA$ {\rm{(}}$\acchi^{(3,22)}_\bullet${\rm{)}}. 
\end{theorem}

\noindent
{\bf{Rank-4-Case-1:}} Let $\chamchi$ be of \cite[Table~3,\,Rows 1-7,10,11,15,16,19]{Hec09}.
Then it is of Cartan-type or of quasi-Caratan-type.
See $\FigFourtyTwoA$. 
If  $\chi$ is of 
$\acchi^{(3,1)}_1$ 
(resp. $\acchi^{(3,2)}_1$, $\acchi^{(3,3)}_1$, $\acchi^{(3,4)}_1$, $\acchi^{(3,5)}_1$), 
then it is of Cartan-type of type 
$A_4$ (resp.~$B_4$, resp.~$C_4$, resp.~$F_4$, resp.~$D_4$),  
a Hamilton circuit of $\chamGammachi$ is
given by 
$\FigEightA$ 
(resp.~$\FigFourA$, resp.~$\FigFourA$,
resp.~$\FigTenA$,  
resp.~$\FigNineA$). 
If  $\chi$ is of 
$\acchi^{(3,1)}_1$ {(resp.~$\acchi^{(3,2)}_1$),
we have $\chi\cong\chi^\prime$
for $\chi^\prime$ of $\acchi^{(3,6)}_1$ and $\acchi^{(3,10)}_1$
(resp.~$\acchi^{(3,3)}_1$, $\acchi^{(3,7)}_1$, $\acchi^{(3,11)}_1$, $\acchi^{(3,15)}_1$, $\acchi^{(3,16)}_1$
and
$\acchi^{(3,19)}_1$),
where $\chi^\prime$ is of quasi-Cartan-type of type 
$A_4$ (resp.~$B_4$).

\begin{equation*}
\begin{array}{l}
$\TopHecFourQUCX$
\\ \quad
\\ \mbox{$\FigFourtyTwoA$: First diagrams of \cite[Table~3,\,Rows 1-7,10,11,15,16,19]{Hec09}}
\end{array}
\end{equation*}

\noindent
{\bf{Rank-4-Case-2:}} Let $\chamchi$ be of \cite[Table~3,\,Rows 8, 12, 13, 20, 21]{Hec09}.
As mentioned by Table~$\tableone$, it is of $(1,2)$-special,
whence a Hamilton circuit map of 
$\chamGammachi$ exists by Proposition~\ref{proposition:maintoolcor}.

\begin{center}
\begingroup
\renewcommand{\arraystretch}{1.5}
Table~$\tableone$: Rank-4  $(1,2)$-special bicharacters \\
  \begin{tabular}{lccccc}
    \hline
$\chamchi$  &   
$\acchi^{(3,8)}_1$ &
$\acchi^{(3,12)}_1$ &
$\acchi^{(3,13)}_1$ &
$\acchi^{(3,20)}_1$ &
$\acchi^{(3,21)}_1$
\\
    \hline \hline
\mbox{Fig.s} & 
\NFigFourtyThreeA, 
\NFigFourtyFourA 
&
\NFigFourtySixA, 
\NFigFourtySevenA 
&
\NFigFourtyNineA, 
\NFigFiftyA, 
&
\NFigFiftyThreeA, 
\NFigFiftyFourA 
&
\NFigFiftySixA, 
\NFigFiftySevenA 
\\
\mbox{Extracted $\al_i$}  &  
$\al_1$ &
$\al_1$ &
$\al_1$ &
$\al_1$ &
$\al_1$
\\
\mbox{Rep.s}  & 
$\{a,b,g \}$ & 
$\{a,c\}$ &
$\{a,e \}$ &
$\{c,{\bar{d}}\}$ &
$\{{\bar{a}},c\}$
\\
\mbox{Fig.s for $\chamchi_{\langle i\rangle}$}  & 
\NFigTwentyFiveA, 
\NFigTwentySevenA, 
\NFigTwentyFiveA, 
&
\NFigTwentySevenA, 
\NFigThirtyTwoA 
&
\NFigThirtyTwoA, 
\NFigTwentyOneA 
&
\NFigThirtySevenA, 
\NFigTwentySevenA 
&
\NFigThirtyFiveA, 
\NFigThirtySevenA 
\\ 
\mbox{Special}  &  
(1,2) & 
(1,2) &
(1,2) &
(1,2) &
(1,2)
\\
\hline
  \end{tabular}
\endgroup
\end{center}

\begin{equation*}
\begin{array}{l}
$\allHecFourEight$
\\ \quad
\\ \mbox{$\FigFourtyThreeA$: Generalized Dynkin diagram of \cite[Table~3, Row~8]{Hec09}}
\end{array}
\end{equation*}

\begin{equation*}
\begin{array}{l}
$\RankFourHecEightTau$
\\ \quad
\\ \mbox{$\FigFourtyFourA$: Changing of diagrams of $\FigFourtyThreeA$} 
\end{array}
\end{equation*}

\begin{equation*}
\begin{array}{l}
${\footnotesize{$\begin{array}{l}
(s^a_3) (s^a_2) s^a_3 s^a_2 s^a_3 
 (s^a_4) s^a_3 s^a_2 s^a_3 s^a_4 s^a_3 s^a_2 
s^a_3 s^a_4 s^a_3 s^a_4 s^a_3 s^a_2 s^a_3 
s^a_4 s^a_3 s^a_4 s^a_3 s^a_4 s^a_3 s^a_2 
s^a_3 s^a_2 s^a_3 s^a_4 s^a_3 s^a_2 s^a_3 \\ \cdot 
s^a_4 s^a_3 s^a_2 s^a_3 s^a_4 s^a_3 s^a_4 \\ \cdot 
s^a_3 s^a_2 s^a_3 s^a_4 s^a_3 s^a_4 s^a_3 
 (s^b_1) (s^c_2) (s^d_3) (s^d_2) (s^e_4)
 (s^e_3) (s^d_4) (s^c_3) (s^c_4) 
(s^b_2) (s^b_4) 
 (s^b_3) s^c_2 s^d_3 s^e_4 s^e_3 s^d_4 s^c_3 
s^b_2 s^b_3 s^b_4 s^c_2 s^c_4 s^d_3 s^e_4 
\\ \cdot
 (s^f_2) (s^f_4) (s^f_3) s^f_4 s^f_3 s^f_4 
 (s^e_2) s^e_3 \\ \cdot s^f_2 s^f_4 s^e_2 s^d_4 
s^d_2 (s^d_1) s^d_2 s^e_4 s^e_3 s^d_4 s^c_3 
s^c_4 s^b_2 s^b_4 s^b_3 s^c_2 s^d_3 s^e_4 
s^e_3 s^d_4 s^c_3 s^b_2 s^b_3 s^b_4 s^c_2 
s^c_4 s^d_3 s^e_4 s^f_2 s^f_4 s^f_3 s^f_4 
s^f_3 s^f_4 s^e_2 s^e_3 s^f_2 s^f_4 s^e_2 
(s^e_1) \\ \cdot s^e_3 s^d_4 s^c_3 s^c_4 s^b_2 
s^b_4 s^b_3 s^c_2 s^d_3 s^e_4 s^e_3 s^d_4 
s^c_3 s^b_2 s^b_3 s^b_4 s^c_2 s^c_4 s^d_3 
s^e_4 s^f_2 s^f_4 s^f_3 s^f_4 s^f_3 s^f_4 
s^e_2 s^e_3 s^f_2 s^f_4 s^e_2 s^e_1 s^e_3 
s^d_4 s^c_3 s^c_4 s^b_2 s^b_4 s^b_3 s^c_2 \\ \cdot 
s^d_3 s^e_4 s^e_3 s^d_4 s^c_3 s^b_2 s^b_3 
s^b_4 s^c_2 s^c_4 s^d_3 s^e_4 s^f_2 s^f_4 
s^f_3 s^f_4 s^f_3 s^f_4 s^e_2 s^e_3 s^f_2 
s^f_4 s^e_2 s^e_1 s^e_3 s^d_4 s^c_3 s^c_4 
s^b_2 s^b_4 s^b_3 s^c_2 s^d_3 s^e_4 s^e_3 
s^d_4 s^c_3 s^b_2 s^b_3 s^b_4 \\ \cdot s^c_2 s^c_4 
s^d_3 s^e_4 s^f_2 s^f_4 s^f_3 s^f_4 s^f_3 
s^f_4 s^e_2 s^e_3 s^f_2 (s^g_1) (s^g_3) 
(s^g_4) 
(s^g_2) s^g_4 s^g_2 s^g_4 s^g_3 s^g_4 s^g_2 
s^g_4 s^g_3 s^g_4 s^g_2 s^g_4 s^g_3 s^g_4 \\ \cdot 
s^g_3 s^g_4 s^g_2 s^g_4 s^g_3 s^g_4 s^g_3 
s^g_4 s^g_3 s^g_4 \\ \cdot s^g_2 
(s^f_1) s^e_2 s^d_4 s^d_2 s^c_3 s^b_2 s^b_4 
s^c_2 s^d_3 s^d_2 s^e_4 s^e_3 s^d_4 s^c_3 
s^c_4 s^b_2 s^b_4 s^b_3 s^c_2 s^d_3 s^e_4 
s^e_3 s^d_4 s^c_3 s^b_2 s^b_3 s^b_4 s^c_2 
s^c_4 s^d_3 s^e_4 s^f_2 s^f_4 s^f_3 s^f_4 
s^f_3 s^f_4 s^e_2 s^e_3 \\ \cdot s^f_2 s^g_1 s^g_2 
s^g_4 s^g_3 s^g_4 s^g_2 s^g_4 s^g_3 s^g_4 
s^g_2 s^g_4 s^g_3 s^g_4 s^g_3 s^g_4 s^g_2 
s^g_4 s^g_3 s^g_4 s^g_3 s^f_1 s^e_2 s^d_4 
s^d_2 s^c_3 s^b_2 s^b_4 s^c_2 s^d_3 s^d_2 
s^d_1 s^d_2 s^c_3 s^b_2 s^b_4 s^c_2 s^d_3 
s^d_2 s^d_1 \\ \cdot s^d_2 s^c_3 s^b_2 s^b_4 s^c_2 
s^d_3 s^d_2 s^d_1 s^d_2 s^c_3 s^b_2 s^b_4 
s^c_2 (s^c_1) s^b_2 
(s^a_1) 
\end{array}$}}$
\\ \quad
\\ \mbox{$\FigFourtyFiveA$: Special Hamilton circuit of \cite[Table~3, Row~8]{Hec09},} \\
\mbox{where $(s^x_i)$'s mean that it is really special
(Length $=336$)
}
\end{array}
\end{equation*}

\begin{equation*}
\begin{array}{l}
$\allHecFourTwelve$
\\ \quad
\\ \mbox{$\FigFourtySixA$: Generalized Dynkin diagram of \cite[Table~3, Row~12]{Hec09}}
\end{array}
\end{equation*}

\begin{equation*}
\begin{array}{l}
$\RankFourHecTwelveTau$
\\ \quad
\\ \mbox{$\FigFourtySevenA$: Changing of diagrams of $\FigFourtySixA$} 
\end{array}
\end{equation*}

\begin{equation*}
\begin{array}{l}
${\footnotesize{$\begin{array}{l}
(s^f_2) (s^f_4) (s^f_3) s^f_4 (s^e_2) 
 (s^g_4) (s^g_2) (s^e_4) (s^c_3) (s^c_4)
  (s^c_2) s^c_4 (s^e_3) s^g_4 
 (s^g_3) s^e_4 s^f_2 s^f_4 s^f_3 s^f_4 s^e_2 
s^g_4 s^g_2 s^e_4 s^c_3 s^c_4 s^c_2 s^c_4 
s^e_3 s^g_4 s^g_3 \\ \cdot
 (s^h_1) (s^h_3) (s^d_4) 
 (s^b_3) (s^a_2) (s^a_3) (s^a_4) (s^b_2)
  (s^b_4) \\ \cdot (s^d_3) (s^h_4) s^h_3 s^d_4 
 (s^d_2) s^b_3 s^a_2 s^a_4 s^b_2 s^d_3 s^d_2 
s^h_4 (s^i_2) (s^i_4) 
 (s^h_2) s^h_3 s^i_2 s^i_4 
 (s^i_3) s^i_4 s^i_3 s^i_4 s^h_2 s^d_4 s^b_3 
s^b_4 s^a_2 s^a_4 
 \\ \cdot  (s^a_1) s^a_4 s^b_2 s^b_4 s^d_3 s^h_4 s^h_3 
s^d_4 s^d_2 s^b_3 s^a_2 s^a_4 \\ \cdot s^b_2 s^d_3 
s^d_2 s^h_4 s^i_2 s^i_4 s^h_2 s^h_3 s^i_2 
s^i_4 s^i_3 s^i_4 s^i_3 s^i_4 s^h_2 s^d_4 
s^b_3 (s^c_1) s^e_3 s^g_4 s^g_3 s^e_4 s^f_2 
s^f_4 s^f_3 s^f_4 s^e_2 s^g_4 s^g_2 s^e_4 
s^c_3 s^c_4 s^c_2 s^c_4 s^e_3 s^g_4 s^g_3 
s^e_4 s^f_2 s^f_4 \\ \cdot s^f_3 s^f_4 s^e_2 
 (s^d_1) s^b_3 s^a_2 s^a_3 s^a_4 s^b_2 s^b_4 
s^d_3 s^h_4 s^h_3 s^d_4 s^d_2 s^b_3 s^a_2 
s^a_4 s^b_2 s^d_3 s^d_2 s^h_4 s^i_2 s^i_4 
s^h_2 s^h_3 s^i_2 s^i_4 s^i_3 s^i_4 s^i_3 
s^i_4 s^h_2 s^d_4 s^b_3 s^b_4 s^a_2 s^a_4 
s^a_1 s^a_4 \\ \cdot s^b_2 s^b_4 s^d_3 s^h_4 s^h_3 
s^d_4 s^d_2 s^b_3 s^a_2 s^a_4 s^b_2 s^d_3 
s^d_2 s^h_4 s^i_2 s^i_4 s^h_2 
 (s^g_1) s^e_4 s^f_2 s^f_4 s^f_3 s^f_4 s^e_2 
s^g_4 s^g_2 s^e_4 s^c_3 s^c_4 s^c_2 s^c_4 
s^e_3 s^g_4 s^g_3 s^e_4 s^f_2 s^f_4 s^f_3 
 (s^f_1) s^e_2 \\ \cdot s^g_4 s^g_2 s^e_4 s^c_3 
s^c_4 s^c_2 s^c_4 s^e_3 s^g_4 s^g_3 s^e_4 
s^f_2 s^f_4 s^f_3 s^f_4 s^e_2 s^g_4 s^g_2 
s^e_4 s^c_3 s^c_4 s^c_2 s^c_4 s^e_3 s^g_4 
s^g_3 s^e_4 s^f_2 s^f_4 s^f_3 s^f_1 s^e_2 
s^g_4 s^g_2 s^e_4 s^c_3 s^c_4 s^c_2 s^c_4 
s^e_3 \\ \cdot s^g_4 s^h_1 s^i_2 s^i_4 s^i_3 s^i_4 
s^i_3 s^i_4 s^h_2 s^d_4 s^b_3 s^b_4 s^a_2 
s^a_4 s^a_3 s^b_2 s^d_3 s^h_4 s^h_3 s^d_4 
s^b_3 s^a_2 s^a_1 s^b_2 s^d_3 s^h_4 s^h_3 
s^g_1 s^g_2 s^e_4 s^c_3 s^c_4 s^c_2 
 (s^b_1) s^a_2 s^a_4 s^a_3 s^b_2 s^d_3 s^h_4 
\\ \cdot s^h_3 s^d_4 s^b_3 s^a_2 s^a_1 s^b_2 s^d_3 
 (s^e_1) 
\end{array}
$}}$
\\ \quad
\\ \mbox{$\FigFourtyEightA$: Hamilton circuit of \cite[Table~3, Row~12]{Hec09} missing $s^i_1$, }
\\ \mbox{where $(s^x_i)$'s mean that $s^u_v$ with $(u,v)\ne (i,1)$ exist
(Length $=288$)}
\end{array}
\end{equation*}

\begin{equation*}
\begin{array}{l}
$\allHecFourThirteen$
\\ \quad
\\ \mbox{$\FigFourtyNineA$: Generalized Dynkin diagram of \cite[Table~3, Row~13]{Hec09}}
\end{array}
\end{equation*}

\begin{equation*}
\begin{array}{l}
$\RankForHecFourThirteenTau$
\\ \quad
\\ \mbox{$\FigFiftyA$: Changing of diagrams of $\FigFourtyNineA$} 
\end{array}
\end{equation*}

\begin{equation*}
\begin{array}{l}
\setlength{\unitlength}{1mm}
\begin{picture}(220,60)(-10,00)
\put(00,0){$\TotalCHHecThreeNineForHecFourThirteen$}
\put(70,0){$\TotalCHHecThreeEightForHecFourThirteen$}
\end{picture}
\\ \quad
\\ \mbox{$\FigFiftyOneA$: Parts of Hamilton circuit of \cite[Table~3, Row~13]{Hec09}} 
\\ \mbox{being $2$-convenient and special}
\end{array}
\end{equation*}


\begin{equation*}
\begin{array}{l}
{\footnotesize{\mbox{$\begin{array}{l}
(s^a_2) (s^a_4) (s^b_3) (s^c_4) (s^c_2)
  (s^b_4) (s^d_2) (s^d_4) (s^d_3) s^d_4 
 (s^b_2) s^c_4 (s^c_3) s^b_4 
 (s^a_3) s^a_4 s^a_2 s^a_4 s^b_3 s^c_4 s^c_2 
s^b_4 s^d_2 s^d_4 s^d_3 s^d_4 s^b_2 s^c_4 
s^c_3 s^b_4 s^a_3 \\ \cdot 
 (s^a_1) s^a_2 s^a_4 s^b_3 s^c_4 s^c_2 s^b_4 
s^d_2 s^d_4 \\ \cdot s^d_3 s^d_4 s^b_2 s^c_4 s^c_3 
s^b_4 s^a_3 s^a_4 s^a_2 s^a_4 s^b_3 s^c_4 
s^c_2 s^b_4 s^d_2 s^d_4 s^d_3 (s^e_1) 
 (s^e_2) (s^e_3) s^e_2 
 (s^e_4) s^e_2 s^e_4 s^e_2 s^e_3 s^e_2 s^e_3 
s^e_2 s^e_4 s^e_2 s^e_3 s^e_2 s^e_4 s^e_2 
s^e_4 s^e_2  \\ \cdot 
 (s^d_1) s^d_4 s^b_2 \\ \cdot s^c_4 s^c_3 s^b_4 
s^a_3 s^a_4 s^a_2 s^a_4 s^b_3 s^c_4 s^c_2 
s^b_4 s^d_2 s^d_4 s^d_3 s^d_4 s^b_2 s^c_4 
s^c_3 s^b_4 s^a_3 s^a_4 s^a_2 s^a_4 s^b_3 
s^c_4 s^c_2 
 (s^c_1) s^c_3 s^b_4 s^a_3 s^a_4 s^a_2 s^a_4 
s^b_3 s^c_4 s^c_2 s^b_4 s^d_2 s^d_4 s^d_3 \\ \cdot 
s^d_4 s^b_2 s^c_4 s^c_3 s^b_4 s^a_3 s^a_4 
s^a_2 s^a_4 s^b_3 s^c_4 s^c_2 s^c_1 s^c_3 
s^b_4 s^a_3 s^a_4 s^a_2 s^a_4 s^b_3 s^c_4 
s^c_2 s^b_4 s^d_2 s^d_4 s^d_3 s^d_4 s^b_2 
s^c_4 s^c_3 s^b_4 s^a_3 s^a_1 s^a_2 s^a_4 
s^b_3 s^c_4 s^c_2 s^b_4 s^d_2 \\ \cdot s^d_4 s^d_3 
s^e_1 s^e_2 s^e_3 s^e_2 s^e_4 s^e_2 s^e_4 
s^e_2 s^e_3 s^e_2 s^e_3 s^e_2 s^e_4 s^e_2 
s^e_3 s^e_2 s^e_4 s^e_2 s^e_4 s^e_2 s^e_3 
s^e_2 s^e_3 s^e_2 s^d_1 s^b_2 s^c_4 s^c_3 
s^b_4 s^a_3 s^a_4 s^a_2 s^a_4 s^b_3 s^c_4 
s^c_2 s^b_4 s^d_2 \\ \cdot s^d_4 s^d_3 s^d_4 s^b_2 
s^c_4 s^c_3 s^b_4 s^a_3 s^a_1 s^a_2 s^a_4 
s^b_3 s^c_4 s^c_2 s^b_4 s^d_2 s^d_4 s^d_3 
s^d_4 s^b_2 
 (s^b_1) s^d_2 s^d_4 s^d_3 s^d_4 s^b_2 s^b_1 
s^d_2 s^d_4 s^e_1 s^e_2 s^e_3 s^e_2 s^d_1 
s^b_2 s^c_4 s^c_3 s^b_4 s^a_3 s^a_1 
\end{array}
$}}}
\\ \quad
\\ \mbox{$\FigFiftyTwoA$: Special Hamilton circuit of \cite[Table~3, Row~13]{Hec09} by 
$\FigFiftyOneA$}, \\ 
\mbox{where $(s^x_i)$'s mean that it is really special
(Length $=240$)
}
\end{array}
\end{equation*}

\begin{equation*}
\begin{array}{l}
$\allHecFourTewnty$
\\ \quad
\\ \mbox{$\FigFiftyThreeA$: Generalized Dynkin diagrams of \cite[Table~3, Row~20]{Hec09}}
\end{array}
\end{equation*}

\begin{equation*}
\begin{array}{l}
$\RankFourHecTwentyTau$
\\ \quad
\\ \mbox{$\FigFiftyFourA$: Changing of diagrams of $\FigFiftyThreeA$}
\end{array}
\end{equation*}

\begin{equation*}
\begin{array}{l}
{\tiny{\mbox{$
\begingroup
\renewcommand{\arraystretch}{1.5}
\begin{array}{l}
(s^e_3) (s^d_2) (s^d_3) (s^e_2) s^e_3 
 (s^c_4) (s^c_3) 
 (s^e_4) s^e_3 s^d_2 s^d_3 s^e_2 s^e_3 s^c_4 
s^c_3 s^e_4 s^e_3 s^c_4 s^c_3 (s^b_2) 
 (s^a_3) (s^a_4) (s^b_3) 
 (s^c_2) s^c_3 s^b_2 s^a_3 s^a_4 s^b_3 s^c_2 
s^c_3 s^b_2 s^a_3 s^a_4 s^b_3 
 (s^d_4) s^d_3 s^e_2 s^e_3 s^c_4 s^c_3
\\ \cdot 
 s^b_2 
s^a_3 (s^a_2) s^b_3 s^d_4 s^d_3 
 (s^b_4) s^a_3 s^a_2 \\ \cdot s^b_3 s^d_4 s^d_3 
s^b_4 s^a_3 s^a_2 s^b_3 s^d_4 (s^{\bar{e}}_1) 
 (s^{\bar{d}}_4) (s^{\bar{c}}_3) (s^{\bar{b}}_2) (s^{\bar{b}}_3) 
 (s^{\bar{a}}_4) (s^{\bar{a}}_3) (s^{\bar{a}}_2) (s^{\bar{b}}_4)
 (s^{\bar{c}}_2) (s^{\bar{d}}_3) 
 (s^{\bar{d}}_2) s^{\bar{c}}_3 s^{\bar{b}}_2 s^{\bar{a}}_4 s^{\bar{a}}_2 s^{\bar{a}}_3 
s^{\bar{b}}_4 s^{\bar{b}}_3 s^{\bar{c}}_2 s^{\bar{d}}_3 s^{\bar{d}}_2 s^{\bar{c}}_3 
 (s^{\bar{c}}_4) s^{\bar{b}}_2 s^{\bar{a}}_4 s^{\bar{a}}_3 s^{\bar{b}}_4 s^{\bar{c}}_2 
s^{\bar{c}}_4 s^{\bar{d}}_3 
\\ \cdot 
(s^{\bar{e}}_4) 
 (s^{\bar{e}}_3) s^{\bar{d}}_4 s^{\bar{d}}_2 s^{\bar{e}}_4 s^{\bar{e}}_3 (s^{\bar{e}}_2)
  (s^d_1) s^e_2 s^e_3 s^d_2 \\ \cdot s^d_3 s^e_2 
s^e_3 s^c_4 s^c_3 s^e_4 s^e_3 s^d_2 s^d_3 
s^e_2 s^e_3 s^c_4 s^c_3 s^e_4 s^e_3 s^c_4 
s^c_3 s^b_2 s^a_3 s^a_4 s^b_3 s^c_2 s^c_3 
s^b_2 s^a_3 s^a_4 s^b_3 s^c_2 s^c_3 s^b_2 
s^a_3 s^a_4 s^b_3 s^d_4 s^d_3 s^e_2 s^e_3 
s^c_4 s^c_3 s^b_2 s^a_3 s^a_2 s^b_3 s^d_4 
s^d_3 s^b_4 s^a_3 s^a_2 s^b_3 s^d_4 \\ \cdot s^d_3 
s^b_4 (s^{\bar{d}}_1) s^{\bar{d}}_2 s^{\bar{c}}_3 s^{\bar{b}}_2 s^{\bar{a}}_4 
s^{\bar{a}}_2 s^{\bar{a}}_3 s^{\bar{b}}_4 s^{\bar{b}}_3 s^{\bar{c}}_2 s^{\bar{d}}_3 s^{\bar{d}}_2 
s^{\bar{c}}_3 s^{\bar{c}}_4 s^{\bar{b}}_2 s^{\bar{a}}_4 s^{\bar{a}}_3 s^{\bar{b}}_4 s^{\bar{c}}_2 
s^{\bar{c}}_4 s^{\bar{d}}_3 s^{\bar{e}}_4 s^{\bar{e}}_3 s^{\bar{d}}_4 s^{\bar{d}}_2 s^{\bar{e}}_4 
s^{\bar{e}}_3 s^{\bar{e}}_2 s^{\bar{e}}_3 s^{\bar{e}}_2 s^{\bar{e}}_3 s^{\bar{d}}_4 s^{\bar{c}}_3 
s^{\bar{b}}_2 s^{\bar{b}}_3 s^{\bar{a}}_4 
 (s^{\bar{a}}_1) s^{\bar{b}}_4 s^{\bar{b}}_3 s^{\bar{c}}_2 s^{\bar{d}}_3 s^{\bar{d}}_2 
s^{\bar{c}}_3 s^{\bar{c}}_4 s^{\bar{b}}_2 s^{\bar{a}}_4 s^{\bar{a}}_3 s^{\bar{b}}_4 \\ \cdot 
s^{\bar{c}}_2 s^{\bar{c}}_4 s^{\bar{d}}_3 s^{\bar{e}}_4 s^{\bar{e}}_3 s^{\bar{d}}_4 s^{\bar{d}}_2 
s^{\bar{e}}_4 s^{\bar{e}}_3 s^{\bar{e}}_2 s^{\bar{e}}_3 s^{\bar{e}}_2 s^{\bar{e}}_3 s^{\bar{d}}_4 
s^{\bar{c}}_3 s^{\bar{b}}_2 s^{\bar{b}}_3 s^{\bar{a}}_4 s^{\bar{a}}_3 s^{\bar{a}}_2 s^{\bar{b}}_4 
s^{\bar{c}}_2 
 (s^a_1) s^a_2 s^b_3 s^d_4 s^d_3 s^e_2 s^e_3 
s^d_2 s^d_3 s^e_2 s^e_3 s^c_4 s^c_3 s^e_4 
s^e_3 s^d_2 s^d_3 s^e_2 s^e_3 s^c_4 s^c_3 
s^e_4 s^e_3 s^c_4 s^c_3 s^b_2 s^a_3 s^a_4 \\ \cdot 
s^b_3 s^c_2 s^c_3 s^b_2 s^a_3 s^a_4 s^b_3 
s^c_2 s^c_3 s^b_2 s^a_3 s^a_4 s^b_3 s^d_4 
s^d_3 s^e_2 s^e_3 s^c_4 s^c_3 s^b_2 s^a_3 
s^a_2 s^b_3 s^d_4 s^d_3 s^b_4 s^{\bar{d}}_1 s^{\bar{d}}_2 
s^{\bar{c}}_3 s^{\bar{b}}_2 s^{\bar{a}}_4 s^{\bar{a}}_2 s^{\bar{a}}_3 s^{\bar{b}}_4 s^{\bar{b}}_3 
s^{\bar{c}}_2 s^{\bar{d}}_3 s^{\bar{d}}_2 s^{\bar{c}}_3 s^{\bar{c}}_4 s^{\bar{b}}_2 s^{\bar{a}}_4 
s^{\bar{a}}_3 s^{\bar{b}}_4 s^{\bar{c}}_2 s^{\bar{c}}_4 s^{\bar{d}}_3 s^{\bar{e}}_4 s^{\bar{e}}_3 
s^{\bar{d}}_4 \\ \cdot s^{\bar{d}}_2 s^{\bar{e}}_4 s^{\bar{e}}_3 s^{\bar{e}}_2 s^d_1 
s^e_2 s^e_3 s^d_2 s^d_3 s^e_2 s^e_3 s^c_4 
s^c_3 s^e_4 s^e_3 s^d_2 s^d_3 s^e_2 s^e_3 
s^c_4 s^c_3 s^e_4 s^e_3 s^c_4 s^c_3 s^b_2 
s^a_3 s^a_4 s^b_3 s^c_2 s^c_3 s^b_2 s^a_3 
s^a_4 s^b_3 s^c_2 s^c_3 s^b_2 s^a_3 s^a_4 
s^b_3 s^d_4 s^d_3 s^e_2 s^e_3 s^c_4 s^c_3 
s^b_2 s^a_3 s^a_2 \\ \cdot s^b_3 s^d_4 s^d_3 s^b_4 
s^{\bar{d}}_1 s^{\bar{d}}_2 s^{\bar{c}}_3 s^{\bar{b}}_2 s^{\bar{a}}_4 s^{\bar{a}}_2 s^{\bar{a}}_3 
s^{\bar{b}}_4 s^{\bar{b}}_3 s^{\bar{c}}_2 s^{\bar{d}}_3 s^{\bar{d}}_2 s^{\bar{c}}_3 s^{\bar{c}}_4 
s^{\bar{b}}_2 s^{\bar{a}}_4 s^{\bar{a}}_3 s^{\bar{b}}_4 s^{\bar{c}}_2 s^{\bar{c}}_4 s^{\bar{d}}_3 
s^{\bar{e}}_4 s^{\bar{e}}_3 s^{\bar{d}}_4 s^{\bar{d}}_2 s^{\bar{e}}_4 s^{\bar{e}}_3 s^{\bar{e}}_2 
s^d_1 s^e_2 s^e_3 s^d_2 s^d_3 s^e_2 s^e_3 
s^c_4 s^c_3 s^e_4 s^e_3 s^d_2 s^d_3 s^e_2 
s^e_3 s^c_4 s^c_3 s^e_4 \\ \cdot s^e_3 s^c_4 s^c_3 
s^b_2 s^a_3 s^a_4 s^b_3 s^c_2 s^c_3 s^b_2 
s^a_3 s^a_4 s^b_3 s^c_2 s^c_3 s^b_2 s^a_3 
s^a_4 s^b_3 s^d_4 s^d_3 s^e_2 s^e_3 s^c_4 
 (s^c_1) s^b_2 s^a_3 s^a_2 s^b_3 s^d_4 s^d_3 
s^b_4 s^a_3 s^a_2 s^b_3 s^d_4 s^d_3 s^b_4 
s^a_3 s^a_2 s^b_3 s^d_4 s^d_3 s^e_2 s^e_3 
s^d_2 s^d_3 s^e_2 s^e_3 s^c_4 \\ \cdot s^c_3 s^e_4 
s^e_3 s^d_2 s^d_3 s^e_2 s^e_3 s^c_4 s^c_3 
s^e_4 s^e_3 s^c_4 s^c_3 s^b_2 s^a_3 s^a_4 
s^b_3 s^c_2 s^c_3 s^b_2 s^a_3 s^a_4 s^b_3 
s^c_2 s^c_3 s^b_2 s^a_3 s^a_4 s^b_3 s^d_4 
s^{\bar{e}}_1 s^{\bar{d}}_4 s^{\bar{c}}_3 s^{\bar{b}}_2 s^{\bar{b}}_3 s^{\bar{a}}_4 s^{\bar{a}}_3 
s^{\bar{a}}_2 s^{\bar{b}}_4 s^{\bar{c}}_2 s^{\bar{d}}_3 s^{\bar{d}}_2 s^{\bar{c}}_3 s^{\bar{b}}_2 
s^{\bar{a}}_4 s^{\bar{a}}_2 s^{\bar{a}}_3 s^{\bar{b}}_4 s^{\bar{b}}_3 s^{\bar{c}}_2 \\ \cdot 
s^{\bar{d}}_3 s^{\bar{d}}_2 s^{\bar{c}}_3 s^{\bar{c}}_4 s^{\bar{b}}_2 s^{\bar{a}}_4 s^{\bar{a}}_3 
s^{\bar{b}}_4 s^{\bar{c}}_2 s^{\bar{c}}_4 s^{\bar{d}}_3 s^{\bar{e}}_4 s^{\bar{e}}_3 s^d_1 
s^d_3 s^e_2 s^e_3 s^d_2 s^d_3 s^e_2 s^e_3 
s^c_4 s^c_3 s^e_4 s^e_3 s^d_2 s^d_3 s^e_2 
s^e_3 s^c_4 s^c_3 s^e_4 s^e_3 s^c_4 s^c_3 
s^b_2 s^a_3 s^a_4 s^b_3 s^c_2 s^c_3 s^b_2 
s^a_3 s^a_4 s^b_3 s^c_2 s^c_3 s^b_2 s^a_3 
s^a_4 \\ \cdot s^b_3 s^d_4 s^d_3 s^e_2 s^e_3 s^c_4 
s^c_3 s^b_2 s^a_3 s^a_2 s^b_3 s^d_4 s^d_3 
s^b_4 s^a_3 s^a_2 s^b_3 s^d_4 s^d_3 s^b_4 
s^{\bar{d}}_1 s^{\bar{d}}_2 s^{\bar{c}}_3 s^{\bar{b}}_2 s^{\bar{a}}_4 s^{\bar{a}}_2 s^{\bar{a}}_3 
s^{\bar{b}}_4 s^{\bar{b}}_3 s^{\bar{c}}_2 s^{\bar{d}}_3 s^{\bar{d}}_2 s^{\bar{c}}_3 s^{\bar{c}}_4 
s^{\bar{b}}_2 s^{\bar{a}}_4 s^{\bar{a}}_3 s^{\bar{b}}_4 s^{\bar{c}}_2 s^{\bar{c}}_4 s^{\bar{d}}_3 
s^{\bar{e}}_4 s^{\bar{e}}_3 s^{\bar{d}}_4 s^{\bar{d}}_2 s^{\bar{e}}_4 s^{\bar{e}}_3 s^{\bar{e}}_2 
s^{\bar{e}}_3 s^{\bar{e}}_2 \\ \cdot s^{\bar{e}}_3 s^{\bar{d}}_4 s^{\bar{c}}_3 s^{\bar{b}}_2 
s^{\bar{b}}_3 s^{\bar{a}}_4 s^{\bar{a}}_1 s^{\bar{b}}_4 s^{\bar{b}}_3 s^{\bar{c}}_2 s^{\bar{d}}_3 
s^{\bar{d}}_2 s^{\bar{c}}_3 s^{\bar{c}}_4 s^{\bar{b}}_2 s^{\bar{a}}_4 s^{\bar{a}}_3 s^{\bar{b}}_4 
s^{\bar{c}}_2 s^{\bar{c}}_4 s^{\bar{d}}_3 s^{\bar{e}}_4 s^{\bar{e}}_3 s^{\bar{d}}_4 s^{\bar{d}}_2 
s^{\bar{e}}_4 s^{\bar{e}}_3 s^{\bar{e}}_2 s^{\bar{e}}_3 s^{\bar{e}}_2 s^{\bar{e}}_3 s^{\bar{d}}_4 
s^{\bar{c}}_3 s^{\bar{b}}_2 s^{\bar{b}}_3 s^{\bar{a}}_4 s^{\bar{a}}_1 s^{\bar{b}}_4 s^{\bar{b}}_3 
s^{\bar{c}}_2 s^{\bar{d}}_3 s^{\bar{d}}_2 s^{\bar{c}}_3 s^{\bar{c}}_4 s^{\bar{b}}_2 s^{\bar{a}}_4 
s^{\bar{a}}_3 s^{\bar{b}}_4 s^{\bar{c}}_2 s^{\bar{c}}_4 \\ \cdot s^{\bar{d}}_3 s^{\bar{e}}_4 
s^{\bar{e}}_3 s^{\bar{d}}_4 s^{\bar{d}}_2 s^{\bar{e}}_4 s^{\bar{e}}_3 s^{\bar{e}}_2 s^{\bar{e}}_3 
s^{\bar{e}}_2 s^{\bar{e}}_3 s^{\bar{d}}_4 s^{\bar{c}}_3 s^{\bar{b}}_2 s^{\bar{b}}_3 s^{\bar{a}}_4 
s^{\bar{a}}_1 s^{\bar{b}}_4 s^{\bar{b}}_3 s^{\bar{c}}_2 s^{\bar{d}}_3 s^{\bar{d}}_2 s^{\bar{c}}_3 
s^{\bar{c}}_4 s^{\bar{b}}_2 s^{\bar{a}}_4 s^{\bar{a}}_3 s^{\bar{a}}_1 s^{\bar{a}}_2 s^{\bar{a}}_3 
s^{\bar{b}}_4 s^{\bar{b}}_3 s^{\bar{c}}_2 s^{\bar{d}}_3 s^{\bar{d}}_2 s^{\bar{c}}_3 s^{\bar{c}}_4 
s^{\bar{b}}_2 s^{\bar{a}}_4 s^{\bar{a}}_3 s^{\bar{b}}_4 s^{\bar{c}}_2 s^{\bar{c}}_4 s^{\bar{d}}_3 
s^{\bar{e}}_4 s^{\bar{e}}_3 s^{\bar{d}}_4 s^{\bar{d}}_2 s^{\bar{e}}_4 s^d_1 \\ \cdot 
s^e_2 s^e_3 s^d_2 s^d_3 s^e_2 s^e_3 s^c_4 
s^c_3 s^e_4 s^e_3 s^d_2 s^d_3 s^e_2 s^e_3 
s^c_4 s^c_3 s^e_4 s^e_3 s^c_4 s^c_3 s^b_2 
s^a_3 s^a_4 s^b_3 s^c_2 s^c_3 s^b_2 s^a_3 
s^a_4 s^b_3 s^c_2 s^c_3 s^b_2 s^a_3 s^a_4 
s^b_3 s^d_4 s^{\bar{e}}_1 s^{\bar{d}}_4 s^{\bar{c}}_3 s^{\bar{b}}_2 s^{\bar{b}}_3 
s^{\bar{a}}_4 s^{\bar{a}}_3 s^{\bar{a}}_2 s^{\bar{b}}_4 s^{\bar{c}}_2 s^{\bar{d}}_3 s^{\bar{d}}_2 
s^{\bar{c}}_3 \\ \cdot s^{\bar{b}}_2 s^{\bar{a}}_4 s^{\bar{a}}_2 s^{\bar{a}}_3 s^{\bar{b}}_4 
s^{\bar{b}}_3 s^{\bar{c}}_2 s^{\bar{d}}_3 s^{\bar{d}}_2 s^{\bar{c}}_3 s^{\bar{c}}_4 s^{\bar{b}}_2 
s^{\bar{a}}_4 s^{\bar{a}}_3 s^{\bar{b}}_4 s^{\bar{c}}_2 s^{\bar{c}}_4 s^{\bar{d}}_3 s^{\bar{e}}_4 
s^{\bar{e}}_3 s^{\bar{d}}_4 s^{\bar{d}}_2 s^{\bar{e}}_4 s^{\bar{e}}_3 s^{\bar{e}}_2 s^{\bar{e}}_3 
s^{\bar{e}}_2 s^d_1 s^e_2 s^e_3 s^c_4 s^c_3 s^b_2 
s^a_3 s^a_2 s^b_3 s^d_4 s^d_3 s^b_4 s^a_3 
s^a_2 s^b_3 s^d_4 s^d_3 s^b_4 s^a_3 s^a_2 
s^b_3 s^d_4 s^{\bar{e}}_1 \\ \cdot s^{\bar{e}}_2 s^{\bar{e}}_3 s^{\bar{e}}_2 
s^{\bar{e}}_3 s^{\bar{d}}_4 s^{\bar{c}}_3 s^{\bar{b}}_2 s^{\bar{b}}_3 s^{\bar{a}}_4 s^{\bar{a}}_3 
s^{\bar{a}}_2 s^{\bar{b}}_4 s^{\bar{c}}_2 s^{\bar{d}}_3 s^{\bar{d}}_2 s^{\bar{c}}_3 s^{\bar{b}}_2 
 (s^{\bar{b}}_1) s^{\bar{c}}_2 s^{\bar{c}}_4 s^{\bar{d}}_3 s^{\bar{e}}_4 s^{\bar{e}}_3 
s^{\bar{d}}_4 s^{\bar{d}}_2 s^{\bar{e}}_4 s^{\bar{e}}_3 s^{\bar{e}}_2 s^{\bar{e}}_3 s^{\bar{e}}_2 
s^{\bar{e}}_3 s^{\bar{d}}_4 s^{\bar{c}}_3 s^{\bar{b}}_2 s^{\bar{b}}_3 s^{\bar{a}}_4 s^{\bar{a}}_3 
s^{\bar{a}}_2 s^{\bar{b}}_4 s^{\bar{c}}_2 s^{\bar{d}}_3 s^{\bar{d}}_2 s^{\bar{c}}_3 s^{\bar{b}}_2 
s^{\bar{a}}_4 s^{\bar{a}}_2 s^{\bar{a}}_1 s^{\bar{a}}_2 s^{\bar{b}}_4 s^{\bar{c}}_2 \\ \cdot 
s^{\bar{d}}_3 s^{\bar{d}}_2 s^{\bar{c}}_3 s^{\bar{b}}_2 s^{\bar{a}}_4 s^{\bar{a}}_2 s^{\bar{a}}_1 
s^{\bar{a}}_2 s^{\bar{b}}_4 s^{\bar{c}}_2 s^{\bar{d}}_3 s^{\bar{d}}_2 s^{\bar{c}}_3 s^{\bar{b}}_2 
s^{\bar{a}}_4 s^{\bar{a}}_2 s^{\bar{a}}_1 s^{\bar{a}}_2 s^{\bar{b}}_4 s^{\bar{c}}_2 s^a_1 
s^a_2 s^b_3 s^{\bar{d}}_1 s^{\bar{d}}_2 s^{\bar{e}}_4 s^{\bar{e}}_3 s^{\bar{e}}_2 
s^{\bar{e}}_3 s^{\bar{e}}_2 s^d_1 s^e_2 s^e_3 s^c_4 s^c_1 
s^b_2 s^a_3 s^a_2 s^b_3 s^d_4 s^d_3 s^b_4 
s^a_3 s^a_2 s^b_3 s^d_4 s^d_3 s^b_4 s^a_3 
s^a_2 \\ \cdot s^b_3 s^d_4 s^{\bar{e}}_1 s^{\bar{e}}_2 s^{\bar{e}}_3 
s^{\bar{d}}_4 s^{\bar{c}}_3 s^{\bar{b}}_2 s^{\bar{b}}_3 s^{\bar{a}}_4 s^{\bar{a}}_3 s^{\bar{a}}_2 
s^{\bar{b}}_4 s^{\bar{c}}_2 s^a_1 s^a_2 s^b_3 s^d_4 s^d_3 
s^b_4 s^a_3 s^a_2 s^b_3 s^d_4 s^{\bar{e}}_1 s^{\bar{e}}_2 
s^{\bar{e}}_3 s^{\bar{d}}_4 s^{\bar{c}}_3 s^{\bar{b}}_2 s^{\bar{b}}_3 s^{\bar{a}}_4 s^{\bar{a}}_3 
s^{\bar{a}}_2 s^{\bar{b}}_4 s^{\bar{c}}_2 s^a_1 s^a_2 s^b_3 s^d_4 
s^d_3 s^b_4 s^{\bar{d}}_1 s^{\bar{d}}_2 s^{\bar{c}}_3 s^{\bar{b}}_2 s^{\bar{a}}_4 
s^{\bar{a}}_2 s^{\bar{a}}_1 s^{\bar{a}}_2 \\ \cdot s^{\bar{b}}_4 s^{\bar{c}}_2 s^a_1 
s^a_2 s^b_3 s^d_4 s^{\bar{e}}_1 s^{\bar{e}}_2 s^d_1 s^e_2 
\end{array}
\endgroup$}}}
\\ \quad
\\ \mbox{$\FigFiftyFiveA$: Hamilton circuit of \cite[Table~3, Row~20]{Hec09} 
missing $s^b_1$, $s^e_1$ and $s^{\bar c}_1$} 
\\ 
\mbox{(In this circuit, there exists $s^x_i$ for every
$(x,i)\in\{a,\ldots,{\bar{e}}\}\times\fkJ_{2,4}$.)}
\\
\mbox{(Length $=960$)}
\end{array}
\end{equation*}

\begin{equation*}
\begin{array}{l}
$\allHecFourTwentyone$
\\ \quad
\\ \mbox{$\FigFiftySixA$: Generalized Dynkin diagram of \cite[Table~3, Row~21]{Hec09}}
\end{array}
\end{equation*}

\begin{equation*}
\begin{array}{l}
$\RankFourHecTwentyoneTau$
\\ \quad
\\ \mbox{$\FigFiftySevenA$: Changing of diagrams of $\FigFiftySixA$} 
\end{array}
\end{equation*}

\begin{equation*}
\begin{array}{l}
{\tiny{\mbox{$
\begingroup
\renewcommand{\arraystretch}{1.5}
\begin{array}{l}
S^\prime:=(s^y_3) (s^x_4) (s^x_3) (s^y_4) s^y_3 
 (s^y_2) s^y_3 s^y_2 s^y_3 s^x_4 s^x_3 
 (s^x_2) s^x_3 s^x_2 s^x_3 s^y_4 s^y_3 s^x_4 
s^x_3 s^x_2 s^x_3 s^y_4 s^y_3 s^x_4 s^x_3 
s^y_4 s^y_3 s^y_2 s^y_3 s^y_2 s^y_3 s^x_4 
s^x_3 s^x_2 s^x_3 s^x_2 s^x_3 s^y_4 s^y_3 
s^x_4 s^x_3 s^x_2 s^x_3 \\ \cdot s^y_4 s^y_3 s^x_4 
s^x_3 s^y_4 s^y_3 s^y_2 
\\ \cdot s^y_3 s^y_2 s^y_3 
s^x_4 s^x_3 s^x_2 s^x_3 s^x_2 s^x_3 s^y_4 
s^y_3 s^x_4 s^x_3 s^x_2 s^x_3 s^y_4 s^y_3 
s^x_4 s^x_3 s^y_4 s^y_3 s^y_2 s^y_3 s^y_2 
s^y_3 s^x_4 s^x_3 s^y_4 s^y_3 s^y_2 s^y_3 
s^y_2 s^y_3 s^y_2 s^y_3 s^x_4 s^x_3 s^x_2 
s^x_3 s^y_4 s^y_3 s^x_4 s^x_3 s^x_2 s^x_3 
\\ \cdot  (s^c_1) (s^c_3) (s^e_4) (s^e_3) (s^d_2)
\\ \cdot (s^d_3) (s^e_2) s^e_3 s^d_2 s^d_3 
 (s^b_4) (s^a_3) (s^a_2) (s^b_3) 
 (s^d_4) s^d_3 s^b_4 s^a_3 s^a_2 s^b_3 s^d_4 
s^d_3 s^b_4 s^a_3 s^a_2 s^b_3 
 (s^c_2) s^c_3 s^e_4 s^e_3 s^d_2 s^d_3 s^b_4 
s^a_3 (s^a_4) s^b_3 s^c_2 s^c_3 
 (s^b_2) s^a_3 s^a_4 s^b_3 s^c_2 s^c_3 s^b_2 
s^a_3 \\ \cdot s^a_4 s^b_3 s^c_2 s^c_3 s^e_4 s^e_3 
 (s^c_4) s^c_3 s^e_4 
\\ \cdot s^e_3 s^d_2 
 (s^d_1) s^b_4 s^a_3 s^a_2 s^b_3 s^c_2 s^c_3 
s^e_4 s^e_3 s^d_2 s^d_3 s^b_4 s^a_3 s^a_4 
s^b_3 s^c_2 s^c_3 s^b_2 s^a_3 s^a_4 s^b_3 
s^c_2 s^c_3 s^b_2 s^a_3 s^a_4 s^b_3 s^c_2 
s^c_3 s^e_4 s^e_3 s^c_4 s^c_3 s^e_4 s^e_3 
s^d_2 s^d_3 s^e_2 s^e_3 s^c_4 s^c_3 s^e_4 
s^e_3 s^d_2 s^d_3 s^e_2 s^e_3 s^d_2 \\ \cdot s^d_3 
s^b_4 s^a_3 s^a_2 s^b_3 s^d_4 s^d_3 s^b_4 
 (s^b_1) s^c_2 s^c_3 s^e_4 s^e_3 s^d_2 s^d_3 
s^b_4 s^a_3 s^a_4 s^b_3 s^c_2 s^c_3 s^b_2 
s^a_3 s^a_4 s^b_3 s^c_2 s^c_3 s^b_2 s^a_3 
s^a_4 s^b_3 s^c_2 s^c_3 s^e_4 s^e_3 s^c_4 
s^c_3 s^e_4 s^e_3 s^d_2 s^d_3 s^e_2 s^e_3 
s^c_4 s^c_3 s^e_4 s^e_3 s^d_2 s^d_3 s^e_2 \\ \cdot 
s^e_3 s^d_2 s^d_3 s^b_4 s^a_3 s^a_2 s^b_3 
s^d_4 s^d_3 s^b_4 s^b_1 s^c_2 s^c_3 s^e_4 
s^e_3 s^d_2 s^d_3 s^b_4 s^a_3 s^a_4 s^b_3 
s^c_2 s^c_3 s^b_2 s^a_3 s^a_4 s^b_3 s^c_2 
s^c_3 s^b_2 s^a_3 s^a_4 s^b_3 s^c_2 s^c_3 
s^e_4 s^e_3 s^c_4 s^c_3 s^e_4 s^e_3 s^d_2 
s^d_3 s^e_2 s^e_3 s^c_4 s^c_3 s^e_4 s^e_3 
s^d_2 \\ \cdot s^d_3 s^e_2 
 (s^y_1) s^x_4 s^x_3 s^x_2 s^x_3 s^y_4 s^y_3 
s^x_4 s^x_3 s^y_4 s^y_3 s^y_2 s^y_3 s^y_2 
s^y_3 s^x_4 s^x_3 s^x_2 s^x_3 s^x_2 s^x_3 
s^y_4 s^y_3 s^x_4 s^x_3 s^x_2 s^x_3 s^y_4 
s^y_3 s^x_4 s^x_3 s^y_4 s^y_3 s^y_2 s^y_3 
s^y_2 s^y_3 s^x_4 s^x_3 s^x_2 s^x_3 s^x_2 \\ \cdot
s^x_3 s^y_4 s^y_3 s^x_4 s^x_3 s^x_2 \\ \cdot s^x_3 
s^y_4 s^y_3 s^x_4 s^x_3 s^y_4 s^y_3 s^y_2 
s^y_3 s^y_2 s^y_3 s^x_4 s^x_3 s^x_2 s^x_3 
s^x_2 s^x_3 s^y_4 s^y_3 s^x_4 s^x_3 s^x_2 
s^x_3 s^y_4 s^y_3 s^x_4 s^x_3 s^y_4 s^y_3 
s^y_2 s^y_3 s^y_2 s^y_3 s^x_4 s^x_3 s^y_4 
s^y_3 s^y_2 s^y_3 s^y_2 
 (s^e_1) s^d_2 s^d_3 s^b_4 s^a_3 s^a_2 s^b_3 
s^d_4 s^d_3 s^b_4 \\ \cdot s^a_3 s^a_2 s^b_3 s^d_4 
s^d_3 s^b_4 s^a_3 s^a_2 s^b_3 s^c_2 s^c_3 
s^e_4 s^e_3 s^d_2 s^d_3 s^b_4 s^a_3 s^a_4 
s^b_3 s^c_2 s^c_3 s^b_2 s^a_3 s^a_4 s^b_3 
s^c_2 s^c_3 s^b_2 s^a_3 s^a_4 s^b_3 s^c_2 
s^c_3 s^e_4 s^e_3 s^c_4 s^c_3 s^e_4 s^e_3 
s^d_2 s^d_1 s^b_4 s^a_3 s^a_2 s^b_3 s^c_2 
s^c_3 s^e_4 s^e_3 s^d_2 \\ \cdot s^d_3 s^b_4 s^a_3 
s^a_4 s^b_3 s^c_2 s^c_3 s^b_2 s^a_3 s^a_4 
s^b_3 s^c_2 s^c_3 s^b_2 s^a_3 s^a_4 s^b_3 
s^c_2 s^c_3 s^e_4 s^e_3 s^c_4 s^c_3 s^e_4 
s^e_3 s^d_2 s^d_3 s^e_2 s^e_3 s^c_4 s^c_3 
s^e_4 s^e_3 s^d_2 s^d_3 s^e_2 s^e_3 s^d_2 
s^d_3 s^b_4 s^b_1 s^c_2 s^c_3 s^e_4 s^e_3 
s^d_2 s^d_3 s^b_4 s^a_3 s^a_4 \\ \cdot s^b_3 s^c_2 
s^c_3 s^b_2 s^a_3 s^a_4 s^b_3 s^c_2 s^c_3 
s^b_2 s^a_3 s^a_4 s^b_3 s^c_2 s^c_3 s^e_4 
s^e_3 s^c_4 s^c_3 s^e_4 s^e_3 s^d_2 s^d_3 
s^e_2 s^e_3 s^c_4 s^c_3 s^e_4 s^e_3 s^d_2 
s^d_3 s^e_2 s^y_1 s^x_4 s^x_3 s^x_2 s^x_3 
s^y_4 s^y_3 s^x_4 s^x_3 s^y_4 s^y_3 s^y_2 
s^y_3 s^y_2 s^y_3 s^x_4 s^x_3 s^x_2 \\ \cdot s^x_3 
s^x_2 s^x_3 s^y_4 s^y_3 s^x_4 s^x_3 s^x_2 
s^x_3 s^y_4 s^y_3 s^x_4 s^x_3 s^y_4 s^y_3 
s^y_2 s^y_3 s^y_2 s^y_3 s^x_4 s^x_3 s^x_2 
s^x_3 s^x_2 s^x_3 s^y_4 s^y_3 s^x_4 s^x_3 
s^x_2 s^x_3 s^y_4 s^y_3 s^x_4 s^x_3 s^y_4 
s^y_3 s^y_2 s^y_3 s^y_2 s^y_3 s^x_4 s^x_3 
s^x_2 s^x_3 s^x_2 s^x_3 s^y_4 s^y_3 s^x_4 \\ \cdot 
s^x_3 s^x_2 s^x_3 s^y_4 s^y_3 s^x_4 s^x_3 
s^y_4 s^y_3 s^y_2 s^y_3 s^y_2 s^y_3 s^x_4 
s^x_3 s^y_4 s^y_3 s^y_2 s^e_1 s^d_2 s^d_3 
s^b_4 s^a_3 s^a_2 s^b_3 s^d_4 s^d_3 s^b_4 
s^a_3 s^a_2 s^b_3 s^d_4 s^d_3 s^b_4 s^a_3 
s^a_2 s^b_3 s^c_2 s^c_3 s^e_4 s^e_3 s^d_2 
s^d_3 s^b_4 s^a_3 s^a_4 s^b_3 s^c_2 s^c_3 
s^b_2 \\ \cdot s^a_3 s^a_4 s^b_3 s^c_2 s^c_3 s^b_2 
s^a_3 s^a_4 s^b_3 s^c_2 s^c_3 s^e_4 s^e_3 
s^c_4 s^c_3 s^e_4 s^e_3 s^d_2 s^d_1 s^b_4 
s^a_3 s^a_2 s^b_3 s^c_2 s^c_3 s^e_4 s^e_3 
s^d_2 s^d_3 s^b_4 s^a_3 s^a_4 s^b_3 s^c_2 
s^c_3 s^b_2 s^a_3 s^a_4 s^b_3 s^c_2 s^c_3 
s^b_2 s^a_3 s^a_4 s^b_3 s^c_2 s^c_3 s^e_4 
s^e_3 s^c_4 \\ \cdot s^c_3 s^e_4 s^e_3 s^d_2 s^d_3 
s^e_2 s^e_3 s^c_4 s^c_3 s^e_4 s^e_3 s^d_2 
s^d_3 s^e_2 s^e_3 s^d_2 s^d_3 s^b_4 s^b_1 
s^c_2 s^c_3 s^e_4 s^e_3 s^d_2 s^d_3 s^b_4 
s^a_3 s^a_4 s^b_3 s^c_2 s^c_3 s^b_2 s^a_3 
s^a_4 s^b_3 s^c_2 s^c_3 s^b_2 s^a_3 s^a_4 
s^b_3 s^c_2 s^c_3 s^e_4 s^e_3 s^c_4 s^c_3 
s^e_4 s^e_3 s^d_2 \\ \cdot s^d_3 s^e_2 s^e_3 s^c_4 
s^c_3 s^e_4 s^e_3 s^d_2 s^d_3 s^e_2 s^y_1 
s^x_4 s^x_3 s^x_2 s^x_3 s^y_4 s^y_3 s^x_4 
s^x_3 s^y_4 s^y_3 s^y_2 s^y_3 s^y_2 s^y_3 
s^x_4 s^x_3 s^x_2 s^x_3 s^x_2 s^x_3 s^y_4 
s^y_3 s^x_4 s^x_3 s^x_2 s^x_3 s^y_4 s^y_3 
s^x_4 s^x_3 s^y_4 s^y_3 s^y_2 s^y_3 s^y_2 
s^y_3 s^x_4 s^x_3 s^x_2 \\ \cdot s^x_3 s^x_2 s^x_3 
s^y_4 s^y_3 s^x_4 s^x_3 s^x_2 s^x_3 s^y_4 
s^y_3 s^x_4 s^x_3 s^y_4 s^y_3 s^y_2 s^y_3 
s^y_2 s^y_3 s^x_4 s^x_3 s^x_2 s^x_3 s^x_2 
s^x_3 s^y_4 s^y_3 s^x_4 s^x_3 s^x_2 s^x_3 
s^y_4 s^y_3 s^x_4 s^x_3 s^y_4 s^y_3 s^y_2 
s^y_3 s^y_2 s^y_3 s^x_4 s^x_3 s^y_4 s^y_3 
s^y_2 s^e_1 s^d_2 s^d_3 s^b_4 \\ \cdot s^a_3 s^a_2 
s^b_3 s^d_4 s^d_3 s^b_4 s^a_3 s^a_2 s^b_3 
s^d_4 s^d_3 s^b_4 s^a_3 s^a_2 s^b_3 s^c_2 
s^c_3 s^e_4 s^e_3 s^d_2 s^d_3 s^b_4 s^a_3 
s^a_4 s^b_3 s^c_2 s^c_3 s^b_2 s^a_3 s^a_4 
s^b_3 s^c_2 s^c_3 s^b_2 s^a_3 s^a_4 s^b_3 
s^c_2 s^c_3 s^e_4 s^e_3 s^c_4 s^c_3 s^e_4 
s^e_3 s^d_2 s^d_1 s^b_4 s^a_3 s^a_2 \\ \cdot s^b_3 
s^c_2 s^c_3 s^e_4 s^e_3 s^d_2 s^d_3 s^b_4 
s^a_3 s^a_4 s^b_3 s^c_2 s^c_3 s^b_2 s^a_3 
s^a_4 s^b_3 s^c_2 s^c_3 s^b_2 s^a_3 s^a_4 
s^b_3 s^c_2 s^c_3 s^e_4 s^e_3 s^c_4 s^c_3 
s^e_4 s^e_3 s^d_2 s^d_3 s^e_2 s^e_3 s^c_4 
s^c_3 s^e_4 s^e_3 s^d_2 s^d_3 s^e_2 s^e_3 
s^d_2 s^d_3 s^b_4 s^b_1 s^c_2 s^c_3 s^e_4 \\ \cdot 
s^e_3 s^d_2 s^d_3 s^b_4 s^a_3 s^a_4 s^b_3 
s^c_2 s^c_3 s^b_2 s^a_3 s^a_4 s^b_3 s^c_2 
s^c_3 s^b_2 s^a_3 s^a_4 s^b_3 s^c_2 s^c_3 
s^e_4 s^e_3 s^c_4 s^c_3 s^e_4 s^e_3 s^d_2 
s^d_3 s^e_2 s^e_3 s^c_4 s^c_3 s^e_4 s^e_3 
s^d_2 s^d_1 s^b_4 s^a_3 s^a_2 s^b_3 s^c_2 
s^c_3 s^e_4 s^e_3 s^d_2 s^d_3 s^b_4 s^a_3 
s^a_4 \\ \cdot s^b_3 s^c_2 s^c_3 s^b_2 s^a_3 s^a_4 
s^b_3 s^c_2 s^c_3 s^b_2 s^a_3 s^a_4 s^b_3 
s^c_2 s^c_3 s^e_4 s^e_3 s^c_4 s^c_3 s^e_4 
s^e_3 s^d_2 s^d_3 s^e_2 s^e_3 s^c_4 s^c_3 
s^e_4 s^e_3 s^d_2 s^d_3 s^e_2 s^e_3 s^d_2 
s^d_3 s^b_4 s^a_3 s^a_2 s^b_3 s^d_4 s^d_3 
s^b_4 s^b_1 s^c_2 s^c_3 s^e_4 s^e_3 s^d_2 
s^d_3 s^b_4 
\end{array}
\endgroup
$}}}
\\ \quad
\\ \mbox{$\FigFiftyEightA^\prime$} 
\end{array}
\end{equation*}

\begin{equation*}
\begin{array}{l}
{\tiny{\mbox{$
\begingroup
\renewcommand{\arraystretch}{1.5}
\begin{array}{l}
S^\prime \cdot s^a_3 s^a_4 s^b_3 s^c_2 s^c_3 
s^b_2 s^a_3 s^a_4 s^b_3 s^c_2 s^c_3 s^b_2 
s^a_3 s^a_4 s^b_3 s^c_2 s^c_3 s^e_4 s^e_3 
s^c_4 s^c_3 s^e_4 s^e_3 s^d_2 s^d_3 s^e_2 
s^e_3 s^c_4 s^c_3 s^e_4 s^e_3 s^d_2 s^d_3 
s^e_2 s^e_3 s^d_2 s^d_3 s^b_4 s^a_3 s^a_2 
s^b_3 s^d_4 s^d_3 s^b_4 s^b_1 s^c_2 s^c_3 
s^e_4 s^e_3 s^d_2 \\ \cdot s^d_3 s^b_4 s^a_3 s^a_4 
s^b_3 s^c_2 s^c_3 s^b_2 s^a_3 s^a_4 s^b_3 
s^c_2 s^c_3 s^b_2 s^a_3 s^a_4 s^b_3 s^c_2 
s^c_3 s^e_4 s^e_3 s^c_4 s^c_3 s^e_4 s^e_3 
s^d_2 s^d_3 s^e_2 s^e_3 s^c_4 s^c_3 s^e_4 
s^e_3 s^d_2 s^d_1 s^b_4 s^a_3 s^a_2 s^b_3 
s^c_2 s^c_3 s^e_4 s^e_3 s^d_2 s^d_3 s^b_4 
s^a_3 s^a_4 s^b_3 s^c_2 \\ \cdot s^c_3 s^b_2 s^a_3 
s^a_4 s^b_3 s^c_2 s^c_3 s^b_2 s^a_3 s^a_4 
s^b_3 s^c_2 s^c_3 s^e_4 s^e_3 s^c_4 s^c_3 
s^e_4 s^e_3 s^d_2 s^d_3 s^e_2 s^e_3 s^c_4 
s^c_3 s^e_4 s^e_3 s^d_2 s^d_3 s^e_2 s^e_3 
s^d_2 s^d_3 s^b_4 s^a_3 s^a_2 s^b_3 s^d_4 
s^d_3 s^b_4 s^b_1 s^c_2 s^c_3 s^e_4 s^e_3 
s^d_2 s^d_3 s^b_4 s^a_3 s^a_4 \\ \cdot s^b_3 s^c_2 
s^c_3 s^b_2 s^a_3 s^a_4 s^b_3 s^c_2 s^c_3 
s^b_2 s^a_3 s^a_4 s^b_3 s^c_2 s^c_3 s^e_4 
s^e_3 s^c_4 s^c_3 s^e_4 s^e_3 s^d_2 s^d_3 
s^e_2 s^e_3 s^c_4 s^c_3 s^e_4 s^e_3 s^d_2 
s^d_3 s^e_2 s^e_3 s^d_2 s^d_3 s^b_4 s^a_3 
s^a_2 s^b_3 s^d_4 s^d_3 s^b_4 s^b_1 s^c_2 
s^c_3 s^e_4 s^e_3 s^d_2 s^d_3 s^b_4 \\ \cdot s^a_3 
s^a_4 s^b_3 s^c_2 s^c_3 s^b_2 s^a_3 s^a_4 
s^b_3 s^c_2 s^c_3 s^b_2 s^a_3 s^a_4 s^b_3 
s^c_2 s^c_3 s^e_4 s^e_3 s^c_4 s^c_3 s^e_4 
s^e_3 s^d_2 s^d_3 s^e_2 s^e_3 s^c_4 s^c_3 
s^e_4 s^e_3 s^d_2 s^d_3 s^e_2 s^y_1 s^x_4 
s^x_3 s^x_2 s^x_3 s^y_4 s^y_3 s^x_4 s^x_3 
s^y_4 s^y_3 s^y_2 s^y_3 s^y_2 s^y_3 s^x_4 \\ \cdot 
s^x_3 s^x_2 s^x_3 s^x_2 s^x_3 s^y_4 s^y_3 
s^x_4 s^x_3 s^x_2 s^x_3 s^y_4 s^y_3 s^x_4 
s^x_3 s^y_4 s^y_3 s^y_2 s^y_3 s^y_2 s^y_3 
s^x_4 s^x_3 s^x_2 s^x_3 s^x_2 s^x_3 s^y_4 
s^y_3 s^x_4 s^x_3 s^x_2 s^x_3 s^y_4 s^y_3 
s^x_4 s^x_3 s^y_4 s^y_3 s^y_2 s^y_3 s^y_2 
s^y_3 s^x_4 s^x_3 s^x_2 s^x_3 s^x_2 s^x_3 
s^y_4 \\ \cdot s^y_3 s^x_4 s^x_3 s^x_2 s^x_3 s^y_4 
s^y_3 s^x_4 s^x_3 s^y_4 s^y_3 s^y_2 s^e_1 
s^d_2 s^d_3 s^b_4 s^a_3 s^a_2 s^b_3 s^d_4 
s^d_3 s^b_4 s^a_3 s^a_2 s^b_3 s^d_4 s^d_3 
s^b_4 s^a_3 s^a_2 s^b_3 s^c_2 s^c_3 s^e_4 
s^e_3 s^d_2 s^d_3 s^b_4 s^a_3 s^a_4 s^b_3 
s^c_2 s^c_3 s^b_2 s^a_3 s^a_4 s^b_3 s^c_2 
s^c_3 s^b_2 \\ \cdot s^a_3 s^a_4 s^b_3 s^c_2 s^c_3 
s^e_4 s^e_3 s^c_4 s^c_3 s^e_4 s^e_3 s^d_2 
s^d_3 s^e_2 s^e_3 s^c_4 s^c_3 s^e_4 s^y_1 
s^x_4 s^x_3 s^x_2 s^x_3 s^y_4 s^y_3 s^x_4 
s^x_3 s^y_4 s^y_3 s^y_2 s^y_3 s^y_2 s^y_3 
s^x_4 s^x_3 s^x_2 s^x_3 s^x_2 s^x_3 s^y_4 
s^y_3 s^x_4 s^x_3 s^x_2 s^x_3 s^y_4 s^y_3 
s^x_4 s^x_3 s^y_4 \\ \cdot s^y_3 s^y_2 s^y_3 s^y_2 
s^y_3 s^x_4 s^x_3 s^x_2 s^x_3 s^x_2 s^x_3 
s^y_4 s^y_3 s^x_4 s^x_3 s^x_2 s^x_3 s^y_4 
s^y_3 s^x_4 s^x_3 s^y_4 s^y_3 s^y_2 s^y_3 
s^y_2 s^y_3 s^x_4 s^x_3 s^x_2 s^x_3 s^x_2 
s^x_3 s^y_4 s^y_3 s^x_4 s^x_3 s^x_2 s^x_3 
s^y_4 s^y_3 s^x_4 s^x_3 s^y_4 s^e_1 s^d_2 
s^d_3 s^b_4 s^a_3 s^a_2 \\ \cdot s^b_3 s^d_4 s^d_3 
s^b_4 s^a_3 s^a_2 s^b_3 s^d_4 s^d_3 s^b_4 
s^a_3 s^a_2 s^b_3 s^c_2 s^c_3 s^e_4 s^e_3 
s^d_2 s^d_3 s^b_4 s^a_3 s^a_4 s^b_3 s^c_2 
s^c_3 s^b_2 s^a_3 s^a_4 s^b_3 s^c_2 s^c_3 
s^b_2 s^a_3 s^a_4 s^b_3 s^c_2 s^c_3 s^e_4 
s^e_3 s^c_4 s^c_3 s^e_4 s^e_3 s^d_2 s^d_1 
s^b_4 s^a_3 s^a_2 s^b_3 s^c_2 \\ \cdot s^c_3 s^e_4 
s^e_3 s^d_2 s^d_3 s^b_4 s^a_3 s^a_4 s^b_3 
s^c_2 s^c_3 s^b_2 s^a_3 s^a_4 s^b_3 s^c_2 
s^c_3 s^b_2 s^a_3 s^a_4 s^b_3 s^c_2 s^c_3 
s^e_4 s^e_3 s^c_4 s^c_3 s^e_4 s^e_3 s^d_2 
s^d_3 s^e_2 s^e_3 s^c_4 s^c_3 s^e_4 s^e_3 
s^d_2 s^d_3 s^e_2 s^e_3 s^d_2 s^d_3 s^b_4 
s^b_1 s^c_2 s^c_3 s^e_4 s^e_3 s^d_2 \\ \cdot s^d_3 
s^b_4 s^a_3 s^a_4 s^b_3 s^c_2 s^c_3 s^b_2 
s^a_3 s^a_4 s^b_3 s^c_2 s^c_3 s^b_2 s^a_3 
s^a_4 s^b_3 s^c_2 s^c_3 s^e_4 s^e_3 s^c_4 
s^c_3 s^e_4 s^e_3 s^d_2 s^d_3 s^e_2 s^e_3 
s^c_4 s^c_3 s^e_4 s^e_3 s^d_2 s^d_3 s^e_2 
s^y_1 s^x_4 s^x_3 s^x_2 s^x_3 s^y_4 s^y_3 
s^x_4 s^x_3 s^y_4 s^y_3 s^y_2 s^y_3 s^y_2 \\ \cdot 
s^y_3 s^x_4 s^x_3 s^x_2 s^x_3 s^x_2 s^x_3 
s^y_4 s^y_3 s^x_4 s^x_3 s^x_2 s^x_3 s^y_4 
s^y_3 s^x_4 s^x_3 s^y_4 s^y_3 s^y_2 s^y_3 
s^y_2 s^y_3 s^x_4 s^x_3 s^x_2 s^x_3 s^x_2 
s^x_3 s^y_4 s^y_3 s^x_4 s^x_3 s^x_2 s^x_3 
s^y_4 s^y_3 s^x_4 s^x_3 s^y_4 s^y_3 s^y_2 
s^y_3 s^y_2 s^y_3 s^x_4 s^x_3 s^x_2 s^x_3 
s^x_2 \\ \cdot s^x_3 s^y_4 s^y_3 s^x_4 s^x_3 s^x_2 
s^x_3 s^y_4 s^y_3 s^x_4 s^x_3 s^y_4 s^y_3 
s^y_2 s^e_1 s^d_2 s^d_3 s^b_4 s^a_3 s^a_2 
s^b_3 s^d_4 s^d_3 s^b_4 s^a_3 s^a_2 s^b_3 
s^d_4 s^d_3 s^b_4 s^a_3 s^a_2 s^b_3 s^c_2 
s^c_3 s^e_4 s^e_3 s^d_2 s^d_3 s^b_4 s^a_3 
s^a_4 s^b_3 s^c_2 s^c_3 s^b_2 s^a_3 s^a_4 
s^b_3 s^c_2 \\ \cdot s^c_3 s^b_2 s^a_3 s^a_4 s^b_3 
s^c_2 s^c_3 s^e_4 s^e_3 s^c_4 s^c_3 s^e_4 
s^e_3 s^d_2 s^d_3 s^e_2 s^e_3 s^c_4 s^c_3 
s^e_4 s^y_1 s^x_4 s^x_3 s^x_2 s^x_3 s^y_4 
s^y_3 s^x_4 s^x_3 s^y_4 s^y_3 s^y_2 s^y_3 
s^y_2 s^y_3 s^x_4 s^x_3 s^x_2 s^x_3 s^x_2 
s^x_3 s^y_4 s^y_3 s^x_4 s^x_3 s^x_2 s^x_3 
s^y_4 s^y_3 s^x_4 \\ \cdot s^x_3 s^y_4 s^y_3 s^y_2 
s^y_3 s^y_2 s^y_3 s^x_4 s^x_3 s^x_2 s^x_3 
s^x_2 s^x_3 s^y_4 s^y_3 s^x_4 s^x_3 s^x_2 
s^x_3 s^y_4 s^y_3 s^x_4 s^x_3 s^y_4 s^y_3 
s^y_2 s^y_3 s^y_2 s^y_3 s^x_4 s^x_3 s^x_2 
s^x_3 s^x_2 s^x_3 s^y_4 s^y_3 s^x_4 s^x_3 
s^x_2 s^x_3 s^y_4 s^y_3 s^x_4 s^x_3 s^y_4 
s^y_3 s^y_2 s^y_3 s^y_2 \\ \cdot s^y_3 s^x_4 s^x_3 
s^y_4 s^y_3 s^y_2 s^y_3 s^y_2 s^y_3 s^y_2 
s^y_3 s^x_4 s^c_1 s^e_4 s^e_3 s^d_2 s^d_3 
s^e_2 s^e_3 s^d_2 s^d_3 s^b_4 s^a_3 s^a_2 
s^b_3 s^d_4 s^d_3 s^b_4 s^a_3 s^a_2 s^b_3 
s^d_4 s^d_3 s^b_4 s^a_3 s^a_2 s^b_3 s^c_2 
s^c_3 s^e_4 s^e_3 s^d_2 s^d_3 s^b_4 s^a_3 
s^a_4 s^b_3 s^c_2 s^c_3 s^b_2 \\ \cdot s^a_3 s^a_4 
  (s^a_1) s^a_2 s^b_3 s^c_2 s^c_3 s^e_4 s^e_3 
s^d_2 s^d_3 s^b_4 s^a_3 s^a_4 s^b_3 s^c_2 
s^c_3 s^b_2 s^a_3 s^a_4 s^b_3 s^c_2 s^c_3 
s^b_2 s^a_3 s^a_4 s^b_3 s^c_2 s^c_3 s^e_4 
s^e_3 s^c_4 s^c_3 s^e_4 s^e_3 s^d_2 s^d_3 
s^e_2 s^e_3 s^c_4 s^c_3 s^e_4 s^e_3 s^d_2 
s^d_1 s^b_4 s^a_3 s^a_2 s^b_3 s^c_2 \\ \cdot s^c_3 
s^e_4 s^e_3 s^d_2 s^d_3 s^b_4 s^a_3 s^a_4 
s^b_3 s^c_2 s^c_3 s^b_2 s^a_3 s^a_4 s^b_3 
s^c_2 s^c_3 s^b_2 s^a_3 s^a_4 s^b_3 s^c_2 
s^c_3 s^e_4 s^e_3 s^c_4 s^c_3 s^e_4 s^e_3 
s^d_2 s^d_1 s^b_4 s^a_3 s^a_2 s^b_3 s^c_2 
s^c_3 s^e_4 s^e_3 s^d_2 s^d_3 s^b_4 s^a_3 
s^a_4 s^b_3 s^c_2 s^c_3 s^b_2 s^a_3 s^a_4 \\ \cdot 
s^b_3 s^c_2 s^c_3 s^b_2 s^a_3 s^a_4 s^b_3 
s^c_2 s^c_3 s^e_4 s^e_3 s^c_4 s^c_3 s^e_4 
s^e_3 s^d_2 s^d_3 s^e_2 s^e_3 s^c_4 s^c_3 
s^e_4 s^e_3 s^d_2 s^d_1 s^b_4 s^a_3 s^a_2 
s^b_3 s^c_2 s^c_3 s^e_4 s^e_3 s^d_2 s^d_3 
s^b_4 s^a_3 s^a_4 s^b_3 s^c_2 s^c_3 s^b_2 
s^a_3 s^a_4 s^b_3 s^c_2 s^c_3 s^b_2 s^a_3 
s^a_4 \\ \cdot s^b_3 s^c_2 s^c_3 s^e_4 s^e_3 s^c_4 
s^c_3 s^e_4 s^e_3 s^d_2 s^d_3 s^e_2 s^e_3 
s^c_4 s^c_3 s^e_4 s^e_3 s^d_2 s^d_3 s^e_2 
s^e_3 s^d_2 s^d_3 s^b_4 s^b_1 s^c_2 s^c_3 
s^e_4 s^e_3 s^d_2 s^d_3 s^b_4 s^a_3 s^a_4 
s^b_3 s^c_2 s^c_3 s^b_2 s^a_3 s^a_4 s^b_3 
s^c_2 s^c_3 s^b_2 s^a_3 s^a_4 s^b_3 s^c_2 
s^c_3 s^e_4 \\ \cdot s^e_3 s^c_4 s^c_3 s^e_4 s^e_3 
s^d_2 s^d_3 s^e_2 s^e_3 s^c_4 s^c_3 s^e_4 
s^e_3 s^d_2 s^d_1 s^b_4 s^a_3 s^a_2 s^b_3 
s^c_2 s^c_3 s^e_4 s^e_3 s^d_2 s^d_3 s^b_4 
s^a_3 s^a_4 s^a_1 s^a_2 s^b_3 s^c_2 s^c_3 
s^e_4 s^e_3 s^d_2 s^d_3 s^b_4 s^a_3 s^a_4 
s^b_3 s^c_2 s^c_3 s^b_2 s^a_3 s^a_4 s^b_3 
s^c_2 s^c_3 s^b_2 \\ \cdot s^a_3 s^a_4 s^b_3 s^c_2 
s^c_3 s^e_4 s^e_3 s^c_4 s^c_3 s^e_4 s^e_3 
s^d_2 s^d_3 s^e_2 s^e_3 s^c_4 s^c_3 s^e_4 
s^e_3 s^d_2 s^d_3 s^e_2 s^e_3 s^d_2 s^d_3 
s^b_4 s^a_3 s^a_2 s^b_3 s^d_4 s^d_3 s^b_4 
s^a_3 s^a_2 s^b_3 s^d_4 s^d_3 s^b_4 s^b_1 
s^c_2 s^c_3 s^b_2 s^a_3 s^a_4 s^b_3 s^c_2 
s^c_3 s^b_2 s^a_3 s^a_4 \\ \cdot s^b_3 s^c_2 s^c_3 
s^e_4 s^e_3 s^c_4 s^c_3 s^e_4 s^e_3 s^d_2 
s^d_3 s^e_2 s^e_3 s^c_4 s^c_3 s^e_4 s^e_3 
s^d_2 s^d_3 s^e_2 s^e_3 s^d_2 s^d_3 s^b_4 
s^a_3 s^a_2 s^b_3 s^d_4 s^d_3 s^b_4 s^a_3 
s^a_2 s^b_3 s^d_4 s^d_1 s^e_2 s^e_3 s^d_2 
s^d_3 s^b_4 s^a_3 s^a_2 s^b_3 s^d_4 s^d_3 
s^b_4 s^a_3 s^a_2 s^b_3 s^d_4 \\ \cdot s^d_3 s^b_4 
s^a_3 s^a_2 s^a_1 s^a_2 s^b_3 s^d_4 s^d_3 
s^b_4 s^a_3 s^a_2 s^b_3 s^d_4 s^d_1 s^e_2 
s^e_3 s^d_2 s^d_3 s^b_4 s^a_3 s^a_2 s^b_3 
s^d_4 s^d_3 s^b_4 s^a_3 s^a_2 s^b_3 s^d_4 
s^d_1 s^e_2 s^e_3 s^c_4 s^c_3 s^e_4 s^e_3 
s^d_2 s^d_3 s^e_2 s^e_3 s^d_2 s^d_3 s^b_4 
s^a_3 s^a_2 s^b_3 s^d_4 s^d_3 s^b_4 \\ \cdot s^a_3 
s^a_2 s^b_3 s^d_4 s^d_1 s^e_2 s^e_3 s^d_2 
s^d_3 s^b_4 s^a_3 s^a_2 s^b_3 s^d_4 s^d_3 
s^b_4 s^a_3 s^a_2 s^b_3 s^d_4 s^d_3 s^b_4 
s^b_1 s^c_2 s^c_3 s^b_2 s^a_3 s^a_4 s^b_3 
s^c_2 s^c_3 s^e_4 s^e_3 s^c_4 s^c_3 s^e_4 
s^e_3 s^d_2 s^d_3 s^e_2 s^e_3 s^c_4 
  (s^x_1) s^x_2 s^x_3 s^y_4 s^e_1 s^d_2 s^d_3 
s^e_2 \\ \cdot s^y_1 s^y_2 s^y_3 s^x_4 s^x_3 s^y_4 
s^y_3 s^y_2 s^y_3 s^y_2 s^y_3 s^y_2 s^y_3 
s^x_4 s^x_3 s^x_2 s^x_3 s^y_4 s^e_1 s^d_2 
s^d_3 s^b_4 s^a_3 s^a_2 s^b_3 s^d_4 s^d_3 
s^b_4 s^a_3 s^a_2 s^b_3 s^d_4 s^d_3 s^b_4 
s^a_3 s^a_2 s^a_1 s^a_2 s^b_3 s^d_4 s^d_3 
s^b_4 s^a_3 s^a_2 s^b_3 s^d_4 s^d_1 s^e_2 
s^e_3 s^c_4 \\ \cdot s^c_3 s^e_4 s^e_3 s^d_2 s^d_3 
s^e_2 s^y_1 s^y_2 s^y_3 s^y_2 s^y_3 s^x_4 
s^x_3 s^y_4 s^y_3 s^y_2 s^y_3 s^y_2 s^y_3 
s^y_2 s^y_3 s^x_4 s^x_3 s^x_2 s^x_3 s^y_4 
s^e_1 s^d_2 s^d_3 s^e_2 s^y_1 s^y_2 s^y_3 
s^x_4 s^x_3 s^y_4 s^y_3 s^y_2 s^y_3 s^y_2 
s^y_3 s^y_2 s^y_3 s^x_4 s^x_3 s^x_2 s^x_3 
s^y_4 s^e_1 s^d_2 \\ \cdot s^d_3 s^b_4 s^a_3 s^a_2 
s^b_3 s^d_4 s^d_3 s^b_4 s^a_3 s^a_2 s^b_3 
s^d_4 s^d_3 s^b_4 s^a_3 s^a_2 s^a_1 s^a_2 
s^b_3 s^d_4 s^d_3 s^b_4 s^a_3 s^a_2 s^a_1 
s^a_2 s^b_3 s^d_4 s^d_1 s^e_2 s^e_3 s^d_2 
s^d_3 s^b_4 s^a_3 s^a_2 s^b_3 s^d_4 s^d_3 
s^b_4 s^a_3 s^a_2 s^b_3 s^d_4 s^d_3 s^b_4 
s^a_3 s^a_2 s^a_1 s^a_2 \\ \cdot s^b_3 s^d_4 s^d_3 
s^b_4 s^a_3 s^a_2 s^a_1 s^a_2 s^b_3 s^d_4 
s^d_1 s^e_2 s^e_3 s^d_2 s^d_3 s^b_4 s^a_3 
s^a_2 s^b_3 s^d_4 s^d_3 s^b_4 s^a_3 s^a_2 
s^b_3 s^d_4 s^d_3 s^b_4 s^a_3 s^a_2 s^a_1 
s^a_2 s^b_3 s^d_4 s^d_3 s^b_4 s^a_3 s^a_2 
s^b_3 s^d_4 s^d_1 s^e_2 s^e_3 s^c_4 s^c_3 
s^e_4 s^e_3 s^d_2 s^d_3 s^e_2 \\ \cdot s^y_1 s^y_2 
s^y_3 s^y_2 s^y_3 s^x_4 s^x_3 s^x_2 s^x_3 
s^y_4 s^e_1 s^d_2 s^d_3 s^b_4 s^a_3 s^a_2 
s^b_3 s^d_4 s^d_3 s^b_4 s^a_3 s^a_2 s^b_3 
s^d_4 s^d_3 s^b_4 s^a_3 s^a_2 s^a_1 s^a_2 
s^b_3 s^d_4 s^d_3 s^b_4 s^a_3 s^a_2 s^b_3 
s^d_4 s^d_1 s^e_2 s^e_3 s^c_4 s^c_3 s^e_4 
s^e_3 s^d_2 s^d_3 s^e_2 s^y_1 s^y_2 \\ \cdot s^y_3 
s^y_2 s^y_3 s^x_4 s^x_3 s^x_2 s^x_3 s^y_4 
s^e_1 s^d_2 s^d_3 s^b_4 s^a_3 s^a_2 s^b_3 
s^d_4 s^d_3 s^b_4 s^a_3 s^a_2 s^b_3 s^d_4 
s^d_3 s^b_4 s^a_3 s^a_2 s^a_1 s^a_2 s^b_3 
s^d_4 s^d_3 s^b_4 s^a_3 s^a_2 s^b_3 s^d_4 
s^d_1 s^e_2 s^e_3 s^c_4 s^c_3 s^e_4 s^e_3 
s^d_2 s^d_3 s^e_2 s^y_1 s^y_2 s^y_3 s^x_4 \\ \cdot 
s^x_3 s^x_2 s^x_3 s^y_4 s^e_1 s^d_2 s^d_3 
s^b_4 s^a_3 s^a_2 s^b_3 s^d_4 s^d_3 s^b_4 
s^a_3 s^a_2 s^b_3 s^d_4 s^d_3 s^b_4 s^a_3 
s^a_2 s^a_1 s^a_2 s^b_3 s^d_4 s^d_3 s^b_4 
s^a_3 s^a_2 s^a_1 s^a_2 s^b_3 s^d_4 s^d_1 
s^e_2 s^e_3 s^y_1 
\end{array}
\endgroup
$}}}
\\ \quad
\\ \mbox{$\FigFiftyEightA$: Special Hamilton circuit of \cite[Table~3, Row~21]{Hec09},
where $S^\prime$ is that} 
\\ \mbox{ of $\FigFiftyEightA^\prime$ and let $x:={\bar{a}}$ and $y:={\bar{b}}$ and let $(s^z_i)$'s mean that}
\\ \mbox{this is exactly special (Length $=2688$)}
\end{array}
\end{equation*}

\noindent
{\bf{Rank-4-Case-3:}} Let $\chamchi$ be of \cite[Table~3,\,Rows 9, 14, 17]{Hec09}.
As mentioned by Table~$\tabletwo$, 
it is of $(i,j)$-convenient
for some $i$ and $j$,
whence a Hamilton circuit map of 
$\chamGammachi$ exists by Proposition~\ref{proposition:maintoolcor}.

\begin{center}
\begingroup
\renewcommand{\arraystretch}{1.5}
Table~$\tabletwo$: Rank-4  $(i,j)$-convenient bicharacters \\
  \begin{tabular}{cccc}
    \hline
$\chamchi$  &  
$\acchi^{(3,9)}_1$ &  
$\acchi^{(3,14)}_1$ &
$\acchi^{(3,17)}_1$  
\\
    \hline \hline
\mbox{Fig.s} &
\NFigFiftyNineA, 
\NFigSixtyA 
&
\NFigSixtyTwoA, 
\NFigSixtyThreeA, 
&
\NFigSixtyFourA, 
\NFigSixtyFiveA 
\\
\mbox{Extracted $\al_i$}  & 
$\al_3$ &
$\al_1$ &
$\al_3$
\\
\mbox{Rep.s}  & 
$\{a,c\}$&
$\{a,h,j\}$ &
$\{g,b,i\}$ 
\\
\mbox{Fig.s for $\chamchi_{\langle i\rangle}$}  &
\NFigSixteenA, 
\NFigThirtyTwoA 
&
\NFigThirtyTwoA, 
\NFigSixteenA, \NFigSixteenA
&
\NFigSeventeenA, 
\NFigThirtyTwoA, 
\NFigSeventeenA 
\\ 
\mbox{Convenient}  &
(3,1) &
(1,3) &
(3,1) 
\\
\hline
  \end{tabular}
\endgroup
\end{center}

\begin{equation*}
\begin{array}{l}
$\allHecFourNine$
\\ \quad
\\ \mbox{$\FigFiftyNineA$: Generalized Dynkin diagram of \cite[Table~3, Row~9]{Hec09}}
\end{array}
\end{equation*}

\begin{equation*}
\begin{array}{l}
$\RankFourHecNineTau$
\\ \quad
\\ \mbox{$\FigSixtyA$: Changing of diagrams of $\FigFiftyNineA$} 
\end{array}
\end{equation*}

\begin{equation*}
\begin{array}{l}
{\footnotesize{\mbox{$\begin{array}{l}
(s^f_2) (s^f_4) (s^f_3) s^f_4 (s^d_2) 
 (s^e_4) (s^e_2) (s^d_4) (s^c_3) (s^c_4)
  (s^c_2) s^c_4 (s^d_3) s^e_4 
 (s^e_3) s^d_4 s^f_2 s^f_4 s^f_3 s^f_4 s^d_2 
s^e_4 s^e_2 s^d_4 s^c_3 s^c_4 s^c_2 s^c_4 
s^d_3 s^e_4 s^e_3 \\ \cdot 
 (s^e_1) s^e_3 s^d_4 s^f_2 s^f_4 s^f_3 s^f_4 
s^d_2 s^e_4 \\ \cdot s^e_2 s^d_4 s^c_3 s^c_4 s^c_2 
s^c_4 s^d_3 s^e_4 s^e_3 s^d_4 s^f_2 s^f_4 
s^f_3 s^f_4 s^d_2 s^e_4 s^e_2 s^d_4 s^c_3 
s^c_4 s^c_2 (s^b_1) (s^a_2) (s^a_4) 
 (s^a_3) s^a_4 s^a_3 s^a_4 (s^b_2) 
 (s^c_1)  \\ \cdot s^d_3 s^e_4 s^e_3 s^d_4 s^f_2 s^f_4 
s^f_3 s^f_4 s^d_2 s^e_4 \\ \cdot s^e_2 s^d_4 s^c_3 
s^c_4 s^c_2 s^c_4 s^d_3 s^e_4 s^e_3 s^d_4 
s^f_2 s^f_4 s^f_3 s^f_4 s^d_2 s^e_4 s^e_2 
s^d_4 s^c_3 s^b_1 s^a_2 s^a_4 s^a_3 s^a_4 
s^a_3 s^a_4 s^b_2 (s^b_4) 
 (s^b_3) s^c_1 s^d_3 s^e_4 s^e_3 s^d_4 s^f_2 
s^f_4 s^f_3 s^f_4 s^d_2 s^e_4 \\ \cdot s^e_2 s^d_4 
s^c_3 s^c_4 s^c_2 s^c_4 s^d_3 s^e_4 s^e_3 
s^d_4 s^f_2 s^f_4 s^f_3 s^f_4 s^d_2 s^e_4 
s^e_2 s^d_4 s^c_3 s^b_1 s^a_2 s^a_4 s^a_3 
s^a_4 s^a_3 s^a_4 s^b_2 s^b_4 s^b_3 s^c_1 
s^d_3 s^e_4 s^e_3 s^d_4 s^f_2 s^f_4 s^f_3 
s^f_4 s^d_2 s^e_4 \\ \cdot s^e_2 s^d_4 s^c_3 s^c_4 
s^c_2 s^c_4 s^d_3 s^e_4 s^e_3 s^d_4 s^f_2 
s^f_4 s^f_3 s^f_4 s^d_2 s^e_4 s^e_2 s^d_4 
s^c_3 s^b_1 s^a_2 s^a_4 s^a_3 s^a_4 s^a_3 
s^a_4 s^b_2 s^c_1 s^d_3 s^e_4 s^e_3 s^d_4 
s^f_2 s^f_4 s^f_3 s^f_4 s^d_2 s^e_4 s^e_2 
s^d_4 \\ \cdot s^c_3 s^c_4 s^c_2 s^c_4 s^d_3 s^e_4 
s^e_3 s^d_4 s^f_2 s^f_4 s^f_3 s^f_4 s^d_2 
s^e_4 s^e_2 s^d_4 s^c_3 s^b_1 s^a_2 s^a_4 
s^a_3 s^a_4 s^a_3 s^a_4 s^b_2 s^b_4 s^b_3 
s^c_1 s^d_3 s^e_4 s^e_3 s^d_4 s^f_2 s^f_4 
s^f_3 s^f_4 s^d_2 s^e_4 s^e_2 s^d_4 \\ \cdot s^c_3 
s^c_4 s^c_2 s^c_4 s^d_3 s^e_4 s^e_3 s^d_4 
s^f_2 s^f_4 s^f_3 s^f_4 s^d_2 s^e_4 s^e_2 
s^d_4 s^c_3 s^b_1 s^a_2 s^a_4 s^a_3 s^a_4 
s^a_3 s^a_4 s^b_2 s^c_1 s^d_3 s^e_4 s^e_3 
s^d_4 s^f_2 s^f_4 s^f_3 s^f_4 s^d_2 s^e_4 
s^e_2 s^d_4 s^c_3 s^c_4 \\ \cdot s^c_2 s^c_4 s^d_3 
s^e_4 s^e_3 s^d_4 s^f_2 s^f_4 s^f_3 s^f_4 
s^d_2 s^e_4 s^e_2 s^d_4 s^c_3 s^b_1 s^a_2 
s^a_4 s^a_3 s^a_4 s^a_3 s^a_4 s^b_2 s^b_4 
s^b_3 s^c_1 s^d_3 s^e_4 s^e_3 s^d_4 s^f_2 
s^f_4 s^f_3 s^f_4 s^d_2 s^e_4 s^e_2 s^d_4 
s^c_3 s^c_4 \\ \cdot s^c_2 s^c_4 s^d_3 s^e_4 s^e_3 
s^d_4 s^f_2 
 (s^f_1) s^d_2 s^e_4 s^e_2 s^d_4 s^c_3 s^c_4 
s^c_2 s^c_4 s^d_3 s^e_4 s^e_3 s^d_4 s^f_2 
s^f_4 s^f_3 s^f_4 s^d_2 s^e_4 s^e_2 s^d_4 
s^c_3 s^c_4 s^c_2 s^c_4 s^d_3 s^e_4 s^e_3 
s^d_4 s^f_2 s^f_1 s^d_2 s^e_4 \\ \cdot s^e_2 s^d_4 
s^c_3 s^c_4 s^c_2 s^c_4 s^d_3 s^e_4 s^e_3 
s^d_4 s^f_2 s^f_4 s^f_3 s^f_4 s^d_2 s^e_4 
s^e_2 s^d_4 s^c_3 s^c_4 s^c_2 s^b_1 s^a_2 
s^a_4 s^a_3 s^a_4 s^a_3 
(s^a_1) s^b_2 s^b_4 s^b_3 s^b_4 s^b_3 s^b_4 
s^a_2 s^a_4 s^a_3 s^a_1 s^b_2 s^b_4 \\ \cdot s^b_3 
s^b_4 s^b_3 s^b_4 s^a_2 s^a_4 s^a_3 s^a_1 
s^b_2 s^b_4 s^b_3 s^b_4 s^b_3 s^b_4 s^a_2 
s^a_1 s^b_2 s^c_1 s^d_3 s^e_4 s^e_3 s^d_4 
s^f_2 s^f_4 s^f_3 s^f_4 s^d_2 s^e_4 s^e_2 
s^d_4 s^c_3 s^b_1 s^a_2 s^a_4 s^a_3 s^a_4 
s^a_3 s^a_4 s^b_2 s^b_4 \\ \cdot s^b_3 s^b_4 s^b_3 
s^c_1 s^c_2 s^c_4 s^d_3 s^e_4 s^e_3 s^d_4 
s^f_2 s^f_4 s^f_3 s^f_4 s^d_2 s^e_4 s^e_2 
s^d_4 s^c_3 s^c_4 s^c_2 s^b_1 s^b_3 s^b_4 
s^b_3 s^b_4 s^a_2 s^a_4 s^a_3 s^a_1 s^b_2 
s^b_4 s^b_3 s^b_4 s^b_3 s^b_4 s^a_2 s^a_1 
s^b_2 s^b_4 \\ \cdot s^b_3 s^b_4 s^b_3 s^b_4 s^a_2 
s^a_1 s^b_2 s^b_4 s^b_3 s^b_4 s^b_3 s^b_4 
s^a_2 s^a_4 s^a_3 s^a_4 s^a_3 s^a_1 s^a_3 
s^a_4 s^a_3 s^a_1 s^a_3 s^a_4 s^a_3 s^a_1 
s^a_3 s^a_1 s^a_3 s^a_4 s^a_3 s^a_1 s^a_3 
s^a_1 s^a_3 s^a_1 s^b_2 s^b_4 s^b_3 s^b_4 \\ \cdot 
s^b_3 s^c_1 s^d_3 s^e_4 s^e_3 s^d_4 s^f_2 
s^f_4 s^f_3 s^f_1 s^f_3 s^f_1 s^f_3 s^f_4 
s^d_2 s^e_4 s^e_2 s^d_4 s^c_3 s^c_4 s^c_2 
s^b_1 s^b_3 s^c_1 s^c_2 s^b_1 s^b_3 s^b_4 
s^b_3 s^c_1 s^c_2 s^b_1 s^b_3 s^c_1 s^c_2 
s^b_1 s^b_3 s^b_4 s^b_3 s^c_1 \\ \cdot s^c_2 s^b_1 
s^b_3 s^c_1 s^c_2 s^b_1 s^b_3 s^c_1 s^c_2 
s^b_1 s^b_3 s^b_4 s^b_3 s^c_1 s^d_3 
 (s^d_1) 
\end{array}$}}}
\\ \quad
\\ \mbox{$\FigSixtyOneA$: Special Hamilton circuit of \cite[Table~3, Row~9]{Hec09}
(Length $=576$)}
\end{array}
\end{equation*}

\begin{equation*}
\begin{array}{l}
$\allHecFourFourteen$
\\ \quad
\\ \mbox{$\FigSixtyTwoA$: Generalized Dynkin diagram of \cite[Table~3, Row~14]{Hec09}}
\end{array}
\end{equation*}

\begin{equation*}
\begin{array}{l}
$\RankFourHecFourteenTau$
\\ \quad
\\ \mbox{$\FigSixtyThreeA$: Changing of diagrams of $\FigSixtyTwoA$} 
\end{array}
\end{equation*}

\begin{equation*}
\begin{array}{l}
{\footnotesize{\mbox{$\begin{array}{l}
(s^b_3) (s^d_4) (s^d_3) (s^f_2) (s^f_3)
  (s^g_4) (s^i_3) (s^i_2) (s^g_3) 
 (s^e_2) (s^e_3) (s^c_4) (s^c_3) (s^b_2)
  (s^a_3) 
 (s^a_4) \\ \cdot s^b_3 s^d_4 s^d_3 s^f_2 s^f_3 s^g_4 
s^i_3 s^i_2 s^g_3 s^e_2 s^e_3 s^c_4 s^c_3 
s^b_2 (s^b_1) (s^c_2) s^c_3 (s^e_4) s^e_3 
 (s^g_2) s^i_3 s^i_2 s^g_3 
 (s^f_4) \\ \cdot s^f_3 (s^d_2) s^d_3 
 (s^b_4) s^a_3 s^a_4 s^b_3 s^c_2 s^c_3 s^e_4 
s^e_3 s^g_2 s^i_3 s^i_2 s^g_3 s^f_4 (s^j_1) 
 (s^j_4) (s^j_3) s^j_4 s^j_3 
 (s^j_2) s^j_3 s^j_4 
 (s^f_1) \\ \cdot s^d_2 s^d_3 s^b_4 s^a_3 s^a_4 s^b_3 
s^c_2 s^c_3 s^e_4 s^e_3 s^g_2 s^i_3 s^i_2 
s^g_3 s^f_4 \\ \cdot s^f_3 s^d_2 s^d_3 s^b_4 s^a_3 
s^a_4 s^b_3 s^c_2 (s^h_1) (s^h_2)
(s^h_3) s^h_2 s^h_3 (s^h_4) s^h_3 s^h_2 
 (s^c_1)  \\ \cdot s^e_4 s^e_3 s^g_2 s^i_3 s^i_2 s^g_3 
s^f_4 s^f_3 s^d_2 s^d_3 s^b_4 s^a_3 s^a_4 
s^b_3 s^c_2 s^c_3 s^e_4 s^e_3 s^g_2 s^i_3 
s^i_2 s^g_3 s^f_4 \\ \cdot s^j_1 s^j_4 s^j_3 s^j_4 
s^j_3 s^j_2 s^j_3 s^j_4 s^f_1 s^d_2 s^d_3 
s^b_4 s^a_3 s^a_4 s^b_3 s^c_2 s^c_3 s^e_4 
s^e_3 s^g_2 s^i_3 s^i_2 s^g_3 s^f_4 s^f_3 
s^d_2 s^d_3 s^b_4 s^a_3 s^a_4 s^b_3 s^c_2 
s^h_1 s^h_2 s^h_3 s^h_2 s^c_1 s^e_4 s^e_3 
s^g_2 \\ \cdot s^i_3 s^i_2 s^g_3 s^f_4 s^f_3 s^d_2 
s^d_3 s^b_4 s^a_3 s^a_4 s^b_3 s^c_2 s^c_3 
s^e_4 s^e_3 s^g_2 s^i_3 s^i_2 
 (s^i_1) s^i_2 s^g_3 s^f_4 s^f_3 s^d_2 s^d_3 
s^b_4 s^a_3 s^a_4 s^b_3 s^c_2 s^c_3 s^e_4 
s^e_3 s^g_2 s^i_3 s^i_2 s^g_3 s^f_4 s^f_3 
s^d_2 \\ \cdot s^d_3 s^b_4 s^a_3 s^a_4 s^b_3 s^c_2 
s^h_1 s^h_2 s^h_3 s^h_2 s^h_3 s^h_4 s^c_1 
s^e_4 s^e_3 s^g_2 s^i_3 s^i_2 s^g_3 s^f_4 
s^j_1 s^j_4 s^f_1 s^d_2 s^d_3 s^b_4 s^a_3 
s^a_4 s^b_3 s^c_2 s^c_3 s^e_4 s^e_3 s^g_2 
s^i_3 s^i_2 s^g_3 s^f_4 s^f_3 s^d_2 \\ \cdot s^d_3 
s^b_4 s^a_3 s^a_4 s^b_3 s^c_2 s^c_3 s^e_4 
s^e_3 s^g_2 s^i_3 s^i_2 s^g_3 s^f_4 s^j_1 
s^j_4 s^j_3 s^j_2 s^j_3 s^j_4 s^j_3 s^j_4 
s^j_3 s^j_2 s^f_1 s^d_2 s^d_3 s^b_4 s^a_3 
s^a_4 s^b_3 s^c_2 s^c_3 s^e_4 s^e_3 s^g_2 
s^i_3 s^i_2 s^g_3 s^f_4 \\ \cdot s^f_3 s^d_2 s^d_3 
s^b_4 s^a_3 s^a_4 s^b_3 s^c_2 s^h_1 s^h_2 
s^h_3 s^h_2 s^h_3 s^h_4 s^c_1 s^e_4 s^e_3 
s^g_2 (s^g_1) s^f_4 s^f_3 s^d_2 s^d_3 s^b_4 
s^a_3 s^a_4 s^b_3 s^c_2 s^h_1 s^h_4 s^h_3 
s^h_2 s^h_3 s^h_2 s^h_3 s^h_4 s^c_1 s^e_4 
s^e_3 s^g_2 \\ \cdot s^i_3 s^i_2 s^g_3 s^f_4 s^j_1 
s^j_4 s^j_3 s^j_2 s^f_1 s^d_2 s^d_3 s^b_4 
s^a_3 s^a_4 s^b_3 s^c_2 s^h_1 s^h_2 s^h_3 
s^h_4 s^c_1 s^e_4 s^e_3 s^g_2 s^i_3 s^i_2 
s^g_3 s^f_4 s^j_1 s^j_4 s^j_3 s^j_2 s^f_1 
s^d_2 s^d_3 s^b_4 s^a_3 s^a_4 
 (s^a_1) s^a_4 
\end{array}$}}}
\\ \quad
\\ \mbox{$\FigSixtyFourA$: Hamilton circuit of \cite[Table~3, Row~14]{Hec09}
missing $s^a_2$, $s^d_1$, $s^e_1$ and $s^i_4$} 
\\ \mbox{(From this Hamilton circuit, by $\FigSixtyTwoA$, 
we obtain another Hamilton circuit}
\\ 
\mbox{of \cite[Table~3, Row~14]{Hec09} missing $s^h_2$, $s^e_3$, $s^d_3$ and $s^j_4$.
Indeed it is obtained by} 
\\ 
\mbox{replacing $1$, $3$ and $x\in\{a,\ldots,j\}$
by $3$, $1$ and $x\circ\left[{{1234}\atop{3214}}\right]$ respectively.)}
\\ 
\mbox{(Length $=360$)}
\end{array}
\end{equation*}

\begin{equation*}
\begin{array}{l}
$\allHecFourSeventeen$
\\ \quad
\\ \mbox{$\FigSixtyFiveA$: Generalized Dynkin diagram of \cite[Table~3, Row~17]{Hec09}}
\end{array}
\end{equation*}

\begin{equation*}
\begin{array}{l}
$\RankFourHecSeventeenTau$
\\ \quad
\\ \mbox{$\FigSixtySixA$: Changing of diagrams of $\FigSixtyFiveA$} 
\end{array}
\end{equation*}

\begin{equation*}
\begin{array}{l}
{\tiny{\mbox{$
\begingroup
\renewcommand{\arraystretch}{1.3}
\begin{array}{l}
(s^b_2) (s^b_1) (s^a_4) (s^a_1) (s^c_2)
  (s^k_1) (s^k_2) (s^c_1) (s^d_4) 
 (s^d_1) (s^h_2) (s^h_1) (s^e_4) (s^l_1)
  (s^l_4) 
 (s^e_1) s^b_2 s^b_1 s^a_4 s^a_1 s^c_2 s^k_1 
s^k_2 s^c_1 s^d_4 s^d_1 s^h_2 s^h_1 s^e_4 
s^l_1 s^l_4 
 (s^l_3) \\ \cdot s^l_4 s^e_1 s^b_2 s^b_1 s^a_4 s^a_1 
s^c_2 s^k_1 s^k_2 s^c_1 s^d_4 s^d_1 s^h_2 
s^h_1 s^e_4 s^l_1 s^l_4 s^e_1 \\ \cdot s^b_2 s^b_1 
s^a_4 s^a_1 s^c_2 s^k_1 s^k_2 s^c_1 s^d_4 
 (s^i_3) (s^j_2) (s^j_1) 
 (s^j_4) s^j_1 s^j_4 s^j_1 (s^i_2) (s^i_1) 
 (s^i_4) 
 (s^d_3) s^h_2 s^h_1 s^e_4 s^l_1 s^l_4 s^e_1 
s^b_2 s^b_1 s^a_4 s^a_1 s^c_2 s^k_1 s^k_2 
s^c_1 s^d_4 s^d_1 s^h_2 s^h_1 s^e_4 s^l_1 
s^l_4 s^e_1 s^b_2 \\ \cdot (s^f_3) (s^g_4) (s^g_1) 
 (s^g_2) s^g_1 s^g_2 s^g_1 \\ \cdot (s^f_4) 
 (s^f_1) (s^f_2) 
 (s^b_3) s^a_4 s^a_1 s^c_2 s^k_1 s^k_2 s^c_1 
s^d_4 s^d_1 s^h_2 s^h_1 s^e_4 s^l_1 s^l_4 
s^e_1 s^b_2 s^b_1 s^a_4 s^a_1 s^c_2 s^k_1 
s^k_2 s^c_1 s^d_4 s^i_3 s^j_2 s^j_1 s^j_4 
s^j_1 s^j_4 s^j_1 s^i_2 s^i_1 s^i_4 s^d_3 
s^h_2 s^h_1 s^e_4 s^l_1 s^l_4 s^e_1 s^b_2 
s^b_1 \\ \cdot s^a_4 s^a_1 s^c_2 s^k_1 \\ \cdot s^k_2 s^c_1 
s^d_4 s^d_1 s^h_2 s^h_1 s^e_4 s^l_1 s^l_4 
s^e_1 s^b_2 s^f_3 s^g_4 s^g_1 s^g_2 s^g_1 
s^g_2 s^g_1 s^f_4 s^f_1 s^f_2 s^b_3 s^a_4 
s^a_1 s^c_2 s^k_1 s^k_2 s^c_1 s^d_4 s^d_1 
s^h_2 s^h_1 s^e_4 s^l_1 s^l_4 s^e_1 s^b_2 
s^b_1 s^a_4 s^a_1 s^c_2 s^k_1 s^k_2 s^c_1 
s^d_4 s^i_3 s^j_2 s^j_1 s^j_4 s^j_1 \\ \cdot s^j_4 
s^j_1 s^i_2 s^d_3 s^h_2 s^h_1 s^e_4 s^l_1 
s^l_4 s^e_1 s^b_2 s^b_1 s^a_4 s^a_1 s^c_2 
s^k_1 s^k_2 s^c_1 s^d_4 s^d_1 s^h_2 s^h_1 
s^e_4 s^l_1 s^l_4 s^e_1 s^b_2 s^f_3 s^g_4 
s^g_1 s^g_2 s^g_1 s^g_2 s^g_1 s^f_4 s^f_1 
s^f_2 s^b_3 s^a_4 s^a_1 s^c_2 s^k_1 s^k_2 
s^c_1 s^d_4 s^d_1 s^h_2 s^h_1 s^e_4 s^l_1 \\ \cdot 
s^l_4 s^e_1 s^b_2 s^b_1 s^a_4 s^a_1 s^c_2 
s^k_1 s^k_2 s^c_1 s^d_4 s^i_3 s^j_2 s^j_1 
s^j_4 s^j_1 s^j_4 s^j_1 s^i_2 s^i_1 s^i_4 
s^d_3 s^h_2 s^h_1 s^e_4 s^l_1 s^l_4 s^e_1 
s^b_2 s^b_1 s^a_4 s^a_1 s^c_2 s^k_1 s^k_2 
s^c_1 s^d_4 s^d_1 s^h_2 s^h_1 s^e_4 s^l_1 
s^l_4 s^e_1 s^b_2 s^f_3 s^g_4 s^g_1 s^g_2 
s^g_1 \\ \cdot s^g_2 s^g_1 s^f_4 s^f_1 s^f_2 s^b_3 
s^a_4 s^a_1 s^c_2 s^k_1 s^k_2 s^c_1 s^d_4 
s^d_1 s^h_2 s^h_1 s^e_4 s^l_1 s^l_4 s^e_1 
s^b_2 s^f_3 s^g_4 s^g_1 s^g_2 s^g_1 s^g_2 
s^g_1 s^f_4 s^b_3 s^a_4 s^a_1 s^c_2 s^k_1 
s^k_2 s^c_1 s^d_4 s^d_1 s^h_2 s^h_1 s^e_4 
s^l_1 s^l_4 s^e_1 s^b_2 s^b_1 s^a_4 s^a_1 
s^c_2 s^k_1 \\ \cdot s^k_2 s^c_1 s^d_4 s^i_3 s^j_2 
s^j_1 s^j_4 s^j_1 s^j_4 s^j_1 s^i_2 s^i_1 
s^i_4 s^d_3 s^h_2 s^h_1 s^e_4 s^l_1 s^l_4 
s^e_1 s^b_2 s^b_1 s^a_4 s^a_1 s^c_2 s^k_1 
s^k_2 s^c_1 s^d_4 s^d_1 s^h_2 s^h_1 s^e_4 
s^l_1 s^l_4 s^e_1 s^b_2 s^f_3 s^g_4 s^g_1 
s^g_2 s^g_1 s^g_2 s^g_1 s^f_4 s^f_1 s^f_2 
s^b_3 s^a_4 s^a_1 \\ \cdot s^c_2 s^k_1 s^k_2 s^c_1 
s^d_4 s^d_1 s^h_2 s^h_1 s^e_4 s^l_1 s^l_4 
s^e_1 s^b_2 s^b_1 s^a_4 s^a_1 s^c_2 s^k_1 
s^k_2 s^c_1 s^d_4 s^i_3 s^j_2 s^j_1 s^j_4 
s^j_1 s^j_4 s^j_1 s^i_2 s^i_1 s^i_4 s^d_3 
s^h_2 s^h_1 s^e_4 s^l_1 s^l_4 s^e_1 s^b_2 
s^b_1 s^a_4 s^a_1 s^c_2 s^k_1 s^k_2 s^c_1 
s^d_4 s^d_1 s^h_2 s^h_1 \\ \cdot s^e_4 s^l_1 s^l_4 
s^e_1 s^b_2 s^f_3 s^g_4 s^g_1 s^g_2 s^g_1 
s^g_2 s^g_1 s^f_4 s^f_1 s^f_2 s^b_3 s^a_4 
s^a_1 s^c_2 s^k_1 s^k_2 s^c_1 s^d_4 s^d_1 
s^h_2 s^h_1 s^e_4 s^l_1 s^l_4 s^e_1 s^b_2 
s^b_1 s^a_4 s^a_1 s^c_2 s^k_1 s^k_2 s^c_1 
s^d_4 s^i_3 s^j_2 s^j_1 s^j_4 s^j_1 s^j_4 
 (s^j_3) s^i_2 s^i_1 s^i_4 s^i_1 \\ \cdot s^i_4 
s^d_3 s^h_2 s^h_1 s^e_4 s^l_1 s^l_4 s^e_1 
s^b_2 s^b_1 s^a_4 s^a_1 s^c_2 s^k_1 s^k_2 
s^c_1 s^d_4 s^d_1 s^h_2 s^h_1 s^e_4 s^l_1 
s^l_4 s^e_1 s^b_2 s^b_1 s^a_4 s^a_1 s^c_2 
s^k_1 s^k_2 
 (s^k_3) s^k_2 s^c_1 s^d_4 s^d_1 s^h_2 s^h_1 
s^e_4 s^l_1 s^l_4 s^e_1 s^b_2 s^b_1 s^a_4 
s^a_1 s^c_2 s^k_1 s^k_2 s^k_3 \\ \cdot s^k_2 s^c_1 
s^d_4 s^d_1 s^h_2 s^h_1 s^e_4 s^l_1 s^l_4 
s^e_1 s^b_2 s^b_1 s^a_4 s^a_1 s^c_2 s^k_1 
s^k_2 s^c_1 s^d_4 s^d_1 s^h_2 s^h_1 s^e_4 
s^l_1 s^l_4 s^e_1 s^b_2 s^f_3 s^g_4 s^g_1 
s^g_2 s^g_1 s^g_2 s^g_1 s^f_4 s^f_1 s^f_2 
s^b_3 s^a_4 s^a_1 s^c_2 s^k_1 s^k_2 s^c_1 
s^d_4 s^d_1 s^h_2 s^h_1 s^e_4 s^l_1 \\ \cdot s^l_4 
s^e_1 s^b_2 s^b_1 s^a_4 s^a_1 s^c_2 s^k_1 
s^k_2 s^c_1 s^d_4 s^i_3 s^j_2 s^j_1 s^j_4 
s^j_1 s^j_4 s^j_1 s^i_2 s^i_1 s^i_4 s^d_3 
s^h_2 s^h_1 s^e_4 s^l_1 s^l_4 s^e_1 s^b_2 
s^b_1 s^a_4 s^a_1 s^c_2 s^k_1 s^k_2 s^c_1 
s^d_4 s^d_1 s^h_2 s^h_1 s^e_4 s^l_1 s^l_4 
s^e_1 s^b_2 s^f_3 s^g_4 s^g_1 s^g_2 s^g_1 \\ \cdot 
s^g_2 (s^g_3) s^f_4 s^b_3 s^a_4 s^a_1 s^c_2 
s^k_1 s^k_2 s^c_1 s^d_4 s^d_1 s^h_2 s^h_1 
s^e_4 s^l_1 s^l_4 s^e_1 s^b_2 s^b_1 s^a_4 
s^a_1 s^c_2 s^k_1 s^k_2 s^c_1 s^d_4 s^i_3 
s^j_2 s^j_1 s^j_4 s^j_1 s^j_4 s^j_3 s^i_2 
s^i_1 s^i_4 s^i_1 s^i_4 s^i_1 s^j_2 s^j_1 
s^j_4 s^j_3 s^i_2 s^i_1 s^i_4 s^i_1 s^i_4 
s^i_1 \\ \cdot s^j_2 s^j_3 s^i_2 s^i_1 s^i_4 s^i_1 
s^i_4 s^d_3 s^h_2 s^h_1 s^e_4 s^l_1 s^l_4 
s^e_1 s^b_2 s^b_1 s^a_4 s^a_1 s^c_2 s^k_1 
s^k_2 s^c_1 s^d_4 s^d_1 s^h_2 s^h_1 s^e_4 
s^l_1 s^l_4 s^e_1 s^b_2 s^f_3 s^g_4 s^g_1 
s^g_2 s^g_1 s^g_2 s^g_1 s^f_4 s^f_1 s^f_2 
s^b_3 s^a_4 s^a_1 s^c_2 s^k_1 s^k_2 s^c_1 
s^d_4 s^d_1 \\ \cdot s^h_2 s^h_1 s^e_4 s^l_1 s^l_4 
s^e_1 s^b_2 s^b_1 s^a_4 s^a_1 s^c_2 s^k_1 
s^k_2 s^k_3 s^k_2 s^c_1 s^d_4 s^d_1 s^h_2 
s^h_1 s^e_4 s^l_1 s^l_4 s^e_1 s^b_2 s^f_3 
s^g_4 s^g_1 s^g_2 s^g_1 s^g_2 s^g_1 s^f_4 
s^f_1 s^f_2 s^f_1 s^f_2 s^b_3 s^a_4 s^a_1 
s^c_2 s^k_1 s^k_2 s^k_3 s^k_2 s^c_1 s^d_4 
s^d_1 s^h_2 s^h_1 \\ \cdot s^e_4 s^l_1 s^l_4 s^e_1 
s^b_2 s^f_3 s^g_4 s^g_1 s^g_2 s^g_1 s^g_2 
s^g_3 s^f_4 s^f_1 s^f_2 s^b_3 s^a_4 s^a_1 
s^c_2 s^k_1 s^k_2 s^c_1 s^d_4 s^d_1 s^h_2 
s^h_1 s^e_4 s^l_1 s^l_4 s^e_1 s^b_2 s^b_1 
s^a_4 s^a_1 s^c_2 s^k_1 s^k_2 s^c_1 s^d_4 
s^i_3 s^j_2 s^j_1 s^j_4 s^j_1 s^j_4 s^j_3 
s^i_2 s^i_1 s^i_4 s^i_1 \\ \cdot s^i_4 s^i_1 s^j_2 
s^j_1 s^j_4 s^j_3 s^i_2 s^i_1 s^i_4 s^i_1 
s^i_4 s^i_1 s^j_2 s^j_1 s^j_4 s^j_1 s^j_4 
s^j_3 s^j_4 s^j_3 s^i_2 s^i_1 s^i_4 s^d_3 
s^h_2 s^h_1 s^e_4 s^l_1 s^l_4 s^e_1 s^b_2 
s^b_1 s^a_4 s^a_1 s^c_2 s^k_1 s^k_2 s^c_1 
s^d_4 s^i_3 s^j_2 s^j_1 s^j_4 s^j_1 s^j_4 
s^j_1 s^i_2 s^d_3 s^h_2 s^h_1 \\ \cdot s^e_4 s^l_1 
s^l_4 s^e_1 s^b_2 s^f_3 s^g_4 s^g_3 s^f_4 
s^f_1 s^f_2 s^f_1 s^f_2 s^f_1 s^g_4 s^g_1 
s^g_2 s^g_1 s^g_2 s^g_3 s^g_2 s^g_1 s^g_2 
s^g_1 s^f_4 s^f_1 s^f_2 s^f_1 s^f_2 s^b_3 
s^a_4 s^a_1 s^c_2 s^k_1 s^k_2 s^c_1 s^d_4 
s^d_1 s^h_2 s^h_1 s^e_4 s^l_1 s^l_4 s^e_1 
s^b_2 s^f_3 s^g_4 s^g_3 s^f_4 s^f_1 \\ \cdot s^f_2 
s^b_3 s^a_4 s^a_1 s^c_2 s^k_1 s^k_2 s^c_1 
s^d_4 s^i_3 s^j_2 s^j_1 s^j_4 s^j_1 s^j_4 
s^j_3 s^i_2 s^i_1 s^i_4 s^i_1 s^i_4 s^i_1 
s^j_2 s^j_3 s^i_2 s^i_1 s^i_4 s^i_1 s^i_4 
s^i_1 s^j_2 s^j_1 s^j_4 s^j_1 s^j_4 s^j_3 
s^j_4 s^j_1 s^j_4 s^j_3 s^i_2 s^i_1 s^i_4 
s^i_1 s^i_4 s^d_3 s^h_2 s^h_1 s^e_4 s^l_1 \\ \cdot 
s^l_4 s^e_1 s^b_2 s^b_1 s^a_4 s^a_1 s^c_2 
s^k_1 s^k_2 s^c_1 s^d_4 s^d_1 s^h_2 s^h_1 
s^e_4 s^l_1 s^l_4 s^e_1 s^b_2 s^f_3 s^f_2 
s^f_1 s^g_4 s^g_1 s^g_2 s^g_1 s^g_2 s^g_3 
s^g_2 s^g_1 s^g_2 s^g_1 s^f_4 s^f_1 s^f_2 
s^f_1 s^f_2 s^b_3 s^a_4 s^a_1 s^c_2 s^k_1 
s^k_2 s^c_1 s^d_4 s^i_3 s^i_4 s^i_1 s^i_4 
s^d_3 \\ \cdot s^h_2 s^h_1 s^e_4 s^l_1 s^l_4 s^e_1 
s^b_2 s^b_1 s^a_4 s^a_1 s^c_2 s^k_1 s^k_2 
s^c_1 s^d_4 s^i_3 s^i_4 s^d_3 s^h_2 s^h_1 
s^e_4 s^l_1 s^l_4 s^e_1 s^b_2 s^f_3 s^f_2 
s^f_1 s^g_4 s^g_3 s^f_4 s^f_1 s^f_2 s^f_1 
s^f_2 s^b_3 s^a_4 s^a_1 s^c_2 s^k_1 s^k_2 
s^c_1 s^d_4 s^d_1 s^h_2 s^h_1 s^e_4 s^l_1 
s^l_4 s^e_1 \\ \cdot s^b_2 s^b_1 s^a_4 s^a_1 s^c_2 
s^k_1 s^k_2 s^c_1 s^d_4 s^d_1 s^h_2 s^h_1 
s^e_4 s^l_1 s^l_4 s^e_1 s^b_2 s^f_3 s^g_4 
s^g_1 s^g_2 s^g_1 s^g_2 s^g_3 s^g_2 s^g_1 
s^g_2 s^g_3 s^f_4 s^f_1 s^f_2 s^f_1 s^f_2 
s^b_3 s^a_4 s^a_1 s^c_2 s^k_1 s^k_2 s^c_1 
s^d_4 s^d_1 s^h_2 s^h_1 s^e_4 s^l_1 s^l_4 
s^e_1 s^b_2 s^b_1 \\ \cdot s^a_4 s^a_1 s^c_2 
 (s^c_3) s^d_4 s^d_1 s^h_2 s^h_1 s^e_4 s^l_1 
s^l_4 s^e_1 s^b_2 s^b_1 s^a_4 s^a_1 s^c_2 
s^c_3 s^d_4 s^d_1 s^h_2 s^h_1 s^e_4 s^l_1 
s^l_4 s^e_1 s^b_2 s^f_3 s^f_2 s^b_3 s^a_4 
s^a_1 s^c_2 s^k_1 s^k_2 s^c_1 s^d_4 s^i_3 
s^j_2 s^j_1 s^j_4 s^j_1 s^j_4 s^j_3 s^j_4 
s^j_1 s^j_4 s^j_3 s^j_4 s^j_3 \\ \cdot s^i_2 s^i_1 
s^i_4 s^i_1 s^i_4 s^d_3 s^h_2 s^h_1 s^e_4 
s^l_1 s^l_4 s^e_1 s^b_2 s^f_3 s^f_2 s^f_1 
s^f_2 s^f_1 s^g_4 s^g_3 s^f_4 s^f_1 s^f_2 
s^f_1 s^f_2 s^f_1 s^g_4 s^g_1 s^g_2 s^g_1 
s^g_2 s^g_3 s^g_2 s^g_1 s^g_2 s^g_3 s^f_4 
s^f_1 s^f_2 s^f_1 s^f_2 s^b_3 s^a_4 s^a_1 
s^c_2 s^k_1 s^k_2 s^c_1 s^d_4 s^i_3 \\ \cdot s^i_4 
s^d_3 s^h_2 s^h_1 s^e_4 s^l_1 s^l_4 s^e_1 
s^b_2 s^f_3 s^f_2 s^b_3 s^a_4 s^a_1 s^c_2 
s^c_3 s^d_4 s^d_1 s^h_2 s^h_1 s^e_4 s^l_1 
s^l_4 s^e_1 s^b_2 s^b_1 s^a_4 s^a_1 s^c_2 
s^c_3 s^d_4 s^i_3 s^j_2 s^j_1 s^j_4 s^j_1 
s^j_4 s^j_3 s^i_2 s^i_1 s^i_4 s^i_1 s^i_4 
s^d_3 s^h_2 s^h_1 s^e_4 s^l_1 s^l_4 s^e_1 \\ \cdot 
s^b_2 s^f_3 s^f_2 s^b_3 s^a_4 s^a_1 s^c_2 
s^k_1 s^k_2 s^c_1 s^d_4 s^i_3 s^i_4 s^d_3 
s^h_2 s^h_1 s^e_4 s^l_1 s^l_4 s^e_1 s^b_2 
s^f_3 s^f_2 s^b_3 s^a_4 s^a_1 s^c_2 s^k_1 
s^k_2 s^c_1 s^d_4 s^i_3 s^i_4 s^d_3 s^h_2 
s^h_1 s^e_4 s^l_1 s^l_4 s^e_1 s^b_2 s^f_3 
s^f_2 s^f_1 s^f_2 s^b_3 s^a_4 s^a_1 s^c_2 
s^k_1 \\ \cdot s^k_2 s^c_1 s^d_4 s^d_1 s^h_2 s^h_1 
s^e_4 s^l_1 s^l_4 s^e_1 s^b_2 s^f_3 s^f_2 
s^b_3 s^a_4 s^a_1 s^c_2 s^k_1 s^k_2 s^c_1 
s^d_4 s^i_3 s^i_4 s^d_3 s^h_2 s^h_1 s^e_4 
s^l_1 s^l_4 s^e_1 s^b_2 s^f_3 s^f_2 s^b_3 
s^a_4 s^a_1 s^c_2 s^k_1 s^k_2 s^c_1 s^d_4 
s^i_3 s^i_4 s^i_1 s^i_4 s^d_3 s^h_2 s^h_1 
s^e_4 s^l_1 \\ \cdot s^l_4 s^e_1 s^b_2 s^f_3 s^f_2 
s^b_3 s^a_4 s^a_1 s^c_2 s^k_1 s^k_2 s^c_1 
s^d_4 s^i_3 s^i_4 s^d_3 s^h_2 s^h_1 s^e_4 
s^l_1 s^l_4 s^e_1 s^b_2 s^f_3 s^f_2 s^b_3 
s^a_4 s^a_1 s^c_2 s^k_1 s^k_2 s^c_1 s^d_4 
s^i_3 s^i_4 s^d_3 s^h_2 s^h_1 s^e_4 
 (s^e_3)
\end{array}
\endgroup
$}}}
\\ \quad
\\ \mbox{$\FigSixtySevenA$: Hamilton circuit of \cite[Table~3, Row~17]{Hec09}} 
\\ \mbox{missing 
$s^b_4$, $s^e_2$, $s^a_2$, $s^a_3$,  
$s^h_3$, $s^h_4$, $s^c_4$, $s^d_2$,
$s^k_4$ and $s^l_2$ 
(Length $=1440$)}
\end{array}
\end{equation*}

\noindent
{\bf{Rank-4-Case-3:}} 
Let $\chamchi$ be of $\acchi^{(3,18)}_1$.
Use the symbols in $\FigSixtyEightA$. 
For $t\in\{a,b,c,d\}$ of $\FigSixtyEightA$, 
let $t^\prime$ be the one obtained from $t$ by changing
$\al_1$ and $\al_4$.
Let $\chamchi_1:=\chamchi\in\chamcalXPi$.
Then  $\champi_\chamchi(\chamchi_1)=a$.
Let $\chamchi_2:=\chamtau_1\chamtau_2\chamtau_3\chamtau_4\chamchi\in\chamcalGchi$ and
$\chamchi_3:=\chamtau_1\chamtau_2\chamtau_3\chamtau_4\chamchi_2\in\chamcalGchi$.
Then $\champichi(\chamchi_2)=e$ and $\champichi(\chamchi_3)=a^\prime$.
By $\FigSixtyNineA$, 
$\{\chamchi_1,\chamchi_2,\chamchi_3\}$ is
a $1$-complete $\chambartau$-representative subset of $\chamcalGchi$.
Moreover Hamilton circuit maps of $\chamGamma((\chamchi_1)_{\langle 1\rangle})$,
$\chamGamma((\chamchi_2)_{\langle 1\rangle})$,
$\chamGamma((\chamchi_3)_{\langle 1\rangle})$
are drawn in by $\FigSeventyA$, 
$\FigSeventyOneA$, 
$\FigSeventyTwoA$ 
respectively. Then by Proposition~\ref{proposition:maintool}
and $\FigSeventyThreeA$, 
we see an existence of 
a Hamilton circuit map for $\chamchi$. In fact, it is written in
$\FigSeventyFourA$. 

\begin{equation*}
\begin{array}{l}
$\allHecFourEighteen$
\\ \quad
\\ \mbox{$\FigSixtyEightA$: Generalized Dynkin diagram of \cite[Table~3, Row~18]{Hec09}}
\end{array}
\end{equation*}


\begin{equation*}
\begin{array}{l}
$\RankFourHecEighteenTau$
\\ \quad
\\ \mbox{$\FigSixtyNineA$: Changing of diagrams of $\FigSixtyEightA$} 
\end{array}
\end{equation*}

\begin{equation*}
\begin{array}{l}
$\spTotalCHHecThreeSixForHecFourEighteen$
\\ \quad
\\ \mbox{$\FigSeventyA$: First part of Hamilton circuit for \cite[Table~3, Row~18]{Hec09}}
\end{array}
\end{equation*}

\begin{equation*}
\begin{array}{l}
$\setlength{\unitlength}{1mm}
\begin{picture}(110,40)(-30,00)
\put(0,0){$\TotalCHHecThreeEightForHecFourEighteen$}
\end{picture}$
\\ \quad
\\ \mbox{$\FigSeventyOneA$: Second part of Hamilton circuit for \cite[Table~3, Row~18]{Hec09}}
\end{array}
\end{equation*}

\begin{equation*}
\begin{array}{l}
$\setlength{\unitlength}{1mm}
\begin{picture}(110,40)(00,00)
\put(0,0){
$\spTotalCHHecThreeFiveForFourEughteen$}
\end{picture}$
\\ \quad
\\ \mbox{$\FigSeventyTwoA$: Third part of Hamilton circuit for \cite[Table~3, Row~18]{Hec09}}
\end{array}
\end{equation*}

\begin{equation*}
\begin{array}{l}
$\JointForHecFourEighteen$
\\ \quad
\\ \mbox{$\FigSeventyThreeA$: Joints of $\FigSeventyA$, 
$\FigSeventyOneA$, 
$\FigSeventyTwoA$} 
\end{array}
\end{equation*}

\begin{equation*}
\begin{array}{l}
{\footnotesize{\mbox{$\begin{array}{l}
(s^a_2) (s^a_3) (s^b_4) (s^c_3) (s^d_2)
  (s^d_3) (s^c_2) (s^b_3) 
 (s^a_4) s^a_3 s^a_2 s^b_4 
 (s^b_2) s^c_3 s^d_2 (s^e_4) (s^e_2) 
 (s^e_3) s^e_2 s^e_3 s^e_2 
 (s^d_4) s^d_3 s^e_4 s^e_2 s^d_4 s^c_2 
 (s^c_4) \\ \cdot s^b_3 s^a_4 s^a_2 s^b_4 s^c_3 s^c_4 
s^d_2 s^d_3 s^c_2 s^b_3 s^b_2 
 (s^b_1) \\ \cdot s^c_3 s^c_4 s^d_2 s^d_3 s^c_2 
s^b_3 s^b_2 s^a_4 s^a_2 s^a_3 s^b_4 s^c_3 
s^d_2 s^d_3 s^c_2 s^b_3 s^a_4 s^a_3 s^a_2 
s^b_4 s^b_2 s^c_3 s^d_2 s^e_4 s^e_2 s^e_3 
s^e_2 s^e_3 s^e_2 s^d_4 s^d_3 s^e_4 s^e_2 
s^d_4 s^c_2 s^c_4 s^b_3 s^b_1 s^c_3 s^c_4 \\ \cdot 
s^d_2 s^d_3 s^c_2 s^b_3 s^b_2 s^a_4 s^a_2 
s^a_3 s^b_4 s^c_3 s^d_2 s^d_3 s^c_2 s^b_3 
s^a_4 s^a_3 s^a_2 s^b_4 s^b_2 s^c_3 s^d_2 
s^e_4 s^e_2 s^e_3 s^e_2 s^e_3 s^e_2 s^d_4 
s^d_3 s^e_4 s^e_2 \\ \cdot (s^g_1) (s^g_3) (s^h_2) 
 (s^h_4) (s^i_3) (s^i_4) (s^i_2) 
 (s^h_3) s^h_4 \\ \cdot s^i_3 s^i_2 s^i_4 s^h_3 
s^h_4 (s^g_2) s^g_3 (s^f_4) (s^f_3) 
 (s^f_2) (s^g_4) s^g_3 
 (s^e_1) s^e_2 s^d_4 s^c_2 s^c_4 s^b_3 s^a_4 
s^a_2 s^b_4 s^c_3 s^c_4 s^d_2 s^d_3 s^c_2 
s^b_3 s^b_2 s^a_4 s^a_2 s^a_3 s^b_4 s^c_3 \\ \cdot 
s^d_2 s^d_3 s^c_2 s^b_3 s^a_4 s^a_3 s^a_2 \\ \cdot 
s^b_4 s^b_2 s^c_3 s^d_2 s^e_4 s^e_2 s^e_3 
s^e_2 s^g_1 s^f_4 s^f_3 s^f_2 s^g_4 s^g_3 
s^f_4 s^f_2 s^f_3 s^g_4 s^g_3 s^h_2 s^h_4 
s^i_3 s^i_4 s^i_2 s^h_3 s^h_4 s^i_3 s^i_2 
\\ \cdot (s^j_1) (s^j_2) (s^j_3) s^j_2 
 (s^j_4) s^j_2 s^j_4 s^j_2 s^j_3 s^j_2 s^j_3 
s^j_2 \\ \cdot s^j_3 s^j_2 s^j_4 s^j_2 s^j_3 s^j_2 
s^j_3 s^j_2 s^j_4 s^j_2 s^j_3 s^j_2 s^j_4 
s^j_2 s^j_3 s^j_2 s^j_4 s^j_2 s^j_4 s^j_2 
s^j_3 s^j_2 s^j_3 s^j_2 s^j_3 s^j_2 s^j_4 
s^j_2 s^j_3 s^j_2 s^j_3 s^j_2 s^j_4 
 (s^i_1) s^i_4 s^h_3 s^h_4 s^g_2 s^g_3 s^f_4 
\\ \cdot s^f_3 s^f_2 s^g_4 s^g_3 s^f_4 s^f_2 s^f_3 
s^g_4 s^g_3 s^h_2 
 (s^h_1) s^g_2 s^g_3 s^f_4 s^f_3 s^f_2 s^g_4 
s^g_3 s^f_4 s^f_2 s^f_3 s^g_4 s^g_3 s^h_2 
s^h_4 s^i_3 s^i_4 s^i_2 s^h_3 s^h_1 s^g_2 
s^g_3 s^f_4 s^f_3 s^f_2 s^g_4 s^g_3 s^f_4 
s^f_2 s^f_3 \\ \cdot s^g_4 s^g_3 s^h_2 s^h_4 s^i_3 
s^i_4 s^i_2 s^h_3 s^h_1 s^g_2 s^g_3 s^f_4 
s^f_3 s^f_2 s^g_4 s^g_3 s^f_4 s^f_2 s^f_3 
s^g_4 s^g_3 s^h_2 s^h_1 s^g_2 s^g_3 s^f_4 
s^f_3 s^f_2 s^g_4 s^g_3 s^e_1 s^e_2 s^d_4 
s^c_2 s^c_4 s^b_3 s^a_4 s^a_2 s^b_4 s^c_3 \\ \cdot 
s^c_4 s^d_2 s^d_3 s^c_2 s^b_3 s^b_2 s^a_4 
s^a_2 s^a_3 s^b_4 s^c_3 s^d_2 s^d_3 s^c_2 
s^b_3 s^a_4 s^a_3 s^a_2 s^b_4 s^b_2 s^b_1 
s^b_2 s^a_4 s^a_2 s^a_3 s^b_4 s^c_3 s^d_2 
s^d_3 s^c_2 s^b_3 s^a_4 s^a_3 s^a_2 s^b_4 
s^b_2 s^c_3 s^d_2 s^e_4 s^e_2 \\ \cdot s^e_3 s^e_2 
s^e_3 s^e_2 s^d_4 s^d_3 s^e_4 s^e_2 s^d_4 
s^c_2 s^c_4 s^b_3 s^a_4 s^a_2 s^b_4 s^c_3 
s^c_4 s^d_2 s^d_3 s^c_2 
 (s^c_1) s^d_2 s^e_4 s^e_2 s^e_3 s^e_2 s^e_3 
s^e_2 s^d_4 s^d_3 
 (s^f_1) s^f_2 s^f_3 s^g_4 s^g_3 s^h_2 s^h_4 
s^i_3 s^i_4 s^i_2 \\ \cdot s^h_3 s^h_1 s^g_2 s^g_3 
s^f_4 s^f_3 s^f_2 s^g_4 s^g_3 s^f_4 s^f_2 
s^f_3 s^g_4 s^g_3 s^h_2 s^h_4 s^i_3 s^i_4 
s^i_2 s^h_3 s^h_4 s^i_3 s^i_2 s^i_4 s^h_3 
s^h_1 s^i_3 s^i_2 s^i_4 s^h_3 s^h_1 s^i_3 
s^i_4 s^i_2 s^h_3 s^h_4 s^i_3 s^i_2 s^i_4 
s^h_3 \\ \cdot s^h_1 s^i_3 s^i_2 s^i_4 s^h_3 s^h_1 
s^i_3 s^i_2 s^i_4 s^h_3 s^h_1 s^i_3 s^i_4 
s^i_2 s^h_3 s^h_4 s^i_3 s^j_1 s^j_3 s^j_2 
s^i_1 s^h_3 s^h_4 s^g_2 s^e_1 s^e_2 s^d_4 
s^d_3 s^f_1 s^f_2 s^f_3 
 (s^d_1) s^c_2 s^c_4 s^b_3 s^a_4 s^a_2 
 (s^a_1) s^a_2 s^a_1
\end{array}
$}}} 
\\ \quad
\\ \mbox{$\FigSeventyFourA$: Special Hamilton circuit of \cite[Table~3, Row~18]{Hec09} obtained} 
\\ \mbox{from 
$\FigSeventyA$, 
$\FigSeventyOneA$, 
and
$\FigSeventyTwoA$ 
by joints mentioned by 
$\FigSixtyEightA$,} 
\\ \mbox{where $g:=d^\prime$, $h:=c^\prime$, $i:=b^\prime$ and $j:=a^\prime$
(Length $=480$)}
\end{array}
\end{equation*}



\noindent
{\bf{Rank-4-Case-4:}} 
Let $\chamchi\in\chamcalXPi$ be
that of $\acchi^{(3,22)}_1$.
Use the symbols in 
$\FigSeventyFiveA$ 
and $ \FigSeventySixA$. 
As shown in 
$\FigSeventyFiveA$, 
for $t\in\{a,b,c,d,{\bar{a}},{\bar{b}},{\bar{c}},{\bar{d}}\}$ of 
$\FigSeventyFiveA$, 
let $t^\prime$ mean the one obtained from $t$ by changing
$\al_2$ and $\al_4$.
Let $\chamchi_1:=\chamchi$.
Then  $\champichi(\chamchi)=a$.
By $\FigSeventySixA$, 
for $\chamchi^\prime$, $\chamchi^{\prime\prime}\in\chamcalGchi$
with $\chamtau_1(\chamchi^\prime)=\chamchi^{\prime\prime}$,
we have
\begin{equation}\label{eqn:bdoneHecFourTewlve}
(\champichi(\chamchi^\prime),\champichi(\chamchi^{\prime\prime}))
\in\{({\bar{a}},{\bar{a}}^\prime),({\bar{c}},{\bar{d}}^\prime),
({\bar{b}},{\bar{b}}^\prime),({\bar{d}},{\bar{c}}^\prime)\}.
\end{equation}
Let $\chamchi_1:=\chamchi$.
Then $\champichi(\chamchi_1)=a$.
Let $\chamchi_2$, (resp. $\chamchi_3$) $\in\chambarcalGchi$
be such that $\champichi(\chamchi_2)={\bar{a}}$
(resp. $\champichi(\chamchi_3)=a^\prime$).
By $\FigSeventySixA$, 
$\{\chamchi_1,\chamchi_2,\chamchi_3\}$ is 
a $1$-complete $\chambartau$-representative subset of $\chamcalGchi$.
Notice that the generalized Dynkin diagram of $\champi_{(\chamchi_1)_{\langle 1\rangle}}((\chamchi_1)_{\langle 1\rangle})$
(resp.~$\champi_{(\chamchi_3)_{\langle 1\rangle}}((\chamchi_3)_{\langle 1\rangle})$
is the same as $a$ of $\FigTwentyEightA$ 
with $\zeta$, $2$, $3$, $4$ (resp.~$\zeta$, $2$, $4$, $3$) in place of 
$\chamq$, $1$, $2$, $3$ respectively,
and that 
the generalized Dynkin diagram of $\champi_{(\chamchi_2)_{\langle 1\rangle}}((\chamchi_2)_{\langle 1\rangle})$
is the same as $a$ of $\FigThirtyOneA$. 
with $\zeta$, $2$, $3$, $4$ in place of 
$-\chamq=-r$, $1$, $2$, $3$ respectively,
By $\FigThirtyA$ 
and
$\FigThirtyTwoA$, 
we see that
for $t\in\fkJ_{1,3}$, there exists a special Hamilton map
of $\chamGamma(\champi_{(\chamchi_t)_{\langle 1\rangle}}((\chamchi_t)_{\langle 1\rangle}))$.
By an argument similar to that for
$\acchi^{(3,18)}_\bullet$,
we see an existence of 
a Hamilton circuit map for $\chamchi$. In fact, it is written in
$\FigEightyOneA$. 
Concerning $\FigSeventySevenA$, $\FigSeventyEightA$
and $\FigSeventyNineA$, the set of the `joints' (compare to $\FigSeventyThreeA$) is 
$\{({\bf 3},a,a,a,a)$, $({\bf 4},a,a,a,a)$, 
$({\bf 3},b,c,{\bar b},{\bar a})$, $({\bf 4},b,b,{\bar b},{\bar b})$,
$({\bf 4},c,d,{\bar a},{\bar d})$, $({\bf 3},d,d,{\bar d},{\bar d})$,   
$({\bf 3},{\bar c},{\bar c},{\bar c},{\bar c})$, $({\bf 3},{\bar c}^\prime,{\bar c}^\prime,{\bar c}^\prime,{\bar c}^\prime)$,  
$({\bf 2},{\bar d}^\prime,{\bar a}^\prime,d^\prime,c^\prime)$, $({\bf 3},{\bar d}^\prime,{\bar d}^\prime,d^\prime,d^\prime)$,  
$({\bf 3},{\bar a}^\prime,{\bar b}^\prime,c^\prime,b^\prime)$, 
\newline
$({\bf 2},{\bar b}^\prime,{\bar b}^\prime,b^\prime,b^\prime)$,   
$({\bf 2},a^\prime,a^\prime,a^\prime,a^\prime)$, $({\bf 3},a^\prime,a^\prime,a^\prime,a^\prime)\}$.

\begin{equation*}
\begin{array}{l}
$\allHecFourTewntytwo$
\\ \quad
\\ \mbox{$\FigSeventyFiveA$: Generalized Dynkin diagrams of \cite[Table~3, Row~22]{Hec09}}
\end{array}
\end{equation*}

\begin{equation*}
\begin{array}{l}
$\RankFourHecTwentytwoTau$
\\ \quad
\\ \mbox{$\FigSeventySixA$: Changing of diagrams of $\FigSeventyFiveA$}
\end{array}
\end{equation*}

\begin{equation*}
\begin{array}{l}
$\setlength{\unitlength}{1mm}
\begin{picture}(200,50)(-15,00)
\put(00,0){$\TotalCHHecThreeSevenForHCHecfirstpartFourTwentytwo$}
\end{picture}$
\\ \quad
\\ \mbox{$\FigSeventySevenA$: First part for Hamilton circuit of \cite[Table~3, Row~22]{Hec09}}
\end{array}
\end{equation*}

\begin{equation*}
\begin{array}{l}
$\setlength{\unitlength}{1mm}
\begin{picture}(140,45)(-5,0)
\put(00,0){$\TotalCHHecThreeNineForHecFourTwetytwo$}
\put(70,0){$\TotalCHHecThreeNineForHecFourTwetytwoDash$}
\end{picture}$
\\ \quad
\\ \mbox{$\FigSeventyEightA$: Second and Third parts for Hamilton circuit of \cite[Table~3, Row~22]{Hec09}}
\end{array}
\end{equation*}

\begin{equation*}
\begin{array}{l}
$\setlength{\unitlength}{1mm}
\begin{picture}(200,50)(-15,00)
\put(00,0){$\TotalCHHecThreeSevenForHCHecfourthpartFourTwentytwo$}
\end{picture}$
\\ \quad
\\ \mbox{$\FigSeventyNineA$: Fourth part for Hamilton circuit of \cite[Table~3, Row~22]{Hec09}}
\end{array}
\end{equation*}

\begin{equation*}
\begin{array}{l}
$\JointForHecTwentyTwo$
\\ \quad
\\ \mbox{$\FigEightyA$: Joints of $\FigSeventySevenA$, 
$\FigSeventyEightA$, 
$\FigSeventyNineA$} 
\end{array}
\end{equation*}

\begin{equation*}
\begin{array}{l}
{\tiny{\mbox{$
\begingroup
\renewcommand{\arraystretch}{1.5}
\begin{array}{l}
S^\prime:=(s^a_3) (s^b_2) (s^b_4) 
 (s^a_2) s^a_3 s^b_2 s^b_4 s^a_2 s^a_3 s^b_2 
s^b_4 s^a_2 s^a_3 s^b_2 s^b_4 s^a_2 s^a_3 
s^b_2 (s^c_3) (s^c_2) (s^d_4) (s^d_2) 
 (s^c_4) s^c_2 (s^b_3) s^b_4 s^c_3 s^d_4 
 (s^d_3) \\ \cdot s^d_2 s^d_3 s^c_4 s^b_3 s^a_2 s^a_3 
s^b_2 s^c_3 s^d_4 s^d_3 s^d_2 s^d_3 s^c_4 
s^b_3 s^a_2 s^a_3 
 (s^a_4) s^a_3 s^b_2 s^b_4 s^a_2 \\ \cdot s^a_3 
s^b_2 s^b_4 s^a_2 s^a_3 s^b_2 s^b_4 s^a_2 
s^a_3 s^b_2 s^b_4 s^a_2 s^a_3 s^b_2 s^c_3 
s^c_2 s^d_4 s^d_3 s^c_4 s^c_2 s^d_4 s^d_3 
s^c_4 s^c_2 s^d_4 s^d_3 s^c_4 s^c_2 s^d_4 
s^d_2 s^c_4 s^c_2 s^d_4 s^d_3 s^c_4 s^c_2 
s^d_4 s^d_3 s^c_4 s^c_2 s^d_4 s^d_3 s^c_4 
s^c_2 \\ \cdot (s^f_1) (s^h_2) (s^h_3) (s^j_4) 
 (s^j_3) (s^k_2) \\ \cdot (s^l_3) (s^l_2) 
 (s^k_3) (s^i_4) (s^i_3) (s^g_2) (s^g_3)
  (s^f_4) (s^e_3) (s^e_4) 
 (s^f_3) s^h_2 s^h_3 s^j_4 s^j_3 s^k_2 s^l_3 
s^l_2 s^k_3 s^i_4 
 (s^i_1) s^g_2 s^g_3 s^f_4 s^e_3 s^e_4 s^f_3 
s^h_2 s^h_3 s^j_4 s^j_3 s^k_2 s^l_3 s^l_2 
s^k_3 s^i_4 \\ \cdot s^i_3 s^g_2 s^g_3 s^f_4 s^e_3 
s^e_4 s^f_3 s^h_2 s^h_3 s^j_4 s^j_3 s^k_2 
s^l_3 s^l_2 \\ \cdot (s^o_1) (s^p_4) (s^p_3) 
(s^o_4) (s^o_2)
s^p_4 s^p_3 s^o_4 s^o_2 s^p_4 s^p_3 s^o_4 
s^o_2 s^p_4 s^p_3 s^o_4 s^o_2 s^p_4 s^p_3 
s^o_4 (s^n_3) (s^n_4)
(s^m_2) (s^m_4) (s^n_2) s^n_4 
(s^o_3) s^o_2 s^n_3 s^m_2 (s^m_3) s^m_4  \\ \cdot s^m_3 
s^n_2 s^o_3 s^p_4 s^p_3 s^o_4 s^n_3 s^m_2 
s^m_3 s^m_4 s^m_3 s^n_2 s^o_3 s^p_4 s^p_3 
(s^p_2) s^p_3 s^o_4 \\ \cdot s^o_2 s^p_4 s^p_3 s^o_4 
s^o_2 s^p_4 s^p_3 s^o_4 s^o_2 s^p_4 s^p_3 
s^o_4 s^o_2 s^p_4 s^p_3 s^o_4 s^n_3 s^n_4 
s^m_2 s^m_3 s^n_2 s^n_4 s^m_2 s^m_3 s^n_2 
s^n_4 s^m_2 s^m_3 s^n_2 s^n_4 s^m_2 s^m_4 
s^n_2 s^n_4 s^m_2 s^m_3 s^n_2 s^n_4 s^m_2 
s^m_3 \\ \cdot (s^j_1) s^j_3 (s^h_4) s^h_3 
 (s^f_2) s^e_3 (s^e_2) s^f_3 
 (s^g_4) s^g_3 \\ \cdot (s^i_2) s^i_3 
 (s^k_4) s^l_3 (s^l_4) s^k_3 
 (s^j_2) s^j_3 s^h_4 s^h_3 s^f_2 s^e_3 s^e_2 
s^f_3 s^g_4 s^g_3 s^i_2 s^i_3 s^k_4 (s^n_1)
s^n_4 s^m_2 s^m_3 s^n_2 s^n_4 s^m_2 s^m_3 
s^n_2 s^n_4 s^m_2 s^m_3 s^n_2 s^n_4 s^m_2 
s^m_4 s^n_2 s^n_4 s^m_2 s^m_3 s^n_2 s^n_4  \\ \cdot 
s^m_2 s^m_3 s^n_2 s^n_4 s^m_2 s^m_3 s^n_2 
s^n_4 s^o_3 \\ \cdot s^p_4 s^p_3 s^o_4 s^o_2 s^p_4 
s^p_3 s^o_4 s^o_2 s^p_4 s^p_3 s^o_4 s^o_2 
s^p_4 s^p_3 s^o_4 s^o_2 s^p_4 s^p_3 s^o_4 
s^n_3 s^n_4 s^m_2 s^m_4 s^n_2 s^n_4 s^o_3 
s^o_2 s^n_3 s^m_2 s^m_3 s^m_4 s^m_3 s^n_2 
s^o_3 s^p_4 s^p_3 s^o_4 s^n_3 s^m_2 s^m_3 
s^m_4 s^m_3 s^n_2 s^o_3 \\ \cdot s^p_4 s^p_3 s^p_2 
s^p_3 s^o_4 s^o_2 \\ \cdot s^p_4 s^p_3 s^o_4 
 (s^l_1) s^k_3 s^i_4 s^i_3 s^g_2 s^g_3 s^f_4 
s^e_3 s^e_4 s^f_3 s^h_2 s^h_3 s^j_4 s^j_3 
s^k_2 s^l_3 s^l_2 s^k_3 s^i_4 s^i_3 s^g_2 
s^g_3 s^f_4 s^e_3 s^e_4 s^f_3 s^h_2 
 (s^h_1) s^j_4 s^j_3 s^k_2 s^l_3 s^l_2 s^k_3 
s^i_4 s^i_3 s^g_2 s^g_3 s^f_4 s^e_3 s^e_4 
s^f_3 s^h_2 s^h_3 s^j_4 s^j_3 s^k_2 \\ \cdot s^l_3 
s^l_2 s^k_3 s^i_4 s^i_3 s^g_2 s^g_3 s^f_4 
s^e_3 s^e_4 
 (s^b_1) s^a_2 s^a_3 s^b_2 s^b_4 s^a_2 s^a_3 
s^b_2 s^b_4 s^a_2 s^a_3 s^b_2 s^b_4 s^a_2 
s^a_3 s^b_2 s^b_4 s^a_2 s^a_3 s^b_2 s^c_3 
s^c_2 s^d_4 s^d_2 s^c_4 s^c_2 s^b_3 s^b_4 
s^c_3 s^d_4 s^d_3 s^d_2 s^d_3 s^c_4 s^b_3 
s^a_2 s^a_3 s^b_2 s^c_3 s^d_4 \\ \cdot s^d_3 s^d_2 
s^d_3 s^c_4 s^b_3 s^a_2 s^a_3 s^a_4 s^a_3 
s^b_2 s^b_4 s^a_2 s^a_3 s^b_2 s^b_4 s^a_2 
s^a_3 s^b_2 s^b_4 s^a_2 s^a_3 s^b_2 s^b_4 
s^a_2 s^a_3 s^b_2 s^c_3 s^c_2 s^d_4 s^d_3 
s^c_4 s^c_2 s^d_4 s^d_3 s^c_4 s^c_2 s^d_4 
s^d_3 s^c_4 s^c_2 s^d_4 s^d_2 s^c_4 s^c_2 
s^d_4 s^d_3 s^c_4 s^c_2 s^d_4 s^d_3 \\ \cdot 
 (s^g_1) s^g_3 s^i_2 s^i_3 s^k_4 s^l_3 s^l_4 
s^k_3 s^j_2 s^j_3 s^h_4 s^h_3 s^f_2 s^e_3 
s^e_2 s^f_3 s^g_4 s^g_3 s^i_2 s^i_3 s^k_4 
s^l_3 s^l_4 s^k_3 s^j_2 s^j_3 s^h_4 s^h_3 
s^f_2 (s^c_1) s^c_2 s^d_4 s^d_3 s^c_4 s^c_2 
s^d_4 s^d_3 s^c_4 s^c_2 s^d_4 s^d_3 s^c_4 
s^c_2 s^d_4 s^d_2 s^c_4 s^c_2 s^d_4 s^d_3 
s^c_4 \\ \cdot s^c_2 s^d_4 s^d_3 s^c_4 s^c_2 s^d_4 
s^d_3 s^c_4 s^c_2 s^b_3 s^a_2 s^a_3 s^b_2 
s^b_4 s^a_2 s^a_3 s^b_2 s^b_4 s^a_2 s^a_3 
s^b_2 s^b_4 s^a_2 s^a_3 s^b_2 s^b_4 s^a_2 
s^a_3 s^b_2 s^c_3 s^c_2 s^d_4 s^d_2 s^c_4 
s^c_2 s^b_3 s^b_4 s^c_3 s^d_4 s^d_3 s^d_2 
s^d_3 s^c_4 s^b_3 s^a_2 s^a_3 s^b_2 s^c_3 
s^d_4 s^d_3 \\ \cdot s^d_2 s^d_3 s^c_4 s^b_3 s^a_2 
s^a_3 s^a_4 s^a_3 s^b_2 s^b_4 s^a_2 s^a_3 
s^b_2 (s^e_1) s^f_3 s^h_2 s^h_3 s^j_4 s^j_3 
s^k_2 s^l_3 s^l_2 s^k_3 s^i_4 s^i_3 s^g_2 
s^g_3 s^f_4 s^e_3 s^e_4 s^f_3 s^h_2 s^h_3 
s^j_4 s^j_3 s^k_2 s^l_3 s^l_2 s^k_3 s^i_4 
s^i_1 s^g_2 s^g_3 s^f_4 s^e_3 s^e_4 s^f_3 
s^h_2 s^h_3 s^j_4 \\ \cdot s^j_3 s^k_2 s^l_3 s^l_2 
s^k_3 s^i_4 s^i_3 s^g_2 s^g_3 s^f_4 s^e_3 
s^e_4 s^f_3 s^h_2 s^h_3 s^j_4 s^j_3 s^k_2 
s^l_3 s^l_2 s^o_1 s^p_4 s^p_3 s^o_4 s^o_2 
s^p_4 s^p_3 s^o_4 s^o_2 s^p_4 s^p_3 s^o_4 
s^o_2 s^p_4 s^p_3 s^o_4 s^o_2 s^p_4 s^p_3 
s^o_4 s^n_3 s^n_4 s^m_2 s^m_4 s^n_2 s^n_4 
s^o_3 s^o_2 s^n_3 s^m_2 \\ \cdot s^m_3 s^m_4 s^m_3 
s^n_2 s^o_3 s^p_4 s^p_3 s^o_4 s^n_3 s^m_2 
s^m_3 s^m_4 s^m_3 s^n_2 s^o_3 s^p_4 s^p_3 
s^p_2 s^p_3 s^o_4 s^o_2 s^p_4 s^p_3 s^o_4 
s^o_2 s^p_4 s^p_3 s^o_4 s^o_2 s^p_4 s^p_3 
s^o_4 s^o_2 s^p_4 s^p_3 s^o_4 s^n_3 s^n_4 
s^m_2 s^m_3 s^n_2 s^n_4 s^m_2 \\ \cdot s^m_3 s^n_2 
s^n_4 s^m_2 s^m_3 s^n_2 s^n_4 \\ \cdot s^m_2 s^m_4 
s^n_2 s^n_4 s^m_2 s^m_3 s^n_2 s^n_4 s^m_2 
s^m_3 s^j_1 s^j_3 s^h_4 s^h_3 s^f_2 s^e_3 
s^e_2 s^f_3 s^g_4 s^g_3 s^i_2 s^i_3 s^k_4 
s^l_3 s^l_4 s^k_3 s^j_2 s^j_3 s^h_4 s^h_3 
s^f_2 s^e_3 s^e_2 s^f_3 s^g_4 s^g_3 s^i_2 
s^i_3 s^k_4 s^n_1 s^n_4 s^m_2  \\ \cdot s^m_3 s^n_2 
s^n_4 s^m_2 s^m_3 s^n_2 s^n_4 s^m_2 \\ \cdot s^m_3 
s^n_2 s^n_4 s^m_2 s^m_4 s^n_2 s^n_4 s^m_2 
s^m_3 s^n_2 s^n_4 s^m_2 s^m_3 s^n_2 s^n_4 
s^m_2 s^m_3 s^n_2 s^n_4 s^o_3 s^p_4 s^p_3 
s^o_4 s^o_2 s^p_4 s^p_3 s^o_4 s^o_2 s^p_4 
s^p_3 s^o_4 s^o_2 s^p_4 s^p_3 s^o_4 s^o_2 
s^p_4 s^p_3 s^o_4 s^n_3 s^n_4 s^m_2 \\ \cdot s^m_4 
s^n_2 s^n_4 s^o_3 s^o_2 s^n_3 s^m_2 s^m_3 \\ \cdot 
s^m_4 s^m_3 s^n_2 s^o_3 s^p_4 s^p_3 s^o_4 
s^n_3 s^m_2 s^m_3 s^m_4 s^m_3 s^n_2 s^o_3 
s^p_4 s^p_3 s^p_2 s^p_3 s^o_4 s^o_2 s^p_4 
s^p_3 s^o_4 s^l_1 s^k_3 s^i_4 s^i_3 s^g_2 
s^g_3 s^f_4 s^e_3 s^e_4 s^f_3 s^h_2 s^h_3 
s^j_4 s^j_3 s^k_2 s^l_3 s^l_2 s^k_3 s^i_4 
s^i_3 s^g_2 s^g_3 s^f_4 s^e_3 s^e_4 s^f_3 
s^h_2 \\ \cdot s^h_1 s^j_4 s^j_3 s^k_2 s^l_3 s^l_2 
s^k_3 s^i_4 s^i_1 s^g_2 s^g_3 s^f_4 s^e_3 
s^e_4 s^f_3 s^h_2 s^h_3 s^j_4 s^j_3 s^k_2 
s^l_3 s^l_2 s^k_3 s^i_4 s^i_3 s^g_2 s^g_3 
s^f_4 s^e_3 s^e_4 s^f_3 s^h_2 s^h_3 s^j_4 
s^j_3 s^k_2 s^l_3 s^l_2 s^k_3 s^i_4 s^i_1 
s^g_2 s^g_3 s^f_4 s^e_3 s^e_4 s^f_3 s^h_2 
s^h_3 s^j_4 \\ \cdot s^j_3 s^k_2 s^l_3 s^l_2 s^k_3 
s^i_4 s^i_3 s^g_2 s^g_3 s^f_4 s^e_3 s^e_4 
s^b_1 s^a_2 s^a_3 s^b_2 s^b_4 s^a_2 s^a_3 
s^b_2 s^b_4 s^a_2 s^a_3 s^b_2 s^b_4 s^a_2 
s^a_3 s^b_2 s^b_4 s^a_2 s^a_3 s^b_2 s^c_3 
s^c_2 s^d_4 s^d_2 s^c_4 s^c_2 s^b_3 s^b_4 
s^c_3 s^d_4 s^d_3 s^d_2 s^d_3 s^c_4 s^b_3 
s^a_2 s^a_3 s^b_2 \\ \cdot s^c_3 s^d_4 s^d_3 s^d_2 
s^d_3 s^c_4 s^b_3 s^a_2 s^a_3 s^a_4 s^a_3 
s^b_2 s^b_4 s^a_2 s^a_3 s^b_2 s^b_4 s^a_2 
s^a_3 s^b_2 s^b_4 s^a_2 s^a_3 s^b_2 s^b_4 
s^a_2 s^a_3 s^b_2 s^c_3 s^c_2 s^d_4 s^d_3 
s^c_4 s^c_2 s^d_4 s^d_3 s^c_4 s^c_2 s^d_4 
s^d_3 s^c_4 s^c_2 s^d_4 s^d_2 s^c_4 s^c_2 
s^d_4 s^d_3 s^c_4 s^c_2 \\ \cdot s^d_4 s^d_3 s^g_1 
s^g_3 s^i_2 s^i_3 s^k_4 s^l_3 s^l_4 s^k_3 
s^j_2 s^j_3 s^h_4 s^h_3 s^f_2 s^e_3 s^e_2 
s^f_3 s^g_4 s^g_3 s^i_2 s^i_3 s^k_4 s^n_1 
s^n_4 s^m_2 s^m_3 s^n_2 s^n_4 s^m_2 s^m_3 
s^n_2 s^n_4 s^m_2 s^m_3 s^n_2 s^n_4 s^m_2 
s^m_4 s^n_2 s^n_4 s^m_2 s^m_3 s^n_2 s^n_4  \\ \cdot 
s^m_2 s^m_3 s^n_2 s^n_4 s^m_2 
\end{array}
\endgroup$}}}
\\ \quad
\\ \mbox{$\FigEightyOneA^\prime$}
\end{array}
\end{equation*}

\begin{equation*}
\begin{array}{l}
{\tiny{\mbox{$
\begingroup
\renewcommand{\arraystretch}{1.5}
\begin{array}{l}
S^\prime\cdot s^m_3 s^n_2 
s^n_4 s^o_3 s^p_4 s^p_3 s^o_4 s^o_2 s^p_4 
s^p_3 s^o_4 s^o_2 s^p_4 s^p_3 s^o_4 s^o_2 
s^p_4 s^p_3 s^o_4 s^o_2 s^p_4 s^p_3 s^o_4 
s^n_3 s^n_4 s^m_2 s^m_4 s^n_2 s^n_4 s^o_3 
s^o_2 s^n_3 s^m_2 s^m_3 s^m_4 s^m_3 s^n_2 
s^o_3 s^p_4 s^p_3 s^o_4 s^n_3 s^m_2 s^m_3  \\ \cdot 
s^m_4 s^m_3 s^n_2 s^o_3 s^p_4 s^p_3 \\ \cdot s^p_2 
s^p_3 s^o_4 s^o_2 s^p_4 s^p_3 s^o_4 s^o_2 
s^p_4 s^p_3 s^o_4 s^o_2 s^p_4 s^p_3 s^o_4 
s^l_1 s^k_3 s^i_4 s^i_3 s^g_2 s^g_3 s^f_4 
s^e_3 s^e_4 s^f_3 s^h_2 s^h_3 s^j_4 s^j_3 
s^k_2 s^l_3 s^l_2 s^k_3 s^i_4 s^i_3 s^g_2 
s^g_3 s^f_4 s^e_3 s^e_4 s^f_3 s^h_2 s^h_1 
s^j_4 s^j_3 s^k_2 s^l_3 s^l_2 s^k_3 s^i_4 \\ \cdot 
s^i_3 s^g_2 s^g_3 s^f_4 s^e_3 s^e_4 s^f_3 
s^h_2 s^h_3 s^j_4 s^j_3 s^k_2 s^l_3 s^l_2 
s^k_3 s^i_4 s^i_1 s^g_2 s^g_3 s^f_4 s^e_3 
s^e_4 s^f_3 s^h_2 s^h_3 s^j_4 s^j_3 s^k_2 
s^l_3 s^l_2 s^k_3 s^i_4 s^i_3 s^g_2 s^g_3 
s^f_4 s^e_3 s^e_4 s^f_3 s^h_2 s^h_1 s^j_4 
s^j_3 s^k_2 s^l_3 s^l_2 s^k_3 s^i_4 s^i_1 
s^g_2 \\ \cdot s^g_3 s^f_4 s^e_3 s^e_4 s^f_3 s^h_2 
s^h_3 s^j_4 s^j_3 s^k_2 s^l_3 s^l_2 s^k_3 
s^i_4 s^i_3 s^g_2 s^g_3 s^f_4 s^e_3 s^e_4 
s^b_1 s^a_2 s^a_3 s^b_2 s^b_4 s^a_2 s^a_3 
s^b_2 s^b_4 s^a_2 s^a_3 s^b_2 s^b_4 s^a_2 
s^a_3 s^b_2 s^b_4 s^a_2 s^a_3 s^b_2 s^c_3 
s^c_2 s^d_4 s^d_2 s^c_4 s^c_2 s^b_3 s^b_4 
s^c_3 s^d_4 \\ \cdot s^d_3 s^d_2 s^d_3 s^c_4 s^b_3 
s^a_2 s^a_3 s^b_2 s^c_3 s^d_4 s^d_3 s^d_2 
s^d_3 s^c_4 s^b_3 s^a_2 s^a_3 s^a_4 s^a_3 
s^b_2 s^b_4 s^a_2 s^a_3 s^b_2 s^b_4 s^a_2 
s^a_3 s^b_2 s^b_4 s^a_2 s^a_3 s^b_2 s^b_4 
s^a_2 s^a_3 s^b_2 s^c_3 s^c_2 s^d_4 s^d_3 
s^c_4 s^c_2 s^d_4 s^d_3 s^c_4 s^c_2 s^d_4 
s^d_3 s^c_4 s^c_2 \\ \cdot s^d_4 s^d_2 s^c_4 s^c_2 
s^d_4 s^d_3 s^c_4 s^c_2 s^d_4 s^d_3 s^g_1 
s^g_3 s^i_2 s^i_3 s^k_4 s^l_3 s^l_4 s^k_3 
s^j_2 s^j_3 s^h_4 s^h_1 s^j_4 s^j_3 s^k_2 
s^l_3 s^l_2 s^k_3 s^i_4 s^i_3 s^g_2 s^g_3 
s^f_4 s^e_3 s^e_4 s^f_3 s^h_2 s^h_3 s^j_4 
s^j_3 s^k_2 s^l_3 s^l_2 s^k_3 s^i_4 s^i_1 
s^g_2 s^g_3 s^f_4 s^e_3 \\ \cdot s^e_4 s^f_3 s^h_2 
s^h_3 s^j_4 s^j_3 s^k_2 s^l_3 s^l_2 s^k_3 
s^i_4 s^i_3 s^g_2 s^g_3 s^f_4 s^e_3 s^e_4 
s^f_3 s^h_2 s^h_1 s^j_4 s^j_3 s^k_2 s^l_3 
s^l_2 s^k_3 s^i_4 s^i_3 s^g_2 s^g_3 s^f_4 
s^e_3 s^e_4 s^f_3 s^h_2 s^h_1 s^j_4 s^j_3 
s^k_2 s^l_3 s^l_2 s^k_3 s^i_4 s^i_3 s^g_2 
s^g_3 s^f_4 s^e_3 s^e_4 s^f_3 \\ \cdot s^h_2 s^h_1 
s^j_4 s^j_3 s^k_2 s^l_3 s^l_2 s^k_3 s^i_4 
s^i_3 s^g_2 s^g_3 s^f_4 s^e_3 s^e_4 s^f_3 
s^h_2 s^h_1 s^j_4 s^j_3 s^k_2 s^l_3 s^l_2 
s^k_3 s^i_4 s^i_3 s^g_2 s^g_3 s^f_4 s^e_3 
s^e_4 s^f_3 s^h_2 s^h_3 s^j_4 s^j_3 s^k_2 
s^l_3 s^l_2 s^k_3 s^i_4 s^i_1 s^g_2 s^g_3 
s^f_4 s^e_3 s^e_4 s^f_3 s^h_2 s^h_3 \\ \cdot s^j_4 
s^j_3 s^k_2 s^l_3 s^l_2 s^o_1 s^p_4 s^p_3 
s^o_4 s^o_2 s^p_4 s^p_3 s^o_4 s^o_2 s^p_4 
s^p_3 s^o_4 s^o_2 s^p_4 s^p_3 s^o_4 s^o_2 
s^p_4 s^p_3 s^o_4 s^n_3 s^n_4 s^m_2 s^m_4 
s^n_2 s^n_4 s^o_3 s^o_2 s^n_3 s^m_2 s^m_3 
s^m_4 s^m_3 s^n_2 s^o_3 s^p_4 s^p_3 s^o_4  \\ \cdot 
s^n_3 s^m_2 s^m_3 s^m_4 s^m_3 s^n_2 s^o_3 \\ \cdot 
s^p_4 s^p_3 s^p_2 s^p_3 s^o_4 s^o_2 s^p_4 
s^p_3 s^o_4 s^o_2 s^p_4 s^p_3 s^o_4 s^o_2 
s^p_4 s^p_3 s^o_4 s^o_2 s^p_4 s^p_3 s^o_4 
s^n_3 s^n_4 s^m_2 s^m_3 s^n_2 s^n_4 s^m_2 
s^m_3 s^n_2 s^n_4 s^m_2 s^m_3 s^n_2 s^n_4 
s^m_2 s^m_4 s^n_2 s^n_4 s^m_2 s^m_3 s^j_1 \\ \cdot 
s^j_3 s^h_4 s^h_3 s^f_2 s^e_3 s^e_2 s^f_3 
s^g_4 \\ \cdot s^g_3 s^i_2 s^i_1 s^g_2 s^g_3 s^f_4 
s^e_3 s^e_4 s^f_3 s^h_2 s^h_3 s^j_4 s^j_3 
s^k_2 s^l_3 s^l_2 s^k_3 s^i_4 s^i_3 s^g_2 
s^g_3 s^f_4 s^e_3 s^e_4 s^f_3 s^h_2 s^h_1 
s^j_4 s^j_3 s^k_2 s^l_3 s^l_2 s^k_3 s^i_4 
s^i_3 s^g_2 s^g_3 s^f_4 s^e_3 s^e_4 s^b_1 
s^a_2 s^a_3 s^b_2 s^b_4 s^a_2 s^a_3 s^b_2 
s^b_4 s^a_2 \\ \cdot s^a_3 s^b_2 s^b_4 s^a_2 s^a_3 
s^b_2 s^b_4 s^a_2 s^a_3 s^b_2 s^c_3 s^c_2 
s^d_4 s^d_2 s^c_4 s^c_2 s^b_3 s^b_4 s^c_3 
s^d_4 s^d_3 s^d_2 s^d_3 s^c_4 s^b_3 s^a_2 
s^a_3 s^b_2 s^c_3 s^d_4 s^d_3 s^d_2 s^d_3 
s^c_4 s^b_3 s^a_2 s^a_3 s^a_4 s^a_3 s^b_2 
s^b_4 s^a_2 s^a_3 s^b_2 s^b_4 s^a_2 s^a_3 
s^b_2 s^b_4 s^a_2 \\ \cdot s^a_3 s^b_2 s^b_4 s^a_2 
s^a_3 s^b_2 s^e_1 s^e_4 s^f_3 s^h_2 s^h_3 
s^j_4 s^j_3 s^k_2 s^l_3 s^l_2 s^k_3 s^i_4 
s^i_1 s^g_2 s^g_3 s^f_4 s^e_3 s^e_4 s^f_3 
s^h_2 s^h_3 s^j_4 s^j_3 s^k_2 s^l_3 s^l_2 
s^k_3 s^i_4 s^i_3 s^g_2 s^g_3 s^f_4 s^e_3 
s^e_4 s^f_3 s^h_2 s^h_3 s^j_4 s^j_3 s^k_2 
s^l_3 s^l_2 s^k_3 s^i_4 \\ \cdot s^i_1 s^g_2 s^g_3 
s^f_4 s^e_3 s^e_4 s^f_3 s^h_2 s^h_3 s^j_4 
s^j_3 s^k_2 s^l_3 s^l_2 s^k_3 s^i_4 s^i_3 
s^g_2 s^g_3 s^f_4 s^c_1 s^c_2 s^d_4 s^d_3 
s^c_4 s^c_2 s^d_4 s^d_3 s^c_4 s^c_2 s^d_4 
s^d_3 s^c_4 s^c_2 s^d_4 s^d_2 s^c_4 s^c_2 
s^d_4 s^d_3 s^c_4 s^c_2 s^d_4 s^d_3 s^g_1 
s^g_3 s^i_2 s^i_3 s^k_4 s^l_3 \\ \cdot s^l_4 s^k_3 
s^j_2 s^j_3 s^h_4 s^h_1 s^j_4 s^j_3 s^k_2 
s^l_3 s^l_2 s^k_3 s^i_4 s^i_3 s^g_2 s^g_3 
s^f_4 s^e_3 s^e_4 s^f_3 s^h_2 s^h_1 s^j_4 
s^j_3 s^k_2 s^l_3 s^l_2 s^k_3 s^i_4 s^i_3 
s^g_2 s^g_3 s^f_4 s^e_3 s^e_4 s^f_3 s^h_2 
s^h_3 s^j_4 s^j_3 s^k_2 s^l_3 s^l_2 s^k_3 
s^i_4 s^i_3 s^g_2 s^g_3 s^f_4 s^e_3 \\ \cdot s^e_4 
s^f_3 s^h_2 s^h_1 s^j_4 s^j_3 s^k_2 s^l_3 
s^l_2 s^k_3 s^i_4 s^i_3 s^g_2 s^g_3 s^f_4 
s^e_3 s^e_4 s^f_3 s^h_2 s^h_1 s^f_2 s^e_3 
s^e_2 s^f_3 s^g_4 s^g_3 s^i_2 s^i_3 s^k_4 
s^l_3 s^l_4 s^k_3 s^j_2 s^j_3 s^h_4 s^h_3 
s^f_2 s^e_3 s^e_2 s^f_3 s^c_1 s^c_2 s^d_4 
s^d_3 s^c_4 s^c_2 s^f_1 s^h_2 s^h_3 s^j_4 \\ \cdot 
s^j_3 s^k_2 s^l_3 s^l_2 s^k_3 s^i_4 s^i_1 
s^g_2 s^g_3 s^f_4 s^e_3 s^e_4 s^f_3 s^h_2 
s^h_3 s^j_4 s^j_3 s^k_2 s^l_3 s^l_2 s^k_3 
s^i_4 s^i_1 s^g_2 s^g_3 s^f_4 s^e_3 s^e_4 
s^f_3 s^h_2 s^h_3 s^j_4 s^j_3 s^k_2 s^l_3 
s^l_2 s^k_3 s^i_4 s^i_3 s^g_2 s^g_3 s^f_4 
s^e_3 s^e_4 s^f_3 s^h_2 s^h_3 s^j_4 s^j_3 
s^k_2 \\ \cdot s^l_3 s^l_2 s^k_3 s^i_4 s^i_1 s^g_2 
s^g_3 s^f_4 s^e_3 s^e_4 s^f_3 s^h_2 s^h_3 
s^j_4 s^j_3 s^k_2 s^l_3 s^l_2 s^k_3 s^i_4 
s^i_1 s^g_2 s^g_3 s^f_4 s^e_3 s^e_4 s^f_3 
s^h_2 s^h_1 s^j_4 s^j_3 s^k_2 s^l_3 s^l_2 
s^k_3 s^i_4 s^i_1 s^k_4 s^l_3 s^l_4 s^k_3 
s^j_2 s^j_3 s^h_4 s^h_3 s^f_2 s^e_3 s^e_2 
s^f_3 s^g_4 \\ \cdot s^g_3 s^i_2 s^i_3 s^k_4 s^l_3 
s^l_4 s^k_3 s^n_1 s^n_4 s^m_2 s^m_3 s^n_2 
s^n_4 s^m_2 s^m_3 s^n_2 s^n_4 
  (s^k_1) s^i_4 s^i_3 s^g_2 s^g_3 s^f_4 s^e_3 
s^e_4 s^f_3 s^h_2 s^h_3 s^j_4 s^j_3 s^k_2 
s^l_3 s^l_2 s^k_3 s^i_4 s^i_1 s^g_2 s^g_3 
s^f_4 s^e_3 s^e_4 s^f_3 s^h_2 s^h_1 s^j_4  \\ \cdot 
s^j_3 s^k_2 s^l_3 s^l_2 s^k_3 \\ \cdot s^i_4 s^i_3 
s^g_2 s^g_3 s^f_4 s^e_3 s^e_4 s^f_3 s^h_2 
s^h_1 s^j_4 s^j_3 s^k_2 s^l_3 s^l_2 s^k_3 
s^i_4 s^i_3 s^g_2 s^g_3 s^f_4 s^e_3 s^e_4 
s^f_3 s^h_2 s^h_1 s^j_4 s^j_3 s^k_2 s^l_3 
s^l_2 s^k_3 s^i_4 s^i_3 s^g_2 s^g_3 s^f_4 
s^e_3 s^e_4 s^f_3 s^h_2 s^h_1 s^j_4 s^j_3 
s^k_2 s^l_3 s^l_2 s^k_3 s^i_4 s^i_1 \\ \cdot s^g_2 
s^g_3 s^f_4 s^e_3 s^e_4 s^f_3 s^h_2 s^h_1 
s^f_2 s^e_3 s^e_2 s^f_3 s^g_4 s^g_3 s^i_2 
s^i_3 s^k_4 s^l_3 s^l_4 s^k_3 s^j_2 s^j_3 
s^h_4 s^h_3 s^f_2 s^e_3 s^e_2 s^f_3 s^c_1 
s^c_2 s^d_4 s^d_3 s^c_4 s^c_2 s^f_1 s^h_2 
s^h_3 s^j_4 s^j_3 s^k_2 s^l_3 s^l_2 s^k_3 
s^i_4 s^i_1 s^g_2 s^g_3 s^f_4 s^e_3 s^e_4 \\ \cdot 
s^f_3 s^h_2 s^h_3 s^j_4 s^j_3 s^k_2 s^l_3 
s^l_2 s^k_3 s^i_4 s^i_3 s^g_2 s^g_3 s^f_4 
s^e_3 s^e_4 s^f_3 s^h_2 s^h_1 s^j_4 s^j_3 
s^k_2 s^l_3 s^l_2 s^k_3 s^i_4 s^i_1 s^g_2 
s^g_3 s^f_4 s^e_3 s^e_4 s^f_3 s^h_2 s^h_1 
s^j_4 s^j_3 s^k_2 s^l_3 s^o_1 s^p_4 s^p_3 
s^o_4 s^l_1 s^l_4 s^k_3 s^j_2 s^j_3 s^h_4 
s^h_3 \\ \cdot s^f_2 s^e_3 s^e_2 s^f_3 s^c_1 s^c_2 
s^d_4 s^d_3 s^c_4 s^c_2 s^f_1 s^h_2 s^h_1 
s^j_4 s^j_3 s^k_2 s^l_3 s^o_1 s^p_4 s^p_3 
s^o_4 s^o_2 s^p_4 s^p_3 s^o_4 s^o_2 s^p_4 
s^p_3 s^o_4 s^l_1 s^l_4 s^k_3 s^n_1 s^n_4 
s^m_2 s^m_3 s^n_2 s^n_4 s^k_1 s^i_4 s^i_1 
s^g_2 s^g_3 s^f_4 s^e_3 s^b_1 s^a_2 s^a_3 
s^b_2 s^b_4 \\ \cdot s^a_2 s^a_3 s^b_2 s^b_4 s^a_2 
s^a_3 s^b_2 s^e_1 s^e_2 s^f_3 s^c_1 s^c_2 
s^d_4 s^d_3 s^c_4 s^c_2 s^f_1 s^h_2 s^h_1 
s^j_4 s^j_3 s^k_2 s^l_3 s^o_1 s^p_4 s^p_3 
s^o_4 s^o_2 s^p_4 s^p_3 s^o_4 s^o_2 s^p_4 
s^p_3 s^o_4 s^l_1 s^l_4 s^k_3 s^n_1 s^n_4 
s^m_2 s^m_3 s^n_2 s^n_4 s^k_1 s^i_4 s^i_1 
s^g_2 s^g_3 s^f_4 \\ \cdot s^e_3 s^e_4 s^b_1 s^a_2
\end{array}
\endgroup$}}}
\\ \quad
\\ \mbox{$\FigEightyOneA$: Hamilton circuit of \cite[Table~3, Row~22]{Hec09}
missing $s^a_1$, $s^d_1$, $s^m_1$ and $s^p_1$,} 
\\ \mbox{where $S^\prime$ is that of $\FigEightyOneA^\prime$, and let}
\\ \mbox{$e:={\bar{b}}$, $f:={\bar{a}}$, $g:={\bar{d}}$, $h:={\bar{c}}$, $i:={\bar{c}}^\prime$,
$j:={\bar{d}}^\prime$, $k:={\bar{a}}^\prime$, $l:={\bar{b}}^\prime$,
$m:=d^\prime$,} 
\\ \mbox{$n:=c^\prime$, $o:=b^\prime$, $p:=a^\prime$
(Length $=2304$)}
\end{array}
\end{equation*}

\newcommand{\chamhatN}{{\hat{N}}}
\newcommand{\chamhatPm}{{{\hat{P}}(\chamhatN,m)}}
\newcommand{\chamhatShatNp}{{\mathfrak{S}}(\chamhatN+1)}
\newcommand{\chamhatsigma}{\sigma}
\newcommand{\chamhatp}{{\hat{p}}}
\newcommand{\chamhatpzero}{\chamhatp_0}
\newcommand{\chamhatfkI}{{\hat{\chamfkI}}}
\newcommand{\chamcheckfkI}{{\check{\chamfkI}}}
\newcommand{\chamlambda}{\lambda}
\newcommand{\chamlambdahatp}{{\chamlambda_\chamhatp}}
\newcommand{\chammu}{\mu}
\newcommand{\chammuhatp}{{\chammu_\chamhatp}}
\newcommand{\chamepsilon}{\epsilon}
\newcommand{\chamcalH}{{\mathcal{H}}}
\newcommand{\chamchihatp}{\chamchi^\chamhatp}
\newcommand{\chamOmega}{\Omega}
\newcommand{\chamhattau}{{\hat{\chamtau}}}

\section{Existence of a Hamilton circuit of $\chamGammachi$
with $|\chamfkI|\geq 5$}\label{section:HCosp}

In this section, analyzing carefully in each case,
we shall see:

\begin{theorem} \label{theorem:existGeqRkFive}
Assume $|\chamfkI|\geq 5$.
For every $\chamchi$ of {\rm{\cite[Table~4]{Hec09}}},
there exists a Hamilton circuit map of 
$\chamGammachi$.
\end{theorem}

\noindent
{\bf{Rank~$\geq 5$-Case-1:}} 
Let $\chamchi$ be one of \cite[Table~4, Rows~1-8, 16, 20, 22]{Hec09}.
Then we can easily see that $\chamchi$ is of Cartan-type or quasi-Cartan-type.
By Theorem~\ref{theorem:Main:CSW},
a Hamilton circuit of $\chamGammachi$ exists.

\begin{equation*}
\begin{array}{l}
$\setlength{\unitlength}{1mm}
\begin{picture}(140,115)(20,0)
\put(00,0){$\allHecHigherquasiCartan$}
\end{picture}$
\\ \quad
\\ \mbox{$\FigEightyTwoA$: One of diagrams of each rows of \cite[Table~4,\,Rows~1-8]{Hec09}}
\end{array}
\end{equation*}

\noindent
{\bf{Rank~$\geq 5$-Case-2:}} Let $\chamq\in\bKt\setminus\{1,-1\}$,
$\chamhatN\in\fkJ_{5,\infty}$, $\chamfkI:=\fkJ_{1,\chamhatN}$
and $m\in\fkJ_{1,\chamhatN-1}$.

Let $\chamcalH$ be a rank-$\chamhatN$
free $\bZ$-module
with a basis $\{\chamepsilon_i|i\in\chamfkI\}$.
Let
$\chamhatfkI:=\fkJ_{1,\chamhatN+1}$. 
Let $\chamhatPm$ be the set formed by all maps 
$\chamhatp:\chamhatfkI\to\fkJ_{0,1}$
with $(-1)^{(\chamhatp(\chamhatN)+1)\chamhatp(\chamhatN+1)}=1$
and $|\{i\in\chamfkI|\chamhatp(i)=0\}|=m$.
For $\chamhatp\in\chamhatPm$, define the $\bZ$-bihomomophism
$\chamlambdahatp:\chamcalH\times\chamcalH\to\bZ$
by $\chamlambdahatp(\chamepsilon_i,\chamepsilon_j):=
\delta_{ij}(-1)^{\chamhatp(i)}$
$(i,j\in\chamfkI)$. For $\chamhatp\in\chamhatPm$
and $i\in\chamfkI$, define $\al^\chamhatp_i\in\chamcalH$ by
\begin{equation*}
\al^\chamhatp_i:=
\left\{\begin{array}{ll}
\chamepsilon_i-\chamepsilon_{i+1} & \mbox{if $i\in\fkJ_{1,\chamhatN-2}$,} \\
\chamepsilon_{\chamhatN-1}-\chamepsilon_\chamhatN & 
\mbox{if $\chamhatp(\chamhatN+1)=0$ and $i=\chamhatN-1$,} \\
2\chamepsilon_\chamhatN & 
\mbox{if $\chamhatp(\chamhatN)=\chamhatp(\chamhatN+1)=1$ and $i=\chamhatN-1$,} \\
\chamepsilon_{\chamhatN-1}+\chamepsilon_\chamhatN & 
\mbox{if $\chamhatp(\chamhatN)=\chamhatp(\chamhatN+1)=0$ and $i=\chamhatN$,} \\
2\chamepsilon_\chamhatN & 
\mbox{if $\chamhatp(\chamhatN)-1=\chamhatp(\chamhatN+1)=0$ and $i=\chamhatN$,} \\
\chamepsilon_{\chamhatN-1}-\chamepsilon_\chamhatN & 
\mbox{if $\chamhatp(\chamhatN)=\chamhatp(\chamhatN+1)=1$ and $i=\chamhatN$.}
\end{array}\right.
\end{equation*}
For $\chamhatp\in\chamhatPm$, define the map
$\chammuhatp:\chamfkI\to\fkJ_{0,1}$ by
$\chammuhatp(i):=\delta_{0,t}$, where $t:=\chamlambdahatp(\al^\chamhatp_i,\al^\chamhatp_i)=0$.
For $\chamhatp\in\chamhatPm$, 
define $\chamchihatp\in\chamcalXPi$
by 
\begin{equation*}
\chamchihatp(\al_i,\al_j):=
\left\{\begin{array}{ll}
(-1)^{\chammuhatp(i)\chammuhatp(j)}
\chamq^{\chamlambdahatp(\al^\chamhatp_i,\al^\chamhatp_j)}
& \mbox{if $i\leq j$}, \\
1 & \mbox{if $i>j$}, \\
\end{array}\right.
\end{equation*} for $i$, $j\in\chamfkI$.
Then we see $|\chamRp(\chamchihatp)|<\infty$ by \cite[Table~4,\,Rows~9,10]{Hec09}. 
For ${\grave{\chamchi}}_1$, ${\grave{\chamchi}}_2\in\chamcalXPi$,
we write ${\grave{\chamchi}}_1{\dot{\equiv}}{\grave{\chamchi}}_2$
if ${\grave{\chamchi}}_1(\al_i,\al_i)={\grave{\chamchi}}_2(\al_i,\al_i)$
and ${\grave{\chamchi}}_1(\al_i,\al_j){\grave{\chamchi}}_1(\al_j,\al_i)
={\grave{\chamchi}}_2(\al_i,\al_j){\grave{\chamchi}}_2(\al_j,\al_i)$
for all $i$, $j\in\chamfkI$.
It follows from \cite[Lemma~4.22]{AYY15} that 
if ${\grave{\chamchi}}_1{\dot{\equiv}}{\grave{\chamchi}}_2$
and $|\chamRp({\grave{\chamchi}}_1)|<\infty$, 
then $\chamRp({\grave{\chamchi}}_2)=\chamRp({\grave{\chamchi}}_1)$.
By \cite[Table~4,\,Rows~9,10]{Hec09}, we have
\begin{equation*}
\begin{array}{l}
\forall\chamhatp\in\chamhatPm,\,\,
\{{\grave{\chamchi}}_1\in\chamcalXPi\,|\,
\exists{\grave{\chamchi}}_2\in\chamcalG(\chamchihatp),\,
{\grave{\chamchi}}_1{\dot{\equiv}}{\grave{\chamchi}}_2\} \\
\quad\quad\quad\quad\quad\quad\quad\quad
=\{{\grave{\chamchi}}_3\in\chamcalXPi\,|\,\exists\chamhatp^\prime\in\chamhatPm,\,
{\grave{\chamchi}}_3{\dot{\equiv}}\chamchi^{\chamhatp^\prime}\}.
\end{array}
\end{equation*}
(For every
${\tilde{\chamchi}}\in\chamcalXPi$ of \cite[Table~4,\,Rows~9,10]{Hec09},
we have 
${\tilde{\chamchi}}{\dot{\equiv}}\chamchihatp$
for some 
$\chamq\in\bKt\setminus\{1,-1\}$,
$\chamhatN\in\fkJ_{5,\infty}$,
$m\in\fkJ_{1,\chamhatN-1}$
and $\chamhatp\in\chamhatPm$. See $\FigEightyThreeA$.)

For $\chamhatp\in\chamhatPm$ and $i$, $j\in\chamfkI$ with $i\ne j$, we have
\begin{equation}\label{eqn:chamNchihatpij}
\chamN^{\chamchihatp}_{i,j}=
\left\{\begin{array}{ll}
0 & \mbox{if $\chamlambdahatp(\al^\chamhatp_i,\al^\chamhatp_j)=0$}, \\
1 & \mbox{if $\chamlambdahatp(\al^\chamhatp_i,\al^\chamhatp_i)
=-2\chamlambdahatp(\al^\chamhatp_i,\al^\chamhatp_j)\ne 0$}, \\
1 & \mbox{if $\chamlambdahatp(\al^\chamhatp_i,\al^\chamhatp_i)
=0$ and $\chamlambdahatp(\al^\chamhatp_i,\al^\chamhatp_j)\ne 0$}, \\
2 & \mbox{if $\chamlambdahatp(\al^\chamhatp_i,\al^\chamhatp_i)
=-\chamlambdahatp(\al^\chamhatp_i,\al^\chamhatp_j)\ne 0$}.
\end{array}\right.
\end{equation}

Define $\chamhatpzero\in\chamhatPm$ by
$\chamhatpzero^{-1}(\{0\})=\fkJ_{1,m}\cup\{\chamhatN+1\}$ and
$\chamhatpzero^{-1}(\{1\})=\fkJ_{m+1,\chamhatN}$.
Let $\chamchi:=\chamchi^{\chamhatpzero}$.
We can define the bijection
$\chamOmega:\chamhatPm\to\chambarcalGchi$
in a way that for every $\chamhatp\in\chamhatPm$,
there exists a unique $\chamOmega(\chamhatp)\in\chambarcalGchi$
with $\chamchi^\prime{\dot{\equiv}}\chamchihatp$
for $\chamchi^\prime\in\chamcalGchi$ satisfying
$\champichi(\chamchi^\prime)=\chamOmega(\chamhatp)$.
For $i\in\chamfkI$, define the bijection
$\chamhattau_i:\chamhatPm\to\chamhatPm$
by $\chamOmega(\chamhattau_i(\chamhatp))=\chambartau_i(\chamOmega(\chamhatp))$
$(\chamhatp\in\chamhatPm)$.

Let $\chamhatShatNp$ be the symmetric group of degree $\chamhatN+1$,
and let $\chamhatsigma_i:=(i,i+1)\in\chamhatShatNp$
$(i\in\fkJ_{1,\chamhatN})$.
Regarding $\chamhattau_i(\chamhatp)$ for $\chamhatp\in\chamhatPm$ and $i\in\chamfkI$,
by \eqref{eqn:chamNchihatpij}, 
we have the following.
\newline\newline
$\bullet$ If $\chamlambdahatp(\al^\chamhatp_i,\al^\chamhatp_i)\ne 0$,
then $\chamhattau_i(\chamhatp)=\chamhatp$.
\newline
$\bullet$ If $i\in\fkJ_{1,\chamhatN-2}$, then $\chamhattau_i(\chamhatp)
=\chamhatp\circ\chamhatsigma_i$.
\newline
$\bullet$ If $\chamlambdahatp(\al^\chamhatp_{\chamhatN-1},
\al^\chamhatp_{\chamhatN-1})=0$
(whence $\chamhatp(\chamhatN+1)=0$), then
$\chamhattau_i(\chamhatp)
=\chamhatp\circ\chamhatsigma_{\chamhatN-1}$.
\newline
$\bullet$ If $\chamlambdahatp(\al^\chamhatp_\chamhatN,
\al^\chamhatp_\chamhatN)=0$, then
$\chamhattau_i(\chamhatp)_{|\chamfkI}
=(\chamhatp)_{{|\chamfkI}}\circ\chamhatsigma_{\chamhatN-1}$
and $(\chamhattau_i(\chamhatp))(\chamhatN+1)=1-\chamhatp(\chamhatN+1)$.
Hence by \eqref{eqn:chamNchihatpij},  for $\chamhatp\in\chamhatPm$ and $i$, $j\in\chamfkI$ with $i\ne j$,
letting $a:=\chamOmega(\chamhatp)$, we have
\begin{equation*}
m^a_{ij}=
\left\{\begin{array}{ll}
0 & \mbox{if $\chamlambdahatp(\al^\chamhatp_i,\al^\chamhatp_j)=0$}, \\
3 & \mbox{if $\chamlambdahatp(\al^\chamhatp_i,\al^\chamhatp_i)
=-2\chamlambdahatp(\al^\chamhatp_i,\al^\chamhatp_j)\ne 0$}, \\
3 & \mbox{if $\chamlambdahatp(\al^\chamhatp_i,\al^\chamhatp_i)
=0$ and $\chamlambdahatp(\al^\chamhatp_i,\al^\chamhatp_j)\ne 0$}, \\
4 & \mbox{if $\chamlambdahatp(\al^\chamhatp_i,\al^\chamhatp_i)
=-\chamlambdahatp(\al^\chamhatp_i,\al^\chamhatp_j)\ne 0$}.
\end{array}\right.
\end{equation*} 
Notice $\chamhatN\geq 5$.
By $\FigEightyThreeA$, $\FigEightyFourA$ and $\FigEightyFiveA$,
using the claim obtained from that of Proposition~\ref{proposition:EdgeSpecial}
by letting $N=\chamhatN$, $i=N-2$
and then changing $1,2,\ldots,\chamhatN$ by
$\chamhatN,\chamhatN-1,\ldots,1$, 
we see that a Hamilton circuit map of $\chamGamma(\chamchihatp)$ exists.
Similarly, by Proposition~\ref{proposition:EdgeSpecial}
for $N=\chamhatN$ and $i=N-3$,
using $\FigThreeA$, $\FigSixA$, $\FigTwentySevenA$ and $\FigThirtyTwoA$, 
we also see that fact.

\begin{equation*}
\begin{array}{l}
$\setlength{\unitlength}{1mm}
\begin{picture}(140,75)(20,0)
\put(00,0){$\allHecHigherNineTen$}
\end{picture}$
\\ \quad
\\ \mbox{$\FigEightyThreeA$: Diagrams of \cite[Table~4,\,Rows~9,10]{Hec09},
where $i\in\fkJ_{1,\chamhatN-2}$, 
$q\in\bKt$, $q^2\ne 1$,}
\\ \mbox{$\chamq_j\in\{-1,\chamq^{\pm 1}\}$
$(j\in\fkJ_{1,\chamhatN-1})$, $r_k\in\{\chamq^{\pm 1}\}$ 
$(k\in\fkJ_{0,\chamhatN-1})$, $r_{i-1}\chamq_i^2 r_i=1$}
\\ \mbox{and $m=|\{i\in\fkJ_{0,\chamhatN}|r_i=\chamq^{-1}\}|$.}
\end{array}
\end{equation*}

\begin{equation*}
\begin{array}{l}
$\setlength{\unitlength}{1mm}
\begin{picture}(90,35)(00,00)
\put(00,0){$\cartanAone$}
\end{picture}$
\\ \quad
\\ \mbox{$\FigEightyFourA$: Special Hamilton circuit of \cite[Table~1,\,Row~1]{Hec09}}
\end{array}
\end{equation*}

\begin{equation*}
\begin{array}{l}
$\setlength{\unitlength}{1mm}
\begin{picture}(90,35)(20,00)
\put(00,0){$\supersltwoone$}
\end{picture}$
\\ \quad
\\ \mbox{$\FigEightyFiveA$: Special Hamilton circuit of \cite[Table~1,\,Row~2]{Hec09}}
\end{array}
\end{equation*}

We emphasis:
\begin{theorem}\label{theorem:HChigherranknineten}
For every $\chamchi$ of 
{\rm{\cite[Table~4,\,Rows~9,10]{Hec09}}},
there exists a Hamilton circuit map of $\chamGammachi$.
\end{theorem}


\noindent
{\bf{Rank~$\geq 5$-Case-3:}} Let $\chamchi$ be of \cite[Table~4, Rows~11-15, 17-19, 21]{Hec09}.
Let $N:=|\chamPi|$. 
Here $N\in\fkJ_{5,7}$ and $\chamfkI=\fkJ_{1,N}$.
Notice Tables~$\tablethree$ and $\tablefour$,
where we show a Hamilton circuit of 
$\chamGamma({\hat{\chamchi}}_{|\bZ\al_{N-3}\oplus
\bZ\al_{N-2}\oplus\bZ\al_{N-1}\oplus\bZ\al_N})$
of ${\hat{\chamchi}}\in\chamcalGchi$.
Recall that 
we have a special Hamilton circuit map of $\chamGamma(\chamchi^\prime)$ for $\chamchi^\prime$ of \cite[Table~3, Rows~11,13,18]{Hec09}
by $\FigSixtyOneA$, $\FigFiftyTwoA$ and $\FigSeventyFourA$.
Hence by Proposition~\ref{proposition:EdgeSpecial}, we have 
a Hamilton circuit map of $\chamGammachi$ 
for $\chamchi$ of \cite[Table~4, Rows~11, 12, 15, 17, 18, 21]{Hec09}. 
Although the Hamilton circuit maps $\FigSixtyFourA$ and $\FigFiftyFiveA$ of 
of $\chamGamma(\chamchi^{\prime\prime})$ for $\chamchi^{\prime\prime}$ of of \cite[Table~3, Rows~14, 20]{Hec09}
are not special, it follows from Proposition~\ref{proposition:EdgeAlmostSpecial} that
we have 
a Hamilton circuit map of $\chamGammachi$ 
for $\chamchi$ of \cite[Table~4, Rows~13, 14, 19]{Hec09}. 
Incidentally, if $\chamhatN\geq 6$, then
by $\FigEightyFourA$ and $\FigEightyFiveA$,
using the claim obtained from that of Proposition~\ref{proposition:EdgeSpecial}
by letting $N=\chamhatN$, $i=N-2$
and then changing $1,2,\ldots,\chamhatN$ by
$\chamhatN,\chamhatN-1,\ldots,1$, 
we see that a Hamilton circuit map of $\chamGammachi$ exists.

\begin{equation*}
\begin{array}{l}
$\TopHecFiveEightQUCX$
\\ \quad
\\ \mbox{$\FigEightySixA$: First diagrams of \cite[Table~4,\,Rows 11-22]{Hec09}}
\end{array}
\end{equation*}

\begin{center}
\begingroup
\renewcommand{\arraystretch}{1.5}
Table~$\tablethree$: Rank-5's rank-4 edge sub-Hamilton circuits \\
  \begin{tabular}{lccccc}
    \hline
$\chamchi$  &   
$\acchi^{(4,11)}_\bullet$ &
$\acchi^{(4,12)}_\bullet$ &
$\acchi^{(4,13)}_\bullet$ &
$\acchi^{(4,14)}_\bullet$ &
$\acchi^{(4,15)}_\bullet$
\\
    \hline \hline
\mbox{Fig.s} & 
\NFigEightA, 
\NFigFiftyTwoA 
&
\NFigFiftyTwoA, 
\NFigSeventyFourA 
&
\NFigEightA, 
\NFigSeventyFourA, 
\NFigFiftyFiveA 
&
\NFigEightA, 
\NFigSixtyFourA 
&
\NFigEightA, 
\NFigSixtyOneA 
\\
\hline
  \end{tabular}
\endgroup
\end{center}

\vspace{0.1cm}

\begin{center}
\begingroup
\renewcommand{\arraystretch}{1.5}
Table~$\tablefour$: Rank-6, 7 and 8's rank-4 edge sub-Hamilton circuits \\
  \begin{tabular}{lccccccc}
    \hline
$\chamchi$  & 
$\acchi^{(4,16)}_\bullet$ &  
$\acchi^{(4,17)}_\bullet$ &
$\acchi^{(4,18)}_\bullet$ &
$\acchi^{(4,19)}_\bullet$ &
$\acchi^{(4,20)}_\bullet$ &
$\acchi^{(4,21)}_\bullet$ &
$\acchi^{(4,22)}_\bullet$
\\
    \hline \hline
\mbox{Fig.s} & 
\NFigEightA &
\NFigEightA, 
\NFigFiftyTwoA 
&
\NFigEightA, 
\NFigFiftyTwoA, 
\NFigSeventyFourA 
&
\NFigEightA, 
\NFigSixtyFourA 
&
\NFigEightA
&
\NFigEightA, 
\NFigFiftyTwoA 
&
\NFigEightA
\\
\hline
  \end{tabular}
\endgroup
\end{center}

\noindent
To obtain the above 
Tables~$\tablethree$ and $\tablefour$,
we can use some diagrams in
\cite{AA17} together with Section~\ref{section:NamesAA17}
below.

After all, by Theorems~\ref{theorem:existRkThree}, \ref{theorem:existRkFour}
and \ref{theorem:existGeqRkFive}, we prove our main theorem:
\begin{theorem}\label{theorem:HChigherrank}
For every $\chamchi\in\chamcalXPi$ with $|\chamRpchi|<\infty$,
there exists a Hamilton circuit map of $\chamGammachi$.
\end{theorem}

\section{Names by \cite{AA17}} \label{section:NamesAA17}
As mentioned in Definition~\ref{definition:genDynkindiagram},
$\acchi^{(x,y)}_\bullet$ means the family of the generalized Dynkin diagrams
belonging to Row-$y$ of Table~$x$ of \cite{Hec09}.
In \cite{AA17}}, a suitable name of each $\acchi^{(x,y)}_\bullet$ has been introduced
with consideration of not only simple Lie superalgerbas
but also modular simple Lie superalgerbas.
The names are RHSs of the following.
\newline\newline
(2) See \cite[Section~7]{AA17}.
\newline\newline
$\acchi^{(3,14)}_\bullet={\mathbf{{\tt{wk}}(4)}}$,
$\acchi^{(1,5)}_\bullet={\mathbf{{\tt{br}}(2)}}$,
$\acchi^{(2,18)}_\bullet={\mathbf{{\tt{br}}(3)}}$.
\newline\newline
(2) See \cite[Section~8]{AA17}.
\newline\newline
$\acchi^{(1,9)}_\bullet={\mathbf{{\tt{brj}}(2;3)}}$,
$\acchi^{(2,13)}_\bullet={\mathbf{g(1,6)}}$,
$\acchi^{(2,15)}_\bullet={\mathbf{g(2,3)}}$,
$\acchi^{(3,18)}_\bullet={\mathbf{g(3,3)}}$,
$\acchi^{(3,20)}_\bullet={\mathbf{g(4,3)}}$,
$\acchi^{(3,21)}_\bullet={\mathbf{g(3,6)}}$,
$\acchi^{(4,11)}_\bullet={\mathbf{g(2,6)}}$,
$\acchi^{(4,12)}_\bullet={\mathbf{{\tt{el}}(5;3)}}$,
$\acchi^{(4,13)}_\bullet={\mathbf{g(8,3)}}$,
$\acchi^{(4,17)}_\bullet={\mathbf{g(4,6)}}$,
$\acchi^{(4,18)}_\bullet={\mathbf{g(6,6)}}$,
$\acchi^{(4,21)}_\bullet={\mathbf{g(8,6)}}$.
\newline\newline
(3) See \cite[Section~9]{AA17}.
\newline\newline
$\acchi^{(1,13)}_\bullet={\mathbf{{\tt{brj}}(2;5)}}$,
$\acchi^{(4,15)}_\bullet={\mathbf{{\tt{el}}(5;5)}}$.
\newline\newline
(3) See \cite[Section~10]{AA17}.
\newline\newline
$\acchi^{(4,14)}_\bullet={\mathbf{{\tt{ufo}}(1)}}$,
$\acchi^{(4,19)}_\bullet={\mathbf{{\tt{ufo}}(2)}}$,
$\acchi^{(2,16)}_\bullet={\mathbf{{\tt{ufo}}(3)}}$,
$\acchi^{(2,17)}_\bullet={\mathbf{{\tt{ufo}}(4)}}$,
$\acchi^{(3,17)}_\bullet={\mathbf{{\tt{ufo}}(5)}}$,
$\acchi^{(3,22)}_\bullet={\mathbf{{\tt{ufo}}(6)}}$,
$\acchi^{(1,7)}_\bullet={\mathbf{{\tt{ufo}}(7)}}$,
$\acchi^{(1,8)}_\bullet={\mathbf{{\tt{ufo}}(8)}}$,
$\acchi^{(1,12)}_\bullet={\mathbf{{\tt{ufo}}(9)}}$,
$\acchi^{(1,14)}_\bullet={\mathbf{{\tt{ufo}}(10)}}$,
$\acchi^{(1,15)}_\bullet={\mathbf{{\tt{ufo}}(11)}}$,
$\acchi^{(1,16)}_\bullet={\mathbf{{\tt{ufo}}(12)}}$.

\hspace{3cm}

\noindent
{\bf{Acknowledgment.}} 
We thank Iv{\`{a}}n Angiono and Masaya Tomie for helpful suggestions.
For the achievement of this paper,
especially rank~4 cases,  the author used 
Wolfram Mathematica~12.1 \cite{Mathematica20}
with his own programs assembled by only easy commands.
This work was partially
supported by JSPS Grant-in-Aid for Scientific Research (C) 19K03420.

\noindent
Hiroyuki Yamane:
Department of Mathematics,
Faculty of Science, Academic Assembly,
University of Toyama,
3190 Gofuku, Toyama-shi, Toyama 930-8555, Japan \newline
E-mail: hiroyuki@sci.u-toyama.ac.jp

\end{document}